\documentclass[11pt]{amsart}                                                         

\usepackage{geometry}

\usepackage{mathtools}	%
\mathtoolsset{showonlyrefs}	%

\usepackage{tikz-cd}
\usetikzlibrary{arrows,
	backgrounds,
	decorations,
	decorations.pathmorphing,
	decorations.pathreplacing,
	decorations.shapes,
	positioning,
	fit,
	petri,
	graphs,
	matrix,
	automata,
	intersections,
	shapes,
	calc,
	tikzmark}
\tikzcdset{arrow style=tikz}
\usepackage{amssymb}
\usepackage{bbm}			%
\usepackage{dsfont}
\usepackage{stmaryrd}		%

\usepackage{afterpage}	%
\usepackage{array}
\usepackage{chngcntr}	%
\usepackage[inline]{enumitem}	%
\usepackage{leftindex}	%
\usepackage{longtable}	%
\usepackage{nicematrix}	%
\usepackage{scalerel}   %
\usepackage{subcaption}	%
\usepackage{xcolor}		%

\usepackage{multirow}	%

\usepackage[unicode,pagebackref=false]{hyperref}	%
\hypersetup{
	bookmarksnumbered=true,
	colorlinks=true,
	linkcolor=blue,
	citecolor=blue,
	filecolor=blue,
	menucolor=blue,
	urlcolor=blue,
	bookmarksopen=true,
	bookmarksdepth=2,
	pageanchor=true}

\numberwithin{equation}{section}

\newcounter{resultcounter}

\theoremstyle{plain}
\newtheorem{result}[resultcounter]{Result}
\newtheorem{thm}[equation]{Theorem}
\newtheorem{lem}[equation]{Lemma}
\newtheorem{prp}[equation]{Proposition}
\newtheorem{cor}[equation]{Corollary}
\newtheorem{dfn}[equation]{Definition}
\newtheorem{ex}[equation]{Example}
\newtheorem{rmk}[equation]{Remark}

\makeatletter
\renewcommand\paragraph{\@startsection{paragraph}{4}{\z@}%
{3.25ex \@plus1ex \@minus.2ex}%
{-1em}%
{\normalfont\itshape}}
\makeatother

\newenvironment{td}[0]		%
{\begin{tikzcd}[ampersand replacement=\&, cells={outer sep=2pt, inner sep=0pt}, 
	]}{\end{tikzcd}}

\newenvironment{cd}[0]		%
{\begin{tikzcd}[ampersand replacement=\&, cells={outer sep=2pt, inner sep=1pt}, 
	row sep=large,
	column sep=large,	
	]}{\end{tikzcd}}

\newenvironment{arq}[0]		%
{\begin{tikzcd}[ampersand replacement=\&, cells={outer sep=2pt, inner sep=2pt}, 
	column sep=0.5cm,
	every cell/.style={font=\small, color=gray},
	arrows={>=stealth, color=gray}
	]}{\end{tikzcd}}

\newenvironment{gd}[0]		%
{\begin{tikzcd}[ampersand replacement=\&, 
	cells={outer sep=0.75pt, inner sep=0.75pt}, 
	row sep={2.5em}, 
	column sep={35pt},
	labels={rectangle, font=\small, inner sep=1.25pt, outer sep=1.5pt, minimum size=0cm},
	baseline=0pt,
	arrows={->,semithick, >=stealth'},
	]}{\end{tikzcd}}

\newenvironment{qrep}[0]		%
{\begin{tikzcd}[ampersand replacement=\&,
	cells={rectangle, inner sep=0.5pt, outer sep=0.5pt,  minimum width=0.0cm, 
		minimum height=0.5cm, 
		text height=0.25cm
	},
	labels={rectangle, font=\small, outer sep=1pt
	}, 	
	arrows={->,semithick, >=stealth'},
	/tikz/bend right/.default=30,
	/tikz/bend left/.default=30,
	column sep=1cm,
	]}{\end{tikzcd}}
\newenvironment{ld}[0]		%
{\begin{tikzcd}[ampersand replacement=\&,
	cells={rectangle, inner sep=0pt, outer sep=1pt,  
	},
	labels={rectangle, font=\small, outer sep=0pt, inner sep=0.5pt
	}, 	
	arrows={->,semithick, >=stealth'},
	/tikz/bend right/.default=30,
	/tikz/bend left/.default=30,
	column sep=1cm,
	]}{\end{tikzcd}}
\newenvironment{gqrep}[7]		%
{\begin{qrep}
	{#3}
	\ar[start anchor=north east, end anchor=north west, bend left,yshift=-5pt]{r}{#5}
	\&
	{#2}   
	\ar[start anchor=south west, end anchor=south east,bend left,yshift=5pt]{l}{#7}
	\ar[start anchor=south east, end anchor=south west, bend right, yshift=5pt]{r}[swap]{#6}
	\& 
	{#1} 
	\ar[start anchor=north west, end anchor=north east, bend right, yshift=-5pt]{l}[swap]{#4}
}{\end{qrep}}

\definecolor{c1}{RGB}{255, 102, 102} %
\definecolor{c2}{RGB}{102, 178, 255} %
\colorlet{c3}{c1!50!c2!50!white}
\colorlet{c4}{green!66!black} %
\definecolor{c5}{RGB}{102, 178, 139} %
\definecolor{green}{rgb}{0.0,0.5,0.0}
\definecolor{red}{rgb}{0.5,0.0,0.0}

\makeatletter
\AtBeginDocument{
\let\c@figure\c@equation
\let\c@table\c@equation
}
\makeatother
\counterwithin{figure}{section}
\counterwithin{table}{section}

\newcommand{\hcinvaut}{\varrho}
\newcommand{\hcinv}{\rho}
\newcommand{\invaut}{\varsigma}
\newcommand{\inv}{\sigma}

\DeclareMathOperator{\diag}{\mathsf{diag}}
\newcommand{\id}{\mathsf{id}}
\newcommand{\JB}{\mathsf{J}}
\newcommand{\Id}{\mathbbm{1}} 

\newcommand{\Ccat}{\mathcal{C}} 
\newcommand{\Dcat}{\mathcal{D}}

\DeclareMathOperator{\add}{\mathsf{add}}
\DeclareMathOperator{\ind}{\mathsf{ind}}
\newcommand{\HC}{\mathsf{HC}}

\DeclareMathOperator{\md}{\mathsf{mod}}
\DeclareMathOperator{\Md}{\mathsf{Mod}}
\DeclareMathOperator{\proj}{\mathsf{proj}}

\DeclareMathOperator{\rep}{\mathsf{rep}}
\newcommand{\Rep}[1]{\rep #1}
\newcommand{\regrep}[1]{\Rep_{\cong}(#1)}
\newcommand{\stmd}[1]{\underline{\md} #1}

\newcommand{\Acat}{\Rep{\A}} 

\newcommand{\bd}{\mathsf{b}}
\newcommand{\fd}{\mathsf{fd}}
\newcommand{\bdmin}{\mathsf{b,m}}

\newcommand{\Com}{\mathsf{Com}}
\newcommand{\Comb}{\Com^{\bd}}
\newcommand{\Combmin}{\Com^{\bdmin}}
\newcommand{\Combminproj}[1]{\Combmin(\proj{#1})}
\newcommand{\Combproj}[1]{\Comb(\proj{#1})}
\newcommand{\Combzeroproj}[1]{\bigoplus_{i \in \Z} \proj{#1}}

\newcommand{\minus}{\mathsf{-}}
\newcommand{\D}{\mathsf{D}}
\newcommand{\Db}[1]{\D^{\bd}(#1)}
\newcommand{\Dbfd}[1]{\D^{\bd}_{\fd}(#1)}
\newcommand{\DbRep}[1]{\Db{\rep{#1}}}
\newcommand{\DbHC}{\D^{\bd}\bigl(\HC_\circ(\lieg, K)\bigr)}
\newcommand{\Dbfdmod}[1]{\Dbfd{\md{#1}}}
\newcommand{\DbfdMod}[1]{\Dbfd{\Md{#1}}}
\newcommand{\Dmmod}[1]{\D^{\minus}(\md{#1})}
\newcommand{\Dbmod}[1]{\Db{\md{#1}}}

\newcommand{\gcat}{\mathsf{gt}}
\newcommand{\Gcat}[1]{\gcat(#1)}
\newcommand{\Gcatfd}[1]{\gcat_{\fd}(#1)}
\newcommand{\Gcatproj}[1]{\gcat_0(#1)}
\newcommand{\GcatCombminproj}[1]{\gcat'(#1)}

\newcommand{\Hot}{\mathsf{Hot}}
\newcommand{\Hotb}{\Hot^{\bd}}
\newcommand{\Hotbfd}{\Hotb_{\fd}}
\newcommand{\Hotbproj}[1]{\Hotb(\proj{#1})}
\newcommand{\Hotbfdproj}[1]{\Hotbfd(\proj{#1})}

\newcommand{\Dtr}[1]{(#1)^t}

\newcommand{\dual}[1]{#1^{\dagger}}
\newcommand{\cdual}[1]{#1^{\ddagger}}
\newcommand{\hcdual}[1]{{#1}^{\vee}}

\newcommand{\EE}{\mathbb{E}}
\newcommand{\FF}{\mathbb{F}}

\renewcommand{\SS}{\mathbb{S}}
\newcommand{\Tw}{\mathbb{T}_{\! S_\star}}

\newcommand{\F}{\mathbb{F}}
 
\newcommand{\Prf}{\mathbb{P}} 
\newcommand{\Gl}{\mathbb{G}} 
\newcommand{\Red}{\mathbb{H}}
\newcommand{\Ind}{\mathbb{I}}
\newcommand{\res}[1]{\leftindex_{\A}{#1}}
\DeclareMathOperator{\syz}{\mathsf{syz}}

\newcommand{\defect}{\delta}
\DeclareMathOperator{\gldim}{\mathsf{gldim}}
\DeclareMathOperator{\prdim}{\mathsf{prdim}}
\DeclareMathOperator{\injdim}{\mathsf{injdim}}

\DeclareMathOperator{\uldim}{\underline{\mathsf{dim}}}

\DeclareMathOperator{\ann}{\mathsf{ann}}
\DeclareMathOperator{\supp}{\mathsf{supp}}
\DeclareMathOperator{\head}{\mathsf{top}}
\DeclareMathOperator{\soc}{\mathsf{soc}}
\DeclareMathOperator{\rad}{\mathsf{rad}}
\renewcommand{\dim}{\mathsf{dim}\,}

\newcommand{\ct}[2]{\beta_{#1#2}}
\newcommand{\ctpr}[2]{\beta^{\mathrm{pr}}_{#1#2}}
\newcommand{\ctdeg}[1]{\beta_{#1}}

\newcommand{\ctalt}[1]{\chi_{#1}}

\newcommand{\C}{C}
\newcommand{\PPr}{P}
\newcommand{\X}{X}
\newcommand{\Y}{Y}
\newcommand{\V}{S}
\newcommand{\W}{T}
\newcommand{\AP}{X}
\newcommand{\BM}{\tilde{M}}
\newcommand{\BP}{\widetilde{P}}
\newcommand{\BS}{\tilde{S}}

\newcommand{\CP}{\AP_{\bu}}
\newcommand{\CPr}{P_{\bu}}
\newcommand{\CL}{L_{\bu}}

\newcommand{\CV}{S_{\bu}}
\newcommand{\CW}{T_{\bu}}
\newcommand{\CX}{X_{\bu}}
\newcommand{\CY}{Y_{\bu}}

\renewcommand{\simeq}{\cong}
\newcommand{\topspan}[1]{\overline{\mathsf{span}}\{{#1}\}}
\newcommand{\colonequals}{\coloneqq}
\newcommand{\equalscolon}{\eqqcolon}

\let\ker\relax 
\DeclareMathOperator{\charac}{\mathsf{char}}
\newcommand{\chr}[1]{\charac{(#1)}}
\DeclareMathOperator{\ker}{\mathsf{ker}}
\DeclareMathOperator{\coker}{\mathsf{coker}}
\DeclareMathOperator{\im}{\mathsf{im}}

\newcommand{\kk}{\mathbbm{k}}
\newcommand{\CC}{\mathbb C}
\newcommand{\N}{\mathbb{N}}

\newcommand{\RR}{\mathbb R}
\newcommand{\Z}{\mathbb{Z}}

\newcommand{\lieg}{\mathfrak{g}}
\newcommand{\idm}{\mathfrak{m}}
\newcommand{\Rx}{R}
\newcommand{\mx}{\mathfrak m}

\newcommand{\A}{\Lambda}
\newcommand{\Aop}{\Lambda^{\op}}
\newcommand{\ACtwo}{\A\Ctwo}
\newcommand{\AI}{\mathsf{A}} 
\newcommand{\B}{\Gamma}
\newcommand{\BI}{\mathsf{B}}
\newcommand{\Ctwo}{\mathsf{C}_2}

\newcommand{\oA}{\mathbin{\otimes_{\A}}}
\newcommand{\oB}{\mathbin{\otimes_{\B}}}
\newcommand{\oAI}{\mathbin{\otimes_{\AI}}}
\newcommand{\oBI}{\mathbin{\otimes_{\BI}}}

\DeclareMathOperator{\GL}{\mathsf{GL}}
\DeclareMathOperator{\Mat}{\mathsf{Mat}}
\DeclareMathOperator{\SL}{\mathsf{SL}}
\DeclareMathOperator{\SO}{\mathsf{SO}}

\DeclareMathOperator{\Aut}{\mathsf{Aut}}
\DeclareMathOperator{\End}{\mathsf{End}}
\DeclareMathOperator{\Ext}{\mathsf{Ext}}
\DeclareMathOperator{\Hom}{\mathsf{Hom}}
\DeclareMathOperator{\RHom}{\mathsf{RHom}}

\newcommand{\fE}{\mathfrak{E}}
\newcommand{\fF}{\mathfrak{F}}
\newcommand{\fX}{\mathfrak{X}}
\newcommand{\fXA}{\mathfrak{X}_{\A}}
\newcommand{\chfXA}{\DbRep{\A}}
\newcommand{\kA}{\mathcal{A}}
\newcommand{\kX}{\mathcal{X}}
\newcommand{\sign}{\mathrm{s}}

\newcommand{\xx}{a}
\newcommand{\yy}{b}
\newcommand{\xe}{e}
\newcommand{\xa}{a_1}
\newcommand{\xb}{a_2}
\newcommand{\xc}{a_3}
\newcommand{\xf}{f}
\newcommand{\ya}{b_1}
\newcommand{\yb}{b_2}
\newcommand{\yc}{b_3}
\newcommand{\yf}{g}

\newcommand{\op}{\mathrm{op}}
\newcommand{\olra}{\overset{\leftrightarrow}}
\newcommand{\ora}[1]{\overset{\rightarrow}{#1}{}}
\newcommand{\ola}[1]{\overset{\leftarrow}{#1}{}}

\def\ol{\overline}

\newcommand{\wh}{\widehat}
\newcommand{\wt}{\widetilde}

\tikzset{rotate/.style={anchor=south, rotate=90, inner sep=.5mm}}
\newcommand{\bt}{\circ}	%
\newcommand{\pullback}{\lrcorner}
\newcommand\bu[1][.8]{\mathbin{\ThisStyle{\vcenter{\hbox{%
				\scalebox{#1}{$\SavedStyle\bullet$}}}}}%
}

\newcommand{\secref}[1]{%
\S\ref{#1}.~\nameref{#1}\\
\hfill{\footnotesize\color{black}p.\,\pageref{#1}}%
}

\begin{document}

\title[Representation theory of the Gelfand quiver]{Representation theory of the Gelfand quiver  \\ and Harish-Chandra modules for $\SL_2(\RR)$}

\author{Igor Burban}
\address{
	Universit\"at Paderborn,
	Institut f\"ur Mathematik,
	Warburger Strasse 100,
	33098 Paderborn,
	Germany
}
\email{burban@math.uni-paderborn.de} 

\author{Wassilij Gnedin}
\address{
	Universit\"at Bielefeld,
	Fakult\"at f\"ur Mathematik,
	Universit\"atsstrasse 25,
	D-33615 Bielefeld,
	Germany
}
\email{wgnedin@math.uni-bielefeld.de}

\begin{abstract}
	In 1970, Gelfand posed the problem of classifying the indecomposable objects in a representation category equivalent to the principal block of Harish-Chandra modules for $\SL_2(\RR)$;
	explicit solutions were obtained by Bondarenko, and, independently, Crawley-Boevey.
	In this article, we give a complete answer to Gelfand's problem from a derived category perspective.
	We classify indecomposable objects in the bounded derived category
	of nilpotent representations of the Gelfand quiver
	in terms of band and string complexes, and determine their images under the derived Auslander-Reiten translation, the sign involution, and the contragredient duality. The four main combinatorial classes are characterized in Lie-theoretic as well as homological terms. For the abelian category of nilpotent representations, we provide 
	projective resolutions, standard homological invariants
	and explicit representation matrices of all indecomposables.
	Our approach can be extended to arrow ideal completions of path algebras of skew-gentle quivers.
\end{abstract}

\maketitle

\section{Introduction}

In his ICM talk in 1970, I.~Gelfand \cite{Gelfand} raised the question of finding
an explicit  description of the indecomposable complex finite-dimensional nilpotent representations of the quiver
\begin{align}
	\label{E:GelfandQuiver}
	\begin{gqrep}{+}{\star}{-}{b_+}{b_-}{a_+}{a_-}\end{gqrep}
	\qquad b_{-} a_{-} = b_{+} a_{+}.
\end{align}
The motivation to study  this problem comes from the fact that the category of such representations of
\eqref{E:GelfandQuiver} is equivalent to the principal block of the category of Harish--Chandra modules 
for the Lie group
$\SL_2(\RR)$.

Already in 1959 Zhelobenko  \cite{Zhelobenko} showed  that the principal block of the category of Harish-Chandra modules
for
$\SL_2(\CC)$ is equivalent to the category of 
complex
fi\-nite-di\-men\-sion\-al
nilpotent representations of the quiver
\begin{align}
	\label{E:ZhelobQuiver}
	\begin{qrep}			
		{\bt}
		\arrow[loop, looseness=7.5, in=150, out=215, "x", ""{name=x1, inner sep=0pt,
			outer sep=0pt, pos=0.175}, ""{name=x2, inner sep=0pt,outer sep=0pt, pos=0.825}]
		\ar[bend left,yshift=0pt, ""{name=b1, inner sep=0, pos=0.33,outer sep=0pt}]{r}{b}
		\&  
		{\bt}   
		\ar[bend left,yshift=0pt, ""{name=a2, inner sep=0pt, outer sep=0pt, pos=0.66}]{l}{a}
	\end{qrep} \qquad bx = xa = 0.
\end{align}
A description of the corresponding indecomposable  representations  was obtained by Gel\-fand and Ponomarev in 1968;  see  \cite{GP}. Their work was very influential for the further development of the representation theory of associative algebras.

However, the case of the quiver \eqref{E:GelfandQuiver} turned out to be more complicated. Quoting \cite{Gelfand}: ``The question of the classification of the objects of this category is apparently solvable but leads to considerable difficulties''.

In 1973 Nazarova and Roiter proved that the category of finite-dimensional nilpotent representations of
\eqref{E:GelfandQuiver} is representation-tame, reducing the problem of description of its indecomposable objects to a certain problem of linear algebra (matrix problem) \cite{NR}. 
However, the correct combinatorics of indecomposables 
were obtained only in the end of the 1980s by
Bondarenko \cite{B1,B2} and Crawley-Boevey \cite{CB}
by substantially different methods.
Further solutions to Gelfand’s problem were given
by Deng \cite{Deng} and Iyama \cite{Iyama}.

In 1981 Khoroshkin obtained a 
quiver description
of the principal blocks of Harish--Chandra modules for the Lorentz groups $\SO^+(n,1)$
proving they  are representation-tame for all $n$ \cite{Khoroshkin}. 
In case $n = 2m$, this quiver is given as follows.
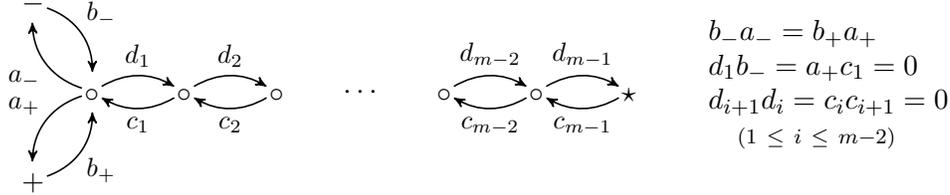
\begin{align}
	\label{eq:HQ}
	\begin{qrep}
		{-}	\ar[bend left=35]{rd}[outer sep = 0pt]{b_-}
		\\
		\&[-0.5cm]	
		{\bt}
		\ar[bend left=35]{lu}[outer sep = 0pt]{a_-}
		\ar[bend left]{r}{d_1}
		\ar[bend right=35]{ld}[swap,outer sep = 0pt]{a_+}
		\&
		{\bt}   
		\ar[bend left]{l}{c_1}
		\ar[bend right, <-]{r}[swap]{c_2}
		\& 
		{\bt}
		\ar[bend right, <-]{l}[swap]{d_2}
		\ar[-, densely dotted, phantom]{rr}[description]{\cdots}
		\&
		\&
		{\bt}
		\ar[bend left]{r}{d_{m-2}}
		\&
		{\bt} 
		\ar[bend left]{l}{c_{m-2}}
		\ar[bend right, <-]{r}[swap]{c_{m-1}}
		\& 
		{\star} 
		\ar[bend right, <-]{l}[swap]{d_{m-1}}	
		\\
		{+}\ar[bend right=35]{ru}[swap,outer sep=0pt]{b_+}
	\end{qrep}
	\qquad
	\begin{array}{l}
		b_- a_- = b_+ a_+
		\\
		d_1 b_- = a_+ c_1 = 0 \\
		d_{i+1} d_{i} =c_{i} c_{i+1} =  0
		\\
		\quad
		{\scriptstyle  (1  \ \leq \   i \ \leq \  m-2)}
	\end{array}
\end{align}
In \cite{Nodal} the first-named author together with Drozd showed that the Khoroshkin quivers are even derived-tame.

Although the indecomposable representations of the Gelfand quiver have been classified via prior approaches, a systematic treatment via the derived category has not been carried out.

In this article, we describe the indecomposable objects in the derived category of nilpotent representations of the Gelfand quiver, and determine their images under the derived Auslander--Reiten translation and contragredient duality.
This approach has two advantages over prior methods:
first, the techniques are amenable to generalizations to a larger class of derived-tame algebras;
second, the results on the derived category yield a conceptual framework 
to study the original, abelian category.
Following this approach, 
we provide a detailed solution to Gelfand's problem, including
projective resolutions, representation matrices, homological invariants and functorial properties 
of all indecomposable nilpotent representations.

\subsection*{Main results}
For the sake of simplicity of the presentation, we fix an algebraically closed field $\kk$ as base field. The category of finite-dimensional nilpotent $\kk$-linear representations of the Gelfand quiver 
will be denoted by $\Acat$, its bounded derived category by $\Db{\Acat}$.
We give approximate descriptions of certain combinatorial notions, called \emph{bands} and \emph{strings of $\Db{\Acat}$}.
\begin{itemize}
	\item
	A \emph{periodic word} $w$ is given by a periodic sequence $(x_i)_{i \in \Z}$, 
	where $x_i$ is a decorated number $\ora{n}_i$ or $\ola{n}_i$ with a positive integer $n_i \in \N$ for each index $i$,
	and a minimal periodic part of $w$ 
	has as many left-oriented as right-oriented numbers.
	Any \emph{band} $\omega$ of $\Db{\Acat}$ has the form $(d,w,m,\lambda)$ consisting of an integer $d\in \Z$,
	a periodic word $w$, a number $m \in \N$ and a scalar $\lambda \in \kk^*$.
	\item
	Any \emph{finite word} $w$ is a sequence $(\alpha,x_1,x_2,\ldots x_{\ell},\beta)$ 
	with ends $\alpha$, $\beta \in \{\star,\diamond \}$, 
	and
	certain 
	decorated numbers
	in between. 
	A finite word is \emph{usual}, \emph{special} or \emph{bispecial} if it has $0$, $1$ or $2$ ends of type $\diamond$, respectively.
	A \emph{string} of $\Db{\Acat}$
	has 
	one of the following forms,
	where $d \in \Z$ is arbitrary. 
	\begin{itemize}
		\item A \emph{usual string} $(d,w)$
		with a usual, 
		asymmetric 
		word $w$.
		\item A \emph{special string} $(d,w,\varepsilon)$ with a special word $w$ and a sign $\varepsilon \in \{+,-\}$.
		\item A \emph{bispecial string}
		$(d,\varepsilon_1,w,\varepsilon_2)$ with a bispecial word 
		$w$ and 
		two signs $\varepsilon_1$, $\varepsilon_2$.
	\end{itemize}
\end{itemize}
In particular, strings are defined using only discrete parameters, while bands depend on a continuous parameter from the base field.
To any band or string $\omega$, 
we associate a `gluing diagram'.
Applying certain combinatorial `gluing rules'
we process such a diagram into an explicit complex $\CP(\omega)$ of $\Db{\Acat}$.
There is a natural notion of equivalence  on the set of bands and strings which reflects when their complexes are isomorphic.
The first main result presented in this article is the following.
\begin{result}[Theorem~\ref{thm:bijection2}]\label{res:A}
The assignment $\begin{td} \omega \ar[mapsto]{r} \& \CP(\omega) \end{td}$ yields a bijection
\begin{align*}
	\begin{cd}
		\{\text{bands and strings of }\Db{\Acat}   \}\big/_{\textstyle \approx}
		\ar[yshift=0pt]{r}{\sim}\&
		\ind \Db{\Acat} 
	\end{cd}
\end{align*}
between
the set of equivalence classes of bands and strings of $\chfXA$ and the set of isomorphism classes of indecomposable 
objects in the category $\Db{\Acat}$.
\end{result}
This result elaborates the classification of indecomposable objects in $\Db{\Acat}$ initiated in \cite{Nodal}, 
whose method we adopt and carry out in full detail.

The second main result involves the following three notions.
\begin{itemize}
\item
For a complex $\CX$ from $\Db{\Acat}$ we
define its \emph{defect} as the number
$$
\defect(\CX) = \sum_{i \in \Z} \dim \Hom_{\Db{\Acat}}(\CX,S_{\star}[i]) \in \N_0
$$
where  $S_\star$ denotes the simple module at vertex $\star$.
In different terms, the defect counts the total number of indecomposable projectives $P_\star$ at all degrees in 
the minimal projective complex isomorphic to $\CX$.
\item
The Gelfand quiver admits a vertical symmetry interchanging the signs $+$ and $-$ in vertices and arrows.
This symmetry induces a $\kk$-algebra automorphism of $\A$,
and thus  an auto-equivalence $\inv$ of $\Db{\Acat}$ such that $\inv^2 \cong \Id$.
\item By known results \cite{IR, vdBergh}, the category $\Db{\Acat}$ admits an auto-equivalence called the \emph{Auslander-Reiten translation $\tau$}, which gives rise to a Serre functor  $\SS = \tau \circ [1]$.
\end{itemize}
It turns out that these homological notions align well with 
the four combinatorial classes of indecomposable objects defined above.
\begin{result}[Proposition~\ref{prp:tau-per} and Theorem~\ref{thm:4cl}]\label{res:B}
Let $\CP$ be a complex of $\DbRep{\A}$.
\begin{enumerate}
	\item \label{B1} The following equivalences hold.
	\begin{align*}
		\begin{array}{lclclcl}
			\delta(\CP) = 0 & \Leftrightarrow & \tau(\CP) \cong \inv(\CP) 
			& \Leftrightarrow & \tau^2(\CP) \cong \CP &\Leftrightarrow & 
			\text{$\CP$ is $\tau$-periodic}.
		\end{array}
	\end{align*}
	\item \label{B2} Assume that $\CP$ is indecomposable and let $\omega$ be a band or string of $\chfXA$ such that $\CP \cong \CP({\omega})$.
	Then  the following equivalences hold.
	\begin{align*}
		\begin{array}{lclclclcl}
			\omega\text{ is a usual string} 
			& \Leftrightarrow & \defect(\CP) = 2 &\Leftrightarrow& \inv(\CP) \cong \CP &\text{and}&\CP\text{ is not $\tau$-periodic}.\\
			\omega\text{ is a special string} 
			& \Leftrightarrow & \defect(\CP) = 1 &\Leftrightarrow& \inv(\CP) \not\cong \CP &\text{and}&\CP\text{ is not $\tau$-periodic}.\\
			\omega\text{ is a bispecial string} & \Leftrightarrow &
			\defect(\CP) = 0&\text{and}& \inv(\CP) \not\cong \CP  
			&\Leftrightarrow&	
			\CP \text{ has $\tau$-period two}.
			\\
			\omega\text{ is a band} & \Leftrightarrow &
			\defect(\CP) = 0&\text{and}& \inv(\CP) \cong \CP  &\Leftrightarrow& 
			\CP \text{ is $\tau$-invariant}.
		\end{array}
	\end{align*} 
\end{enumerate}
\end{result}
The first statement motivates the terminology `defect' and indicates that the functor $\tau$ is essentially determined by the Betti invariants of the complex $\CX$. 
The intuition for the first set of equivalences in \eqref{B2} comes from the fact that the defect corresponds to the number of ends of type $\star$ of the underlying word $w$, while $\inv$ simply flips the signs in $\omega$ if it has any sign data. 
The first set of equivalences yields a homological characterization of the four combinatorial classes. If $\kk = \CC$, this characterization can be interpreted in terms of the derived category of the principal blocks of Harish-Chandra modules over $\SL(2,\RR)$
as the defect as well as the involution admit a Lie-theoretic origin.
The second set of equivalences in \eqref{B2} follows from the first statement and describes the four classes in terms of natural functors on $\DbRep{\A}$.
This result implies that bispecial strings form homogeneous tubes of rank two, while bands form homogeneous tubes of rank one in the Auslander-Reiten quiver of $\Db{\Acat}$. 

In addition to these functors, 
the category $\DbRep{\A}$ admits a contravariant involutive auto-equivalence $\cdual{(-)}$, which originates from contragredient duality on blocks of Harish-Chandra modules equivalent to the  category $\Rep{\A}$.
The next result describes the action of the aforementioned functors on indecomposable objects in $\DbRep{\A}$.
\begin{result}[Proposition~\ref{prp:inv}, Theorems~\ref{thm:tau} and \ref{thm:duals}]\label{res:fun}
Let $\F$ be an auto-equivalence of $\Db{\Acat}$ given by the involution $\inv$, Auslander-Reiten translation $\tau$, or contragredient duality $\cdual{(-)}$. For any band or string $\omega$ we determine the band respectively string $\F(\omega)$ 
such that there is an isomorphism $\F(\CX(\omega))  \cong \CX(\F(\omega))$ in $\Db{\Acat}$.
\end{result}
For instance, if $\omega$ is a band $(d,w,m,\lambda)$,
then $\inv(\omega) = \tau(\omega) = \omega$ and $\cdual{\omega} = (d^t,w^t, m,\frac{1}{\lambda})$ where $d^t = 1-d$ and $w^t$ is given by flipping all orientations in the sequence $w$.

In the last part of this article, we shift focus to the abelian category $\Acat$. 
At first, we describe which band or string complexes are projective resolutions of objects in $\Acat$ (Theorem~\ref{thm:proj-res}).
The parametrizing data has the following approximate description.
\begin{itemize}
\item
Any band of $\Acat$ has the form $(w,m,\lambda)$
with a periodic sequence $w= 
(\ora{p}_{i},\ola{q}_i)_{i \in \Z}$ where $p_i, q_i \in \N$, 
$m \in \N$ and $\lambda \in \kk^*$.
\item
Any string of $\Acat$ is given by a sequence $w = (\alpha, x_1,x_2,\ldots, x_{\ell}, \beta)$
with ends $\alpha,\beta \in \{\hat{\star},\star,\diamond\}$
and certain decorated numbers,
together with the same amount of signs as ends of type $\diamond$ in $w$.
\end{itemize}	
These notions form a compact notation for certain bands and strings $\omega$ of $\DbRep{\A}$, and are motivated by the following result.

\begin{result}[Theorem~\ref{thm:prbij}]
The assignment $\begin{td}
\omega \ar[mapsto]{r} \& M(\omega) = H_0(\CP(\omega))
\end{td}$ yields a bijection
\begin{align*}
\begin{cd}
\{\text{bands and strings of }\Acat  \}\big/_{\textstyle \approx}
\ar[yshift=0pt]{r}{\sim}\&
\ind \Acat
\end{cd}
\end{align*}
on 
the set of equivalence classes of bands and strings of $\Acat$ and the set of isomorphism classes of indecomposable 
objects in the category $\Acat$.
\end{result}

The involution $\inv$ and the duality endofunctor $\cdual{(-)}$ preserve the abelian category $\Acat$, in contrast to the derived Auslander-Reiten translation.
Therefore, one may specialize part of Result~\ref{res:fun} to describe these functors in terms of band and string representations (Theorem~\ref{thm:dualities}).

Next, we describe certain homological invariants of indecomposable objects in $\Acat$.
\begin{result}[Theorem~\ref{thm:inv} and Corollary~\ref{cor:chi}]
For any band or string $\omega$ of $\Acat$ we determine the projective and the injective dimension, the top-, the socle- and the Jordan--Hölder-multiplicities, and the Euler characteristic of the indecomposable representation $M(\omega)$.
\end{result}
These invariants can be computed using only the numerical data underlying $\omega$.
For instance, for a band $\omega = (w,m,\lambda)$, 
the indecomposable representation $M(\omega)$ has projective and injective dimension one, $\head M \cong \soc M \cong (S_+\oplus S_-)^{mk}$
where $2k$ denotes the period of the word $w$, $\uldim M = (mn,mn,mn)$
where $n$ is the sum of all numbers in a minimal periodic part of $w$, and vanishing Euler characteristic.

Finally, we give an explicit description of the indecomposable representations in the category $\Acat$, addressing Gelfand's original question. 
In this course, the description of 
the quiver representation of
$M(\omega)$ can be built from certain cyclic representations determined by the combinatorial data of its band or string $\omega$.
\begin{result}[Propositions~\ref{prp:sh-band}, \ref{prp:bands}, \ref{prp:str-base}]
For any band or string $\omega$ of $\Acat$
we provide an explicit description of vector spaces and matrices of the quiver representation of $M(\omega)$.
\end{result}

There are analogues of the results above for an arbitrary base field $\kk$.
We note that the main results of this article extend to any Khoroshkin quiver, and to principal blocks of the Lorentz groups $\SO^+(n,1)$ for any $n$.
The underlying approach applies more generally to the arrow ideal completion of the path algebra of any skew-gentle quiver.
\subsection*{Outline and approach}
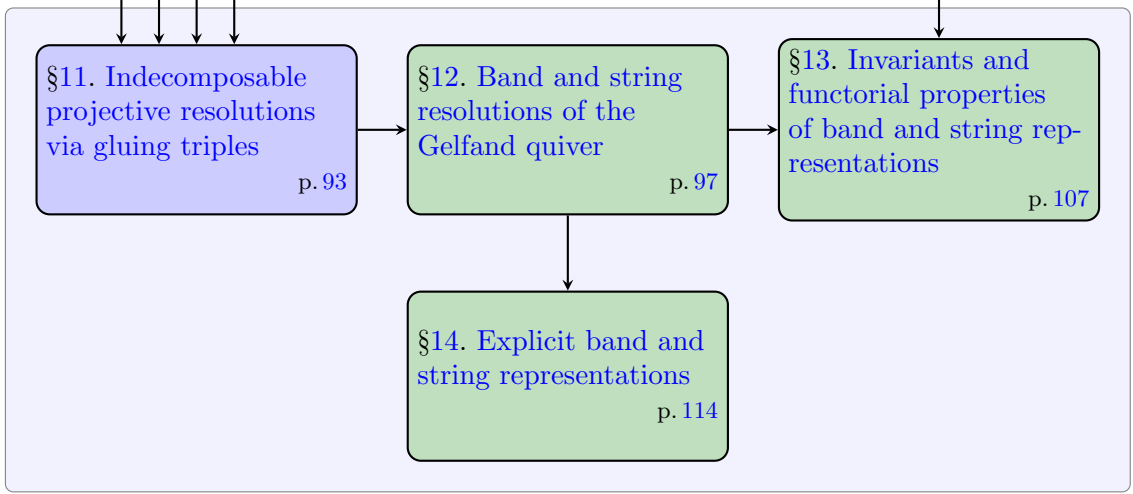
\begin{figure}
\begin{tikzpicture}[
node distance=1cm and 0.65cm,
box/.style={rectangle, draw, thick, minimum width=4cm, text width=4cm, minimum height=2.25cm, align=left, fill=blue!20, rounded corners=5pt},
part/.style={rectangle, draw,  gray, inner sep=0.4cm, rounded corners=3pt, fill=blue!5},
arrow/.style={->, thick, >=stealth, rounded corners=5pt}
]
\node[box] (n1) {\secref{sec:HC_modules}};
\node[box, right=of n1] (n2) {\secref{sec:basic_rep}}; 
\node[box, right=of n2] (n3) {\secref{sec:basic_homological}};
\node[box, below=of n1, yshift=-1.25cm] (n4) {\secref{sec:glue}};
\node[box, right=of n4] (n5) {\secref{sec:reduce}};
\node[box, right=of n5] (n6) {\secref{app:mp}};
\node[box, below=of n4, fill=green!25] (n7) {\secref{sec:complex}};
\node[box, right=of n7] (n8) {\secref{sec:proofs}};
\node[box, right=of n8, fill=green!25] (n9) {\secref{sec:derived_fun}};
\node[box, below=of n7, yshift=-1.25cm] (n10) {\secref{sec:H0-triples}};
\node[box, right=of n10, fill=green!25] (n11) {\secref{sec:resolutions}};
\node[box, fill=green!25, right=of n11] (n12) {\secref{sec:invariants}};
\node[box, fill=green!25, below=of n11] (n13) {\secref{sec:bases}};
\begin{scope}[on background layer]
\node[part, fit=(n1)(n2)(n3), label=above:{Part 1: Motivation and conceptual observations}] (part1) {};
\node[part, fit=(n4)(n5)(n6)(n7)(n8)(n9), label=above:{Part 2: Derived category of the Gelfand quiver}] (part2) {};
\node[part,  fit=(n10)(n11)(n12)(n13), label=above:{Part 3: Abelian category of the Gelfand quiver}] (part3) {};
\draw[arrow, rounded corners=6pt] (n3.south) -- ($(n6.90)+(90:1cm)$) -- ++(2.4cm, 0) -- ++(0, -3.75cm) -- ++(-2.4cm, 0) -- ($(n9.north)$);
\draw[arrow] (n4) -- ($(n7.90)+(90:0.5cm)$) -- ++(-2.4cm, 0) -- ++(0, -3.5cm) -- ++(1.4cm, 0) -- ([xshift=-1cm]n10.north);
\draw[arrow] (n5) -- (n8);
\draw[arrow] ([xshift=-0.5cm]n6.south) -- ++(0, -0.5cm)
-| ([xshift=0.5cm]n8.north);
\draw[arrow] (n7) -- (n8);
\draw[arrow] (n7) -- (n9);
\draw[arrow] ([xshift=-0.5cm]n7.south) -- ([xshift=-0.5cm]n10.north);
\draw[arrow] ([xshift=0.5cm]n4.south) -- ++(0, -0.5cm)
-| ([xshift=-0.5cm]n8.north);
\draw[arrow] (n9.south)
-- (n12.north);
\draw[arrow] (n8) -- (n9);
\draw[arrow] (n10) -- (n11);
\draw[arrow] (n11) -- (n12);
\draw[arrow] (n11) -- (n13);
\draw[arrow] (n8.south) 
-- ++(0, -0.75cm)
-|
([xshift=0cm]n10.north);
\draw[arrow] ([xshift=-0.5cm]n9.south) 
-- ++(0, -1cm)
-|
([xshift=0.5cm]n10.north);	 
\end{scope}
\end{tikzpicture}
\caption{Outline of contents and their proof-dependencies}
\label{fig:dep}
\end{figure}

In the following, we briefly describe the structure and the main methods of this article.
The proof-dependencies between the sections are depicted in Figure~\ref{fig:dep}.

Sections \ref{sec:HC_modules} to \ref{sec:basic_homological} form the first part of this work, which provides
an overview of the representation theory of the Gelfand quiver
without using heavy combinatorial machinery.
In Section~\ref{sec:HC_modules}, we give an account of the Lie-theoretic origins of the category $\Acat$ of nilpotent representations of the Gelfand quiver, as well as the involution $\inv$ and 
contragredient duality $\cdual{(-)}$.
Section~\ref{sec:basic_rep} 
yields natural examples of indecomposable representations of the Gelfand quiver, in particular, a classification of its cyclic representations.
In Section~\ref{sec:basic_homological}, we consider  homological properties of the ring $\A$, introduce the main functors on the derived category $\Db{\Acat}$ 
and relate the derived Auslander-Reiten translation $\tau$ to the defect, proving Result~\ref{res:B}~\eqref{B1}.

The primary goal of the second part of the paper, Sections~\ref{sec:glue} to \ref{sec:derived_fun}, is to provide a classification of indecomposable objects in the category $\Db{\Acat}$.
Following \cite{Nodal}, our method is based on the use of 
the \emph{category of gluing triples} 
$\Gcatfd{\A}$ and  the \emph{category of regular representations}
$\regrep{\fXA}$ of a certain combinatorial datum $\fXA$.
These categories arise together with certain natural functors
\begin{align}\label{eq:reduc}
\begin{cd}
\Db{\Acat} \ar{r}{\F} \&  \Gcatfd{\A} \ar{r}{\Red} \& \regrep{\fXA}
\end{cd}
\end{align}
which induce  a bijection and 
a surjective map on sets of isomorphism classes of indecomposable objects, respectively,
\begin{align}
\label{eq:reduc2}
\begin{cd}
\ind \Db{\Acat} \ar[yshift=4pt]{r}{\F}[swap]{\sim} 
\ar[yshift=-4pt,<-]{r}[swap]{\Gl}
\& \ind \Gcatfd{\A} \ar[yshift=4pt,twoheadrightarrow]{r}{\Red} 
\ar[yshift=-4pt,<-, hookleftarrow]{r}[swap]{\Ind}
\& \ind \regrep{\fXA}
\end{cd}
\end{align}
where $\Gl$ is the inverse of $\F$ and $\Ind$ the left inverse of $\Red$.
The complement of the image of the embedding $\Ind$
can be described explicitly and is given by countably many isomorphism classes of indecomposable objects. 
In summary, this approach allows to reduce the classification problem 
of the category $\Db{\Acat}$
to that of the category $\regrep{\fXA}$.
The latter can be formulated as a  matrix problem, which has been solved by Bondarenko~\cite{B1,B2}. His solution provides explicit canonical forms for the indecomposable objects 
in the category $\regrep{\fXA}$.

Section~\ref{sec:glue}
introduces the category $\Gcatfd{\A}$, the functor $\FF$ and the inverse $\Gl$ on indecomposable objects in diagram \eqref{eq:reduc}.
Section~\ref{sec:reduce} contains the definitions of the category $\regrep{\fXA}$, the functor $\Red$ and the assignment $\Ind$ and ends with a refined categorical picture underlying \eqref{eq:reduc}.
Section~\ref{app:mp} forms an interlude, in which
we recall Bondarenko's canonical forms for a class of matrix problems including $\regrep{\fXA}$ without proof.
The sole purpose of 
Section~\ref{sec:complex} is to 
formulate the resulting classification of indecomposable objects in $\Db{\Acat}$ (Result~\ref{res:A}).
The proof of this classification is based on the conversion of Bondarenko's canonical forms into the corresponding indecomposable objects in $\Db{\Acat}$ along the composition $\Gl\circ \Ind$ in \eqref{eq:reduc2}, which is carried out in Section~\ref{sec:proofs}.

In the third and last part of this article, we return our focus on the abelian category $\Acat$.
Section~\ref{sec:H0-triples} addresses the question how to identify the band or string complexes which are projective resolutions of objects in $\Acat$.
Some necessary conditions are provided using the language of gluing triples. They are verified to be sufficient 
by going through the list of candidate band and string complexes.
In Section~\ref{sec:resolutions}, 
we introduce the bands and strings of $\Acat$ in a compact notation, and describe their associated projective resolutions from scratch. 
An intermediate step in the gluing process of a band or string resolution is its `cyclification', a simpler projective resolution with the same Betti numbers.
The cyclification already determines the main homological invariants computed in Section~\ref{sec:invariants}.
Finally, in Section~\ref{sec:bases} we compute a quiver representation of the homology of each band and string resolution.
A basis of the latter can be identified with a natural basis of the homology of the associated cyclification.

The main results of this article are stated in a self-contained manner.
In particular, readers interested only in results on the derived category 
$\Db{\Acat}$ may start with Section~\ref{sec:complex}, followed by Section~\ref{sec:derived_fun} skipping the proofs.
In the same spirit, readers interested only in results on the abelian category $\Acat$ may begin with Section~\ref{sec:resolutions}, and then consult Section~\ref{sec:invariants} 
on invariants and functors,
or Section~\ref{sec:bases} on the final description of indecomposable representations.
Both reading paths presume familiarity with the notions below.

\subsection*{Main notation}
Throughout this article, we will use the following notation without 
frequent back references.
\begin{longtable}{cll}
$\kk$ & an algebraically closed field & \\
$\Rx$ & the ring of formal power series $\kk \llbracket t \rrbracket$  \\
$\mx $ & the maximal ideal  of $\Rx$ \\
$\A$ & the Gelfand order & \eqref{E:GelfandOrder} \\
$\Rep{\A}$ & the category of  finite-dimensional left $\A$-modules & Definition~\ref{dfn:GO}\\
$S_\star$, $S_+$, $S_-$ & simple representations of the Gelfand quiver & Remark~\ref{rmk:sim} \\
$P_\star$, $P_+$, $P_-$, $I_\star$ & indecomposable  $\A$-lattices &  \eqref{eq:lattices}\\
$\inv$ & involution functor & \eqref{E:reflection},
\S\ref{subsec:fun3}\\
$\cdual{(-)}$ & the contragredient duality & \eqref{E:LieDuality}, \S\ref{subsec:fun3}\\
$\dual{(-)}$ & duality functor preserving each simple $\A$-module &  \eqref{E:dual}, \S\ref{subsec:fun3} \\
$\DbRep{\A}$ & the bounded derived category of $\Rep{\A}$ & \S\ref{subsec:basic-der} \\
\end{longtable}
Moreover, we will denote by $\N$ the set of positive integers  and by $\N_0$ 
the set of non-negative integers.

\bigskip
\noindent
\emph{Acknowledgement}. We are grateful to Yu.~Drozd for numerous discussions of the results of this paper. 
This work grew out of the PhD thesis of the second-named author \cite{Gnedin}.
Our work was partially supported by the  German Research Foundation SFB-TRR 358/1 2023 -- 491392403.

\clearpage

\clearpage

\newpage
\section{\texorpdfstring{Harish-Chandra modules for $\SL_2(\RR)$ and the Gelfand quiver}{Harish-Chandra modules for SL(2,R) and the Gelfand quiver}}
\label{sec:HC_modules}
This section provides some Lie-theoretic background motivating the  study of the representation theory of the Gelfand quiver.
First, we recall that the category of Harish-Chandra modules for $\SL_2(\RR)$ 
admits a block decomposition
such that each block is equivalent to the category of $\CC$-linear nilpotent representations of the one-loop quiver, the cyclic quiver with two vertices
or the Gelfand quiver.
Then we consider a natural involution 
and the contragredient duality on Harish-Chandra modules, as well as their incarnations
$\inv$ and  $\cdual{(-)}$ on the category $\Rep{\A}$ of nilpotent representations of the Gelfand quiver.

\subsection{\texorpdfstring{Blocks of Harish-Chandra modules for $\SL_2(\RR)$}{Blocks of Harish-Chandra modules for SL(2,R)}}
\label{subsec:HCblocks}
\noindent
Let  $G = \SL_2(\RR)$, $K = \SO_2(\RR)$, $\lieg \colonequals \mathsf{Lie}(G) \otimes_{\RR} \CC \cong \mathfrak{sl}_2(\CC)$ and $U(\lieg)$ be the universal enveloping algebra of $\lieg$. We denote
\begin{align*}
	h = \begin{pmatrix}
		0 & -i \\
		i & 0
	\end{pmatrix}, &&
	x = \frac{1}{2}
	\begin{pmatrix}	1 & i \\
		i & -1
	\end{pmatrix},&&
	y = \frac{1}{2}
	\begin{pmatrix}
		1 & -i \\
		-i & -1
	\end{pmatrix}.
\end{align*}
Then $\lieg = \bigl\langle h, x, y\bigr\rangle_\CC$ and 
\begin{align}
	\label{E:Standard}
	[h, x] = 2x,&& [h, y] = -2y, &&[x, y] = h.
\end{align}

\noindent
Recall the following standard results.
\begin{enumerate}
	\item[(a)] $K$ is a maximal compact subgroup of $G$ and any other such subgroup is conjugate to $K$. 
	\item[(b)] Any continuous topologically irreducible representation of $K$ is 
	one-dimensional and isomorphic to the map
	\begin{align}\label{E:KDual}
		\begin{td}\chi_n \colon K \ar{r} \& \CC^\ast, \& k_\theta \colonequals  
			\begin{pmatrix}
				\cos(\theta) & \sin(\theta) \\
				-\sin(\theta) & \cos(\theta)
			\end{pmatrix}
			\ar[mapsto]{r} \& \exp(2\pi in\theta) \end{td}
	\end{align}
	for some $n \in \Z$. 
\end{enumerate}

\begin{dfn} A complex vector space $M$ is a Harish-Chandra $(\lieg, K)$-module if it has a structure $(M, \circ)$ of a finitely generated $U(\lieg)$-module as well as a structure $(M, \,\cdot\,)$ of an admissible algebraic representation of $K$, that is, $M$ is a direct sum of irreducible representations \eqref{E:KDual} of $K$ with finite multiplicities, such that 
	$$
	\left.\frac{d}{dt}\right|_{t = 0} \exp(t z)\cdot  v  = z \circ v 
	$$
	for any $z \in \mathsf{Lie}(K)$ and $v \in M$.
\end{dfn}

\begin{rmk}\label{rmk:HC}
	Let $M$ be a left $U(\lieg)$-module. For any $n\in \Z$ we denote $M_n \colonequals \bigl\{v \in M \mid h \circ v = nv \bigr\}$. Then $M$ is a Harish-Chandra $(\lieg, K)$-module if and only if
	\begin{enumerate}
		\item[(a)] $M$ is finitely generated over $U(\lieg)$,
		\item[(b)] $M \cong \bigoplus_{n \in \mathbb{Z}} M_n$ as a complex vector space and $\dim_{\CC}(M_n) < \infty$ for any $n \in \Z$. 
	\end{enumerate}
	In what follows, we shall omit $\circ$  and $\cdot$ when speaking about the action of $\lieg$  and $K$ on a Harish-Chandra $(\lieg, K)$-module.  We shall denote by $\HC(\lieg, K)$ the category of Harish-Chandra $(\lieg, K)$-modules. It is clear that $\HC(\lieg, K)$ is abelian. 
\end{rmk}

\smallskip
\noindent
Let $M$ be an object of $\HC(\lieg, K)$.
It follows from \eqref{E:Standard} that $x(M_{p-1}) \subseteq M_{p+1}$ and 
$y(M_{p+1}) \subseteq M_{p-1}$ for any $p \in \Z$. Consider the full subcategories 
$\HC^\pm(\lieg, K)$ of $\HC(\lieg, K)$ consisting of those $M$ for which
$M_n = 0$ for all odd $n$ (in the case $+$), respectively for all even $n$ (in the case $-$). It is clear that $\HC^\pm(\lieg, K)$ are mutually orthogonal abelian subcategories of $\HC(\lieg, K)$ such  that we have a decomposition 
$$
\HC(\lieg, K) = \HC^{+}(\lieg, K) \oplus \HC^{-}(\lieg, K).
$$
However, the category $\HC(\lieg, K)$ admits a much finer splitting. 
Let 
\begin{equation}\label{E:Casimir}
	c \colonequals h^2 - 2h + 4xy  = h^2 + 2h + 4yx \in U(\lieg)
\end{equation} be the Casimir element of $U(\lieg)$. Then $\CC[c]$ is the center of $U(\lieg)$.

\begin{dfn}\label{dfn:HCb}
	For any $\gamma \in \CC$ let $\HC_\gamma(\lieg, K)$ be the full subcategory of 
	$\HC(\lieg, K)$ consisting of those modules $M$ for which 
	there exists $m \in \N$ (depending on $M$) such that $(c- \gamma \mathbbm{1})^m \cdot M = 0$. It is clear that $\HC_\gamma(\lieg, K)$ is an abelian subcategory of $\HC(\lieg, K)$.
\end{dfn}

\begin{thm}\label{T:blocks}
	The following
	results are true.
	\begin{itemize}
		\item[(a)] Any indecomposable object of $\HC(\lieg, K)$ belongs to some subcategory $\HC_\gamma(\lieg, K)$ for a uniquely determined $\gamma \in \CC$. 
		\item[(b)] For $\gamma' \ne \gamma''$ the corresponding subcategories $\HC_{\gamma'}(\lieg, K)$ and $\HC_{\gamma''}(\lieg, K)$ are $\Ext$-orthogonal to each other and the  category $\HC(\lieg, K)$ splits into a sum of blocks
		\begin{equation}\label{E:blocks}
			\HC(\lieg, K) = \bigoplus_{\gamma  \in \CC } \HC_{\gamma}(\lieg, K).
		\end{equation}
		Moreover, every object of the category $\HC_{\gamma}(\lieg, K)$ has finite length for any $\gamma \in \CC$. 
		\item[(c)] If $\gamma \in \CC \setminus \bigl\{l^2-1  \, \big| \,  l \in \N_0\bigr\}$ then we have a further splitting $$\HC_{\gamma}(\lieg, K) = \HC_{\gamma}^+(\lieg, K) \oplus \HC_{\gamma}^-(\lieg, K).$$
		Moreover, each category $\HC_{\gamma}^{\pm}(\lieg, K)$ is equivalent 
		to the category of finite dimensional nilpotent representations of the one-loop quiver 
		$$\begin{qrep}
			\bt \arrow[loop, looseness=7.5, in=35, out=-35]
		\end{qrep}.$$
		\item[(d)] For any $l \in \N_0$ let $\gamma(l) = l^2 -1$. Then  the block $\HC_{\gamma(0)}(\lieg, K)$ is equivalent to the category of finite dimensional nilpotent representations of the cyclic quiver
		\begin{equation}\label{E:CyclicQuiver}
			\begin{qrep}			
				{\bt} \ar[start anchor=north east, end anchor=north west, bend left,yshift=-3pt]{r}{c_+}
				\& 
				{\bt}   
				\ar[start anchor=south west, end anchor=south east,bend left,yshift=3pt]{l}{c_-}
			\end{qrep}
		\end{equation}
		whereas  the block $\HC_{\gamma(l)}(\lieg, K)$ is equivalent to the category of finite dimensional nilpotent representations of the Gelfand quiver \eqref{E:GelfandQuiver} for any $l \in \N$.
	\end{itemize}
\end{thm}

\begin{proof}
	The statements (a) and (b) are true in a much greater generality. For the considered case of $\lieg = \mathfrak{sl}_2(\CC)$ we refer to \cite[Section 3.6]{MazorchukBook}. Note that $\Hom$-orthogonality of $\bigl(\HC_{\gamma}(\lieg, K)\bigr)_{\gamma \in \CC}$ implies their $\Ext$-orthogonality. 
	
	\smallskip
	\noindent
	(c) For any object $M$ of $\HC(\lieg, K)$ and $p \in \Z$, consider the following diagram of vector spaces and linear maps
	\begin{equation}\label{E:Fragment}
		\begin{qrep}			
			{M_{p-1}} \ar[start anchor=north east, end anchor=north west, bend left,yshift=-3pt]{r}{x_p}
			\& 
			{M_{p+1}}   
			\ar[start anchor=south west, end anchor=south east,bend left,yshift=3pt]{l}{y_p}
		\end{qrep}
	\end{equation}
	where $x_p := x\Big|_{M_{p-1}}$ and $y_p := y\Big|_{M_{p+1}}$.
	It follows using \eqref{E:Casimir} that
	\begin{equation}\label{E:KeyIdentities}
		\left\{
		\begin{array}{l}
			x_p y_p = \frac{1}{4}\bigl(c - h^2 + 2h\bigr)\Big|_{M_{p+1}} = 
			\frac{1}{4}\bigl(c_{p+1} -(p^2-1)\mathbbm{1}\bigr) \\ 
			y_p x_p = \frac{1}{4}\bigl(c - h^2 - 2h\bigr)\Big|_{M_{p-1}} = 
			\frac{1}{4}\bigl(c_{p-1} -(p^2-1)\mathbbm{1}\bigr)
		\end{array}
		\right.
	\end{equation}
	where $c_{p\pm 1} = c\big|_{M_{p\pm 1}}$. 
	As a consequence, if $\gamma \in \CC \setminus \bigl\{l^2-1  \, \big| \,  l \in \N_0\bigr\}$ and $M \in \HC_\gamma^\pm(\lieg, K)$ then the morphisms $x_p$ and $y_p$ are isomorphisms for each fragment 
	\eqref{E:Fragment} of  $M$. We can choose bases of the weight spaces of $M$ in such a way that all maps $y_p$ in \eqref{E:Fragment}  are given by the identity matrices. Then it follows from \eqref{E:KeyIdentities} that $w:= c_{p-1} = c_{p+1}$. 
	Moreover, from \eqref{E:KeyIdentities}  we obtain:  
	$$w = c_{p+1+2k}, \;  y_{p+2k} = \mathbbm{1} \; \mbox{\rm and}\; x_{p+2k} = \frac{1}{4}\bigl(w-((p+2k)^2-1) \mathbbm{1}\bigr)\quad \mbox{\rm for all} \quad k \in \Z.$$
	Assigning to $M$ the matrix $w$, we obtain the stated equivalence of categories. 
	
	\smallskip
	\noindent
	(d) Let $M \in \HC_{-1}(\lieg, K)$. Then $M \in \HC^-(\lieg, K)$. Moreover, it follows from \eqref{E:KeyIdentities} that $x_{2k}$ and $y_{2k}$ are isomorphisms for all $k \in \Z \setminus \{0\}$, whereas both  compositions $x_0 y_0$ and $y_0 x_0$ are nilpotent operators. Assigning to $M$ the  following representation
	\begin{equation}\label{E:ReprCyclic}
		\begin{qrep}			
			{M_{-1}} \ar[start anchor=north east, end anchor=north west, bend left,yshift=-3pt]{r}{x_0}
			\& 
			{M_{+1}}   
			\ar[start anchor=south west, end anchor=south east,bend left,yshift=3pt]{l}{y_0}
		\end{qrep}
	\end{equation}
	of the cyclic quiver \eqref{E:CyclicQuiver}, we obtain an equivalence of categories; see \cite[Section 3.8]{MazorchukBook} or \cite[Proposition 3.6]{ABR}.
	
	\smallskip
	\noindent
	Now, consider $M \in \HC_{\gamma(l)}(\lieg, K)$ with $l \in \N$. Then the $\lieg$-module structure on $M$ is determined by the following diagram 
	
	\begin{align}\label{E:DiagramHCmodule}
		\begin{qrep}
			\vphantom{M_l}
			{\ldots}
			\ar[start anchor=north east, end anchor=north west, bend left,yshift=-3pt]{r}{x}
			\&
			{M_{-l-1}}
			\ar[start anchor=south west, end anchor=south east,bend left,yshift=3pt]{l}{y}
			\ar[start anchor=north east, end anchor=north west, bend left,yshift=-3pt]{r}{x_{-l}}
			\&
			{M_{-l+1}}
			\ar[start anchor=south west, end anchor=south east,bend left,yshift=3pt]{l}{y_{-l}}
			\ar[start anchor=north east, end anchor=north west, bend left,yshift=-3pt]{r}{x}
			\&
			{\ldots}
			\vphantom{M_l}
			\ar[start anchor=south west, end anchor=south east,bend left,yshift=3pt]{l}{y}
			\ar[start anchor=north east, end anchor=north west, bend left,yshift=-3pt]{r}{x}
			\&
			{M_{l-1}}
			\ar[start anchor=south west, end anchor=south east,bend left,yshift=3pt]{l}{y}
			\ar[start anchor=north east, end anchor=north west, bend left,yshift=-3pt]{r}{x_l}
			\&
			{M_{l+1}}
			\ar[start anchor=south west, end anchor=south east,bend left,yshift=3pt]{l}{y_l}
			\ar[start anchor=north east, end anchor=north west, bend left,yshift=-3pt]{r}{x}
			\&
			{\ldots}\vphantom{M_l}
			\ar[start anchor=south west, end anchor=south east,bend left,yshift=3pt]{l}{y}				
		\end{qrep}
	\end{align}
	where we omit certain indices to simplify our notation. 
	It follows from \eqref{E:KeyIdentities} that the compositions $x_{\pm l} y_{\pm l}$ and 
	$ y_{\pm l} x_{\pm l}$ are nilpotent, whereas the remaining linear maps $x$ and $y$ are isomorphisms. 
	\\
	\noindent
	Set $\begin{td} \tilde{x}:= x_{l-2} \dots x_{-l+2}: M_{-l+1} \ar{r} \&  M_{l-1} \end{td}$ and 
	$\begin{td} \tilde{y}:= y_{-l+2} \dots y_{l-2}: M_{l-1} \ar{r} \& M_{-l+1} \end{td}$. An easy calculation (see \cite[Lemma 1.7]{ABR}) shows that
	\begin{equation}
		\tilde{x} x_{-l} y_{-l} = y_l x_l \tilde{x} \quad \mbox{\rm and} \quad x_{-l} y_{-l} \tilde{y} = \tilde{y} y_l x_l.
	\end{equation}
	We can  attach  to $M$ (given by the diagram \eqref{E:DiagramHCmodule})  the following representation
	\begin{equation}\label{E:HCtoRepsgelfand}
		\begin{gqrep}{M_{l+1}}{M_{-l+1}}{M_{-l-1}}{\tilde{x}^{-1} y_l}{x_{-l}}{x_l \tilde{x}}{y_{-l}}\end{gqrep}
	\end{equation}
	of the Gelfand quiver \eqref{E:GelfandQuiver}. It can be shown that this assignment defines 
	the asserted equivalence of categories; see \cite[Theorem 1.7]{ABR} or \cite[Section 3.9]{MazorchukBook} for further details. 
\end{proof}

\begin{dfn}\label{dfn:GO}
	Recall that	 $\kk$ is an 
	algebraically closed
	field,
	$\Rx  = \kk\llbracket t\rrbracket$, $\mx = (t)$ and
	\begin{equation}\label{E:GelfandOrder}
		\A = 
		\begin{pmatrix}
			\Rx & \mx & \mx\\
			\Rx & \Rx & \mx \\
			\Rx & \mx & \Rx
		\end{pmatrix}\subset
		{\Mat}_{3 \times 3}(\Rx).
	\end{equation}
	Note that $\A$ is isomorphic to the arrow ideal completion of the path algebra of the Gelfand quiver \eqref{E:GelfandQuiver}. This isomorphism is given by  the following identifications.
	\begin{align*}
		\begin{array}{lllllll}
			e_\star = 
			\begin{pNiceArray}{ccc}
				1 & 0 & 0 \\
				0 & 0 & 0 \\
				0 & 0 & 0
			\end{pNiceArray}, &
			e_+ = 
			\begin{pNiceArray}{ccc}
				0 & 0 & 0 \\
				0 & 1 & 0 \\
				0 & 0 & 0
			\end{pNiceArray}, & 
			e_- = 
			\begin{pNiceArray}{ccc}
				0 & 0 & 0 \\
				0 & 0 & 0 \\
				0 & 0 & 1
			\end{pNiceArray}
			\\
			b_- = \begin{pNiceArray}{ccc}
				0 & 0 & t \\
				0 & 0 & 0 \\
				0 & 0 & 0
			\end{pNiceArray},
			&
			a_- = \begin{pNiceArray}{ccc}
				0 & 0 & 0 \\
				0 & 0 & 0 \\
				1 & 0 & 0
			\end{pNiceArray}&
			b_+= \begin{pNiceArray}{ccc}
				0& t & 0 \\
				0 & 0 & 0 \\
				0 & 0 & 0
			\end{pNiceArray}, &
			a_+= \begin{pNiceArray}{ccc}
				0 & 0 & 0 \\
				1 & 0 & 0 \\
				0 & 0 & 0
			\end{pNiceArray}.
		\end{array}
	\end{align*}
	Note that $\A$ is an order (see for instance \cite{CurtisReiner} for the corresponding definition). Hence, in what follows $\A$ will be called the \emph{Gelfand order}. The algebra $\A$ belongs to the class of so-called nodal orders; see \cite{NodalFirst} as well as 
	\cite{Nodal, NodalOrders}.
	Next, $\Rep{\A}$ denotes the category of finite dimensional left $\A$-modules. It is equivalent to the category of
	finite dimensional nilpotent representations of \eqref{E:GelfandQuiver}.
\end{dfn}

\begin{rmk}\label{rmk:sim}
	For any $l \in \N$ we have an equivalence of categories 
	$$
	\begin{cd}
		\EE_l \colon \HC_{\gamma(l)}(\lieg, K) \ar{r}{\sim} \& \Rep{\A}
	\end{cd}
	$$
	given by \eqref{E:HCtoRepsgelfand} with $\kk = \CC$. Indecomposable finite dimensional objects of $\HC(\lieg, K)$ are determined by their dimensions up to isomorphism. Let $F_l$ be such an object of dimension $l \in \N$. Then it belongs to the exceptional block $\HC_{\gamma(l)}(\lieg, K)$ and 
	\begin{equation}\label{E:FiniteDimSimple}
		\EE_l(F_l) = S_\star = 
		\left[\begin{gqrep}{0}{\kk}{0}{}{}{}{}
		\end{gqrep}
		\right]
	\end{equation}
	is the simple $\A$-module associated to the vertex $\star$ of the Gelfand quiver \eqref{E:GelfandQuiver}. For the remaining two simple $\A$-modules $S_\pm$ we have:
	\begin{itemize}
		\item[(a)] $S_-$ is the image of the highest weight module $U(\lieg)/U(\lieg)\langle x, h +(l+1)\rangle$.
		\item[(b)] $S_+$ is the image of the lowest weight module \, $U(\lieg)/U(\lieg)\langle y, h-(l+1)\rangle$.
	\end{itemize}
\end{rmk}

\begin{rmk} Consider the following matrix algebra 
	\begin{equation}\label{E:Hereditary}
		O
		= 
		\begin{pmatrix}
			\Rx & \mx \\
			\Rx &  \Rx
		\end{pmatrix} \subset {\Mat}_{2 \times 2}(\Rx).
	\end{equation}
	Then $O$ is isomorphic to the arrow ideal completion of the path algebra of the cyclic quiver \eqref{E:CyclicQuiver}. In these terms we have an equivalence of categories
	$$
	\begin{cd}	\EE_0 \colon \HC_{\gamma(0)}(\lieg, K) \ar{r}{\sim} \&  \Rep{O} \end{cd}
	$$
	given by the assignment \eqref{E:ReprCyclic}. 
\end{rmk}

\subsection{The involution on Harish-Chandra modules}
\label{subsec:HC-inv}
\noindent
There is a 
natural involution
\begin{align}
	\hcinvaut = \mathsf{Ad}_v: \; (\lieg, K) \longrightarrow (\lieg, K), &&
	\text{where }v = \begin{pmatrix} 1 & 0 \\ 0 & -1\end{pmatrix} \in \GL_2(\RR)
\end{align} 
Note that  $\hcinvaut^2 = \mathsf{id}$, $\hcinvaut(h) = -h$, $\hcinvaut(x) = y$ and $\hcinvaut(y) = x$.
We get an induced involutive  auto-equivalence $\begin{td} \hcinv \colon \HC(\lieg, K) \ar{r}{\sim} \& \HC(\lieg, K) \end{td}$. Since $\hcinvaut(c) = c$, the functor $\hcinv$ preserves each block $\HC_\gamma(\lieg, K)$ for any $\gamma \in \CC$.

\smallskip
\noindent
On the other hand, the algebra $\A$  has an obvious involution $\begin{td} \invaut \colon \A \ar{r}{\sim} \& \A \end{td}$, 
which is uniquely determined by the prescriptions
\begin{align}\label{E:Involution}
	\invaut(e_\pm) = e_\mp, \; 
	\invaut(a_\pm) = a_\mp \; \; \mbox{\rm and} \; \; 
	\invaut(b_\pm) = b_\mp.
\end{align} 
We obtain the induced involutive auto-equivalence 
$
\inv \colon \begin{td} \Rep{\A} \ar{r}{\sim} \&  \Rep{\A}\end{td}
$
given by the rule
\begin{equation}\label{E:reflection}
	\begin{gqrep}{V_+}{V_\star}{V_-}{B_+}{B_-}{A_+}{A_-}\end{gqrep}
	\quad \stackrel{\inv}\longmapsto \quad
	\begin{gqrep}{V_-}{V_\star}{V_+}{B_-}{B_+}{A_-}{A_+}\end{gqrep}.
\end{equation}
Let $\HC_\circ(\lieg, K) = \HC_{\gamma(1)}(\lieg, K)$ be the principal block (i.e.~the block containing the one-dimensional representation of $\lieg$). Then the proof of the following result is straightforward.

\begin{lem} The diagram of categories and functors 
	\begin{align}
		\begin{cd}
			\HC_\circ(\lieg, K) \ar{r}{\EE_1} \ar{d}{\hcinv} \& \Rep{\A} \ar{d}{\inv} \\
			\HC_{\circ}(\lieg, K) \ar{r}{\EE_1} \& 
			\Rep{\A}
		\end{cd}
	\end{align}
	commutes in the sense that there is an isomorphism of functors $\EE_1 \circ \hcinv \cong \inv \circ \EE_1$.
	\qed
\end{lem}  
\begin{rmk}
	With slightly more effort (e.g.~using \cite[Lemma 1.8]{ABR}) one can show that for any $l \in \N$ and  $M \in \HC_{\gamma(l)}(\lieg, K)$ we have an isomorphism $ \inv( \EE_l(M)) \cong \EE_l(\hcinv(M))$.
\end{rmk}

\subsection{Contragredient duality on Harish-Chandra modules}
\noindent
Recall that the \emph{contragredient dual} $(\hcdual{M}, \hcdual{\pi})$ of
an object  $(M, \pi)$  of $\HC(\lieg, K)$ is defined as follows.
Let
$M^\ast := \Hom_{\CC}(M, \CC)$ and $ \begin{td} (\pi(z))^\ast \colon {M^\ast} \ar{r} \&
	{M^\ast} \end{td}$ be the dual of the linear map of $\pi(z)$ for any $z \in \lieg$. Then we pose
$$
\hcdual{M} := \left\{l \in  M^\ast \,\big|\, \dim_\CC(K\cdot l) < \infty \right\} \;
\mbox{and}\; \hcdual{\pi}(z) = -\bigl(\pi(z)\bigr)^\ast\Big|_{\hcdual{M}} \; \mbox{for}\;  z \in \lieg.
$$
If $M \cong \bigoplus_{n \in \mathbb{Z}} M_n$ is the decomposition of $M$ into a direct sum of eigenspaces of $h$ then $\hcdual{M} \cong \bigoplus_{n \in \mathbb{Z}} M_n^\ast$.
This gives  a duality functor 
$$
\begin{td} \hcdual{(-)} \colon \HC(\lieg, K) \ar{r} \&  \HC(\lieg, K), \& (M, \pi) \ar[mapsto]{r} \& (\hcdual{M}, \hcdual{\pi}),\end{td}
$$
mapping each category  $\HC_\gamma(\lieg, K)$ to itself for any $\gamma \in \CC$. 

\begin{lem}\label{L:twoduals}
	For any $l \in \N$ consider the following diagram of categories and functors
	$$
	\begin{cd}
		\HC_{\gamma(l)}(\lieg, K) \ar{d}{{\hcdual{(-)}}} \ar{r}{\EE_l} \& 				\Rep{\A} \ar{d}{\cdual{(-)}} \\
		\HC_{\gamma(l)}(\lieg, K) \ar{r}{\EE_l}  \& \Rep{\A}
	\end{cd}
	$$
	where the functor $\cdual{(-)}$ on $\Rep{\A}$ is given by the rule
	\begin{align}\label{E:LieDuality}
		\begin{gqrep}{V_+}{V_\star}{V_-}{B_+}{B_-}{A_+}{A_-}\end{gqrep}
		\quad \stackrel{\cdual{(-)}}\longmapsto \quad
		\begin{gqrep}{V_-^*}{V^*_\star}{V^*_+}{A^*_-}{A^*_+}{B^*_-}{B^*_+}\end{gqrep}.
	\end{align}
	Then there exists an isomorphism of functors $\cdual{(-)} \circ \EE_l \cong \EE_l \circ  \hcdual{(-)}$. 
\end{lem}

\begin{proof} Note the following easy facts.
	
	\smallskip
	\noindent
	\underline{Claim 1}. For any  $\lambda_+, \lambda_- \in \left\{-1, 1\right\}$ the assignment
	\begin{align}\label{E:LieDualite}
		\begin{gqrep}{V_+}{V_\star}{V_-}{B_+}{B_-}{A_+}{A_-}
		\end{gqrep}
		\quad \longmapsto \quad
		\begin{gqrep}{V_+}{V_\star}{V_-}{\lambda_+ B_+}{\lambda_- B_-}{\lambda_+ A_+}{\lambda_- A_-}
		\end{gqrep}
	\end{align}
	defines an endofunctor of the category $\Rep{\A}$ (whose action on morphisms is the identity), which is isomorphic to the identity functor. 
	
	\smallskip
	\noindent
	\underline{Claim 2}. Let $M$ be an object of $\HC_{\gamma(l)}(\lieg, K)$ given by its ``diagrammatic presentation'' \eqref{E:DiagramHCmodule}.  Then we have the following representation of the Gelfand quiver \eqref{E:GelfandQuiver}:
	\begin{equation}\label{E:HCtoRepsgelfand2}
		\begin{gqrep}
			{M_{l+1}}
			{M_{l-1}}
			{M_{-l-1}} 
			{y_l}
			{\tilde{x} x_{-l}}
			{x_l}
			{y_{-l} \tilde{x}^{-1}}
		\end{gqrep}
	\end{equation}
	This assignment defines another equivalence of categories 
	$\begin{td} \EE'_l \colon \HC_{\gamma(l)}(\lieg, K) \ar{r}{\sim} \&  \Rep{\A} \end{td}$. Note that the following diagram in the category of vector spaces is commutative: 
	\begin{align}
		\begin{qrep}
			{M_{-l-1}} 
			\ar[start anchor=north east, end anchor=north west, bend left,yshift=-3pt]{r}{x_{-l}} \ar{dd}{\Id}
			\& 
			{M_{-l+1}}
			\ar[start anchor=south west, end anchor=south east,bend left,yshift=3pt]{l}{y_{-l}}
			\ar{dd}{\tilde{x}}
			\ar[start anchor=south east, end anchor=south west, bend right, yshift=3pt]{r}[swap]{x_{l} \tilde{x}}
			\&
			{M_{l+1}}
			\ar[start anchor=north west, end anchor=north east, bend right, yshift=-3pt]{l}[swap]{\tilde{x}^{-1} y_{l}} \ar{dd}{\Id}
			\\
			\\
			{M_{-l-1}} 
			\ar[start anchor=north east, end anchor=north west, bend left,yshift=-3pt]{r}{\tilde{x} x_{-l}}
			\&
			{M_{l-1}} 
			\ar[start anchor=south west, end anchor=south east,bend left,yshift=3pt]{l}{y_{-l} \tilde{x}^{-1}}
			\ar[start anchor=south east, end anchor=south west, bend right, yshift=3pt]{r}[swap]{x_{l}}
			\&
			{M_{l+1}}
			\ar[start anchor=north west, end anchor=north east, bend right, yshift=-3pt]{l}[swap]{y_{l}}
		\end{qrep}
	\end{align}
	It defines an isomorphism of functors $\EE_l \cong \EE'_l$. 
	
	\smallskip
	\noindent
	The 
	proof
	follows from the above claims combined with the description of the duality functor $\hcdual{(-)}$ in terms of diagrammatic presentations \eqref{E:DiagramHCmodule} of objects of $\HC_{\gamma(l)}(\lieg, K)$. 
\end{proof}
For further background on Harish-Chandra modules, we refer to Schmid \cite{Schmid}
and Wallach \cite{Wallach}.

\section{Elementary representation theory of the Gelfand quiver}
\label{sec:basic_rep}
In this section, we gather a few basic observations on finitely generated modules over the Gelfand order $\A$. 
In particular, we recall the description of lattices over the Gelfand order, introduce natural bases of its indecomposable projective modules,
determine all cyclic representations as well as all 
bricks in $\Rep{\A}$ up to isomorphism, and consider the smallest example of a one-parameter family $(M_{\lambda})_{\lambda \in \kk}$ of indecomposable representations.
Finally, we sketch the relationship of the Gelfand order to an associated gentle order first considered by Khoroshkin \cite{Khoroshkin}.

Cyclic representations will reappear as elementary ingredients of a gluing procedure in Section~\ref{sec:bases}, where their classification will be recovered in terms of projective resolutions.

\subsection{Lattices over the Gelfand order}
As before,
let $\kk$ be an algebraically closed
field, $R = \kk\llbracket t\rrbracket$, $\mx = R t$ and $\A$ be the Gelfand order \eqref{E:GelfandOrder}. Note that
the ring $R$ is isomorphic to the center of the matrix algebra $\A$
via the assignment $\begin{td} r \ar[mapsto]{r} \& \diag(r,r,r) \end{td}$ for any element $r \in R$.
\begin{dfn}
	A finitely generated $\A$-module $L$ is called a \emph{lattice} if it is free as $R$-module. 
\end{dfn}
We will frequently use the indecomposable $\A$-lattices 
\begin{align}\label{eq:lattices}
	P_\star = 
	\begin{pmatrix}
		R \\
		R \\
		R
	\end{pmatrix},&& 
	P_+ = 
	\begin{pmatrix}
		\idm \\
		R \\
		\idm
	\end{pmatrix}
	&&
	P_- = 
	\begin{pmatrix}
		\idm \\
		\idm \\
		R
	\end{pmatrix}
	\quad
	\text{and} \quad 
	I_\star = \rad P_\star = 
	\begin{pmatrix}
		\mx \\
		\Rx \\
		\Rx
	\end{pmatrix}.
\end{align}
For any vertex $v \in Q_0 \colonequals \{\star,+,-\}$ the projective module $P_v$ can be identified with $\A e_v$.

The natural inclusion of $I_\star$ into $P_\star$ will be denoted by $\iota$.
The modules $P_\star$, $P_+$ and $P_-$ are precisely the indecomposable projective $\A$-modules up to isomorphism.

\begin{rmk}
	Any indecomposable $\A$-lattice is isomorphic to one of the four lattices above.
	More precisely,
	the Gelfand order $\A$ admits a hereditary envelope $\B$, see \eqref{eq:emb2}, 
	such that any indecomposable $\A$-lattice is a direct summand of $\A \oplus \B$; see \cite[Theorem B and Example 4]{RingelRoggenkamp}.
	More generally,	any nodal order admits a hereditary envelope 
	(see \cite{Nodal, NodalOrders} for the definition),
	and a similar classification of indecomposable lattices holds true by 
	\cite{RingelRoggenkamp, MCM}.
\end{rmk}

\subsection{Bases of indecomposable projective modules}
\label{subsec:ind-proj}

We recall that the power series ring $\Rx$ has the topological basis $\{t^n \mid n\in \N_0 \}$ with respect to the $\mx$-adic topology. Based on this observation, the projective $\A$-module $P_{\star}$ admits a natural topological $\kk$-linear basis 
\begin{align*}
	P_\star 
	= \topspan{ t^n_v \mid v \in Q_0, n\in \N_0}
	\quad
	\text{ with }
	\quad
	t_\star^n = \begin{pmatrix}
		t^n \\ 0 \\ 0
	\end{pmatrix},
	\quad
	t_+^n = \begin{pmatrix}
		0 \\ t^n \\ 0
	\end{pmatrix},
	\quad
	t_-^n = \begin{pmatrix}
		0 \\ 0 \\ t^n
	\end{pmatrix}.
\end{align*}
where the closure of the $\kk$-linear span is taken with respect to the Jacobson radical of the Gelfand order $\A$. 
As the modules $P_+$ and $P_-$ embed into $P_\star$, each one admits a similar basis
\begin{align*}
	P_{\pm} = 
	\topspan{  t^0_{\pm}, t_v^p \mid v \in Q_0, p\in \N }.
\end{align*}
The quiver representations of the modules $P_-$, $P_\star$ and $P_+$ can be depicted as follows.
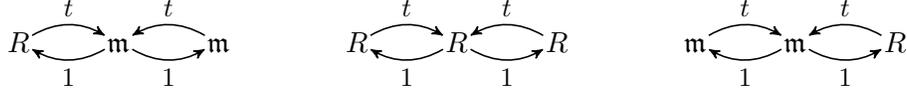
\begin{align}
	\begin{gqrep}{\mx}{\mx}{\Rx}{t}{t}{1}{1}\end{gqrep}&&
	\begin{gqrep}{\Rx}{\Rx}{\Rx}{t}{t}{1}{1}\end{gqrep}&&
	\begin{gqrep}{\Rx}{\mx}{\mx}{t}{t}{1}{1}\end{gqrep}
\end{align}
The $\kk$-linear bases of $P_-$, $P_\star$ and $P_+$ can be visualized by the infinite diagrams
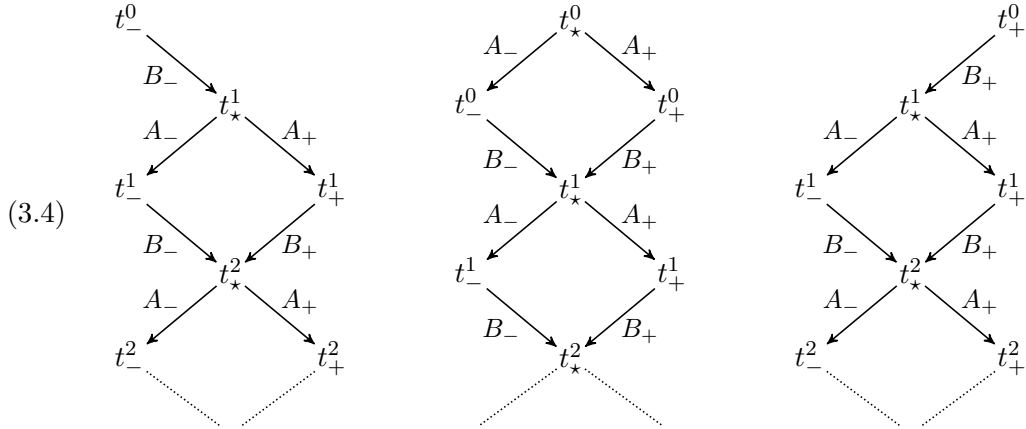
\begin{align}
	\label{eq:proj-diag}
	\begin{ld}
		t_-^0 \ar{rd}[swap]{B_-} \& \& \\
		\& t_\star^1 \ar{ld}[swap]{A_-} \ar{rd}{A_+} \& \\
		t_-^1 \ar{rd}[swap]{B_-} \& \& \ar{ld}{B_+} t_+^1 \\
		\& t_\star^2 \ar{ld}[swap]{A_-} \ar{rd}{A_+} \& \\
		t_-^2 \ar[densely dotted,-]{rd }\& \& t_+^2 \ar[densely dotted,-]{ld}\\
		\&  	\phantom{\bt} \& 
	\end{ld}
	&&
	\begin{ld}
		\& t_\star^0 \ar{ld}[swap]{A_-} \ar{rd}{A_+} \& \\
		t_-^0 \ar{rd}[swap]{B_-} \& \& \ar{ld}{B_+} t_+^0 \\
		\& t_\star^1 \ar{ld}[swap]{A_-} \ar{rd}{A_+} \& \\
		t_-^1 \ar{rd}[swap]{B_-} \& \& t_+^1 \ar{ld}{B_+}  \\
		\& t_\star^2  \ar[densely dotted,-]{rd} \ar[densely dotted,-]{ld} \&  \\
		\phantom{\bt} \& \& \phantom{\bt}
	\end{ld}
	&&
	\begin{ld}
		\& \& t_+^{0} \ar{ld}{B_+}\\
		\& t_\star^1 \ar{ld}[swap]{A_-} \ar{rd}{A_+} \& \\
		t_-^1 \ar{rd}[swap]{B_-} \& \& \ar{ld}{B_+} t_+^1 \\
		\& t_\star^2 \ar{ld}[swap]{A_-} \ar{rd}{A_+} \& \\
		t_-^2  \ar[densely dotted,-]{rd}  \& \& t_+^2 \ar[densely dotted,-]{ld} \\
		\& \phantom{\bt} \&
	\end{ld}
\end{align}
where $B_-$, $A_-$, $B_+$ and $A_+$ denote the matrices of the associated quiver representation.

\subsection{Cyclic modules and their resolutions}
\label{subsec:cyc-mod-res}
For a finitely generated module $M$ the intersection $\rad M$ of its maximal left submodules can be identified with $(\rad \A ) M$, where $\rad \A$ denotes the Jacobson radical of the ring $\A$.

Let us recall that a non-zero $\A$-module $M$ is called \emph{cyclic} if any of the following equivalent conditions holds.
\begin{enumerate}
	\item There is an element $v \in M$ such that $M = \A v$.
	\item The $\kk$-linear space $\head M = M/\rad M $ is one-dimensional.
	\item $M$ is isomorphic to a quotient of an indecomposable projective $\A$-module.
\end{enumerate}
To construct specific examples of cyclic $\A$-modules, for any symbol $\alpha 
\in Q_0 \cup \{\hat{\star}\}$, number $p \in \N_0$ and vertex $\beta \in Q_0$
we consider a complex of projective $\A$-modules
\begin{align}\label{eq:cyc-res}
	\CX(\alpha,p,\beta) = 
	\left\{
	\begin{array}{rl}
		\begin{cd} P_\star \ar{r}{
				\begin{psmallmatrix}
					\phantom{-}t \\ -t
				\end{psmallmatrix}
			} \& P_+ \oplus P_- \ar{r}{
				\begin{psmallmatrix}
					{t^{p}}&
					{t^{p}}	
			\end{psmallmatrix}	} \& P_{\beta}
		\end{cd} & \text{if }\alpha =\hat{\star}, \\
		\begin{cd} \mathstrut \& P_{\alpha} \ar{r}{t^{p}} \& P_{\beta}
		\end{cd}
		& \text{if }\alpha \in Q_0.
	\end{array} 
	\right.
\end{align}
In particular, we introduce the symbol $\hat{\star}$ to indicate a projective $P_\star$ at degree $2$ in $\CX(\alpha,p,\beta)$.
To ensure that the first differential is minimal, that is, that its image is contained in $\rad P_{\beta}$, we have also to require that $p \neq 0$ if $(\alpha,\beta) \notin \{(\hat{\star},\star), (+,\star),(-,\star)\}$.
It can be verified that the complex $\CX(\alpha,p,\beta)$ is a minimal projective resolution and its homology $M(\alpha,p,\beta) = H_0 (\CX(\alpha,p,\beta))$ is a finite-dimensional cyclic $\A$-module.

The projective resolution $\CX(\alpha,p,\beta)$ is quasi-isomorphic to 
the two-term complex 
$$\CL(\alpha,p,\beta) = 
\begin{cases}
	(
	\begin{cd} I_\star \ar{r}{t^{p}}\& P_{\beta}
	\end{cd})&	\text{if }\alpha = \hat{\star}\\
	(\begin{cd}  P_{\alpha} \ar{r}{t^{p}} \& P_{\beta}
	\end{cd}) &  \text{if }\alpha \in Q_0.
\end{cases}
$$
	The only non-trivial differential of the complex
	$\CL(\alpha,p,\beta)$  is given by the embedding of the first syzygy of $M(\alpha,p,\beta)$ into its projective cover $P_{\beta}$.
	It
	is called a \emph{syzygy resolution} of the module $M(\alpha,p,\beta)$ in Table~\ref{tab:cyc-base} on page~\pageref{tab:cyc-base}.

The monomial basis of $P_{\beta}$ induces 
a `standard' basis
of the cyclic $\A$-module $M(\alpha,p,\beta)$, which can be visualized by a finite subdiagram of the respective
basis diagram of $P_{\beta}$ in 
\eqref{eq:proj-diag}. 
Varying the possible triples $(\alpha,p,\beta)$, we obtain the twelve families of cyclic representations  $M(\alpha,p,\beta)$
listed in Table~\ref{tab:cyc-base}.

					In fact, the families of Table~\ref{tab:cyc-base} provide
					a complete classification of finite-dimensional cyclic $\A$-modules by the next statement.
						\begin{prp}\label{prp:fd-cyc}
							Any finite-dimensional cyclic $\A$-module is isomorphic to $M(\alpha,p,\beta)$ 
							for 
							certain $\alpha \in Q_0 \cup \{\hat{\star}\}$,  $\beta \in Q_0$ and $p \in \N_0$ such that $p \neq 0$ if $(\alpha,\beta) \notin \{(\hat{\star},\star), (+,\star),(-,\star)\}$.
						\end{prp}
						\begin{proof}
							Given a cyclic finite-dimensional $\A$-module $M$,
							there is a short exact sequence
							$$
							\begin{cd} 0 \ar{r} \& L 
								\ar{r} \& P_{\beta} \ar{r}{\pi} \& M \ar{r} \& 0\end{cd}
							$$
							where $\beta \in Q_0$ and $L = \ker \pi$ is a non-zero submodule of $P_{\beta}$.
							Since $\A$ is Noetherian and $P_{\beta}$ is finitely generated, so is $L$. Let $x$ be a non-zero element in a minimal set of generators of $L$. Let $a \in Q_0$ such that $e_a x \neq 0$. 
							Since $e_a L \subseteq e_a P_{\beta} \subseteq R$ and any power series $p \in R$ with $p(0) \neq 0$ is invertible, it holds that $e_a x = u t^n_a$ for a unit $u \in R$ and a number $n \in \N_0$.
							In particular, $L$ can be generated by monomials of the form $t^p_a$ with $a \in Q_0$ and $p \in \N_0$, where $p \neq 0$ if $\beta = \pm$. Given any two such monomials, the respective diagram of $P_{\beta}$ in \eqref{eq:proj-diag} shows that
							these are either given 
							by $t_+^p$ and $t_-^p$, or that 
							one of them is a right divisor of the other. 
							It follows that 
							$L = \A (t_+^p, t_-^p)$
							or $L = \A t_a^p$ for a vertex $a\in Q_0$. Thus, $M \cong M(\alpha,p,\beta)$.
						\end{proof}
						For a finite-dimensional $M$ we denote by $\ell \ell(M)$ its \emph{Loewy length}, that is, the smallest number $p \in \N$ such that $\rad^p M = 0$. The Loewy length can also be expressed as the nilpotency degree of the operator $t$ on $M$.
						The \emph{socle of $M$} will be denoted by $\soc M$ and is the submodule of  elements of $M$ annihilated by $\rad \A$.
						\begin{cor}
							Any cyclic finite-dimensional $\A$-modules $M$ and $M'$ are   isomorphic if and only if
							$\head M \cong \head M'$ and $\uldim M = \uldim M'$ if and only if 
							$\head M \cong \head M'$, $\soc M \cong \soc M'$ 
							and $\ell\ell(M) = \ell \ell(M')$.
						\end{cor}
						\begin{proof}
							Note that $\head M(\alpha,p,\beta) \cong S_\beta$ is determined by $\beta$.
							The first equivalence follows from the computations in Table~\ref{tab:cyc-base}.
							A closer look at the representations shows also that
							$\soc M(\alpha,p,\beta)$ is determined by $\alpha$, and $\ell\ell\bigl(M(\alpha,p,\beta)\bigr)$ by $p$, which implies the last equivalence.
						\end{proof}
						We recall that any cyclic $\A$-module is indecomposable by Nakayama's Lemma. 
						In particular, Table~\ref{tab:cyc-base} provides first examples of indecomposable objects of the category $\Rep{\A}$.

						\begin{rmk}
							Cyclic finite-dimensional $\A$-modules and finite-dimensional $O$-modules correspond precisely to those objects of $\HC(\lieg, K)$, which arise as Harish-Chandra modules associated to vector-valued polyharmonic weak Maa\ss{} forms; see \cite[Theorem 3.7]{ABR}.
						\end{rmk}
						
						\subsection{A preprojective quotient and bricks}
					The quotient $\bar{\A} \colonequals \A/(a_- b_-,a_+ b_+)$
					of the Gelfand order $\A$
					is the path algebra of the quiver with relations
					\begin{align*}
						\begin{qrep}
							{-}
							\ar[start anchor=north east, end anchor=north west, bend left,yshift=-5pt,
							""{name=a1, inner sep=0pt, outer sep=0pt, pos=0.66,swap}
							]{r}{b_-}
							\&
							{\star}   
							\ar[start anchor=south west, end anchor=south east,bend left,yshift=5pt,
							""{name=a2, inner sep=0pt, outer sep=0pt, pos=0.33,swap}
							]{l}{a_-}
							\ar[start anchor=south east, end anchor=south west, bend right, yshift=5pt,
							""{name=b2, inner sep=0pt, outer sep=0pt, pos=0.33}
							]{r}[swap]{a_+}
							\& 
							{+} 
							\ar[start anchor=north west, end anchor=north east, bend right, yshift=-5pt,
							""{name=b1, inner sep=0pt, outer sep=0pt, pos=0.66}
							]{l}[swap]{b_+}
						\end{qrep}
						&& 
						\begin{array}{l}
							b_- a_- = b_+ a_+,\\
							a_- b_- = a_+ b_+ = 0.
						\end{array}
					\end{align*}
					The quotient $\bar{\A}$ is precisely 
					\emph{the preprojective algebra of Dynkin type $\mathbb{A}_3$}.
					Its category $\Rep{\bar{\A}}$ of finite-dimensional modules 
					is well-known.
					\begin{lem}\label{lem:preproj}
						The Auslander-Reiten quiver of the category $\Rep{\bar{\A}}$
						is given by 
						
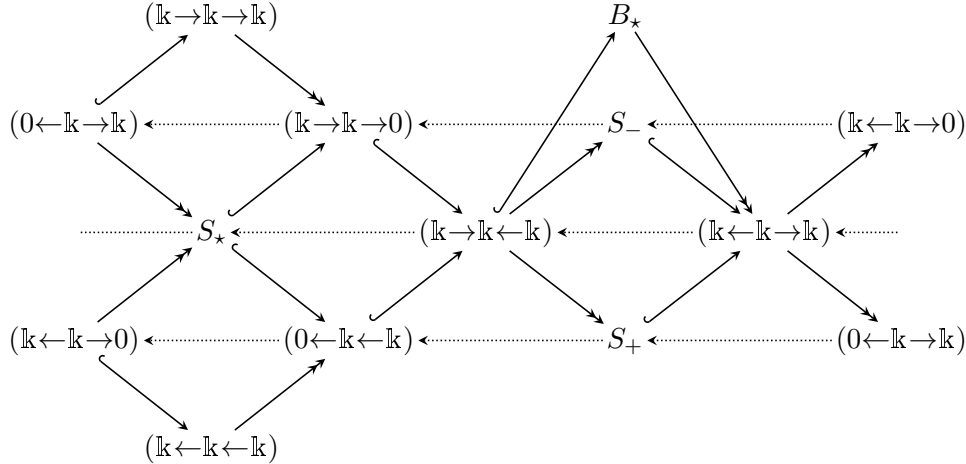
\begin{align*}
							\begin{tikzcd}[ampersand replacement=\&, outer sep=3em, cells={outer sep=0.75pt, inner sep=0.75pt}, 
								labels={rectangle, font=\small, inner sep=1.25pt, outer sep=1.5pt, minimum size=0cm},
								baseline=0pt,column sep = 0pt, row sep=1cm,
								arrows={->,semithick,>=stealth},
								]
								\&		    (\kk {\rightarrow} \kk { \rightarrow} \kk) \ar[twoheadrightarrow]{rd} \&  \& \& 
								B_\star
								\ar[twoheadrightarrow]{rdd}  \&  \\
								(0{\leftarrow}\kk {\rightarrow} \kk) \ar[hook]{ru} \ar[twoheadrightarrow]{rd} \&   \& 	\ar[densely dotted]{ll}	
								(\kk {\rightarrow} \kk {\rightarrow} 0)
								\ar[hook]{rd} \& \& S_- \ar[hook]{rd} \ar[densely dotted]{ll} \& \& (\kk {\leftarrow} \kk {\rightarrow} 0) \ar[densely dotted]{ll}
								\\
								\mathstrut \& S_\star \ar[densely dotted,-]{l} \ar[hook]{rd} \ar[hook]{ru} \& \& (\kk {\rightarrow} \kk {\leftarrow} \kk) \ar[densely dotted]{ll} \ar[twoheadrightarrow]{rd} \ar[twoheadrightarrow]{ru}
								\ar[hook]{ruu} \& 	\phantom{(\kk {\leftarrow} \kk {\leftarrow} \kk)} \& (\kk { \leftarrow} \kk { \rightarrow} \kk) \ar[twoheadrightarrow]{rd} \ar[twoheadrightarrow]{ru} \ar[densely dotted]{ll} \& \ar[densely dotted]{l} \mathstrut
								\\
								(\kk {\leftarrow} \kk  {\rightarrow} 0) \ar[hook]{rd} \ar[twoheadrightarrow]{ru} \&   \& 	
								(0 {\leftarrow} \kk {\leftarrow} \kk)
								\ar[hook]{ru} \ar[densely dotted]{ll} \& \& S_+ \ar[hook]{ru} \ar[densely dotted]{ll} \& \& (0 {\leftarrow} \kk {\rightarrow} \kk) \ar[densely dotted]{ll}
								\\
								\&		    (\kk { \leftarrow} \kk {\leftarrow} \kk) \ar[twoheadrightarrow]{ru} \&  \& \&   \& 
							\end{tikzcd}
						\end{align*}
						where we use abbreviations 
						\begin{align*}
										B_\star = \left[
										\begin{gqrep}{\kk}{\kk^2}{\kk}{\begin{psmallmatrix} 0 \\ 1 \end{psmallmatrix}}{\begin{psmallmatrix} 0 \\ 1 \end{psmallmatrix}}{\begin{psmallmatrix} 1 & 0 \end{psmallmatrix}}{\begin{psmallmatrix}1 & 0 \end{psmallmatrix}}
										\end{gqrep} \right],
										&&	
										(\kk {\rightarrow} \kk {\leftarrow} \kk) = \left[
										\begin{gqrep}{\kk}{\kk}{\kk}{\scriptstyle 1}{\scriptstyle 1}{\scriptstyle 0}{\scriptstyle 0}
										\end{gqrep}\right]
									\end{align*}
									and similar notation for the remaining modules.%
									\qedhere
								\end{lem}
								In particular, the number of isomorphism classes of indecomposable $\bar{\A}$-modules is precisely the number of roots of $\mathbb{A}_3$.
								We omit the proof of the last statement as it is a standard application of the covering technique in representation theory of quivers going back to work by Bongartz and Gabriel \cite{BG}.
							\begin{dfn}An object $X$ of 
								an additive category 
								$\mathcal{C}$
								is
								a \emph{brick} if
								$\End_{\mathcal{C}}(X)$ is a division ring. 
							\end{dfn}

							\begin{prp}\label{P:Quasi-simple} Up to isomorphism, the 
								bricks of the category $\Rep{\A}$ are
								given by the representations of Lemma~\ref{lem:preproj} without the representation $B_\star$.
							\end{prp}
							\begin{proof}
								It is straightforward to check that each of the eleven representations in question is a brick.
								Vice versa, let $M \in \Rep{\A}$ be 
								a brick.
								Then the left multiplication $t \cdot \colon M \to M $ defines a non-invertible endomorphism, which implies that $t M = 0$.
								Hence, $M$ is a module over the finite-dimensional $\kk$-algebra $\A/t\A$.
									Since $\A/t \A
									\cong
									\bar{\A}/\bar{\A}(b_- a_-)\bar{\A}$,
									Lemma~\ref{lem:preproj} implies the claim.
							\end{proof}

							\smallskip
							\noindent
							Theorem \ref{T:blocks} and Proposition \ref{P:Quasi-simple} recover the
							following statement.
							\begin{cor}
								An object $M$ of the $\CC$-linear category $\HC(\lieg, K)$ 
								has one-dimensional endomorphism ring
								if and only if the following hold:
									$M$ belongs to one of the blocks from Theorem~\ref{T:blocks}, 
									$\dim_{\CC} M_n \leq 1$ for any $n \in \Z$, and 
									for any $l \in \Z$ with $x_l =y_l =0$ it follows that $M_{l-1} = 0$ or $M_{l+1} = 0$.
								\qed
							\end{cor}

							\subsection{A one-parameter family of indecomposable representations}

							The simplest classes of indecomposable $\A$-modules, 
							cyclic modules and bricks,
							are described by discrete invariants.
							Next, we consider a family of modules depending on a parameter from the base field.
							The following construction produces 
							a Gelfand quiver representation
							from  any finite-dimensional representation of the quiver of type 
							$\widetilde{\mathbb{D}}_4$ with two sources and two sinks.
							\begin{align}
								\nonumber
								&
								\begin{ld}
									U_1 \ar{rd}[swap]{B} \& \& W_1 \ar{ld}{D}\\
									\&		V  \ar{rd}{C}  \ar{ld}[swap]{A} \\	
									U_2  \& \& W_2
								\end{ld}
								\qquad
								\longmapsto \qquad
								\left[
								\begin{gqrep}{W_1\oplus W_2}{V}{U_1 \oplus U_2}{
										\begin{psmallmatrix} D& 0
										\end{psmallmatrix}
									}{	
										\begin{psmallmatrix}
											B & 0
										\end{psmallmatrix}
									}{	
										\begin{psmallmatrix}
											0 \\ C
										\end{psmallmatrix}
									}{	\begin{psmallmatrix}
											0
											\\ A
										\end{psmallmatrix}
								}\end{gqrep}\right]
								\\
							\intertext{We may use this to construct a family of representations of the Gelfand quiver}
							\label{eq:sh-band}
							&
							\begin{ld}
								\kk \ar{rd}[swap]{
									\begin{psmallmatrix}	
										0 \\ 1
									\end{psmallmatrix}
								} \& \& \kk \ar{ld}{
									\begin{psmallmatrix}	
										1 \\0 
									\end{psmallmatrix}
								}\\
								\&		\kk^2 \ar{rd}{
									\begin{psmallmatrix}	
										1 & 1 
									\end{psmallmatrix}
								}  \ar{ld}[swap]{
									\begin{psmallmatrix}	
										\lambda & 1
									\end{psmallmatrix}
								} \\	
								\kk  \& \& \kk
							\end{ld}
							\qquad
							\longmapsto \qquad
							M_\lambda = \left[
							\begin{gqrep}{\kk^2}{\kk^2}{\kk^2}{\begin{psmallmatrix}
										1 & 0 \\
										0 & 0
								\end{psmallmatrix}}{\begin{psmallmatrix}
										0 & 0 \\
										1 & 0
								\end{psmallmatrix}}{\begin{psmallmatrix}
										0 & 0 \\
										1 & 1
								\end{psmallmatrix}}{\begin{psmallmatrix}
										0 & 0\\
										\lambda & 1 
							\end{psmallmatrix}}\end{gqrep}\right]
						\end{align}
						indexed by a scalar $\lambda \in \kk$.
						Straightforward computations show that
						\begin{enumerate}
							\item $\End_{\A}(M_\lambda) \cong \kk[\phi, \psi]/(\phi^2, \phi \psi, \psi^2)$ is a local $\kk$-algebra, 
							\item $M_{\lambda} \not\cong M_{\mu}$ for any $\mu \in \kk$ with $\mu \neq \lambda$. 
						\end{enumerate}
						In particular, $(M_\lambda)_{\lambda \in \kk}$ forms a family of pairwise non-isomorphic indecomposable representations.
						Hence,
						the category $\Rep{\A}$ is not representation-discrete. 
						
						Moreover, the following holds.
						\begin{itemize}
							\item 
							If $\lambda \neq 0$, 
							there is an isomorphism $\begin{td}\inv(M_\lambda) \ar{r}{\sim} \& M_{\lambda} \end{td}$
							given by $(\lambda\Id, \begin{psmallmatrix} 0 & 1 \\ \lambda & 0 \end{psmallmatrix},  \Id)$.
							\item
							If $\lambda = 0$, $A_+ B_- = 0$ is not similar to $A_- B_+ \neq 0$, and thus $\inv(M_{0}) \not\cong M_0$. 	
							\item It holds that $\lambda = 1$ if and only if 
							the socle of $M_{\lambda}$ contains $S_\star$.
						\end{itemize}
						These observations show that
						the family $(M_\lambda)_{\lambda \in \kk}$ of representations splits naturally into three homologically distinct types:
						the generic case $M_{\lambda}$ with $\lambda \in \kk\backslash\{0,1\}$, and the two degenerate cases
						$M_0$ and $M_1$.

						\subsection{The gentle order associated to the Gelfand order}
						
						Let $B$ denote the arrow ideal completion of the $\kk$-linear path algebra of the  quiver with relations
						\begin{align}\label{E:Auslander}
							\begin{qrep}
								{-}
								\ar[start anchor=north east, end anchor=north west, bend left,yshift=-5pt,
								""{name=a1, inner sep=0pt, outer sep=0pt, pos=0.66},<-
								]{r}{d_-}
								\&
								{\bt}   
								\ar[start anchor=south west, end anchor=south east,bend left,yshift=5pt,
								""{name=a2, inner sep=0pt, outer sep=0pt, pos=0.33},<-
								]{l}{c_-}
								\ar[start anchor=south east, end anchor=south west, bend right, yshift=5pt,
								""{name=b2, inner sep=0pt, outer sep=0pt, pos=0.33,swap}
								]{r}[swap]{d_+}
								\& 
								{+} 
								\ar[start anchor=north west, end anchor=north east, bend right, yshift=-5pt,
								""{name=b1, inner sep=0pt, outer sep=0pt, pos=0.66,swap}
								]{l}[swap]{c_+}
							\end{qrep}
							&&
							d_- c_+ = d_+ c_- = 0.
						\end{align}
						The algebra $B$ can be identified with the matrix algebra
						$$
						\left\{\left.
						\left(
						\left(
						\begin{array}{cc}
							a_{11} & a_{12} \\
							a_{21} & a_{22}
						\end{array}
						\right), 
						\left(
						\begin{array}{cc}
							b_{11} & b_{12} \\
							b_{21} & b_{22}
						\end{array}
						\right)
						\right) \in O \times O \, \right| \,  a_{22}(0) = b_{22}(0) 
						\right\},
						$$
						where $O$ is the algebra defined by \eqref{E:Hereditary}. It follows that $B$ is  a nodal order of \emph{gentle type}. 
						A classification of indecomposable objects in $\Rep{B}$ can be obtained in a similar way to the case of the Zhelobenko quiver \eqref{E:ZhelobQuiver}; see \cite{GP, Ringel}.
						
						For the remainder of this section, we assume that  $\chr{\kk} \neq 2$.
						In this case, the representation theory of the algebra $B$ is closely related to that of the Gelfand order via a skew group algebra as intermediate link.
						
						Recall that the cyclic group $\Ctwo = \bigl\langle \invaut \,\big| \, \invaut^2 = \id \bigr\rangle$ acts on the Gelfand order $\A$ by the rules~\eqref{E:Involution}.
						This gives rise to the  skew group algebra
						$ \ACtwo = \{ a\,\id  + b\,\invaut \mid a,b \in \A\}$
						with the multiplication defined 
						by $a\,g \cdot b\,h = a \,g(b)\,(g\circ h)$
						for any $a,b \in \A$ and any $g,h \in \Ctwo$.
						\begin{lem}\label{L:SkewProduct}
							In the setup above,
							the algebra $\ACtwo$ is Morita equivalent to the algebra $B$.
						\end{lem}
						
						\begin{proof} The algebra $\ACtwo$ has the primitive idempotents $e_\pm$ and 
							$e_\star^\pm = \frac{1}{2} \cdot (e_\star \id \pm e_\star \invaut)$. Since $\invaut(e_\pm) = e_\mp$, we have the following identity in $\ACtwo$:
							$
							\invaut e_\pm \id = e_\mp \invaut.
							$
							Hence, the idempotents $e_+$ and $e_-$ are conjugated in $\ACtwo$. Let $e = e_+ + e_{\star}^+ +  e_\star^-$. Then the algebras $\ACtwo$ and $\ACtwo^b =  e (\ACtwo) e$ are Morita equivalent
							and satisfy 	$\rad (\ACtwo^b) = e (\rad \ACtwo) e$	 by \cite[Proposition~5.13]{CurtisReiner}.
							It follows that there is an isomorphism $\begin{td} B  \ar{r}{\sim} \& \ACtwo^b \end{td}$
							such that
							the primitive idempotents of $B$ at vertices $-$, $\bt$ and $+$ are mapped to $e_\star^-$, $e_+$ and $e_\star^+$, respectively, the arrows
							$c_\pm$ to $a_+ e_\star^{\pm}$, and $d_\pm$ to $e_\star^{\pm} b_+$.
						\end{proof}

						The natural embedding $\begin{td} \A \ar{r} \& \ACtwo \end{td}$ gives rise to a restriction functor
						\begin{align*}
							\begin{cd}
								\res{(-)} \colon \Rep{\ACtwo}  \ar{r} \& \Rep{\A}
							\end{cd}
						\end{align*}
						Moreover, there is an involution on $\ACtwo$ interchanging 
						signs in its paths, which induces an involutive auto-equivalence
						$\widetilde{\inv}$ on $\Rep{\ACtwo}$. 
						Work by Khoroshkin \cite{Khoroshkin} and Reiten-Riedtmann \cite{ReitenRiedtmann}
						imply the following result.
						\begin{thm}
							In the setup above, the following statements hold.
							\begin{enumerate}
								\item Any indecomposable object $M$ in $\Rep{\A}$	is isomorphic to a direct summand of the restriction $\res{N}$ of an indecomposable object $N$ from $\Rep{\ACtwo}$.
								\item
								For any indecomposable object $N$ from $\Rep{\ACtwo}$ the following holds.
								\begin{itemize}
									\item If $N$ is not $\widetilde{\inv}$-invariant, then $\res{N}$ is indecomposable, $\inv$-invariant and $\res{N} \cong \res{(\widetilde{\inv}(N))}$.
									\item If $N$ is $\widetilde{\inv}$-invariant, then $\res{N} \cong M \oplus \inv(M)$ for an indecomposable non-$\inv$-invariant object $M$ from $\Rep{\A}$ .
								\end{itemize}
							\end{enumerate}
						\end{thm}
						Together with the previously stated classification results on $\Rep{B}$, this theorem  
						can be used
						to
						conclude that the category $\Rep{\A}$ is representation-tame if 
						$\chr{\kk} \neq 2$.
						While
						the decomposition of $\res{N}$
						can be determined for any individual $\widetilde{\inv}$-invariant indecomposable object $N$ in $\Rep{\ACtwo}$ using the methods  in \cite{Khoroshkin}, describing the decompositions uniformly for all such $N$ is a non-trivial task.

						\begin{rmk} 
							If $\kk = \CC$,
							the  category $\Rep{B}$, and thus the category $\Rep{\ACtwo}$, is equivalent to the principal block of the category of Harish-Chandra modules for the reductive group $\mathsf{PGL}_2(\RR)$; see \cite[Section 6]{BGG}.
						\end{rmk}

\section{Basic homological algebra of the Gelfand quiver}
\label{sec:basic_homological}

This section collects some homological properties of the derived category $\DbRep{\A}$ of nilpotent representations of the Gelfand order.
We observe that  $\A$ is a self-opposite ring of  global dimension two, introduce an elementary duality endofunctor $\dual{(-)}$ on $\DbRep{\A}$, and relate it to the derived versions of the involution $\inv$ and the contragredient duality $\cdual{(-)}$ introduced before.
Then we provide a factorization of the derived Auslander-Reiten translation on $\DbRep{\A}$ in terms of a spherical 
auto-equivalence $\Tw$ associated to the simple module $S_\star$ (Proposition~\ref{P:twist}). 
This result is crucial to 
prove that $\tau$-periodicity, $\tau^2$-invariance and vanishing defect are all equivalent conditions (Proposition~\ref{prp:tau-per},  Result~\ref{res:B}~\eqref{B1}). 
This section prepares the ground for the explicit description of the functors on the derived category in Section~\ref{sec:derived_fun}.

\subsection{A duality preserving simple modules}
\label{subsec:dual-sim}
As in the previous sections, let
$\Rep{\A}$ denote the category of finite-dimensional  left modules of the Gelfand order $\A$.

Next, we observe 
that the Gelfand order $\A$ is \emph{self-opposite},
that is,
we describe an isomorphism 
$(\Aop,\ast)  \cong (\A,\,\cdot\,)$ of $\kk$-algebras.
This observation will lead to the definition of a 
duality functor
$\begin{td} \dual{(-)}\colon \Rep{\A} \ar{r}{\sim} \& \Rep{\A} \end{td}$.
The opposite algebra $\Aop$ 
of the Gelfand order
may be identified with 
arrow ideal completion of the path algebra of the quiver with relations
\begin{align}\label{E:GelfandQuiverOpp}
	\begin{gqrep}{+}{\star}{-}{a_+^\circ}{a_-^{\circ}}{b_+^\circ}{b_-^\circ}
	\end{gqrep} && a^\circ_{-}\ast b^\circ_{-}  =  a^\circ_{+}\ast b^\circ_{+}.
\end{align}
The last quiver can be identified with the Gelfand quiver
via $\begin{td} b^\circ_{\pm} \ar[mapsto]{r} \& a_{\pm} \end{td}$ and $\begin{td} a^{\circ}_{\pm} \ar[mapsto]{r} \& b_{\pm} \end{td}$. 
This identification can be used to define
an isomorphism
$\begin{td} \A \ar{r}{\sim} \& \Aop \end{td}$.
In equivalent terms, we may consider the
$\kk$-linear anti-automorphism $\begin{td} \varphi \colon \A \ar{r}{\sim} \& \A \end{td}$ which 
is uniquely determined by the prescriptions
\begin{align}
	\label{eq:anti-inv}
	\varphi(e_\star) = e_\star,
	\quad	\varphi(e_\pm) = e_\pm, 
	\quad
	\varphi(a_\pm) = b_\pm \quad \mbox{\rm and} \quad
	\varphi(b_\pm) = a_\pm.
\end{align}
The restriction along the isomorphism $\varphi$ gives rise to an equivalence
\label{eq:Phi}
$\phi \colon
\begin{td}\rep{\Aop} \ar{r}{\sim} \&\rep{\A}\end{td}$.
Composing this equivalence 
with the $\kk$-linear duality 
$\begin{td} \Hom_{\kk}(\,-\,, \kk)\colon  \Rep{\A} \ar{r}{\sim} \& \Rep{\Aop} \end{td}$
yields a contravariant equivalence
\begin{align}
	\dual{(-)} \colonequals \phi \circ \Hom_{\kk}(\,-\,, \kk) \colon \begin{cd}\rep{\A} \ar{r}{\sim} \& \rep{\A}\end{cd}
\end{align}
Given an object in the category $\Rep{\A}$
its image under this duality is obtained by taking the $\kk$-linear dual of each vector space and each linear map.

\begin{align}\label{E:dual}
\begin{gqrep}{V_+}{V_\star}{V_-}{B_+}{B_-}{A_+}{A_-}
\end{gqrep}
\quad
\overset{\dual{(-)}}{\longmapsto} \quad
\begin{gqrep}{V_+^*}{V_\star^*}{V_-^*}{A_+^*}{A_-^*}{B_+^*}{B_-^*}
\end{gqrep}
\end{align}
In particular, the functor $\dual{(-)}$ is a contravariant, involutive equivalence which preserves each simple $\A$-module up to isomorphism.
For any objects $X, Y$ of $\Rep{\A}$ and
$p \in \Z$ 
we have: 
$$\Ext^p_{\A}(X, Y) \cong \Ext^p_{\A}(\dual{Y}, \dual{X}).$$
The duality functor $\dual{(-)}$
will be useful to determine socle and injective dimension of $X$ using the observations
\begin{align*}
\begin{array}{c}
	\soc(X) \cong \head(\dual{X}),
\end{array}
&&
\begin{array}{c}
	\injdim (X)  = \prdim (\dual{X}).
\end{array}
\end{align*}
\begin{rmk}
There is an isomorphism of functors
$$
\begin{cd}  \Hom_{\kk}(\,-\,, \kk) \cong
	\Hom_{\Rx}(-,E_{\kk}) \colon 
	\Rep{\A} \ar{r}{\sim} \& \Rep{\Aop} \end{cd}
$$
where $E_{\kk}$ denotes the injective hull of the field $\kk$ considered as an $\Rx$-module.
Therefore, $\dual{(-)}$ may be viewed as an `endofunctor version' of Matlis duality for $\rep \A$.
\end{rmk}

\subsection{Basic facts on the derived category of the Gelfand order}
\label{subsec:basic-der}

\begin{lem}\label{L:dimensions}
The simple objects of the category $\Rep{\A}$ have the projective and injective dimensions
\begin{equation}
	\prdim(S_\pm) = 1 = \injdim(S_\pm) \quad \mathrm{and} \quad \prdim(S_\star) = 2 = \injdim(S_\star).
\end{equation}
In particular, the global dimension of each of the abelian categories $\Rep{\A}$, $\md\A$ and $\Md\A$ is equal to two.
\end{lem}
\begin{proof}
A minimal projective resolution of $S_\pm$ 
is given by 
$(\begin{td}P_\star \ar{r}{t} \& P_{\pm}
\end{td})$
and one of $S_\star$ by
\begin{align}\label{E:projresSst}
	\begin{cd}
		P_\star \ar{r}{
			\left(\begin{smallmatrix}
				\phantom{-}t\\
				-t
			\end{smallmatrix}\right)	
		} \& P_+ \oplus P_- \ar{r}{\left(\begin{smallmatrix} \iota & \iota \end{smallmatrix}\right)}\& P_\star
	\end{cd}
\end{align}
This implies the claims on the projective dimension. Since the duality functor $\dual{(-)}$ \eqref{E:dual} preserving each simple module up to isomorphism, we have:
$\prdim(S_v) = \injdim(S_v)$ for any $v  \in \{\star,+, -\}$. The result on the global dimensions follows from the equalities $\gldim\bigl(\Rep{\A}\bigr) = \gldim(\md \A) = \gldim(\Md\A) =
\prdim_{\A}(\A/\rad \A)$, 
see for example  \cite[Proposition 2.2]{IR}.
\end{proof}

Let $\Md \A$  denote the category 
of left $\A$-modules,
and  $\md \A$ its full subcategory given by finitely generated modules.

\begin{prp}
There are exact equivalences  of triangulated categories
\begin{align*}
	\begin{cd}
		\DbRep{\A} \ar{r}{\sim} \& \Dbfdmod{\A} \ar{r}{\sim} \& \DbfdMod{\A}
	\end{cd} 
\end{align*}
where $\Dbfdmod{\A}$  denotes the full subcategory
of the derived category $\Dbmod{\A}$  consisting of complexes with finite-dimensional homology, and 
$\DbfdMod{\A}$ is defined similarly.
\end{prp}

\noindent
\begin{proof}[Comment on the proof]
The result follows from the fact that $\A$ is finitely generated as module over its center $R$, see for example \cite[Lemma 2.5]{IR}.
\end{proof}

Let $\proj \A$ be the category of finitely generated projective $\A$-modules and
$\Hotbproj{\A}$ 
the bounded homotopy category of complexes of modules in $\proj \A$.
Since $\A$ has finite global dimension,
the category $\Hotbproj{\A}$ is equivalent to the derived category $\Dbmod{\A}$.
The full subcategory $\Hotbfdproj{\A}$
given by complexes in $\Hotbproj{\A}$ with finite-dimensional homology 
is equivalent to $\DbRep{\A}$.

\subsection{Relationship of the involution and the dualities}
\label{subsec:fun3}
We recall that there 
is an involution on the Gelfand quiver
interchanging $+$ and $-$ in vertices and arrows,
which induces an auto-equivalence $\begin{td} \inv \colon \Rep{\A} \ar{r}{\sim} \& \Rep{\A} \end{td}$ of order two.
The latter extends further to an auto-equivalence
of $\DbRep{\A}$
which will be denoted again by $\inv$.
\begin{rmk}
The action of $\inv$ on a projective complex from $\Dbfdmod{\A}$ is determined by $\inv(P_{\pm}) \cong P_{\mp}$, $\inv(P_\star) \cong P_\star$
and $\inv(\begin{td} P_u \ar{r}{t^n} \& P_v \end{td}) \cong (\begin{td} \inv(P_{u}) \ar{r}{t^n} \& \inv(P_v) \end{td})$ for any $u,v \in \{\star,+,-\}$ and $n \in \N_0$
with $n \neq 0$ if $(u,v) \neq (\pm,\star)$.
\end{rmk}

We recall that the \emph{contragredient duality} on an appropriate block of Harish-Chandra modules
$\HC_{\gamma(l)}(\lieg, K)$
gives rise to a contravariant equivalence $\begin{td} \cdual{(-)} \colon \Rep{\A} \ar{r}{\sim} \& \Rep{\A} \end{td}$
whose extension to $\DbRep{\A}$ will again be denoted by the same notation.
Similarly, we extend the duality $\dual{(-)}$ from $\Rep{\A}$ to 
$\DbRep{\A}$ as well.

The explicit descriptions of the involutive functors $\cdual{(-)}$, $\inv$,  and $\dual{(-)}$ 
on representations
in \eqref{E:reflection}, \eqref{E:LieDuality} and \eqref{E:dual}
imply that 
the composition of any two of these functors is equal to the third one.
In particular, it holds that 
\begin{align}\label{eq:duals}
\begin{cd}
	\cdual{(-)} = \dual{(-)} \circ \inv = \inv \circ \dual{(-)} \colon 
	\DbRep{\A} \ar{r}{\sim} \& \DbRep{\A}
\end{cd}.  
\end{align}
Therefore, it will be sufficient to give explicit descriptions of the functors $\inv$ and $\dual{(-)}$ 
on indecomposable objects in order to describe the contragredient duality on
$\DbRep{\A}$.

\subsection{Derived Auslander-Reiten translation}
\label{subsec:der-AR}
Recall that the dualizing $\A$-bimodule  $\Omega$ has the following description:
\begin{equation}
\begin{cd}
	\Omega \colonequals \Hom_{\Rx}(\A, \Rx) \cong
	\begin{pmatrix}
		\Rx & \Rx &  \Rx \\
		\Rx t^{-1} & \Rx & \Rx t^{-1} \\
		\Rx t^{-1}& \Rx t^{-1}& \Rx
	\end{pmatrix}
	\ar{r}{\cdot t}[swap]{\sim}\&
	\begin{pmatrix}
		\Rx t & \Rx t & \Rx t \\
		\Rx & \Rx t & \Rx \\
		\Rx & \Rx & \Rx t
	\end{pmatrix},
\end{cd}
\end{equation}
where the left and right actions of $\A$ are given by the matrix multiplication.
The next result is a special case of \cite[Lemma 6.4.1]{vdBergh} or
\cite[Theorem 3.7]{IR}.
To simplify notation of morphism spaces, we abbreviate $\D(\A) = \D(\Md \A)$ below.
\begin{thm} The functor $\tau = \Omega \stackrel{\mathbb{L}}\otimes \,-\,: \begin{td} \DbRep{\A} \ar{r} \&
	\DbRep{\A} \end{td}$ {\rm{(}}called \emph{Auslander--Reiten translation}{\rm{)}} is an auto-equivalence. Moreover, for any pair of objects
$\CX, Y_{\bu}$ of $\DbRep{\A}$ we have a bifunctorial isomorphism
\begin{equation}\label{E:SerreDual}
	\bigl(\Hom_{\D(\A)}(\CX, Y_{\bu})\bigr)^* \cong
	\Hom_{\D(\A)}(Y_{\bu}, \tau(\CX)[1]),
\end{equation}
where $(-)^\ast$ denotes the $\kk$-linear dual $\Hom_{\kk}(-,\kk)$.
In other words, $\SS = \tau \circ [1]$ is a  Serre functor of the category $\DbRep{\A}$.
\end{thm}

\begin{rmk}\label{R:tauaction}
The action of $\tau$ on projective complexes can be described very explicitly. In fact, $\tau$ is also an auto-equivalence
of $\Dbmod{\A}$ with the following action on the indecomposable projective $\A$--modules:
$\tau(P_\pm) \cong  P_\mp$ and 
$\tau(P_\star) \cong  I_\star$.
Moreover, the action of $\tau$ on morphisms between projective $\A$-modules is determined by   the rules:
\begin{itemize}
	\item $\bigl( \begin{td} P_{\pm} \ar{r}{t^n} \& P_{\mp} \end{td} \bigr) \longmapsto \bigl( \begin{td} P_{\mp} \ar{r}{t^n} \&  P_{\pm} \end{td}\bigr)$ for any $n \in \N$,
	\item $\bigl( \begin{td} P_\star  \ar{r}{t^n} \& P_\pm \end{td} \bigr) \longmapsto \bigl( \begin{td} I_\star \ar{r}{t^n} \& P_\mp \end{td} \bigr)$ for any $n \in \N$,
	\item $\bigl( \begin{td} P_\pm \ar{r}{t^n} \&  P_\star \end{td} \bigr) \longmapsto \bigl( \begin{td} P_\mp \ar{r}{t^n} \&  I_\star \end{td} \bigr)$ for any $n \in \N_0$.
\end{itemize}
\end{rmk}

\begin{lem}\label{L:SimpleSpherical}
For the simple module $S_\star$ we have
\begin{equation}
	\Ext^p_{\A}(S_\star, S_\star) \cong
	\begin{cases}
		\kk & \text{if } p \in \{0,2\}, \\
		0 & \text{otherwise},
	\end{cases}
	\quad \mbox{\rm and} \quad
	\tau(S_\star) \cong S_\star[1].
\end{equation}
In other words, $S_\star$ is a \emph{two-spherical object} in the triangulated category  $\DbRep{\A}$.
\end{lem}
\begin{proof}
The result on the $\Ext$ spaces follows immediately from the projective resolution \eqref{E:projresSst}. To show the second result note that the complex $\tau(S_\star)$ 
is given by
$$
\begin{cd}
	I_\star  \ar{r}{
		\left(\begin{smallmatrix}
			\phantom{-}t\\
			-t
		\end{smallmatrix}\right)	
	} \& P_- \oplus P_+ \ar{r}{\left(\begin{smallmatrix} \iota & \iota \end{smallmatrix}\right)}\& I_\star 
\end{cd} 
$$ 
and has only one non-vanishing homology, 
which is given by 
$S_\star$ at degree one.
\end{proof}

\begin{prp}\label{P:BoundsDimensions}
For any object  $X$  of $\Rep{\A}$  we have: $\prdim(X), \injdim(X) \in \{1, 2\}$. Moreover, 
\begin{itemize}
	\item[(a)] $\prdim(X) = \,2$  if and only if $\Hom_{\A}(S_\star, X) \ne 0$.
	\item[(b)] $\injdim(X) = 2$ if and only if $\Hom_{\A}(X, S_\star) \ne 0$.
\end{itemize}
\end{prp}
\begin{proof}
Let
$
(\begin{td}
	F_2 \ar{r} \&   F_1 \ar{r} \&  F_0  \end{td})
$
be a minimal projective resolution of $X$. Clearly, $F_2 \in \add(P_\star \oplus P_+ \oplus P_-)$. However,
according to Lemma~\ref{L:dimensions} we have: $\Hom_{\A}(F_2, S_\pm) \cong \Ext^2_{\A}(X, S_\pm) = 0$. Therefore,
$F_2 \cong P_\star^{m}$ for $m = \dim \Ext^2_{\A}(X, S_\star)$. By Lemma~\ref{L:SimpleSpherical} and Serre duality
\eqref{E:SerreDual} we have:
$
\Ext^2_{\A}(X, S_\star)^\ast \cong \Hom_{\A}(S_\star, X),
$
hence the result on the projective dimension of $X$. The claim about the injective dimension follows from  the facts that    $\dual{S_\star} = S_\star$ and $\injdim(X) = 
\prdim\bigl(\dual{X}\bigr)$.
\end{proof}
The last proposition implies the following purely Lie-theoretic consequence.
\begin{cor} \label{cor:hc-dim} For any   object  $M$ of $\HC(\lieg, K)$ we have: $\prdim(M), \injdim(M) \in \{1, 2\}$. Moreover, 
\begin{itemize}
	\item[(a)] $\prdim(M) = \,2$  if and only if $M$ contains a finite-dimensional subobject.
	\item[(b)] $\injdim(M) = 2$ if and only if $M$ admits a finite-dimensional quotient.
\end{itemize}
\end{cor}

\begin{proof} The restrictions on the projective and injective dimensions of $M$ follow from the block decomposition \eqref{E:blocks} and the facts that $\mathsf{gl.dim}(R) = 1 = \mathsf{gl.dim}(O)$ and $\mathsf{gl.dim}(\A) = 2$. Recall that indecomposable finite-dimensional objects of $\HC(\lieg, K)$ are determined by their dimensions. Let $F_l$ be an indecomposable object of $\HC(\lieg, K)$ of dimension $l$. Then 
$F_l$ belongs to the block $\HC_{\gamma(l)}(\lieg, K)$ and $\EE_l(F_l) = S_\star$; see \eqref{E:FiniteDimSimple}. Hence, the statements about projective and injective dimension of $M$ follow from Proposition~\ref{P:BoundsDimensions}.
\end{proof}

According to work by Seidel and Thomas \cite{SeidelThomas}, the spherical object  $S_\star$ induces
an auto--equivalence $\Tw$ of $\DbRep{\A}$ called \emph{dual spherical twist}  given by
\begin{equation}
\Tw(\CX) \colonequals
\mathsf{Cone}\Bigl(
\begin{td}
	\CX  \ar{r}{\mathsf{ev}} \&
\Hom^{\bu}_{\kk} \bigl(\Hom^{\bu}_{\A}(\CX, S_{\star}),S_{\star} \bigr)  \end{td} \Bigr) [-1].
\end{equation}

\begin{prp}\label{P:twist}
We have the following isomorphisms of endo-functors of $\DbRep{\A}$:
$$
\tau 
\simeq \inv \circ 	\Tw \simeq 	\Tw \circ \inv.
$$
\end{prp}

\begin{proof}
According to \cite[Theorem 6.14]{IR}, the dual twist functor admits the following description:
$\Tw \simeq \Upsilon_{\bu} \stackrel{\mathbb{L}}\otimes \,-\,$, where $\Upsilon_{\bu}$ is the complex
of $\A$--bimodules
$$
\Upsilon_{\bu} = \bigl(\begin{cd} \dots \ar{r} \&  0\ar{r} \& \A
\ar{r}{\gamma} \& \End_{\kk}(S_\star) \ar{r} \&  0 \ar{r} \&  \dots \end{cd}\bigr).
$$
Here, $\A$ is located at  the zero degree and $\gamma$ is the canonical map. It is clear that $\Upsilon_{\bu} \cong \inv(\Omega)[0]$,
which implies  the statement.
\end{proof}
\begin{rmk}
Proposition~\ref{P:twist} provides a description of the Auslander--Reiten translation $\tau$ in other incarnations
of the category $\DbRep{\A}$ like $\DbHC$.
Moreover,
an analogous result holds
for the
derived category of the principal block of Harish-Chandra modules over the 
Lorentz group $\SO(n,1)$ or its identity component $\SO^+(n,1)$, for any $n \geq 1$ \cite[Theorem~2.2.3]{Gnedin}. 
\end{rmk}

\subsection{The defect and $\tau$-periodic objects}

\begin{dfn}
For any object $\CX$ of $\DbRep{\A}$ we define its \emph{defect} $\delta(\CX)$ by the formula
\begin{equation}\label{E:DefectComplex}
\delta(\CX) \colonequals
\sum\limits_{p \in \mathbb{Z}} \dim \Hom_{\D(\A)}(\CX, S_\star[p]).
\end{equation}
\end{dfn}
We note that the defect
is a well-defined number from $\N_0$ since the complex is perfect. If the complex $\CX$ is minimal, that is, if $\im \partial_p \subseteq \rad \X_{p-1}$ for any $p \in \Z$, its defect $ \delta(\CX)$ is equal to the number $
\sum_{p\in \Z} [\X_p : P_\star]
$,
where $[X_p:P_\star]$ denotes the multiplicity of the indecomposable projective $P_\star$ as direct summand  of the projective $\X_p$ at degree $p$.

\begin{lem}\label{lem:defect}
For any object $\CX$ of $\DbRep{\A}$ we have:
\begin{equation}\label{E:DefectComplexTwo}
\delta(\CX)  =
\sum\limits_{p \in \mathbb{Z}} \dim  \Hom_{\D(\A)}(S_\star, \CX[p]).
\end{equation}
Moreover, it holds that
$\delta(\CX) = \delta\bigl(\tau(\CX)\bigr) =  \delta\bigl(\inv(\CX)\bigr) =
\delta\bigl(\dual{\CX}\bigr)=
\delta\bigl(\cdual{\CX}\bigr)$.
\end{lem}
\begin{proof}
The first formula \eqref{E:DefectComplexTwo}  follows from the Serre duality \eqref{E:SerreDual} and the fact that  $\tau(S_\star) \cong S_\star[1]$; see Lemma~\ref{L:SimpleSpherical}. Since $\tau$ is an auto-equivalence of $\DbRep{\A}$, the latter fact also implies that $\delta(\CX) = \delta\bigl(\tau(\CX)\bigr)$. Since $\inv(P_\pm) = P_\mp$ and $\inv(P_\star) = P_\star$, the formula
$\delta(\CX) = \delta\bigl(\inv(\CX)\bigr)$ is true for any object $\CX$ of $\Dbmod{\A}$. Finally, since $\dual{S_\star} \cong S_\star$ and $\dual{(-)}$ is an involutive contravariant auto-equivalence of $\DbRep{\A}$, we conclude that 
$$
\dim \Hom_{\D(\A)}(\CX, S_\star[p]) = 
\dim \Hom_{\D(\A)}(S_\star, \dual{\CX}[p]).
$$
Hence, the forelast equality $\delta(\CX) = \delta\bigl(\dual{\CX}\bigr)$ is true by the formula  \eqref{E:DefectComplexTwo}, and the last one by \eqref{eq:duals}.
\end{proof}
\begin{rmk}\label{rmk:def-lie}
For a Harish-Chandra module $M$ from the principal block $\HC_\circ(\lieg, K)$ and $p \in \N_0$, the vector space
$\Ext^p_{(\lieg, K)}(F_1, M)$ is the \emph{$p$-th relative Lie cohomology of $M$}, where $F_1$ denotes the one-dimensional
$U(\lieg)$-module, see \cite{Borel-Wallach, Knapp-Vogan}. Thus, the defect of an object of $\Rep{\A}$  is 
the total dimension of the  relative Lie cohomology of the corresponding Harish-Chandra module.
\end{rmk}
In what follows, we 
denote by $\add(P_{+-})$
the full subcategory
of $\proj \A$ consisting of direct summands of projective modules $P_+$ and $P_-$.
Moreover, we denote by $S_\star^\perp$ the full subcategory of $\DbRep{\A}$ given by
objects $\CX$ 
such that $\Hom_{\D(\A)}(\CX,S_\star[i])=0$ for all integers $i \in \Z$.
Similarly, $\leftindex^{\perp}{{S_\star}}$
denotes the left orthogonal subcategory of $S_\star$ in $\DbRep{\A}$.

\begin{prp}\label{prp:tau-per}
For any object $\CX$ from $\DbRep{\A}$ the following statements are equivalent.
\begin{enumerate}
\item \label{tau-per1}$\delta(\CX)=0$, that is, $\CX$ is isomorphic to a complex from $\Hotbfd\bigl(\add(P_{+-})\bigr)$.
\item \label{tau-per2a}
There is an isomorphism $\tau(\CX) \cong \inv(\CX)$. 
\item \label{tau-per2b} 		The complex $\CX$ is $\tau^2$-invariant, that is, $\tau^2(\CX) \cong \CX$. 
\item \label{tau-per3} The complex $\CX$ is $\tau$-periodic, that is, 
there exists $p \in \N$ with $\tau^p(\CX) \cong \CX$.
\end{enumerate}
\end{prp}
\begin{proof}
\eqref{tau-per1} $\Rightarrow$ \eqref{tau-per2a}:
Let $\delta(\CX)=0$. Then $\CX \in \mathstrut^\perp S_\star$ and thus 
$\tau(\CX) \cong \inv( \Tw(\CX)) \cong \inv(\CX)$ by Proposition~\ref{P:twist}, which shows~\eqref{tau-per2a}.

The implications \eqref{tau-per2a} $\Rightarrow$ \eqref{tau-per2b} $\Rightarrow$ \eqref{tau-per3} are true by definitions.

To show \eqref{tau-per3} $\Rightarrow$ \eqref{tau-per1} suppose that
$\tau^p(\CX) \cong \CX$ for some $p > 0$.
Assume that $\CX$ is not isomorphic to a complex from $\Hotbfd\bigl(\mathsf{add}(P_{+-})\bigr)$. Then there exists
a maximal integer $m \in \mathbb{Z}$ such that $\Hom_{\D(\A)}(\CX, S_\star[m]) \ne 0$. Using Lemma~\ref{L:SimpleSpherical}
it follows that
$$
\Hom_{\D(\A)}\bigl(\CX, S_\star[m]\bigr) \cong
\Hom_{\D(\A)}\bigl(\tau^{p}(\CX), \tau^{p}(S_\star)[m]\bigr) \cong
\Hom_{\D(\A)}\bigl(\CX, S_\star[m+p]\bigr) = 0,
$$
a contradiction. This shows the remaining implication.
\end{proof}
\begin{rmk}
Note that we have an isomorphism of functors 
\begin{equation}
\tau\Big|_{\Hotb\bigl(\add(P_{+-})\bigr)} \simeq \inv\Big|_{\Hotb\bigl(\add(P_{+-})\bigr)}.
\end{equation} Moreover, the above proof and Lemma~\ref{L:SimpleSpherical} imply that
$
\Hotbfd\bigl(\add(P_{+-})\bigr) = S_\star^{\perp} = \leftindex^{\perp}{{S_\star}}$.
\end{rmk}

\begin{cor} An object   $\CX$ of $\DbRep{\A}$ 
is $\tau$-invariant if and only if $\dual{\CX}$ is $\tau$-invariant. 
\end{cor}
\begin{proof} 
Since $\dual{(-)}$ is an involution, it is sufficient to prove the `only if'-implication.
Assume that $\CX$ is $\tau$-invariant. 	Lemma~\ref{lem:defect} and Proposition~\ref{prp:tau-per} imply that $\delta(\CX) = \delta(\dual{\CX}) = 0$, and thus $\tau(\CX) \cong \inv(\CX)$ as well as $\tau(\dual{\CX}) \cong \inv(\dual{\CX})$.
It follows that $\tau(\dual{\CX}) \cong \inv(\dual{\CX}) = \dual{\bigl(\inv(\CX)\bigr)} \cong \dual{\bigl(\tau(\CX)\bigr)} \cong \dual{\CX}$.
\end{proof}

\section{The gluing method}
\label{sec:glue}

Following \cite{Nodal}, we 
recall a categorical framework  in this section, the 
\emph{category of gluing triples} $\Gcat{\A}$.
Its construction is based on 
a
realization of the Gelfand order as pullback ring
in a diagram of certain structurally simpler rings $\AI$, $\BI$ and $\B$, as shown on the left.
\begin{align*}
	\begin{cd}
		\A \ar[hookrightarrow]{r} \ar{d} \ar[twoheadrightarrow]{d} 
		\arrow[phantom]{dr}[very near start,description]{\pullback}
		\& \B \ar[twoheadrightarrow]{d}{\pi} \\
		\AI \ar[hookrightarrow]{r}{\iota} \& \BI 
	\end{cd}
	&&
	\begin{cd}
		\CX \ar[mapsto]{rr}{\B  \otimes_{\A} -}  
		\ar[mapsto]{d}[swap]{\AI \oA -} 
		\& \& \CY 
		\ar[mapsto]{d}
		\\
		\CV 
		\ar[mapsto]{r}
		\& \BI \oAI \CV \ar{r}{\vartheta}[swap]{\sim} \& \BI \oB \CY
	\end{cd}
\end{align*}
The pullback diagram gives rise to a diagram of homotopy categories and functors.
In particular, to any complex $\CX$ of projective $\A$-modules we may associate a complex $\CY$ 
of projective $\B$-modules, a complex $\CV$ of projective $\AI$-modules and an isomorphism
$\vartheta$ of complexes of projective $\BI$-modules.
The triple $\F(\CX) = (\CV,\CY,\vartheta)$ yields an example of an object in the category
$\Gcat{\A}$. In fact, the complex $\CX$ can be reconstructed from this triple. The general notion of a triple is obtained relaxing the condition that $\CY$, $\CV$ and $\vartheta$ originate from a complex of projective $\A$-modules.
A slight variation of this construction
allows to relate $\DbRep{\A}$ to $\Gcatfd{\A}$ in Corollary~\ref{cor:dbrepgl}.

Given an arbitrary gluing triple $\gamma$, the main novelty of this section is an explicit construction of a complex $\CX'$ with $\FF(\CX') \cong \gamma$, which will be useful to prove the gluing rules in the later combinatorics of band and string complexes.

\subsection{Categories of gluing triples}
\label{subsec:cat-glue}
As before, the category $\DbRep{\A}$ is equivalent to the category $\Hotbfdproj{\A}$.
At first, we will focus on a categorical construction $\Gcat{\A}$ for the larger category  
$\Hotbproj{\A}$, and at the end pass to the relevant subcategories.
Our starting point is 
the embedding of
the Gelfand order into a hereditary order
\begin{align*}
	\begin{cd}
		\A
		= 
		\begin{pmatrix}
			\Rx & \idm & \idm\\
			\Rx & \Rx & \idm \\
			\Rx & \idm & \Rx
		\end{pmatrix}
		\ar[hookrightarrow]{r} \&
		\B
		= 
		\begin{pmatrix}
			\Rx & \idm & \idm\\
			\Rx & \Rx & \Rx \\
			\Rx & \Rx & \Rx
		\end{pmatrix}
	\end{cd}
\end{align*}
The annihilator ideal $I =\ann_{\A}(\B/\A) =  \A e_{\star} \A$ is mapped onto the ideal $\B e_{\star} \B$ under this embedding. Factoring out these ideals yields a monomorphism of semisimple $\kk$-algebras
\begin{align}\label{eq:emb2}
	\begin{cd} \AI \colonequals \A/I \cong 
		\begin{pmatrix} \kk & 0  \\ 
			0 & \kk   \end{pmatrix} \ar[hookrightarrow]{r}{\iota}
		\& 
		\begin{pmatrix}  
			\kk & \kk \\  \kk & \kk \end{pmatrix}
		\cong \B/I \equalscolon \BI\,.
	\end{cd}
\end{align}
\begin{rmk}
	The choices above are natural in the following sense.
	\begin{enumerate}
		\item There is an $\Rx$-algebra isomorphism
		$\begin{td} \B \ar{r}{\sim} \&\bigl(\End_{\A}(\rad \A)\bigr)^{\op} \end{td}$
		mapping an element $\gamma \in \B$ to the right multiplication $\varrho_\gamma$ on $\rad \A$.
		This shows that the overring $\B$ has
		an intrinsic definition obtained by a `normalization process' from $\A$.
		\item 
		The ideal $I$ is the \emph{conductor ideal}
		of the pair $(\A,\B)$ and may be viewed as their largest common two-sided ideal. In particular, the semisimple algebras $\AI$ and $\BI$ are as small as possible.
	\end{enumerate}
\end{rmk}
In particular, the Gelfand order can be realized as pullback ring in the diagram on the left
\begin{align}
	\label{eq:pb}
	\begin{cd}
		\A \ar[hookrightarrow]{r} \ar{d} \ar[twoheadrightarrow]{d} 
		\arrow[phantom]{dr}[very near start,description]{\pullback}
		\& \B \ar[twoheadrightarrow]{d}{\pi} \\
		\AI \ar[hookrightarrow]{r}{\iota} \& \BI 
	\end{cd}
	&&
	\begin{cd}
		\Hotbproj{\A} \ar{r}{\B  \otimes_{\A} -} \ar{d}{\AI  \otimes_{\A} -} \& \Hotbproj{\B} \ar{d}{\BI \otimes_{\B} -} \\
		\Hotbproj{\AI} \ar{r}{\BI  \otimes_{\AI} -} \& \Hotbproj{\BI}
	\end{cd}
\end{align} 
where $\pi$ denotes the natural projection. 
The pullback diagram 
gives rise to a diagram of bounded homotopy categories  and tensor product functors on the right.
This diagram commutes in the sense that there is an isomorphism of compositions of functors 
\begin{align}\label{eq:can}
	\begin{cd}
		\kappa \colon \BI \oAI (\AI \oA -)
		\ar{r}{\sim}
		\& \BI \oBI (\B \oA -)
	\end{cd}
\end{align}
The basic idea of the definition of a gluing triple is 
that an object $\CX$ from $\Hotbproj{\A}$ 
should be represented by its associated natural transformation $\kappa_{\CX}$.
The next definition yields the central construction to realize this idea.
\begin{dfn} \label{dfn:triples}
	The category $\Gcat{\A}$  \emph{of gluing triples} associated to the pair of functors 
	\begin{align}
		\label{eq:gt-func}
		\begin{cd}
			\Hotbproj{\AI} \ar{r}{\BI  \otimes_{\AI} -} \& \Hotbproj{\BI}
			\ar[<-]{r}{\BI  \otimes_{\B} -} \& \Hotbproj{\B}
		\end{cd}
	\end{align}
	is defined as follows.
	\begin{enumerate}
		\item An object of $\Gcat{\A}$ is given by a triple
		$(\CV,\CY,\vartheta)$ comprised from  a complex 
		$\CV$ from $\Hotbproj{\AI}$, a complex
		$\CY$ from $\Hotbproj{\B}$ 
		and an isomorphism $$\begin{td} \vartheta \colon 
			\BI \otimes_{\AI} \CV \ar{r}{\sim} \& \BI \otimes_{\B} \CY \end{td}$$ 
		in the category $\Hotbproj{\BI}$.
		\item A morphism 
		from an object $(\CV,\CY,\vartheta)$ to an object $(\CV',\CY',\vartheta')$ in $\Gcat{\A}$ is given by a pair
		$(\phi,\psi)$ which consists of a morphism $\begin{td} \phi \colon \CV \ar{r} \& \CV' \end{td}$ in $\Hotbproj{\AI}$
		and a morphism $\begin{td} \psi\colon \CY \ar{r} \& \CY' \end{td}$ in $\Hotbproj{\B}$ such that the  diagram 
		\begin{align*}
			\begin{cd}
				\BI \otimes_{\AI} \CV
				\ar{d}[swap]{
					\id  \otimes_{\AI} \phi
				}
				\ar{r}{\vartheta}[swap]{\sim} \& 
				\BI \otimes_{\B} \CY 
				\ar{d}{ \id \otimes_{\B} \psi}
				\\
				\BI \otimes_{\AI} \CV' \ar{r}{\vartheta'}[swap]{\sim} \& \BI \otimes_{\B} \CY'
			\end{cd}
		\end{align*}
		commutes 	 in $\Hotbproj{\BI}$, that is, both compositions of morphisms are homotopic.
	\end{enumerate}
\end{dfn}
It is straightforward to define the direct sum of gluing triples. In particular, the category $\Gcat{\A}$ is 
$\kk$-linear 
and 
additive.
\begin{rmk}
	\begin{enumerate}
		\item
		There is a notion of 
		\emph{pullback category} 
		or \emph{comma category}   associated to the pair of functors in \eqref{eq:gt-func},
		which is defined in the same way as the category of gluing triples except that
		the morphism $\vartheta$ of an object $(\CV,\CY,\vartheta)$
		does not have to be an isomorphism.
		\item
		In fact, the category $\Gcat{\A}$ depends on the choice of the overring $\B$
		and the choice of a common two-sided ideal $I$ of $\A$ and $\B$,  but we suppress these dependencies in the notation.
	\end{enumerate}
\end{rmk}

In the next remark, we recall certain useful features of semiperfect rings.
\begin{rmk} \label{rmk:minimal}
	For any semiperfect ring $S$, the following statements are true.
	\begin{enumerate}
		\item Any complex in $\Hotbproj{S}$ is isomorphic to a \emph{minimal complex}, that is, a complex $(P_{\bu},\partial)$ such that
		$\im \partial_i \subseteq \rad P_{i-1}$ for each degree $i \in \Z$.
		\item \label{rmk:minimal2} A morphism between minimal complexes in $\Hotbproj{S}$ is a  homotopy equivalence if and only if it is an  isomorphism of complexes.
	\end{enumerate}
\end{rmk}
The notion of a minimal complex 
has a counterpart in the category of gluing triples introduced in the next statement.
\begin{lem}\label{lem:min-tri}
	Any object in $\Gcat{\A}$ is isomorphic to a \emph{minimal triple}, that is, 
	a gluing triple $(\CV,\CY, \vartheta)$ satisfying any of the following, equivalent conditions.
	\begin{enumerate}
		\item \label{mt1} The complex $\CV$ has zero differentials, the complex $\CY$ is minimal and the morphism $\vartheta$ is an isomorphism of complexes. 
		\item \label{mt2} The
		complexes $\CV$ and $\CY$ are minimal.
		\item \label{mt3} The complex $\CY$ is minimal and $\vartheta$ is an isomorphism of complexes.
	\end{enumerate}	
\end{lem}
\begin{proof}
	Let $(\CV',\CY',\vartheta')$ be a gluing triple.
	Since $\AI$ is semisimple there is a 
	complex $\CV$ with zero differentials and an isomorphism 
	$\begin{td} \phi \colon \CV' \ar{r}{\sim} \&  \CV \end{td}$  in $\Hotbproj{\AI}$.
	As $\B$ is semiperfect, there is a minimal complex $\CY$ and an isomorphism
	$\begin{td} \psi\colon \CY' \ar{r}{\sim} \&  \CY \end{td}$  in $\Hotbproj{\B}$. 
	Setting $\vartheta \colonequals (\BI \otimes\psi) \cdot \vartheta' \cdot (\BI \otimes\phi)^{-1}$ it follows that $\begin{td} (\phi,\psi)\colon (\CV',\CY',\vartheta') \ar{r} \& (\CV,\CY,\vartheta) \end{td}$ defines an isomorphism of gluing triples, which shows the first claim.
	
	To show the implication \eqref{mt2} $\Rightarrow$ \eqref{mt1}, let $(\CV,\CY, \vartheta)$ be a gluing triple with minimal complexes $\CV$ and $\CY$. 
	Since $\AI$ is semisimple and
	$\rad \B \subseteq \iota(I)$, both complexes $\CV$ and $\BI \oB \CY$ have zero differentials, and thus
	$\vartheta$ is an isomorphism of complexes.
	
	To show that \eqref{mt3} implies \eqref{mt1}, let $(\CV,\CY, \vartheta)$ be a triple with a minimal complex $\CY$ and $\vartheta$ an isomorphism of complexes.
	As there is an embedding $\begin{td} \AI \ar[hookrightarrow]{r} \& \BI \end{td}$, the  unit $\begin{td} \V_i  \ar[hookrightarrow]{r} \& \BI \oAI \V_i \end{td}$ is injective for each degree $i \in \Z$. As $\BI \oB \CY \cong \BI \oAI \CV$ has zero differentials, so does $\CV$. The remaining implications are trivial.
\end{proof}
There is a natural functor
$\begin{td} \F\colon \Hotbproj{\A} \ar{r} \& \Gcat{\A} \end{td}$
defined as follows.
\begin{itemize}
	\item For any object $\CX$ from $\Hotbproj{\A}$ we set $\F(\CX) \colonequals (\AI \otimes_{\A} \CX, \B \otimes_{\A} \CX, \kappa_{\CX} )$
	where $\kappa_{\CX}$ is the natural isomorphism induced from \eqref{eq:can}, that is,
	the composition of the following canonical isomorphisms:
	\begin{align*}
		\begin{td}
			\BI \otimes_{\AI} \AI \otimes_{\A} \CX \ar{r}{\sim} \&
			\BI \otimes_{\A} \CX \ar{r}{\sim} \&
			\BI \otimes_{\B} \B \otimes_{\A} \CX
		\end{td} 
	\end{align*}
	\item For any morphism $\begin{td}\xi \colon \CX \ar{r} \& \CX' \end{td}$ we set $\F(\xi) \colonequals (\AI \otimes_{\A} \xi, \B \otimes_{\A} \xi)$.
\end{itemize}
A straightforward verification shows
$\F$ 
defines  a functor. 
\begin{lem}\label{lem:minimal}
	A complex $\CX$
	in $\Hotbproj{\A}$ is minimal if and only if the complex $\B \oA \CX$ is minimal if and only if  
	the gluing triple $\F(\CX)$ is minimal.
\end{lem}
\begin{proof}
	The complex $\CX$ is minimal if and only if the complex $(\head \A) \otimes_{\A} \CX$ has zero differentials if and only if 
	$(\head \B) \otimes_{\head \A} \head \A \otimes_{\A} \CX \cong  (\head \B) \otimes_{\B} \B \otimes_{\A} \CX$ has
	zero differentials
	if and only if $\B \otimes_{\A} \CX$ is minimal, which is equivalent to the minimality of $\F(\CX)$ by Lemma~\ref{lem:min-tri} using that $\kappa_{\CX}$ is an isomorphism of complexes.
\end{proof}
\begin{rmk}\label{rmk:support}
	For a minimal complex $\CX$ we denote 
	$\supp \CX = \{ i \in \Z \mid X_i \neq 0 \}$.
	For any minimal complex $\CX$ from $\Hotbproj{\A}$
	Lemma~\ref{lem:minimal} implies that $\supp \CX = \supp (\B \oA \CX)$.
\end{rmk}
Vice versa, we may assign a complex $\CX \colonequals \Gl(\gamma)$ of $\A$-modules to any object $\gamma = (\CV,\CY,\vartheta)$ from $\Gcat{\A}$ as follows, setting $\CW \colonequals \BI \oB \CY$ below.
\begin{itemize}
	\item 
	Note that there is a natural projection $\begin{td} \eta_{\CY} \colon \CY \ar[twoheadrightarrow]{r} \& \CW\end{td}$.
	\item
	The isomorphism $\begin{td} \vartheta \colon \BI \otimes_{\AI} \CV \ar{r}{\sim} \& \CW \end{td}$ has an adjoint $\begin{td} \vartheta^* \colon \CV \ar[hookrightarrow]{r} \& \CW \end{td}$ in the category $\Hotbproj{\AI}$  given by the composition $\begin{td} \vartheta \cdot \eta_{\CV} \colon \CV \ar[hookrightarrow]{r} \& \BI \oAI \CV \ar{r}{\sim} \& 
		\CW \end{td}$. 
	\item  
	We take the pullback of 
	$\eta_{\CY}$ and $\vartheta^*$  
	in the category of complexes of $\A$-modules
	\begin{align}\label{eq:pb-complex}
		\begin{cd}
			\CX \ar[twoheadrightarrow]{d}[swap]{\beta} \ar[hookrightarrow]{r}{\alpha}
			\arrow[phantom]{dr}[very near start,description]{\pullback}
			\& \CY \ar[twoheadrightarrow]{d}{\eta_{\CY}} \\
			\CV 
			\ar[hookrightarrow]{r}{\vartheta^*}
			\&
			\CW
		\end{cd}  
	\end{align}
\end{itemize}
A priori, the pullback $\CX$ is a complex of finitely generated $\A$-modules.
We will see shortly 
that these $\A$-modules are actually projective, which allows us to view
the complex 
$\CX$ as an
object of the category $\Hotbproj{\A}$.
\begin{rmk}
	Any morphism $(\phi,\psi)$  in the category $\Gcat{\A}$ is given by a pair of homotopy classes $(\ol{g},\ol{h})$ 
	of morphisms of complexes.
	There is a natural definition of a morphism  $\Gl(g,h)$
	of complexes of $\A$-modules using the universal property of the pullback.
	However, if $g$  and $h$ are homotopic to zero, 
	$\Gl(g,h)$ may not be homotopic to zero, in general.
	Therefore, $\Gl$ does not extend to a functor in a natural way.
\end{rmk}
We recall a weakened version of the notion of an equivalence functor \cite[Definition 1.1.4]{Baues}.
\begin{dfn}
	A functor $\begin{td} F\colon \Ccat \ar{r} \& \Dcat \end{td}$
	is \emph{detecting}
	if $F$ is essentially surjective, full and reflects isomorphisms, that is, 
	a morphism $\xi$ in $\Ccat$ is an isomorphism if $F(\xi)$ is.
\end{dfn}
\begin{rmk}
	Any detecting functor $\begin{td} F \colon \Ccat \ar{r} \& \Dcat \end{td}$ induces an equivalence of categories
	$\begin{td} \Ccat/\ker F \ar{r}{\sim} \& \Dcat \end{td}$
	where the kernel of $F$ is given by all morphisms $\xi$ in $\Ccat$ such that $F(\xi) =0$
	and $\Ccat/\ker F$ denotes the quotient category with respect to the ideal $\ker F$.
\end{rmk}
The next statement shows the utility of the category $\Gcat{\A}$.
\begin{thm}\label{thm:detect}
	In the setup above, the functor $\begin{td}\F \colon \Hotbproj{\A} \ar{r} \& \Gcat{\A}\end{td}$ is detecting.
	Moreover, the following statements hold.
	\begin{enumerate}
		\item \label{det1} Any complex $\CX$ from $\Hotbproj{\A}$
		is isomorphic to $\Gl \F (\CX)$.
		\item \label{det2} Any gluing triple $\gamma$ from $\Gcat{\A}$
		is isomorphic to $\F\Gl(\gamma)$.
		\item \label{det3} It holds that $\ker F = \ker \bigl( \begin{td} \B \oA - \colon \Hotbproj{\A} \ar{r}\& \Hotbproj{\B} \end{td}\bigr)$.
	\end{enumerate}
\end{thm}
The above statement follows from \cite{Nodal}. We give a brief alternative proof with minimal complexes as key ingredient.
\begin{proof} 
	For a semiperfect ring $S$  let $\Combminproj{S}$ denote the full subcategory of the category
	$\Combproj{S}$
	given by minimal complexes.
	If $S$ is semisimple, this category will be identified with the category $\bigoplus_{i \in \Z} \proj S$ whose objects are viewed as bounded complexes of finitely generated projective $S$-modules with zero differentials.

	The pullback diagram \eqref{eq:pb} gives rise to diagrams of categories and functors
	\begin{align}\label{eq:cats}
		\begin{cd}
			\proj{\A} \ar{r}{\B \oA -} \ar{d}{\AI \oA -} \& \proj \B \ar{d}{\BI \otimes_{\B} -} \\
			\proj{\AI} \ar{r}{\BI  \otimes_{\AI} -} \& \proj \BI
		\end{cd}
		&&
		\begin{cd}
			\Combminproj{\A} \ar{r}{\B \oA -} \ar{d}{\AI \oA -} \& \Combminproj{\B} \ar{d}{\BI \otimes_{\B} -} \\
			\Combzeroproj{\AI} \ar{r}{\BI  \otimes_{\AI} -} \& \Combzeroproj{\BI}
	\end{cd} \end{align}
	where the top functor $\B \oA -$ in the right diagram is well-defined by Lemma \ref{lem:minimal}.
	
	We focus on the left diagram first.
	Let $\Gcatproj{\A}$
	denote the category of gluing triples associated to the pair of functors
	$$
	\begin{cd}
		\proj{\AI} \ar{r}{\BI  \otimes_{\AI} -} \& \proj \BI \ar[<-]{r}{\BI \otimes_{\B} -}  \&
		\proj \B. 
	\end{cd}
	$$
	A straightforward adaptation of the definitions of $\F$ and $\Gl$
	yields a pair of functors 
	\begin{align*}
		\begin{cd}
			\proj{\A} = \add \A \ar[yshift=3pt]{r}{\F_0} \ar[yshift=-3pt,<-]{r}[swap]{\Gl_0} \& 
			\add \F_0(\A) \subseteq 
			\Gcatproj{\A}
		\end{cd} 
	\end{align*}
	For any object $X$ from $\proj \A$ there is a natural transformation
	$\begin{td} \eta_X \colon X \ar{r} \& \Gl_0\F_0(X) \end{td}$ by the universal property of the pullback.
	Similarly, for any object $\gamma = (\V,\Y,\vartheta)$ from $\Gcatproj{\A}$ there is a natural transformation $\begin{td} (\phi_{\gamma},\psi_{\gamma}) \colon \F_0 \Gl_0(\gamma) \ar{r} \& \gamma \end{td}$
	which is given by 
	the $\AI$-linear counit morphism and the $\B$-linear counit morphism
	$$
	\begin{td} \phi_{\gamma}\colon   \AI \oA \V \ar{r} \& \V,\
		a \oA s \ar[mapsto]{r} \& a \beta(s) 
		\& \psi_{\gamma} \colon \B \oA \Y \ar{r} \& \Y, \
		x \oA y \ar[mapsto]{r} \& x \alpha(y)\end{td}$$
	with $\alpha$ and $\beta$ as in \eqref{eq:pb-complex}.
	As $\eta_{X}$ and $(\phi_\gamma,\psi_{\gamma})$ are isomorphisms for $X = \A$
	respectively $\gamma= \F_0(\A)$, it follows that the functors
	$\F_0$ and $\Gl_0$ are quasi-inverse equivalences.
	By a result due to Milnor \cite[Theorem 2.3]{Milnor} (see also \cite[Lemma 2.6]{Nodal}),
	any object from $\Gcatproj{\A}$
	is isomorphic to an object from
	$\add \F_0(\A)$.
	
	Shifting focus to the right diagram in \eqref{eq:cats},
	we denote by $\GcatCombminproj{\A}$
	the category of gluing triples associated to the pair of functors
	$$
	\begin{cd}
		\Combzeroproj{\AI} \ar{r}{\BI  \otimes_{\AI} -} \& \Combzeroproj{\BI} \ar[<-]{r}{\BI \otimes_{\B} -}  \&
		\Combminproj{\B}
	\end{cd}
	$$
	The functors $\F_0$ and $\Gl_0$ extend directly to a pair of quasi-inverse equivalences of the categories in the top row of the diagram
	\begin{align*}
		\begin{cd}
			\Combminproj{\A} \ar{d}{\Prf_1} \ar[yshift=5pt]{r}{\F'}[swap]{\sim} \ar[yshift=-3pt,<-]{r}[swap]{\Gl'}  \& 
			\GcatCombminproj{\A}
			\ar{d}{\Prf_2} 
			\\  
			\Hotbproj{\A} \ar{r}{\F} \& 
			\Gcat{\A}
		\end{cd}
	\end{align*}
	The functor $\Prf_1$ in the above diagram is defined by viewing a minimal complex $\CX$ from $\Combminproj{\A}$ as an object in $\Hotbproj{\A}$, and assigning a
	morphism $f$ of minimal complexes 
	to its homotopy class $\overline{f}$.
	Similarly, the functor $\Prf_2$ maps an object
	$(\CV,\CY,p)$ from $\GcatCombminproj{\A}$ to $(\CV,\CY,\overline{p})$ 
	and a morphism $(g,h)$ of minimal triples to the pair $(\overline{g}, \overline{h})$ of their homotopy classes.
	
	According to Remark~\ref{rmk:minimal}, the functor $\Prf_1$ is detecting.
	On the other hand, the functor $\Prf_2$ is detecting as well:
	it is essentially surjective because of 
	Lemma~\ref{lem:min-tri}, full by its definition 
	and that of the category $\GcatCombminproj{\A}$,
	and reflects isomorphisms, because
	any pair $(\overline{g}, \overline{h})$ of homotopy equivalences between minimal triples is already a pair of isomorphisms of complexes due to Remark~\ref{rmk:minimal}~\eqref{rmk:minimal2}.
	Since $\F \circ \Prf_1 = \Prf_2 \circ \F'$ 
	and the functors $\F'$, $\Prf_1$ and $\Prf_2$ are detecting, so is the functor $\F$.
	
	The last diagram of categories gives rise to a diagram of bijections on sets of isomorphism classes of objects.
	\begin{align*}
		\begin{cd}
			\text{[}\Combminproj{\A}\text{]}
			\arrow[d,  "\rotatebox{90}{$\sim$}" swap, "\Prf_1"]
			\ar[yshift=5pt]{r}{\F'}[swap]{\sim} \ar[yshift=-3pt,<-]{r}[swap]{\Gl'}  \& 
			\text{[}\GcatCombminproj{\A}\text{]}
			\arrow[d,  "\rotatebox{90}{$\sim$}" swap, "\Prf_2"]
			\\
			\text{[}\Hotbproj{\A}\text{]}
			\ar[yshift=5pt]{r}{\F} \ar[yshift=-3pt,<-]{r}[swap]{\Gl}  \& 
			\text{[}\Gcat{\A}\text{]}
		\end{cd}
	\end{align*}
	The equalities $\Prf_2\circ \F' = \F \circ \Prf_1$ and $\Prf_1\circ \Gl' = \Gl \circ \Prf_2$ yield claims \eqref{det1} and \eqref{det2}.
	To show the remaining claim \eqref{det3} we note that 
	for any morphism $\xi$ in $\Hotbproj \A$ with $\B \oA \xi = 0$ it follows that
	$\BI \oB \B \oA \xi = \BI \oA \AI \oA \xi = 0$, and thus $\AI \oA \xi = 0$, because $\BI \oAI -$ is faithful.
	Therefore $\ker (\B \oA -) = \ker (\B \oA -) \cap \ker (\AI \oA -) = \ker F$.
\end{proof}

\begin{cor}\label{cor:KRS}
	The category $\Gcat{\A}$ has the Krull-Remak-Schmidt property.
\end{cor}
\begin{proof}
	We recall that the category $\DbRep{\A}$ has the Krull-Remak-Schmidt property, see~\cite[Corollary~A.6]{Nodal}.
	Since the functor $\F$ is full and reflects the zero object, the endomorphism ring of any indecomposable object in $\Gcat{\A}$ is local.
\end{proof}
To adapt the last theorem to the category $\DbRep{\A}$, we observe the following.
\begin{lem}\label{lem:tri-fd}
	Let $\gamma = (\CV,\CY,\vartheta)$ be an object in $\Gcat{\A}$ and $\CX \colonequals \Gl(\gamma)$ its pullback complex.
	Then the sum $
	\sum_{i \in \Z} \dim H_i(\CX)$ is finite if and only if so is the sum $
	\sum_{i \in \Z} \dim H_i(\CY)$.
\end{lem}
\begin{proof}
	According to \eqref{eq:pb-complex}
	there is a short exact sequence of complexes of $\A$-modules
	$$
	\begin{cd}
		0 \ar{r} \& \CX \ar{r} \& \CV \oplus \CY \ar{r} \& \CW \ar{r} \& 0
	\end{cd}
	$$
	with $\CW = \BI \oB \CY$.
	In particular, there is a long exact homology sequence
	$$
	\begin{cd}
		\ldots
		H_{i+1}(\CW) \ar{r} \&
		H_i(\CX) \ar{r}\& H_i(\CV) \oplus H_i(\CY) \ar{r} \& H_i(\CW) \ldots
	\end{cd}
	$$
	As $\AI$ and $\BI$ are finite-dimensional algebras, 
	the vector spaces $H_{i}(\CW)$ and $H_i(\CV)$ are finite-dimensional at any degree $i \in \Z$. This implies the claim.
\end{proof}

We denote by $\Gcatfd{\A}$
the full subcategory of $\Gcat{\A}$ given by gluing triples $(\CV,\CY,\vartheta)$ such that $\CY$ has finite-dimensional homology at each degree.
\begin{cor} \label{cor:dbrepgl}
	There are  mutually inverse bijections on the sets of isomorphism classes of indecomposable objects
	\begin{align*}
			\begin{cd}\ind \DbRep{\A}    \ar[yshift=5pt]{r}{\F}[swap]{\sim}
				\ar[<-,yshift=-3pt]{r}[swap]{\Gl} \&  \ind \Gcatfd{\A}
			\end{cd}
	\end{align*}
\end{cor}
\begin{proof}
	The functors $\F$ and $\Gl$ 
	induce quasi-inverse bijections between the set of isomorphism classes of indecomposable objects in $\Hotbproj{\A}$ and that of $\Gcat{\A}$ by Theorem~\ref{thm:detect}.
	The map $\Gl$ restricts further to a bijection of the subsets in the claim
	by 
	Lemma~\ref{lem:tri-fd}.
\end{proof}

\subsection{Gluing of differentials}
\label{subsec:glu}
In order to describe the indecomposable objects of the category $\DbRep{\A}$, it will be crucial to 
construct the complex $\Gl(\gamma)$ explicitly for any gluing triple $\gamma$ from $\Gcat{\A}$ up to isomorphism. 
This is the goal of the present subsection.
Throughout this subsection, 
let $\gamma=(\CV, \CY, \vartheta)$ be a gluing triple in $\Gcat{\A}$, which we may assume to be minimal according to Lemma~\ref{lem:minimal}.
We set $\CX = \Gl(\gamma)$ as in \eqref{eq:pb-complex}.
\subsubsection{Problem and motivation}
To simplify notation, we focus first on the definition of the differential $\partial_1$ of $\CX$.
The $\A$-linear map
$\partial_1$ is obtained using the universal property of the kernel
in the diagram
$$
\begin{cd}
	0 \ar{r} \& \X_1 \ar{r}{\begin{psmallmatrix} \beta_1 \\ \alpha_1 \end{psmallmatrix}} \ar[dashed]{d}{\partial_1} \& \V_1 \oplus \Y_1 \ar{d}{0 \oplus \wt{\partial}_1 } \ar{r} \& \W_1 \ar{d}{0} \ar{r} \& 0
	\\	
	0 \ar{r} \& \X_0 \ar{r}{\begin{psmallmatrix} \beta_0 \\ \alpha_0 \end{psmallmatrix}} \& \V_0 \oplus \Y_0 \ar{r} \& \W_0  \ar{r} \& 0
\end{cd}
$$
in which the maps $\alpha_1$ and $\alpha_0$ are monomorphisms of $\A$-modules.
In more detail, the differential
$\wt{\partial}_1$ of the complex $\CY$ restricts to the map $\wt{\partial}_1|_{\im \alpha_1}$ in the diagram
$$
\begin{cd}
	\X_1 \ar[dashed]{d}{\partial_1} \ar{r}{\check{\alpha}_1}[swap]{\sim} \& \im \alpha_1 \ar[dashed]{d}{\wt{\partial}'_1} \ar[hookrightarrow]{r} \& \Y_1 \ar{d}{\wt{\partial}_1} \\
	\X_0 \ar{r}{ \check{\alpha}_0}[swap]{\sim} \& \im \alpha_0 \ar[hookrightarrow]{r} \& \Y_0
\end{cd}
$$
and the differential $\partial_1$ is precisely the composition $\check{\alpha}_{0}^{-1} \cdot \wt{\partial}_1|_{\im \alpha_1} \cdot \check{\alpha}_1$.
However, to  determine $\partial_1$ explicitly we would need to choose a set of $\A$-linear generators $\{x_1,x_2,\ldots, x_n\}$ 
of $\im \alpha_{1}$, a set of $\A$-linear generators $\{y_1,y_2,\ldots, y_m\}$ of $\im \alpha_0$ 
and express $\wt{\partial}'_1(x_i)$ of each generator $x_i$ as a $\A$-linear combination $\sum_{j=1}^m \lambda_{ij} y_j$ of the generators of $\im \alpha_0$.
Instead of performing these computations, we will replace $(\CV,\CY,\vartheta)$ by an isomorphic triple $(\CV',\CY',\vartheta')$  
whose pullback complex $\CX'$ admits
simpler embeddings $\begin{td} \alpha'_1 \colon \X'_1 \ar[hookrightarrow]{r} \& \Y'_1 \end{td}$ and $\begin{td} \alpha'_0 \colon \X'_0 \ar[hookrightarrow]{r} \& \Y'_0 
\end{td}$ with natural choices of generators. This will come at the expense that the differential $\wt{\partial}'_1$ of the complex $\CY'$ 
will be more complicated than the differential $\wt{\partial}_1$ of $\CY$.

\subsubsection{Passing to a gluing triple in the essential image}

Next, we fix some notation.
Up to isomorphism the indecomposable projective $\B$-modules and the only simple $\BI$-module are given by
\begin{align*}
	\BP_\star = \begin{pmatrix} \Rx \\ \Rx \\ \Rx \end{pmatrix},
	&&
	\BP_\diamond = \begin{pmatrix} \mx \\ \Rx \\ \Rx \end{pmatrix},
	&&
	\wt{S}_\diamond = \begin{pmatrix} 0 \\ \kk \\ \kk \end{pmatrix}.
\end{align*}

\paragraph{Lifting automorphisms of semisimple $\BI$-modules}
Let $\wt{S}$ be a semisimple $\BI$-module and $\BP$ a projective $\B$-module such that $\BI \oB \BP \cong \wt{S}$.
We lift an automorphism $M$ of $\wt{S}$ to a specific automorphism $\wh{M}$ of $\BP$ using the following considerations.
\begin{enumerate}
	\item Any $\BI$-linear morphism $\begin{td} \wt{S}_{\diamond} \ar{r} \& \wt{S}_{\diamond}\end{td}$ is given by right multiplication  with $\lambda$ for certain $\lambda \in \kk$.
	Since $\Rx$ is a $\kk$-algebra, such a $\BI$-linear endomorphism  can be lifted to an endomorphism of the $\B$-module $\BP_{\diamond}$ given by right multiplication with $\lambda$ as well.
	\item We assume that $\wt{S} = \wt{S}_\diamond^n$ for some $n \in \N_0$
	and that $\BP = \BP_\diamond^n \oplus \BP_\star^r$ with $r \in \N_0$.
	Then the automorphism $M$ of $\wt{S}$ is given by a matrix with entries in $\kk$.
	\item The first observation yields the commutative diagram 	of $\kk$-algebras on the left
	\begin{align}
		\label{eq:lift}
		\begin{cd}
			\Aut_{\B}(\BP_{\diamond}^n)  \ar[twoheadrightarrow]{d} \ar{r}{\sim} \& \GL_{n}(\Rx) \ar[xshift=0pt,twoheadrightarrow]{d}[swap,xshift=-5pt]{p} \&[-0.5cm] M \\
			\Aut_{\BI}(\wt{S}_{\diamond}^n) \ar{r}{\sim} \& \GL_{n}(\kk) \ar[xshift=5pt, hookrightarrow,dashed]{u} \& M					\ar[mapsto,dashed]{u}
		\end{cd}
		&&
		\begin{cd}
			\BP_{\diamond}^n \oplus \BP_\star^r  \ar[dashed]{r}{\wh{M}
				=  \begin{psmallmatrix}M & 0 \\ 0 & \Id \end{psmallmatrix} }[swap]{\sim} 
			\ar[->>]{d} \&[0.5cm] \BP_{\diamond}^n  \oplus \BP_{\star}^r \ar[->>]{d} \\
			\wt{S}_{\diamond}^n \ar{r}{M}[swap]{\sim} \& \wt{S}_{\diamond}^n
		\end{cd}
	\end{align}	
	where the surjective ring morphism $p$ admits a left inverse.
	The image of the matrix $M$ under the left inverse is denoted by $M$ again.
	The diagram on the right shows the
	choice of the lift $\wh{M}$
	of
	the automorphism of $\wt{S}^n$ given by the matrix $M$.
\end{enumerate}
Summarized, we may view any automorphism of $\wt{S}^n$ as an automorphism of $\wt{P}_\diamond^n$ and take its direct sum with an identity of $\wt{P}_\star^r$.

\paragraph{Notation for embeddings}

Let  $\begin{td}\iota \colon \AI \ar[hookrightarrow]{r} \& \BI \end{td}$
denote the embedding of semisimple $\kk$-algebras from \eqref{eq:emb2}.
We denote by 
$\begin{td}\iota_{\pm} \colon S_{\pm} \ar[hookrightarrow]{r} \& \wt{S}_{\diamond} \end{td}$
the $\AI$-linear map
which maps $a e_{\pm}$ to $\iota(a e_{\pm})$. Its 
adjoint $\BI$-linear map
is denoted by $\begin{td}\tilde{\iota}_{\pm} \colon \BI \otimes_{\AI} S_{\pm} \ar{r}{\sim} \& \wt{S}_{\diamond}\end{td}$
and maps $b \otimes s$ to $b \iota_{\pm} (s)$.

\paragraph{Construction of $\CX'$}
Given the minimal triple $\gamma = (\CV,\CY,\vartheta)$ in $\Gcat{\A}$ as above, we define a new complex $\CX'$ isomorphic to $\CX = \Gl(\gamma)$ by the following steps for each $i \in \Z$.
\begin{enumerate}
	\item
	We may assume that there are decompositions $\V_i = \bigoplus_{a=1}^{n_i} S{(a)}$
	with	$n_i \in \N_0$ and  each $S{(a)}$ given either by $S_+$ or $S_-$,
	and $\Y_i = \BP_{\diamond}^{n_i} \oplus \BP_{\star}^{r_i}$
	with $r_i \in \N_0$.
	\item Next, we choose certain bases for the $\BI$-modules $\BI \otimes_{\AI} \V_i$ and $\BI \otimes_{\B} \Y_i$.
	\begin{enumerate}
		\item
		There is a commutative diagram in the category of  $\AI$-modules 
		$$
		\begin{cd} 
			S_i = \bigoplus_{a=1}^{n_i} S{(a)} \ar[hookrightarrow]{d}{\eta_{S_i}} \ar[hookrightarrow]{r}{\jmath_i =  \bigoplus \iota_{S{(a)}}} \& \wt{S}_{\diamond}^{n_i} = \bigoplus_{a=1}^{n_i} \wt{S}_{\diamond} \\ 
			\BI \otimes_{\AI} \left(\bigoplus_{a=1}^{n_i} S{(a)}\right) \ar{r}{c}[swap]{\sim}  \ar[dashed]{ru}{\zeta_i} \& \bigoplus_{a=1}^{n_i} \left(\BI \otimes_{\AI} S{(a)}\right) 
			\arrow[u, "\bigoplus \widetilde{\iota}_{S{(a)}}"  swap, "\rotatebox{90}{$\sim$}"]
		\end{cd}
		$$
		We define a new isomorphism $\zeta_i = \bigoplus \widetilde{\iota}_{S{(a)}} \cdot c$.
		\item 
		There exists a $\BI$-linear isomorphism $\begin{td} \nu_i \colon \BI \otimes_{\B} \Y_i \ar{r}{\sim} \& \wt{S}_{\diamond}^{n_i} \end{td}$.
		We define a $\B$-linear morphism $\begin{td} \wt{\pi}_i \colon \Y_i \ar[twoheadrightarrow]{r} \& \wt{S}_{\diamond}^{n_i} \end{td}$ of $\B$-linear modules as the composition $\nu_i \cdot \eta_{Y_i}$
		with the unit morphism $\begin{td} \eta_{Y_i} \colon	
			\Y_i \ar[twoheadrightarrow]{r} \& \BI \otimes_{\B} \Y_i \end{td}$.
	\end{enumerate}
	Since $\BI \cong \Mat_{2 \times 2}(\kk)$ according to~\eqref{eq:emb2}, having chosen the isomorphisms $\zeta_i$ and $\nu_i$ 
	allows us to depict $\vartheta_i$ by an invertible matrix  $\begin{td} M_i = \nu_i \vartheta_i \, \zeta_i^{-1} \, \colon \wt{S}_{\diamond}^{n_i} \ar{r}{\sim} \& \wt{S}_{\diamond}^{n_i}\end{td}$.
	\item We define a new complex $\CY'$ of projective $\B$-modules as follows.
	For each $i \in \Z$ we set  $\Y'_i  = \Y_i$, $\wt{\pi}'_i = \wt{\pi}_i$,
	lift $M_i$ to an invertible matrix $\wh{M}_i$
	as described in~\eqref{eq:lift}
	and set
	$\begin{td} {\wt{\partial}}'_i = {\wh{M}}^{-1}_{i-1} \cdot {\wt{\partial}}_i \cdot {\wh{M}}_{i} \colon \Y'_i \ar{r} \& \Y'_{i-1} \end{td}$,
	where $\wt{\partial}_i$ denotes the differential of $\CY$ at degree $i$.
	In different terms, we define morphisms making the following diagrams commute.
	\begin{align*}
		\begin{cd}
			\BI \otimes_{\AI} \V_i  
			\arrow[d, "\rotatebox{90}{$\sim$}"  swap, "\vartheta_i"]
			\ar{r}{\zeta_i}[swap]{\sim} \& 
			\wt{S}_{\diamond}^{n_i}  \arrow[d, dashed, "\rotatebox{90}{$\sim$}"  swap, "M_i"] \& \Y'_i = \BP_{\diamond}^{n_i} \oplus \BP_{\star}^{r_i} \ar[twoheadrightarrow]{l}[swap]{\wt{\pi}'_i}  
			\arrow[d, dashed, "\rotatebox{90}{$\sim$}"  swap, "{\wh{M}_i = \begin{psmallmatrix} M_i & 0 \\ 0 & \Id \end{psmallmatrix}}"]
			\\ 
			\BI \otimes_{\AI} \Y_i \ar{r}{\nu_i}[swap]{\sim} \&  		  \wt{S}_{\diamond}^{n_i} 
			\&  \ar[twoheadrightarrow]{l}[swap]{\wt{\pi}_i} \Y_i
			= \BP_{\diamond}^{n_i} \oplus \BP_{\star}^{r_i}  
		\end{cd}
		&&
		\begin{cd}
			Y'_i 
			\arrow[d, "\rotatebox{90}{$\sim$}"  swap, "\wh{M}_{i}"] \ar[dashed]{r}{\wt{\partial}'_i}   
			\& Y'_{i-1} 
			\arrow[d, <-, "\rotatebox{90}{$\sim$}"  swap, "\wh{M}^{-1}_{i-1}"]\\
			Y_i \ar{r}{\wt{\partial}_i} 
			\& Y_{i-1}
		\end{cd}
	\end{align*}
	This yields a complex $\CY'$ of finitely generated projective $\B$-modules,
	which is isomorphic to $\CY$. Thus $\CY'$ defines an object in $\Hotbfdproj{\B}$.
	\item 
	Let	${\wt{S}}_{\bu}$ denote the complex of semisimple $\BI$-modules with trivial differentials and $\wt{S}_i = \wt{S}_{\diamond}^{n_i}$ for each integer $i \in \Z$.
	Let $\begin{td} \jmath\colon \CV \ar[hookrightarrow]{r} \& {\wt{S}}_{\bu} \end{td}$ denote the monomorphism of complexes given by 
	$\begin{td} \jmath_i = \bigoplus \iota_{S{(a)}} \colon  S_i \ar[hookrightarrow]{r} \& \wt{S}_i \end{td}$
	for each $i \in \Z$.
	The complex $\CX'$ is defined as the pullback in the category of complexes of $\A$-modules with respect to the morphisms $\jmath$ and $\wt{\pi}'$.
	$$
	\begin{cd}
		\CX' 				\arrow[phantom]{dr}[very near start,description]{\pullback} \ar[dashed,hookrightarrow]{r}{\alpha'} \ar[dashed,twoheadrightarrow]{d}{\beta'}\&[0.5cm] \CY' \ar[twoheadrightarrow]{d}{\wt{\pi}'} 
		\\
		\CV \ar[hookrightarrow]{r}{\jmath} \& {\wt{S}}_{\bu} 
	\end{cd}
	$$
\end{enumerate}
This completes the definition of the complex $\CX'$.
\begin{prp}
	In the notations above, the following statements hold.
	\begin{enumerate}
		\item	The complexes $\CX'$ and $\CX$ are isomorphic. In particular, $\CX'$ defines an object in $\Hotbfdproj{\A}$.
		\item At each degree $i \in \Z$ the differential of the complex $\CX'$ is given by 
		\begin{align}\label{eq:diff}
			\partial'_i = \wt{\partial}'_i|_{X'_i} = 
			(\wh{M}_{i-1}^{-1} \cdot \wt{\partial}_i \cdot \wh{M}_i)|_{X'_i}
		\end{align}
		where the restriction to $\X'_i$ is to be understood along the 
		$\A$-linear monomorphism $\begin{td}\alpha'_i \colon \X'_i \ar[hookrightarrow]{r} \& \Y'_i \end{td}$ of $\A$-modules, that is, along the embedding
		$$
		\begin{cd} \bigoplus_{a=1}^{n_i} P{(a)} \oplus P_{\star}^{r_i} \ar[hookrightarrow]{r}{\bigoplus \iota_{P(a)} \oplus \id } \&[0.5cm] \BP_{\diamond}^{n_i} \oplus \BP_{\star}^{r_i} \end{cd}$$
		where $\begin{td} \iota_{P(a)} \colon P{(a)} \ar[hookrightarrow]{r} \& \BP_{\diamond} \end{td}$ denotes the natural lift of $\begin{td} \iota_{S(a)} \colon S{(a)} \ar[hookrightarrow]{r} \& \wt{S}_{\diamond} \end{td}$.
	\end{enumerate}
\end{prp}
\begin{proof}
	We note that the collection $(M_i)_{i \in \Z}$ defines an endomorphism $M$ of the complex ${\wt{S}}_{\bu}$.
	The complex $\CX'$ is defined by the short exact sequence of the first row
	\begin{align*}
		\begin{cd}
			0 \ar{r} \& \CX' \ar[dashed]{d} \ar{r} \& \CV \oplus \CY' 
			\arrow[d, "\rotatebox{90}{$\sim$}"  swap, "\id \oplus
			\wh{M}"]
			\ar{r}{\begin{psmallmatrix} \jmath & \wt{\pi}'\end{psmallmatrix}} \& {\wt{S}}_{\bu} 
			\arrow[d, "\rotatebox{90}{$\sim$}"  swap, "M"]
			\ar{r} \& 0 \\
			0 \ar{r} \& \CX''  \ar{r} \& \CV \oplus \CY  \ar[equal]{d} \ar[yshift=3pt]{r}{\begin{psmallmatrix} M \jmath & \wt{\pi} \end{psmallmatrix}} \& {\wt{S}}_{\bu} 					\arrow[d, <-, "\rotatebox{90}{$\sim$}"  swap, "\nu"]  \ar{r} \& 0 \\
			0 \ar{r} \& \CX  \ar[dashed]{u} \ar{r} \& \CV \oplus \CY   \ar[yshift=3pt]{r}{\begin{psmallmatrix} \vartheta^* & \eta_{\CY} \end{psmallmatrix}} \& \BI \otimes_{\B} \CY  \ar{r} \& 0 
		\end{cd}
	\end{align*}
	As the squares on the right commute, the universal property of the pullback implies that there are morphisms towards $\CX''$ making the left squares commute. As both middle and both right vertical maps are isomorphisms, so are both left vertical maps, which shows the first claim.
	The second claim follows from the construction of $\CX$.
\end{proof}
Formula~\eqref{eq:diff} will be central to prove the gluing rules in Subsection~\ref{subsec:glurul}.

\begin{rmk}
		The complex $\CX'$ is also isomorphic to the pullback complex $\Gl(\gamma')$ of the gluing triple $\gamma' = (\CV,\CY',\vartheta')$,
		where for each index $i \in \Z$ the morphism $\vartheta'_i$ corresponds to the identity matrix in the previously chosen bases, that is, there is a commutative diagram
		\begin{align*}
			\begin{cd}
				\BI \otimes_{\AI} \V_i  \ar[dashed]{d}{\vartheta'_i}  \ar{r}{\zeta_i}[swap]{\sim} \& 
				\wt{S}_{\diamond}^{n_i}  \ar{d}{\Id} \\ 
				\BI \otimes_{\AI} \Y'_i \ar{r}{\nu_i}[swap]{\sim} \&  		  \wt{S}_{\diamond}^{n_i} 
			\end{cd}
		\end{align*}
		In fact, the gluing triple $\gamma'$ can be identified with $\F(\CX)$, $\vartheta'$ with $\kappa_{\CX}$, 
		and each embedding $\alpha'_i$ with the unit morphism
		$\begin{td} \X_i \ar[hookrightarrow]{r} \& \B \oA \X_i \end{td}$.
		In this context, the matrices $(\wh{M}_i)_{i \in \Z}$ encode the transformations used to obtain a decomposition of $\B \otimes_{\A} \CX$.
\end{rmk}

\section{Reduction to matrix problem}
\label{sec:reduce}
In this section, we first recall the combinatorics of the category $\DbRep{\B}$ of the hereditary envelope 
$\B$ of the Gelfand order $\A$.
These combinatorics motivate the definition of the category of regular representations $\regrep{\fXA}$.
We describe the relationship of the former to the category $\Gcatfd{\A}$, and conclude the categorical reductions for the category $\DbRep{\A}$. The overall picture is summarized in Theorem~\ref{thm:picture}, its essential consequences in Corollary~\ref{cor:reduction}.
Using the category $\regrep{\fXA}$, we reduce the classification problem for the category $\Gcatfd{\A}$ 
to an explicit problem of linear algebra,  which is spelled out in terms of partitioned matrices in Subsection~\ref{subsec:mp} and visualized in Figure~\ref{fig:der-mp}.

\subsection{The derived representation theory of the envelope of the Gelfand order}
The overring $\B$ of the Gelfand order has two indecomposable projective modules which will be denoted by $\BP_\star$ and $\BP_\diamond$, and two simple modules
$\wt{S}_\star$ and	$\wt{S}_{\diamond}$ up to isomorphism.
\begin{align*}
	\B
	= 
	\begin{pmatrix}
		R & \idm & \idm\\
		R &  R & R \\
		R & R & R
	\end{pmatrix}
	&&
	\BP_{\star} = \begin{pmatrix}
		R \\ R \\ R \end{pmatrix}
	&&
	\BP_{\diamond} = \begin{pmatrix}
		\mx \\ R \\ R
	\end{pmatrix}
	&&
	\wt{S}_\star = 
	\begin{pmatrix} \kk \\ 0 \\ 0 \end{pmatrix}
	&&
	\wt{S}_\diamond =
	\begin{pmatrix} 0 \\ \kk \\ \kk \end{pmatrix}
\end{align*}
For vertices $u,v \in \{\diamond,\star\}$ we denote $\mathcal{N}_{uv} = \N_0$ if $(u,v) =(\diamond,\star)$ and $\mathcal{N}_{uv} = \N$ otherwise.
\begin{lem}\label{lem:ind-env} 
	The ring $\B$ is hereditary and uniserial. In particular, the
	following statements hold.
	\begin{enumerate}
		\item A $\B$-module $N$ is finite-dimensional and indecomposable
		if and only if $N$ is isomorphic to 
		the quotient $\BP_v/\rad^\ell(\BP_v)$
		for a vertex $v \in \{\star,\diamond\}$ 
		and a number $\ell \in \N$.
		\item An object $\CY$ from $\DbRep{\B}$ is indecomposable
		if and only if $\CY$ is isomorphic to a shift of a projective resolution of an indecomposable finite-dimensional module, that is,  $\CY$ is isomorphic to a complex 
		\begin{align}
			\label{eq:ind-H}
			(\begin{gd}\BP_{u} \ar{r}{\cdot t^n} \& \BP_v \end{gd})[d] \text{ with }
			u,v \in \{\star,\diamond\}, n \in \mathcal{N}_{uv} \text{ and }d \in \Z.
		\end{align}
	\end{enumerate}
\end{lem}
\begin{proof}
	It is straightforward to check that the 
	homology of a complex of the form \eqref{eq:ind-H} is a quotient of the form in the first claim.
	As each quotient is a cyclic module, each of these complexes is indecomposable.
	
	Minimal projective resolutions of simple $\B$-modules are given by the complexes in \eqref{eq:ind-H} with $d=0$, $u \neq v$, and minimal $n \in \mathcal{N}_{uv}$. This shows that $\B$ is hereditary.
	By a well-known fact (see e.g. \cite[Proposition~4.4.15]{Krause1}), 
	any indecomposable 
	object in $\Dbmod{\B}$ is isomorphic to the shift of an indecomposable projective resolution.
	Note that the differential in any projective resolution $(\begin{td} \wt{\partial}_1 \colon \Y_1 \ar{r} \& \Y_0 \end{td})$ 
	is given by a matrix in $\Mat_{m \times n}(\Rx)$ with $m,n \in \N_0$. The application of the algorithm for the Smith normal form
	to such a matrix yields a decomposition of $\CY$ into the indecomposable complexes of the form \eqref{eq:ind-H} with $d=0$, stalk complexes $\BP_{v}$  
	and acyclic complexes $\begin{td} \BP_v \ar{r}{\id} \& \BP_v \end{td}$ with $v \in \{\star,\diamond\}$.
	As $\Hotbfdproj{\B}$ is equivalent to $\DbRep{\B}$, the `only if'-implications in both claims follow.
\end{proof}
	For each of the resolutions in \eqref{eq:ind-H}, 
	we will denote its only non-trivial homology 
	by $\BM_{uv}^{(n)}$. As $n$ runs through $\mathcal{N}_{uv}$ and $u,v$ through $\{\star,\diamond\}$, these modules yield a complete list of indecomposable finite-dimensional $\B$-modules up to isomorphism.

	\begin{rmk}
		Given a $\B$-module $\BM_{uv}^{(n)}$ as above,
		its socle is isomorphic to $\wt{S}_{u}$, its top to $\wt{S}_v$ and
		its length is equal to $2n$ if $u=v$, to $2n-1$ if $(u,v) = 
		(\star,\diamond)$ and to $2n+1$ if $(u,v) = (\diamond,\star)$.
	\end{rmk}

	\begin{rmk}
		The ring $\B$ is Morita equivalent to the arrow ideal completion of the path algebra of the cyclic quiver with two vertices~\eqref{E:CyclicQuiver}.
		In particular, the results of this subsection yield a complete description of the derived category
		$\Db{\HC_{\gamma(0)}(\lieg, K)}$.
	\end{rmk}

	The Auslander-Reiten quiver of $\DbRep{\B}$ is given by countably many homogeneous tubes of rank two.
	
	Figures~\ref{fig:ray} and \ref{fig:coray} show two tubes, in which we highlight a certain ray and a certain coray which will be relevant below.

	\begin{lem}\label{lem:row-trans}
		Let $\CY$ and $\CY'$ be indecomposable complexes of the form \eqref{eq:ind-H}.
		Then the image of the map
		$$
		\begin{cd}
			\Hom_{\DbRep{\B}}(\CY,\CY') \ar{r}{\BI \oB - } \& 
			\Hom_{\DbRep{\BI}}(\BI \oB \CY,\BI \oB \CY')
		\end{cd}
		$$
		has a $\kk$-linear basis given by morphisms $\BI \oB \psi'$ 
		such that 
		$\psi'_i = \id$ for certain $i \in \Z$ and $\psi'$ is one of the following morphisms.
		\begin{enumerate}
			\item $\psi'$ is an identity morphism, that is, 
			$\psi'$ is given by 
			\begin{align*}
				\begin{cd} \CY \ar{d}{\psi'} \\ \CY' 
				\end{cd} \qquad
				\begin{cd}
					\BP_{\diamond} \ar[equal]{d}{\psi'_i} \ar{r}{t^{p}} \& \BP_u \ar[equal]{d} \\
					\BP_{\diamond} \ar{r}{t^{p}} \& \BP_u
				\end{cd}
				\qquad\text{or}\qquad
				\begin{cd}
					\BP_v \ar[equal]{d} \ar{r}{t^{q}} \& \BP_\diamond \ar[equal]{d}{\psi'_i} \\			\BP_v \ar{r}{t^q} \& \BP_\diamond
				\end{cd}
			\end{align*}		
			for certain $p \in \mathcal{N}_u$ and $q \in \N$.
			\item	$\psi'$ is a finite composition of successive morphisms in the ray of the Auslander-Reiten quiver shown in Figure~\ref{fig:ray}. More explicitly, $\psi'$ is given by
			\begin{align*}
				\begin{cd} \CY \ar{d}{\psi'} \\ \CY' 
				\end{cd} \qquad
				\begin{cd}
					\&		\BP_{\diamond} \ar[equal]{d} \ar{r}{t^p} \& \BP_u \ar{d}{t^{q-p}} \\
					\&		\BP_{\diamond} \ar{r}{t^{q}} \& \BP_v 
				\end{cd}
			\end{align*}
			with $p \in \mathcal{N}_u$, 
			$q \in \mathcal{N}_v$
			such that
			$q  \geq p$ if $(u,v) = (\diamond,\star)$
			and $q  > p$ otherwise.
			\item 
			$\psi'$ corresponds to the shift of an element in $\Ext^1(N,N')$ for certain indecomposable finite-dimensional modules $N, N'$.
			More explicitly, $\psi'$ is given by
			\begin{align*}
				\begin{cd} \CY 
					\ar{d}{\psi'} \\ \CY' 
				\end{cd} \qquad
				\begin{cd}	
					\&		\BP_{\diamond}  \ar[equal]{d} \ar{r}{t^{m}} \& \BP_{u}  \\
					\BP_{v}  \ar{r}{t^n} \& \BP_{\diamond}  
				\end{cd}
			\end{align*}
			with $u,v \in
			\{\diamond,\star\}$, $m \in \mathcal{N}_u$ and
			$n \in \N$.
			\item $\psi'$ is
			a finite composition of successive morphisms 
			in the coray of the Auslander-Reiten quiver shown 
			in Figure~\ref{fig:coray}.
			More explicitly, $\psi'$ is given by
			\begin{align*}
				\begin{cd} \CY \ar{d}{\psi'} \\ \CY' 
				\end{cd} \qquad
				\begin{cd}
					\BP_{u} \ar{d}{t^{q-p}}  \ar{r}{t^{q}} \& \BP_{\diamond} \ar[equal]{d}
					\\
					\BP_{v} \ar{r}{t^{p}} \& \BP_{\diamond}
				\end{cd}
			\end{align*}			
			with $p,q\in \N$ such that $q \geq p$ if $(u,v) = (\diamond,\star)$
			and $q > p$ otherwise.
		\end{enumerate}

						\end{lem}
						\begin{proof}
							The `if'-implication is straightforward. 
							To show the `only if'-implication let $\psi$ be a morphism with $\BI \oB \psi_i \neq 0$.
							Then $\im \psi_i = \Y_i = \Y'_i = \BP_{\diamond}$ and $\psi_i = \lambda r$ for 
							a scalar $\lambda \in \kk^*$ and
							a unit $r \in R^*$.
							As $\psi$ is a morphism of complexes, 
							there are the following cases.
							\begin{align*}
								\begin{cd} \CY \ar{d}{\psi} \\ \CY' 
								\end{cd} \qquad
								\begin{cd}
									\BP_{\diamond}\arrow[d, "\lambda r" swap, "\rotatebox{90}{$\sim$}"] \ar{r}{t^{p}} \& \BP_u \ar{d}[swap]{\lambda s} \\
									\BP_{\diamond} \ar{r}{t^{q}} \& \BP_v
								\end{cd}
								\qquad
								&
								\begin{cd}	
									\&		\BP_{\diamond} \arrow[d, "\lambda r" swap, "\rotatebox{90}{$\sim$}"] \ar{r}{t^{p}} \& \BP_{u}  \\
									\BP_{v}  \ar{r}{t^n} \& \BP_{\diamond}  
								\end{cd}
								&
								\begin{cd}
									\BP_u  \ar{d}[swap]{\lambda s}  \ar{r}{t^{m}} \& \BP_\diamond \arrow[d, "\lambda r" swap, "\rotatebox{90}{$\sim$}"] \\			\BP_v \ar{r}{t^n} \& \BP_\diamond
								\end{cd}\\
							\end{align*}
							where $p \in \mathcal{N}_u$, $q \in \mathcal{N}_v$, $m, n \in \N$, $\BP_{\diamond}$ is located at degree $i$, $u,v\in \{\star,\diamond\}$ and $s \in R^*$ makes the respective diagram commute. If $\psi$ is an automorphism, the first and third diagram lead to the first case.
							Otherwise one of the remaining three cases occurs.
						\end{proof}
						
						\afterpage{
							\begin{figure}
								\caption{A ray in the tube of the Auslander-Reiten quiver 
								}
								\label{fig:ray}
								\begin{subfigure}[t]{0.5\textwidth}
									$
									\begin{td}
										\\[-0.25cm]
										\&[1cm]	\BP_{\diamond} \ar[equal]{d} \ar{r}{1} \& \BP_{\star} \ar{d}{t} \& \BM_{\diamond\star}^{(1)}[i-1] \quad \ar[xshift=-15pt,hookrightarrow]{d} \\
										\&	\BP_{\diamond} \ar[-,densely dotted]{dd} \ar{r}{t} \& \BP_{\diamond} \ar[densely dotted,-]{dd} \& \BM_{\diamond\diamond}^{(1)}[i-1] \quad  \ar[xshift=-15pt,densely dotted,-]{dd} \\
										\\
										\&	\BP_{\diamond} \ar[equal]{d} \ar{r}{t^{n-1}} \& \BP_{\star} \ar{d}{t} \& \BM_{\diamond\star}^{(n)}[i-1] \ar[xshift=-15pt,hookrightarrow]{d} \quad \\
										\&	\BP_{\diamond} \ar[equal]{d} \ar{r}{t^n} \& \BP_{\diamond} \ar{d}{1} 
										\& \BM_{\diamond\diamond}^{(n)}[i-1] \quad  \ar[xshift=-15pt,hookrightarrow]{d}
										\\
										\&	\BP_{\diamond}  \ar[densely dotted,-]{d} \ar{r}{t^n} \& \BP_{\star}  \ar[densely dotted,-]{d}
										\& \BM_{\diamond\star}^{(n+1)}[i-1] \ar[xshift=-15pt,densely dotted,-]{d} \\
										\&	\mathstrut \& \mathstrut \& \mathstrut
									\end{td}$
								\end{subfigure}
								\qquad
								\begin{subfigure}[t]{0.33\textwidth}
									$					
									\underbrace{
										\begin{arq}
											\mathstrut \ar[densely dotted, -]{ddd}	\&   \& \& \mathstrut \ar[densely dotted, -]{dd} \\ \\
											\& {\color{black} \BM^{(3)}_{\diamond \diamond}} \arrow[rd, two heads]            \&                                                            \& \BM^{(3)}_{\star\star} \arrow[llld, dashed,two heads]                                        \\
											{\color{black} \BM^{(3)}_{\diamond \star}} \arrow[rd, two heads] \arrow[ru,hook, color=black]            \&                                                               \&  \BM^{(3)}_{\star \diamond} \arrow[rd, two heads] \arrow[ru,hook] \&                                                                                   \\
											\& \BM^{(2)}_{\star\star}\arrow[rd, two heads] \arrow[ru,hook]                   \&                                                            \& {\color{black} \BM^{(2)}_{\diamond \diamond}} \arrow[llld, two heads, dashed] \arrow[color=black, lllu, dashed,hook] \\
											\BM^{(2)}_{\star \diamond}  \arrow[rd, two heads] \arrow[ru,hook] \&                                                               \& {\color{black}  \BM^{(2)}_{\diamond \star}} \arrow[rd, two heads] \arrow[ru,hook, color=black]            \&                                                                                   \\
											\& {\color{black} \BM^{(1)}_{\diamond \diamond}} \arrow[rd, two heads] \arrow[ru,hook, color=black] \&                                                            \& \BM^{(1)}_{\star\star} \arrow[llld, dashed, two heads] \arrow[lllu, dashed,hook]                   \\
											{\color{black} \BM^{(1)}_{\diamond \star}} \arrow[color=black, ru,hook]                       \&                                                               \& \BM^{(1)}_{\star \diamond} \arrow[ru, hook]                       \&                                                                                  
									\end{arq}}_{\text{shifted by $[i-1]$}}$
								\end{subfigure}
							\end{figure}
							\begin{figure}
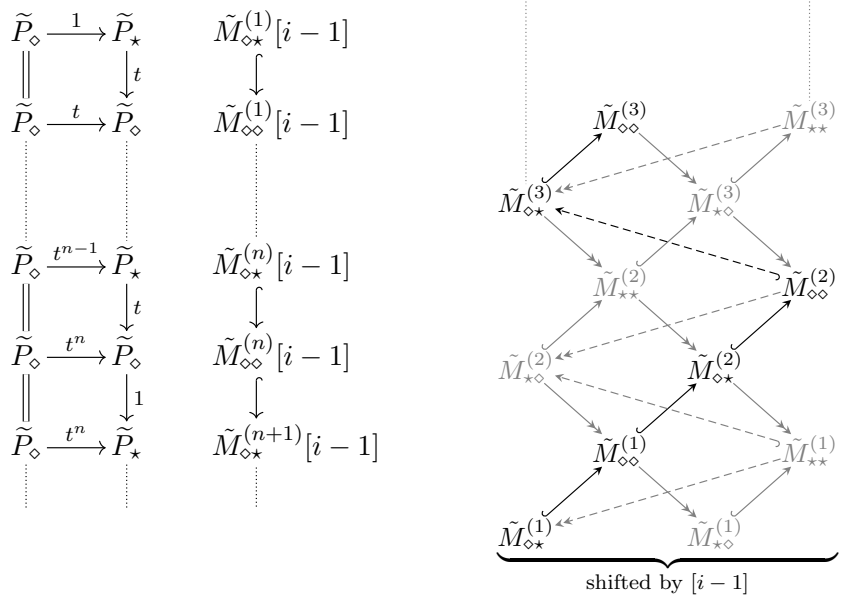

								\caption{A coray in the tube of the Auslander-Reiten quiver
								}
								\label{fig:coray}
								\begin{subfigure}[t]{0.5\textwidth}	
									$\begin{td}
										\mathstrut \ar[-,densely dotted]{d} \& \mathstrut \ar[-,densely dotted]{d} \&[0.75cm]\& \mathstrut \ar[xshift=-10pt,-,densely dotted]{d}  \\
										\BP_{\star}  \ar{d}{t} \ar{r}{t^{n+1}} \& \BP_{\diamond}  \ar[equal]{d} \& \& \BM_{\star \diamond}^{(n+1)}[i] \ar[xshift=-10pt,two heads]{d}\\
										\BP_{\diamond}   \ar{d}{1} \ar{r}{t^{n}} \& \BP_{\diamond}  \ar[equal]{d}  
										\& \& \BM_{\diamond \diamond}^{(n)}[i] \ar[xshift=-10pt,two heads]{d} \quad
										\\
										\BP_{\star}\ar[-,densely dotted]{dd} \ar{r}{t^{n}} \& \BP_{\diamond} \ar[densely dotted,-]{dd} 
										\& \& \BM_{\star \diamond}^{(n)}[i]  \quad \ar[-,densely dotted,xshift=-10pt]{dd}
										\\
										\\
										\BP_{\diamond} \ar{d}{1} \ar{r}{t} \& \BP_{\diamond} \ar[equal]{d} 
										\& \& \BM_{\diamond \diamond}^{(1)}[i] \ar[xshift=-10pt,two heads]{d} \quad \\
										\BP_{\star} \ar{r}{t} \& \BP_{\diamond}
										\& \& \BM_{\star \diamond}^{(1)}[i] \quad
									\end{td}$				
								\end{subfigure}	
								\begin{subfigure}[t]{0.33\textwidth}
									\	
									$\underbrace{
										\begin{arq}
											\mathstrut \ar[densely dotted, -]{ddd}	\&  \& \& \mathstrut \ar[densely dotted, -]{dd} \\ \\
											\& {\color{black} \BM^{(3)}_{\diamond \diamond}} \arrow[rd, two heads,color=black]            \&                                                            \& \BM^{(3)}_{\star\star} \arrow[llld, dashed,two heads]                                        \\
											\BM^{(3)}_{\diamond \star} \arrow[rd, two heads] \arrow[ru,hook]            \&                                                               \& {\color{black} \BM^{(3)}_{\star \diamond}} \arrow[rd, two heads,color=black] \arrow[ru,hook] \&                                                                                   \\
											\& \BM^{(2)}_{\star\star}\arrow[rd, two heads] \arrow[ru,hook]                   \&                                                            \& {\color{black} \BM^{(2)}_{\diamond \diamond}} \arrow[llld, two heads, dashed,color=black] \arrow[lllu, dashed,hook] \\
											{\color{black} \BM^{(2)}_{\star \diamond}} \arrow[rd, two heads,color=black] \arrow[ru,hook] \&                                                               \& \BM^{(2)}_{\diamond \star} \arrow[rd, two heads] \arrow[ru,hook]            \&                                                                                   \\
											\& {\color{black} \BM^{(1)}_{\diamond \diamond}} \arrow[rd, two heads,color=black] \arrow[ru,hook] \&                                                            \& \BM^{(1)}_{\star\star} \arrow[llld, dashed, two heads] \arrow[lllu, dashed,hook]                   \\
											\BM^{(1)}_{\diamond \star} \arrow[ru,hook]                       \&                                                               \& {\color{black} \BM^{(1)}_{\star \diamond}} \arrow[ru, hook]                       \&                                                                                  
									\end{arq}}_{\text{shifted by $[i]$}}$
								\end{subfigure}
							\end{figure}
							\clearpage
						}
						
						\begin{lem}\label{lem:kernel}
							The kernel of the functor 
							$
							\begin{td} \BI \oB - \colon \Hotbfd(\proj \B) \ar{r} \& \Hotbfd(\proj \BI) \end{td}
							$
							is given by the full subcategory
							$
							\Hotbfd(\add \BP_{\star}) = \add \bigl( \begin{td} \BP_{\star} \ar{r}{t^n} \& \BP_{\star}[i] \end{td} \mid n \in \N, i \in \Z \bigr)
							$.
						\end{lem}
						\begin{proof}
							We claim that a non-invertible morphism $\begin{td}\psi\colon \Y \ar{r}\& \Y' \end{td}$ of minimal indecomposable complexes satisfies 
							$\im \psi_i \subseteq I \Y'_i$ for each $i \in \Z$ only if it factors through an object in $\Hotbfd(\add \BP_{\star})$.	
							\begin{itemize}
								\item 
								If $\supp \Y = \supp \Y'$, the morphism $\psi$ is a finite composition of irreducible morphisms.
								The Auslander-Reiten quiver in Figure~\ref{fig:ray}
								shows that any irreducible morphism satisfies the claim, which implies the claim for $\psi$.
								\item Otherwise $\supp \Y \cap \supp \Y' = \{ i\}$ for certain  $i \in \Z$, 
								and $\psi$ has only one non-zero component at degree $i$ 
								and factors through $\BP_{\star}[i]$.
								As $\rad^n (H_i(\Y')) = 0$ for certain $n \in \N$ any morphism $\begin{td}\BP_{\star} [i]\ar{r} \& \Y' \end{td}$, and thus $\psi$, factors through 
								$(\begin{td} \BP_{\star} \ar{r}[yshift=-1pt]{t^n} \& \BP_{\star} \end{td})[i]$.
							\end{itemize}
							As the kernel of a functor is an ideal, it follows that the kernel of $\BI \oB -$ is contained in the full subcategory $\Hotbfd(\add \BP_{\star})$. The other inclusion is trivial.
						\end{proof}

						\subsection{Bunch of semichains for the  derived category}

						The category $\Gcatfd{\A}$ can be defined using the functors
						\begin{align*}
							\begin{cd}
								\Hotbproj{\AI} \ar{r}{\BI  \otimes_{\AI} -} \& \Hotbproj{\BI}
								\ar[<-]{r}{\BI  \otimes_{\B} -} \& \Hotbfdproj{\B}
							\end{cd}
						\end{align*}
						We introduce a notion closely related to the 
						behaviour of indecomposable objects and their morphisms
						in the categories $\Hotbfdproj{\AI}$ and $\Hotbfdproj{\B}$ under the respective tensor product functors.		
						\begin{dfn}\label{def:chkX}
							Let $(\fXA,\geq,\approx)$, abbreviated to $\fXA$ later on, be the following data.
							\begin{itemize}
								\item Let $\fXA = \fE \cup \fF$,
								$\fE = \bigcup_{i \in \Z} \fE_{i}$, 
								$\fF = \bigcup_{i \in \Z} \fF_i$,
								where for each $i \in \Z$ the set
								$\fF_i$ is the set $\{ f_i^+, f_i^- \}$ with trivial partial order, and
								the set $\fE_i$ is the chain 
								\begin{align*}
									\fE_i = \{
									\xx_i^n, \yy^n \mid n \in \N
									\}
									\quad
									\text{with}
									\quad
									\xx^p_i < \xx^{p+1}_i < \yy_i^{q+1} < \yy_i^q \text{ for any }p,q \in \N.
								\end{align*}
								Two elements $e$, $e'$ in $\fXA$ are comparable if 
								and only if $e$, $e' \in \fE_i$ for some $i\in\Z$. 
								\item We consider the smallest equivalence relation $\approx$ on $\fXA$ satisfying
								$\xx_{i+1}^{2n} \approx \yy^{2n}_{i}$
								for any $n \in \N$ and any $i \in \Z$.
							\end{itemize}
						\end{dfn}
						The final purpose of the datum  $\fXA$	
						is to encode a matrix problem which is visualized in Figure~\ref{fig:der-mp} and explained in Subsection~\ref{subsec:mp}.
						Next, we motivate the definition of the elements in $\fXA$ and their equivalence classes by the following features.
						\begin{enumerate}
							\item
							Elements of $\fF$ have trivial equivalence classes and correspond to isomorphism classes of stalk complexes of simple $\AI$-modules, that is, there is a bijection
							\begin{align*}
								\begin{cd}
									\fF = \fF/\!\approx \ar{r}{\sim} \& \ind \Hotbproj{\AI}
									\& f_i^{\pm} \ar[mapsto]{r} \& S_{\pm}[i]
								\end{cd}
							\end{align*}
							\item 
							Equivalence classes of elements of $\fE$ correspond to isomorphism classes of indecomposable objects in $\Hotbfdproj{\B}$ which are not annihilated by the functor $\BI \oB -$.
							In other terms, there is a bijection
							\begin{align*}
								\begin{cd} 	 \fE/\!\approx \ar{r}{\sim} \& \ind \Hotbfdproj{\B}\backslash \ind \ker (\BI \oB - ),
									\& \ol{e} \ar[mapsto]{r} \& \CY(\ol{e}) \end{cd}
							\end{align*}
							where $\CY(\ol{e})$ is  given by the second row of the following table.
							\begin{align}
								\label{eq:Yx}
								\begin{array}{|c||c|c|c|}
									\hline
									\ol{e}
									&
									\{ a_{i+1}^{2n}, 
									b_i^{2n}\} &
									\{b_i^{2n-1}\} & 	\{a_{i+1}^{2n+1}\} \\ \hline \hline
									\CY(\ol{e}) &		(\begin{gd}
										\BP_{\diamond} \ar{r}{\cdot t^{n}} \& \BP_\diamond \end{gd})[i]
									&
									(\begin{gd}
										\BP_{\star} \ar{r}{\cdot t^{n}} \& \BP_\diamond \end{gd})[i] 
									&
									(\begin{gd}
										\BP_{\diamond} \ar{r}{\cdot t^{n}} \& \BP_\star \end{gd})[i] \\ \hline
									\BI \oB \CY(\ol{e})  &
									(\begin{gd}
										\BS_{\diamond} \ar{r}{0} \& \BS_{\diamond}
									\end{gd})[i]
									&
									(\begin{gd}
										0 \ar{r} \& \BS_{\diamond}
									\end{gd})[i]
									&
									(\begin{gd}
										\BS_{\diamond} \ar{r} \& 0
									\end{gd})[i] \\ \hline
								\end{array}
							\end{align}
							The computation in the last row indicates two further properties.
							\begin{enumerate}
								\item The cardinality of each equivalence class $\ol{e}$ is equal to the number of simple summands
								of 	$\BI \oB \Y(\ol{e})$.
								\item Elements of the set $\fE$ correspond to the simple direct summands of $\BI \oB \CY$
								where $\CY$ runs through all indecomposable objects in $\Hotbfdproj{\B}$ up to isomorphism.
							\end{enumerate}
						\end{enumerate}
						The definition of the order relations in $\fXA$ is motivated
						by the next statement, which follows from  
						Lemma~\ref{lem:row-trans}.
						\begin{rmk}\label{rmk:order}
							For any $i \in \Z$, 
							each relation $x < a$ in $\fE_i$ corresponds to a morphism $\begin{td}\psi \colon \CY(x) \ar{r} \& \CY(a)\end{td}$ such that $\BI \oB \psi_{i}$ is the identity of $\wt{S}_{\diamond}[i]$.
							In fact, the relation $x<a$ is irreducible if and only if the morphism $\psi$ is irreducible in the sense of Auslander-Reiten theory.
						\end{rmk}

						\begin{dfn}\label{dfn:repX-A}
							The category $\regrep{\fXA}$ of \emph{regular representations} of $\fXA$ is defined as follows.
							\begin{enumerate}
								\item
								An object of the category $\regrep{\fXA}$
								is given by a family of $\kk$-linear isomorphisms
								$(\begin{td} \mu_i \colon V_i \ar{r}{\sim} \& W_i \end{td})$
								such that 
								$V_i= V_{f_i^+} \oplus V_{f_i^-} $
								and $W_i = \bigoplus_{e \in \fE_i} W_e$ with 
								$\dim W_a = \dim W_e$ for any $a,e \in \fE$ with $a \approx e$, 
								and there are only finitely many integers $i\in \Z$ with $W_i \neq 0$.
								\item
								A morphism between two objects $(\mu_i)_{i \in \Z}$ and $(\mu'_i)_{i \in \Z}$
								is given by families of morphism $(\tau_i)_{i \in \Z}$ and $(\sigma_i)_{i \in \Z}$
								such that 
								the diagram
								$$
								\begin{cd}
									V_i = V_{f_i^+} \oplus V_{f_i^-} \ar{d}[swap]{\tau_{i}=
										\begin{psmallmatrix}
											\tau_i^+ & 0
											\\
											0 &  \tau_i^-
										\end{psmallmatrix}
									} \ar{r}{\mu_i} \& W_i = \bigoplus_{x \in \fE_{i}} W_x 	\ar{d}{\sigma_{i} = (\sigma_{ax})} \\
									V'_i = V'_{f_i^+} \oplus V'_{f_i^-} \ar{r}{\mu'_\iota}
									\& 
									W'_{i}  = \bigoplus_{a \in \fE_{i}}W'_a
								\end{cd}
								$$
								commutes
								and the following conditions are satisfied.
								\begin{itemize}
									\item For any elements $a, x \in \fE$ with $a \approx x$ it holds that $\sigma_{aa} = \sigma_{xx}$.
									\item For any $i \in \Z$ and any elements $a, x \in \fE_i$ with $x \not\leq a$
									it holds that $\sigma_{a x} = 0$. 
								\end{itemize}
							\end{enumerate}
						\end{dfn}

						\subsection{Relationship of gluing triples and matrix representations}
						The next goal is to define a functor $\begin{td} \Red\colon \Gcatfd{\A} \ar{r} \& \regrep{\fX} \end{td}$ which shall be used to reduce the classification problem of the first category to a matrix problem.
						Since $\BI$ is the matrix algebra $\Mat_{2\times 2}(\kk)$,
						it will be convenient to use the Morita equivalence functor
						\begin{align}\label{eq:Morita}
							\begin{cd}
								(-)^b \colon
								\md \BI
								\ar{r}{\sim} \& \md \kk
							\end{cd}
						\end{align}

						\begin{rmk}\label{rmk:triple}
							Any object in $\Gcatfd{\A}$ is isomorphic to a minimal triple $(\CV,\CY,\vartheta)$ of the following form.
							\begin{enumerate}
								\item The complex $\CV$ is a finite direct sum of stalk complexes of simple $\AI$-modules, in particular, for any degree $i \in \Z$ it holds that $\V_i = S_+^{m_{+,i}} \oplus S_-^{m_{-,i}}$ with $m_{\pm,i} \in \N_0$.
								\item  The complex $\CY$ is a finite direct sum of indecomposable complexes described in \eqref{eq:ind-H}, in particular, $$\CY = \bigoplus_{
									\begin{smallmatrix}
										u,v \in \{ \star, \diamond\}\\ n\in \mathcal{N}_{uv}, i \in \Z \end{smallmatrix}}
								\bigl((\begin{gd} \BP_u \ar{r}{t^{n}} \& \BP_v
								\end{gd}) [i]\bigr)^{m_{u,v,n,i}}
								$$
								for certain multiplicities $m_{u,v,n,i} \in \N_0$, where, as before, $\mathcal{N}_{uv} = \N_0$ if $(u,v)=(\diamond,\star)$ respectively $\mathcal{N}_{uv} = \N$ otherwise. 
								\item The morphism $\vartheta$ of complexes of $\BI$-modules
								is given by a family of $\BI$-linear isomorphisms
								$(\begin{td}\vartheta_i \colon \BI \oAI \V_i \ar{r}{\sim} \& \BI \oB \Y_i \end{td})_{i \in \Z}$.
								Because of the Morita equivalence \eqref{eq:Morita}, each $\BI$-linear morphism $\vartheta_i$ is determined by an invertible matrix
								$M_i \in \GL_{m(i)}(\kk)$
								such that the sum $\sum_{i \in \Z}m(i)$ is finite  and for each $i \in \Z$ it holds that
								\begin{align} \label{eq:const}
									m(i) =  m_{+,i} + m_{-,i} = \sum_{n\in \N
									} m_{\star,\diamond,n,i} + m_{\diamond,\diamond,n,i} + m_{\diamond,\star,n,i-1}
									+ m_{\diamond,\diamond,n,i-1}.
								\end{align}
							\end{enumerate}
						\end{rmk}
						Next, we define prescriptions on gluing triples in such a form and their morphisms.
						\begin{dfn}\label{dfn:H}
							Let $\gamma$ be a minimal triple $(\CV,\CY,\vartheta)$ such that $\CV$ and $\CY$
							are given in the form as in the last remark.
							\begin{itemize}
								\item 
								We define $\Red(\gamma) = (\begin{td} \mu_i\colon V_i \ar{r} \& W_i \end{td})_{i \in \Z}$
								with $\mu_i$ given by 
								$$\begin{cd}  \vartheta_i^b \colon\underset{\varepsilon \in \{+,-\}}{\bigoplus} (\BI \oAI \V_{\varepsilon}^{m_{\varepsilon,i}})^b \ar{r}{\sim} \&
									\underset{e \in \fE_i}{\bigoplus} (\BI \oAI \Y_i(\ol{e}))^b
								\end{cd}.$$
								\item 
								Given another minimal triple $\gamma' = (\CV', \CY', \vartheta')$ with $\CV'$ and $\CY'$ of the form as in Remark~\ref{rmk:triple} and a morphism
								$\begin{td}(\phi,\psi) \colon \gamma \ar{r} \& \gamma' \end{td}$ of gluing triples, 
								we set
								$\Red(\phi,\psi) = 
								\bigl((\tau_i)_{i\in\Z}, (\sigma_i)_{i\in\Z}\bigr)$
								with $\tau_i   = 
								(\BI \oAI \phi_i)^b$ and $\sigma_i = (\BI \oB \psi_i)^b$.
							\end{itemize}

						\end{dfn}
						In different terms, to obtain $\Red(\gamma)$ we apply the Morita equivalence functor to each $\BI$-linear morphism $\vartheta_i$ and regard it as a partitioned matrix with two vertical blocks corresponding to summands $S_{+}^{m_{+,i}}$ respectively $S_-^{m_{-,i}}$ of $\V_i$,
						and horizontal blocks indexed by indecomposable objects of $\DbRep{\B}$ which have projectives of type $\BP_{\diamond}$ at degree $i$.
						We note also that each of the matrices $\tau_i$ and $\sigma_i$ can be viewed as a block matrix compatible with these partitions.
						\begin{rmk}\label{rmk:mat-gt}
							Vice versa, given an object $\mu$ in $\regrep{\fX}$ 
							we define an object
							$\Ind(\mu)$ in $\Gcatfd{\A}$ 
							as follows.
								\begin{enumerate}
									\item Set $\CV = \bigoplus_{i \in \Z} (S_{+}^{m_{+,i}} \oplus S_-^{m_{-,i}} )[i]$
									with  $m_{\pm,i} = \dim V_{f_i^{\pm}}$.
									\item Set $\CY = \bigoplus_{\ol{e} \in \fE/\approx} \CY(\ol{e})^{m_{e}}$ with $m_{e} = \dim W_e$
									and $\CY(\ol{e})$ defined by~\eqref{eq:Yx}.
									\item By applying the quasi-inverse 
									of the Morita equivalence in \eqref{eq:Morita}
									for each $i \in \Z$ the $\kk$-linear map $\mu_i$ gives rise to a $\BI$-linear map $
									\begin{td} \vartheta_i \colon \BI \oAI \V_i \ar{r}{\sim} \& \BI \oB \Y_i \end{td}$.
								\end{enumerate}
							It is straightforward to check that $\Ind(\mu)$ is indeed an object of $\Gcatfd{\A}$ and
							that 
							$\Red\Ind(\mu) = \mu$.
						\end{rmk}

				\begin{prp}\label{prp:H}
					The prescriptions in Definition~\ref{dfn:H} yield
					a 
					full and essentially surjective
					functor $\begin{td}\Red \colon \Gcatfd{\A} \ar{r} \& \regrep{\fX} \end{td}$
					such that  $\ker \Red$ is given by the full subcategory of objects in $\Gcatfd{\A}$ 
					which are isomorphic to objects of the form
					$(0,\CY,0)$.  
				\end{prp}
				\begin{proof}
					For any object $\gamma$ of $\Gcatfd{\A}$ we fix an isomorphism to a minimal triple $(\CV,\CY,\vartheta)$
					with $\CV$ and $\CY$ decomposed as described in Remark~\ref{rmk:triple}.
					This allows to extend the prescriptions of Definition~\ref{dfn:H} to all objects and morphisms of the category $\Gcatfd{\A}$.
					It is straightforward to check that $\Red(\gamma)$ is indeed an object of $\regrep{\fX}$.
					Lemma~\ref{lem:row-trans} implies  that 
					$\Red(\phi,\psi)$ is a morphism in $\regrep{\fX}$. It follows that $\Red$ defines a functor.
					
					The functor $\Red$ is essentially surjective by Remark~\ref{rmk:mat-gt}.
					To show that $\Red$ is full, 
					let $(\tau,\sigma)$  be a morphism 
					$\begin{td} \Red(\gamma) \ar{r} \& \Red(\gamma')\end{td}$ in $\regrep{\fX}$.
					We may assume that $\gamma$ and $\gamma'$ are of the form as in Remark~\ref{rmk:triple}. 
					\begin{itemize}
						\item 				Since the functor $\begin{td}\BI \oAI - \colon \add S_{\pm} \ar{r} \& \add \wt{S}_{\diamond} \end{td}$ is fully faithful,
						there is a morphism $\phi$ such that $\BI \oAI \phi = \tau$.
						Since $\BI \oAI -$ reflects isomorphisms, 
						the morphism $\tau$ is invertible if and only if so is $\phi$. 
						\item 
						For any distinct elements $a$, $x \in \fE$ with $a \approx x$ the pair
						of identical matrices $(\sigma_{aa},\sigma_{xx})$
						lifts to an endomorphism $\psi_{aa}$ of $\CY(\ol{a})^{m_a} = \CY(\ol{x})^{m_x}$
						by the first case in Lemma~\ref{lem:row-trans}
						We note that $\sigma$ is an isomorphism if and only if $\sigma_{aa}$ is invertible for each $a \in \fE$. In this case, the lifts $\psi_{aa}$ can be chosen to be isomorphisms as well.
						Remark~\ref{rmk:order} and
						the remaining cases in Lemma~\ref{lem:row-trans} show that
						for each $i\in\Z$ and any $a,x \in \fE_i$ with $x < a$,
						the block matrix $\begin{td}\sigma_{ax} \colon W_x \ar{r} \& W'_a \end{td}$ of the block matrix $\sigma_i$ 
						can be lifted to a morphism of complexes $\begin{td} \psi_{ax} \colon {\CY(x)}^{m_x} \ar{r} \& {\CY(a)}^{m_a} \end{td}$.
						As the remaining blocks $\sigma_{ax}$ are zero, this shows that there is a morphism $\psi$ such that $\BI \oB \psi = \sigma$.
						
					\end{itemize}
					The pair $(\phi,\psi)$ is a morphism of gluing triples since $\sigma_i \mu_i = \mu'_i \tau_i$ for each $i \in \Z$. Therefore, the functor $\Red$ is full.
					
					Let $(\phi,\psi)$ be a morphism in $\Gcatfd{\A}$ with $\Red(\phi,\psi)=0$.					
					Then $\BI \oB \phi= 0$ implies that $\phi = 0$, since $\BI \oAI -$ is faithful.
					Therefore, the kernel of $\Red$ is given by morphisms of gluing triples $(0,\psi)$ with $\BI \oB \psi = 0$.
					Lemma~\ref{lem:kernel} completes the description of $\ker \Red$.
				\end{proof}

				\subsection{Final categorical picture}

				We return our focus towards the category $\DbRep{\A}$ and consider the composition of the functors $\Red$ and $\F$.

				\begin{thm}\label{thm:picture}
					There is a commutative diagram of categories and functors
					\begin{align*}
						\begin{cd}
							\&[-0.25cm] 		 \& \ker(\F)\ar[hookrightarrow]{d}  \ar[hookrightarrow]{ld} \& \\
							\&\ker(\Red \circ \F)  \ar[hookrightarrow]{r}	 \arrow[dd, "\rotatebox{90}{$\sim$}"{pos=0.5,xshift=-8pt}  swap, "\F_{\star}" ] \arrow[ld, "\F'_{\star}" swap, "\rotatebox{30}{$\sim$}"{pos=0.6,yshift=3pt}] \&	\DbRep{\A} \ar{dd}{\F} \ar{r}{\Red \circ \F} \& \regrep{\fXA} \ar[equal]{dd} \\[-0.5cm]
							\mathcal{K} \arrow[rd, "\F''_{\star}" swap, "\rotatebox{-20}{$\sim$}"{pos=0.3,yshift=-5pt}] \&\\[-0.5cm]
							\& \ker(\Red)	\ar[hookrightarrow]{r} \&	\Gcatfd{\A} \ar{r}{\Red} \& \regrep{\fXA}
						\end{cd}
					\end{align*}
					where the functor $\F$ is detecting, its restriction $\F_{\star}$ factors through equivalences, the functor $\Red$ is full and essentially surjective, and the following holds.
					\begin{align*}
						\ker(\F) &=  \ker\bigl(\begin{cd} \B \oA - \colon \Hotbfdproj{\A} \ar{r}\& \Hotbfdproj{\B} \end{cd}\bigr). \\
						\ker(\Red \circ \F) &=  \ker\bigl(\begin{cd} \AI \oA - \colon \Hotbfdproj{\A} \ar{r}\& \Hotbfdproj{\AI} \end{cd}) = \Hotbfd(\add P_{\star}\bigr).\\
						\mathcal{K} &= \ker\bigl(\begin{cd} \BI \oB - \colon \Hotbfdproj{\B} \ar{r}\& \Hotbfdproj{\BI} \end{cd}) = \Hotbfd(\add \BP_{\star}\bigr). \\
						\ker(\Red) &
						\cong 
						\{ (0,\CY,0) \mid \CY \in \mathcal{K} \}.
					\end{align*}

				\end{thm}
				\begin{proof}
					Properties of the functors $\F$ and $\Red$ were shown in Theorem~\ref{thm:detect} and Proposition~\ref{prp:H}.
					The functor $\F'_{\star}$ is given by the restriction of $\B \oA -$ to $\ker(\Red \circ \F)$.
					Since the embedding $\begin{td} \A \ar[hookrightarrow]{r} \& \B \end{td}$
					restricts to an isomorphism $e_\star \A e_{\star} \cong  e_{\star} \B e_{\star}$, the map 
					$$
					\begin{cd}\End_{\A}(P_{\star}) \ar{r}
						\&
						\End_{\A}(\BP_{\star})	\& \psi \ar[mapsto]{r} \& \B \oA \psi
					\end{cd}
					$$
					is an isomorphism, and thus $\F'_\star$
					is an equivalence of categories. 
					
					The functor $\F''_\star$ maps any morphism $\begin{td}\psi \colon \CY \ar{r} \& \CY' \end{td}$ from $\mathcal{K}$ to a morphism of gluing triples $\begin{td}(0,\psi) \colon (0,\CY,0) \ar{r} \& (0,\CY',0)\end{td}$. The definitions of $\ker(\Red)$ and $\F$ imply that
					$\F''_{\star}$ is an equivalence of categories and $\F_{\star} = \F''_{\star} \circ \F'_{\star}$.
				\end{proof}
				
				To derive some consequences of the last theorem, we observe a simple categorical fact.
				\begin{lem}\label{lem:Cker}
					Let $\mathcal{C}$ be a Krull-Remak-Schmidt category, $\mathcal{D}$ an additive category and
					$\begin{td} F \colon \mathcal{C} \ar{r} \& \mathcal{D} \end{td}$ be an additive, full and essentially surjective functor
					such that its kernel is given by a full subcategory $\mathcal{K}$ of $\mathcal{C}$.
					Then $F$ restricts to a detecting functor
					$$ \begin{cd} \mathcal{C}_{\mathcal{K}} \ar{r} \& \mathcal{D}, \& X \ar[mapsto]{r} \& F(X) \end{cd} $$
					where $\mathcal{C}_{\mathcal{K}}$ denotes the full subcategory of objects in $\mathcal{C}$ which do not have non-zero direct summands in ${\mathcal{K}}$.
				\end{lem}
				The statement is well-known in case $F$ is the quotient functor $\begin{td} \md A \ar{r} \& \stmd{A} \end{td}$ for a finite-dimensional algebra $A$ \cite[Lemma 4.3]{SY}. The proof of this special case extends to the situation above. Nevertheless, we give a proof for the convenience of the reader.
				\begin{proof}
					The induced functor is full, and can be shown to be essentially surjective using the Krull-Remak-Schmidt property of $\mathcal{C}$.
					Assume that 
					$\begin{td}\alpha \colon X \ar{r} \& X' \end{td}$ is a morphism in $\mathcal{C}_{\mathcal{K}}$ such that
					$F(\alpha)$ is an isomorphism. Since $F$ is full, there is a morphism $\begin{td}\beta \colon Y \ar{r} \& X \end{td}$ such that $\id_X - \beta \alpha$ factors through an object $K$ from $\mathcal{K}$, that is $\id_X - \beta \alpha = vu$ for certain morphisms $\begin{td}u \colon X \ar{r}\& K \end{td}$ and $\begin{td}v \colon K \ar{r}\& X\end{td}$. Since the indecomposable summands of $X$ and those of $K$ are pairwise non-isomorphic, it follows that 
					$v \in \rad(K,X)$,
					$vu \in \rad(X,X)$ and thus
					$\id_X - vu = \beta \alpha$ is invertible. A similar argument shows that $\alpha \beta$ is invertible. Therefore $\alpha$ is an isomorphism.
				\end{proof}
				
				\begin{cor}\label{cor:reduction}
					The functors $\F$ and $\Red$ restrict to 
					detecting functors
					\begin{align*}
						\begin{cd}
							\DbRep{\A}_{\ker(\Red \circ \F)} \ar{r} 	\&	\Gcatfd{\A}_{\ker(\Red)} \ar{r}  \& \regrep{\fXA} 
						\end{cd}
					\end{align*}
					In particular, there is a pair of quasi-inverse bijections of sets of isomorphism classes of indecomposable objects
					\begin{align*}
						\ind  \DbRep{\A} \backslash 
						\ind \ker(\Red \circ \F) 
						\begin{cd} \mathstrut \ar[yshift=3pt]{r}{\Red \circ \F}[swap]{\sim} 
							\ar[yshift=-3pt,<-]{r}[swap]{\Gl \circ  \Ind}
							\& 
							\ind \regrep{\fXA} 
						\end{cd}
					\end{align*}
					where $\ind \ker(\Red \circ \F) = 
					\ind \Hotbfd(\add P_{\star}) = 
					\{ (\begin{gd} P_{\star} \ar{r}{t^n} \& P_\star \end{gd})[i] \mid n\in \N, i \in \Z \}.$
				\end{cor}
				\begin{proof}
					Since $\End_\A(P_{\star}) \cong \Rx = \kk \llbracket t\rrbracket$, 
					there is an equivalence 
					$\begin{td} \DbRep{\Rx} \ar{r}{\sim}  \&
						\Hotbfd(\add {P}_{\star})\end{td}$ which maps
					$(\begin{gd}\Rx \ar{r}{t^n} \& \Rx\end{gd})[i]$
					to
					$(\begin{gd}P_{\star} \ar{r}{t^n} \& P_{\star}\end{gd})[i]$ for any $n \in \N$ and $i\in \Z$.
					This observation yields a description of $\ind \Hotbfd(\add P_{\star})$.
					The remaining statements follow from
					Theorem~\ref{thm:picture} and Lemma~\ref{lem:Cker}.
				\end{proof}
				
				In more concrete terms, the corollary states that
				a representation
				$\mu$ is indecomposable if and only if so is its associated complex $\Gl\Ind(\mu)$
				and that 
				two representations $\mu$ and $\mu'$ are isomorphic in $\regrep{\fXA}$ if and only if 
				the complexes $\Gl\Ind(\mu)$ and 
				$\Gl\Ind(\mu')$  are isomorphic in $\DbRep{\A}$.
				In particular, the problem to classify the indecomposable objects in the category $\DbRep{\A}$
				is reduced to the problem to classify the 
				indecomposable objects in $\regrep{\fXA}$. 
				The latter is formulated in purely linear algebraic language in the next subsection.

				\subsection{Explicit description of matrix problem}\label{subsec:mp}
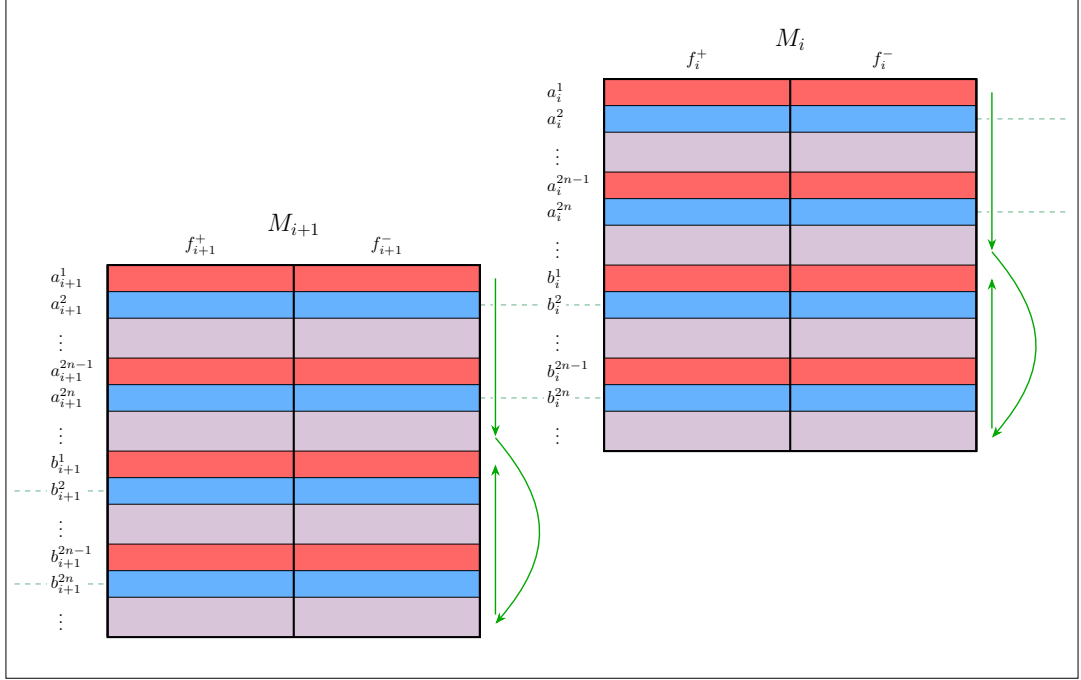
\begin{figure}
					\caption{Matrix problem of the bunch of skew-gentle type $\fXA$}\label{fig:der-mp}
					\centering
					\scalebox{0.66}{$
						\begin{tikzpicture}[every node/.style={font=\small}, >=Stealth, yscale=6/14*1.25, xscale=1.25, show background rectangle]
							\def\shift{-7}

							\draw[dashed,  color=c5] (-2,-1.5+ \shift) -- (0,-8.5);  %
							\draw[dashed, color=c5] (-2,-5+ \shift) -- (0,-12);     %
							\draw[dashed,  color=c5] (-2+8,-1.5) -- (7.5,-1.5); 
							\draw[dashed, color=c5] (-2+8,-5) -- (7.5,-5);   
							\draw[dashed,  color=c5] (-1.5-8,-15.5) -- (-8,-15.5);  
							\draw[dashed, color=c5] (-1.5-8,-19) -- (-8,-19);     
						\def\blockdata{
							1/c1/0/1, 2/c2/1/2, 3/c3/2/3.5,
							4/c1/3.5/4.5, 5/c2/4.5/5.5, 6/c3/5.5/7,
							7/c1/7/8, 8/c2/8/9, 9/c3/9/10.5,
							10/c1/10.5/11.5, 11/c2/11.5/12.5, 12/c3/12.5/14}

						\foreach \i/\c/\y/\z in \blockdata {
							\fill[\c] (-8,-\y+ \shift) rectangle (-2,-\z+ \shift);
							\draw[thin] (-8,-\y+ \shift) -- (-2,-\y+ \shift);
							\draw[very thick] (-8,-\y+ \shift) -- (-8,-\z+ \shift);
							\draw[very thick] (-2,-\y+ \shift) -- (-2,-\z+ \shift);}

						\foreach \i/\c/\y/\z in \blockdata {
							\fill[\c] (0,-\y) rectangle (6,-\z);
							\draw[thin] (0,-\y) -- (6,-\y);
							\draw[very thick] (0,-\y) -- (0,-\z);
							\draw[very thick] (6,-\y) -- (6,-\z);
						}

						\node[anchor=west, fill=white, inner sep=2pt, outer sep=1pt] at (-9,-0.5 + \shift) {$\xx^1_{i+1}$};
						\node[anchor=west, fill=white, inner sep=2pt, outer sep=1pt] at (-9,-1.5 + \shift) {$\xx^2_{i+1}$};
						\node[anchor=west, fill=white, inner sep=2pt, outer sep=1pt] at (-9,-4 + \shift) {$\xx^{2n-1}_{i+1}$};
						\node[anchor=west, fill=white, inner sep=2pt, outer sep=1pt] at (-9,-5 + \shift) {$\xx^{2n}_{i+1}$};
						
						\node[anchor=west, fill=white, inner sep=2pt, outer sep=1pt] at (-9,-12+ \shift) {$\yy^{2n}_{i+1}$};
						\node[anchor=west, fill=white, inner sep=2pt, outer sep=1pt] at (-9,-11+ \shift) {$\yy^{2n-1}_{i+1}$};
						\node[anchor=west, fill=white, inner sep=2pt, outer sep=1pt] at (-9,-8.5+ \shift) {$\yy^2_{i+1}$};
						\node[anchor=west, fill=white, inner sep=2pt, outer sep=1pt] at (-9,-7.5+ \shift) {$\yy^1_{i+1}$};
						
						\node at (-0.75-8, -2.75+\shift) { \vdots}; %
						\node at (-0.75-8, -6.25+\shift) { \vdots}; %
						\node at (-0.75-8, -9.75+\shift) { \vdots}; %
						\node at (-0.75-8, -13.25+\shift) { \vdots}; %
						\node[anchor=west, fill=white, inner sep=2pt, outer sep=1pt] at (-1,-0.5) {$\xx^1_i$};
						\node[anchor=west, fill=white, inner sep=2pt, outer sep=1pt] at (-1,-1.5) {$\xx^2_i$};
						\node[anchor=west, fill=white, inner sep=2pt, outer sep=1pt] at (-1,-4) {$\xx^{2n-1}_i$};
						\node[anchor=west, fill=white, inner sep=2pt, outer sep=1pt] at (-1,-5) {$\xx^{2n}_i$};
						
						\node[anchor=west, fill=white, inner sep=2pt, outer sep=1pt] at (-1,-12) {$\yy^{2n}_i$};
						\node[anchor=west, fill=white, inner sep=2pt, outer sep=1pt] at (-1,-11) {$\yy^{2n-1}_i$};
						\node[anchor=west, fill=white, inner sep=2pt, outer sep=1pt] at (-1,-8.5) {$\yy^2_i$};
						\node[anchor=west, fill=white, inner sep=2pt, outer sep=1pt] at (-1,-7.5) {$\yy^1_i$};

						\node at (-0.75, -2.75) { \vdots}; %
						\node at (-0.75, -6.25) { \vdots}; %
						\node at (-0.75, -9.75) { \vdots}; %
						\node at (-0.75, -13.25) { \vdots}; %
						\draw[->, thick, color=c4, preaction={draw=white, line width=3pt}] (6.25-8,-0.5+\shift) -- (6.25-8,-6.45+\shift); %
						\draw[->, thick, bend left=20, color=c4] (6.25-8,-6.5+\shift) to (6.25-8,-13.5+\shift); %
						\draw[->, thick, color=c4] (6.25-8,-13.15+\shift) -- (6.25-8,-7.5+\shift); %

						\draw[->, thick, color=c4, preaction={draw=white, line width=3pt}] (6.25,-0.5) -- (6.25,-6.45); %
						\draw[->, thick, bend left=20, color=c4] (6.25,-6.5) to (6.25,-13.5); %
						\draw[->, thick, color=c4] (6.25,-13.15) -- (6.25,-7.5); %

						\draw[very thick] (3,0) -- (3,-14); %
						\draw[very thick] (-5,0+ \shift) -- (-5,-14+ \shift); %

						\draw[very thick] (0,0) rectangle (6,-14);
						\draw[very thick] (-8,0+ \shift) rectangle (-2,-14+ \shift);

						\node at (1.5,0.75) {$\xf_i^+$};
						\node at (4.5,0.75) {$\xf_i^-$};
						\node at (3,1.5) {\Large $M_i$};
						
						\node at (-6.5,0.75+ \shift) {$\xf_{i+1}^+$};
						\node at (-3.5,0.75+ \shift) {$\xf_{i+1}^-$};
						\node at (-5,1.5+ \shift) {\Large $M_{i+1}$};
						\node at (0,2.5) {};
						\node at (0,-22) {};
					\end{tikzpicture}
					$} 
				\vspace{0.5em}
				\begin{flushleft}
					{\small 
						Each matrix $M_i$ is invertible with entries in $\kk$.
						Vertical blocks are indexed by the set $\fE_i$, horizontal blocks by $\fF_i$.
						Dotted lines indicate tied labels with equal numbers of rows. Arrows indicate order relations.
					}
				\end{flushleft}
			\end{figure}

				The classification problem for the category $\regrep{\fXA}$ can be stated as follows.
				We are given a collection of matrices $(M_i)_{i \in \Z}$ such that $M_i \in \GL_{m(i)}(\kk)$ such that $m(i) =0 $ for all but finitely many integers $i$,
				and each matrix $M_i$ is partitioned into vertical and horizontal stripes as indicated in Figure~\ref{fig:der-mp}.
			More precisely, the following holds.
			\begin{enumerate}
				\item The columns of each matrix $M_i$ are divided into two vertical stripes labeled by the symbols $f_i^+$ and $f_i^-$. 
				\item The rows of each matrix $M_i$ are divided into horizontal stripes labeled by finitely  many elements from the chain $\fE_i$
				such that for any number $n \in \N$
				the number of rows in the horizontal block labeled by $\xx^{2n}_{i+1}$
				is equal to the number of rows in the horizontal block labeled by $\yy^{2n}_{i}$.
				We will say that the blocks $\xx^{2n}_{i+1}$ and $\yy^{2n}_{i}$ are \emph{tied}.
			\end{enumerate}
			
			Two collections of matrices $(M_i)_{i \in \Z}$
			and $(M'_i)_{i \in \Z}$ are isomorphic if and only if one can be obtained from the other by a sequence of transformations of the following types:
			\begin{enumerate}
				\item \emph{Independent row or column transformations within blocks}. 
				Let $i \in \Z$.
				\begin{enumerate}
					\item For each sign $\varepsilon \in \{ + , - \}$ we may perform any elementary transformations of columns within the vertical block in $M_i$ labeled by $f_i^{\varepsilon}$.
					\item For any
					$e \in \{ \xx_{i}^{2n-1}, \yy_i^{2n-1} \mid n \in \N \}$ 
					we may perform any elementary transformations of rows within the horizontal blocks in $M_i$ labeled by $e$. 
				\end{enumerate}
				\item \emph{Simultaneous row transformations of tied blocks}. For any $i\in \Z$ and any element $e \in \{ \xx_i^{2n}, \yy_i^{2n} \mid n \in \N \}$
				we may perform any elementary transformation of rows in the horizontal block with label $e$ in the matrix $M_i$ together with the \emph{same} transformations of rows in the horizontal block which is tied to $e$, that is, 
				the horizontal block $\xx_{i+1}^{2n}$ in matrix $M_{i+1}$ if $e = \yy_i^{2n}$, respectively
				the horizontal block $\yy_{i-1}^{2n}$ in the matrix $M_{i-1}$ if $e = \xx_i^{2n}$.
				\item \emph{Row transformations between different blocks}. 
				For each integer $i \in \Z$ 
				and for any elements $e, e' \in \fE_i$ with $e < e'$ 
				we may add the multiple by $\lambda \in \kk$ of any row of the horizontal stripe with the smaller label $e$
				to any row of the horizontal stripe with the larger label $e'$.
			\end{enumerate}
			The solution to this problem is provided by work of Bondarenko which is summarized and adapted to a broader class of matrix problems in the next section. 

\section{Bondarenko's canonical forms for bunches of skew-gentle type}
\label{app:mp}

In \cite{B1, B2}, Bondarenko described the canonical normal forms for certain tame matrix problems.
Here, we specialize Bondarenko's classification results to a more narrow class of matrix problems over combinatorial structures called \emph{bunches of skew-gentle type}.
A prototypical example of such a matrix problem is the four subspace problem. 
In this section, we restrict our focus to matrix representations with invertible matrices, introduce bands and strings $\omega$ of $\fX$ and recall Bondarenko's canonical forms $M(\omega)$, which provide a complete classification of matrix representations of $\fX$ (Theorem~\ref{thm:Bondarenko}).
In the next sections, this result is  applied to the  bunch of skew-gentle type $\fXA$ and thus forms the foundation for the classification of indecomposable objects in $\DbRep{\A}$.

\subsection{Bunches of skew-gentle type and matrix problems}

A \emph{chain} is a partially ordered set such that any two elements are comparable.
By a \emph{two-point link} we mean a partially ordered set given by two incomparable elements.
\begin{dfn}\label{dfn:bunch}
	A \emph{bunch of skew-gentle type} $(\fX, \geq, \approx)$ is given by the following data.
	\begin{enumerate}
		\item A disjoint union $\fX = \fE \cup \fF$ called \emph{the set of labels}	such that \emph{the set of horizontal labels} $\fE = \bigcup_{\imath \in I} \fE_\imath$ is a disjoint union of chains indexed by a finite or countable index set $I$,
		and \emph{the set of vertical labels} $\fF = \bigcup_{\imath \in I} \fF_\imath$ is a disjoint union indexed by the same index set, where each set $\fF_\imath$ is a singleton set or a two-point link.
		\item The set of $\fX$ is endowed with a partial order $\geq$ 
		setting $x \geq y$ in $\fX$ if and only if there exists an index $\imath \in I$ such that $x,y \in \fE_\imath$
		and $x \geq y$ in $\fE_{\imath}$.
		\item The relation $\approx$ is an equivalence relation on 
		the set $\fX$ such that
		\begin{itemize}
			\item each equivalence class has cardinality one or two, 
			\item any element from a two-point link 
			is only equivalent to itself,
			\item no horizontal label is equivalent to a vertical label.
		\end{itemize}
	\end{enumerate}
	We will abbreviate a bunch of skew-gentle type $(\fX, \geq, \approx)$ by $\fX$.
	
\end{dfn}

\begin{rmk}\label{rmk:red-gen}
	The reduction method of the present article allows to reduce the classification problem of the right-bounded derived category 
	$\Dmmod{A}$
	of the arrow ideal completion of the path algebra $A$ of any skew-gentle quiver, which may have oriented cycles,
	to a matrix problem of a  bunch of skew-gentle type.
	This motivates the terminology above.
\end{rmk}

The following datum
will serve as running example for this section.
\begin{ex}\label{ex:run}
	Assume that $\fX_0 = \fE \cup \fF = (\fE_1 \cup \fE_2) \cup (\fF_1 \cup \fF_2)$ is
	given by the disjoint union of the chains
	$\fE_1 = \{
	\xa < \xb < \xc  \}$ 
	and $\fE_2 = \{\ya > \yb > \yc \}$, and
	two-point links
	$\fF_1 = \{ \xf_+, \xf_- \}$ and $\fF_2 = \{ \yf_+, \yf_- \}$.
	Let $\approx$ be the unique equivalence relation on the set $\fX$ satisfying $\xa \approx \ya$ and $\xb \approx \yb$.
	Then $\fX_0$ is a bunch of skew-gentle type.
\end{ex}

For the remainder of this section, we fix a bunch of skew-gentle type $\fX$ and a 
field $\kk$. For simplicity of the presentation, we assume that $\kk$ is algebraically closed. The case of an arbitrary field will be addressed in Remark~\ref{rmk:field2}.
\begin{dfn}\label{dfn:repX}
	\begin{enumerate}
		\item 
		A \emph{$\kk$-linear representation} of $\fX$ is given by a collection of  $\kk$-linear maps $\mu = (\begin{td} \mu_\imath\colon V_{\imath} \ar{r} \& W_{\imath} \end{td} )_{\imath \in I}$
		such that $\sum_{\imath \in I} \dim V_{\imath} + \dim W_{\imath}$ is finite, 
		$V_{\imath} = \bigoplus_{y \in \fF_{\imath}} V_{y}$ and 
		$W_{\imath} = \bigoplus_{x \in \fE_{\imath}} W_x$ with $V_b = V_y$ for any $b,y \in \fF$ with $b \approx y$ and, similarly, $W_a = W_{x}$ for any $a,x \in \fE$ with $a \approx x$. 
		\item A \emph{morphism $\begin{td} \mu \ar{r} \& \mu' \end{td}$ of $\kk$-linear representations of $\fX$} is given by 
		a collection of pairs of $\kk$-linear maps $(\sigma, \tau) = (\sigma_{\imath}, \tau_{\imath})_{\imath \in I}$
		such that for each index $\imath \in I$ the diagram
		$$
		\begin{cd}
			V_{\imath} = \bigoplus_{y \in \fF_{\imath}} V_{y} \ar{r}{\mu_\imath} \ar{d}{\tau_{\imath}=(\tau_{by})} \& W_{\imath}  = \bigoplus_{x \in \fE_{\imath}} W_x
			\ar{d}{\sigma_{\imath} = (\sigma_{ax})}
			\\
			V'_{\imath} = \bigoplus_{b \in \fF_{\imath}} V'_{b} \ar{r}{\mu'_\imath}
			\& 
			W'_{\imath}  = \bigoplus_{a \in \fE_{\imath}}W'_a
		\end{cd}
		$$
		commutes
		and the following conditions are satisfied.
		\begin{itemize}
			\item For any index $\imath \in I$ and any elements $a, x \in \fE_\imath$ with $ x \not\leq a$  it holds that $\sigma_{a x} = 0$.
			\item For any index $\imath \in I$ and any elements $b,y \in \fF_{\imath}$ with $b \neq y$ it holds that $\tau_{by} = 0$.
			\item For any elements $b, y \in \fF$ with $b \approx y$ it holds that $\tau_{bb} = \tau_{yy}$.
			\item For any elements $a, x \in \fE$ with $a \approx x$ it holds that $\sigma_{aa} = \sigma_{xx}$.
		\end{itemize}
		\item 
		The \emph{direct sum of $\kk$-linear representations} $\mu$ and $\mu'$ of $\fX$ is defined by $(\mu_\imath \oplus \mu'_{\imath})_{\imath \in I}$. 
		\item  
		A $\kk$-linear representation $\mu$ of $\fX$ is called \emph{trivial}
		if all of its vector spaces are zero, 
		\emph{indecomposable} if it is not isomorphic to a direct sum 
		of non-trivial representations, 
		and \emph{regular} if each  $\kk$-linear map $\mu_{\imath}$ is an isomorphism.	
	\end{enumerate}
\end{dfn}
The $\kk$-linear representations of a bunch of skew-gentle type $\fX$ and their morphisms form 
an additive, $\kk$-linear category $\rep{\fX}$, which has the Krull-Remak-Schmidt property by \cite{Kleiner-Roiter}.
We are interested in the problem to classify the indecomposable objects in 
the full subcategory of $\rep{\fX}$
given by regular representations, which will be denoted by 
$\regrep{\fX}$.
Choosing bases in each $\kk$-linear representation of a bunch of skew-gentle type $\fX$ allows to rephrase this classification problem for $\regrep{\fX}$ in terms of matrices.

\begin{rmk}
	A regular $\kk$-linear matrix representation of $\fX$ can be described by a family of invertible block matrices $M = (M_\imath)_{\imath \in I}$ of the following form.
	\begin{itemize}
		\item For each index $\imath \in I$ the matrix $M_\imath$ has entries from the field $\kk$ and finitely many rows and columns.
		\item For each index $\imath \in I$ the rows of the matrix $M_\imath$ are divided into horizontal stripes which are labeled by elements from $\fE_\imath$.	
		Similarly, the columns of the matrix $M_\imath$ are divided into vertical stripes labeled by elements from $\fF_\imath$.
		\item For any distinct elements $x,y \in \fX$ 
		such that $x \approx y$ the horizontal stripes labeled by $x$ and $y$ have the same number of rows.
	\end{itemize}
	Two matrix representations $M$ and $M'$ of the bunch of skew-gentle type $\fX$ are isomorphic if and only if $M$ can be obtained from $M'$ by a sequence of transformations of any of the following types.
	\begin{enumerate}
		\item \emph{Independent transformations within stripes}.
		Let $x \in \fX$ be an element which is only equivalent to itself.
		Let $\imath \in I$ denote the index such that $x \in \fE_\imath \cup \fF_\imath$.
		\begin{itemize}
			\item If $x \in \fE_\imath$, we may perform any elementary transformation of rows in the horizontal stripe labeled by $x$.
			\item If $x \in \fF_{\imath}$, we may perform any elementary transformation of columns in the vertical stripe corresponding to $x$.
		\end{itemize}
		\item \emph{Simultaneous transformations of equivalent stripes}.
		Assume that $x,y \in \fX$ are two distinct, equivalent elements. Then either 
		$x \in \fE_{\imath}$ and $y \in \fE_{\jmath}$ or $x \in \fF_{\imath}$ and $y \in \fF_{\jmath}$ for certain indices $\imath, \jmath \in I$.
		\begin{itemize}
			\item In the first case,
			we may perform any elementary transformation of rows in the horizontal stripe with label $x$ in the matrix $M_\imath$ and the \emph{same} transformation of rows in the horizontal stripe with label $y$ in the matrix $M_\jmath$.   
			\item Similarly, in the second case, we may perform 
			any elementary column transformation in the vertical stripe with label $x$ in the matrix $M_\imath$ and the \emph{same} column transformation in the vertical stripe with label $y$ in the matrix $M_\jmath$.   
		\end{itemize}
		\item \emph{Independent transformations between stripes}.
		Let $x,y \in \fX$ with $x < y$. Then $x,y \in \fE_\imath$ with $\imath \in I$. 
		For any $\lambda \in \kk$ we may add 
		the $\lambda$-multiple of any row of the horizontal stripe with label $x$ to any row of the horizontal stripe with label $y$.	
	\end{enumerate}
	The notions of direct sum and indecomposability of
	representations translate into matrix language in a straightforward way.
	Using this translation, the problem to classify the indecomposable objects in $\regrep{\fX}$ 
	can be formulated as a matrix problem, that is, a problem of linear algebra.
\end{rmk}

\begin{ex}
	An isomorphism of $\kk$-linear representations $\begin{td} (\mu_1,\mu_2) \ar{r}{\sim} \& (\mu'_1,\mu'_2) \end{td}$ of
	the bunch $\fX_0$ from Example~\ref{ex:run}
	is given 
	by the four invertible vertical maps
	\begin{align*}
		\begin{cd}
			V_{\xf_+} \oplus V_{\xf_-} \ar{r}{\mu_1} \ar{d}{\tau_{1}=
				\left(
				\begin{smallmatrix}
					* & 0 \\
					0 & *
				\end{smallmatrix}\right)}
			\& W_{\xa} \oplus W_{\xb} \oplus W_{\xc} 
			\ar{d}{\sigma_{1} = 
				\begin{psmallmatrix}
					\sigma'_{1} & 0 & 0 \\
					* & \sigma''_{1} & 0 \\
					* & * & *	\mathstrut
			\end{psmallmatrix}}
			\\
			V'_{\xf_+} \oplus V'_{\xf_-} \ar{r}{\mu'_1}
			\& 
			W_{\xa} \oplus W_{\xb} \oplus W_{\xc} 
		\end{cd}
		&&
		\begin{cd}
			V_{\yf_+} \oplus V_{\yf_-} \ar{r}{\mu_2} \ar{d}{\tau_{2}=
				\begin{psmallmatrix}
					* & 0 \\
					0 & *
			\end{psmallmatrix}}
			\& W_{\xa} \oplus W_{\xb} \oplus W_{\yc} 
			\ar{d}{\sigma_{2} = 
				\begin{psmallmatrix}
					\sigma'_{1} & * & * \\
					0 & \sigma''_{1} & * \\
					0 & 0 & *	\mathstrut
			\end{psmallmatrix}}
			\\
			V'_{\yf_+} \oplus V'_{\yf_-} \ar{r}{\mu'_2}
			\& 
			W_{\xa} \oplus W_{\xb} \oplus W_{\yc} 
		\end{cd}
	\end{align*}
	where each symbol $*$ can be an arbitrary $\kk$-linear map.
	Choosing bases for the vector spaces allows to depict the classification problem 
	for $\rep{\fX_0}$ by the following diagram.
	\begin{align}\label{eq:mp-run}
		\begin{array}{c}
			\begin{tikzpicture}
				\node () at (-3,0) {
					$\begin{NiceArray}{|c|c|}[first-row,first-col,rules={width=1pt}]
						& \xf_+ &   \xf_-  \\
						\hline
						\xa & \rowcolor{c2} &    \\ \hline
						\xb &
						\rowcolor{c2}  &  \\ \hline
						\xc &
						\rowcolor{c1} &   \\ \hline
					\end{NiceArray}$
				};
				\node () at (3,0) {	
					$\begin{NiceArray}{|c|c|}[first-row,last-col,rules={width=1pt}]
						\yf_+ &  \yf_-  \\
						\hline
						\rowcolor{c2}   &  & \ya  \\ \hline
						\rowcolor{c2} &  &  \yb \\ \hline
						\rowcolor{c1}  &&   \yc
						\\ 
						\hline
					\end{NiceArray}$
				};
				\draw[-, dashed,color=c4](-2, 0.25) -- (2, 0.25);
				\draw[-, dashed,color=c4](-2, -0.25) -- (2, -0.25);
				\draw[->, thick, color=c4](-4.5, 0.25) -- (-4.5, -1);
				\draw[<-, thick, color=c4](4.5, 0.25) -- (4.5, -1);
			\end{tikzpicture}
		\end{array}
	\end{align}
	If $\xa$ and $\ya$ are the only horizontal stripes with non-zero rows, we obtain a representation of the four subspace quiver.
	In particular, the matrix problem of $\fX_0$, that is, the problem to classify its indecomposable representations up to isomorphism, generalizes the four subspace problem.
	On the other hand, the problem to classify all regular indecomposable representations of $\fX_0$ is still representative
	for the complexity of the matrix problem of the bunch of skew-gentle type $\fXA$ from Definition~\ref{def:chkX}.
\end{ex}
\begin{rmk}
	Definition~\ref{def:chkX} yields an example of a bunch of skew-gentle type $\fXA$.
	In particular, Definition~\ref{dfn:repX-A}
	is the specialization of Definition~\ref{dfn:repX} to
	$\fXA$.
	Subsection~\ref{subsec:mp}
	provides a description for the classification problem of $\regrep{\fXA}$ in terms of matrices
	which is possible for any bunch of skew-gentle type.
	In particular, one may perform the same transformations with the matrices in~\eqref{eq:mp-run}
	as described in Subsection~\ref{subsec:mp}.
\end{rmk}

\subsection{Alphabet of a bunch of skew-gentle type}\label{subsec:abc}

We recall that we have fixed a bunch of skew-gentle type $\fX$.
\begin{dfn}
	We associate a new 
	datum $(\kA_{\fX},\geq,\sim,-)$
	to $\fX$ where $\kA_{\fX}$ is called the \emph{set of letters of $\fX$} endowed with a partial order $\geq$ and symmetric relations $\sim$ and $-$. 
	\begin{itemize}
		\item The set $\kA_{\fX}$ is obtained by `collapsing each two-point link' in the set $\fX$ as follows. 
		Each two-point link $\fF_{\imath}= \{f_{\imath,+}, f_{\imath,-}\}$ in the union $\fF = \bigcup_{\imath} \fF_{\imath}$
		is replaced by a singleton set $\bar{\fF}_{\imath}=\{f_{\imath}\}$.
		This yields new sets $\bar{\fF}$ and $\kA_\fX = \fE \cup \bar{\fF}$.
		\item The set $\kA_{\fX}$ inherits the partial order from $\fX$, that is, two distinct elements $x,y \in \kA_{\fX}$ are comparable if and only if $x,y \in \fE_{\imath}$ for an index $\imath \in I$. 
		\item 
		An element $f \in \kA_{\fX}$ obtained from collapsing a two-point link is called \emph{special}. For any special element $f $ we set $f \sim f$. 
		Any distinct elements $x,y\in \fX$ with $x \approx y$ in $\fX$
		exist also in $\kA_{\fX}$, and we set $x \sim y$.
		In this case, each of the symbols $x$ and $y$ is called \emph{tied}.
		\item For any $x, y \in \kA_{\fX}$ we set $x - y$ if there is an index $\imath \in I$ such that $x \in \fE_{\imath}$ and $y \in \bar{\fF}_{\imath}$ or vice versa.	
	\end{itemize}
	We will abbreviate $(\kA_{\fX},\geq,\sim,-)$ by $\kA_{\fX}$.
\end{dfn}

\begin{ex}
	In case of the bunch $\fX_0$ from Example~\ref{ex:run}
	the set of letters is given by $\kA_{\fX_0} = \fE_1 \cup \fE_2 \cup \{ \xf \} \cup \{ \yf \}$,
	the symmetric relation $\sim$ on $\kA_{\fX_0}$ is generated by
	$\xa \sim \ya$, $\xb \sim \yb$, $\xf \sim \xf$, $\yf \sim \yf$, and the symmetric relation $-$ by $\xf - e_1$
	for any $e_1 \in \fE_1$ and $\yf - e_2$ for any $e_2 \in \fE_2$.
	In particular, the elements $\xc$ and $\yc$ are neither tied nor special. 
\end{ex}

\begin{dfn}\label{def:cyc-word}
	Let $w = (x_i \rho_i)_{i \in \Z}$ be a sequence 	such that 
	\begin{itemize}
		\item each symbol $x_i$ is given by a letter from	$\kA_{\fX}$, 
		\item each symbol $\rho_i$ is given by a relation $-$ 
		if $i$ is even, respectively by $\sim$ if $i$ is odd, 
		\item the	relation $x_{i}\rho_i x_{i+1}$ is satisfied in $\kA_{\fX}$ for each index $i \in \Z$.
	\end{itemize}
	For the sequence $w$
	we make the following definitions.
	\begin{enumerate}
		\item For any integer $j \in 2\Z$ the \emph{$j$-th shift of the sequence $w$} is defined by 
		$w_j = (x_{i+j}\rho_{i+j})_{i \in \Z}$.
		\item The \emph{opposite of the sequence $w$} is set to $w^{\op} = (x_{1-i} \rho_{-i})_{i \in \Z}$.
		\item The word $w$ is \emph{symmetric} if 
		there is an integer $j\in 2\Z$ such that $w_j = w^{\op}$.
		\item The sequence $w$ is a \emph{periodic word} if there exists a number $\ell \in 2\N$ such that $w = w_\ell$. 
		Assuming $\ell$ to be minimal, it is called the \emph{period of $w$} and we denote
		$\dot{w} = (x_1 \rho_1 x_2 \rho_2 \ldots x_{\ell} \rho_{\ell} )$.
	\end{enumerate}	
\end{dfn}

\begin{ex}	\label{ex:cyc1}
	An example of a periodic word of the bunch of skew-gentle type $\fX_0$ from Example~\ref{ex:run} is given by
	the $8$-periodic sequence
	\begin{align*}
		w&= 
		{\color{blue}
			\ldots
			\xf \sim \xf - \xa \sim \ya - \yf \sim \yf - \ya \sim \xa \, \underset{\rho_0}{-} \,} \xf \sim \xf - \xa \sim \ya - \yf \sim \yf - \ya \sim \xa -  {\color{blue} \ldots }.
		\intertext{
			Shifting $w$ by $2$ corresponds to moving $w$ two steps to the left, while
			taking the opposite of $w$ corresponds to a vertical flip with axis of symmetry at $\rho_0$.
			On the other hand,
		}
		w^{\op}&= {\color{blue}\ldots 
			\xa \sim \ya - \yf \sim \yf - \ya \sim \xa - \xf \sim \xf  \, \underset{\rho_0}{-} \, }
		\xa \sim \ya - \yf \sim \yf - \ya \sim \xa - \xf \sim \xf -  {\color{blue} \ldots }
	\end{align*}
	is equal to $w_2$,
	and thus the periodic word $w$ is symmetric.
\end{ex}

\begin{dfn}\label{def:band}
	A \emph{band} $\omega$ of $\fX$ is a triple $(w,m,\lambda)$
	given by a cyclic word $w$, a positive integer $m \in \N$, called the \emph{multiplicity of the band $\omega$}, and a scalar $\lambda \in \kk^*$ such that $\lambda \neq (-1)^{ \frac{t(w)}{2} +1}$ if $w$ is symmetric of period $\ell$,
	where $t(w)$
	denotes the number of subwords $x_i \sim x_{i+1}$ with tied letters, that is, $x_i \neq x_{i+1}$, in the periodic part $\dot{w}$.
\end{dfn}
We note that the number $t(w)$ is even for any symmetric periodic word $w$.

\begin{ex}\label{ex:cyc2}
	For any $m \in \N$ and $\lambda \in \kk^*$
	the triple $(w,m,\lambda)$ defines a band of $\fX_0$ assuming
	\begin{enumerate}
		\item $w$ is the periodic word with $
		\dot{w} = 
		(\xf \sim \xf - \xa \sim \ya - \yf \sim \yf - \ya \sim \xa -
		\xf \sim \xf - \xa \sim \ya - \yf \sim \yf - \yb \sim \xb \, -),$ or
		\item $w$ is the symmetric word from Example~\ref{ex:cyc1} and $\lambda \neq -1$.
	\end{enumerate}	
	In fact, for any $\ell$-periodic word $w$ of the bunch $\fX_0$ it holds that $t(w) = \frac{\ell}{2}$.
\end{ex}

\begin{dfn}\label{dfn:fin-word}
	A \emph{finite word}  of ${\fX}$ is given by
	$w = (x_1 \rho_1 x_2  \ldots \rho_{\ell-1} x_\ell)$
	such that
	\begin{itemize}
		\item each symbol	$x_i$ is a letter from $\kA_{\fX}$, 
		\item each symbol $\rho_i$ is one of the relations $\sim$ or $-$, and $\rho_{i-1} \neq \rho_{i}$ for each index $1 < i < \ell$,
		\item the relation $x_i \rho_i x_{i+1}$ is satisfied for each index $1 \leq i < \ell$,
		\item \emph{regularity conditions:} neither $x_1$ nor $x_{\ell}$ is tied and $\ell \geq 2$.
	\end{itemize}
	The number $\ell$ is the \emph{length} of the word $w$. Furthermore, we use the following terminology.
	\begin{enumerate}
		\item The \emph{opposite of $w$} is defined by the word $w^{\op} = (x_{\ell} \rho_{\ell-1} \ldots x_2 \rho_1 x_1)$. 
		\item 
		The word $w$ is \emph{symmetric} if $w = w^{\op}$.
		\item
		The left end of the word $w$ is called \emph{special} if $x_1$ is special and $\rho_1$ is given by $-$.
		Similarly, the right end of $w$ is called \emph{special} if $x_\ell$ is special and $\rho_{\ell-1}$ is given by $-$. The word $w$ is \emph{usual}, \emph{special} respectively \emph{bispecial} if it has zero, one respectively two special ends.
	\end{enumerate}
\end{dfn}
Similar to periodic words, taking the opposite of a finite word corresponds to a vertical flip.
Any special word is already asymmetric.

\begin{dfn}\label{dfn:fin-word2}
	A string of $\fX$ is given by one of the following.
	\begin{enumerate}
		\item
		A \emph{usual string} is given by a usual, asymmetric word $w$. 
		\item A \emph{special string} $(w,\varepsilon)$ is given by a special word $w$ and 
		a sign $\varepsilon \in \{+,-\}$. 
		\item
		A \emph{bispecial string}
		$(\varepsilon_1, w, \varepsilon_2)$ is given by a bispecial word $w$ and signs $\varepsilon_1,\varepsilon_2 \in \{+,-\}$. 
	\end{enumerate}
\end{dfn}

\begin{ex}\label{ex:str}
	The following data define a usual, a special and a bispecial string  of the bunch $\fX_0$ from Example~\ref{ex:run}.
	\begin{enumerate}
		\item $w= (\xc - \xf \sim \xf - \xa \sim \ya - \yf \sim \yf - \yc)$. 
		\item $(w,\varepsilon)$
		with $w =
		(\yc - \yf \sim \yf - \ya \sim \xa - \xf \sim \xf - \xa \sim \ya - \yf)$.
		\item \label{bspl-1a} $(-,w,+)$ with 
		$w = (\xf - \xa \sim \ya - \yf \sim \yf - \yb \sim \xb - \xf 
		\sim \xf - \xb \sim \yb - \yf \sim \yf - \ya \sim \xa - \xf)$.
	\end{enumerate}
\end{ex}

\subsection{Matrix representations of bands and strings}
\label{subsec:strbdrep}
Throughout this subsection, let $\omega$ be a string or band of $\fX$ and let $w$ denote its underlying word. 
The next goal is to define a $\kk$-linear representation
$M(\omega)$
of $\fX$.
This construction requires several preparatory steps involving mostly the word $w$.
\begin{enumerate}
	\item Each special letter $f$ is changed into $f_+$ or $f_-$.
	\item Each subword $f_+ \sim f_-$ or $f_- \sim f_+$ in $w$ is decorated by a horizontal arrow.
	\item In the final step, each horizontal arrow
	and each relation $-$ in $w$ give rise to non-zero entries in the matrix representation $M(\omega)$.
\end{enumerate}
Each of these steps is carried out in one of the next subsections.
The first step determines the vector spaces in the representation
$M(\omega)$, while the second determines its matrices.
In fact, the first two steps can be interchanged and the second one is the more technical.

\subsubsection{Signs of special letters}\label{subsub:signs}
In case the word $w$ is bispecial, we need first to identify its  `smallest asymmetric' part for the first step.
\begin{dfn}
	\begin{enumerate}
		\item
		Given a bispecial word $v$ of length $\ell$ of $\fX$ and a number $m \in \N$,
		the \emph{$m$-th power of $v$} is defined as the word $v^m = v \sim v^{\op} \sim v \sim v^{\op} \ldots \sim $
		of length $\ell m$	which ends with $v$ if $m$ is odd, respectively with $v^{\op}$ if $m$ is even.
		\item The \emph{primitive root} of a bispecial word $w$ of $\fX$ is the shortest possible subword $v$ 
		satisfying $w = v^m$ for a number $m \in \N$, which is called the \emph{multiplicity} of 
		$v$ in 
		$w$.
	\end{enumerate}
\end{dfn}
\begin{rmk}
	As the formulation suggests, any bispecial word of $\fX$ has a unique primitive root with a well-defined multiplicity. A bispecial word $w$ is symmetric if and only if 
	its primitive root has even multiplicity in $w$.
	In this sense, primitive roots may be viewed as bispecial words satisfying a stronger notion of asymmetry.
\end{rmk}
With the above preparations, each special letter of $w$ is given a sign as follows.
\begin{itemize}
	\item If $\omega$ is not a bispecial string, any special subword $f \sim f$ is changed to $f_+ \sim f_-$.
	\item If $\omega$ is a special string $(w,\pm)$ then the special end $f$ of $w$ is changed to $f_{\pm}$.
	\item 
	If $\omega$ is a bispecial string $(\varepsilon_1, w,\varepsilon_2)$, then $w = v \sim v^{\op} \sim v \ldots v_m$ is 	the $m$-th power of 
	its primitive root $v$ for certain $m \in \N$.
	\begin{itemize}
		\item The left end $x_1$ of any copy of $v$ in $w$ is decorated by $\varepsilon_1$, any special subword $f\sim f$ in $v$ is changed into $f_+ \sim f_-$
		and the right end $x_{\ell}$ of $v$ is decorated by $\varepsilon_2$. 
		\item 	In every copy of $v^{\op}$ in $w$, we choose the opposite sign $\ol{\varepsilon}_2$ of $\varepsilon_2$ at its left end $x_\ell$, $f_- \sim f_+$ for special subwords in between, and $\ol{\varepsilon}_1$ at its right end $x_1$. 
	\end{itemize}
	Thus, the right end of $w$ obtains sign $\varepsilon_2$ if $m$ is odd, respectively $\ol{\varepsilon}_1$ if $m$ is even.
\end{itemize}
The resulting `polarized' word 
is a sequence of relations and labels from the set underlying the bunch $\fX$. 

\begin{ex}\label{ex:bispl}
	The primitive root of the bispecial word $w$ from Example~\ref{ex:str}
	is given by 
	$v =  
	\xf - \xa \sim \ya - \yf \sim \yf - \yb \sim \xb - \xf$
	and has multiplicity two in $w$, that is, $w = v \sim v^{\op}$.
	Because of the signs in the bispecial string $(+,w,-) $  the polarized version of $w$ is given by
	\begin{align*} 
		\xf_+ - \xa \sim \ya - \yf_+ \sim \yf_- - \yb \sim \xb - \xf_- 
		\sim \xf_+ - \xb \sim \yb - \yf_- \sim \yf_+ - \ya \sim \xa - \xf_-.
	\end{align*}
\end{ex}

\subsubsection{Orienting special subwords}
\label{subsec:A-orient}
As the next step does not depend on the signs of special letters determined in the previous one, we suppress these in the notation of this subsection. 
To orient special subwords, we need to consider a canonical extension of a word.
\begin{dfn}\label{dfn:ambient}
	For a finite or periodic word $w$ 
	\emph{its ambient word} 
	is defined by
\begin{align*}
	\overline{w} = \begin{cases}
		\begin{array}{ll}
			\phantom{w^{\op} \sim \ }w
			& \text{if $w$ is a periodic or usual word,}\\
			\phantom{w^{\op} \sim \ } w \sim w^{\op}
			& \text{if $w$ is a special word with special right end,}\\
			w^{\op} \sim w
			& \text{if $w$ is a special word with special left end,}\\
			\phantom{w^{\op} \sim \ }w_\infty
			& \text{if $w$ is bispecial,}
		\end{array}
	\end{cases}
\end{align*}
where $w_{\infty}$ is the periodic word of $\fX$ with periodic part $\dot{w}_\infty = w$.
\end{dfn}
In the following, we denote by  $(x_i)_{i \in I_{\overline{w}}}$ the sequence of letters in $\overline{w}$, where the index set $I_{\overline{w}}$
of $\overline{w}$ is the same index set as that of $w$ in the first case, the integer interval $\{1,2, \ldots 2\ell\}$ in the second, $\{-\ell-1, \ldots \ell\}$ in the third respectively $\Z$ in the last case.

To decorate each subword $x_i \rho_i x_{i+1} = f \sim f$ with $f \in \fF$ in $w$ with a left or right orientation $\overleftrightarrow{x_i \sim x_{i+1}}$, we consider the maximal symmetric subword of the form $v^{\op} \rho_i v$ in $\overline{w}$. Then one of the following occurs.
\begin{itemize}
\item It holds that $v^{\op} \rho_i v = \overline{w}$, that is, $x_{i-p} = x_{i+1+p}$ for any $p \in \N$.
This may only occur if $\overline{w}$ is a periodic symmetric word with an axis of symmetry at $\tilde{r}_i$, that is,
$$
\overline{w} =(
\underbrace{\quad \ldots \quad -x_{i-2} \sim x_{i-1} - x_i}_{v^{\op}} \sim \underbrace{x_{i+1} - x_{i+2} \sim  x_{i+3} - \quad \ldots \quad}_{v}).
$$
In this case, we set
$\overleftarrow{x_i \sim x_{i+1}}$ in $w$.
\item Otherwise, there exists a minimal number $q = q(i) \in \N$ such that
$x_{i-q} \neq x_{i+q+1}$. 
$$
\overline{w} =(
\ldots
{x_{i-q}} - 
\underbrace{x_{i-q+1} \ldots \quad  \sim x_{i-1} - x_i}_{v^{\op}} \sim \underbrace{x_{i+1} - x_{i+2} \sim  \quad \ldots x_{i+q}}_{v} - {x_{i+q+1}} \ldots).
$$ 
Then there is an index $i \in I$ such that
the letters $l_i = x_{i-q}$ and $r_i = x_{i+q+1}$ belong to the same chain $\fE_\imath$.
In this case, we set
$$
\overleftarrow{x_i \sim x_{i+1}}\quad \text{if }l_i < r_i \qquad
\text{respectively} \qquad
\overrightarrow{x_i \sim x_{i+1}} \quad \text{if }l_i > r_i.
$$
Summarized in terms of a mnemonic, 
we orient the arrow \emph{towards the smaller letter outside the symmetric core}.
More precisely, we ignore the maximal symmetric subword $u = v^{\op} \rho_i v$, compare the first letter $l_i$ left from $u$
to the first letter $r_i$ right from $u$
and choose the orientation in $\overleftrightarrow{x_i \rho_i x_{i+1}}$ towards the smaller of both letters.
\end{itemize}
The outcome is denoted by $\olra{w}$.
We note that for a periodic word $w$ it is sufficient to describe the signs and oriented arrows in the periodic part $\dot{w}$, while for bispecial word $w$ signs and arrows are determined by its primitive root.

\begin{ex}\label{ex:or}
We describe the signs and oriented arrows for the words underlying the bands and strings 
from Examples~\ref{ex:cyc2} and Examples~\ref{ex:str}.
\begin{enumerate}
\item 
For the previous periodic symmetric word $w$ we obtain using the symmetry rule that
$
\olra{w} =  (\overleftarrow{\xf_+ \sim \xf_-} - \xa \sim \ya - \overleftarrow{\yf_+ \sim \yf_-} - \ya \sim \xa - )$.
\item In case of the previous asymmetric periodic word $w$,
because of $\xb>\xa$ and $\ya > \yb$ 
it follows that $$\olra{w} =
(\overrightarrow{\xf_+ \sim \xf_-} - \xa \sim \ya - \overrightarrow{\yf_+ \sim \yf_-} - \ya \sim \xa -
\overrightarrow{\xf_+ \sim \xf_-} - \xa \sim \ya - \overrightarrow{\yf_+ \sim \yf_-} - \yb \sim \xb -
).$$
\item Similarly, for the referenced usual word $w$ 
its oriented version with signs is given by
$\olra{w} = (\xc - \overrightarrow{\xf_+ \sim \xf_-} - \xa \sim \ya - \overrightarrow{\yf_+ \sim \yf_-} - \yc)$. 
\item For the special string $(w,+)$ from above,
we consider the ambient word
\begin{align*}
	\overline{w} =
	\underbrace{\yc - \yf \sim \yf - \ya \sim \xa - \xf \sim \xf - \xa \sim \ya 
		- \yf}_{w} \ {\color{blue} \sim \underbrace{\yf - \ya \sim \xa \quad \ldots \quad  - \yc}_{w^{\op}}}.
\end{align*}	
Each special subword in $w$ is decorated by a left arrow because of $\yc < \ya$, that is,
$\olra{w} = 
(\yc - \overleftarrow{\yf_+ \sim \yf_-} - \ya \sim \xa - \overleftarrow{\xf_+ \sim \xf_-} - \xa \sim \ya - \yf_+)$.
\item Continuing Example~\ref{ex:bispl} 
with  the bispecial string $(+,w,-)$, we obtain
\begin{align}  
	\olra{w} =  \xf_+ - \xa \sim \ya - \overrightarrow{\yf_+ \sim \yf_-} - \yb \sim \xb - \overleftarrow{\xf_- 
		\sim \xf_+} - \xb \sim \yb - \overleftarrow{\yf_- \sim \yf_+} - \ya \sim \xa - \xf_-.
\end{align}
\end{enumerate}
\end{ex}

\subsubsection{Construction of canonical forms}\label{subsec:can}

Given the version of $w$ with oriented arrows and signs, next we construct a diagram of $\kk$-vector spaces and $\kk$-linear maps
together with certain unoriented edges, which is then `folded' into  
a matrix representation $M(\omega)$ of $\fX$.
We distinguish between strings and bands in this process.
\begin{itemize}
\item Assume that $\omega$ is a string.
\begin{enumerate}
\item \label{diag} \emph{Construction of the diagram}.
Any
non-tied 
label $a$ of $w$ gives rise to a one-dimensional vector space denoted by $\kk_{a}$.
Any tied subword $x \sim x'$ of $w$ yields an unoriented edge $\begin{td} \kk_{x}\ar[-,densely dotted]{r} \& \kk_{x'} \end{td}$ in the diagram.
Any subword of the form $y -x$, $x - y$, $\overrightarrow{y \sim z} - x$
or  $x - \overleftarrow{z \sim y}$
with $x \in \fE$ and $y \in \fF$ gives rise to a map $\begin{td} \kk_y \ar{r}{1} \& \kk_x \end{td}$
in the diagram.
\item \emph{Folding the diagram}.
For each label $x \in \fF$ the vector space $V_x$ of $\mu(\omega)$ is given by the direct sum of one-dimensional vector spaces corresponding to all letters in $w$ equal to $x$. Similarly, for each label $y \in \fE$ the vector space $W_y$ is the direct sum of 
copies of $\kk$ corresponding to letters $y$ in $w$.
The diagram arrows of the previous step yield all
non-zero $\kk$-linear maps
$\begin{td} \mu_{yx} \colon V_x \ar{r} \& W_y \end{td}$ in the representation $\mu(\omega)$.
The remaining edges indicate the vector spaces which have to be identified.
\end{enumerate}
This completes the definition of the string representation $\mu(\omega)$.
\item Assume that $\omega$ is a band $(w,m,\lambda)$.
\begin{enumerate}
\item \emph{Construction of the diagram}.
We carry out step~\eqref{diag} above for the periodic part $\dot{w} = (x_i \rho_i)_{i=1}^{\ell}$ together with the following changes.
Each $1$-dimensional space is replaced by an $m$-dimensional vector space and
each identity map by the identity matrix $\Id$ of size $m$.
Moreover, we view $\rho_\ell$ in $\dot{w}$ as a relation connecting the letters $x_{\ell}$ and $x_1$, which gives rise to two or one more arrows in the diagram as follows.
\begin{itemize}
	\item 
	If there is a subword $\overrightarrow{x_{\ell-1} \sim x_{\ell}} - x_1$ 
	or $x_{\ell} - \overleftarrow{x_1 \sim x_2}$ 	
	we use	the Jordan block $\JB_\lambda$ of size $m$ with eigenvalue $\lambda$
	for the induced $\kk$-linear map 
	$
	\begin{td} \kk^m_{x_{\ell-1}} \ar{r} \& \kk^m_{x_{1}} \end{td}$ 
	respectively
	$\begin{td} \kk^m_{x_{\ell}} \ar{r} \&  \kk^m_{x_{2}} \end{td}$,
	and the identity matrix $\Id$ for the $\kk$-linear map $\begin{td} \kk^m_{x_{\ell}} \ar{r} \& \kk^m_{x_1} \end{td}$.
	\item Otherwise 
	the Jordan block $\JB_\lambda$ defines the $\kk$-linear map  $\begin{td} \kk^m_{x_{\ell}} \ar{r} \& \kk^m_{x_1} \end{td}$.
\end{itemize}
\item \emph{Folding the diagram}. The second step is carried out as for strings.
\end{enumerate}
This completes the definition of the band representation $\mu(\omega)$.
\end{itemize}
In both cases, the datum $\mu(\omega)$ is indeed a representation of $\fX$ because 
we may identify the vector spaces with labels $x$ and $x'$ for any tied subword $x \sim x'$ of $\olra{w}$.

\subsection{Examples of canonical forms}
We describe the matrix representations of the bands and strings from Example~\ref{ex:or}.
\paragraph{A usual string}  For the usual string
with oriented word $\olra{w} = (\xc - \overrightarrow{\xf_+ \sim \xf_-} - \xa \sim \ya - \overrightarrow{\yf_+ \sim \yf_-} - \yc)$
we obtain the diagram
$$
\begin{td}
\kk_{\xc} \ar[<-]{r}{1} \& \kk_{\xf_+} \ar[bend right]{rr}[swap]{1} \& \kk_{\xf_-}  \ar{r}{1} \& \kk_{\xa} \ar[densely dotted,-]{r}\& \kk_{\ya}  \ar[<-]{r}{1} \&  \kk_{\yf_+} 
\ar[bend right]{rr}[swap]{1}
\& \kk_{\yf_-} \ar{r}{1} \& \kk_{\yc}
\end{td}
$$
Folding the diagram yields the representation 
$M(\omega)$
given by 
$$
\begin{array}{cc}
\begin{array}{c}	
\begin{cd}
\kk_{\xf_+} \oplus \kk_{\xf_-} \ar{r}{M_1} \& \kk_{\xa} \oplus \kk_{\xc}
\end{cd} \\
\begin{cd}
\kk_{\yf_+} \oplus \kk_{\yf_-} \ar{r}{M_2} \& \kk_{\xa} \oplus \kk_{\yc} \end{cd}
\end{array}
&
\begin{tikzpicture}[baseline=0pt,yshift=0.25cm]
\node () at (-4,-0.25)	{$(M_1,M_2) = $};
\node (C) at (-1.5,0) {		
$\begin{pNiceArray}{c|c}[first-row,first-col]
&\xf_+ & \xf_-  \\
\xa & 1 & 1   \\
\hline
\xc & 1 & 0  
\end{pNiceArray}$};
\node (D) at (1.5,0) {	
$
\begin{pNiceArray}{c|c}[first-row,last-col]
\yf_+ & \yf_-& \\
1 & 0 & \ya \\
\hline
1 & 1 & \yc
\end{pNiceArray}$
};
\draw[-, dashed, color=c4](-0.5, 0) -- (0.5, 0);
\draw[->, thick, color=c4](-2.75, 0) -- (-2.75, -0.75);
\draw[<-, thick, color=c4](2.75, 0) -- (2.75, -0.75);
\end{tikzpicture}
\end{array}.
$$
\paragraph{A symmetric band}
For the band $(\omega,1,\lambda)$ with
$\olra{w} =  (\overleftarrow{\xf_+ \sim \xf_-} - \xa \sim \ya - \overleftarrow{\yf_+ \sim \yf_-} - \ya \sim \xa - )$ we obtain the cyclic diagram
$$
\begin{td}
\mathstrut  \ar[-,densely dotted]{r} \&		\kk_{\xa} \&
\kk_{\xf_+} \ar{l}{1}
\&
\kk_{\xf_-} \ar[bend left]{ll}{\lambda} 
\ar{r}{1} \&
\kk_{\xa} \ar[-,densely dotted]{r}\& \kk_{\ya} 
\& \kk_{\yf_+} \ar{l}[swap]{1} \&
\kk_{\yf_-} \ar[bend left]{ll}{1}
\ar{r}{1}
\&
\kk_{\ya}  \ar[-,densely dotted]{r} \& \mathstrut 
\end{td}.
$$
Folding this diagram yields the following matrix representation.
\begin{align}\label{eq:bd-sym}
\begin{array}{cc}
\begin{array}{c}
\begin{cd}
\kk_{\xf_+} \oplus \kk_{\xf_-} \ar{r}{M_1} \& \kk_{\xa}^2
\end{cd}
\\
\begin{cd}
\kk_{\yf_+} \oplus \kk_{\yf_-} \ar{r}{M_2} \& \kk_{\xa}^2
\end{cd}
\end{array}&
\begin{tikzpicture}[baseline=0pt,yshift=0.25cm]
\node () at (-4,-0.25)	{$(M_1,M_2) = $};		
\node () at (-1.5,0) {		
$\begin{pNiceArray}{c|c}[first-row,first-col]
	&\xf_+ & \xf_-  \\
	\xa & 0 & 1   \\
	\xa & 1 & \lambda  
\end{pNiceArray}$};
\node () at (1.5,0) {	
$
\begin{pNiceArray}{c|c}[first-row,last-col]
	\yf_+ & \yf_-& \\
	1 & 1 & \ya \\
	0 & 1 & \ya
\end{pNiceArray}$
};
\draw[-, dashed, color=c4](-0.5, -0.25) -- (0.5, -0.25);
\end{tikzpicture}
\end{array}
\end{align}

\paragraph{An asymmetric band}
For the band $(w,m,\lambda)$ with any $m\in \N$, $\lambda \in \kk^*$ and
$$\olra{w} =
(\overrightarrow{\xf_+ \sim \xf_-} - \xa \sim \ya - \overrightarrow{\yf_+ \sim \yf_-} - \ya \sim \xa -
\overrightarrow{\xf_+ \sim \xf_-} - \xa \sim \ya - \overrightarrow{\yf_+ \sim \yf_-} - \yb \sim \xb -
)$$
the
matrix representation is given by the following square matrices of size $4m$.
\begin{align*}
\begin{tikzpicture}[baseline=0pt,yshift=0.25cm]
\node () at (-6,-0.25)	{$(M_1,M_2) = $};		
\node () at (-2.5,0) {		
$\begin{pNiceArray}{cc|cc}[first-row,first-col]
&\xf_+ & \xf_+ & \xf_- & \xf_- \\
\xa& \Id & 0 & \Id & 0  \\
\xa&0 & \Id & 0 & 0  \\
\xa&0 & \Id & 0 & \Id \\
\hline
\xb&\JB_{\lambda} & 0 & 0 & 0 
\end{pNiceArray}$};
\node () at (2.5,0) {	
$\begin{pNiceArray}{cc|cc}[first-row,last-col]
\yf_+ & \yf_+ & \yf_- & \yf_- \\
\Id & 0 & 0 & 0 & \ya \\
\Id & 0 & \Id & 0 & \ya \\
0 & \Id & 0 & 0 & \ya \\
\hline
0 & \Id & 0 & \Id & \yb 
\end{pNiceArray}$};
\draw[-, dashed, color=c4](-0.75, -0) -- (0.75, -0);
\draw[-, dashed, color=c4](-0.75, -1.2) -- (0.75, -1.2);
\draw[->, thick, color=c4](-4.5, 0.6) -- (-4.5, -1.25);
\draw[<-, thick, color=c4](4.5, 0.6) -- (4.5, -1.25);
\end{tikzpicture}
\end{align*}

\paragraph{A special string}
For the special string $(w,+)$ with the oriented word
$\olra{w} = 
(\yc - \overleftarrow{\yf_+ \sim \yf_-} - \ya \sim \xa - \overleftarrow{\xf_+ \sim \xf_-} - \xa \sim \ya - \yf_+)$ yields the representations.
\begin{align*}
\begin{tikzpicture}[baseline=0pt]
\node () at (-4,-0.25)	{$(M_1,M_2) = $};		
\node () at (-1.5,0) {		
$\begin{pNiceArray}{c|c}[first-row,first-col]
&\xf_+ & \xf_-  \\
\xa & 1 & 0  \\
\xa & 1 & 1  
\end{pNiceArray}$};
\node () at (2.2,-0.225) {	
$
\begin{pNiceArray}{cc|c}[first-row,last-col]
\yf_+ & \yf_+ & \yf_- & \\
0 & 0 & 1 & \ya \\
0 & 1 & 0 & \ya \\
\hline
1 & 0 & 1& \yc
\end{pNiceArray}$};
\draw[-, dashed, color=c4](-0.5, -0.25) -- (0.8, -0.25);
\draw[<-, thick, color=c4](3.75, 0.25) -- (3.75, -1.25);
\end{tikzpicture}
\end{align*}

\paragraph{A bispecial string}
For the bispecial string $(+,w,-)$ with oriented word
\begin{align}  
\olra{w} =  \xf_+ - \xa \sim \ya - \overrightarrow{\yf_+ \sim \yf_-} - \yb \sim \xb - \overleftarrow{\xf_- 
\sim \xf_+} - \xb \sim \yb - \overleftarrow{\yf_- \sim \yf_+} - \ya \sim \xa - \xf_-
\end{align}
we obtain the matrix representation below.
\begin{align*}
\begin{tikzpicture}[baseline=0pt,yshift=0.25cm]
\node () at (-6,-0.25)	{$(M_1,M_2) = $};		
\node () at (-2.5,0) {		
$\begin{pNiceArray}{cc|cc}[first-row,first-col]
&\xf_+ & \xf_+ & \xf_- & \xf_- \\
\xa& 1 & 0 & 0 & 0  \\
\xa&0 &0  & 0 & 1  \\
\hline
\xb&0 & 1 & 1 & 0 \\
\xb&0 & 1 & 0 & 0 
\end{pNiceArray}$};
\node () at (2.5,0) {	
$\begin{pNiceArray}{cc|cc}[first-row,last-col]
\yf_+ & \yf_+ & \yf_- & \yf_- \\
1 & 0 & 0 & 0 & \ya \\
0 & 1 & 0 & 0 & \ya \\
\hline
1 & 0 & 1 & 0 & \yb \\
0 & 1 & 0 & 1 & \yb 
\end{pNiceArray}$};
\draw[-, dashed, color=c4](-0.75, 0.25) -- (0.75, 0.25);
\draw[-, dashed, color=c4](-0.75, -0.85) -- (0.75, -0.85);
\draw[->, thick, color=c4](-4.5, 0.6) -- (-4.5, -1.25);
\draw[<-, thick, color=c4](4.5, 0.6) -- (4.5, -1.25);
\end{tikzpicture}
\end{align*}

\subsection{Bondarenko's classification}
Some bands and strings yield isomorphic representations. This motivates the next notion.
For simplicity of presentation, throughout this subsection we make an additional assumption on the bunch of skew-gentle type $\fX$.
\begin{align}\label{eq:even}\tag{*}
\text{\emph{for any periodic word $w$ of $\fX$, 
the number of oriented arrows in $\olra{w}$ is even.}}
\end{align} 
This assumption is satisfied by the example $\fX_0$ as well as the later case of interest  $\fXA$.
\begin{dfn}\label{def:equiv}
\emph{Equivalence} on the set of bands and strings of a bunch of skew-gentle type $\fX$ satisfying \eqref{eq:even} is defined as the equivalence relation generated by the following.
\begin{enumerate}
\item Any band $(w,m,\lambda)$ is equivalent to the bands $(w_2,m,\frac{1}{\lambda})$ and $(w^{\op}, m,\lambda^{\pm})$,
where $\lambda^{\pm} = \lambda$ if $w$ is asymmetric and $\lambda^{\pm} = \frac{1}{\lambda}$ otherwise.
\item Any usual string $w$ is equivalent to the usual string $w^{\op}$. 
\item Any special string $(w,\varepsilon)$ is equivalent to the special string $(w^{\op}, \varepsilon)$.
\item Any bispecial string $(\varepsilon_1,w, \varepsilon_2)$ 
is equivalent to the bispecial string $(\varepsilon_2, (v^{\op})^m, \varepsilon_1)$,
where $v$ denotes the primitive root of the bispecial word $w$ and $m$ its multiplicity.
\end{enumerate}
\end{dfn}

The following statement follows from work by Bondarenko~\cite{B1, B2}.
\begin{thm}\label{thm:Bondarenko}
Let $\fX$ be any bunch of skew-gentle type satisfying~\eqref{eq:even} and $\kk$ an algebraically closed field.
Then the following statements hold.
\begin{enumerate}
\item For any string or band $\omega$ the corresponding representation $M(\omega)$ is a regular indecomposable $\kk$-linear  representation of $\fX$.
\item Any indecomposable regular $\kk$-linear representation $M$ is isomorphic to the representation $M(\omega)$ for a suitable string or band $\omega$ of $\fX$.
\item Any two strings or bands $\omega$ and $\omega'$ of $\fX$ are equivalent if and only if their corresponding representations $M(\omega)$ and $M(\omega')$ are isomorphic.
\end{enumerate}
\end{thm}
In brief terms, the theorem states that there is a bijection between the set of equivalence classes of strings and bands of $\fX$ and the set of isomorphism classes of regular indecomposable $\kk$-linear representations of $\fX$.

\begin{rmk} \label{rmk:sym-isom}
In the setup of Theorem~\ref{thm:Bondarenko},
the canonical representations of some bands and strings of $\fX$ admit the following isomorphic presentations.
\begin{enumerate}
\item For any band $\omega = (w,m,\lambda)$ the representation $M(\omega)$ is isomorphic to the representation with the Jordan block $\JB_\lambda$ in $M({\omega})$ replaced by any similar matrix.
\item Let $\omega$ be a symmetric band or a bispecial string
with multiplicity larger than one.
Let $M'(\omega)$ denote the representation
defined in the same way as $M(\omega)$ except that we use the opposite orientation $\overrightarrow{x_i \sim x_{i+1}}$ at each axis of symmetry in $w$.

\begin{enumerate}
\item If $\omega$ is a symmetric band $(w,m,\lambda)$, 
the representation $M'(\omega)$ 
is isomorphic to the representation $M(\omega')$
for the band $\omega' = (w,m,\frac{1}{\lambda})$.
\item If $\omega$ is bispecial string $(\varepsilon_1, w,\varepsilon_2)$
with multiplicity $m > 1$, 
the representation $M'(\omega)$
is isomorphic to $M(\omega')$
for the bispecial string 
$\omega' = \omega$ if $m$ is odd and $\omega' = (\ol{\varepsilon}_1, w,\ol{\varepsilon}_2)$ if $m$ is even.
\end{enumerate} 
\end{enumerate}
\end{rmk}

\begin{rmk}\label{rmk:field2}
Assume that $\kk$ is an arbitrary field. 
Then Definitions~\ref{def:band} and \ref{def:equiv} have to be adapted as follows.
\begin{itemize}
\item In Definition~\ref{def:band}, a band is a triple $(w,m,p)$ where $w$ is a cyclic word, $m\in \N$ and
$p \in \kk [x]$ is a monic irreducible polynomial such that $p(x)\neq x$ and also $p(x) \neq x + (-1)^{t(w)}$ if $w$ is symmetric. 
\item 
In Definition~\ref{def:equiv}, any band $(w,m,p)$ is equivalent to the bands $(w_2, m, \wt{p})$,
where 
$\wt{p}(x) = \frac{1}{p(0)} \, x^{\deg p}\,  p(\frac{1}{x}) $ denotes the reciprocal monic polynomial of $p$, as well as  $(w^{\op},m,p^\pm)$ with  
$p^\pm = p $  if $w$ is asymmetric
and $p^\pm = \wt{p}$ otherwise. 
\item In Subsection~\ref{subsec:can}, we replace the Jordan block $\JB_{\lambda}$ of size $m$
in the band representation $M(\omega)$ with the Frobenius block $\Phi_{p^m}$ with minimal polynomial $p^m$.
\end{itemize}
With these changes, Theorem~\ref{thm:Bondarenko} remains valid over the arbitrary field $\kk$.
\end{rmk}

We give a few more remarks for context, which are not essential later on.
\begin{enumerate}
\item	A \emph{bunch of gentle type} $\fX$ is a bunch of skew-gentle type 
together with the requirement that 
there are no two-point links.
A similar statement as Remark~\ref{rmk:red-gen} is true for gentle quivers.
In this case, the steps 
described in Subsections~\ref{subsub:signs} and \ref{subsec:A-orient}
are void and the construction of canonical forms $M(\omega)$ is considerably simpler.
\item The notion of a \emph{bunch of semichains} $\fX$ can be obtained from the notion of a bunch of skew-gentle type via the following modifications.
\begin{itemize}
\item Horizontal labels are allowed to be equivalent to vertical labels.
\item For any index $\imath$ each of the set $\fE_\imath$ or $\fF_\imath$ is allowed to be a subposet of the 
`infinite semi-chain' $(\{i_+,i_- \mid i \in \Z\},\leq)$  with $i_{\pm} < j_{\pm}$ and  $i_{\pm} <  j_\mp$ if $i < j$.
\end{itemize}
\item There are extensions of the notions of finite words, strings and their associated canonical forms which provide a complete classification of indecomposable, not necessarily regular representations of any bunch of semichains $\fX$. 
\end{enumerate}

\subsection{Indecomposable representations of the four subspace problem}
In this final subsection, we apply 
Theorem~\ref{thm:Bondarenko} to recover certain series of indecomposable  representations of the four subspace quiver $Q$, which is depicted below. 
$$
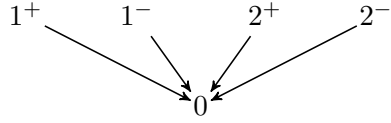
\begin{qrep}
1^+ \ar{rrrd}\&[-0.5cm] \&[-0.5cm] 1^- \ar{rd} \&[-0.5cm] \&[-0.5cm] 2^+ \ar{ld} \&[-0.5cm]\&[-0.5cm]  2^- \ar{llld} \\
\& \& 	\& 0
\end{qrep}
$$
The classification of indecomposable representations of $Q$ was obtained by Gelfand and Ponomarev \cite{GP2}
for an algebraically closed field and by Brenner~\cite{Brenner} for an arbitrary field, see also \cite{MZ} for an elementary approach.
As before, we consider $\kk$-linear representations over an algebraically closed field $\kk$.

Let
$\fX_Q = \fE \cup \fF = (\fE_1 \cup \fE_2) \cup (\fF_1 \cup \fF_2)$ with the singleton sets
$\fE_1 = \{
\xa \}$ 
and $\fE_2 = \{\ya \}$, and
two-point links
$\fF_1 = \{ \xf_+, \xf_- \}$ and $\fF_2 = \{ \yf_+, \yf_- \}$,
and the only non-reflexive equivalence relation given by  $\xa\approx \ya$.
Any representation of $\fX_Q$ is given by a pair of partitioned matrices
$$
\begin{array}{c}
\begin{tikzpicture}
\node () at (-3,0) {
$\begin{NiceArray}{|c|c|}[first-row,first-col,rules={width=1pt}]
& \xf_+ &   \xf_-  \\
\hline
\xa \ & \rowcolor{c2}A_+  &  A_-   \\ \hline
\end{NiceArray}$
};
\node () at (3,0) {	
$\begin{NiceArray}{|c|c|}[first-row,last-col,rules={width=1pt}]
\yf_+ &  \yf_-  \\
\hline
\rowcolor{c2} B_+ & B_- & \ \ya  \\ \hline
\end{NiceArray}$
};
\draw[-, dashed, color=c4](-1.75, -0.25) -- (1.75, -0.25);
\end{tikzpicture}
\end{array}
$$
and defines a representation $V$ of the quiver $Q$ via the matrices
$(A_+,A_-,B_+,B_-)$.
In fact, the category of $\kk$-linear
representations $\rep \fX_Q$ is isomorphic to the category $\rep Q$
of finite-dimensional $\kk$-linear representations of the four subspace quiver.

The next examples describe all indecomposable regular representations of $\fX_Q$ up to isomorphism.
For a  
representation $V$ of $Q$
we set $n_i = \dim V_i$ for each vertex $i \in Q_0$,
and denote the dimension vector of $V$ by $(n_0, n_{1^{+}}, n_{1^-}, n_{2^{+}}, n_{2^-})$.
\begin{ex}
There is only one family of bands  $\omega = (w,m,\lambda)$ up to equivalence, where $\dot{w} = 
(\xf_+ \sim \xf_- - \xa \sim \ya - \yf_+ \sim \yf_- - \ya \sim \xa - )$, $m \in \N$ and $\lambda \in \kk^*$ has to satisfy $\lambda \neq 1$ as $w$ is symmetric and has exactly two tied subwords.
The band representation $M(\omega)$ is
given by adapting the matrices in \eqref{eq:bd-sym} to multiplicity $m$, and yields:
\begin{align}
\label{eq:bigband}
\begin{tikzpicture}[baseline=0pt,yshift=0.25cm]
\node () at (-4,-0.25)	{$(M_1,M_2) = $};		
\node () at (-1.5,0) {		
$\begin{pNiceArray}{c|c}[first-row,first-col]
&\xf_+ & \xf_-  \\
\xa & 0 & \Id   \\
\xa & \Id & \JB_\lambda  
\end{pNiceArray}$};
\node () at (1.5,0) {	
$
\begin{pNiceArray}{c|c}[first-row,last-col]
\yf_+ & \yf_-& \\
\Id & \Id & \ya \\
0 & \Id & \ya
\end{pNiceArray}$
};
\draw[-, dashed, color=c4](-0.5, -0.25) -- (0.5, -0.25);
\end{tikzpicture}
\end{align}
The corresponding representation of $Q$ is a $\tau$-invariant indecomposable representation
$R^{(m)}_\lambda$  with dimension vector $(2m,m,m,m,m)$,
where $\tau$ denotes the Auslander-Reiten translation.
\end{ex}
\begin{ex}
Any bispecial string is equivalent to $(\varepsilon_1, v^n, \varepsilon_2)$ with $\varepsilon_1, \varepsilon_2 \in \{+,-\}$, $n \in \N$ with $v = 
\xf - \xa \sim \ya - \yf$.

\begin{itemize}
\item If $n=5$, the matrix representation $M(\omega)$
has the form:
\begin{align*}
(M_1,M_2) = 
\begin{pNiceArray}{c|c|c|c|c}[first-row,first-col]
&\xf_{\varepsilon_1} & \xf_{\ol{\varepsilon}_1} & \xf_{{\varepsilon_1}} & \xf_{\ol{\varepsilon}_1} & \xf_{{\varepsilon_1}}  \\
\xa & 1 & 0 & 0 & 0 & 0 \\
\xa & 0 & 1 & 1 & 0 & 0  \\
\xa & 0 & 0 & 1 & 0 & 0  \\
\xa & 0 & 0 & 0 & 1 & 1  \\
\xa & 0 & 0 & 0 & 0 & 1 
\end{pNiceArray}
&&
\begin{pNiceArray}{c|c|c|c|c}[first-row,last-col]
\yf_{\varepsilon_2} & \yf_{\ol{\varepsilon}_2} & \yf_{\varepsilon_2} & \yf_{\ol{\varepsilon}_2} & \yf_{\varepsilon_2} &
\\
1 & 1 & 0 & 0 & 0 & \ya \\
0 & 1 & 0 & 0 & 0 & \ya \\
0 & 0 & 1 & 1 & 0 & \ya \\
0 & 0 & 0 & 1 & 0 & \ya \\
0 & 0 & 0 & 0 & 1 & \ya \\
\end{pNiceArray}
\end{align*}
For any higher $n \in \N$ the matrices can be extended
by continuing the diagonal to have unit entries and the superdiagonal to form alternating binary sequences.
\item Assume that $n =2m$ with $m \in \N$ and $(\varepsilon_1,\varepsilon_2) = (-,+)$. 
Then $M(\omega)$ is isomorphic 
to the matrix representation~\eqref{eq:bigband} 
with multiplicity $m$ and the block $\JB_\lambda$ replaced by the upper triangular, nilpotent Jordan block $\JB_0$.
\item More generally, if $n=2m$ with $m \in \N$ and the signs $\varepsilon_1, \varepsilon_2$ are arbitrary, we obtain 
a $\tau^2$-invariant indecomposable representation $R^{(2m)}_{\varepsilon_1,\varepsilon_2}$
with dimension vector $(2m,m,m,m,m)$.
\item If $n =2m-1$ with $m \in \N$, then the quiver representation of $M(\omega)$ is isomorphic to the $\tau^2$-invariant indecomposable representation $R^{(2m-1)}_{\varepsilon_1,\varepsilon_2}$
with dimension vector $(2m, m - \delta_{\varepsilon_1,-}, m-\delta_{\varepsilon_1,+}, 
m - \delta_{\varepsilon_2,-}, m-\delta_{\varepsilon_2,+})$.
\end{itemize}
\end{ex}
There are no usual or special strings of $\fX_Q$ corresponding to regular representations.
Summarized, the regular representations of $\fX_Q$ yield all homogeneous tubes of rank one, which have indecomposables $(R^{(m)}_\lambda)_{m\in\N}$ with $\lambda \in \kk\backslash\{0,1\}$,
and two homogeneous tubes of rank two 
which have indecomposables
$(R^{(m)}_{\varepsilon_1, \varepsilon_2})_{m \in \N}$ with $\varepsilon_1, \varepsilon_2 \in \{+,-\}$.
However, the remaining homogeneous tube of rank two, which contains indecomposable modules with dimension vectors $(2,1,1,0,0)$ and $(2,0,0,1,1)$, is missing.
In Bondarenko's description of canonical forms for the whole category $\rep{\fX_Q}$, the indecomposables of this tube correspond to certain usual non-regular strings.
This indicates that the classification of $\rep{\fX_Q}$ in terms of bands and strings is not intrinsic: there exist $\tau$-periodic indecomposable objects in $\rep Q$ which are neither bispecial strings nor bands.

\section{Band and string complexes of the Gelfand order}\label{sec:complex}

In this section, we introduce 
sequences of numerical data called periodic and finite words, which are the main ingredients to define bands and strings of $\chfXA$.
Words are visualized in terms of \emph{gluing diagrams}.
For any band or string $\omega$ we construct a complex $\CP(\omega)$ in $\DbRep{\A}$.
The gluing diagram associated to $\omega$ is processed in three steps:  vertices are replaced by projective $\A$-modules, 
gluing edges are oriented via a technically involved symmetry rule in an ambient word, oriented edges then induce additional $\A$-linear maps in the gluing diagram.
The outcome is summarized in the compact form of a cycle or line diagram in Subsection~\ref{subsec:hor-not}.
Such a diagram yields directly
the complex $\CP(\omega)$.
We give examples of the construction of band and string complexes in Subsection~\ref{subsec:examples},
and formulate one of the main results of this part of the article in Theorem~\ref{thm:bijection2}.
Its proof is based on the categorical reductions of the preceding three sections
and
deferred until the next section. The latter also clarifies the connection between the notion of bands and strings of the bunch $\fXA$ in the sense of 
Section~\ref{app:mp} and the notions of the present, self-contained section.

\subsection{Bands and strings of indecomposable projective complexes} \label{subsec:strands}

As before, we assume that the base field $\kk$ is algebraically closed. Arbitrary fields are addressed in Remark~\ref{rmk:field}.
We are going to define 
\emph{periodic words} / \emph{bands} / \emph{finite words} / \emph{strings of $\chfXA$}.
For brevity, we will omit that these notions are meant to be associated to $\chfXA$, and use the more precise wording when it
will become relevant in the next section.

We introduce a set of formal symbols called \emph{oriented numbers}
\begin{align*}
	\kX =  \ora{\N} \cup \ola{\N}
	\quad	\text{ with } \quad \ora{\N} = \{ \ora{n} \mid n \in \N \}
	\quad	\text{ and } \quad
	\ola{\N} = \{ \ola{n} \mid n \in \N \}.
\end{align*}
For any oriented number $x$, we denote by $x^t$ the oriented number obtained by flipping the arrow in $x$.
We will use the notation $\olra{n}$ as a placeholder 
for $\ora{n}$ or $\ola{n}$, just as $\pm$ stands for $+$ or $-$.
\begin{dfn}\label{dfn:per-wo}
	Let $w = (x_i)_{i \in \Z}$ be a sequence of positive oriented numbers.
	\begin{enumerate}
		\item For any integer $j \in \Z$  the \emph{$j$-th shift of $w$} is given by 
		$w_j \colonequals (x_{i+j})_{i \in \Z}$.
		\item The sequence $w$ has period $\ell$ if $w=w_{\ell}$ and $\ell$ is the minimal number from $2\N$ 
		with this property,
		and its \emph{periodic part} $\dot{w} = (x_1,x_2,\ldots x_{\ell-1}, x_{\ell})$ has as many left-oriented numbers as right-oriented ones. 
		In this case, the sequence $w$ is called a \emph{periodic word}.
	\item The \emph{opposite word of $w$} is defined by 
	$
	w^{\op} \colonequals (x^{t}_{1-i})_{i \in \Z}$.
	\item The word $w$ is \emph{symmetric} if there is an integer $j \in \Z$ such that $w_j = (w_j)^{\op}$, that is, $x_{j-p} = x^{t}_{j+1+p}$ for any $p \geq 0$.
\end{enumerate}
\end{dfn}
For a periodic word $w$ we will depict its periodic part $\dot{w} = (\olra{n}_1, \olra{n}_2, \ldots, \olra{n}_{\ell-1}, \olra{n}_{\ell})$ 
by a  \emph{closed gluing diagram}
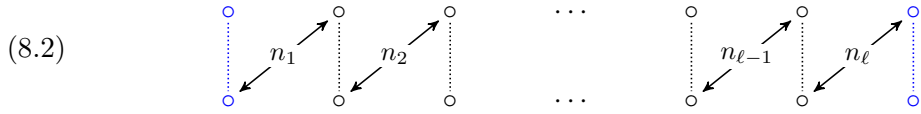
\begin{align} \label{eq:glu-cl}
\begin{gd}
	{\color{blue} \bt} \ar[color=blue,densely dotted,-]{d}\& \bt \ar[densely dotted,-]{d} \&\bt  \ar[densely dotted,-]{d} \&  \cdots   \&
	\bt \ar[densely dotted,-]{d}  \& \bt \ar[densely dotted,-]{d}  \&  {\color{blue} \bt} \ar[densely dotted,-,color=blue]{d}   \\
	{\color{blue} \bt}  \ar[<->]{ru}[description]{n_1} \& \bt  \ar[<->]{ru}[description]{n_2} \& \bt \&  \cdots  \& \bt \ar[<->]{ru}[description]{n_{\ell-1}} \& \bt \ar[<->]{ru}[description]{n_\ell} 	\& {\color{blue} \bt}
\end{gd}
\end{align}
in which the first and the last vertical edge should be identified, and for each index $1 \leq i \leq \ell$ the diagonal arrow decorated by $n_i$ is oriented towards the top right vertex if $\olra{n}_i$ is right-oriented,
respectively towards the bottom left vertex if $\olra{n}_i$ is left-oriented. 
The diagram for a periodic part of the opposite of $w$ can be obtained by a rotation by 180° of the diagram above.
\begin{dfn}\label{dfn:band}
A \emph{band} $\omega$ 
is given by data  $(d, w,m, \lambda)$ consisting of  an integer $d \in \Z$, a periodic word $w$, 
a number $m \in \N$, called the \emph{multiplicity} of the band,
and		a parameter $\lambda \in \kk^*$
satisfying 
$\lambda 
\neq -1$
if $w$ is symmetric.	
A band $\omega$ is called \emph{symmetric} if its underlying word $w$ is symmetric.
\end{dfn}
For the next notion we define $\kX_0  =\kX \cup \{ \ora{0}, \ola{0} \}$. Moreover, we set $(\ora{0})^t = \ola{0}$, and vice versa.
\begin{dfn}\label{dfn:ow-ab}
A \emph{finite word}
is given by  a sequence
$(\alpha, x_1, x_2, \ldots x_{\ell}, \beta)$ 
with ends $\alpha, \beta \in \{\star,\diamond\}$, oriented numbers $x_2,\ldots, x_{\ell-1} \in \kX$
and $x_1, x_{\ell} \in \kX_0$  such that 
the following conditions hold.
\begin{itemize}
\item
$x_i = \ola{0}$ only if $i=1$, $\alpha = \star$ and $w \neq (\star,x_1,\star)$. 
\item $x_i = \ora{0}$ only if $i = \ell$, $\beta = \star$ and $w \neq (\star,x_1,\star)$.
\end{itemize}  
Any word $w$ of the form above gives rise to the following notions.
\begin{enumerate}
\item The \emph{opposite word} of the word $w$ is defined as 
$
w^{\op} \colonequals (\beta,  x^t_{\ell},  \ldots x_2^t, x^t_1, \alpha)$.
\item The word $w$ is \emph{asymmetric} if $w \neq w^{\op}$.
\item  An end $\alpha$ or $\beta$ of the word $w$ is called \emph{special}
if it is given by $\diamond$.
\item The word $w$ is \emph{usual}, \emph{special} respectively \emph{bispecial} if it has zero, one respectively two special ends.	
\end{enumerate}
\end{dfn}
A finite word $(\alpha,\olra{n}_1,\olra{n}_2, \ldots, \olra{n}_{\ell}, \beta)$ will be depicted by a gluing diagram
\begin{align}
\label{eq:glu-op}
\begin{gd}
\& \bt \ar[densely dotted,-]{d} \&\bt  \ar[densely dotted,-]{d} \&  \cdots   \&
\bt \ar[densely dotted,-]{d}  \& \bt \ar[densely dotted,-]{d}  \&  \beta \\
\alpha  \ar[<->]{ru}[description]{n_1} \& \bt  \ar[<->]{ru}[description]{n_2} \& \bt \&  \cdots  \& \bt \ar[<->]{ru}[description]{n_{\ell-1}} \& \bt \ar[<->]{ru}[description]{n_\ell} 	\& 
\end{gd}
\end{align}
where the diagonal arrow with decoration $n_i$ is oriented in the same direction as its underlying oriented number $\olra{n}_i$.
Again, the opposite of the word can be visualized by a half-turn of the above diagram.

\begin{rmk}
The possibility in Definition~\ref{dfn:ow-ab} that a finite word $w$ may begin with $(\star,\ola{0})$ or end with $(\ora{0},\star)$
has to do with the fact that
a map 
$\begin{td} P_{\pm} \ar{r}{t^n} \& P_\star \end{td}$
is non-invertible for $n=0$ as it yields the natural inclusion $\begin{td}P_{\pm} \ar[hookrightarrow]{r} \& P_\star \end{td}$ with two-dimensional cokernel.
\end{rmk}

\begin{dfn}\label{dfn:strings}
There are three classes of
\emph{strings}
defined as follows.
\begin{enumerate}
\item A \emph{usual string} $(d,w)$ is given by an integer $d\in\Z$ and a usual asymmetric word $w$. 
\item  A \emph{special string} $(d,w, \varepsilon)$ is given by an integer $d\in \Z$ , any special word $w$, and a sign $\varepsilon\in \{+,-\}$.
\item A \emph{bispecial string} $(d,\varepsilon_1,w, \varepsilon_2)$
is given by an  integer $d \in \Z$,
any bispecial word $w$, and two signs $\varepsilon_1, \varepsilon_2 \in \{+,-\}$.		
\end{enumerate}
\end{dfn}

\subsection{Definition of band and string complexes}\label{sec:gluing}

Throughout this subsection, we fix a band or string $\omega$,
that is, $\omega$ is a 
a band $(d,w,m,\lambda)$, or a string of the form $(d,w)$, $(d,w,\varepsilon)$ or $(d,\varepsilon_1, w,\varepsilon_2)$.
The next goal is to describe the construction of a bounded complex $\CP(\omega)$
in the category $\DbRep{\A}$.
In terms of the gluing diagram of $w$, this construction is divided into the following steps.
\begin{enumerate}
\item \label{eq:polar} 
Each vertex in the gluing diagram is attributed a certain degree $d \in \Z$
and gives rise to a projective $\A$-module.
\item \label{eq:orient} Each edge in the gluing diagram obtains an orientation.
\item \label{eq:arrows} The oriented edges induce new arrows in the gluing diagram.
\end{enumerate}
Each step corresponds to one of the next subsections.
Step~\eqref{eq:polar} 
determines already the projective module $\AP_i(\omega)$ in $\CP(\omega)$ at each degree $i \in \Z$.
Step~\eqref{eq:orient} is technically involved and depends only on the sequence of oriented numbers in the word $w$. Steps~\eqref{eq:orient} and~\eqref{eq:arrows}
determine the differentials of the complex $\CP(\omega)$.

\subsubsection{Vertices and projectives}
\label{subsec:polar}
Before we can carry out the first step, we need a few notions concerning degrees and bispecial words.
\begin{dfn}\label{dfn:deg}
Let $(d,w)$ be a pair, where $d \in \Z$
and $w$ is	a periodic word with $\dot{w} = (x_1,x_2,\ldots, x_{\ell})$
or $w$ is a finite word $(\alpha,x_1,x_2,\ldots, x_{\ell},\beta)$.
\begin{enumerate}
\item The \emph{degree sequence} associated to $(d,w)$ 
is the sequence $(d_0, d_1, \ldots, d_{\ell})$ defined inductively by $d_0 =d$
and $d_{i} = d_{i-1} - 1$ if $x_i$ is right oriented, respectively $d_{i} = d_{i-1} + 1$ if $x_i$ is left oriented for each index $1 \leq i \leq \ell$.
\item	The \emph{opposite degree} of the pair $(d,w)$ is set to  $d^{\op} \colonequals d_{\ell}$.
\item  The \emph{support} of $(d,w)$ is the set of integers $
\supp(d,w) = 
\{d_0,d_1, \ldots, d_{\ell}\}$.
\end{enumerate}
\end{dfn}

\begin{dfn}
\label{dfn:root}
Below, we consider bispecial words
in the sense of Definition~\ref{dfn:ow-ab}.
\begin{enumerate}
\item For a bispecial word $v$ we denote by $\dot{v} =  (x_1, x_2, \ldots x_{\ell})$ its underlying sequence of oriented numbers and set 
$\dot{v}^{\op} \colonequals (x_{\ell}, \ldots x_2,x_1)$.
For any  number $m \in \N$ 
the \emph{$m$-th power of $\dot{v}$} is defined as the concatenation of sequences $\dot{v}^m = \dot{v} \circ \dot{v}^{\op} \circ \dot{v} \circ  \dot{v}^{\op} \ldots \circ \dot{v}_m$
which ends with $\dot{v}_m = \dot{v}$ if $m$ is odd, respectively $\dot{v}_m = \dot{v}^{\op}$ if $m$ is even.
\item 	The \emph{primitive root of a bispecial word $w$} is given by the shortest bispecial word $v$
such that 	$\dot{w} = \dot{v}^m$ for a number $m \in \N$.
In this case, the number $m$ is called the \emph{multiplicity} of the primitive root in $w$.
\end{enumerate}
\end{dfn}
In this step, vertices in the gluing diagram are replaced by
projective $\A$-modules.
\begin{itemize}
\item Assume that $\omega$ is a band.
Denoting the periodic part by $\dot{w} = (\olra{n}_1, \olra{n}_2,\ldots, \olra{n}_{\ell-1}, \olra{n}_{\ell})$
of its underlying word,
the diagram~\eqref{eq:glu-cl} is changed into the following.
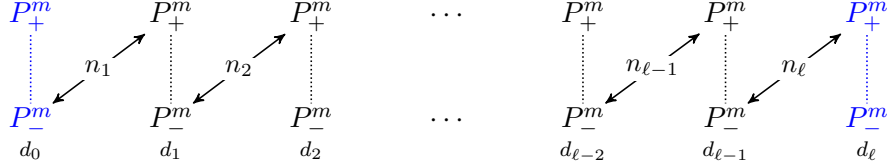
\begin{align}
\label{eq:glu-cl2}
\begin{gd}
	{{\color{blue} P_{+}^m}} \ar[color=blue,densely dotted,-]{d}\& P_{+}^m \ar[densely dotted,-]{d} \& P_{+}^m  \ar[densely dotted,-]{d} \&  \cdots   \&
	P_{+}^m \ar[densely dotted,-]{d}  \& P_{+}^m \ar[densely dotted,-]{d}  \&  {\color{blue} P_{+}^m} \ar[densely dotted,-,color=blue]{d}   \\
	\underset{d_0}{{\color{blue} P_{-}^m}}  \ar[<->]{ru}[description]{n_1} \& \underset{d_1}{P_{-}^m}  \ar[<->]{ru}[description]{n_2} \& \underset{d_2}{P_{-}^m} \&  \cdots  \& \underset{d_{\ell-2}}{P_{-}^m} \ar[<->]{ru}[description]{n_{\ell-1}} \& \underset{d_{\ell-1}}{P_{-}^m} \ar[<->]{ru}[description]{n_\ell} 	\& \underset{d_{\ell}}{{\color{blue} P_{-}^m}}
\end{gd}
\end{align}
Again, the left and the right column have to be identified, which is possible since $d_\ell = d_{0}$.
The sequence of degrees is given by Definition~\ref{dfn:deg}. 
\item Assume that $\omega$ is a usual string $(d,w)$. Denoting $w = (\star, \olra{n}_1, \olra{n}_2, \ldots, \olra{n}_{\ell},\star)$,
the diagram~\eqref{eq:glu-op} is changed into the following.
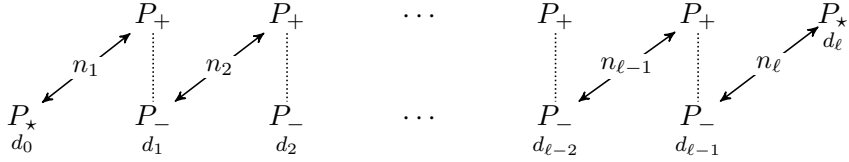
\begin{align*}
\begin{gd}
	\& P_{+} \ar[densely dotted,-]{d} \& P_{+}  \ar[densely dotted,-]{d} \&  \cdots  
	\&
	P_{+} \ar[densely dotted,-]{d}  \& P_{+} \ar[densely dotted,-]{d}  \&  \smash[b]{\underset{d_{\ell}}{P_{\star}}} \\
	\underset{d_0}{P_{\star}}  \ar[<->]{ru}[description]{n_1} \& \underset{d_1}{P_{-}}  \ar[<->]{ru}[description]{n_2} \& \underset{d_2}{P_{-}} \&  \cdots  \& \underset{d_{\ell-2}}{P_{-}} \ar[<->]{ru}[description]{n_{\ell-1}} \& \underset{d_{\ell-1}}{P_{-}} \ar[<->]{ru}[description]{n_\ell} 	\& 
\end{gd}
\end{align*}
Again, the degree sequence is determined by Definition~\ref{dfn:deg}.
\item If $\omega$ is  a special string $(d,w,\varepsilon)$
the diagram is given by the previous diagram except that $P_\varepsilon$ replaces $P_\star$ at degree $d_0$ if $\alpha = \diamond$, respectively $P_\star$ at degree $d_{\ell}$ if $\beta = \diamond$.
\item Assume that $\omega$ is bispecial $(d,\varepsilon_1, w,\varepsilon_2)$	 with $w = (x_1,x_2,\ldots, x_{\ell})$. 
Let $v$ be the primitive root of $w$ and $m$ its multiplicity.
Each copy of $v$ or $v^{\op}$ in $w$ corresponds to a full subdiagram of the gluing diagram of $w$.
The two types of copies are considered separately.
\begin{itemize}
\item
In the subdiagram corresponding to any copy of $v$, we replace its left bottom vertex by 
$P_{\varepsilon_1}$ and its right top vertex by $P_{\varepsilon_2}$,
and use the same convention for intermediate vertices as for 
usual or special strings, that is,
any other vertex in the top row by $P_+$ and any remaining vertex in the bottom row by $P_-$.
\item In the subdiagram of each copy of $v^{\op}$ we follow precisely the opposite convention;
the leftmost vertex becomes $P_{\ol{\varepsilon}_1}$, any gluing edge connects $P_-$ at the top to $P_+$ at the bottom, and the rightmost vertex is changed to $P_{\ol{\varepsilon}_2}$, 
where $\ol{\varepsilon}_1$ and  $\ol{\varepsilon}_2$
denote the opposite signs of $\varepsilon_1$ and $\varepsilon_2$, respectively.
\end{itemize}
As in the case of usual and special strings, we place the degrees $d_0,d_1, \ldots, d_{\ell}$ below the columns of the gluing diagram of $w$. If the primitive root has length $k$, we obtain a diagram which begins with
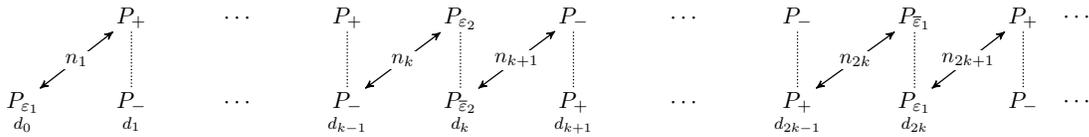
\begin{align*}
\scalebox{0.8}{$
	\begin{gd}
		\& P_{+} \ar[densely dotted,-]{d} \&  \cdots  
		\&
		P_{+} \ar[densely dotted,-]{d}  \&  {P_{\varepsilon_2}} \ar[densely dotted,-]{d}
		\& P_{-} \ar[densely dotted,-]{d} \&  \cdots  
		\&
		P_{-} \ar[densely dotted,-]{d} \&  
		{P_{\ol{\varepsilon}_1}}  
		\ar[densely dotted,-]{d} 
		\& P_+  \ar[densely dotted,-]{d} \ar[phantom]{r}[description]{\cdots}	\&  
		\mathstrut 
			\\
			\underset{d_0}{P_{\varepsilon_1}}  \ar[<->]{ru}[description]{n_1} \& \underset{d_1}{P_{-}} \&  \cdots  \&  \underset{d_{k-1}}{P_{-}} \ar[<->]{ru}[description]{n_k} 	\& \underset{d_k}{P_{\ol{\varepsilon}_2}} \ar[<->]{ru}[description]{n_{k+1}} \&
			\underset{d_{k+1}}{P_{+}} \&  \cdots  \&  \underset{d_{2k-1}}{P_{+}} 
			\ar[<->]{ru}[description]{		n_{2k}	} 	
			\& 
			\underset{d_{2k}}{P_{\varepsilon_1}} \ar[<->]{ru}[description]{n_{2k+1}} \& P_- \ar[phantom]{r}[description]{\cdots} \& \mathstrut
		\end{gd}$}
\end{align*}
and ends with the projective $P_{\varepsilon_2}$ if $m$ is odd, respectively $P_{\ol{\varepsilon}_1}$ if $m$ is even.
\end{itemize}		
In any of the diagrams above, 
both projectives connected by a vertical edge are attributed the degree below the bottom projective.
For each degree $i \in \Z$, 
the projective $\A$-module $\AP_i(\omega)$ of the complex $\CP(\omega)$ is defined to be
the direct sum of all projectives with degree $i$.
The gluing diagram obtained at this stage will be called \emph{polarized gluing diagram of $\omega$}.
\subsubsection{Orienting edges} \label{subsec:orient}
In the next step, each vertical edge in the gluing diagram of $w$ is given an orientation.
For this,  we need first to carry out certain preparations.
\begin{dfn}
Below, we consider periodic and finite words as in Definitions~\ref{dfn:per-wo} and \ref{dfn:ow-ab}.
\begin{enumerate}
	\item For a periodic word $w$, its ambient word $\ol{w}$ is $w$ itself.
	\item Assume that $w$ is a finite word $(\alpha, x_1, x_2, \ldots, x_{\ell-1}, x_{\ell}, \beta)$
	with $\ell > 1$. Its ambient word $\ol{w}$ is defined as follows.
	\begin{enumerate}
		\item If $w$ is bispecial, the word $\ol{w}$ is the symmetric periodic word $(y_i)_{i \in \Z}$
		such that
		$ (y_{1}, y_2,  \ldots y_{2 \ell})
		= 
		(x_1,x_2, \ldots x_{\ell}, x_\ell,  \ldots x_2,  x_{1})$
		and 					$y_{i+2\ell} = y_i$ for any integer $i \in \Z$.
		\item If $w$ is usual,  we set $\ol{w} = (\star,y_1, y_2, \ldots, y_{\ell},\star) = (\star, \tilde{x}_1,x_2,\ldots, x_{\ell-1}, \tilde{x}_{\ell},\star)$
		where 	
		$\tilde{x}_1$ and $\tilde{x}_{\ell}$ are obtained from $x_1$ and $x_{\ell}$  by replacing the underlying numbers $n_1$ respectively $n_{\ell}$
		by the positive half-integers
		\begin{align*}
			&&
			\tilde{n}_1 = \begin{cases} 
				n_1 +\frac{1}{2} & \text{if } (\alpha,x_1)=(\star,{\ola{n}_1}),
				\\
				n_1 - \frac{1}{2} & \text{if }  (\alpha,x_1)=(\star,{\ora{n}_1}),
			\end{cases}
			&& 
			\tilde{n}_\ell = \begin{cases} 
				n_\ell -\frac{1}{2} & \text{if }
				(x_{\ell},\beta) = ({\ola{n}_\ell},\star),
				\\
				n_\ell + \frac{1}{2} & \text{if }(x_{\ell},\beta) = ({\ora{n}_\ell},\star).
			\end{cases}
		\end{align*}
		\item If only the right end of $w$ is special, set
		\begin{align*} &&\ol{w} = 	(\star,y_1,y_2, \ldots  y_{2\ell},\star)  = (\star, \wt{x}_1,x_2 \ldots x_{\ell-1}, x_{\ell},  
			{\color{blue!75}
				x_{\ell},x_{\ell-1}, \ldots x_2, \wt{x}_1, \star}).\end{align*}
		\item If only the left end of $w$ is special, set 
		\begin{align*}
			&&\ol{w} = (\star,y_{-\ell-1},y_{-\ell-2}, \ldots  y_{\ell},\star) = 
			({\color{blue!75} \star, \wt{x}_\ell,x_{\ell-1}, \ldots x_2, x_{1},} \,
			x_1,x_2,\ldots x_{\ell-1}, \wt{x}_{\ell}, \star).\end{align*}
	\end{enumerate}
\end{enumerate}
\end{dfn}
In any case, the ambient word $\ol{w}$ is symmetric.
\paragraph{Orientation rule in combinatorial terms}

For each index $1 \leq i < \ell$ in case $w$ is finite respectively for
each index $i \in \Z$ if $w$ is periodic,
we define a new symbol $\kappa_i \in \{ \uparrow, \downarrow \}$,
which will be placed between $x_i$ and $x_{i+1}$ in $w$.
The orientation of the arrow $\kappa_i$ is determined by a subsequence of oriented numbers in the ambient word $\ol{w}$ as follows.
\begin{enumerate}
\item If $\ol{w}$ is periodic with 
$y_{i-k} = y^{t}_{i+k+1}$ for any $k \in \N_0$, set $\kappa_i$ to be $\downarrow$.
\item Otherwise, 
there exists a minimal number
$k \in \N_0$ with $y_{i-k} \neq y^{t}_{i+k +1}$.
Let ${p},{q} \in \{ \frac{n}{2} \mid n \in \N\}$ denote the numbers such that
$y_{i-k} = \olra{p}$ and $y_{i+k+1} = \olra{q}$. 
In particular, we have the situation
\begin{align*}
	\overline{w}= 
	{\color{blue!75}\ldots
		\underbrace{y_{i-k}}_{\olra{p}}, 
		\	
	}
	{
		\underbrace{ y_{i-k+1},  \ldots  
			y_{i}, y_{i+1},  
			\ldots  y_{i+k}}_{\text{maximal symmetric}}}, \
	{\color{blue!75}
		\underbrace{y_{i+k+1}}_{\olra{q}} \ldots}
\end{align*}
where the maximal symmetric subsequence is empty if and only if $y_i \neq y_{i+1}$.
With these notations we set
\begin{align}
	\label{eq:updownarrow}
	\kappa_i = \begin{cases}
		\downarrow &
		\text{if }
		(\olra{p}, \olra{q}) = (\ola{p}, \ola{q})\text{ or }
		(\ora{p},\ola{q})\text{ with } {p} > {q} \text{ or }
		(\ola{p},\ora{q})\text{ with } {p} < {q} ,\\
		\uparrow & 
		\text{if }
		(\olra{p}, \olra{q}) = (\ora{p}, \ora{q})\text{ or }
		(\ora{p},\ola{q})\text{ with } {p} < {q} \text{ or }
		(\ola{p},\ora{q})\text{ with } {p} > {q}.
	\end{cases}
\end{align}
\end{enumerate}
We record the outcome in a new combinatorial datum $w^{\updownarrow}$ as follows.
\begin{itemize}
\item In case the word $w$ is finite, 
we apply the rule for $\kappa_i$ to each index $1 \leq i \leq \ell$ and obtain a  sequence of the form $	w^{\updownarrow} = 
(\alpha, x_1
\kappa_1 
x_2  \kappa_2 \ldots 
x_\ell, \beta)$.
\item 
If $w$ is periodic, it is sufficient to apply the rule above to each index $1 \leq i \leq \ell$  and to consider 
$w^{\updownarrow} = (x_i \kappa_i)_{1 \leq i \leq \ell}$ identifying $\kappa_\ell$ with $\kappa_0$ in what follows.
\end{itemize}
In the corresponding gluing diagram, the arrow $\kappa_i \in \{\uparrow, \downarrow \}$ corresponds to orienting the $i$-th gluing edge in the same direction as $\kappa_i$.

\paragraph{Orientation mnemonic in terms of gluing diagrams}
The orientation rule \eqref{eq:updownarrow} for the gluing edge $\kappa_i$ at degree $d_i$ has the following diagrammatic interpretation.
\begin{enumerate}
\item The word $\ol{w}$ admits a gluing diagram of the form \eqref{eq:glu-cl} or \eqref{eq:glu-op}
but possibly with half integers as decorations.
\item Starting from the gluing edge $\kappa_i$ we read the diagonal arrows one step to the left, and one step to the right.
If these arrows have the same numerical values and opposite orientations we contract the subdiagram 
containing the two squares containing these to a single gluing edge denoted again by $\kappa_i$.
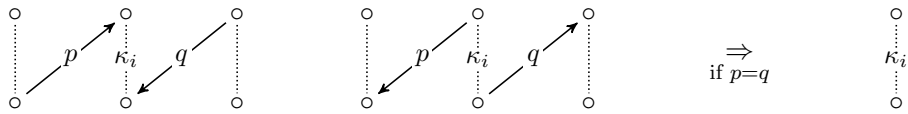
\begin{align*}
\begin{gd}
	\bt \ar[densely dotted,-]{d} \& \bt  \ar[densely dotted,-]{d}[description]{\kappa_i} \& \bt \ar[densely dotted,-]{d} \\
	\bt \ar[->]{ru}[description]{p}  \& \bt \ar[<-]{ru}[description]{q} \& \bt
\end{gd}
&&
\begin{gd}
	\bt \ar[densely dotted,-]{d} \& \bt  \ar[densely dotted,-]{d}[description]{\kappa_i} \& \bt \ar[densely dotted,-]{d} \\
	\bt \ar[<-]{ru}[description]{p}  \& \bt \ar[->]{ru}[description]{q} \& \bt
\end{gd}
&& 
\underset{\text{if }p = q}
{\Rightarrow} &&
\begin{gd}
	\bt \ar[densely dotted,-]{d}[description]{\kappa_i} \\
	\bt
\end{gd}
\end{align*}
\item 
The previous step is iterated until one of the following occurs.
\begin{itemize}
\item If the gluing diagram of $\ol{w}$ is periodic and we end up with a single gluing edge by iterated contractions, the edge $\kappa_i$ in the original diagram of $w$ is oriented downwards.
\item Otherwise we obtain a pair of arrows incident to gluing edge 
$\kappa_i$
which have different values $p$ or $q$ or identical orientations.
\begin{itemize}
	\item If the arrows have the same orientations, we orient $\kappa_i$ in such a way that it creates a source or a sink in the diagram.
	\begin{align*}
		{
			\setlength{\extrarowheight}{3pt}
			\begin{array}{|r||c|c|}
				\hline
				\text{case} &		(\ola{p},\ola{q}) & 		(\ora{p},\ora{q}) \\ \hline
				& & \\[-10pt]
				\text{diagram} &
				\begin{gd}
					\bt \ar[densely dotted,-]{d} \& \bt  \ar[densely dotted,->]{d}[description]{\kappa_i} \& \bt \ar[densely dotted,-]{d} \\
					\bt \ar[<-]{ru}[description]{p}  \& \bt \ar[<-]{ru}[description]{q} \& \bt
				\end{gd}
				&
				\begin{gd}
					\bt \ar[densely dotted,-]{d} \& \bt  \ar[densely dotted,<-]{d}[description]{\kappa_i} \& \bt \ar[densely dotted,-]{d} \\
					\bt \ar[->]{ru}[description]{p}  \& \bt \ar[->]{ru}[description]{q} \& \bt
				\end{gd}
				\\[20pt]
				\hline
			\end{array}
		}
	\end{align*}
	A \emph{transition vertex} is a vertex which is neither a sink nor a source with respect to any kind, diagonal or vertical, arrows.
	Note that in both cases a flip of the orientation would create a transition vertex incident to the oriented edge.
	\item If the arrows have different orientations, it holds that $p \neq q$. We orient $\kappa_i$ in such a way that it creates a sink or a source incident to the arrow with the smaller value.
	\begin{align*}
		{
			\setlength{\extrarowheight}{3pt}
			\begin{array}{|r||c|c|c|c|}
				\hline
				\text{case} &
				(\ora{p},\ola{q}), p > q
				&
				(\ola{p},\ora{q}), p < q  \\
				\hline
				& & \\[-10pt]
				\text{diagram} &
				\begin{gd}
					\bt \ar[densely dotted,-]{d} \& \bt  \ar[densely dotted,->]{d}[description]{\kappa_i} \& \bt \ar[densely dotted,-]{d} \\
					\bt \ar{ru}[description]{p}  \& \bt \ar[<-]{ru}[description]{q} \& \bt
				\end{gd}
				&
				\begin{gd}
					\bt \ar[densely dotted,-]{d} \& \bt  \ar[densely dotted,->]{d}[description]{\kappa_i} \& \bt \ar[densely dotted,-]{d} \\
					\bt \ar[<-]{ru}[description]{p}  \& \bt \ar[->]{ru}[description]{q} \& \bt
				\end{gd}
				\\[20pt]
				\hline
				\hline
				\text{case} & 		(\ora{p},\ola{q}), p < q&
				(\ola{p},\ora{q}), p > q \\
				\hline
				& & \\[-10pt]
				\text{diagram} &
				\begin{gd}
					\bt \ar[densely dotted,-]{d} \& \bt  \ar[densely dotted,<-]{d}[description]{\kappa_i} \& \bt \ar[densely dotted,-]{d} \\
					\bt \ar{ru}[description]{p}  \& \bt \ar[<-]{ru}[description]{q} \& \bt
				\end{gd}
				&
				\begin{gd}
					\bt \ar[densely dotted,-]{d} \& \bt  \ar[densely dotted,<-]{d}[description]{\kappa_i} \& \bt \ar[densely dotted,-]{d} \\
					\bt \ar[<-]{ru}[description]{p}  \& \bt \ar[->]{ru}[description]{q} \& \bt
				\end{gd}
				\\[20pt]												\hline
		\end{array}}
	\end{align*}		
	In each case, the flip of the orientation would destroy the desired sink or source, and create a transition vertex. Therefore the orientation is unique.
	In each case, there is a source and sink incident to the diagonal arrow with the smaller value.
\end{itemize}
The edge $\kappa_i$ in the gluing diagram of $w$ is oriented in the same direction as the processed edge above.
\end{itemize}
\end{enumerate}
This diagrammatic procedure recovers the rules~\eqref{eq:updownarrow}.
A brief  mnemonic for the diagrammatic orientation rules 
is that the orientation of $\kappa_i$ \emph{avoids transitions through smaller values}.	
\subsubsection{Induced arrows and homomorphisms}
\label{subsec:ind-arr}
So far, we have replaced the vertices by projectives and oriented all gluing edges in the gluing diagram of the band or string $\omega$.

The next rules provide new arrows depending on the particular orientations of the gluing edges.
As these rules do not depend on the isomorphism class of the projective at each vertex, any projective $P_\star$, $P_+$ or $P_-$ will be denoted simply by $P$ in the diagrams below.

\paragraph{Induced arrows in string diagrams}
Assume that $\omega$ is a string.
For any degree $d_i$ we set $\sign_i = (-1)^{d_i}$ and denote by $\ol{\sign}_i$ the opposite sign.
\begin{enumerate}
	\item \emph{Source rules}. 
	Whenever an oriented vertical edge and a diagonal arrow start in a common source, we add a horizontal arrow.
	The diagonal and the induced arrow are changed into multiplication maps with monomials $\pm t^{n_i}$, where the power $n_i$ 
	is determined by the decoration of the diagonal arrow.
	
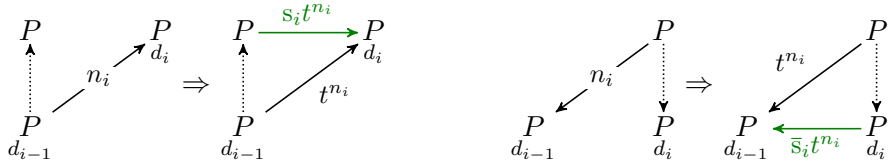
\begin{align}
		\label{eq:open-ind}
		\begin{gd}
			P  \ar[densely dotted,<-]{d} \& \smash[b]{\underset{d_i}{P}} \\
			\underset{d_{i-1}}{P} \ar{ru}[description]{n_i}
		\end{gd}
		\Rightarrow 
		\begin{gd}
			P  \ar[densely dotted,<-]{d} \ar[color=green]{r}{ \sign_i t^{n_i}} \& \smash[b]{\underset{d_i}{P}} \\
			\underset{d_{i-1}}{P} \ar{ru}[swap]{t^{n_i}}
		\end{gd}
		&& 
		\begin{gd}
			\& P \ar[densely dotted]{d} \\
			\underset{d_{i-1}}{P} \ar[<-]{ru}[description]{n_i} \& \underset{d_{i}}{P}
		\end{gd}
		\Rightarrow 
		\begin{gd}
			\& P \ar[densely dotted]{d} \\
			\underset{d_{i-1}}{P} \ar[<-]{ru}{t^{n_i}} \ar[<-, color=green]{r}[swap]{\ol{\sign}_i t^{n_i}} \& \underset{d_{i}}{P}
		\end{gd}			
	\end{align}
	\item \emph{Sink rules}.
	Similarly, there is an additional arrow whenever an oriented vertical edge and an arrow have a common sink.
	
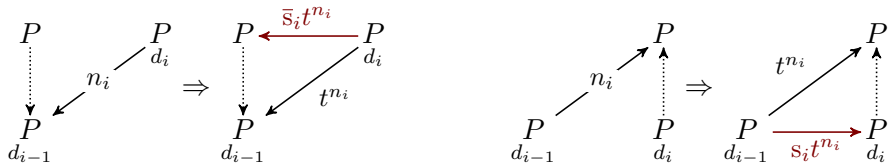
\begin{align}
		\begin{gd}
			P  \ar[densely dotted]{d} \& \smash[b]{\underset{d_i}{P}} \\
			\underset{d_{i-1}}{P} \ar[<-]{ru}[description]{n_i}
		\end{gd}
		\Rightarrow 
		\begin{gd}
			P  \ar[densely dotted]{d} \ar[color=red,<-]{r}{\ol{\sign}_i t^{n_i}} \& \smash[b]{\underset{d_i}{P}} \\
			\underset{d_{i-1}}{P} \ar[<-]{ru}[swap]{t^{n_i}}
		\end{gd}
		&&
		\begin{gd}
			\& P \ar[densely dotted,<-]{d} \\
			\underset{d_{i-1}}{P} \ar{ru}[description]{n_i} \& \underset{d_{i}}{P}
		\end{gd}
		\Rightarrow 
		\begin{gd}
			\& P \ar[densely dotted,<-]{d} \\
			\underset{d_{i-1}}{P} \ar{ru}{t^{n_i}} \ar[color=red]{r}[swap]{\sign_i t^{n_i}} \& \underset{d_{i}}{P}
		\end{gd}			
	\end{align}
	\item \emph{Iterated application}. The sink and source rules have to be applied once again taking the induced arrows into account.
	
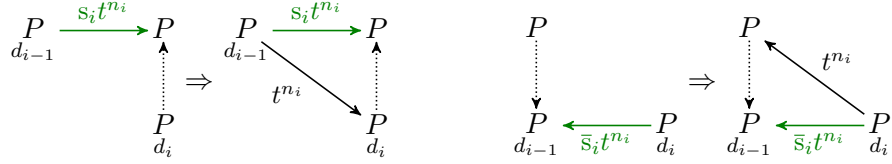
\begin{align}
		\begin{gd}
			\smash[b]{\underset{d_{i-1}}{P}}    \ar[color=green]{r}{\sign_i t^{n_i}}  \& P \ar[densely dotted,<-]{d}  \\
			\& \underset{d_{i}}{P} 
		\end{gd}
		\Rightarrow 
		\begin{gd}
			\smash[b]{\underset{d_{i-1}}{P}}    \ar[color=green]{r}{\sign_i t^{n_i}} \ar{rd}[swap]{t^{n_i}} \& P \ar[densely dotted,<-]{d}  \\
			\& \underset{d_{i}}{P} 
		\end{gd}
		&&
		\begin{gd}
			P \ar[densely dotted]{d}  \& \\
			\underset{d_{i-1}}{P}  \ar[color=green,swap,<-]{r}{\ol{\sign}_i t^{n_i}} \& \underset{d_{i}}{P}
		\end{gd}
		\Rightarrow 
		\begin{gd}
			P \ar[densely dotted]{d}  \ar[<-]{rd}{t^{n_i}} \& \\
			\underset{d_{i-1}}{P}  \ar[color=green,swap,<-]{r}{\ol{\sign}_i t^{n_i}} \& \underset{d_{i}}{P}
	\end{gd}		\end{align}
	Any arrow induced in this step has a positive sign. 
\end{enumerate}

\paragraph{Induced arrows in band diagrams}	
Assume that $\omega$ is a band $(d,w,m,\lambda)$.
\begin{enumerate}
	\item  First, the parameter $\lambda$ is placed below the vertex in the bottom row incident to the arrow decorated by $n_1$.
	The induced arrows involving the first gluing edge involve a Jordan block $\JB_{\lambda}$ or $\JB_{\frac{1}{\lambda}}$ of size $m$
	with eigenvalue $\lambda$ respectively $\frac{1}{\lambda}$.
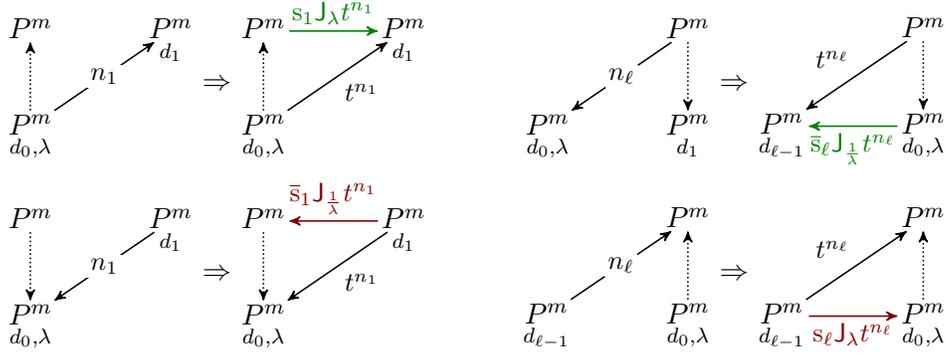
	\begin{align}
		\begin{gd}
			P^m  \ar[densely dotted,<-]{d} \& \smash[b]{\underset{d_1}{P^m}} \\
			\underset{d_0, \lambda}{P^m} \ar{ru}[description]{n_1}
		\end{gd}
		\Rightarrow 
		\begin{gd}
			P^m  \ar[densely dotted,<-]{d} \ar[color=green]{r}{ \sign_1  \JB_{\lambda}  t^{n_1} } \& \smash[b]{\underset{d_1}{P^m}} \\
			\underset{d_0, \lambda}{P^m} \ar{ru}[swap]{t^{n_1}}
		\end{gd}
		&& 
		\begin{gd}
			\& P^m \ar[densely dotted]{d} \\
			\underset{d_0, \lambda}{P^m} \ar[<-]{ru}[description]{n_\ell} \& \underset{d_1}{P^m}
		\end{gd}
		\Rightarrow 
		\begin{gd}
			\& P^m \ar[densely dotted]{d} \\
			\underset{d_{\ell-1}}{P^m} \ar[<-]{ru}{t^{n_\ell}} \ar[<-, color=green]{r}[swap]{\ol{\sign}_\ell  \JB_{\frac{1}{\lambda}} t^{n_\ell} } \& \underset{d_0,\lambda}{P^m}
		\end{gd}		
		\\	
		\begin{gd}
			P^m  \ar[densely dotted]{d} \& \smash[b]{\underset{d_1}{P^m}} \\
			\underset{d_0, \lambda}{P^m} \ar[<-]{ru}[description]{n_1}
		\end{gd}
		\Rightarrow 
		\begin{gd}
			P^m  \ar[densely dotted]{d} \ar[color=red,<-]{r}{\ol{\sign}_1   \JB_{\frac{1}{\lambda}} t^{n_1} } \& \smash[b]{\underset{d_1}{P^m}} \\
			\underset{d_0, \lambda}{P^m} \ar[<-]{ru}[swap]{t^{n_1}}
		\end{gd}
		&&
		\begin{gd}
			\& P^m \ar[densely dotted,<-]{d} \\
			\underset{d_{\ell-1}}{P^m} \ar{ru}[description]{n_\ell} \& \underset{d_0,\lambda}{P^m}
		\end{gd}
		\Rightarrow 
		\begin{gd}
			\& P^m \ar[densely dotted,<-]{d} \\
			\underset{d_{\ell-1}}{P^m} \ar{ru}{t^{n_\ell}} \ar[color=red]{r}[swap]{{\sign}_\ell  \JB_{\lambda}  t^{n_\ell} } \& \underset{d_0,\lambda}{P^m}
		\end{gd}	
	\end{align}
	\item We apply the same rules as in the case of an open gluing diagram.
	The induced maps with monomials $\pm t^{n_i}$ are actually diagonal matrices of size $m$ with $\pm t^{n_i}$ on the diagonal.
	Moreover, an iterated application of source and sink rules may give rise to induced maps with Jordan blocks.
	
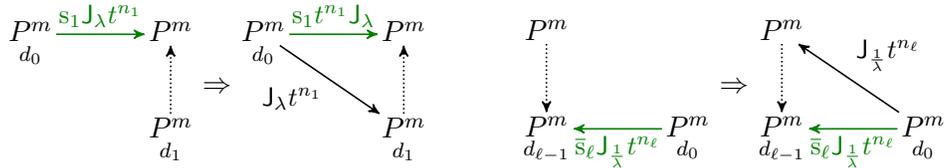
\begin{align}
		\begin{gd}
			\smash[b]{\underset{d_{0}}{P^m}}    \ar[color=green]{r}{\sign_1  \JB_{\lambda}  t^{n_1}}  \& P^m \ar[densely dotted,<-]{d}  \\
			\& \underset{d_{1}}{P^m} 
		\end{gd}
		\Rightarrow 
		\begin{gd}
			\smash[b]{\underset{d_{0}}{P^m}}    \ar[color=green]{r}{\sign_1   t^{n_1} \JB_{\lambda}} \ar{rd}[swap]{  \JB_{\lambda} t^{n_1}} \& P^m \ar[densely dotted,<-]{d}  \\
			\& \underset{d_{1}}{P^m} 
		\end{gd}
		&&
		\begin{gd}
			P^m \ar[densely dotted]{d}  \& \\
			\underset{d_{\ell-1}}{P^m}  \ar[color=green,swap,<-]{r}{\ol{\sign}_\ell   \JB_{\frac{1}{\lambda}} t^{n_\ell}  } \& \underset{d_{0}}{P^m}
		\end{gd}
		\Rightarrow 
		\begin{gd}
			P^m \ar[densely dotted]{d}  \ar[<-]{rd}{\JB_{\frac{1}{\lambda}} t^{n_\ell}  } \& \\
			\underset{d_{\ell-1}}{P^m}  \ar[color=green,swap,<-]{r}{\ol{\sign}_\ell \JB_{\frac{1}{\lambda}} t^{n_\ell} } \& \underset{d_{0}}{P^m}
	\end{gd}		\end{align}
\end{enumerate}

\begin{rmk}
	By the last step, any two consecutive gluing edges with the same orientation give rise to three induced arrows.
	
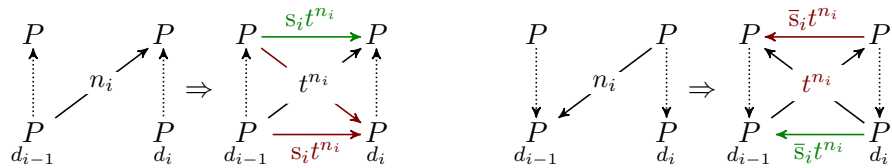
\begin{align}
		\begin{gd}
			P  \ar[densely dotted,<-]{d}     \& P \ar[densely dotted,<-]{d}  \\
			{\underset{d_{i-1}}{P}} \ar{ru}[description]{n_i} \& \underset{d_{i}}{P} 
		\end{gd}
		\Rightarrow 
		\begin{gd}
			P  \ar[densely dotted,<-]{d}   		   \ar[color=green]{r}{\sign_i t^{n_i}} \ar[color=red]{rd} \& P \ar[densely dotted,<-]{d}  \\
			{\underset{d_{i-1}}{P}} \ar{ru}[description, crossing over, inner sep=0pt, inner sep=1pt]{t^{n_i}}	\ar[color=red]{r}[swap]{\sign_i t^{n_i}}	\& \underset{d_{i}}{P} 
		\end{gd}
		&&
		\begin{gd}
			P  \ar[densely dotted]{d}     \& P \ar[densely dotted]{d}  \\
			{\underset{d_{i-1}}{P}} \ar[<-]{ru}[description]{n_i} \& \underset{d_{i}}{P} 
		\end{gd}
		\Rightarrow 
		\begin{gd}
			P  \ar[densely dotted]{d}   		   \ar[color=red,<-]{r}{\ol{\sign}_i t^{n_i}} \ar[<-]{rd} \& P \ar[densely dotted]{d}  \\
			{\underset{d_{i-1}}{P}} \ar{ru}[color=red,description, crossing over, inner sep=0pt,<-]{t^{n_i}}	\ar[color=green,<-]{r}[swap]{\ol{\sign}_i t^{n_i}}	\& \underset{d_{i}}{P} 
		\end{gd}
	\end{align}
	In any case, all arrows between two consecutive gluing edges are oriented in the same direction.
	In particular, each square gives rise to a map $\begin{td} P^2 \ar{r}{M_i t^{n_i}} \& P^{2} \end{td}$
	or $\begin{td} P^2 \ar[<-]{r}{M_i t^{n_i}} \& P^{2} \end{td}$.
\end{rmk}

\paragraph{Folding}
At last, we `fold' the gluing diagram in the following sense.
\begin{enumerate}
	\item All vertical gluing edges are ignored. 
	\item For any degree $i \in \Z$ we define a projective 
	$\AP_i$ as the direct sum of all projectives which were attributed degree $i$ in the diagram.
	\item  For each degree $i \in \Z$ the decorations of the arrows in the gluing diagram 
	starting at projectives of degree $i$ 
	yield all non-zero entries of a map $\partial_{i} \colon \begin{td} \AP_i \ar{r} \& \AP_{i-1} \end{td}$.
\end{enumerate}
This yields a sequence 		of projective $\A$-modules and $\A$-linear homomorphisms 
\begin{align*}
	\CP(\omega) = (			\begin{cd}
		\AP_{b} \ar{r}{\partial_{b}} \& \AP_{b-1} \ar{r} \& \ldots \& \AP_{a+1 }\ar{r}{\partial_{a+1}} \& \AP_a
	\end{cd}	)
\end{align*}
such that $[a,b] = \supp(d,w)$.
It is straightforward to check that $\partial_{i+1} \partial_i = 0$ for each degree $i \in \Z$, that is, $\CP(\omega)$ is a complex in $\DbRep{\A}$.

\subsubsection{Band and string complexes as cycle and line diagrams}
\label{subsec:hor-not}
In this subsection, we present a description of the complex $\CP(\omega)$ 
`before folding'.
Aside from its brevity, this notation will be useful to describe certain morphisms of complexes.

\paragraph{Band complexes via cycle diagrams}
Assume that $\omega$ is a band $(d,w,m,\lambda)$. Denoting $\dot{w} = (x_1, x_2, \ldots, x_{\ell-1}, x_{\ell})$, the \emph{cycle diagram} has the form 
\begin{align}
	\label{eq:cyc-bd-diag}
	\begin{gd}
		\underset{d_0}{{P_{+-}^{(m)} }}
		\ar[<->]{r}{M_1  t^{n_1}}  \& \underset{d_1}{P_{+-}^{(m)}} 
		\ar[densely dotted,-]{r}\& \underset{d_{i-1}}{P_{+-}^{(m)}}  \ar[<->]{r}{M_i   t^{n_i}} 
		\&
		\underset{d_{i}}{P_{+-}^{(m)}}
		\ar[densely dotted,-]{r}
		\&
		\underset{d_{\ell-1}}{P_{+-}^{(m)}}
		\ar[ <->,bend left=20]{llll}[inner sep=0.5pt]{ M_{\ell}   t^{n_{\ell}}}  
	\end{gd}
\end{align}
where 
$P_{+-}^{(m)} = P_+^m \oplus P_-^m$, $M_i$ are certain matrices explained below,
the horizontal map below each decoration $M_i t^{n_i}$ 
has the same orientation as the oriented number $x_i = \olra{n}_i$ of $\dot{w}$ if $1 \leq i < \ell$,
and the orientation opposite to $\olra{n}_{\ell}$ in case $i = \ell$ due to the depiction above.
The matrices $M_i$ are completely determined
by the sequence $w^{\updownarrow} = (x_i \kappa_i)_{i=1}^{\ell}$
and described in Table~\ref{tab:bd-mat}.
\begin{longtable}{|c|c||c|c||c|c|}
	\caption{Coefficient matrices for cycle diagrams of bands} 	\label{tab:bd-mat} \\
	\hline
	$(\kappa_{\ell}, x_1, \kappa_1)$ & $M_1$ &
	$(\kappa_{i-1}, x_i, \kappa_i)$ & $M_i$ &
	$(\kappa_{\ell-1}, x_\ell, \kappa_\ell)$ & $M_\ell$ \\
	\hline
	\hline
	& & & & & \\[-5pt]
	$(\downarrow, \ora{n}_1, \downarrow)$ &
	$\begin{pNiceArray}{cc}[columns-width=15pt] 0 & \Id \\ 0 & 0 \end{pNiceArray}$ &
	$(\downarrow, \ora{n}_i, \downarrow)$ &
	$\begin{pNiceArray}{cc}[columns-width=15pt] 0 & \Id \\ 0 & 0 \end{pNiceArray}$ &
	$(\downarrow, \ora{n}_\ell, \downarrow)$ &
	$\begin{pNiceArray}{cc}[columns-width=15pt] 0 & \Id \\ 0 & 0 \end{pNiceArray}$ 
	\\[15pt]
	$(\uparrow, \ora{n}_1, \downarrow)$ &
	$\begin{pNiceArray}{cc}[columns-width=15pt] \sign_1 \JB_{\lambda} & \Id \\ 0 & 0 \end{pNiceArray}$ &
	$(\uparrow, \ora{n}_i, \downarrow)$ &
	$\begin{pNiceArray}{cc}[columns-width=15pt]  \sign_i  \Id & \Id \\ 0 & 0 \end{pNiceArray}$ &
	$(\uparrow, \ora{n}_\ell, \downarrow)$ &
	$\begin{pNiceArray}{cc}[columns-width=15pt]  \sign_\ell  \Id & \Id \\ 0 & 0 \end{pNiceArray}$ \\[15pt]
	$(\downarrow, \ora{n}_1, \uparrow)$ &
	$\begin{pNiceArray}{cc}[columns-width=15pt] 0 & \Id \\ 0 & \sign_1 \Id \end{pNiceArray}$ &
	$(\downarrow, \ora{n}_i, \uparrow)$ &
	$\begin{pNiceArray}{cc}[columns-width=15pt] 0 & \Id \\ 0 & \sign_i \Id \end{pNiceArray}$ &
	$(\downarrow, \ora{n}_\ell, \uparrow)$ &
	$\begin{pNiceArray}{cr}[columns-width=15pt] 0 & \Id \\ 0 & \sign_\ell \JB_{\lambda} \end{pNiceArray}$ \\[15pt]
	$(\uparrow, \ora{n}_1, \uparrow)$ &
	$\begin{pNiceArray}{cc}[columns-width=15pt] \sign_1 \JB_{\lambda} & \Id \\ \JB_{\lambda} & \sign_1 \Id \end{pNiceArray}$ &
	$(\uparrow, \ora{n}_i, \uparrow)$ &
	$\begin{pNiceArray}{cc}[columns-width=15pt] \sign_i \Id & \Id \\ \Id &  \sign_i \Id \end{pNiceArray}$ &
	$(\uparrow, \ora{n}_\ell, \uparrow)$ &
	$\begin{pNiceArray}{cc}[columns-width=15pt] \sign_\ell \Id & \Id \\ \JB_{\lambda} &  \sign_\ell \JB_{\lambda} \end{pNiceArray}$ \\[10pt]
	\hline
	& & & & & \\[-5pt]
	$(\uparrow, \ola{n}_1, \uparrow)$ &
	$\begin{pNiceArray}{cc}[columns-width=15pt] 0 & 0 \\ \Id & 0 \end{pNiceArray}$ &
	$(\uparrow, \ola{n}_i, \uparrow)$ &
	$\begin{pNiceArray}{cc}[columns-width=15pt] 0 & 0 \\ \Id & 0 \end{pNiceArray}$ &
	$(\uparrow, \ola{n}_\ell, \uparrow)$ &
	$\begin{pNiceArray}{cc}[columns-width=15pt] 0 & 0 \\ \Id & 0 \end{pNiceArray}$ \\[15pt]
	$(\uparrow, \ola{n}_1, \downarrow)$ &
	$\begin{pNiceArray}{cc}[columns-width=15pt] 0 & 0 \\ \Id & \ol{\sign}_1 \Id \end{pNiceArray}$ &
	$(\uparrow, \ola{n}_i, \downarrow)$ &
	$\begin{pNiceArray}{cc}[columns-width=15pt] 0 & 0 \\ \Id & \ol{\sign}_i \Id \end{pNiceArray}$ &
	$(\uparrow, \ola{n}_\ell, \downarrow)$ &
	$\begin{pNiceArray}{cc}[columns-width=15pt] 0 & 0 \\ \Id & \ol{\sign}_{\ell} \JB_{\frac{1}{\lambda}} \end{pNiceArray}$ \\[15pt]
	$(\downarrow, \ola{n}_1, \uparrow)$ &
	$\begin{pNiceArray}{cc}[columns-width=15pt] \ol{\sign}_1 \JB_{\frac{1}{\lambda}} & 0 \\ \Id & 0 \end{pNiceArray}$ &
	$(\downarrow, \ola{n}_i, \uparrow)$ &
	$\begin{pNiceArray}{cc}[columns-width=15pt] \ol{\sign}_i \Id & 0 \\ \Id & 0 \end{pNiceArray}$ &
	$(\downarrow, \ola{n}_\ell, \uparrow)$ &
	$\begin{pNiceArray}{cc}[columns-width=15pt] \ol{\sign}_\ell \Id & 0 \\ \Id & 0 \end{pNiceArray}$ \\[15pt]
	$(\downarrow, \ola{n}_1, \downarrow)$ &
	$\begin{pNiceArray}{cc}[columns-width=15pt] \ol{\sign}_1 \JB_{\frac{1}{\lambda}} & \JB_{\frac{1}{\lambda}} \\ \Id & \ol{\sign}_1 \Id \end{pNiceArray}$ &
	$(\downarrow, \ola{n}_i, \downarrow)$ &
	$\begin{pNiceArray}{cc}[columns-width=15pt] \ol{\sign}_i \Id & \Id \\ \Id & \ol{\sign}_i \Id \end{pNiceArray}$ &
	$(\downarrow, \ola{n}_\ell, \downarrow)$ &
	$\begin{pNiceArray}{cc}[columns-width=15pt] \ol{\sign}_{\ell} \Id & \JB_{\frac{1}{\lambda}} \\ \Id & \ol{\sign}_{\ell} \JB_{\frac{1}{\lambda}} \end{pNiceArray}$ \\
	\hline
	\multicolumn{6}{|c|}{} \\[-10pt]
	\multicolumn{6}{|c|}{
		\parbox{0.95\linewidth}{
			where $\Id$ denotes the identity matrix of size $m$, $\JB_{\lambda}$ and $\JB_{\frac{1}{\lambda}}$ are Jordan blocks of size $m$ with eigenvalues $\lambda$ respectively $\frac{1}{\lambda}$,  
			in the middle two columns we assume that $1 < i < \ell$, and
			for any index $1 \leq j \leq \ell$ we set
			$\sign_j = (-1)^{d_j}$ and $\ol{\sign}_j = - \sign_j$.
		}
	} \\
	\hline
\end{longtable}
The matrices of Table~\ref{tab:bd-mat}
were obtained by the application of the sink-and-source rules of the previous section to the oriented gluing diagram of the band $\omega$.

			\paragraph{String complexes with trivial multiplicity via line diagrams}
			Let $\omega$ be a usual, special or bispecial primitive string $(d,w)$ of $\chfXA$. 
			Then $\CP(\omega)$ can be represented by the diagram
			\begin{align}
				\label{eq:lindiag2}
				\begin{gd}
					\underset{d_0}{P_{\alpha}} 
					\ar[<->]{r}{M_1 t^{n_1}}  \& \underset{d_1}{P_{+-}}
					\ar[densely dotted,-]{r}\& \underset{d_{i-1}}{P_{+-}} \ar[<->]{r}{M_i t^{n_i}} 
					\&
					\underset{d_i}{P_{+-}}
					\ar[densely dotted,-]{r}
					\&
					\underset{d_{\ell-1}}{P_{+-}}
					\ar[<->]{r}{ M_{\ell} t^{n_{\ell}}} \& \underset{d_{\ell}}{P_{\beta}}
				\end{gd}
			\end{align}
			where $P_{+-} = P_+ \oplus P_-$, the orientation of each horizontal map is the same as that of the oriented number $x_i$ 
			in $w$, and the matrices $M_i$ are given as follows:
			for each intermediate index $1 < i < k$ 
			the matrix $M_i$ depends on the subword $(\kappa_{i-1}, x_i,\kappa_i)$ of $w^{\updownarrow}$ and determined using the middle two columns of Table~\ref{tab:bd-mat} setting $m=1$;
			the first and the last matrix $M_1$ and $M_{\ell}$ are determined
			by the beginning and the end of
			the sequence $w^{\updownarrow}$ via the following table.
			\begin{align*}
				\begin{array}{|c|c|c|c||c|c|c|c|}
					\hline
					(x_1,\kappa_1) & M_1 & 	(x_1,\kappa_1) & M_1	&	(\kappa_{\ell-1},x_{\ell})	& 	M_\ell &(\kappa_{\ell-1},x_{\ell})		& 	M_\ell
					\\ 
					\hline
					\hline
					& & &  &&&& \\[-5pt] 
					({\ora{n}_1},\downarrow) & 
					\begin{pNiceArray}[columns-width=0pt]{c}
						1 \\
						0
					\end{pNiceArray} 
					&
					({\ora{n}_1},\uparrow) &
					\begin{pNiceArray}[columns-width=0pt]{c}
						1 \\
						\sign_1
					\end{pNiceArray} 
					&		
					(\uparrow, {\ola{n}_\ell}) &
					\begin{pNiceArray}[columns-width=0pt]{c}
						\ol{\sign}_\ell   \\ 	1
					\end{pNiceArray}
					&
					(\downarrow,{\ola{n}_\ell}) & 
					\begin{pNiceArray}[columns-width=0pt]{c}
						0 \\
						1
					\end{pNiceArray} 
					\\[15pt]
					({\ola{n}_1},\downarrow) & 
					\begin{pNiceArray}[columns-width=0pt]{cc}
						1 & \ol{\sign}_1
					\end{pNiceArray} 
					&
					({\ola{n}_1},\uparrow) 
					&
					\begin{pNiceArray}[columns-width=5pt]{cc}
						1 &
						0
					\end{pNiceArray}	
					&
					(\uparrow, {\ora{n}_\ell}) & 
					\begin{pNiceArray}[columns-width=5pt]{cc}
						0 &
						1
					\end{pNiceArray}
					&
					(\downarrow,{\ora{n}_\ell}) &
					\begin{pNiceArray}[columns-width=5pt]{cc}
						{\sign}_\ell   & 1
					\end{pNiceArray}\\[5pt]
					\hline
				\end{array}
			\end{align*}	
			As for bands, these matrices originate from the induced arrows of the previous subsection.

	\paragraph{Compressed notation for bispecial strings}
	Assume that $\omega$ is a non-primitive bispecial string $(d,\varepsilon_1, w, \varepsilon_2)$.
	Let $m$ denote the multiplicity of the 
	primitive root $v$ of $w$.
	On the one hand, the complex $\CP(\omega)$ admits a presentation like \eqref{eq:lindiag2}
	but with certain changes of signs in the intermediate projectives as described in Subsection~\ref{subsec:ind-arr}.
	On the other hand, 
	the complex $\CP(\omega)$ admits
	a more compact notation involving block matrices of size $m$, which is similar to those of bands and described next.

	For any sign $\varepsilon \in \{ +,-\}$
	we define a projective $\A$-module $P^{(m,\varepsilon)}$ as the alternating direct sum
	$P_{\varepsilon} \oplus P_{\ol{\varepsilon}} \oplus P_{\varepsilon} \oplus 
	P_{\ol{\varepsilon}} 
	\ldots    
	$
	with $m$ indecomposable summands; the last summand of 
	$P^{(m,\varepsilon)}$
	is given by $P_\varepsilon$ if $m$ is odd and $P_{\ol{\varepsilon}}$ if $m$ is even.
	As in the case of bands, we set $P_{+-}^{(m)} = P_+^m \oplus P_-^m$.
	
	Let $\ell$ be the length of the primitive root of $w$.
	We refer to Subsection~\ref{subsec:bisp-str1} for the case $\ell = 1$, and assume $\ell > 1$ below.
	Then $\CP(\omega)$ can be depicted by the diagram
	\begin{align}
		\label{eq:bisp-lin}
		\begin{gd}
			\underset{d_0}{P^{(m,\varepsilon_1)}} 
			\ar[<->]{r}{M_1  t^{n_1}}  \& \underset{d_1}{P_{+-}^{(m)}}
			\ar[densely dotted,-]{r}\& \underset{d_{i-1}}{P_{+-}^{(m)}}  \ar[<->]{r}{M_i   t^{n_i}} 
			\&
			\underset{d_{i}}{P_{+-}^{(m)}}
			\ar[densely dotted,-]{r}
			\&
			\underset{d_{\ell-1}}{P_{+-}^{(m)}}
			\ar[<->]{r}{ M_{\ell}   t^{n_{\ell}}}
			\& \underset{d_\ell}{P^{(m,\varepsilon_2)} }
		\end{gd}
	\end{align}
	with the following matrices:
	for any index $1 < i < \ell$ the matrix $M_i$ is a square matrix of size $2m$ and determined using the middle two columns of Table~\ref{tab:bd-mat};
	the remaining matrices $M_1$ and $M_{\ell}$ at the ends are given by the next table.

	\begin{align*}
		\begin{array}{|c|c|c|c||c|c|c|c|}
			\hline
			(x_1,\kappa_1) & M_1 & 	(x_1,\kappa_1) & M_1	&	(\kappa_{\ell-1},x_{\ell})	& 	M_\ell &(\kappa_{\ell-1},x_{\ell})		& 	M_\ell
			\\ 
			\hline
			\hline
			& & &  &&&& \\[-5pt] 
			({\ora{n}_1},\downarrow) & 
			\begin{pNiceArray}[columns-width=0pt]{c}
				A \\
				0
			\end{pNiceArray} 
			&
			({\ora{n}_1},\uparrow) &
			\begin{pNiceArray}[columns-width=0pt]{c}
				A \\
				\sign_1 A
			\end{pNiceArray} 
			&		
			(\uparrow, {\ola{n}_\ell}) &
			\begin{pNiceArray}[columns-width=0pt]{c}
				0 \\
				\bar{B}
			\end{pNiceArray} 
			&
			(\downarrow,{\ola{n}_\ell}) & 
			\begin{pNiceArray}[columns-width=0pt]{c}
				\ol{\sign}_\ell \bar{B}  \\ 	\bar{B}
			\end{pNiceArray}
			\\[15pt]
			({\ola{n}_1},\downarrow) & 
			\begin{pNiceArray}[columns-width=0pt]{cc}
				\bar{A} & \ol{\sign}_1 \bar{A}
			\end{pNiceArray} 
			&
			({\ola{n}_1},\uparrow) 
			&
			\begin{pNiceArray}[columns-width=0pt]{cc}
				\bar{A} &
				0
			\end{pNiceArray}	
			&
			(\uparrow, {\ora{n}_\ell}) & 
			\begin{pNiceArray}[columns-width=5pt]{cc}
				{\sign}_\ell  {B} & {B}
			\end{pNiceArray}
			&
			(\downarrow,{\ora{n}_\ell}) &
			\begin{pNiceArray}[columns-width=0pt]{cc}
				0 &
				{B}
			\end{pNiceArray}
			\\[5pt]
			\hline
			\multicolumn{8}{|c|}{}
			\\[-10pt] 
			\multicolumn{8}{|c|}{
				\parbox{0.95\linewidth}{
					\text{where the matrices $A$ and $B$  are given by the bidiagonal matrices of size $m$ }\\
					\quad 	$A = 
					\Id + \smashoperator{\sum_{1 \leq i < \frac{m}{2}}} \sign_1   E_{2i,2i+1}
					=
					\begin{psmallmatrix}
						1& 0 \\
						& 1 & \sign_1 \\
						& & 1 & 0\\
						& & & 1 & \scalebox{0.5}{$\ddots$} \\
						&&&& \scalebox{0.5}{$\ddots$}
					\end{psmallmatrix}, \quad
					B = 
					\Id + \smashoperator{\sum_{1\leq i \leq \frac{m}{2}}} \sign_\ell E_{2i-1,2i}
					=
					\begin{psmallmatrix}
						1& \sign_{\ell} \\
						& 1 & 0 \\
						& & 1 & \sign_{\ell }\\
						& & & 1 & \scalebox{0.5}{$\ddots$} \\
						&&&&\scalebox{0.5}{$\ddots$}
					\end{psmallmatrix}$,\\
					\text{
						the matrices $\bar{A}$ and $\bar{B}$ are defined by replacing all signs $\sign_1$ in $A$ and $\sign_\ell$ in $B$ }\\
					by $\ol{\sign}_1$ and $\ol{\sign}_{\ell}$, respectively. 
				}
			}
			\\
			\hline
		\end{array}
	\end{align*}

In particular, the matrices of $\CP(\omega)$ are determined by the multiplicity $m$ of the underlying word $w$
and the oriented version
$v^{\updownarrow}$ of its primitive root.

\begin{rmk}
	It can be verified that for any bispecial string $\omega$ the line diagram of the bispecial string complex $\CP(\omega)$ 
	in the previous subsection is
	the same as the presentation in the present subsection in case $\ell >1$ respectively Subsection~\ref{subsec:bisp-str1} in case $\ell =1$ up to renumbering of columns and rows. 
\end{rmk}

\subsection{Examples}\label{subsec:examples}
We recall that 
the indecomposable projective
$\A$-modules admit non-invertible morphisms
which are given by right multiplication with monomials
\begin{align*}
	\begin{cd}
		P_{u} \ar{r}{t^n} \& P_{v}
	\end{cd}
	&&
	\text{ where }
	P_\star = \begin{pmatrix} 
		\Rx \\
		\Rx \\
		\Rx
	\end{pmatrix}, \
	P_+ = 
	\begin{pmatrix} 
		\mx \\
		\Rx \\
		\mx
	\end{pmatrix}, \
	P_- = 
	\begin{pmatrix} 
		\mx \\
		\mx \\
		\Rx
	\end{pmatrix}
\end{align*}
with $u,v \in Q_0$ and $n \in \N_0$ if $(u,v)=(+,\star)$ or $(-,\star)$ and $n \in \N$ otherwise.

\subsubsection{A  usual string of width three}
\label{subsec:usual}
Let
$\omega$ be the usual string $(2,w)$
with $w= (\star,\ora{p}, \ora{q},\star)$
with any $p\in \N, q\in \N_0$.
Then $w^{\updownarrow} = (\star, \ora{p} \uparrow \ora{q}, \star)$, and the initial and the final diagrams of $\omega$, 
and the string complex $\CP(\omega)$ are given as follows.
\begin{align*}
	\begin{gd}
		\& \bt \ar[densely dotted,-]{d} \& \star \\
		\star \ar{ru}{p} 		\& \bt \ar{ru}{q} 
	\end{gd}
	&\Rightarrow&
	\begin{gd}
		\& P_+ \ar[color=green]{r}{t^q} \ar[densely dotted,<-]{d} \& \smash[b]{\underset{0}{P_\star}} \\
		\underset{2}{P_\star} \ar{ru}{t^p} \ar[color=red]{r}[swap]{-t^p} \& \underset{1}{P_-} \ar{ru}[swap]{t^q} 
	\end{gd}
	&\Rightarrow &  
	(\begin{cd}
		\underset{2}{P_{\star}} \ar{r}{
			\begin{psmallmatrix}
				\phantom{-}t^p \\ -t^p \end{psmallmatrix}		
		}\&
		\underset{1}{P_{+-}} \ar{r}{
			\begin{psmallmatrix}  t^{q} & t^{q} \end{psmallmatrix}	
		}
		\& \underset{0}{P_{\star}}
	\end{cd}).
\end{align*}
If $(p,q)=(1,0)$, the complex
$\CP(\omega)$ is a projective resolution of the simple module $S_{\star}$.

\subsubsection{An equioriented special string}
Let $\omega$ be the special string $(0,w,\varepsilon)$ with 
any sign $\varepsilon\in \{+,-\}$ and
$w =(\diamond,{\ola{p}_1},{\ola{q}_1},{\ola{p}_2},{\ola{q}_2},\star)$, 
where	$p_1 \in \N_0$ and $q_1,p_2,q_2 \in \N$.
Then
$w^{\updownarrow} = 
(\diamond,{\ola{p}_1} \downarrow {\ola{q}_1} \downarrow {\ola{p}_2} \downarrow {\ola{q}_2},\star)$, and the diagrams of $\omega$ and $\CP(\omega)$ are given below.
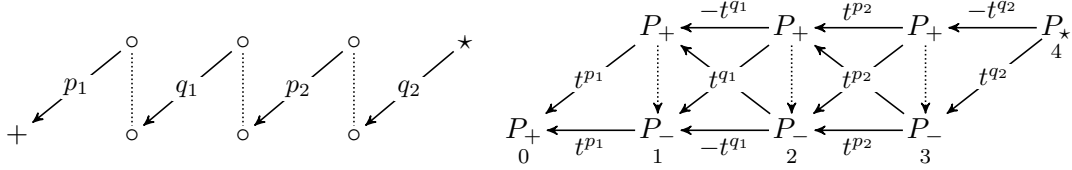
\begin{align*}
	\begin{gd}
		\& \bt \ar[densely dotted,-]{d} \&  \bt \ar[densely dotted,-]{d} \&\bt \ar[densely dotted,-]{d} \& \star \\
		+ \ar[<-]{ru}[description]{p_1} 		\& \bt \ar[<-]{ru}[description]{q_1} \&  \bt \ar[<-]{ru}[description]{p_2} \& \bt \ar[<-]{ru}[description]{q_2} 
	\end{gd}
	&&
	\begin{gd}
		\& P_+ \ar[densely dotted]{d} \ar[<-]{r}{-t^{q_1}} \ar[<-]{rd} \&  P_+ \ar[densely dotted]{d} \ar[<-]{r}{t^{p_2}} \ar[<-]{rd} \& P_+ \ar[densely dotted]{d}  \ar[<-]{r}{-t^{q_2}} \& \smash[b]{\underset{4}{P_\star}} \\
		\underset{0}{P_+} \ar[<-]{ru}[description]{t^{p_1}}  \ar[<-]{r}[swap]{t^{p_1}}		\& \underset{1}{P_-} 
		\ar[<-]{ru}[description]{t^{q_1}}  \ar[<-]{r}[swap]{-t^{q_1}}
		\&  \underset{2}{P_-} 
		\ar[<-]{ru}[description]{t^{p_2}}  \ar[<-]{r}[swap]{t^{p_2}}
		\& \underset{3}{P_-} 	 	\ar[<-]{ru}[description]{t^{q_2}}  
	\end{gd}
\end{align*}
In terms of notation of Subsection~\ref{subsec:hor-not}
we obtain the line diagram
\begin{align*}
	\begin{gd}
		\underset{0}{P_+} \ar[<-]{r}{
			\begin{pmatrix}
				t^{p_1} & t^{p_1}
			\end{pmatrix}	
		} \&[10pt] \underset{1}{P_{+-}} \ar[<-]{r}{ 
			\begin{pmatrix}
				-t^{q_1} &  t^{q_1} \\
				t^{q_1} & -t^{q_1}
			\end{pmatrix}	
		} \&[10pt] \underset{2}{P_{+-}} \ar[<-]{r}{
			\begin{pmatrix}
				t^{p_2} & t^{p_2} \\
				t^{p_2} & t^{p_2}
			\end{pmatrix}	
		} \&[10pt] \underset{3}{P_{+-}} \ar[<-]{r}{
			\begin{pmatrix}
				-t^{q_2} \\ t^{q_2} 
			\end{pmatrix}	
		} \&[10pt] \underset{4}{P_{\star}}
	\end{gd}
\end{align*}
Reading the line diagram from right to left yields the string complex $\CP(\omega)$.
In particular, the above matrices happen to coincide with the differentials of the complex $\CP(\omega)$.

\subsubsection{A usual string with explicit numbers}
\label{subsec:ex-us2}
We consider an explicit usual string $\omega = (2,w)$
with
$w = (\star,\ora{2},\ora{1},\ola{2},\ora{3},\ora{0},\star)$.
By considering pairs of successive oriented numbers, 
a straightforward application of the orientation rule \eqref{eq:updownarrow} yields that $w^{\updownarrow} = 
(\star,\ora{2} \uparrow \ora{1} \uparrow \ola{2} \downarrow \ora{3} \uparrow \ora{0},\star)$.
In particular, the initial gluing diagram of $w$ 
has the following form.

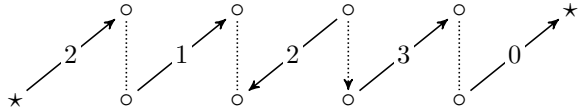
\begin{align*}
		\begin{gd}
			\& \bt \ar[densely dotted,-]{d}    \& \bt \ar[densely dotted,-]{d}  \& \bt \ar[densely dotted]{d} \& \bt  \ar[densely dotted,-]{d} \& \star \\
			\star \ar{ru}[description]{2} 
			\& \bt \ar{ru}[description]{1} \&  \bt \ar[<-]{ru}[description]{2}  \& \bt \ar{ru}[description]{3} \& \bt \ar{ru}[description]{0} 
		\end{gd}
	\end{align*}

	Its final version is given as follows.
	
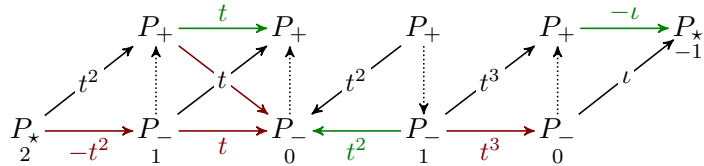
\begin{align*}
		\begin{gd}
			\& P_+ \ar[densely dotted,<-]{d} \ar[red]{rd} \ar[green]{r}{t}  \& P_+ \ar[densely dotted,<-]{d}  \& P_+ \ar[densely dotted]{d} \& P_+ \ar[green]{r}{-\iota} \ar[densely dotted,<-]{d} \& \smash[b]{\underset{-1}{P_{\star}}} \\
			\underset{2}{P_{\star}} \ar{ru}[description]{t^2} \ar[red]{r}[swap]{-t^2}
			\& \underset{1}{P_-} \ar{ru}[description]{t} \ar[red]{r}[swap]{t}\&  \underset{0}{P_-} \ar[<-]{ru}[description]{t^2} \ar[<-,green]{r}[swap]{t^2} \& \underset{1}{P_-} \ar{ru}[description]{t^3} \ar[red]{r}[swap]{t^3}\& \underset{0}{P_-} \ar{ru}[description]{\iota} 
		\end{gd}  
\end{align*}
Relocating the projectives in the last diagram such that
projectives with a common degree are located in the same column yields
the diagram on the left. It can be used as an intermediate step when `folding' the diagram above in order to obtain the string complex $\CP(\omega)$ which is described below.
\begin{align*}
	\begin{gd}
		P_\star \ar{r}{t^2} \ar[color=red]{rd}[swap]{-t^2}
		\& P_+ \ar[densely dotted,<-]{d} \ar[color=red]{rdd} \ar[color=green]{rd}{t} \& \& \\
		\& P_- \ar{r}[description]{t} \ar[color=red]{rd}[swap]{t} \& P_+ \ar[densely dotted,<-]{d} \& \\
		\& P_+ \ar{r}[description]{t^2} \ar[densely dotted,->]{d} \& P_-  \& \\
		\& P_- \ar{r}[description]{t^3} \ar[color=green]{ru}[description]{t^2} \ar[color=red]{rd}[swap]{t^3} \& P_+ \ar[densely dotted,<-]{d} \ar[color=green]{rd}{-\iota} \& \\
		\underset{2}{\mathstrut} \& \underset{1}{\mathstrut} \& \underset{0}{P_-} \ar{r}[swap]{\iota} \& \underset{-1}{P_\star}
	\end{gd}
	\begin{array}{cl}
		\Rightarrow& 
		(\begin{cd}
			P_{\star} \ar{r}{\partial_2} 
			\& P_{+-}^2 \ar{r}{\partial_1} \&
			P_{+-}^2 \ar{r}{\partial_0} \& P_{\star}
		\end{cd})
		\\[10pt]
		&\text{where $P_{+-}^2 = P_+ \oplus P_- \oplus P_+ \oplus P_-$} \\[10pt]
		&
		\partial_2  =\begin{pNiceArray}{c}
			t^2   \\
			-t^2  \\
			0  \\
			0  
		\end{pNiceArray}
		\quad
		\partial_1 = \begin{pNiceArray}[columns-width=auto]{cccc}
			t & t & 0&0  \\
			t & t & t^2 & t^2  \\ 
			0 & 0 & 0 & t^3\\
			0& 0 & 0 & t^3 
		\end{pNiceArray} \\
		\\
		&\quad \partial_0 = \begin{pNiceArray}{cccc}
			0 & 0 & -\iota & \iota 
		\end{pNiceArray},
		\\ 
		\\
		&
		\text{and 
			$\begin{td} \iota \colon P_{\pm} \ar[hookrightarrow]{r} \& P_\star\end{td}$
			denotes the natural inclusion}.
	\end{array}
\end{align*}
Replacing the sequence $(2,1,2,3,0)$ in $w$ by $(p_1,q_1,p_2,q_2,p_3) \in \N^{\times 4} \times \N_0$
such that $q_1 < p_2 < q_2$ yields the same diagrams. 
The differentials of the corresponding complex can be obtained by adjusting the monomial powers.
\subsubsection{A special string of width two}
\label{subsec:sp-str}
Let $\omega$ be the special string $(2,w,-)$ with underlying word
$w=  (\star,{\ora{q}_0},{\ora{p}_1},{\ola{q}_1},{\ora{p}_2}, {\ola{q}_2},\diamond)$
such that $q_0 \in \N$  and $p_1 = q_1 = p_2 = q_2 \in \N$.
Denoting $\tilde{q}_0 = q_0 - \frac{1}{2}$,
the ambient word of $w$ is given by 
$$\overline{w} =  (\star,\ora{\tilde{q}}_0,\underbrace{{\ora{p}_1},{\ola{q}_1},{\ora{p}_2}, {\ola{q}_2},
	{\color{blue}{\ora{q}_2},{\ola{p}_2}, {\ora{q}_1}, {\ola{p}_1},}}_{\text{equal numbers}} {\color{blue}\ola{\tilde{q}}_0},\star).$$
The ambient word admits a diagrammatic description
\begin{align*}
	\scalebox{0.9}{$
		\begin{gd}
			\& \bt \ar[densely dotted,-]{d}    \& \bt \ar[densely dotted,-]{d}  \& \bt \ar[densely dotted]{d} \& \bt  \ar[densely dotted,-]{d} \& \bt  \ar[densely dotted,-]{d} 
			\& \bt \ar[densely dotted,-]{d}    \& \bt \ar[densely dotted,-]{d}  \& \bt \ar[densely dotted]{d} \& \bt  \ar[densely dotted,-]{d} \& \star 	\\
			\star \ar{ru}[description]{{\tilde{q}_0}} 
			\& \bt \ar{ru}[description]{{p_1}} \&  \bt \ar[<-]{ru}[description]{{q_1}}  \& \bt \ar{ru}[description]{{p_2}} \& \bt \ar[<-]{ru}[description]{{q_2}} 
			\&
			\bt \ar{ru}[description]{{q_2}} 
			\& \bt \ar{ru}[description]{{p_2}} \&  \bt \ar[<-]{ru}[description]{{q_1}}  \& \bt \ar{ru}[description]{{p_1}} \& \bt \ar[<-]{ru}[description]{{\tilde{q}_0}} 
		\end{gd}$}
\end{align*}
To determine the first orientation arrow in $w^{\updownarrow}$ 
we compare the arrow $\wt{q}_0$ with the arrow $\wt{p}_1$.
As both are oriented to the right, the orientation avoiding transition vertices is upwards.
Determining any other orientation arrow, we need to ignore a maximal symmetric subdiagram at that arrow, which boils down to
comparing the right-oriented arrow  with $\wt{q}_0$ to some other arrow.
Therefore, it follows that
$w=  (\star,{\ora{q}_0} \uparrow {\ora{p}_1} \uparrow {\ola{q}_1} \uparrow {\ora{p}_2} \uparrow {\ola{q}_2}, -)$.
Applying the sink and source rules yields the following diagram.
\begin{align}
	\label{eq:sp-str-d}
			\begin{gd}
				\& P_+ \ar[densely dotted,<-]{d} \ar{r}{t^{p_1}}  \ar{rd} \& P_+ \ar[densely dotted,<-]{d}  \& P_+ \ar[densely dotted,<-]{d} \ar{rd} \ar{r}{t^{p_2}} \& P_+  \ar[densely dotted,<-]{d} \& P_- \\
				\underset{2}{P_\star} \ar{ru}[description]{t^{q_0}} \ar{r}[swap]{-t^{q_0}}
				\& \underset{1}{P_-} \ar{ru}[description]{t^{p_1}} \ar{r}[swap]{t^{p_1}}\&  \underset{0}{P_-} \ar[<-]{ru}[description]{t^{q_1}}  \& \underset{1}{P_-} \ar{ru}[description]{t^{p_2}} \ar{r}[swap]{t^{p_2}}\& \underset{0}{P_-} \ar[<-]{ru}[description]{t^{q_2}} 
			\end{gd}
		\end{align}
		In different terms, the diagram can be subsumed by the line diagram.
		\begin{align*}
			\begin{gd}
				\underset{2}{P_\star} \ar{r}{\begin{pmatrix}
						t^{q_0} \\
						-t^{q_0} \end{pmatrix}} \& \underset{1}{P_{+-}} \ar{r}{ 			\begin{pmatrix}
						t^{p_1} & t^{p_1}\\
						t^{p_1} & t^{p_1}\end{pmatrix}} \& \underset{0}{P_{+-}} \ar[<-]{r}{ 			\begin{pmatrix}
						0 & 0  \\
						t^{q_1} & 0 \end{pmatrix}} \& \underset{1}{P_{+-}} \ar{r}{ 			\begin{pmatrix}
						t^{p_2} & t^{p_2}\\
						t^{p_2} & t^{p_2}\end{pmatrix}} \&
				\underset{0}{P_{+-}} \ar[<-]{r}{ 			\begin{pmatrix}
						0 \\
						t^{q_2} \end{pmatrix}} \& \underset{1}{P_{-}}
			\end{gd} 
\end{align*}
Folding the line diagram yields that the string complex $\CP(\omega)$ is the three-term complex
$(\begin{td}
P_\star \ar{r}{\partial_2} \& P^2_{+-} \oplus P_- \ar{r}{\partial_1} \&  P^2_{+-}
\end{td})$
with differentials
\begin{align}
\partial_2 = \begin{pNiceArray}{c}
	t^{q_0} \\
	-t^{q_0}\\
	0 \\
	0 \\
	0 
\end{pNiceArray},
\quad
\partial_1 = 
\begin{pNiceArray}{ccccc}
	t^{p_1} & t^{p_1} &0&0&0 \\
	t^{p_1}&t^{p_1}&t^{q_1}&0&0\\
	0&0&t^{p_2}&t^{p_2}&0\\
	0&0&t^{p_2}&t^{p_2}&t^{q_2} \\
\end{pNiceArray}.
\end{align}
It can be checked that we obtain the same complex assuming merely  $p_1 \leq q_1 \geq p_2 \leq q_2$.
\subsubsection{A bispecial string with multiplicity three}
Let $\omega$ be the bispecial string $(1,\varepsilon_1,w,\varepsilon_2)$ with 
arbitrary signs $\varepsilon_1, \varepsilon_2 \in \{+,-\}$
and
$w = (\diamond, \ora{1}, \ola{2}, \ola{3}, \ora{3}, \ora{2},\ola{1}, 
\ora{1}, \ola{2}, \ola{3},\diamond)$.
The primitive root of the latter is given by  $v=(\diamond,  \ora{1}, \ola{2}, \ola{3}, \diamond)$ and has multiplicity three in $w$.
The orientation rules yield $v^{\updownarrow} = 
(\diamond, 
\ora{1} 
\uparrow \ola{2} \downarrow 
\ola{3},\diamond)$.
The gluing diagram and the complex of the primitive root are then given as follows.
\begin{align*}
\begin{array}{ccc}
	\begin{gd}
		\& P_+ \ar[densely dotted,<-]{d}   \& P_+ \ar[densely dotted,->]{d} \ar[<-]{r}{-t^3} \& \smash[b]{\underset{2}{P_{\varepsilon_2}}} \\
		\underset{1}{P_{\varepsilon_1}} \ar{ru}[description]{t} \ar{r}[swap]{t}
		\& \underset{0}{P_-} \ar[<-]{ru}[description]{t^{2}} \ar[<-]{r}[swap]{t^{2}}\&  \underset{1}{P_-} \ar[<-]{ru}[description]{t^{3}}  
	\end{gd}
	&\qquad&
	\begin{cd}
		\\[0.4cm]
		\underset{2}{P_{\varepsilon_2}} \ar{r}[yshift=5pt]{\left(\begin{smallmatrix}
				\phantom{-}0 \\
				-t^3 \\
				\phantom{-}t^3
			\end{smallmatrix}\right)} \& \underset{1}{P_{\varepsilon_1} \oplus P_{+-}} \ar{r}[yshift=5pt]{\left(\begin{smallmatrix}
				\phantom{-}t & 0 & 0 \\
				-t & t^2 & t^2
			\end{smallmatrix}\right)} \& \underset{0}{P_{+-}}
	\end{cd} 
\end{array}
\end{align*}
From here, we have two possibilities to determine the complex $\CP(\omega)$.

\begin{enumerate}
\item To describe $\CP(\omega)$ via the compressed notation, we use 
the arrow orientations in $v^{\updownarrow}$,
and the tables of Subsection~\ref{subsec:hor-not} 
to determine the coefficient matrices 
for the line diagram~\eqref{eq:bisp-lin}.
This diagram is given by
\begin{align*}
	\begin{array}{cl}
		\begin{gd}
			P^{(3,\varepsilon_1)}
			\ar{r}{
				\begin{pmatrix}
					t A \\
					t A
				\end{pmatrix}	
			} \&
			P_{+-}^3
			\ar[<-]{r}{
				\begin{pmatrix}
					0 & 0 		
					\\
					t^2 & t^2
				\end{pmatrix}
			} \&
			P_{+-}^3
			\ar[<-]{r}{
				\begin{pmatrix}
					t^3 B \\ t^3 B
				\end{pmatrix}	
			} \&
			P^{(3,\varepsilon_2)}
		\end{gd}
		&
		\begin{array}{ll}
			\text{with }
			&
			A = \begin{pNiceArray}[small]{ccc}
				1 & 0 & 0 \\
				0 & 1 & 1 \\
				0 & 0 & 1  \\
			\end{pNiceArray},\quad
			B = \begin{pNiceArray}[small]{ccc}
				1 & 1 & 0 \\
				0 & 1 & 0 \\
				0 & 0 &1 
			\end{pNiceArray},
			\\
			&P^{(3, \varepsilon)}
			=P_{\varepsilon} \oplus P_{\ol{\varepsilon}} \oplus P_{\varepsilon},\\
			&	P_{+-}^3 = P_+^3 \oplus P_-^3.
		\end{array}
	\end{array}
\end{align*}
where all block matrices have size three due to the multiplicity of $v$ in $w$.
Folding the diagram yields the complex $\CP(\omega) = (
\begin{td}
	P^{(3,\varepsilon_2)} \ar{r}{\partial_2} \& 
	P^{(3,\varepsilon_1)} \oplus
	P_{+-}^3 \ar{r}{\partial_1} \& P_{+-}^3
\end{td}
)$
with the differentials:
\begin{align*}
			\partial_2 =
			{
				\begin{pNiceArray}[small, columns-width=auto,first-row,last-col]{ccc}
					\varepsilon_2 & \ol{\varepsilon}_2 & \varepsilon_2 & \\[3pt]
					0 &0&0  & \varepsilon_1\\
					0 &0&0  & \ol{\varepsilon}_1\\
					0 &0&0 & \varepsilon_1\\[3pt]
					t^3 & t^3 & 0 & +\\
					0 & t^3 & 0 & +\\
					0 & 0 & t^3 & +\\[3pt]
					t^3 & t^3 & 0 & -\\
					0 & t^3 & 0 & -\\
					0 & 0 & t^3 & -
			\end{pNiceArray}}\quad
			\partial_1 
			=
			{
				\begin{pNiceArray}[small, first-row, last-col, columns-width=auto]{ccccccccc}
					\varepsilon_1 & \ol{\varepsilon}_1 & \varepsilon_1 & + & + & + & - & - & - \\[3pt]
					t & 0 & 0 & 0 & 0 & 0 & 0 & 0 & 0 & +\\
					0 & t & t & 0 & 0 & 0 & 0 & 0 & 0 & +\\
					0 & 0 & t & 0 & 0 & 0 & 0 & 0 & 0 & +\\[3pt]
					t & 0 & 0 & t^2 & 0 & 0 & t^2 & 0 & 0 & -\\
					0 & t & t & 0 & t^2 & 0 & 0 & t^2 & 0 & -\\
					0 & 0 & t & 0 & 0 & t^2 & 0 & 0 & t^2 & -
			\end{pNiceArray}}
\end{align*} 
\item Now, we reconsider the bispecial string $\omega$ via the more conventional method.
Using the downward rule at axes of symmetry we obtain 
$$w^{\updownarrow} = 
(\diamond, 
\ora{1} 
\uparrow \ola{2} \downarrow 
\ola{3} 
\, {\color{blue}\downarrow}\,
\ora{3} 
\uparrow \ora{2} \downarrow 
\ola{1} \,{ \color{blue} \downarrow} \,
\ora{1} 
\uparrow \ola{2} \downarrow 
\ola{3},
\diamond)$$
and thus the gluing diagram
\begin{align*}
	\scalebox{0.85}{$
		\begin{gd}
			\& P_+ \ar[densely dotted,<-]{d}   \& P_+ \ar[densely dotted,->]{d} \ar[<-]{r}{-t^3} \ar[color=blue,<-]{rd} \& {P_{\varepsilon_2}} \ar[color=blue,densely dotted]{d}
			\& P_- \ar[densely dotted,<-]{d}  \ar{r}{t^2} \& P_- \ar[densely dotted]{d}  \ar[<-]{r}{t}  \ar[color=blue,<-]{rd} \& P_{\ol{\varepsilon}_1} \ar[color=blue,densely dotted]{d}
			\& P_+ \ar[densely dotted,<-]{d}   \& P_+ \ar[densely dotted,->]{d} \ar[<-]{r}{-t^3}  \& {P_{\varepsilon_2}}
			\\
			\underset{1}{P_{\varepsilon_1}} \ar{ru}[description]{t} \ar{r}[swap]{t}
			\& \underset{0}{P_-} \ar[<-]{ru}[description]{t^{2}} \ar[<-]{r}[swap]{t^{2}}\&  \underset{1}{P_-} \ar[<-]{ru}[description]{t^{3}} \ar[color=blue,<-]{r}[swap]{-t^{3}} \&
			\underset{2}{P_{\ol{\varepsilon}_2}} \ar{ru}[description]{t^3} \ar{r}[swap]{-t^3}
			\& \underset{1}{P_+} \ar{ru}[description]{t^{2}}  \&  \underset{0}{P_+} \ar[<-]{ru}[description]{t} \ar[color=blue,<-]{r}[swap]{t} \&
			\underset{1}{P_{\varepsilon_1}} \ar{ru}[description]{t} \ar{r}[swap]{t}
			\& \underset{0}{P_-} \ar[<-]{ru}[description]{t^{2}} \ar[<-]{r}[swap]{t^{2}}\&  \underset{1}{P_-} \ar[<-]{ru}[description]{t^{3}}
		\end{gd}
		$}
\end{align*}
Folding this diagram 
and reordering the projectives suitably yields the same complex $\CP(\omega)$ as above.
\end{enumerate}

\subsubsection{Bispecial strings with primitive roots of length one}\label{subsec:bisp-str1}
Assume that $\omega$ is a bispecial string $(d,\varepsilon_1, w, \varepsilon_2)$
with any $d \in \Z$, and signs $\varepsilon_1, \varepsilon_2 \in \{+,-\}$.
The primitive root $v$ of $w$ has length one if and only if
$w$ is given by a word $(\diamond, \olra{n}, \olra{n},\ldots \olra{n}, \diamond)$ with $n \in \N$ and alternating orientations.
We assume that this to be the case. Then the length of $w$ is precisely the multiplicity $m$ of $v$ in $w$.
There are the following two possibilities.
\begin{itemize}
\item Assume that $w$ begins with a right-oriented number, or, equivalently, $v=(\diamond, \ora{n}, \diamond)$.
The diagram of $w$ is given by
$$
\scalebox{0.85}{$
	\begin{gd}
		\& P_{\varepsilon_2} \ar[densely dotted]{d}  \ar[<-]{r}{t^n} \ar[<-]{rd} \& P_{\ol{\varepsilon}_1} \ar[densely dotted]{d}  \& P_{\varepsilon_2} \ar[densely dotted]{d} \ar[<-]{r}{t^n} \ar[<-]{rd} \& P_{\ol{\varepsilon}_1}  \ar[densely dotted]{d}  \& P_{\varepsilon_2}  \ar[<-]{r}{t^n} \ar[densely dotted]{d}  \& P_{\ol{\varepsilon}_1}  \ar[densely dotted,-]{rr} \& \& P_{\varepsilon_2}  \ar[<-, color=gray]{r}{t^n} \ar[color=gray,densely dotted]{d} \& {\color{gray} \underset{d}{P_{\ol{\varepsilon}_1}}} \\
		\underset{d}{P_{\varepsilon_1}} \ar{ru}[description]{t^n} 
		\& \underset{d-1}{P_{\ol{\varepsilon}_2}} \ar[<-]{ru}[description]{t^n} \ar[<-]{r}[swap]{t^n} \&  \underset{d}{P_{\varepsilon_1}} \ar{ru}[description]{t^n}  \& \underset{d-1}{P_{\ol{\varepsilon}_2}} \ar[<-]{ru}[description]{t^n} \ar[<-]{r}[swap]{t^n} \&  \underset{d}{P_{\varepsilon_1}} \ar{ru}[description]{t^n} \& \underset{d-1}{P_{\ol{\varepsilon}_2}} \ar[<-]{ru}[description]{t^n} \ar[densely dotted,-]{rr} \& \& \underset{d}{P_{\varepsilon_1}} \ar[->]{ru}[description]{t^n} \&
		{\color{gray} \underset{d-1}{P_{\ol{\varepsilon}_2}}} \ar[<-, color=gray]{ru}[description]{t^n}
	\end{gd}$}
$$
where the shaded arrows and projectives exist if and only if $m$ is even.
It follows that  $\CP(\omega) = 
(\begin{td}
	\underset{d}{P^{(m,\varepsilon_1)}}
	\ar{r}{M \,t^{n}}\&   \underset{d-1}{P^{(m,\varepsilon_2)}}
\end{td})
$ with the matrix
\begin{align*}
	M = 
	\Id + \smashoperator{\sum_{1 \leq i <  m}} E_{i,i+1} +
	\smashoperator{\sum_{1 \leq i < \frac{m}{2}}} E_{2i-1,2i+1} =
	\begin{psmallmatrix}
		1 & 1 & 1 &  &   \\
		& 1 & 1 & 0 &  \\
		&  & 1 & 1 & 1 &    \\
		&&  & 1 & 1 & 0 & \\
		&&&  & 1 & 1& \scalebox{0.5}{$\ddots$}  \\
		&&&& & 1 &  \scalebox{0.5}{$\ddots$} \\
		&&&&&  &  \scalebox{0.5}{$\ddots$}   \\
	\end{psmallmatrix}.
\end{align*}
\item Now assume that the $w$ begins with a left-oriented number.
In this case its gluing diagram can be obtained from the previous one by transposing all solid arrows.
It follows that
$\CP(\omega) = 
(\begin{td}
	\underset{d+1}{P^{(m,\varepsilon_2)}}
	\ar{r}{M^T t^{n}}\&   \underset{d}{P^{(m,\varepsilon_1)}}
\end{td})$, where $M^T$ denotes the transposition of the previous matrix $M$.
\end{itemize}

\subsubsection{A family of symmetric bands}
\label{subsec:sym-band}
Let
$\omega$ be the band $(0,w,m,\lambda)$
such that its word has periodic part
$\dot{w} = ( \ola{1},\ora{1},\ora{2},\ola{1},\ora{1},\ola{2})$, $m=1$ and $\lambda \in \kk^*$.
The diagram of $\dot{w}$ is given by 
\begin{align*}
\begin{gd}
\circ \ar[densely dotted,-]{d} 		\& \circ \ar[densely dotted,-]{d}   \& \circ \ar[densely dotted,-]{d}  \& \circ \ar[densely dotted,-]{d} 
\& \circ \ar[densely dotted,-]{d} 
\& \circ \ar[densely dotted,-]{d} 
\& \circ \ar[densely dotted,-]{d} 
\\
\circ  \ar[<-]{ru}[description]{1} \&
\circ \ar{ru}[description]{1} \&
\circ \ar{ru}[description]{2} \&
\circ  \ar[<-]{ru}[description]{1} \&
\circ  \ar{ru}[description]{1} \&
\circ  \ar[<-]{ru}[description]{2} \&\circ
\end{gd}
\end{align*} 
where the first and the last vertial edge are identified.
Since the diagram has point symmetry about the midpoint of the second vertical edge,
the periodic word $w$ is symmetric and we need to assume also that $\lambda \neq -1$.
By the orientation rules, it holds that $w^{\updownarrow} =
( \ola{1} \downarrow \ora{1} \uparrow \ora{2} \downarrow \ola{1} \downarrow \ora{1} \uparrow \ola{2} \downarrow )$
and thus the final gluing diagram has the form
\begin{align}
\label{eq:symm-bd-glu}
\begin{gd}
P_+ \ar[densely dotted]{d} 	\ar[<-]{rd} \ar[<-]{r}{\frac{1}{\lambda}t }	\& P_+ \ar[densely dotted]{d}   \& P_+ \ar[densely dotted,<-]{d} \ar{r}{-t^2} \& P_+  \ar[densely dotted]{d}   		\ar[<-]{r}{-t} 		\ar[<-]{rd}
\& P_+ \ar[densely dotted]{d} 
\& P_+ \ar[densely dotted,<-]{d} 
\& P_+ \ar[densely dotted]{d} 
\\
\underset{0,\lambda}{P_-}  \ar[<-]{ru}[description]{t} \ar[<-]{r}[swap]{t}\&
\underset{1}{P_-} \ar{ru}[description]{t} \ar{r}[swap]{t} \&
\underset{0}{P_-} \ar{ru}[description]{t^2} \&
\underset{-1}{P_-} \ar[<-]{ru}[description]{t} 
\ar[<-]{r}[swap]{-t}
\&
\underset{0}{P_-}  \ar{ru}[description]{t} 		 		 		
\ar{r}[swap]{-t}
\&
\underset{-1}{P_-} \ar[<-]{ru}[description]{t^2} 
\ar[<-]{r}[swap]{-{\frac{1}{\lambda}} t^2}
\&
\underset{0,\lambda}{P_-}
\end{gd}
\end{align}
Folding the diagram yields 
that the band complex $\CP(\omega)$ is given by
the three-term complex
$(\begin{td}
{P_{+-}} \ar{r}{\partial_1} \& {P_{+-}^3} \ar{r}{\partial_0} \& {P_{+-}^2}
\end{td})$
with the differentials
\begin{align*}
\partial_1=
\begin{pNiceArray}[first-row,last-col]{cc}
+ & - \\
\frac{1}{\lambda} t & t & +\\
t & t & - \\
0 & t & +\\
0 & t & - \\
0 & 0 & +\\
0 & 0 & -
\end{pNiceArray}
&&
\partial_0
=
\begin{pNiceArray}[first-row,last-col,columns-width=auto]{cccccc}
+ & - & + & - & + & - \\
-t^2 & t^2 & -t & t & 0 &0& +\\
0  & 0 & t & -t & 0 & 0&-\\
0  & 0 & 0 & t & 0 & 0&+\\
0  & 0 & 0 & -t &  t^2  & 	-\frac{1}{\lambda} t^2 &-\\	
\end{pNiceArray}
\end{align*}
Increasing the multiplicity from $1$ to an arbitrary number $m \in \N$
we need to enlarge $P_{+-}$ to $P_{+-}^{(m)}$, interpret zero entries as zero square matrices of size $m$,
view any entry $\pm t^n$ as a diagonal matrix $\pm t^n \Id$,
and replace $\pm \frac{1}{\lambda}$ by the Jordan block $\pm \JB_{\frac{1}{\lambda}}$.

\subsubsection{A family of asymmetric bands} \label{subsec:as-bd}

Let $\omega$ be the band $(0,w,m,\lambda)$ with any  $m \in \N$, $\lambda \in \kk^*$, and
$$
\dot{w} = (\ora{p}_1,\ola{q}_1,\ola{q}_2, \ola{q}_3, \ora{p}_2,\ora{p}_3, \ola{q}_4, \ora{p}_4)
$$
where all underlying numbers are arbitrary positive integers.
As the sequence of oriented arrows admits no axis of symmetry, the word $w$ is not symmetric.
As the numerical values in $w$ are kept arbitrary, only orientations between equioriented numbers are uniquely determined and we obtain
$$
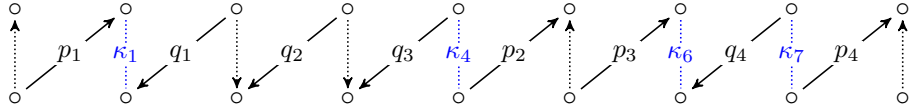
\begin{gd}
\bt \ar[densely dotted,<-]{d}	\& \bt \ar[densely dotted,-, color=blue]{d}[description]{\kappa_1}    \& \bt \ar[densely dotted,->]{d}  \& \bt \ar[densely dotted,->]{d} \& \bt  \ar[densely dotted,-, color=blue]{d}[description]{\kappa_4} \& \bt  \ar[densely dotted,<-]{d} 
\& \bt \ar[densely dotted,-, color=blue]{d}[description]{\kappa_6}    \& \bt \ar[densely dotted,-, color=blue]{d}[description]{\kappa_7}  \& \bt \ar[densely dotted,<-]{d} 	\\
\bt \ar{ru}[description]{{p_1}} 
\& \bt \ar[<-]{ru}[description]{{q_1}} \&  \bt \ar[<-]{ru}[description]{{q_2}}  \& \bt \ar[<-]{ru}[description]{{q_3}} 
\&
\bt \ar{ru}[description]{p_2} 
\& \bt \ar{ru}[description]{{p_3}} \&  \bt \ar[<-]{ru}[description]{q_4}  \& \bt \ar{ru}[description]{p_4} \& \bt
\end{gd}
$$
where $\kappa_1, \kappa_4, \kappa_6, \kappa_7 \in \{\downarrow,\uparrow\}$ may take on arbitrary orientations.
\begin{itemize}
\item $\kappa_6 = \downarrow$ if and only if $p_3 > q_4$, and similarly $\kappa_7 = \downarrow$ if and only if $q_4 < p_4$.
\item $\kappa_1 = \downarrow$ if and only if $p_1 > q_1$ or ($p_1 = q_1$ and $p_4 > q_2$) or $p_4 = p_1 = q_1 = q_2$.
\item $\kappa_4 = \downarrow$ if and only if $q_3 < p_2$ or ($q_3 = p_2$ and $q_2 < p_3$) or $q_2 = q_3 = p_2 = p_3$. 
\end{itemize} 

We define a pair of binary coefficients by $(\delta_i,\ol{\delta}_i) = (1,0)$ if $\kappa_i = \downarrow$ and $(\delta_i,\ol{\delta}_i) = (0,1)$ if $\kappa_i = \uparrow$ for each index $i \in \{1,4,6,7\}$.
An application of the source and sink rules yields the diagram
\begin{align*}
\begin{gd}
P_+ \ar[densely dotted,<-]{d} \ar{r}{- \lambda t^{p_1}}	\ar[dashed, color=blue]{rd}\& 
P_+ \ar[color=blue,densely dotted,<->]{d}  \ar[color=blue, dashed,<-]{r}{-\delta_1 t^{q_1}} \ar[<-,color=blue,dashed]{rd} \& 
P_+ \ar[densely dotted,->]{d}  \ar[<-]{r}{t^{q_2}} \ar[<-]{rd} \&
P_+ \ar[densely dotted,->]{d} \ar[<-]{r}{-t^{q_3}}	\ar[<-,color=blue,dashed]{rd} \& 
P_+  \ar[color=blue,densely dotted,<->]{d} \ar[dashed,color=blue]{r}{-\ol{\delta}_4 t^{p_2}}	\ar[color=blue,dashed]{rd} \&
P_+  \ar[densely dotted,<-]{d} \ar{r}{t^{p_3}}	\ar{rd} \& 
P_+ \ar[color=blue,densely dotted,<->]{d}   \ar[<-,color=blue,dashed]{r}{\delta_6 t^{q_4}} \ar[<-,color=blue,dashed]{rd} \& 
P_+ \ar[color=blue,densely dotted,<->]{d}  \ar[color=blue,dashed]{r}{\ol{\delta}_7 t^{p_4}} \ar[color=blue,dashed]{rd}
\& P_+ \ar[densely dotted,<-]{d} \\
\underset{0,\lambda}{P_-} \ar{ru}[description]{t^{p_1}} \ar[dashed, color=blue,swap]{r}{-\ol{\delta}_1 t^{p_1}}
\& \underset{-1}{P_-} \ar[<-]{ru}[description]{t^{q_1}} \ar[<-]{r}[swap]{-t^{q_1}} \&  \underset{0}{P_-} \ar[<-]{ru}[description]{t^{q_2}} \ar[<-]{r}[swap]{t^{q_2}}  \& \underset{1}{P_-} \ar[<-]{ru}[description]{t^{q_3}} \ar[<-,dashed,color=blue]{r}[swap]{-\delta_4 t^{q_3}}
\&
\underset{2}{P_-} \ar{ru}[description]{t^{p_2}} \ar{r}[swap]{-t^{p_2}}
\& \underset{1}{P_-} \ar{ru}[description]{t^{p_3}} \ar{r}[swap]{t^{p_3}} \&  \underset{0}{P_-} \ar[<-]{ru}[description]{t^{q_4}} 
\ar[<-,dashed,color=blue]{r}[swap]{\delta_7 t^{q_4}}
\& \underset{1}{P_-} \ar{ru}[description]{t^{p_4}} 
\ar{r}[swap]{ \frac{1}{\lambda} t^{q_4}}
\& \underset{0,\lambda}{P_-}
\end{gd}
\end{align*}
In the more compact notation we obtain a cycle diagram
\begin{align*}
\begin{gd}
\underset{0}{P_{+-}} \ar{r}{M_1 t^{p_1}} \&
\underset{-1}{P_{+-}} \ar[<-]{r}{M_2 t^{q_1}} \&
\underset{0}{P_{+-}} \ar[<-]{r}{M_3 t^{q_2}} \&
\underset{1}{P_{+-}} \ar[<-]{r}{M_4 t^{q_3}} \&
\underset{2}{P_{+-}} \ar{r}{M_5 t^{p_2}} \&
\underset{1}{P_{+-}} \ar{r}{M_6 t^{p_3}} \&
\underset{0}{P_{+-}} \ar[<-]{r}{M_7 t^{q_4}} \&
\underset{1}{P_{+-}} 
\ar[->,bend left=15]{lllllll}[inner sep=0.5pt]{ M_8   t^{p_4} }  
\end{gd}
\end{align*}
where the coefficient matrices are given as follows:
\begin{align*}
{
\setlength{\arraycolsep}{1pt}
\begin{array}{cccccccc}
	M_1 & M_2 &
	M_3 & M_4 &
	M_5 & M_6 &
	M_7 & M_8\\
	\begin{pNiceArray}[columns-width=15pt]{cc}
		-\lambda & 1 \\
		\ol{\delta}_1  & 			-	\ol{\delta}_1
	\end{pNiceArray}
	&
	\begin{pNiceArray}[columns-width=15pt]{cc}
		-	\delta_1 & \delta_1 \\
		1 & -1
	\end{pNiceArray}			
	&
	\begin{pNiceArray}[columns-width=15pt]{cc}
		1 & 1 \\
		1 & 1
	\end{pNiceArray}			
	&
	\begin{pNiceArray}[columns-width=15pt]{cc}
		-	1 & \delta_4  \\
		1 & -\delta_4 
	\end{pNiceArray}			
	&
	\begin{pNiceArray}[columns-width=15pt]{cc}
		-\ol{\delta}_4  & 1\\
		\ol{\delta}_4  & -1
	\end{pNiceArray}
	&
	\begin{pNiceArray}[columns-width=15pt]{cc}
		1 & 1 \\
		1 & 1
	\end{pNiceArray}			
	&
	\begin{pNiceArray}[columns-width=15pt]{cc}
		\delta_6  & \delta_6 \delta_7   \\
		1 & \delta_7  
	\end{pNiceArray}			
	&
	\begin{pNiceArray}[columns-width=15pt]{cc}
		\ol{\delta}_7  & 1 \\
		\ol{\delta}_7  & \frac{1}{\lambda}
	\end{pNiceArray}			
\end{array}
}
\end{align*}
In the following we abbreviate $\tilde{M}_1 = M_1 t^{p_1}$, $\tilde{M}_2 = M_2 t^{q_2}$ and so on.
An intermediate step for folding is given by placing projectives with shared degrees in the same columns.
$$
\begin{gd}
P_{+-} \ar{r}{\tilde{M}_4} \ar{rd}[description]{\tilde{M}_5} \& P_{+-} \ar{rd}[description]{\tilde{M}_3} \& P_{+-} \ar{r}{\tilde{M}_1} \& P_{+-}\\
\& P_{+-} \ar{rd}[description]{\tilde{M}_6} \& P_{+-} \ar{ru}[description]{\tilde{M}_2} \\
\& P_{+-} \ar{r}[description]{\tilde{M}_7} \ar[color=blue,densely dotted]{ruu}[description]{\tilde{M}_8}  \& P_{+-} 
\end{gd}
\begin{array}{c}
\Rightarrow 
\begin{cd}
\underset{2}{P_{+-}} \ar{r}{
	\partial_2  
} \& \underset{1}{P_{+-}^3} \ar{r}{\partial_1}
\& \underset{0}{P_{+-}^3} \ar{r}{
	\partial_0
} \& \underset{-1}{P_{+-}}
\end{cd}
\\[10pt]
\text{with }
\partial_2  = 
\begin{pmatrix}
\tilde{M}_4 \\
\tilde{M}_5 \\
0
\end{pmatrix},
\quad
\partial_1=
\begin{pmatrix}
0 & 0 & \tilde{M}_8 \\
\tilde{M}_3 & 0 & 0 \\
0 & \tilde{M}_6 & \tilde{M}_7
\end{pmatrix},
\\[10pt]
\partial_0 =
\begin{pmatrix}
\tilde{M}_1 & \tilde{M}_2 & 0
\end{pmatrix}.
\end{array}
$$
At last, we obtain the band complex $\CP(\omega)$ on the right.
This example indicates that band complexes can be described more generally for periodic words with arbitrary numerical values
using coefficient matrices which reflect all possible orientations of gluing edges. A similar observation holds true for string complexes.

\subsection{Main gluing result}
By the last subsection, any band and any string provides a complex in $\Db{\Acat}$.
However, some bands or strings may yield  complexes that are isomorphic in the derived category. This motivates the next notion.
\begin{dfn}\label{dfn:equiv}
We consider the smallest equivalence relation on the set of bands and strings of $\chfXA$ such that the following holds.
\begin{enumerate}
\item  Any usual  string $(d,w)$ is equivalent to $(d^{\op},w^{\op})$.
\item  Any special  string $(d,w,\varepsilon)$ is equivalent to $(d^{\op},w^{\op},\varepsilon)$.
\item Any bispecial string 
$(d,\varepsilon_1,w,\varepsilon_2)$ 
is equivalent to $(d^{\op},\varepsilon_2,(v^{\op})^m,\varepsilon_1)$,
where $v$ denotes the primitive root of $w$ and $m$ its multiplicity.
\item \label{dfn:equiv3} Any band $(d,w,m,\lambda)$ is equivalent to 
$(d_1,w_1, m, \lambda)$
as well as
$(d,w^{\op}, m, \lambda^{\mp})$ where
$\lambda^{\mp} = \frac{1}{\lambda}$ if $w$ is asymmetric, and $\lambda^{\mp} = \lambda$ otherwise.
\end{enumerate}
\end{dfn}
In particular, any string is not equivalent to any band or a string of another of the three classes.
The main result of this section is the following.
\begin{thm} \label{thm:bijection2}
The assignment $\begin{td} \omega \ar[mapsto]{r} \& \CP(\omega) \end{td}$ 
with $\CP(\omega)$ defined by 
Subsection~\ref{sec:gluing}
yields a bijection
\begin{align}
\label{eq:bij-perf-fd2}
\begin{cd}
	\{\text{bands and strings of }\chfXA  \}\big/_{\textstyle \approx}
	\ar[yshift=0pt]{r}{\sim}\&
	\ind \Db{\Acat} 
\end{cd}
\end{align}
from
the set of equivalence classes of bands and strings of $\chfXA$ to the set of isomorphism classes of indecomposable 
objects in the category $\Db{\Acat}$.
\end{thm}

\begin{rmk}
There exists a bijection of the form \eqref{eq:bij-perf-fd2} if and only if the following statements hold.
\begin{enumerate}	
\item For any band or string $\omega$ of $\DbRep{\A}$, the complex $\CP(\omega)$ defines an indecomposable object in $\DbRep{\A}$.
\item Any two bands or strings $\omega$ and $\omega'$  are equivalent if and only if there is an isomorphism $\CP(\omega) \cong \CP(\omega')$ in 
$\DbRep{\A}$.
\item For any indecomposable complex $\CP$ in $\DbRep{\A}$ there is a band or string $\omega$ such that there is an isomorphism $\CP \cong \CP(\omega)$
in  $\DbRep{\A}$.
\end{enumerate}
\end{rmk}

For technical reasons, we need to single out a series of strings.
\begin{dfn}\label{dfn:trivial}
Any string $(d,w)$ with $d\in\Z$
and $w = (\star,x_1,\star)$ will be called \emph{trivial}.
\end{dfn}

\begin{rmk}
The bijection in \eqref{eq:bij-perf-fd2} restricts to a bijection
\begin{align*}
		\begin{cd}
			\{\text{bands and bispecial strings of }\chfXA  \}\big/_{\textstyle \approx}
			\ar[yshift=0pt]{r}{\sim}\&
			\ind \Hotbfd\bigl(\add (P_{+-})\bigr) \end{cd}.
\end{align*}
In particular, the subcategory 
$\Hotbfd\bigl(\add (P_{+-})\bigr)$
constitutes a large part of $\DbRep{\A}$, whose Auslander-Reiten quiver is given by homogeneous tubes of ranks one and two.\\
The bijection in \eqref{eq:bij-perf-fd2} restricts to another bijection
$$
\begin{cd}
	\{\text{trivial strings of }\chfXA  \}\big/_{\textstyle \approx}
	\ar[yshift=0pt]{r}{\sim}\&
	\ind \Hotbfd(\add P_\star)
\end{cd}.
$$
However,  the subcategory $\Hotbfd(\add P_\star)$
describes only a small part of $\DbRep{\A}$.
\end{rmk}

\begin{rmk}\label{rmk:field}
In case of an arbitrary field $\kk$, Theorem~\ref{thm:bijection2} is also true after 
the following changes in the definitions involving bands.
\begin{enumerate}
\item 
In Definition~\ref{dfn:band}, the parameter $\lambda$ in a band $(d,w,m,\lambda)$ has to be replaced by a monic irreducible polynomial $p \in \kk [x]$ such that $p(x) \neq x$.
If $w$ is symmetric,
it is also required that $p(x) \neq x + 1$.
\item 
Definition~\ref{dfn:equiv}~\eqref{dfn:equiv3} should be changed to the statement that
any band $(d,w,m,p)$ is equivalent to $(d_1,w_1, m, p)$ as well as $(d,w^{\op},m, p_{\mp})$ with $p_-(x) \colonequals  
\frac{1}{p(0)} \, x^{\deg p}\,  p(\frac{1}{x})$
if $w$ is asymmetric and $p_+(x) = p(x)$ otherwise.
\item \label{rmk:field-mat} Instead of Jordan blocks $\JB_{\lambda}$
and $\JB_{\lambda^{-1}}$ of size $m$
in the construction of a band complex
one may use any square matrix of size $m$ with minimal polynomial $p^m$ respectively $p_-^m$
such as the Frobenius block with the respective minimal polynomial.
\end{enumerate}
\end{rmk}

\section{Proof of the Gluing Theorem}\label{sec:proofs}
The goal of this section is to prove Theorem~\ref{thm:bijection2} on the correspondence of bands and strings to indecomposable complexes. 
The proof has three steps.
\begin{enumerate}
	\item
	First, we observe that any band or non-trivial string $\omega$ of $\fXA$ 
	in the sense of Subsection~\ref{sec:gluing}
	admits an expansion to
	a band or string of $\fXA$ as defined in Subsection~\ref{subsec:abc}.
\end{enumerate}
So far, the previous results combine to a diagram of bijections of sets
\begin{align}\label{eq:bij-diag}
	\begin{cd} 
		\{\text{bands and strings of }\chfXA  \}\!\big/_{\textstyle \approx} 		 \arrow[d, "" swap, "\rotatebox{90}{$\sim$}"]  \ar[dashed]{r}{\sim} \&[0.25cm] \ind\DbRep{\A} \arrow[d,<-, "\text{Cor.~\ref{cor:reduction}}" swap, "\rotatebox{90}{$\sim$}"] 
		\\
		\{\text{bands and strings of }{\fXA}  \}\!\big/_{\textstyle \approx} \cup \mathcal{I}_\star
		\ar{r}{\sim}[swap]{\text{Thm.~\ref{thm:Bondarenko}}} \& 	\ind \regrep{\fXA}  \cup \mathcal{I}_\star
	\end{cd}
\end{align}
where $\mathcal{I}_\star = \ind \Hotbfd(\add P_\star)$,
the left vertical bijection follows from the first step, the bottom horizontal bijection is a consequence of Bondarenko's classification result,
and the right vertical bijection is the main result of the previous reduction via the category of gluing triples.
For any band or string $\omega$, we define the `pullback complex' $\CX(\omega)$ to be the image under the top horizontal bijection in the diagram. We prove that
the complex $\CX(\omega)$ is isomorphic to the band or string complex $\CP(\omega)$ as defined by the combinatorial rules in Subsection~\ref{sec:gluing} as follows.
\begin{enumerate} \setcounter{enumi}{1}
	\item We define a complex $\CP'(\omega)$ by changing some signs in the differentials of the complex $\CP(\omega)$, and show that
	$\CP'(\omega) \cong \CP(\omega)$.
	\item 
	In the essential step of the proof,  we show that  $\CP'(\omega) \cong \CX(\omega)$.
	The most technical part is to
	describe the differentials of the complex $\CX(\omega)$ explicitly, which makes use of the construction of Subsection~\ref{subsec:glu}.
\end{enumerate}
The three steps correspond to the three subsections below.

\subsection{Expansion of non-trivial strings and bands}
\label{subsec:exp}
We describe the expansion of non-trivial strings of $\chfXA$,
and then that of bands separately.
We recall that trivial strings have been introduced in Definition~\ref{dfn:trivial}.

\subsubsection{Expanded notation for non-trivial strings}
\label{subsec:exp-str}
Let $(d,w)$ be a pair given by a degree $d \in \Z$ and a finite word $w$ of $\chfXA$, that is, 
$w = (\alpha, x_1, x_2, \ldots, x_{\ell}, \beta)$ with ends $\alpha, \beta \in \{\star,\diamond\}$
and oriented numbers $x_1$, $x_2$ $\ldots$, $x_{\ell} \in \kX$. We assume that $w \neq (\star,x_1,\star)$, which ensures that $w$ does not appear in a trivial string.

We define an \emph{expanded word} $(d,w)_e$ of the bunch of skew-gentle type $\fXA$ from Definition~\ref{def:chkX} in the following.

Let $(d_0,d_1,\ldots, d_{\ell})$ denote the sequence of integers
associated to the pair $(d,w)$ from Definition~\ref{dfn:deg}, that is,
$d_0 = d$ and $d_i = d_{i-1} + 1$ if $x_i
= {\ola{n}_i}$ respectively $d_i = d_{i-1} -1$ if $x_i = {\ora{n}_i}$ for a number $n_i \in \N$.
We interlace the word $w$ with that sequence of integers
setting
$w' = (\alpha, d_0, x_1, d_1, x_2, d_2 \ldots d_{\ell-1}, x_{\ell}, d_{\ell}, \beta)$,
and view this expression as a formal composition
\begin{align}
	\label{eq:comp}
	(\alpha, d_0, x_1, d_1) \circ (d_1, x_2, d_2) \circ \ldots \circ (d_{\ell-2}, x_{\ell-1}, d_{\ell-1} ) 
	\circ (d_{\ell-1}, x_{\ell}, d_{\ell}, \beta)
\end{align}
The first factor in this expression is replaced by an expression $w_e^{(1)}$ defined by the next table, in which we have included gluing diagrams for the sake of completeness.
					\begin{longtable}{|c|c|c|}
						\hline
						$\text{gluing diagram}$ & $(\alpha, d_{0}, x_1, d_1)$ & $w_e^{(1)}$ \\ \hline\hline
						\endfirsthead
						\hline
						$\text{gluing diagram}$ & $(\alpha, d_{0}, x_1, d_1)$ & $w_e^{(1)}$ \\ \hline\hline
						\endhead
						\hline
						\endfoot
						$\begin{gd}
							\& \bt \ar[densely dotted,-]{d} \\
							\bt \ar{ru}[description]{n_1}\& \bt
						\end{gd}$
						&
						$(\diamond, d_{0}, {\ora{n}_1}, d_1)$
						&
						$f_{d_{0}} - a_{d_{0}}^{2 n_1}  \sim b_{d_1}^{2 n_1} - f_{d_1}$
						\\ \hline
						$\begin{gd}
							\& \bt \ar[densely dotted,-]{d} \\
							\bt \ar[<-]{ru}[description]{n_1}\& \bt
						\end{gd}$
						&
						$(\diamond, d_{0}, {\ola{n}_1}, d_1)$
						&
						$f_{d_{0}} - b_{d_{0}}^{2 n_1}  \sim a_{d_1}^{2 n_1} - f_{d_1}$
						\\ \hline
						$\begin{gd}
							\& \bt \ar[densely dotted,-]{d} \\
							\star \ar{ru}[description]{n_1}\& \bt
						\end{gd}$
						&
						$(\star, d_{0}, {\ora{n}_1}, d_1)$
						&
						$b_{d_{1}}^{2 n_1-1}  - f_{d_1}$
						\\ \hline
						$\begin{gd}
							\& \bt \ar[densely dotted,-]{d} \\
							\star \ar[<-]{ru}[description]{n_1}\& \bt
						\end{gd}$
						&
						$(\star, d_{0}, {\ola{n}_1}, d_1)$
						&
						$a_{d_{1}}^{2 n_1+1}  - f_{d_1}$
						\\ \hline
					\end{longtable}
					\addtocounter{equation}{-1}

					For each index $2 \leq i < \ell$ we expand the $i$-th factor in \eqref{eq:comp} as follows.

					\begin{align} \label{eq:main-expand}
						\begin{array}{|c|c|c|}
							\hline
							\text{gluing diagram}&
							(d_{i-1}, x_i, d_i) & w_e^{(i)}   \\ \hline \hline
							\begin{gd}
								\bt \ar[densely dotted,-]{d} \& \bt \ar[densely dotted,-]{d} \\
								\bt \ar{ru}[description]{n_i}\& \bt
							\end{gd}
							&
							(d_{i-1}, {\ora{n}_i}, d_i) & f_{d_{i-1}} - a_{d_{i-1}}^{2 n_i}  \sim b_{d_i}^{2 n_i} - f_{d_i}
							\\ \hline
							\begin{gd}
								\bt \ar[densely dotted,-]{d} \& \bt \ar[densely dotted,-]{d} \\
								\bt \ar[<-]{ru}[description]{n_i}\& \bt
							\end{gd}		&
							(d_{i-1}, {\ola{n}_i}, d_i) & f_{d_{i-1}} - b_{d_{i-1}}^{2 n_i}  \sim a_{d_i}^{2 n_i} - f_{d_i}  
							\\ \hline
						\end{array}
					\end{align}

					The last factor in \eqref{eq:comp} is expanded as follows.

					\begin{longtable}{|c|c|c|}
						\hline
						$\text{gluing diagram}$ & $(d_{\ell-1}, x_\ell, d_\ell, \beta)$ & $w^{(\ell)}_e$ \\ \hline\hline
						\endfirsthead
						\hline
						$\text{gluing diagram}$ & $(d_{\ell-1}, x_\ell, d_\ell, \beta)$ & $w^{(\ell)}_e$ \\ \hline\hline
						\endhead
						\hline
						\endfoot
						$\begin{gd}
							\bt \ar[densely dotted,-]{d} \& \bt \\
							\bt \ar{ru}[description]{n_\ell} \&
						\end{gd}$
						&
						$(d_{\ell-1}, {\ora{n}_\ell}, d_\ell, \diamond)$
						&
						$f_{d_{\ell-1}} - a_{d_{\ell-1}}^{2 n_\ell}  \sim b_{d_\ell}^{2 n_\ell} - f_{d_\ell}$
						\\ \hline
						$\begin{gd}
							\bt \ar[densely dotted,-]{d} \& \bt \\
							\bt \ar[<-]{ru}[description]{n_\ell} \&
						\end{gd}$
						&
						$(d_{\ell-1}, {\ola{n}_\ell}, d_\ell, \diamond)$
						&
						$f_{d_{\ell-1}} - b_{d_{\ell-1}}^{2 n_\ell}  \sim a_{d_\ell}^{2 n_\ell} - f_{d_\ell}$
						\\ \hline
						$\begin{gd}
							\bt \ar[densely dotted,-]{d} \& \star \\
							\bt \ar{ru}[description]{n_\ell} \&
						\end{gd}$
						&
						$(d_{\ell-1}, {\ora{n}_\ell}, d_\ell, \star)$
						&
						$f_{d_{\ell-1}} - a_{d_{\ell-1}}^{2 n_\ell+1}$
						\\ \hline
						$\begin{gd}
							\bt \ar[densely dotted,-]{d} \& \star \\
							\bt \ar[<-]{ru}[description]{n_\ell} \&
						\end{gd}$
						&
						$(d_{\ell-1}, {\ola{n}_\ell}, d_\ell, \star)$
						&
						$f_{d_{\ell-1}} - b_{d_{\ell-1}}^{2 n_\ell-1}$
						\\ \hline
					\end{longtable}
					\addtocounter{equation}{-1}
					The expanded word of $(d,w)$ is set to
					$(d,w)_e = w^{(1)}_e \sim w^{(2)}_e \sim \ldots 
					\sim w^{(\ell)}_e$.

					\begin{lem}\label{lem:str-bij}
						There is a bijection between
						the
						set of equivalence classes of non-trivial strings of $\chfXA$ and  the set of equivalence classes of strings	of $\fXA$
						$$
						\begin{cd}
							\{ \text{non-trivial strings of }\chfXA  \}\!\big/_{\textstyle \approx} \ar{r}{\sim} \& 
							\{ \text{strings of }{\fXA}  \}\!\big/_{\textstyle \approx} \& \omega \ar[mapsto]{r} \& \omega_e
						\end{cd}
						$$
						with $\omega_e$ defined as follows.
						$$
						\begin{array}{|c||c|c|c|}
							\hline
							\text{non-trivial string} & \text{usual string} & \text{special string} & \text{bispecial string}  \\
							\omega & (d,w) & (d,w,\varepsilon) & (d,\varepsilon_1,w ,\varepsilon_2)  \\ \hline
							\omega_e & 
							(d,w)_e & ((d,w)_e, \varepsilon) &
							(\varepsilon_1,(d,w)_e, \varepsilon_2) \\
							\hline
						\end{array}
						$$
					\end{lem}
					\begin{proof}
						With the definitions above
						and Definitions~\ref{dfn:fin-word}, \ref{dfn:fin-word2} and \ref{def:equiv}
						for $\fXA$, the claim follows from the following observations.
						\begin{enumerate}
							\item The assignment $(d,w) \mapsto (d,w)_e$ yields a bijection between 
							the set of pairs $(d,w)$ with $d \in \Z$ and $w$ a finite word  of $\chfXA$ and the set of finite 
							words of $\fXA$, introduced in Definition~\ref{dfn:fin-word}.		
							\item The word $w$ has as many special ends as $(d,w)_e$.
							\item 
							We set $(d,w)^{\op} = (d^{\op}, w^{\op})$. 
							It can be checked that
							$((d,w)^{\op})_e = ((d,w)_e)^{\op}$. In particular, $w$ is symmetric if and only if $(d,w)_e$ is symmetric. \qedhere
						\end{enumerate}
					\end{proof}

					\subsubsection{Expanded notation for bands}
					\label{subsec:exp-bnd}

					Let $(d,w)$ be a pair given by a degree $d \in \Z$ 
					and periodic word 	 $w = (x_i)_{i \in \Z}$ of $\chfXA$.
					We define an expanded word $(d,w)_e$ of $\fXA$ by a similar procedure as for finite words as follows.
					
					Let $\ell \in 2\N$ denote the period of $w$,
					and $(d_i)_{i \in \Z}$ the degree sequence with $d_0 = d$, $d_{i} = d_{i-1} \pm 1$ with $-1$ if $x_i = {\ola{n}_i}$ and $+1$ if $x_i = {\ora{n}_i}$  for each $i \in \Z$.
					Interlacing this integer sequence with $w$ and `decomposing' yields
					a periodic expression
					\begin{align*}
						{\color{blue}\ldots (d_{-1}, x_0, d_0) \ \circ \ } 
						(d_0, x_1, d_1) \circ (d_1, x_2, d_2) \circ \ldots  
						\circ \ (d_{\ell-1}, x_{\ell}, d_{\ell}) \circ {\color{blue} (d_{\ell}, x_{\ell+1}, d_{\ell+1}) \circ \ldots }.
					\end{align*}
					For each $i \in \Z$ each factor is translated in the same way as in~\eqref{eq:main-expand}. 
					As previously, we
					set $(d,w)_e = \ldots w^{(0)}_e \sim w^{(1)}_e \sim w^{(2)}_e \sim \ldots \sim w^{(\ell-1)}_e \sim w^{(\ell)}_e \ldots \ $, which
					defines a periodic word $(a_i r_i)_{i \in \Z}$ by identifying
					$r_1$ with the unique relation $\sim$ in $w^{(1)}_e$.
					
					When passing from bands of $\chfXA$ to bands of $\fXA$, we need to rescale the continuous parameters in certain cases.
					For a periodic word $w$ of $\chfXA$ and 
					$\lambda \in \kk^*$ we set $s_w
					= (-1)^{\frac{\ell}{2}} $ if $w$ is symmetric and has period $\ell$, respectively $s_w = 1$ otherwise.
					\begin{lem}\label{lem:bnd-bij}
						The
						assignment 
						$\begin{td}
							\omega =  (d,w,m,\lambda) \ar[mapsto]{r} \& \omega_e =((d,w)_e,m,s_w \lambda) \end{td}$
						defines a bijection on the sets of equivalence classes
						$$
						\begin{cd}
							\{ \text{bands of }\chfXA  \}\!\big/_{\textstyle \approx} \ar{r}{\sim} \& 
							\{ \text{bands of }{\fXA}  \}\!\big/_{\textstyle \approx}
						\end{cd} 
						$$
						which restricts to a bijection on equivalence classes of symmetric bands.
					\end{lem}
					\begin{proof}
						As in the proof of the previous lemma, we make a few observations.
						\begin{enumerate}
							\item For a pair $(d,w)$ with a word $w$ of $\chfXA$ of period $\ell$, the word $(d,w)_e$ has the form 
							\begin{align*}
								\ldots
								f_{d_0} - 
								\underbrace{e'_1 \sim e''_1}_{a_1 r_1 a_2} - f_{d_1} \sim
								f_{d_1} - e'_2 \sim e''_2 - f_{d_2} \sim
								\ldots
								f_{d_{\ell-1}} - e'_{\ell} \sim e''_{\ell} - f_{d_\ell} {\color{blue} \sim \underbrace{f_{d_\ell} -}_{a_{4\ell} r_{4\ell}} \ldots}
							\end{align*}
							with certain $e'_i, e''_i \in \fE$ for each $i \in \Z$. The period of $(d,w)_e$ divides $4\ell$.
							\item 
							The assignment $(d,w) \mapsto (d,w)_e$ yields a bijection between 
							the set of pairs $(d,w)$ with $d \in \Z$ and $w$ a periodic word of $\chfXA$ and the
							set $\mathcal{P}$ of periodic words $(a_i r_i)_{i \in \Z}$ of $\fXA$ with $a_1 \in \fE$.
							\item For any $j \in \Z$ it holds that $(d_j,w_j)_e = ((d,w)_e)_{4j}$.
							Vice versa,
							for any periodic word $v$ from $\mathcal{P}$ 
							and 
							any $p \in \Z$ with $v_{p} = v$ it follows that $p \in 4 \Z$. 
							\item The last fact implies that the period of $(d,w)_e$ is precisely $4 \ell$.
							Thus a periodic part of $(d,w)_e$ has $t((d,w)_e) = \frac{\ell}{2}$ tied subwords.
							\item \label{sym-op}
							It holds that 
							$((d,w)_e)^{\op} = ((d,w^{\op})_e)_2$.  		In particular,	the word $w$ is symmetric if and only if $(d,w)_e$ is symmetric.
							\item For a band $\omega = (d,w,m,\lambda)$ of $\chfXA$ and its expansion $\omega_e = ((d,w)_e,m,\lambda)$ the following holds.
							\begin{enumerate}
								\item If $w$ is non-symmetric, neither is $(d,w)_e$. So $\omega_e$ is a non-symmetric band.
								\item \label{sym-step} If $w$ is symmetric of period $\ell$, then $\lambda \neq -1$ and
								$s_w \lambda =(-1)^{\frac{\ell}{2}}\lambda \neq (-1)^{\frac{\ell}{2}+1} \equalscolon \ol{s}_w$.
								As $t((d,w)_e) = \frac{\ell}{2}$, the parameter $\lambda$ in the symmetric band $\omega_e$ is admissible.
								Moreover, assigning $\begin{td} \lambda \ar[mapsto]{r}\&s_w\lambda \end{td}$
								defines a bijection $\begin{td} \kk^*\backslash\{-1\} \ar{r}{\sim} \& \kk^*\backslash \{ \ol{s}_w \} \end{td}$.
							\end{enumerate}
							In both cases it follows that $\omega_e$ is a well-defined band of $\fXA$.
							\item 
							Expanding the equivalent bands
							$\omega=(d,w,m,\lambda)$ and $(d_1,w_1,m,\lambda)$
							yields equivalent bands 
							$\omega_e = ((d,w)_e,m,s_w \lambda)$ and $(((d,w)_e)_4,m,s_w \lambda)$.
							The band $\omega$ is also equivalent to the band $\omega^{\op} \colonequals (d,w^{\op},m,\lambda^{\mp})$.
							Then 
							$(\omega^{\op})_e = ((d,w^{\op})_e,m,s_w \lambda^{\mp})$ is equivalent to $(((d,w)_e)^{\op}, m, s_w \lambda^{\pm})$ by \eqref{sym-op}, and thus to
							$\omega_e$ as well.
							It follows that any two bands $\omega$ and $\omega'$ are equivalent if and only if $\omega_e$ and $\omega'_e$ are equivalent.
							\item Finally, given any band of $\fXA$ we can pass to an equivalent band such that its underlying word lies in $\mathcal{P}$. 
							Together with the
							previous properties it follows that 
							the assignment $\begin{td}\omega \ar[mapsto]{r} \& \omega_e \end{td}$ is surjective.
						\end{enumerate}
						The last arguments complete the proof that the assignment in question is a well-defined bijection.
						The bijections restricts to sets of equivalence of symmetric bands because of~\eqref{sym-step}.
					\end{proof}

					\subsubsection{Main statement on expansion of bands and strings}
					\begin{prp} 		\label{prp:bijection}
						Mapping any band or non-trivial string $\omega$ of $\chfXA$ to its expansion $\omega_e$ 
						and any trivial string $\omega$ of $\chfXA$ to the complex $\CX(\omega)$ 
						yields
						a well-defined bijection 
						\begin{align*}
							\begin{cd} 
								\{\text{bands and strings of }\chfXA  \}\!\big/_{\textstyle \approx} \ar{r}{\sim} \& 
								\{\text{bands and strings of }{\fXA}  \}\!\big/_{\textstyle \approx}\cup
								\mathcal{I}_\star
							\end{cd}
					\end{align*}
					where $\mathcal{I}_\star =	\ind\Hotbfd(\add P_\star)$.
				\end{prp}
				\begin{proof}
					As strings cannot be equivalent to bands neither in compact nor in extended notation, the bijection on equivalence classes on strings in Lemma~\ref{lem:str-bij}
					and the bijection on equivalence classes of bands in Lemma~\ref{lem:bnd-bij} provide a bijection 
					\begin{align*}
						\begin{cd}
							\{ \text{bands and non-trivial strings of }\chfXA  \}\!\big/_{\textstyle \approx} \ar{r}{\sim} \& 
							\{ \text{bands and strings of }{\fXA}  \}\!\big/_{\textstyle \approx}
						\end{cd}
					\end{align*}
					Since $\End_\A(P_{\star}) \cong \Rx = \kk \llbracket t\rrbracket$, 
					there is also a bijection
					\begin{align*}
							\begin{cd}
								\{ \text{trivial strings of }\chfXA  \}\!\big/_{\textstyle \approx} 
								\ar{r}{\sim} \& 
								\mathcal{I}_\star
							\end{cd} 
							&&
							\begin{td}
								\bigl(d, (\star,\olra{n},\star)\bigr) \ar[mapsto]{r} 
								\& 
								({P}_{\star} \ar{r}{t^n} \& {P}_{\star})[d] \end{td} 
						\end{align*}	
							A combination of these bijections yields the claim.
						\end{proof}

						\subsubsection{Translation of bands and strings with oriented words}
						\label{subsec:od-mat}
						Subsection~\ref{subsec:exp-str} has provided a translation between pairs $(d,w)$ given by an integer $d \in \Z$ and words $w$ of $\chfXA$ and gluing diagrams.
						We need to extend this translation to  words
						with oriented arrows and gluing diagrams with oriented edges.
						
						By the translation in \ref{subsec:exp-str},  a pair $(d,w)$ of an integer $d \in \Z$ and  a finite word $w$ of length $\ell$
						gives rise to a sequence
						$w' = (\alpha,d_0,x_1,d_1,x_2,d_2, \ldots x_{\ell},d_{\ell},\beta)$
						such that $\|d_{i} - d_{i-1}\| =1$ for each $1 \leq i \leq \ell$.
						Each degree $d_i$ corresponds 
						to a vertical edge in the gluing diagram of $w$ on the one hand, and special subwords $f_d \sim f_d$ of the expanded word $(d,w)_e$ on the other hand.
						Going further, the word $w$ gives rise 
						to a sequence $(x_i \kappa_i)_{i \in \Z}$ with additional parameters $\kappa_i \in \{\uparrow, \downarrow\}$ defined in \eqref{eq:updownarrow}, while in Subsection~\ref{subsec:A-orient}
						each special subword is decorated by a left or right arrow.
						The left and right oriented arrows above special subwords correspond to downward and upward arrows in the gluing diagram, irrespective of the surrounding elementary gluing diagrams.
						\begin{align} \label{eq:ornt-exp}
							\begin{array}{|c|c|c|}
								\hline
								\text{gluing diagram}&
								\kappa_i & (d,w)_e   \\ \hline \hline
								\begin{gd}
									\bt \ar[densely dotted,->]{d}  \\
									\bt
								\end{gd}
								&
								\downarrow 
								&
								\overleftarrow{f_{d_{i}} \sim f_{d_i}} 
								\\ \hline
							\end{array}
							&&
							\begin{array}{|c|c|c|}
								\hline
								\text{gluing diagram}&
								\kappa_i & (d,w)_e   \\ \hline 
								\hline
								\begin{gd}
									\bt \ar[densely dotted,<-]{d}  \\
									\bt
								\end{gd}
								&
								\uparrow 
								&
								\overrightarrow{f_{d_{i}} \sim f_{d_i}} \\ \hline
							\end{array}
						\end{align}
						The same translations are carried out in case of pair $(d,w)$ with a periodic word $w$ of $\chfXA$.
						This yields a datum denoted by $((d,w)^\kappa)_e$.

						The translation of oriented gluing diagrams is compatible with expansion of bands and strings in the following sense.
						\begin{lem}\label{lem:ornt-ext}
							Let $(d,w)$ be the pair underlying a band or string $\omega$ of $\chfXA$.
							Then 
							$({(d,w)}^\kappa)_e$ coincides with the oriented version of  the word $(d,w)_e$.
						\end{lem}
						\begin{proof}
							In the notations above, let	$\ol{w}$ be the ambient word 
							of $w$ 
							and $\ol{v}$ the ambient word of $v = (d,w)_e$ and $\olra{v}$ its oriented version.
							\begin{itemize}
								\item The word $\ol{w}$ is periodic and symmetric with $\ol{w}_i = (\ol{w}_i)^{\op}$ if and only if  the word $\ol{v}$ is periodic and  symmetric with 
								$\ol{v}_{4i} = (\ol{v}_{4i} )^{\op}$.
								In this case, $\kappa_i$ is given by $\downarrow$, and matches with the left-oriented arrow above the $i$-th special subword in $\olra{v}$
								via the translation rule \eqref{eq:ornt-exp}.
								\item Assume that the previous case above does not occur at $\kappa_i$.
								\begin{itemize}
									\item Assume that $\kappa_i = \downarrow$.
									Then one of the following cases occurs.
									\begin{itemize}
										\item If
										$(a_i, b_i) = (\ora{p},\ora{q})$, then $x_i = \xx_d^{2p} < \yy_d^{2q} = y_i$,
										\item 		If $(a_i, b_i) = (\ola{p},\ora{q})$ with $p < q$
										then $x_i = \xx_d^{2p} < \xx^{2q}_d = y_i$,
										\item	If $(a_i,b_i) = (\ora{p},\ola{q})$ with $p > q$, then $x_i = \yy_d^{2p} < \yy_d^{2q} = y_i$, 
									\end{itemize}
									where throughout $p,q \in \{ \frac{m}{2} \mid m \in \N\}$ and $d \in \Z$.
									In each case, it follows that $x_i < y_i$, or, equivalently, the arrow above the $i$-th special subword $\xf_{d_i} \sim \xf_{d_i}$ is left-oriented. 
									\item 
									Similarly, it is straightforward to show that $\kappa_i = \uparrow$ leads to three cases each of which implies that $x_i > y_i$. Thus the $i$-th special subword of $\olra{v}$ is right-oriented.
								\end{itemize}
							\end{itemize} 
							This completes the proof of  the claim.
						\end{proof}

						\begin{rmk}\label{rmk:sgn-ext}
							Let $\omega$ be a bispecial string $(d,\varepsilon_1,w, \varepsilon_2)$ of $\chfXA$. Let 
							$v$ denote the primitive root of $w$ and $m \in \N$  its multiplicity.
							Then the bispecial word $(d,w)_e$ of $\fXA$
							has primitive root $(d,v)_e$ and multiplicity $m$ as well.
						\end{rmk}

						The additional translations will be used to
						translate a band or string $\omega$ of $\chfXA$ into an explicit matrix representation
						$\mu(\omega_e)$ of $\regrep{\fXA}$.

						\subsection{Band or string complexes with a degree-independent sign convention}
						For any band or string $\omega$ be a band or string $\omega$ of $\chfXA$, we define $\CP'(\omega)$ as follows.
						\begin{enumerate}
							\item Starting with the gluing diagram of the underlying word of $\omega$, we carry out the vertex replacement 
							as in Subsection~\ref{subsec:polar}
							for the complex $\CP(\omega)$.
							\item Similarly, 
							the gluing edges are oriented using Subsection \ref{subsec:orient} as before.
							\begin{itemize}
								\item If $\omega$ is a band $(d,w,m,\lambda)$ with a word of period $\ell$ we use $\lambda' = (-1)^{\frac{\ell}{2}} \lambda$ as decorating eigenvalue below the projective at degree zero.
							\end{itemize}
							\item In the third step, we use different signs on induced arrows
							than those in Subsection~\ref{subsec:ind-arr}.
							\begin{itemize}
								\item Whenever we apply a source rule, the induced arrow carries a positive sign.
								\item Whenever we apply a sink rule, the induced arrow obtains a negative sign.
							\end{itemize}
							The remaining step via folding is the same as for the complex $\CP(\omega)$.
						\end{enumerate}
						It is not hard to see that this procedure defines a complex $\CP'(\omega)$ of $\DbRep{\A}$. 
						\begin{rmk}
							With this sign convention, any two consecutive gluing edges with the same orientation 
							give rise to a diagram with two positive and two negative arrows:
							
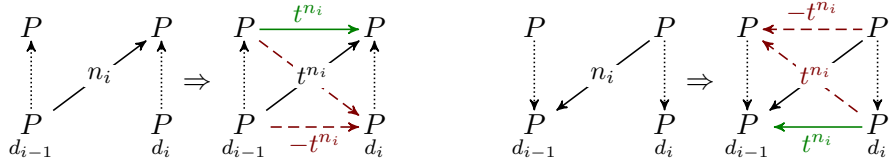
\begin{align}
								\begin{gd}
									P  \ar[densely dotted,<-]{d}     \& P \ar[densely dotted,<-]{d}  \\
									{\underset{d_{i-1}}{P}} \ar{ru}[description]{n_i} \& \underset{d_{i}}{P} 
								\end{gd}
								\Rightarrow 
								\begin{gd}
									P  \ar[densely dotted,<-]{d}   		   \ar[color=green]{r}{t^{n_i}} \ar[color=red,dashed]{rd} \& P \ar[densely dotted,<-]{d}  \\
									{\underset{d_{i-1}}{P}} \ar{ru}[description, outer sep=0pt]{t^{n_i}}	\ar[color=red,dashed]{r}[swap]{- t^{n_i}}	\& \underset{d_{i}}{P} 
								\end{gd}
								&&
								\begin{gd}
									P  \ar[densely dotted]{d}     \& P \ar[densely dotted]{d}  \\
									{\underset{d_{i-1}}{P}} \ar[<-]{ru}[description]{n_i} \& \underset{d_{i}}{P} 
								\end{gd}
								\Rightarrow 
								\begin{gd}
									P  \ar[densely dotted]{d}   		   \ar[color=red,<-,dashed]{r}{- t^{n_i}} \ar[<-,dashed,color=red]{rd} \& P \ar[densely dotted]{d}  \\
									{\underset{d_{i-1}}{P}} \ar[<-]{ru}[color=red,description,outer sep=0pt]{t^{n_i}}
									\ar[color=green,<-]{r}[swap]{t^{n_i}}	\& \underset{d_{i}}{P} 
								\end{gd}
							\end{align}
							Here, solid arrows have positive signs, and dashed arrows negative ones.
							The same rules apply to arrows with Jordan blocks.
							
						\end{rmk}

						\begin{lem}\label{lem:sgns}
							For any band or string $\omega$ of $\chfXA$ there is an isomorphism
							$
							\CP'(\omega) \cong \CP(\omega)$ of complexes.
						\end{lem}
						\begin{proof}
							First we assume that $\omega$ is a band $(d,w,m,\lambda)$
							such that $d$ is odd. We denote the periodic part of $w$ by $\dot{w} = (x_1, \ldots, x_{\ell})$.
							Replacing any negative Jordan block of the form $-J_{\lambda^{\pm 1}}$ by its similar matrix $J_{-\lambda^{\pm 1}}$
							in the differentials of $\CP'(\omega)$
							yields an isomorphic complex.
							
							Depending on the orientation of $x_1$ the 
							local gluing diagram around its number
							may have one of two forms 
							
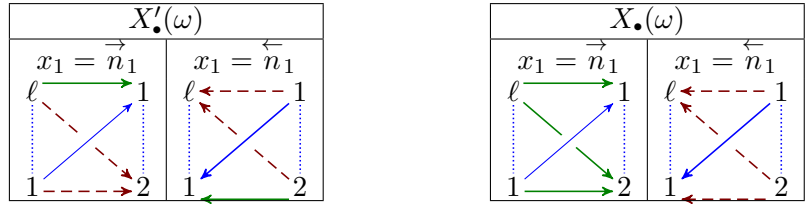
\begin{align*}
								\begin{array}{|c|c|}
									\hline
									\multicolumn{2}{|c|}{\CP'(\omega)}
									\\
									\hline
									x_1 = \ora{n}_1 & x_1 = \ola{n}_1 \\
									\begin{gd}
										\ell \ar[color=green,yshift=4pt, solid]{r} 
										\ar[color=red,dashed]{rd} 
										\ar[color=blue,-,densely dotted]{d}
										\& 
										1
										\ar[color=blue,-,densely dotted]{d} \\
										1
										\ar[color=red,yshift=-1pt,dashed]{r} 
										\ar[color=blue,->,thin, crossing over]{ru}
										\& 
										2
									\end{gd}
									&
									\begin{gd}
										\ell 
										\ar[color=red,yshift=1pt,<-,dashed]{r} 
										\ar[color=red,<-,dashed]{rd} 
										\ar[color=blue,-,densely dotted]{d}
										\& 
										1
										\ar[color=blue,-,densely dotted]{d} \\
										1 \ar[color=green,yshift=-4pt,<-,solid]{r} 
										\ar[color=blue,<-, crossing over]{ru}
										\& 
										2
									\end{gd}\\
									\hline
								\end{array}
								&&
								\begin{array}{|c|c|}
									\hline
									\multicolumn{2}{|c|}{\CP(\omega)}
									\\
									\hline
									x_1 = \ora{n}_1 & x_1 = \ola{n}_1 \\
									\begin{gd}
										\ell \ar[color=green,yshift=4pt, solid]{r} 
										\ar[color=green,solid]{rd} 
										\ar[color=blue,-,densely dotted]{d}
										\& 
										1
										\ar[color=blue,-,densely dotted]{d} \\
										1
										\ar[color=green,yshift=-1pt,solid]{r} 
										\ar[color=blue,->,thin, crossing over]{ru}
										\& 
										2
									\end{gd}
									&
									\begin{gd}
										\ell 
										\ar[color=red,yshift=1pt,<-,dashed]{r} 
										\ar[color=red,<-,dashed]{rd} 
										\ar[color=blue,-,densely dotted]{d}
										\& 
										1
										\ar[color=blue,-,densely dotted]{d} \\
										1 \ar[color=red,yshift=-4pt,<-,dashed]{r} 
										\ar[color=blue,<-, crossing over]{ru}
										\& 
										2
									\end{gd}\\
									\hline
								\end{array}
							\end{align*}
							where solid  arrows have positive signs, and dashed arrows have negative signs.
							The diagrams of $\CP'(\omega)$ and $\CP(\omega)$ can be depicted as follows.
							\begin{align*}
								\CP'(\omega)\colon 
								{
									\tikzcdset{crossing over clearance=0.5ex}
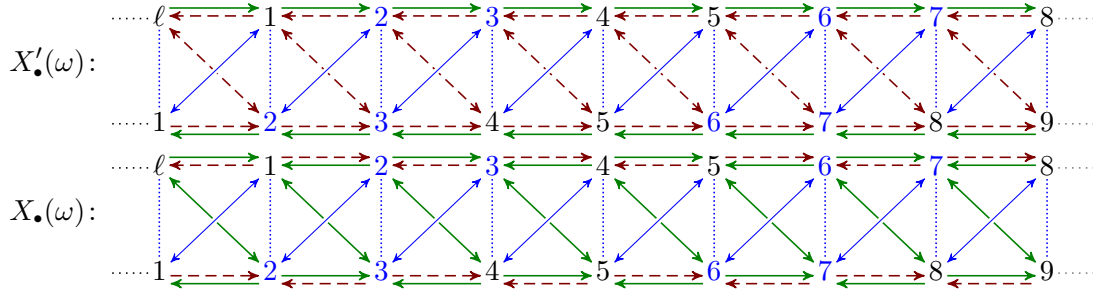
									\begin{gd}
										\mathstrut 	\ar[-, dotted]{r} \&[-0.7cm]	\ell
										\ar[color=green,yshift=4pt,solid]{r} \ar[color=red,yshift=1pt,<-,dashed]{r} 
										\ar[color=red,<->,dashed]{rd}  
										\ar[color=blue,-,densely dotted]{d}
										\&
										1 \ar[color=green,yshift=4pt,solid]{r} \ar[color=red,yshift=1pt,<-,dashed]{r} \ar[color=red,<->,dashed]{rd}  
										\ar[color=blue,-,densely dotted]{d}
										\& 
										{\color{blue}	2 }
										\ar[color=green,yshift=4pt,solid]{r} \ar[color=red,yshift=1pt,<-,dashed]{r} \ar[color=red,<->,dashed]{rd}  
										\ar[color=blue,-,densely dotted]{d}
										\& 
										{\color{blue}	3 }	 
										\ar[color=green,yshift=4pt,solid]{r} \ar[color=red,yshift=1pt,<-,dashed]{r} \ar[color=red,<->,dashed]{rd}  
										\ar[color=blue,-,densely dotted]{d}
										\& 
										4
										\ar[color=green,yshift=4pt,solid]{r} \ar[color=red,yshift=1pt,<-,dashed]{r} \ar[color=red,<->,dashed]{rd}  
										\ar[color=blue,-,densely dotted]{d}
										\& 
										5
										\ar[color=green,yshift=4pt,solid]{r} \ar[color=red,yshift=1pt,<-,dashed]{r} \ar[color=red,<->,dashed]{rd}  
										\ar[color=blue,-,densely dotted]{d}
										\& 
										{	\color{blue} 6}
										\ar[color=green,yshift=4pt,solid]{r} \ar[color=red,yshift=1pt,<-,dashed]{r} \ar[color=red,<->,dashed]{rd}  
										\ar[color=blue,-,densely dotted]{d}
										\& 
										{	\color{blue} 7}	 
										\ar[color=green,yshift=4pt,solid]{r} \ar[color=red,yshift=1pt,<-,dashed]{r} 
										\ar[color=red,<->,dashed]{rd}  
										\ar[color=blue,-,densely dotted]{d}
										\& 
										8
										\ar[color=blue,-,densely dotted]{d}
										\ar[-,dotted, gray]{r}
										\&[-0.7cm]  \mathstrut
										\\
										\mathstrut 	\ar[-, dotted]{r} \&
										1 
										\ar[color=green,yshift=-4pt,<-,solid]{r} \ar[color=red,yshift=-1pt,dashed]{r} 
										\ar[color=blue,<->,thin, crossing over]{ru}\&
										{\color{blue} 2} \ar[color=green,yshift=-4pt,<-,solid]{r} \ar[color=red,yshift=-1pt,dashed]{r} 
										\ar[color=blue,<->,thin, crossing over]{ru}
										\& 
										{\color{blue} 3} \ar[color=green,yshift=-4pt,<-,solid]{r} \ar[color=red,yshift=-1pt,dashed]{r} 
										\ar[color=blue,<->,thin, crossing over]{ru}
										\& 
										4 \ar[color=green,yshift=-4pt,<-,solid]{r} \ar[color=red,yshift=-1pt,dashed]{r} 
										\ar[color=blue,<->,thin, crossing over]{ru}
										\& 
										5 \ar[color=green,yshift=-4pt,<-,solid]{r} \ar[color=red,yshift=-1pt,dashed]{r} 
										\ar[color=blue,<->,thin, crossing over]{ru}
										\& 
										{\color{blue} 6} \ar[color=green,yshift=-4pt,<-,solid]{r} \ar[color=red,yshift=-1pt,dashed]{r} 
										\ar[color=blue,<->,thin, crossing over]{ru}
										\& 
										{\color{blue} 7} \ar[color=green,yshift=-4pt,<-,solid]{r} \ar[color=red,yshift=-1pt,dashed]{r} 
										\ar[color=blue,<->,thin, crossing over]{ru}
										\& 
										8 \ar[color=green,yshift=-4pt,<-,solid]{r} \ar[color=red,yshift=-1pt,dashed]{r} 
										\ar[color=blue,<->,thin, crossing over]{ru}
										\& 9 \ar[-,dotted, gray]{r}\&    \mathstrut
									\end{gd}
								}
								\\
								\CP(\omega)\colon 
								{
									\tikzcdset{crossing over clearance=0.5ex}
									\begin{gd}
										\mathstrut 	\ar[-, dotted]{r} \&[-0.7cm]	\ell
										\ar[color=green,yshift=4pt,solid]{r} \ar[color=red,yshift=1pt,<-,dashed]{r} \ar[color=green,<->,solid]{rd}  
										\ar[color=blue,-,densely dotted]{d}
										\&
										1 \ar[color=red,yshift=4pt,dashed]{r} \ar[color=green,yshift=1pt,<-,solid]{r} \ar[color=green,<->,solid]{rd}   
										\ar[color=blue,-,densely dotted]{d}
										\& 
										{\color{blue} 2} 
										\ar[color=green,yshift=4pt,solid]{r} \ar[color=red,yshift=1pt,<-,dashed]{r} \ar[color=green,<->,solid]{rd}   
										\ar[color=blue,-,densely dotted]{d}
										\& 
										{\color{blue} 3}
										\ar[color=red,yshift=4pt,dashed]{r} \ar[color=green,yshift=1pt,<-,solid]{r} \ar[color=green,<->,solid]{rd}   
										\ar[color=blue,-,densely dotted]{d}
										\& 
										4
										\ar[color=green,yshift=4pt,solid]{r} \ar[color=red,yshift=1pt,<-,dashed]{r} \ar[color=green,<->,solid]{rd}   
										\ar[color=blue,-,densely dotted]{d}
										\& 
										5
										\ar[color=red,yshift=4pt,dashed]{r} \ar[color=green,yshift=1pt,<-,solid]{r} \ar[color=green,<->,solid]{rd}   
										\ar[color=blue,-,densely dotted]{d}
										\& 
										{\color{blue} 6}
										\ar[color=green,yshift=4pt,solid]{r} \ar[color=red,yshift=1pt,<-,dashed]{r} \ar[color=green,<->,solid]{rd}   
										\ar[color=blue,-,densely dotted]{d}
										\& 
										{\color{blue} 7}
										\ar[color=red,yshift=4pt,dashed]{r} \ar[color=green,yshift=1pt,<-,solid]{r}
										\ar[color=green,<->,solid]{rd}   
										\ar[color=blue,-,densely dotted]{d}
										\& 
										8
										\ar[color=blue,-,densely dotted]{d}
										\ar[-,dotted, gray]{r}
										\&[-0.7cm]   \mathstrut
										\\
										\mathstrut 	\ar[-, dotted]{r} \& 
										1 
										\ar[color=green,yshift=-4pt,<-,solid]{r} \ar[color=red,yshift=-1pt,dashed]{r} 
										\ar[color=blue,<->,thin, crossing over]{ru}\&
										{\color{blue} 2} \ar[color=red,yshift=-4pt,<-,dashed]{r} \ar[color=green,yshift=-1pt,solid]{r} 
										\ar[color=blue,<->,thin, crossing over]{ru}
										\& 
										{\color{blue} 3} \ar[color=green,yshift=-4pt,<-,solid]{r} \ar[color=red,yshift=-1pt,dashed]{r} 
										\ar[color=blue,<->,thin, crossing over]{ru}
										\& 
										4 \ar[color=red,yshift=-4pt,<-,dashed]{r} \ar[color=green,yshift=-1pt,solid]{r} 
										\ar[color=blue,<->,thin, crossing over]{ru}
										\& 
										5 \ar[color=green,yshift=-4pt,<-,solid]{r} \ar[color=red,yshift=-1pt,dashed]{r} 
										\ar[color=blue,<->,thin, crossing over]{ru}
										\& 
										{\color{blue} 6} \ar[color=red,yshift=-4pt,<-,dashed]{r} \ar[color=green,yshift=-1pt,solid]{r} 
										\ar[color=blue,<->,thin, crossing over]{ru}
										\& 
										{\color{blue} 7} \ar[color=green,yshift=-4pt,<-,solid]{r} \ar[color=red,yshift=-1pt,dashed]{r} 
										\ar[color=blue,<->,thin, crossing over]{ru}
										\& 
										8 \ar[color=red,yshift=-4pt,<-,dashed]{r} \ar[color=green,yshift=-1pt,solid]{r} 
										\ar[color=blue,<->,thin, crossing over]{ru}
										\& 9 \ar[-,dotted, gray]{r}\&   \mathstrut
									\end{gd}
								}
							\end{align*}
							In each square of each diagram the solid diagonal arrow between the bottom left and the top right number is an arrow  decorated by a non-induced entry $t^{n_i}$.
							Each square of each diagram subsumes both possibilities for the orientation of a letter $x_i$, and either all arrows in such a square are oriented to the left, or all arrows are oriented to the right.
							
							In the course of the construction of the complexes, the vertices in the diagrams are replaced by projectives.
							There is an isomorphism $\begin{td} \CP'(\omega) \ar{r}{\sim} \& \CP(\omega) \end{td}$ 
							given by multiplying 
							the $m$ rows as well as the $m$ columns 
							containing $t^{n_i}$ as non-induced entry 
							with $-1$
							for each index $1 \leq i \leq \ell$ such that $i \equiv 2 \text{ or } 3$ modulo $4$.
							
							Assume that $\omega$ is a band $(d,w,m,\lambda)$ with an even integer $d\in \Z$. Then $\omega[1] = (d+1,w,m,\lambda)$ 
							is a band starting with an odd degree.
							It follows that
							$\CP'(\omega)[1] \cong \CP'(\omega[1]) \cong \CP(\omega[1]) \cong \CP(\omega) [1]$, where the second isomorphism is due to the first part of the proof. Thus, $\CP'(\omega) \cong \CP(\omega)$.
							
							The same arguments apply to any string $\omega$ with an underlying word of length $\ell$.
						\end{proof}
						\begin{rmk}
							The advantage of the presentation $\CP'(\omega)$ over $\CP(\omega)$ is that the sign rules do not depend on the degrees associated to the diagram vertices.
							The disadvantage of $\CP'(\omega)$ is that we lose the `signless' property of differentials of $\CP(\omega)$ at odd degrees.
							We will work with the presentation $\CP'(\omega)$ to prove Theorem~\ref{thm:bijection2}, and with $\CP(\omega)$ starting from Section~\ref{sec:derived_fun} when computing functors and restricting our focus to projective resolutions of finite-dimensional $\A$-modules.
						\end{rmk}

						\subsection{Construction of pullback complexes of bands and strings}
						\label{subsec:glurul}	
						\begin{dfn}
							For a band or string $\omega$  of $\chfXA$, the \emph{pullback complex}
							$\CPr(\omega)$ is the image under the composition of bijections in~\eqref{eq:bij-diag}.
						\end{dfn}
						In more detail, the pullback complex is defined by the following steps.
						\begin{itemize}
							\item If $\omega$ is a trivial string, then $\omega$ is equivalent to $(d,(\star,\ola{n},\star))$ and $\CPr(\omega)$ is the complex
							$(\begin{td} P_\star \ar{r}{t^n} \& P_\star \end{td})[d]$.
							\item Assume that $\omega$ is non-trivial.
							\begin{enumerate}
								\item The datum $\omega$ is expanded into a band respectively string $\omega_e$ of $\fXA$ via the translations in Subsection~\ref{subsec:exp-bnd} respectively~\ref{subsec:exp-str}.
								\item $\omega_e$ yields a regular matrix representation $\mu(\omega_e)$ in $\regrep{\fXA}$ by the three-step construction of canonical forms in Subsection~\ref{subsec:strbdrep}.
								\item The matrix representation $\mu(\omega_e)$ determines a gluing triple $\Ind(\mu(\omega_e))$ in $\Gcatfd{\A}$ via Remark~\ref{rmk:mat-gt}.
								\item Then $\CPr(\omega)$ is precisely the pullback complex $\Gl\Ind(\mu(\omega_e))$ as defined in \eqref{eq:pb-complex}.
							\end{enumerate}
						\end{itemize}
						In particular, $\CPr(\omega)$ is an indecomposable projective complex of $\DbRep{\A}$.
						
						\begin{rmk}\label{rmk:ind-triple}
							Assume that $\omega$ is band with multiplicity one or a string.
							The complexes of the gluing triple $\Ind(\mu(\omega_e))=(\CV,\CY, \vartheta)$ can be described using the polarized gluing diagram of $\omega$.
							\begin{enumerate}
								\item The complex of $\AI$-modules $\CV$ is the direct sum $\bigoplus_{i \in \Z} (S_+^{m_{+,i}} \oplus S_-^{m_{-,i}}) [i]$ with $m_{\pm,i}$ given by the multiplicity of projectives of type $P_{\pm}$ with underlying degree $i$ in the gluing diagram. 
								\item  The complex of $\B$-modules $\CY$ is the direct sum of indecomposable complexes $$\CY = \bigoplus_{
									\begin{smallmatrix}
										u,v \in \{ \star, +,-\}\\ n\in \mathcal{N}_{uv}, i \in \Z \end{smallmatrix}}
								\bigl((\begin{gd} \BP_{u'} \ar{r}{t^{n}} \&  \BP_{v'}
								\end{gd}) [i]\bigr)^{m_{u,v,n,i}}
								$$
								where we set $\mathcal{N}_{uv} = \N_0$ if $(u,v)=(\pm,\star)$ respectively $\mathcal{N}_{uv} = \N$ otherwise,
								$u' \colonequals {\diamond}$ if $u =\pm$ and $u'  = \star$ otherwise, the vertex $\tilde{v}$ is defined similarly, and $m_{u,v,n,i}$ is the total number of 
								diagonal morphisms of the form 
								$$\begin{gd} \&  \smash[b]{\underset{i}{P_v}}\\
									\underset{i+1}{P_{u}} \ar{ru}[description]{t^n} \&  \end{gd} \qquad\text{or}\qquad
								\begin{gd} \&  \smash[b]{\underset{i+1}{P_u}}\\
									\underset{i}{P_{v}} \ar[<-]{ru}[description]{t^n} \&  \end{gd}
								$$
								in the polarized gluing diagram of $\omega$.
								If $\omega$ is band with multiplicity $m$, such diagonal morphisms appear as $m$-fold direct sums.
							\end{enumerate}
							The gluing morphism $\vartheta$ of the triple $\Ind({\mu(\omega_e)})$
							is given by a family $(\vartheta_i)_{i \in \Z}$ of isomorphisms of stalk complexes of $\BI$-modules and
							can be determined using the oriented gluing diagram of $\omega$ as follows.
							\begin{enumerate}\setcounter{enumi}{2}
								\item 
								\begin{enumerate}
									\item
									By the preceding observations, any entry $P_{\pm}$ at a degree $i$ 
									without decoration by an eigenvalue in the oriented gluing diagram corresponds to an $\AI$-module
									$\V_\pm[i]$ as well as a $\B$-module $\BP_\diamond[i]$.
									Such a projective contributes an identity entry
									$\begin{td} \BI \oAI S_{\pm}[i] \ar{r}{1} \& \BI \oB \BP_\diamond[i] \end{td}$
									in the morphism $\vartheta_i$,
									where the domain as well as the codomain are identified with the stalk complex of the simple $\BI$-module at degree $i$.
									\item
									An oriented gluing edge at degree $i$ decorated by an eigenvalue $\lambda$ translates into a $2$-by-$2$--minor of $\vartheta_i$  as follows.
									\begin{align*}
										\begin{gd}
											P_{\pm} \ar[->,densely dotted]{d}\\
											\underset{i,\lambda}{P_{\mp}}
										\end{gd}
										\Rightarrow
										\begin{td}
											\BI \oAI S_{\pm}[i] \ar{r}{1} \& \BI \oB \BP_{\diamond}[i] \\
											\underset{\mathstrut}{\BI \oAI S_{\mp}[i]}  \ar{r}{\lambda}  \ar{ru}[description]{1}\& 
											\BI \oB \BP_{\diamond}[i]
										\end{td}
										&&
										\begin{gd}
											P_{\pm} \ar[<-,densely dotted]{d}\\
											\underset{i,\lambda}{P_{\mp}}
										\end{gd}
										\Rightarrow
										\begin{td}
											\BI \oAI S_{\pm}[i] \ar{r}{1}\ar{rd}[description]{\lambda}  \& \BI \oB \BP_{\diamond}[i] \\
											\underset{\mathstrut}{\BI \oAI S_{\mp}[i]}  \ar{r}{1}  \& 
											\BI \oB \BP_{\diamond}[i]
										\end{td}
									\end{align*}
									If there is no decoration $\lambda$, the same rules apply with $\lambda$ replaced by $1$.
									All remaining entries of $\vartheta_i$ at any degree $i \in \Z$ are zero.
								\end{enumerate}
							\end{enumerate}
							If $\omega$ is a band with multiplicity $m > 1$, the same considerations apply with $m$-fold direct sums, units replaced by identity matrices and $\lambda$ by the Jordan block $\JB_{\lambda}$ of size $m$.
						\end{rmk}

						Let $\omega$ be a band or non-trivial string of $\chfXA$.
						On the one hand, we
						define a complex $\CPr'(\omega)$ 
						by changing certain signs in the band respectively string complex $\CPr(\omega)$.
						On the other hand, we will show that $\CPr'(\omega)$ is precisely the complex
						obtained from $\omega$ by following the inverse direction in the categorical reductions of the previous section.

						\paragraph{Gluing diagrams of width two}
						
						Next, we show that the definition of $\CPr'(\omega)$ 
						matches the applications of the map on objects
						from $\regrep{\fXA}$ to $\DbRep{\A}$ in a special case, which is given by a gluing diagram corresponding to a complex with only two non-trivial terms.
						This case will provide a local template which is extendable to the general situation.

						\begin{lem}\label{lem:spc-cpx}	
							Assume that $\omega$ is a usual word of $\chfXA$ with gluing diagram
							\begin{align}
								\label{eq:rnd-ornt}
								\begin{gd}
									\& P_{+} \ar[densely dotted,<->]{d}[swap]{g_1} \& P_{+}  \ar[densely dotted,<->]{d}[swap]{h_1} \&
									P_{+} \ar[densely dotted,<->]{d}[swap]{g_2} \& P_{+}  \ar[densely dotted,<->]{d}[swap]{h_2} \&
									\cdots
									\&  
									P_{+} \ar[densely dotted,<->]{d}[swap]{g_k}  \& \smash[b]{\underset{0}{P_{\star}}} \\
									\underset{1}{P_{\star}}  \ar[->]{ru}[description]{t^{p_1}} \& \underset{0}{P_{-}}  \ar[<-]{ru}[description]{t^{q_1}} \& 
									\underset{1}{P_{-}}  \ar[->]{ru}[description]{t^{p_2}} \& \underset{0}{P_{-}}  \ar[<-]{ru}[description]{t^{q_2}}  \& \underset{1}{P_{-}}\&   \cdots
									\& \underset{1}{P_{-}} \ar[<->]{ru}[description]{t^{p_{k}}} 	\& 
								\end{gd}
							\end{align}
							For each $1 \leq i \leq k$ with $i < k$ in case of $\ol{\gamma}_i$
							we denote binary coefficients 
							\begin{align}
								\label{eq:rnd-diff}
								\begin{array}{ll}
									(\alpha_i,\ol{\alpha}_i) = \begin{cases}
										(-1,0) & \text{if }g_i = \downarrow, \\
										(0,-1) & \text{if }g_i = \uparrow,
									\end{cases}
									&
									\gamma_i = \alpha_i \beta_i = 
									\begin{cases}
										-1 & \text{if } (g_i,h_i) = (\downarrow,\downarrow), \\
										\phantom{-}0 & \text{otherwise},
									\end{cases}
									\\
									\\
									(\beta_i,\ol{\beta}_i) = \begin{cases}
										(1,0) & \text{if } h_i = \downarrow, \\
										(0,1) & \text{if }h_i = \uparrow,
									\end{cases}     
									&
									\ol{\gamma}_i = \ol{\beta}_i \ol{\alpha}_{i+1} = 
									\begin{cases}
										1 & \text{if } (h_i,g_{i+1}) = (\uparrow,\uparrow), \\
										0 & \text{otherwise}.
									\end{cases}   
								\end{array}
							\end{align}
							Then the pullback complex $\CPr(\omega)$ 
							is isomorphic to the string complex $\CP(\omega)$ with degree-independent sign convention,
							that is, to the two-term complex
							$(\begin{td}
								\AP_1 \ar{r}{\partial_1}  \& \AP_0
							\end{td})$
							with the differential $\partial_1$ given by the pentadiagonal matrix
							\begin{align}
								\label{rnd:diff}
								\begin{pNiceArray}{ccccccccc}[first-row, last-col
									, columns-width = auto
									]
									\CodeBefore
									\chessboardcolors{red!10}{blue!10}
									\foreach \i in {1,...,9}{
										\pgfmathtruncatemacro{\jmin}{max(1,\i-2)}
										\pgfmathtruncatemacro{\jmax}{min(9,\i+2)}
										\foreach \j in {\jmin,...,\jmax}{
											\pgfmathtruncatemacro{\sum}{mod(\i-\j,2)}
											\ifnum\sum=0
											\cellcolor{red!25}{\i-\j}
											\else
											\cellcolor{blue!25}{\i-\j}
											\fi
										}
									}
									\Body
									P_{\star} & P_+ & P_- & P_+ &&& P_- & P_+ & P_- 
									\\
									t^{p_1} & \alpha_1 t^{q_1} & \gamma_{1} t^{q_1} & 0 &\Cdots &\Cdots&\Cdots&\Cdots&0&P_+ \\
									\ol{\alpha}_1 t^{p_1} & t^{q_1} & \beta_1 t^{q_1} & 0 &\Ddots&&& &\Vdots & P_-  \\
									0 & \ol{\beta}_1 t^{p_2} & t^{p_2 }& \alpha_2 t^{q_2} & \Ddots &&&&\Vdots&P_+  \\
									0 & \ol{\gamma}_1 t^{p_2} & \ol{\alpha}_2 t^{p_2} & t^{q_2} &\Ddots &\Ddots&&& \Vdots& P_- \\
									\Vdots		&\Ddots &\Ddots &\Ddots &\Ddots & \Ddots &\Ddots&  \\
									\Vdots	&&\Ddots &\Ddots &\Ddots &\Ddots & \Ddots & \Ddots & 0 &  \\
									\Vdots &&&&		\Ddots & \Ddots &		
									t^{p_{k-1}}
									& 
									*
									& *
									&  P_+  \\
									\Vdots &&	&&& \Ddots & 	
									*
									& 
									t^{q_{k-1}} 
									& 
									*
									&  P_-  \\
									0 &	\Cdots 	&&&& 0&		* & 
									*
									& t^{p_{k}}&  P_\star	\end{pNiceArray}
							\end{align}

									\end{lem}
									\begin{proof} 
										To simplify the presentation we prove the claim in the case $k = 4$, the general case follows by the same argument.
										
										By the translations in Subsection \ref{subsec:exp-str} for gluing diagrams
										together with the translation of oriented gluing edges in Subsection~\ref{subsec:od-mat}
										the oriented gluing diagram in \eqref{eq:rnd-ornt} encodes the decorated word
										\begin{align*}
											w = \yy_0^{2p'_1} - 
											&\overleftrightarrow{\xf_0^+ \sim \xf_0^-} - \yy_0^{2q_1} \sim \xx_1^{2 q_1} -
											\overleftrightarrow{\xf_1^+ \sim \xf_1^-} - \xx_1^{2p_2} \sim   \\
											\yy_0^{2 p_2} - &\overleftrightarrow{\xf_0^+ \sim \xf_0^-} - \yy_0^{2q_2} \sim \xx_1^{2 q_2} -
											\overleftrightarrow{\xf_1^+ \sim \xf_1^-} - \xx_1^{2p_3} \sim 
											\\ 
											\yy_0^{2 p_3} - &\overleftrightarrow{\xf_0^+ \sim \xf_0^-} - \yy_0^{2q_3} \sim \xx_1^{2 q_3} -
											\overleftrightarrow{\xf_1^+ \sim \xf_1^-} - \xx_1^{2p''_4}
										\end{align*}
										where
										$2p'_1= 2p_1 - 1$ and $2p''_4= 2p_4 + 1$,
										the arrow in the $i$-th subword of the form  $\overleftrightarrow{f_0^+ \sim f_0^-}$
										is left oriented if and only if 
										$\alpha_i = 1$ or, equivalently, $g_i = \downarrow$, and, similarly, 
										the arrow in the $i$-th  subword  $\overleftrightarrow{f_1^+ \sim f_1^-}$
										is left oriented if and only if $\beta_i =1$ or $h_i = \downarrow$
										for each index $1 \leq i \leq 3$.
										
										The word $w$ gives rise 
										to two triagonal matrices
										$M_1, M_0 \in \Mat_{6 \times 6}(\kk)$
										by the construction in \ref{subsec:can}.
										These matrices will be described below and constitute
										a representation $\mu$ of $\regrep{\fXA}$.
										By Remark~\ref{rmk:mat-gt}
										the representation $\mu$ gives rise to the following object
										$\Ind(\mu)=(\CV, \CY, \vartheta)$ in $\Gcatfd{\A}$.
										\begin{enumerate}
											\item
											The complex $\CV$ is given by the two-term complex $(\begin{td} \V_1  \ar{r} \&   \V_0 \end{td})$ with zero differential
											and $\V_1 = \V_0 = (S_+ \oplus S_-)^3$. For each index $i \in \{0,1\}$
											the module $\BI \oAI \V_i$
											has a set of generators indexed by elements of $\fF_i$ appearing in $w$.
											\item
											The complex
											$\CY$ is the two-term complex $(\begin{td} \Y_1 \ar{r} \& \Y_0 \end{td})$ given by the direct sum 
											$$
											\begin{cd}
												\BP_{\star} \oplus \BP_{\diamond}^6 
												\ar{r}{\tilde{\partial}_1} \&
												\BP_{\diamond}^6 \oplus \BP_{\star} 
											\end{cd} 
											\quad
										\text{ with }\tilde{\partial}_1 = \mathrm{diag}(t^{p_1},t^{q_1},t^{p_2},t^{q_2},t^{p_3},t^{q_3},t^{p_4}).$$
										The projective module
										$\Y_1$ 
										has a set of generators given 
										by the union of the formal symbol
										$\xx^{\star}_1$ representing the generator of $\BP_{\star}$ together 
										with letters from $\fE_1$ appearing in $w$. Similarly, $\Y_0$ admits a set of generators given by letters from $\fE_0$ appearing in $w$ together with a symbol $\yy^{\star}_0$.
										\item
										The matrices $\vartheta_1,\vartheta_0 \colon \begin{td} \wt{S}_\diamond^6 \ar{r} \&  \wt{S}_\diamond^6 \end{td}$ are given by application of the Morita equivalence $\begin{td} \md \kk \ar{r}{\sim} \& \md \BI \end{td}$
										to the $\kk$-linear maps of the matrices $M_1$ and $M_0$.
									\end{enumerate}
									Applying $\Gl$ to $\Ind({\mu})$ 
									yields a 
									two-term complex
									$\CPr = (\begin{td} 
										\PPr_1 \ar{r}{\partial_1^{\CPr}} \&  \PPr_0 \end{td})$
									of projective $\A$-modules
									with 
									$\PPr_1 = 
									P_{\star} \oplus (P_+ \oplus P_-)^3$	
									and $\PPr_0 	= (P_+ \oplus P_-)^3 \oplus P_{\star}$.
									In particular,  $\PPr_1$ has a set of generators given by a symbol $\xf_1^{\star}$ together with lifts of the generators of $\V_1$,
									and, similarly, there are generators of $\PPr_0$ which are lifts of generators of $\V_0$ together with a symbol $\xf_0^{\star}$.
									Moreover, there are natural 
									$\A$-linear embeddings
										$\begin{td} \PPr_1 
											\ar[hookrightarrow]{r} \& \Y_1 \end{td}$
										and
										$\begin{td} \PPr_0  \ar[hookrightarrow]{r} \& \Y_0 \end{td}$.
										According to formula \eqref{eq:diff}
										the composition
										\begin{align*}\begin{cd}
												\Y_1 \ar{r}{\hat{M}_1} \& \Y_1 \ar{r}{\wt{\partial}_1} \& \Y_0 \ar{r}{\hat{M}_0^{-1}} \& \Y_0
											\end{cd}
										\end{align*}
										can be restricted to a map $\begin{td} \PPr_1 \ar{r} \& \PPr_0 \end{td}$
										which  yields
										the differential $\partial_1$ of $\CPr$,
										where
										the matrices 
										$\hat{M}_1$ and $\hat{M}_0$ 
										are lifts of $M_1$ and $M_0$ obtained using Remark~\eqref{eq:lift}.
										More precisely, 
										the three factors of $\partial_1$ are 
										given as follows.
										\begin{align*}
											\scalebox{0.6}{$
												\underbrace{		\begin{pNiceArray}{|cc:cc:cc|c|}[first-row, last-col, columns-width=auto]
														\yy^{2p'_1}_0 &\yy^{2q_1}_0& \yy^{2p_2}_0& \yy^{2q_2}_0 & \yy^{2p_3}_0& \yy^{2q_3}_0 & \yy_0^{\star}\\
														1 & \alpha_1 & \Block{2-2}{0} & & \Block{2-2}{0} &   & \Block{2-1}{0} &\xf_0^+ \\
														{\ol{\alpha}_1}
														& 1 &  &  &   &  & &\xf_0^- \\ 
														\hdottedline
														\Block{2-2}{0} &  & 1 & \alpha_2 &\Block{2-2}{0} & & \Block{2-1}{0} &\xf_0^+\\
														&  & \ol{\alpha}_2 & 1 &&& &\xf_0^- \\
														\hdottedline
														\Block{2-2}{0}  &&\Block{2-2}{0}& & 
														1 & \alpha_{3} & \Block{2-1}{0} &\xf_0^+\\
														& && & {\ol{\alpha}_{3}} & 1 && \xf_0^-\\
														\hline
														\Block{1-2}{0} & & \Block{1-2}{0} && \Block{1-2}{0} & & {\color{blue} 1 }& \xf_0^\star
												\end{pNiceArray}}_{\displaystyle \hat{M}_0^{-1}}
												\cdot
												\underbrace{
													\begin{pNiceArray}{ccccccc}[first-row,last-col, columns-width=auto]
														\xx_1^\star & \xx^{2q_1}_1&\xx^{2p_2}_1 & \xx^{2q_2}_1&\xx^{2p_3}_1 & \xx^{2q_3}_1 & \xx^{2p''_4}_1\\
														t^{p_1} & 0 & \Cdots   &\Cdots   &\Cdots   & \Cdots & 0 & \yy^{2p'_1}_0\\
														0 & t^{q_1} & \Ddots &  & &  & \Vdots &\yy^{2q_1}_0\\
														\Vdots 	&  \Ddots & t^{p_2} &\Ddots  &&&  \Vdots&\yy^{2p_2}_0 \\
														\Vdots 			&  & \Ddots& t^{q_2} &\Ddots&  &  \Vdots&\yy^{2q_2}_0\\
														\Vdots &&&\Ddots& t^{p_3} & \Ddots & \Vdots& \yy^{2p_k}_0\\
														\Vdots & & &  & & t^{q_3}& 0		&  \yy^{2q_3}_0 \\
														0 &\Cdots &\Cdots & \Cdots && 0& t^{p_4}		&  \yy_0^*
													\end{pNiceArray}
												}_{ \displaystyle \tilde{\partial}_1}
												\cdot
												\underbrace{
													\begin{pNiceArray}{|c|cc:cc:cc|}[first-row, last-col, columns-width=auto]
														\xf_1^\star & \xf_1^+ & \xf_1^- &\xf_1^+ & \xf_1^- &\xf_1^+ & \xf_1^-  \\
														{\color{blue} 1 } & \Block{1-2}{0} & & \Block{1-2}{0} & & \Block{1-2}{0} &  & \xx^{\star}_1 \\ \hline
														\Block{2-1}{0} &	1 & \beta_1 & \Block{2-2}{0} & & \Block{2-2}{0} &   & \xx^{2q_1}\\
														&	{\ol{\beta}_1}
														& 1 &  &  & &&    \xx^{2p_2}\\ 
														\hdottedline
														\Block{2-1}{0} &	\Block{2-2}{0} &  & 1  &  \beta_2 &\Block{2-2}{0} & & \xx^{2q_2}\\
														&	&  & \ol{\beta}_2&1&&& \xx^{2p_3} \\
														\hdottedline
														\Block{2-1}{0} & \Block{2-2}{0}	  &&\Block{2-2}{0}& & 1 & \beta_{3} &  \xx^{2q_3} \\ 
														&	&&& &  \ol{\beta}_{3}  & 1 &  \xx^{2p''_4}\\
														\hline
													\end{pNiceArray}
												}_{{\displaystyle \hat{M}_1}}$
											}
										\end{align*}
										
										We note that the matrices $M_1$ and $M_0^{-1}$ are maximal minors of $\hat{M}_1$ and $\hat{M}_0^{-1}$, respectively, and that each lifted matrix was obtained by adding an identity matrix as direct summand to $M_i$ with appropriate labeling.
										A computation of the composition above
										yields that
										\begin{align*}
											\partial_1 
											=
											\begin{pNiceArray}{ccccccc}[first-row, last-col
												, columns-width = auto,name=mymatrix
												]
												\CodeBefore
												\chessboardcolors{red!10}{blue!10}
												\foreach \i in {1,...,7}{
													\pgfmathtruncatemacro{\jmin}{max(1,\i-2)}
													\pgfmathtruncatemacro{\jmax}{min(7,\i+2)}
													\foreach \j in {\jmin,...,\jmax}{
														\pgfmathtruncatemacro{\sum}{mod(\i-\j,2)}
														\ifnum\sum=0
														\cellcolor{red!25}{\i-\j}
														\else
														\cellcolor{blue!25}{\i-\j}
														\fi
													}
												}
												\Body
												P_{\star} & P_+ & P_- & P_+ & P_- & P_+ & P_- & \\
												t^{p_1} & \alpha_1 t^{q_1} & \gamma_{1} t^{q_1} & 0 &\Cdots &\Cdots&0&P_+ \\
												\ol{\alpha}_1 t^{p_1} & t^{q_1} & \beta_1 t^{q_1} & 0 &\Ddots&&\Vdots & P_-  \\
												0 & \ol{\beta}_1 t^{p_2} & t^{p_2 }& \alpha_2 t^{q_2} & \gamma_2 t^{q_2}  & \Ddots & \Vdots &P_+  \\
												0 & \ol{\gamma}_1 t^{p_2} & \ol{\alpha}_2 t^{p_2} & t^{q_2} &\beta_2 t^{q_2} &0& 0 & P_- \\
												\Vdots		&\Ddots & 0 & \ol{\beta}_2 t^{p_3} &t^{p_3} & \alpha_3 t^{q_3}  & \gamma_3 t^{q_3} & P_+ \\
												\Vdots	& & \Ddots &\ol{\gamma}_2 t^{p_3} &\ol{\alpha}_3 t^{p_3} &t^{q_3} & \beta_3 t^{q_3} & P_- \\
												0 &\Cdots& \Cdots & 0 &		0 &  \ol{\beta}_3 t^{p_4}&		
												t^{p_{4}}
												& P_\star  \\
											\end{pNiceArray}
									\end{align*}
									which is a special case of the matrix $\partial_1$ from \eqref{rnd:diff}.
									It is straightforward to verify that $\partial_1$ is precisely the differential obtained by application of the gluing rules to the oriented gluing diagram in~\eqref{eq:rnd-ornt}.
								\end{proof}

							\begin{lem}\label{lem:spc-band}
								Let $m\in \N$ and assume that $\lambda \in \kk^*$ 
								decorates a vertex incident to one of the gluing edges in \eqref{eq:rnd-ornt}, that is,
								there is an index $1 \leq i < k$  such that
								one of the following cases occurs.
								
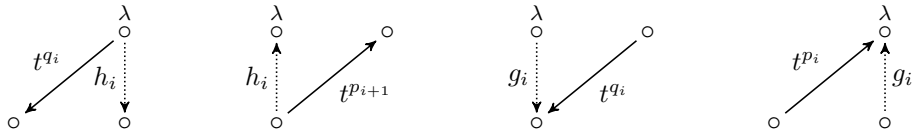
\begin{align*}
									\begin{gd}
										\& 			\overset{\lambda}{\bt} \ar[densely dotted]{d}[swap]{h_i}  \\
										\bt \ar[<-]{ru}{t^{q_i}}
										\&  			\bt
									\end{gd}
									&&
									\begin{gd}
										\overset{\lambda}{\bt} \ar[densely dotted,<-]{d}[swap]{h_i}  \& \bt \\
										{\bt} \ar{ru}[swap]{t^{p_{i+1}}}
										\& 
									\end{gd}
									&&
									\begin{gd}
										\overset{\lambda}{\bt} \ar[densely dotted]{d}[swap]{g_i}  \& \bt \\
										\bt \ar[<-]{ru}[swap]{t^{q_i}}
										\& 
									\end{gd}
									&&
									\begin{gd}
										\& 	\overset{\lambda}{\bt} \ar[densely dotted,<-]{d}{g_i}\\
										{\bt} \ar{ru}{t^{p_i}} \& \bt
									\end{gd}
								\end{align*}
								Then the conclusion of Lemma~\ref{lem:spc-cpx}
								remains valid with the following modifications.
								\begin{itemize}
									\item Each row label and each column label has to be replaced by its $m$-fold direct sum.
									\item Each entry has to be viewed as a square matrix of size $m$ with the original entry on the diagonal.
									\item One of the coefficients in the matrix~\eqref{rnd:diff}
									has to be redefined in each of the four cases as follows.
									\begin{align*}
										\beta_i = \JB_\lambda && \ol{\beta}_i = \JB_{\frac{1}{\lambda}} && \alpha_i = - \JB_\lambda &&
										\ol{\alpha}_i = -\JB_{\frac{1}{\lambda}}
									\end{align*}
									where $\JB_\lambda$ denotes a Jordan block of size $m$ and eigenvalue $\lambda$.
								\end{itemize}
							\end{lem}
							\begin{proof}
								\begin{enumerate}
									\item
									In the first two cases,
									to simplify notation we denote $p =p_{i+1}$ and $q= q_i$.
									\begin{enumerate}
										\item In the first case, 
										the subword
										$
										\xx^{2 q} - \overleftarrow{\xf_1^+ \sim \xf_1^-} - \xe''
										$ of $w$
										contributes
										the coefficient matrix $\beta_i = \JB_\lambda$ 
										in the horizontal stripe with label
										$\xx^{2q}$ and the vertical stripe labeled $\xf_1^-$ of $M_1$.
										In the following, we will only depict the relevant minors of the matrices. In case of the matrix $M_1$ this yields the following abbreviated notation.
										\begin{align*}
											M_1 = 
											\begin{pNiceArray}{c:cc:c}[first-row, last-col, columns-width=auto]
												& \xf_1^{+} & \xf_1^{-} & & \\
												*& 0 &0  & *&\\
												\hline
												0&\Id & \JB_\lambda &0 & \xx^{2q}\\
												0&0 & \Id & 0& \xe'' \\
												\hline
												*		&0 &0 & * & 
											\end{pNiceArray}
											\qquad
											\Rightarrow
											\qquad
											M_1: \
											\begin{bNiceArray}{cc}[first-row, last-col, columns-width=auto]
												\xf_1^{+} & \xf_1^{-}  & \\
												\Id & \JB_\lambda  & \xx^{2q}\\
												0 & \Id  & \xe'' \\
											\end{bNiceArray}
										\end{align*}

										\item 
										In the second case,
										the parameter $\lambda$ decorates the relation in $\xe' - \xf_1^+$
										inside of the subword $\xe' - \overrightarrow{\xf_1^+ \sim \xf_1^-} - \xe'$.
										In this situation, we may 
										replace the matrix 
										$\JB_{\lambda}$ in the block $(\xe',\xf_1^+)$ of matrix $M_1$
										by the similar matrix
										$\JB^{-1}_{\lambda^{-1}}$ 
										according to Remark~\ref{rmk:sym-isom},
										and pass to an isomorphic representation $(M'_1, M'_0)$ of $\fXA$
										by multiplying the vertical stripe with index $\xf_1^+$ and block matrix $\JB^{-1}_{\lambda^{-1}}$ in the horizontal stripe $\xe'$ of matrix $M_1$ with $\JB_{\frac{1}{\lambda}}$.
										In this case, 
										\begin{align*}
											M_1: \
											\begin{bNiceArray}{cc}[first-row, last-col, columns-width=auto]
												\xf_1^{+} & \xf_1^{-} &  \\
												\JB_{\frac{1}{\lambda}}^{-1} & 0 &  \xe' \\
												\Id & \Id &  \xx^{2 p} \\
											\end{bNiceArray}
											\cong
											\begin{bNiceArray}{cc}[first-row, last-col, columns-width=auto]
												\xf_1^{+} & \xf_1^{-} &  \\
												\Id & 0 &  \xe'\\
												\JB_{\frac{1}{\lambda}} & \Id &  \xx^{2 p} 
											\end{bNiceArray}
										\end{align*}
									\end{enumerate}
									We obtain the claim by the following computation summarizing the first two cases.
									\begin{align*}
										\begin{bNiceArray}{cc}[first-row, last-col, columns-width=auto]
											\yy^{2q}  & \yy^{2 p} & \\
											\Id & 0 &  \xf_0^-\\
											0 & \Id &  \xf_0^+ \\
										\end{bNiceArray}
										\cdot
										\begin{bNiceArray}{cc}[first-row, last-col, columns-width=auto]
											\xe' & \xe''  & \\
											t^{q} & 0  & \yy^{2q}\\
											0 & t^{p}  & \yy^{2 p} \\
										\end{bNiceArray}
										\cdot
										\begin{bNiceArray}{cc}[first-row, last-col, columns-width=auto]
											\xf_1^{+} & \xf_1^{-} &  \\
											\Id & \beta_i &  \xe'\\
											\ol{\beta}_i & \Id &  \xe'' \\
										\end{bNiceArray}
										= 
										\begin{bNiceArray}{cc}[first-row, last-col, columns-width=auto]
											\xf_1^{+} & \xf_1^{-}  & \\
											t^{q} & \beta_i t^q &  \xf_0^-\\
											\ol{\beta}_i t^{p} & t^{p}  & \xf_0^+ \\
										\end{bNiceArray}
									\end{align*}
									\item In the last two cases, we abbreviate $q= q_i$ and $p = p_i$, respectively.
									\begin{enumerate}
										\item In the third case, the parameter $\lambda$ decorates the relation 
										of the long elementary subword in $\xe' - \overleftarrow{\xf_0^+ \sim \xf_0^-} - \yy^{2q}$.
										\item
										In the last case, 
										the parameter $\lambda$ decorates the first relation in the subword $\yy^{2p} - \overrightarrow{\xf_0^+ \sim \xf_0^-} -\xe''$.
										Again, we use the matrix $\JB_{\frac{1}{\lambda}^{-1}}$ instead of $\JB_{\lambda}$ and pass to an isomorphic representation of $\fXA$ 
										by multiplying the vertical stripe $\xf_0^-$ with  $\JB_{\frac{1}{\lambda}}$. 
									\end{enumerate}
									In particular, the relevant matrices in these cases are given as follows.
									\begin{align*}
										M_0 : \
										\begin{bNiceArray}{cc}[first-row, last-col, columns-width=auto]
											\xf_0^+ & \xf_0^-&  \\
											\Id & \JB_{\lambda} &  \xe'\\
											0 & \Id  &  \yy^{2 q} \\
										\end{bNiceArray}
										\quad \text{\large respectively} \quad
										M_0 : \ 
										\begin{bNiceArray}{cc}[first-row, last-col, columns-width=auto]
											\xf_0^+ & \xf_0^- &  \\
											\JB^{-1}_{\frac{1}{\lambda}} & 0  & \yy^{2p}\\
											\Id & \Id &  \xe'' 
										\end{bNiceArray}
										\cong
										\begin{bNiceArray}{cc}[first-row, last-col, columns-width=auto]
											\xf_0^+ & \xf_0^-  & \\
											\Id & 0  & \yy^{2p}\\
											\JB_{\frac{1}{\lambda}} & \Id &  \xe'' 
										\end{bNiceArray}
									\end{align*}
									In these cases, the claim follows from the computation
									\begin{align*}
										\begin{bNiceArray}{cc}[first-row, last-col, columns-width=auto]
											\xe' & \xe'' &  \\
											\Id & \alpha_i &  \xf_0^+\\
											\ol{\alpha}_i & \Id & \xf_0^- 
										\end{bNiceArray}
										\cdot
										\begin{bNiceArray}{cc}[first-row, last-col, columns-width=auto]
											\xx^{2p} & \xx^{2 q}  & \\
											t^{p} & 0 &  \yy^{2p}\\
											0 & t^{q} &  \yy^{2 q} 
										\end{bNiceArray}
										\cdot
										\begin{bNiceArray}{cc}[first-row, last-col, columns-width=auto]
											\xf_1^- & \xf_1^+ &  \\
											\Id & 0 &  \xe'\\
											0 & \Id &  \xe'' 
										\end{bNiceArray}
										= 
										\begin{bNiceArray}{cc}[first-row, last-col, columns-width=auto]
											\xf_1^- & \xf_1^+  & \\
											t^{p} &  \alpha_i t^q  & \xf_0^+\\
											\ol{\alpha}_i t^{p} & t^{q}  & \xf_0^- 
										\end{bNiceArray}
									\end{align*}
									\qedhere
								\end{enumerate}
							\end{proof}

						\begin{prp}\label{prp:gluing}
							Let $\omega$ be a band or string of $\chfXA$. 
							Then the pullback complex $\CPr(\omega)$ is isomorphic to the complex $\CP'(\omega)$.
						\end{prp}
						\begin{proof}
							Let $i \in \Z$. The brutal truncation $\varsigma_i(\CPr') = (\begin{td} \PPr'_i \ar{r}{\partial_i} \&  \PPr'_{i-1} \end{td})$ of the complex 
							$\CPr' = \CPr'(\omega)$ at degrees $i$ and $i-1$
							corresponds to a subdiagram of the gluing diagram of $\omega$
							given by projectives at degrees $i$ and $i-1$.
							This subdiagram is given by the union of $n_i \geq 1$ connected subdiagrams, which gives rise to 
							decompositions $\varsigma_i(\CPr')=	\bigoplus_{j=1}^{n_i} \CPr'{(j)}$ and $\partial_i = \bigoplus_{j=1}^{n_i} \partial'_i(j)$
							of certain (possibly decomposable) two-term complexes $\CPr'{(j)}$.
							Each summand 
							$\CPr'{(j)}$
							has an underlying gluing diagram.
							
							\begin{itemize}
								\item Assume that $\omega$ is a string.
								Note the gluing diagram of $\CPr'{(j)}$
								can be obtained from a
								gluing diagram of the form~\eqref{eq:rnd-ornt} by relabeling of vertices or iterated removal of end vertices and their incident arrows.
								Each differential
								$\partial_i(j)$ 
								is given by a certain composition
								$\hat{M}_{i-1}^{-1}(j) \,\wt{\partial}(j)\, \hat{M}_i(j)$.
								The factors of this composition can be obtained by appropriate relabeling or removal of rows of columns of the 
								matrices
								$\hat{M}_{i-1}^{-1}$, $\wt{\partial}$,
								and $\hat{M}_i$ with appropriate coefficients 
								from the proof of Lemma~\ref{lem:spc-cpx}.
								In particular,
								each differential $\partial'_i(j)$
								can be obtained from the matrix in \eqref{rnd:diff}
								by relabeling or removal of rows or columns, where the power of any monomial of the form $t^{n}$ 
								corresponding to an arrow
								$\begin{gd} \bt \ar{r}{n}  \& \bt
								\end{gd}$
								is allowed to be zero  
								if that arrow was changed into $\begin{gd} \bt \ar{r}{n} \& \star \end{gd}$.
								Together with Lemma~\ref{lem:spc-cpx}, it follows that each map $\partial_i(j)$ can be obtained by application of the gluing rules to the connected subdiagram underlying $\CPr'(j)$.
								\item Assume that $\omega$ is a band. 
								Then the underlying diagram of $\varsigma_i(\CPr')$ is either closed and has width two, or it is a union of open gluing diagrams as in the string case.
								In any case, one of the vertices of the diagram is decorated by $\lambda$, 
								and Lemma~\ref{lem:spc-cpx} together with Lemma~\ref{lem:spc-band}
								imply that the differential $\partial_i$ can be obtained via gluing rules.
							\end{itemize}
							In each case, the main claim follows.
						\end{proof}

						\subsubsection{Proof of the gluing theorem}
						Next, we recall and prove  
						Theorem~\ref{thm:bijection2} which is formulated more briefly below.
						\newcounter{backup}
						\setcounter{backup}{\value{equation}}
						\renewcommand{\theequation}{\ref{thm:bijection2}'}
						\begin{thm} \label{thm:bijection2'}
								The assignment $\begin{td} \omega \ar[mapsto]{r} \& \CP(\omega) \end{td}$ yields a bijection
								\begin{align}
									\label{eq:bij-perf-fd3}
									\begin{cd}
										\{\text{bands and strings of }\chfXA  \}\!\big/_{\textstyle \approx}
										\ar[yshift=0pt]{r}{\sim}\&
										\ind \Db{\Acat} 
									\end{cd}
								\end{align}
						\end{thm}
						\setcounter{equation}{\value{backup}}
						\renewcommand{\theequation}{\thesection.\arabic{equation}}
						\begin{proof}
							An application of Theorem~\ref{thm:Bondarenko} to $\fXA$ combined with Corollary~\ref{cor:reduction} yields bijections
							\begin{align*} 
								\begin{td}
									\{ \text{bands and strings of }{\fXA}  \}\!\big/_{\textstyle \approx} \ar{r}{\sim} \&
									\ind \regrep{\fX} \ar{r}{\sim} \&  \ind\DbRep{\A}\backslash
									\ind 
									\Hotbfd(\add P_{\star})\end{td}
							\end{align*}
							Together with Proposition~\ref{prp:bijection}, this establishes all bijections
							in diagram~\eqref{eq:bij-diag}.
							The composition of these bijections coincides with the one in question by Proposition~\ref{prp:gluing} and Lemma~\ref{lem:sgns}.
						\end{proof}
						\begin{rmk}
							For an arbitrary field $\kk$, band complexes are defined as in Remark~\ref{rmk:field}.
							The bijection~\eqref{eq:bij-perf-fd3} then follows by the same proof, replacing Bondarenko's Theorem~\ref{thm:Bondarenko} with its generalization from
							Remark~\ref{rmk:field2}.
						\end{rmk}
						
						For later use, we record an alternative presentation of symmetric band complexes.
						\begin{lem}\label{lem:sym-bd}
							For any symmetric band $\omega = (d,w,m,\lambda)$ of $\chfXA$,
							the band complex $\CX(\omega)$ is isomorphic to the complex $\CX^{\uparrow}(\omega')$ for the band $\omega' = (d,w,m,\frac{1}{\lambda})$, which defined in the same way as $\CX(\omega')$ however setting the opposite orientation $\uparrow$ at each axis of symmetry in $w^{\updownarrow}$.
						\end{lem}
						\begin{proof}
							Using the expanded notation from Subsection~\ref{subsec:exp-bnd}
							and~\eqref{eq:ornt-exp},
							this statement follows from  Remark~\ref{rmk:sym-isom}
							for the matrix representation $\mu(\omega)$ underlying $\PPr(\omega) \cong \CX(\omega)$. 
						\end{proof}

\section{Action of natural functors on band and string complexes}
\label{sec:derived_fun}

The derived category $\Db{\Acat}$
of nilpotent Gelfand quiver representations
admits three natural auto-equivalences:
the involution $\inv$, the Auslander-Reiten translation $\tau$ and the 
contragredient duality $\cdual{(-)}$.
In the present section, we 
express these functors
in terms of 
basic operations on the gluing diagram of a band or string $\omega$ of $\Db{\Acat}$ such as
the half-turn, the transpose of each arrow, and the replacement of each projective $P_{\star}$ by the gluing diagram of the projective resolution of $I_{\star}$.
The diagrammatic interpretation allows
to describe the action of the three natural auto-equivalences $\tau$, $\inv$ and $\cdual{(-)}$
on bands and strings of $\Db{\Acat}$.
Moreover, we give a homological characterization of the four classes of band and string complexes.
The statements of the results of this section presume only familiarity with the notions of bands and strings from Section~\ref{sec:complex}.

\subsection{The derived involution}

We recall that
the involution functor $\inv$ 
is 
uniquely determined by the prescription on the right
\begin{align*}
	\begin{td}
		\inv \colon 
		\DbRep{\A} \ar{r}{\sim} \& \DbRep{\A}
	\end{td}
	&&
	\inv(\begin{cd} \underset{d+1}{P_u} \ar{r}{t^n} \& \underset{d}{P_v} \end{cd})
	=(\begin{td} \underset{d+1}{P_{\invaut(u)}} \ar{r}{t^n} \& \underset{d}{P_{\invaut(v)}} \end{td})
\end{align*}
for any integer $d \in \Z$, vertices $u,v \in \{+,\star,-\}$ and number
$n \in \N$ respectively $n \in \N_0$ if $(u,v) = (\pm,\star)$, where $\invaut(\pm) = \mp$ and $\invaut(\star) = \star$. 
The functor $\inv$ is induced by a $\kk$-algebra involution of $\A$, which was introduced in \eqref{E:Involution} and which has a Lie-theoretic origin.
Given the construction of band and string complexes, there is 
the following straightforward description of the involution functor $\inv$ in terms of bands and strings.

\begin{prp}\label{prp:inv}
	Let $\omega$ be a band or string of $\chfXA$ and $\CP(\omega)$ its corresponding indecomposable complex.
	Then $\inv(\CP(\omega)) \cong \CP(\inv(\omega))$, where $\inv(\omega)$ is defined by the following table.
	\begin{align*}
		\begin{array}{|c|c|c|c|c|}
			\hline
			& \text{usual string} & \text{special string} & \text{bispecial string}  & \text{band} \\
			\omega & (d,w) & (d,w,\varepsilon) & (d,\varepsilon_1, w, \varepsilon_2) & (d,w,m,\lambda) \\ \hline
			\inv(\omega) & 
			(d,w) & (d,w,\ol{\varepsilon}) &
			(d,\ol{\varepsilon}_1, w, \ol{\varepsilon}_2) & (d,w,m,\lambda) \\
			\hline
		\end{array}
	\end{align*}
\end{prp}
\begin{proof}
	Let $\omega$ be a band or string of $\chfXA$. 
	We claim that
	$\inv(\CX(\omega)) \cong \CX(\omega_\inv^{\op})$
	for $\omega_\inv^{\op} \colonequals (\inv(\omega))^{\op}$.
	\begin{itemize}
		\item Assume that $\omega$ is a usual string, a special string or an asymmetric band. Switching the signs in all vertices of the final gluing diagram of $\CX(\omega)$
		and rotating the diagram by $180^\circ$ yields precisely the final gluing diagram underlying $\CX(\omega^{\op})$.
		This shows that there is an isomorphism $\inv(\CX(\omega)) \cong \CX(\omega^{\op}_{\inv})$ given by reordering of projectives.
		\item If $\omega$ is a bispecial string $(d,\varepsilon_1, w,\varepsilon_2)$  or a symmetric band $(d,w,m,\lambda)$, 
		we may adapt the previous argument using Remark~\ref{rmk:sym-isom}.
	\end{itemize}
	In both cases, the claim follows.
	It holds that $\CX(\omega_\inv^{\op}) \cong \CX(\sigma(\omega))$ by 
	Theorem~\ref{thm:bijection2} and Definition~\ref{def:equiv}.
\end{proof}

\subsection{Homological characterization of bands and strings}
The previous result allows to deduce a purely homological characterization 
of the four combinatorially defined classes of bands and strings in terms of 
the involution and the defect. 
As introduced in \eqref{E:DefectComplex}, the defect of a complex $\CX$ from $\DbRep{\A}$ is given by
$$
\defect(\CX) = \sum_{p \in \Z} \dim \Hom_{\D(\A)}(\CX,S_{\star}[p]).
$$
Similar to the involution,
the defect has an interpretation in Lie-theoretic terms \ref{rmk:def-lie}.

\begin{thm}	\label{thm:4cl}
	Let $\omega$ be a band or string of $\chfXA$ and $\CP = \CP({\omega})$ its corresponding indecomposable complex.
	Then  the following equivalences hold.
	\begin{align*}
		\begin{array}{lclclclcl}
			\omega\text{ is a usual string} 
			& \Leftrightarrow & \defect(\CP) = 2 &\Leftrightarrow& \inv(\CP) \cong \CP &\text{and}&\CP\text{ is not $\tau$-periodic}.\\
			\omega\text{ is a special string} 
			& \Leftrightarrow & \defect(\CP) = 1 &\Leftrightarrow& \inv(\CP) \not\cong \CP &\text{and}&\CP\text{ is not $\tau$-periodic}.\\
			\omega\text{ is a bispecial string} & \Leftrightarrow &
			\defect(\CP) = 0&\text{and}& \inv(\CP) \not\cong \CP  
			&\Leftrightarrow&	
			\CP \text{ has $\tau$-period two}.
			\\
			\omega\text{ is a band} & \Leftrightarrow &
			\defect(\CP) = 0&\text{and}& \inv(\CP) \cong \CP  &\Leftrightarrow& 
			\CP \text{ is $\tau$-invariant}.
		\end{array}
	\end{align*} 
\end{thm}
\begin{proof}
	By the definitions, $\defect(\CP)$ is given by the number of non-special ends in the word underlying the datum $\omega$. Therefore,
	$\omega$
	is a usual string if and only if $\defect(\CP) = 2$, 
	a special string if and only if $\defect(\CP) = 1$,
	and a bispecial string or a band if and only if $\defect(\CP) = 0$.
	The equivalences follow then from Propositions~\ref{prp:inv} and~\ref{prp:tau-per}.
\end{proof}
We note 
that any complex $\CX$ of $\DbRep{\A}$ with defect greater than two must be decomposable.

\begin{ex}
	We reconsider the family of low-dimensional indecomposable representations $(M_{\lambda})_{\lambda \in \kk}$ with $M_{\lambda}$ given by \eqref{eq:sh-band}. 
	The observations below \eqref{eq:sh-band} yield
	that
	$\delta(M_{\lambda}) \neq 0$ if and only if $\lambda = 1$,
	and that
	$\inv(M_{\lambda}) \not\cong M_{\lambda}$ if and only if $\lambda = 0$.
	
	By the last theorem, the projective resolution 
	$\CX(M_\lambda)$ of $M_\lambda$ with $\lambda \in \kk \backslash \{0,1\}$ should be isomorphic to a band complex,
	$\CX(M_1)$ to a usual string complex,
	and
	$\CX(M_0)$ to a bispecial string complex.
	For any $\lambda \in \kk$ a minimal projective resolution of $M_{\lambda}$
	is given by 
	\begin{align}\label{eq:bd-res}
		\CX(M_\lambda) = 
		\begin{cases}
			\begin{array}{rl}
				\begin{cd} P_- \oplus P_+ \ar{r}{
						\left(
						\begin{smallmatrix}
							t& -t \\
							-\lambda t & t		
						\end{smallmatrix}\right)
					} \& P_+ \oplus P_- \end{cd} &
				\text{if }\lambda \neq 1,\\
				\begin{cd} P_\star \ar{r}{
						\left(\begin{smallmatrix}		
							0 \\
							-t \\
							t
						\end{smallmatrix}\right)
					} \&
					P_\star \oplus P_- \oplus P_+
					\ar{r}{
						\left(
						\begin{smallmatrix}
							0 & t& -t \\
							t^2&-t & t		
						\end{smallmatrix}\right)		
					} \&
					P_+ \oplus P_-
				\end{cd} & 
				\text{if }\lambda = 1.
			\end{array}
		\end{cases}
	\end{align}
	Considering the projective resolutions, it follows 
	that there are isomorphisms
	\begin{enumerate}
		\item $\CX(M_1) \cong \CP(\omega)$ for the usual string $\omega = (2,w)$ with $w
		= (\star,\ora{1},\ora{1},\ola{2},\star)$,
		\item 
		$\CX(M_0) \cong \CP(\omega)$ for the bispecial string $(1,-,w,+)$
		with $w = (\ora{1},\ola{1})$, and
		\item $\CX(M_{\lambda}) \cong \CP(\omega)$ for the band $
		\omega = (1,w,m,\lambda')$ with $m=1$, and $\lambda' = \frac{\lambda-1}{\lambda}$.
	\end{enumerate}
	This confirms the prediction of the theorem.
\end{ex}

\subsection{Resolving radicals and the Auslander-Reiten translation}
\label{subsec:res-rad}

Assume that $w$ is a finite word of $\chfXA$, that is, 
$w = (\alpha, x_1, x_2, \ldots x_{\ell}, \beta)$, $\alpha,\beta \in \{\star,\diamond \}$ and $x_i = \ola{n}_i$ or $\ora{n}_i$ for each $1 \leq i \leq \ell$.
For a number $n \in \N$ we denote $n^+ = n + 1$
and $n^- = n-1$ in what follows. 
\begin{dfn} \label{dfn:rad-words}
	In the notations above,we define sequences $\ell(w)$, $r(w)$ and a finite word $w_\star$ of $\chfXA$ as follows.
	\begin{enumerate}
		\item 
		If $\alpha \neq \star$, set $\ell(w) = w$. Otherwise, $\ell(w)$ is defined according to the table
		\begin{align*}
			\begin{array}{|r||r|rl|r|r|r|r|}
				\hline
				w &  (\star, \ora{n}_1, v) & (\star, {\color{blue}\ola{n}_1}, v) & n_1 > 0 & (\star, {\color{blue} \ola{0}, \ora{n}_2}, v) & (\star, {\color{blue} \ola{0}}, \ola{n}_2, v) & (\star, {\color{blue}\ola{0}},\varepsilon)\\ 
				\hline
				\ell(w) & 
				(\star, {\color{blue} \ora{1} }, \ora{n}_1, v) &
				(\star, {\color{blue} \ora{1}}, \ola{n}_1, v)& &
				(\star, { \ora{n}^+_2}, v) &
				(\star, \ola{n}^-_2,v) & 
				(\star, {\color{blue}\ora{1}},\varepsilon)
				\\ \hline 
			\end{array}
		\end{align*}
		\item We set $r(w) = (\ell(w^{\op}))^{\op}$, that is, $r(w)$ is defined via rules dual to $\ell(w)$.
		\item \label{tau-w}
		We set $w_{\star} = \ell(r(w))$, which can be verified to coincide with $r(\ell(w))$.
	\end{enumerate}
\end{dfn}
The gluing diagrams of $\ell(w)$ with the assumptions as above 
are given by the following.
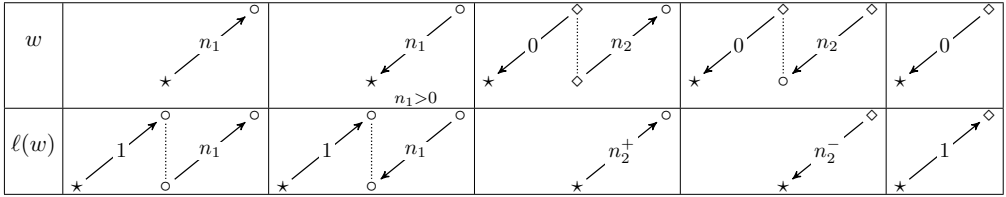
\begin{align*}
	\scalebox{0.8}{$
		{
			\setlength{\arraycolsep}{3pt}
			\begin{array}{|c|r|r|r|r|r|}
				\hline
				w&
				\begin{gd}
					\phantom{\bt} \& \bt \\
					\star \ar{ru}[description]{n_1} \& 
				\end{gd} 
				&
				\underset{n_1 > 0}{
					\begin{gd}
						\&  \bt  \\
						\star \ar[<-]{ru}[description]{n_1}\&
				\end{gd}  }
				&
				\begin{gd}
					\&	\diamond \ar[-,densely dotted]{d}  \& \bt \\
					\star \ar[<-]{ru}[description]{0} \&	\diamond \ar{ru}[description]{n_2} \& 
				\end{gd} &
				\begin{gd}
					\&	\diamond \ar[-,densely dotted]{d}  \&  \diamond \\
					\star \ar[<-]{ru}[description]{0} \&	\bt \ar[<-]{ru}[description]{n_2} 
				\end{gd}
				&
				\begin{gd}
					\& \diamond\\
					\star \ar[<-]{ru}[description]{0} \&
				\end{gd}
				\\
				\hline
				\ell(w) &
				\begin{gd}
					\& \bt \ar[-,densely dotted]{d} \& \bt \\
					\star \ar{ru}[description]{1} \&	\bt \ar{ru}[description]{n_1} \& 
				\end{gd} 
				&
				\begin{gd}
					\&  \bt \ar[-,densely dotted]{d} \& \bt \\
					\star \ar{ru}[description]{1}\& \bt \ar[<-]{ru}[description]{n_1}
				\end{gd}  
				&
				\begin{gd}
					\&	  \& \bt \\
					\&	\star \ar{ru}[description]{n_2^+} \& 
				\end{gd} &
				\begin{gd}
					\&	  \&  \diamond \\
					\&	\star \ar[<-]{ru}[description]{n_2^-} 
				\end{gd}
				&
				\begin{gd}
					\& \diamond\\
					\star \ar{ru}[description]{1} \&
				\end{gd} \\ 		\hline
		\end{array}}
		$}
\end{align*}

\begin{enumerate}
	\item We define $w^{\vee}$ by reversing the oriented arrow above each number in the word $w_{\star}$.
\end{enumerate}
\begin{rmk}
	The word $w = (\star,\ola{0},\ola{1},\star)$ is the only case when 
	the sequence $\ell(w)$ is not a word of $\chfXA$. 
	In this case, $\ell(w) = (\star, \ola{0},\star)$, but $w_{\star}$ is still defined and equal to $w$ itself.
	Similar considerations apply to $r(w^{\op})$.
\end{rmk}

\begin{dfn}
	We will call a special word of the form $(\star,\ola{0},\diamond)$
	or $(\diamond,\ora{0},\star)$
	\emph{exceptional}.
\end{dfn}

\begin{lem}\label{lem:exc-T}
	For any band or string $\omega$ of $\chfXA$ the following conditions are equivalent.
	\begin{enumerate}
		\item \label{sp-emb1}
		$\omega$ is a special string $(d, w,\varepsilon)$
		for some integer $d \in \Z$, an exceptional word $w$ and some sign $\varepsilon$.
		\item \label{sp-emb2}
		$\CX(\omega) \cong (\begin{td} P_{\varepsilon} \ar{r}{\iota}\& P_{\star} \end{td}) [d]$
		for some integer $d \in \Z$ and sign $\varepsilon$.
		\item \label{sp-emb3} $\tau(\CX(\omega)) \cong S_\varepsilon [d]$ for some integer $d \in \Z$ and sign $\varepsilon$.
	\end{enumerate}
\end{lem}
\begin{proof}
	The equivalence \eqref{sp-emb1} $\Leftrightarrow$ \eqref{sp-emb2} holds true by the construction of $\CX(\omega)$.
	The equivalence
	$\eqref{sp-emb2} \Leftrightarrow \eqref{sp-emb3}$ follows from
	$
	\tau((\begin{td} {P_{\pm}} \ar{r}{\iota}\& {P_{\star}} \end{td})[d])
	\cong 
	(\begin{td} {P_{\mp}} \ar{r}{\iota}\& {I_{\star}} \end{td})[d]
	\cong
	S_{\pm}[d]
	$.
\end{proof}
The remaining goal of this subsection is to show the following result.
\begin{thm}\label{thm:tau}
	Let $\omega$ be a band or string of $\chfXA$.
	Then there are isomorphisms of complexes 
	$\Tw(\CP(\omega)) \cong \CP(\omega_{\star})$
	and
	$\tau(\CP(\omega)) \cong \CP(\tau(\omega))$
	with $\omega_{\star}$ and $\tau(\omega)$ defined as follows.
	\begin{itemize}
		\item If $\omega$ 
		is not a special string with an exceptional word,
		then $\omega_{\star}$ and $\tau(\omega)$
		are defined by the table
		\begin{align}
			\label{tab:tau/dual}
			\begin{array}{|c|c|c|c|c|}
				\hline
				& \text{usual string} & \text{special string} & \text{bispecial string} & \text{band} \\
				\omega & (d,w) & (d,w,\varepsilon) & (d,\varepsilon_1, w, \varepsilon_2) & (d,w,m,\lambda) \\ \hline
				\omega_{\star} & 
				(d_{\star},w_{\star}) & (d_{\star},w_{\star}, {\varepsilon}) &
				(d,{\varepsilon}_1,w, {\varepsilon}_2) & (d,w,m,\lambda) 
				\\
				\hline 
				\tau(\omega) & 
				(d_{\star},w_{\star}) & (d_{\star},w_{\star}, \ol{\varepsilon}) &
				(d,\ol{\varepsilon}_1,w, \ol{\varepsilon}_2) & (d,w,m,\lambda) 
				\\
				\hline
			\end{array}
		\end{align}
		where $d_{\star } = d+1$ if $w$ begins with $\star$ and
		$d_{\star} = d $ otherwise, 
		$w_{\star}$ is the word defined in \eqref{tau-w} above.
		\item
		For a special string $(d,w,\varepsilon)$ with an exceptional word $w$ 
		it holds that  $\omega_{\star} = (d_\star, w_{\star},\ol{\varepsilon})$ and  $\tau(\omega)  = 
		(d_\star, w_{\star},\varepsilon)$.	
	\end{itemize}
\end{thm}
In particular, 
$\Tw$ preserves the sign data in all strings but the exceptional special strings, which motivates the terminology.
The proof of the last theorem requires a few preparations. 
The basic idea of the proof is already present in the next example.
\begin{ex}
We reconsider the usual string $\omega = (2,w)$
with 
$w = (\star,\ora{2},\ora{1},\ola{2},\ora{3},\ora{0},\star)$ from
Subsection~\ref{subsec:ex-us2}.
An application of $\Tw$ to  the usual string complex $\CX(\omega)$ 
yields a complex $\CL(\omega)$ given by the middle row below, 
where $I_{\star} = \rad P_{\star}$.
\begin{align*}
	\begin{gd}
		\CX(\omega) \&\&
		{P_\star} \ar{r}{\left(\begin{smallmatrix}
				t^{2} \\
				-t^{2} \end{smallmatrix}\right)} \& {P_{+-}} \ar{r}{M_2} \& {P_{+-}} \ar[<-]{r}{M_3} \& {P_{+-}} \ar{r}{ 			\left(\begin{smallmatrix}
				0 & t^3\\
				0 & t^3 \end{smallmatrix}\right)} \&
		{P_{+-}} \ar{r}{ 			\left(\begin{smallmatrix}
				-\iota & \iota
			\end{smallmatrix}\right)} \& {P_{\star}}
		\\
		\CL(\omega)\& 0\ar{r} \&
		I_\star \ar{r}{\left(\begin{smallmatrix}
				t^{2} \\
				-t^{2} \end{smallmatrix}\right)} \& {P_{+-}} \ar{r}{M_2} \& {P_{+-}} \ar[<-]{r}{M_3} \& {P_{+-}} \ar{r}{ 			\left(\begin{smallmatrix}
				0 & t^3\\
				0 & t^3 \end{smallmatrix}\right)} \&
		{P_{+-}} \ar{r}{ 			\left(\begin{smallmatrix}
				-\iota & \iota
			\end{smallmatrix}\right)} \& {I_{\star}}
		\\
		\CX' \ar[dashed]{u} \&
		P_\star \ar{u} \ar{r}{\left(\begin{smallmatrix}
				t \\
				t \end{smallmatrix}\right)}
		\&
		{P_{+-}} \ar{r}{\left(\begin{smallmatrix}
				-t^{2} &  t^2 \\
				t^{2} & -t^2 \end{smallmatrix}\right)} 
		\ar{u}{
			\left(\begin{smallmatrix} -\iota & \iota \end{smallmatrix}\right)	
		}
		\& {P_{+-}} 
		\ar{u}[swap]{\id}
		\ar{r}{M_2} \& {P_{+-}} \ar{u}{\id} \ar[<-]{r}{M_3} \& {P_{+-}} \ar{u}[swap]{\id} \ar{r}{ 			\left(\begin{smallmatrix}
				0 & t^2 \end{smallmatrix}\right)} \&
		{P_{\star}} \ar{u}[swap]{\left(\begin{smallmatrix} t \\ t \end{smallmatrix}\right)} \ar{r} \& 0 \ar{u}
	\end{gd} 
\end{align*}			
`Folding' the bottom diagram yields a quasi-isomorphism of complexes 
$\begin{td} \CL(\omega) \ar{r}\& \CX' \end{td}$.
The crucial point is that 
the complex $\CX'$ coincides with the string complex of the 
string $\omega_{\star} = (3,w_{\star})$ 
with 
$w_{\star} = (\star,\ora{1},\ora{2},\ora{1},\ola{2},\ora{3},\star)$.
These observations imply that $\Tw(\CX(\omega)) \cong \CX(\omega_{\star})$ as claimed in Theorem~\ref{thm:tau}.
\end{ex}

\begin{lem}\label{lem:qis}
Let $\CX'$ and $\CL$ be the complexes given by folding
the top respectively bottom row of a commutative diagram
of $\A$-modules
$$
\begin{gd}
	P_2 \ar[hookrightarrow]{r}{M_0
	} \ar{d} \& P_1  \ar{r}{M_1
	} 
	\ar[->>]{d}{
		\pi
}
\& P_0 \ar[equal]{d} \ar[<->]{r}{M_2} \& X(3) \ar[equal]{d}
\ar[densely dotted,-]{rr} \& \& X({\ell}) \ar[equal]{d} \ar[<->]{r}{M_\ell} \& X(\ell+1) \ar[equal]{d}
\\
0 \ar{r}
\&
L  \ar[hookrightarrow]{r}{
	M'_1	
} \& P_0 \ar[<->]{r}{M_2} \& X(3)
\ar[densely dotted,-]{rr}  \& \& X({\ell}) \ar[<->]{r}{M_{\ell}} \& X({\ell+1}) 
\end{gd}
$$
in which $\im M_0 = \ker \pi$, 
each pair of arrows below a matrix $M_2, M_3, \ldots, M_\ell$ is oriented 
in the same direction,
and any well-defined composition of horizontal arrows is zero.
Then folding the vertical arrows yields a quasi-isomorphism of complexes $
\begin{td} \CX' \ar[twoheadrightarrow]{r} \& \CL \end{td}$.
\end{lem}
\begin{proof}Folding yields the morphism of complexes
$$
\begin{cd}
\mathstrut \ar[densely dotted,-]{r}\&
X_3 \ar[equal]{d}
\ar{r}{
	\left(
	\begin{smallmatrix}
		0
		\\
		\partial'_3
	\end{smallmatrix}
	\right)	
} \& 
P_2 \oplus X'_2 
\ar[->>]{d}[swap]{\left(\begin{smallmatrix}
		0 & \id 
	\end{smallmatrix}
	\right)
} 
\ar{r}{
	\left(
	\begin{smallmatrix}
		M_0 & 0 \\
		0 & \partial'_2
	\end{smallmatrix}
	\right)
}
\&
P_1 \oplus X'_1 \ar[->>]{d}[swap]{
	\left(
	\begin{smallmatrix}
		\pi & 0 \\
		0 & \id
	\end{smallmatrix}
	\right)
}
\ar{r}{
	\left(
	\begin{smallmatrix}
		M_1 & \partial'_1 \\
		0 & \partial''_1
	\end{smallmatrix}
	\right)
}
\&
P_0 \oplus X'_0 \ar[equal]{d}
\ar{r}{\partial_0}\&
X_{-1} \ar[equal]{d} \ar[densely dotted,-]{r}\& \mathstrut
\\
\mathstrut \ar[densely dotted,-]{r} \&
X_3 \ar{r}{\partial'_3} \&  X'_2
\ar{r}{
	\left(
	\begin{smallmatrix}
		0 \\ \partial'_2
	\end{smallmatrix}
	\right)	
} \&
L \oplus X_1'
\ar{r}{
	\left(
	\begin{smallmatrix}
		M'_1 & \partial'_1 \\
		0 & \partial''_1
	\end{smallmatrix}
	\right)
}
\&
P_0 \oplus
X_0' \ar{r}{\partial_0} \& X_{-1}
\ar[densely dotted,-]{r}\& \mathstrut
\end{cd}
$$
which can be checked to be a quasi-isomorphism.
\end{proof}

\begin{lem}\label{lem:rad-res}
Let $\omega$ be a usual or special string of $\chfXA$. Set $\CL(\omega) \colonequals \Tw(\CX(\omega))$.
Let $\CX'(\omega)$ denote the complex obtained from the string complex $\CX(\omega)$
by replacing each end of the form $P_{\star}$ in its gluing diagram with the projective resolution of $\rad P_\star$. 
Then there is a quasi-isomorphism  $\begin{td} \CX'(\omega) \ar[twoheadrightarrow]{r} \& \CL(\omega) \end{td}$ of complexes
of $\A$-modules.
\end{lem}
\begin{proof}
Assume that $\omega$ is a usual string with an underlying word of length $\ell$.
We depict the oriented gluing diagram of $\CX(\omega)$  in the vertical notation
$$
\begin{gd}
\CX(\omega) \&
\mathstrut \& 
P_\star  \ar[<->]{r}{M_1} \& 
P_{+-}  \ar[<->]{r}{M_2} \& P_{+-} \ar[densely dotted,-]{r}  \& P_{+-} \ar[<->]{r}{M_\ell}  \& P_\star \& \mathstrut
\end{gd}
$$
where the arrow below $M_i$ is oriented to the left if $x_i$ is left-oriented, and otherwise oriented to the right for each index $1 \leq i\leq \ell$.

For the purposes of this proof, it will convenient 
to depict the gluing diagram in vertical notation, described by top of the next diagram.
\begin{align*}
\begin{gd}
	\CX'(\omega) \ar[dashed]{d} \&
	P_{\star} \ar{d} \ar{r}{\left(\begin{smallmatrix}
			\phantom{-}t \\
			-t
		\end{smallmatrix}\right)} \& 
	P_{+-} \ar[twoheadrightarrow]{d}[swap]{\pi} \ar[<->]{r}{M'_1} \& 
	P_{+-} \ar[equal]{d} \ar[densely dotted,-]{rr} \& \& P_{+-} \ar[<->]{r}{M'_\ell} \ar[equal]{d}  \& P_{+-}  \ar[<-]{r}{\left(\begin{smallmatrix}
			\phantom{-}t \\
			-t
		\end{smallmatrix}\right)} \ar[twoheadrightarrow]{d}{\pi} \& P_{\star} \ar{d} \\
	\CL(\omega) \& 0 \ar{r} \& I_{\star}  \ar[<->]{r}{M_{1}} \& P_{+-} 
	\ar[densely dotted,-]{rr} \&  \& P_{+-} \ar[<->]{r}{M_{\ell}} \& I_{\star} \ar[<-]{r} \& 0
\end{gd}
\end{align*}
If the arrow below $M_1$ is oriented to the right, then $M'_1$ is defined as the composition $M_1 \circ \pi$. If the arrow is left-oriented and $M_1 = \begin{psmallmatrix}A & B \end{psmallmatrix}$, we set $M'_1 = A \oplus B$.
Similar considerations are carried out for the matrix $M'_\ell$.

Since each square commutes, `folding' both rows together with the vertical morphisms yields a morphism of complexes $			\begin{td}\CX'(\omega)\ar{r} \&   
\CL(\omega) \end{td}$. 
In case the arrows below $M_1$ are right-oriented and the arrows below $M_{\ell}$ are left-oriented, repeated applications of Lemma~\ref{lem:qis}
yield that the morphism is a quasi-isomorphism.
By similar 
arguments for the other orientations of the arrows below $M_1$ and $M_{\ell}$, it can be verified that this morphism is a quasi-isomorphism. 

The case of a special string $\omega$ requires resolving only one radical and is therefore similar and simpler.
\end{proof}

Next, we need to show that the  diagram obtained by 
the visual changes described in   Subsection~\ref{subsec:res-rad}
to the ends of the oriented gluing diagram of $\omega$ results
indeed in the oriented gluing diagram of  $\omega_{\star}$.
\begin{prp}\label{prp:orient-tau}
Let $\omega$ be a usual string $(d,w)$ or a special string $(d,w,\varepsilon)$.
Let $\omega_{\star}$ be the datum defined in Table~\ref{tab:tau/dual}. 
Then each common gluing edge of the diagrams of the words $w$ and $w_{\star}$ 
is oriented in the same direction.
\end{prp}
\begin{proof}
Assume that $w$ and $w_{\star}$ are given in their expanded form as introduced in Subsection~\ref{subsec:exp-str}, the expansion of the first word denoted again by $w$, and the second by $w'$ for simplicity.
Let $\tilde{w}$ denote the ambient word of $w$ and $\tilde{w}'$ that of $w'$. 
By definition the ends of $\tilde{w}$ and $\tilde{w}'$ are given by elements from $\fE$. 
In this proof, for elements $e'$, $e \in \fE$ we write $e' \prec e$ if $e$ is not a minimal element in $\fE$ and $e'$ is the maximal element with $e' < e$. We write $e' \preceq e$ if $e' \prec e$ or $e' = e$.
Moreover, we call an end of $\tilde{w}$ \emph{extending}
if it is given by a non-minimal element of $\fE$ and \emph{contracting} otherwise.
The terminology is motivated by the following observations on the form of $w'$. 
\begin{align*}
\begin{array}{|r||r|r|}
	\hline
	& \multicolumn{1}{c|}{\text{extending left end} } & \multicolumn{1}{c|}{\text{contracting left end}} \\ \hline \hline
	\tilde{w} & \xe -  \ldots & \xx^1_{d} - \xf_{d} \sim \xf_d - \ol{e} \sim e-
	\ldots  \\ \hline  
	\tilde{w}' & \yy^1_d - \xf_d \sim \xf_d - \ol{\xe}' \sim \xe' - \ldots & e' - \ldots\\ 
	\hline
\end{array}
\end{align*}
In both cases, it holds that $\xe' \prec \xe$. 
More precisely,  there are $n \in \N$ and $d \in \Z$
such that
$(\ol{\xe}', \xe', \xe)
= ( \yy_{d}^{2n}, \xx_{d+1}^{2n},\xx_{d+1}^{2n+1})$
or 
$(\xx_{d}^{2n}, \yy_{d-1}^{2n}, \yy_{d-1}^{2n-1}
)
$
in case the left end of $\tilde{w}$ is extending, 
respectively
$(\ol{\xe}, \xe, \xe')
= ( \yy_{d}^{2n}, \xx_{d+1}^{2n},\xx_{d+1}^{2n-1})$
or 
$(\xx_{d}^{2n}, \yy_{d-1}^{2n}, \yy_{d-1}^{2n+1}
)
$ in the contracting case.

Assume that we determine the oriented arrow above the $i$-th common special subword $\xf \sim \xf$ of $w$ and $w'$.
By Subsection~\ref{subsec:A-orient} the orientation in $w$ is determined by two letters $l_i$ and $r_i$  appearing directly before and after the maximal symmetric subword $v \sim v^{\op}$ of a certain length $2q = 2q(i) \in 2\N$ with $\xf \sim \xf$ in the middle, which is considered 
in the ambient word $\tilde{w}$.
In case of the word $w'$ the corresponding length of the maximal symmetric word shall be denoted by $2q'$, and the orientation-determining letters shall be denoted by $l'_i$ and $r'_i$, respectively.
These letters are well-defined, comparable and distinct, since $w'$ is a usual or special string.

It is sufficient to show that either $l_i < r_i$ and $l'_i < r'_i$, or $l_i > r_i$ and $l'_i > r'_i$, which is done by distinguishing the following cases.

\begin{enumerate}
\item Assume that $q = q'$. Equivalently, $l'_i$ and $l_i$ are located at the same index, and $r'_i$ and $r_i$ share the same index was well.
In particular, it holds that $l'_i \preceq l_i$ and $r'_i \preceq r_i$.
\begin{itemize}
	\item If $l_i < r_i$, it holds that $l'_i \preceq l_i \leq r'_i \preceq r_i$.
	Since $l'_i \neq r'_i$, thus $l'_i < r'_i$.	
	\item If $l_i > r_i$ it holds that $l_i \succ l'_i \geq r_i \succeq r'_i$, and similarly as above, $l'_i \neq r'_i$ implies that $l'_i > r'_i$.
\end{itemize}
\item Otherwise it holds that $q \neq q'$. This may happen only if $l_i$ is the left end of $\tilde{w}$ and $l'_i$ is the left end of $\tilde{w}'$.
We distinguish further between the type of end of $\tilde{w}$.
\begin{enumerate}
	\item	
	\label{L2}
	Assume that the left end $l_i$ of $\tilde{w}$ is extending.
	Then $l'_i$ is a maximal element of $\fE$, thus $l'_i > r'_i$, and the following situation occurs.
	\begin{align*}
		\begin{array}{rl}
			\tilde{w} &  \phantom{l'_i - f_d \sim f_d - \ol{e}' \sim \ } l_i - v \sim v^{\op} - r_i \ldots \\
			\tilde{w}' & l'_i - f_d \sim f_d - \ol{e}' \sim e' - v \sim v^{\op} - e' \sim \ol{e}' - f_d \sim f_d -r_i'\ldots		\end{array}
	\end{align*}
	Since $l_i' \neq r'_i$,  it follows that $r_i$ is not an extending end of $\tilde{w}$. Therefore, $r_i$ is given by the letter in $\tilde{w}'$ at the same position, that is, $r_i = e'$.
	So $l_i = e \succ e' = r_i$.
	\item \label{L3} Otherwise the left end $l_i$ of $\tilde{w}$ is contracting. Then $l_i < r_i$, and we have the situation
	\begin{align*}
		\begin{array}{rl}
			\tilde{w} &  l_i - f \sim f - \ol{e} \sim e - v \sim v^{\op} - e \sim \ol{e} - f \sim f -r_i \ldots		 \\
			\tilde{w}' & \phantom{l_i - f \sim f - \ol{e} \sim} l'_i - v \sim v^{\op} - r'_i \ldots 
		\end{array}
	\end{align*}
	Since $l_i \neq r_i$, $r_i$ is not a contracting end of $\tilde{w}$, and thus $r'_i = e$.
	It follows that $l'_i = e' \prec e =  r'_i$.
\end{enumerate}
\end{enumerate}
Summarized, the above case-analysis  shows that $(l_i, r_i)$ and $(l'_i, r'_i)$ satisfy the same order relation, and thus the arrows above the special subword in $w$ and $w'$ have the same orientation. 
\end{proof}
\begin{proof}[Proof of Theorem~\ref{thm:tau}]
To show the first claim, assume that $\omega$ is a band or string of $\chfXA$.
The complex	$\CL(\omega) = \Tw(\CX(\omega))$ 
is quasi-isomorphic to the complex $\CX'(\omega)$ defined in Lemma~\ref{lem:rad-res}.
\begin{itemize}
\item If $\omega$ is a special string with an exceptional word, the corresponding claims are true by Lemma~\ref{lem:exc-T}. 
\item Excluding the previous case, the complex $\CX'(\omega)$ may be viewed as the complex obtained from an oriented gluing diagram.
The crucial observation is that this gluing diagram can be identified with the gluing diagram of the 
string $\omega_{\star}$ through the use of Proposition~\ref{prp:orient-tau}. The formula for $d_{\star}$ is obtained by inspecting the possible beginnings of $\ell(w)$ in Definition~\ref{dfn:rad-words}~\eqref{tau-w}.
This shows that $\CX'(\omega) = \CX(\omega_{\star}) \cong \CL(\omega) = \Tw(\CX(\omega))$. 
\end{itemize} 
The remaining claim concerning the Auslander-Reiten translation 
is a consequence of the first claim, combined with
the isomorphism of functors $\Tw^{-1} \cong \inv \circ \tau$ from Proposition~\ref{P:twist} and
the description of the involution on bands and strings from Proposition~\ref{prp:inv}.
\end{proof}

\subsection{The transpose}
The next result follows from work by Iyama and Reiten \cite{IR}. It will be useful to compute the derived duality $\dual{(-)}$ below.
\begin{prp}\label{prp:IR}
There are isomorphisms of functors
\begin{align*}
\RHom_{\kk}(-,\kk) \circ [-1]
&
\cong
\RHom_{R}(-,R) 
\\
&
\cong
\RHom_{\A}(-,\A) \circ \tau^{-1} 
\colon 
\begin{cd}\Dbfdmod{\A} \ar{r}{\sim} \& \Dbfdmod{\Aop}
\end{cd}.
\end{align*}
\end{prp}
\begin{proof}[Comment on the proof]
The first isomorphism follows from
\cite[Lemma 3.6]{IR} and \cite[Proof of Theorem 3.8]{IR}, 
the second from 
\cite[Proof of Proposition 3.6(2)]{IR}.
\end{proof}

We consider the contravariant equivalence 
\begin{align*}
\begin{cd}
\Dtr{-} \colonequals  \mathsf{R}\phi \circ \RHom_{\A}(-,\A) \circ [1] \colon
\DbRep{\A} \ar{r}{\sim} \& \DbRep{\A}
\end{cd}
\end{align*}
where $\phi\colon
\begin{td}\rep{\Aop} \ar{r}{\sim} \&\rep{\A}\end{td}$ denotes the equivalence from
Subsection~\ref{subsec:dual-sim}. 
In notation of paths, the functor $\Dtr{-}$
is uniquely determined by the prescriptions
\begin{align*}
(\begin{td} \underset{d+1}{P_{\pm}} \ar{r}{\cdot a_{\pm}} \& \underset{d}{P_{\star}} \end{td})
\ \overset{\Dtr{-}}{\longmapsto} \
(\begin{td} \underset{-d}{P_{\pm}} \ar[<-]{r}{b_{\pm}} \& \underset{1-d}{P_{\star}}\end{td})
&&
(\begin{td} \underset{d+1}{P_{\star}} \ar{r}{b_\pm} \& \underset{d}{P_{\pm}} \end{td})
\ \overset{\Dtr{-}}{\longmapsto} \
(\begin{td} \underset{-d}{P_{\star}} \ar[<-]{r}{a_{\pm}} \& \underset{1-d}{P_{\pm}}
\end{td})
\end{align*}
for any $d \in \Z$. In particular, the functor $\Dtr{-}$ preserves the path length.
Viewing the paths as monomials in $\Rx$, the functor operates via
\begin{align*}
\begin{cases}
(\begin{td} \underset{d+1}{P_{\pm}} \ar{r}{t^{n}} \& \underset{d}{P_{\star}} \end{td})
\ \overset{\Dtr{-}}{\longmapsto} \
(\begin{td} \underset{-d}{P_{\pm}} \ar[<-]{r}{t^{n+1}} \& \underset{1-d}{P_{\star}}\end{td}) 
& \text{for any }n \in \N_0,
\\
(\begin{td} \underset{d+1}{P_{\star}} \ar{r}{t^{n}} \& \underset{d}{P_{\pm}} \end{td})
\ \overset{\Dtr{-}}{\longmapsto} \
(\begin{td} \underset{-d}{P_{\star}} \ar[<-]{r}{t^{n-1}} \& \underset{1-d}{P_{\pm}}
\end{td}) & \text{for any }n \in \N,\\
(\begin{td} \underset{d+1}{P_{u}} \ar{r}{t^n} \& \underset{d}{P_{v}} \end{td})
\ \overset{\Dtr{-}}{\longmapsto} \
(\begin{td} \underset{-d}{P_{u}} \ar[<-]{r}{t^n} \& \underset{1-d}{P_{v}} \end{td})
& \text{for any }n \in \N\text{ if $(u,v) \notin \{(\pm,\star),(\star,\pm)\}$},
\end{cases}
\end{align*}
where
$d \in \Z$ and $u,v \in \{\star,+,-\}$.
In terms of diagrams, the application of this functor corresponds to taking the transpose of each arrow and adjusting the monomial power by the rules above. 
This motivates the next definition.
\begin{dfn}\label{dfn:w-tr}
The \emph{transpose} of a word $w$ of $\chfXA$ is defined as follows.
\begin{enumerate}
\item If $w = (x_i)_{i \in \Z}$ is periodic, 
we set $\check{w}^t = (x_{i}^{t})_{i \in \Z}$.
\item
If $w$ is a finite word of the form
$(\alpha, x_1, x_2, \ldots x_{\ell},\beta)$
with $w \neq (\star,x_1,\star)$
we set  $\check{w}^{t} = (\alpha,\check{x}_1^{t},x_2^t, \ldots, x_{\ell-1}^t, \check{x}_\ell^{t}, \beta)$
with 
\begin{align*}
	\check{x}_1^t = \begin{cases}
		\ola{n}_1^-&
		\text{if } (\alpha,x_1)=(\star, \ora{n}_1), \\
		\ora{n}_1^+&
		\text{if } (\alpha,x_1)=(\star, \ola{n}_1), \\
		x_1^t & \text{otherwise},
	\end{cases}
	&&
	\check{x}_\ell^t = \begin{cases}
		\ora{n}_\ell^-&
		\text{if } (x_\ell,\beta) = (\ola{n}_\ell,\star), \\
		\ola{n}_\ell^+&
		\text{if } (x_\ell,\beta) = (\ora{n}_\ell,\star), \\
		x_\ell^t & \text{otherwise}.
	\end{cases} 
\end{align*}
If $w = (\star,x_1,\star)$
we set $\check{w}^t =w$.
\end{enumerate}
\end{dfn}
It is straightforward to verify that the transpose of a word is again a periodic respectively finite word in the sense of Definitions~\ref{dfn:per-wo}, \ref{dfn:ow-ab}.
\begin{lem}\label{lem:tr}
Assume that $w$ is the word underlying a band or string $\omega$ of $\chfXA$.
Then the sequence of arrows in $(\check{w}^t)^{\updownarrow}$ 
can be obtained by inverting the sequence of arrows in
$w^{\updownarrow}$ except those at axes of symmetry.
\end{lem}
\begin{proof}
We assume first that $\omega$ is a string
with an underlying word
$w = (\alpha, x_1,x_2, \ldots x_{\ell},\beta)$. 
On the one hand,
$w^{\updownarrow} 
= (\alpha, x_1 \kappa_1 x_2 \kappa_2 \ldots  x_{\ell},\beta)$
with each arrow $\kappa_i$ determined via the ambient word $\ol{w}$.
On the other hand,
$(\check{w}^t)^{\updownarrow} = 
(\alpha,y_1 \kappa'_{1} y_2 \kappa'_2 \ldots y_\ell,\beta)$
with $y_1 = \check{x}_1^t$, $y_i = x_i^t$ for any $1 < i < \ell$ and $y_\ell = \check{x}_\ell^t$,
and each arrow $\kappa'_i$ determined via the ambient word $\ol{\check{w}^t}$.
Going through the possible ends of $w$, it can be verified that $\ol{\check{w}^t} = (\ol{w})^t$.
Let $1 \leq i< \ell$. 
It follows that
for any index $p \in \N_0$
it holds that
$x_{i-p} = x_{i+1+p}^{t}$
for the letters in $\ol{w}$
if and only if  $y_{i+1+p} = y^t_{i-p}$ for the letters in $\ol{\check{w}^t}$.
One of the following two cases occurs.
\begin{itemize}
	\item  
	Assume that 
	$x_{i-p} = x_{i+1+p}^{t}$
	for any $p \in \N_0$
	for the letters in $\ol{w}$.
	Then  $y_{i-p}^t = y_{i+1+p}$ for any $p \in \N_0$ for the letters in $\ol{\check{w}^t}$.
	In this case it holds that $\kappa_i = \downarrow = \kappa'_i$ by the symmetry rule in Subsection~\ref{subsec:orient}.
	\item Otherwise, there is a minimal number $k \in \N_0$ with  $x_{i-k} \neq x_{i+k+1}^t$. Since 
	$\ol{\check{w}^t} =(\ol{w})^{t}$, $k$ is also the minimal number in $\N_0$ with
	$y_{i-k} \neq y_{i+k+1}^t$. Moreover, $y_{i-k} = x_{i+k+1}^t$ and
	$y_{i+1+k} =  x_{i-k}^t$.
	Note that interchanging $\olra{p}$ and $\olra{q}$
	and changing their horizontal orientations
	in rule \eqref{eq:updownarrow} flips the arrow $\kappa_i$.
	Therefore, $\kappa_i \neq \kappa'_i$. 
\end{itemize} 
This shows the claim in the case that $\omega$ is a string.

The proof for the case that $\omega$ is a band 
is the same as above except for some formal changes of notation.
\end{proof}

\begin{ex}
	As in Subsection~\ref{subsec:sym-band},
	let $\omega$ be the band $(0,w,1,\lambda)$
	such that $w$ has periodic part
	$\dot{w} = ( \ola{1},\ora{1},\ora{2},\ola{1},\ora{1},\ola{2})$, and $\lambda \in \kk^*\backslash\{-1\}$.
	On the one hand, we apply 
	$\Dtr{-}$
	to the  diagram  \eqref{eq:symm-bd-glu}
	underlying the band complex $\CX(\omega)$
	and obtain:
	$$
	\begin{gd}
		P_+ \ar[densely dotted,-]{d} \ar{rd} \ar{r}{\frac{1}{\lambda}t}
		\& P_+ \ar[densely dotted,-]{d}
		\& P_+ \ar[densely dotted,-]{d} \ar[<-]{r}{-t^2}
		\& P_+ \ar[densely dotted,-]{d} \ar{r}{-t} \ar{rd}
		\& P_+ \ar[densely dotted,-]{d}
		\& P_+ \ar[densely dotted,-]{d}
		\& P_+ \ar[densely dotted,-]{d}
		\\
		\underset{1}{P_-} \ar{ru}[description]{t} \ar{r}[swap]{t}
		\&
		\underset{0}{P_-} \ar[<-]{ru}[description]{t} \ar[<-]{r}[swap]{t}
		\&
		\underset{1}{P_-} \ar[<-]{ru}[description]{t^2}
		\&
		\underset{2}{P_-} \ar{ru}[description]{t} \ar{r}[swap]{-t}
		\&
		\underset{1}{P_-} \ar[<-]{ru}[description]{t} \ar[<-]{r}[swap]{-t}
		\&
		\underset{2}{P_-} \ar{ru}[description]{t^2} \ar{r}[swap]{-{\frac{1}{\lambda}} t^2}
		\&
		\underset{1}{P_-}
	\end{gd}
	$$
	On the other hand, we consider
	the band $\omega^t = (1,\check{w}^t,1,\frac{1}{\lambda})$.
	Its oriented gluing diagram has the form:
	$$
	\begin{gd}
		P_+ \ar[densely dotted,<-]{d} 		\& P_+ \ar[color=blue,densely dotted,-]{d}[description]{\kappa'_1}   \& P_+ \ar[densely dotted,->]{d}  \& P_+ \ar[densely dotted,<-]{d} 
		\& P_+ \ar[color=blue,densely dotted,-]{d}[description]{\kappa'_4} 
		\& P_+ \ar[densely dotted,->]{d} 
		\& P_+ \ar[densely dotted,<-]{d} 
		\\
		\underset{1,\frac{1}{\lambda}}{P_-}  \ar[->]{ru}[description]{1} \&
		\underset{0}{P_-} \ar[<-]{ru}[description]{1} \&
		\underset{1}{P_-} \ar[<-]{ru}[description]{2} \&
		\underset{2}{P_-}  \ar[->]{ru}[description]{1} \&
		\underset{1}{P_-}  \ar[<-]{ru}[description]{1} \&
		\underset{2}{P_-}  \ar[->]{ru}[description]{2} \&
		\underset{1,\frac{1}{\lambda}}{P_-}
	\end{gd}
	$$
	Using the symmetry rule
	the gluing edges $\kappa'_1$ and $\kappa'_4$
	have to be oriented downwards, and we obtain
	\begin{align*}
		w^{\updownarrow} =(
		\ola{1} \underset{\kappa_1}{\downarrow} \ora{1} \uparrow \ora{2} \downarrow \ola{1} \underset{\kappa_4}{\downarrow} \ora{1} \uparrow \ola{2} \downarrow 
		)
		\quad\text{and}\quad
		{(\check{w}^t)}^{\updownarrow} =( 
		\ora{1} \underset{\kappa_1'}{\downarrow} \ola{1} \downarrow \ola{2} \uparrow
		\ora{1} \underset{\kappa_4'}{\downarrow} \ola{1} \downarrow \ora{2} \uparrow )
	\end{align*}
	in accordance with Lemma~\ref{lem:tr}.
	By Lemma~\ref{lem:sym-bd}
	using the upward orientations at $\kappa'_1$ and $\kappa'_4$
	yields an complex $\CX^{\uparrow}(\omega^t)$ isomorphic to the complex $\CX(\omega^t)$ obtained via the usual convention.
	The complex $\CX^{\uparrow}(\omega^t)$ 
	has precisely the underlying diagram given by the first diagram above,
	which shows that $\Dtr{\CX(\omega)} \cong \CX(\omega^t)$.
\end{ex}	

\begin{prp}\label{prp:rev}
	Let $\omega$ be a band or string of $\chfXA$.
	Then $\Dtr{\CX(\omega)} \cong \CX(\omega^t)$ in $\DbRep{\A}$,
	where $\omega^t$ is defined as follows, setting 
	$d^t = 1 - d$.
	\begin{align*}
		\begin{array}{|c||c|c|c|c|}
			\hline
			& \text{usual string} & \text{special string} & \text{bispecial string} & \text{band} \\
			\omega & (d,w) & (d,w,\varepsilon) & (d,\varepsilon_1,w, \varepsilon_2) & (d,w,m,\lambda) \\ 
			\hline 
			\omega^t&
			(d^t,\check{w}^t)
			& (d^t,\check{w}^t,{\varepsilon})
			& (d^t, \varepsilon_1, \check{w}^t, \varepsilon_2) &
			(d^t,\check{w}^t,m,\lambda^{\mp}  ) \\
			\hline
		\end{array}
	\end{align*}	
\end{prp}
\begin{proof}
	We prove the claim in the case that $\omega$ is a band $(d,w,m,\lambda)$, the proof for strings 
	uses similar, simpler considerations.	
	We denote $\dot{w} = (x_1,x_2,\ldots,x_{\ell})$
	and $w^{\updownarrow} = (x_i\kappa_i)_{i=1}^{\ell}$.
	In the construction of the differentials of $\CX(\omega)$, we may use transposed Jordan blocks
	because of Remark~\ref{rmk:field}~\eqref{rmk:field-mat}.
	Applying  $ \Dtr{-}$ to the closed gluing diagram 
	\eqref{eq:cyc-bd-diag}
	of $\omega$
	yields a diagram of the form
	\begin{align*}
		\begin{gd}
			\underset{d^t_0}{{P_{+-}^{(m)} }}
			\ar[<->]{r}{M_1^T  t^{n_1}}  \& \underset{d^t_1}{P_{+-}^{(m)}} 
			\ar[densely dotted,-]{r}\& \underset{d^t_{i-1}}{P_{+-}^{(m)}}  \ar[<->]{r}{M_i^T   t^{n_i}} 
			\&
			\underset{d^t_{i}}{P_{+-}^{(m)}}
			\ar[densely dotted,-]{r}
			\&
			\underset{d^t_{\ell-1}}{P_{+-}^{(m)}}
			\ar[ <->,bend left=20]{llll}[inner sep=0.5pt]{ M_{\ell}^T   t^{n_{\ell}}}  
		\end{gd}
	\end{align*}
	with transposed matrices and altered degrees.
	We note that the transposed matrices can be viewed as matrices of oriented subwords.
	\begin{itemize}
		\item Let $1 < i < \ell$. We denote $M_i = M(\kappa_{i-1},x_i,\kappa_i,d_i)$ to indicate that the matrix $M_i$ of Table~\ref{tab:bd-mat}
		is determined by these parameters.
		Note that the transposed matrix $M_i^T$ is precisely
		$M(\ol{\kappa}_{i-1},x_i^t,\ol{\kappa}_i, d_i^t)$
		for any values of $(\kappa_{i-1}, x_i, \kappa_i,d_i)$.
		\item At index $1$ it
		holds that $M_1 = M(\kappa_\ell, x_1, \kappa_1, d_1, \JB^T_\lambda)$ and $M_1^T= 
		M(\ol{\kappa}_\ell, x_1^t, \ol{\kappa}_1, d^t_1, \JB_{\frac{1}{\lambda}})$.
		\item The previous arguments are also true
		at index $\ell$. 
	\end{itemize}
	It follows that the cyclic diagram above
	is precisely the diagram associated to the sequence $ (x_1^t \ol{\kappa}_{1}x_2^t \ol{\kappa}_2 x_3^t \ldots,  x^t_{\ell} \ol{\kappa}_\ell)$
	together with multiplicity $m$ and eigenvalue $\frac{1}{\lambda}$.
	By Lemma~\ref{lem:tr},
	the last sequence is precisely 
	$(\check{w}^t)^{\updownarrow}$
	using upward oriented arrows at axes of symmetry.
	Together with Lemma~
	\ref{lem:sym-bd}
	it follows that $\Dtr{\CX(\omega)})
	= \CX^\uparrow(\omega^t) \cong \CX(\omega^t)$.
\end{proof}

\subsection{Dualities} \label{subsec:duals}
As in Subsection~\ref{subsec:res-rad}
let $w$ be a word of the form $(\alpha,x_1,x_2,\ldots x_{\ell},\beta)$.

\begin{enumerate}
	\item If $\alpha = \diamond$, we set
	$\ell^{-1}(w) = w$,
	and otherwise $\ell^{-1}(w)$ is defined via the table:
	\begin{align*}
		\begin{array}{|r||r|rl|r|r|r|r|}
			\hline
			w & 
			(\star, \ola{n}_1,v)
			&			
			(\star, {\color{blue} {\ora{n}_1}}, v) 
			& n_1 > 1
			&
			(\star, {\color{blue} \ora{1}}, \ora{n}_2, v) &
			(\star, {\color{blue} \ora{1}, \ola{n}_2 }, v)& 
			(\star, {\color{blue}\ora{1}},\diamond)
			\\  
			\hline
			\ell^{-1}(w) & 
			(\star, {\color{blue} \ola{0}}, \ola{n}^+_1, v) &
			(\star, {\color{blue} \ola{0}}, \ora{n}_1^-, v)  & &
			(\star, \ora{n}_2, v) & (\star, {\color{blue}\ola{n}_2}, v)  &  (\star, {\color{blue}\ola{0}},\diamond)
		\\
		\hline
		\end{array}
	\end{align*}
	\item We set $r^{-1}(w) = (\ell^{-1}(w^{\op}))^{\op}$ and $w^{-1}_\star = \ell^{-1}(r^{-1}(w))$.
\end{enumerate}
It is straightforward to verify that the operators $\ell^{-1}$, $r^{-1}$
and $(-)_{\star}^{-1}$ define operators inverse to $\ell$, $r$ and $(-)_{\star}$, respectively.

\begin{dfn}
	We will call a special string $\omega = (d,w,\varepsilon)$ of $\chfXA$
	\emph{simple} if it satisfies any of the following equivalent conditions.
	\begin{enumerate}
		\item $w = (\star,\ora{1},\diamond)$ or $(\diamond,\ola{1},\star)$.
		\item The complex $\CX(\omega)$ is isomorphic to a simple $\A$-module up to shift.
		\item $\CX(\omega) \cong S_\varepsilon [d]$ for some integer $d\in \Z$.
	\end{enumerate}
\end{dfn}

\begin{thm}\label{thm:duals}
	Let $\omega$ be a band or string of $\chfXA$.
	Then there are isomorphisms of complexes 
	$\dual{(\CX(\omega))} \cong \CX(\dual{\omega})$ and $\cdual{(\CX(\omega))} \cong \CX(\cdual{\omega})$
	with $\dual{\omega}$ and $\cdual{\omega}$ defined as follows.
	\begin{itemize}
		\item If $\omega$ is not a simple special string,
		then $\dual{\omega}$ and $\cdual{\omega}$ are given by the table
		\begin{align*}
			\begin{array}{|c|c|c|c|c|}
				\hline
				& \text{usual string} & \text{special string} & \text{bispecial string} & \text{band} \\
				\omega & (d,w) & (d,w,\varepsilon) & (d,\varepsilon_1,w, \varepsilon_2) & (d,w,m,\lambda) \\ 
				\hline
				\dual{\omega} &
				(d^t_\star,\dual{w})
				& (d^t_\star,\dual{w},\ol{\varepsilon})
				& (d^t_{\star}, \ol{\varepsilon}_1, \check{w}^t,   \ol{\varepsilon}_2) &
				(d^t_{\star},\check{w}^t,m,\lambda^{\mp}  ) \\
				\hline
				\cdual{\omega} &
				(d^t_\star,\dual{w})
				& (d^t_\star,\dual{w},{\varepsilon})
				& (d^t_{\star}, \varepsilon_1, \check{w}^t,   \varepsilon_2) &
				(d^t_{\star},\check{w}^t,m,\lambda^{\mp}) \\
				\hline
			\end{array}
		\end{align*}
		$d^t_{\star} = 2 - d$ if $w$ begins with $\star$ and
		$d^t_\star  =  1 - d$ otherwise,  and $\dual{w} = (\Tw^{-1}(w))^{t}$.
		\item 
		If $\omega$ is a simple special string $(d,w,\varepsilon)$, then
		$\dual{\omega} = 
		(d^t_\star,\dual{w},{\varepsilon})$
		and $\cdual{\omega}=(d^t_\star,\dual{w},\ol{\varepsilon})$.
	\end{itemize}
\end{thm}
\begin{proof}
	Propositions~\ref{prp:IR} and~\ref{P:twist}
	imply that $\dual{(-)} = \Dtr{-} \circ \Tw^{-1} \circ \inv$.
	By \eqref{eq:duals} it holds that $\cdual{(-)} = \dual{(-)} \circ \inv = \inv \circ \dual{(-)}$.
	Therefore the claims follow from Propositions~\ref{prp:rev} and \ref{prp:inv} combined with  Theorem~\ref{thm:tau}.
\end{proof}

\section{Indecomposable projective resolutions via gluing triples}
\label{sec:H0-triples}
In the present section, we address the question which band and string complexes are actually projective resolutions.
Throughout the rest of this article, the term \emph{projective resolution} means always a projective resolution of a module.

The answer is formulated in Theorem~\ref{thm:proj-res} in homological and combinatorial terms, as well as an enumeration of 
the relevant bands and strings in Table~\ref{tab:res-str-bands}.
The proof of the theorem uses the realization from Subsection~\ref{subsec:cat-glue}
of a complex $\CX$ in $\DbRep{\A}$ as pullback of a certain
gluing triple $(\CV,\CY,\vartheta)$.
In the first subsection, we formulate simple necessary conditions for $\CX$ to be a projective resolution in terms of the normalization complex $\CY$.
This restricts the list of candidates for projective resolutions severely. It turns out that any candidate complex qualifies as a projective resolution, because it is quasi-isomorphic to a `syzygy resolution' by which we mean a two-term complex of the form $(\begin{td} L\ar[hookrightarrow]{r} \& P_0 \end{td})$ 
with projective module $P_0$.

The last subsection derives a formula for the Jordan--H\"{o}lder-multiplicities of the non-trivial homology 
of a projective resolution in terms of the complexes $\CV$ and $\CY$.

\subsection{Criteria for projective resolutions}

Throughout this subsection, let $\CX$ 
denote a minimal complex from $\DbRep{\A}$.
We recall that the Gelfand order $\A$ can be viewed as the pullback 
of certain ring morphisms on the left,  and  $\CX$ as the pullback complex on the right
\begin{align}
	\label{eq:pb2}
	\begin{cd}
		\A \ar[hookrightarrow]{r} \ar{d} \ar[twoheadrightarrow]{d} 
		\arrow[phantom]{dr}[very near start,description]{\pullback}
		\& \B \ar[twoheadrightarrow]{d}{\pi} \\
		\AI \ar[hookrightarrow]{r}{\iota} \& \BI 
	\end{cd}
	&&
	\begin{cd}
		\CX \ar[twoheadrightarrow]{d}[swap]{\alpha} \ar[hookrightarrow]{r}{\beta}
		\arrow[phantom]{dr}[very near start,description]{\pullback}
		\& \CY \ar[twoheadrightarrow]{d}{\eta_{\CY}} \\
		\CV 
		\ar[hookrightarrow]{r}{\vartheta^*}
		\&
		\CW
	\end{cd}  
\end{align} 
where in the left diagram $\B$ is the hereditary envelope of $\A$,
and
the rings $\AI$ and $\BI$ are semisimple,
and in the right diagram the 
morphism
$\eta_{\CY}$  is the natural projection from $\CY= \B \oA \CX$ onto $\CW = \BI \oB \CY$,
and $\vartheta^*$
is 
the composition of 
the natural isomorphism $\begin{td}\vartheta \colon \BI \oAI \CV \ar{r}{\sim} \& \BI \oB \CY \end{td}$ with the unit morphism $\begin{td}\eta_{\CV}\colon \CV \ar{r} \& \BI \oAI \CV \end{td}$.
Further details can be found in Subsection~\ref{subsec:cat-glue}.

\begin{lem}\label{lem:H+}
	In the notations above, assume that $H_n(\CY) = 0$ for all $n > 0$. Then $H_n(\CX) = 0$ and $\X_{n+1} = 0$ for all $n > 0$.
\end{lem}
\begin{proof}
	By assumption,
	$\CY$ has the form $(\begin{td} Y_1 \ar[hookrightarrow]{r} \& Y_0\ldots
	\end{td})$.
	Using that there is a monomorphism 
	$\begin{td}\beta \colon \CX \ar[hookrightarrow]{r}\& \mathstrut_{\A} \CY \end{td}$
	of complexes of $\A$-modules, it follows that $\CX = (\begin{td} X_1 \ar[hookrightarrow]{r}\& X_0\ldots
	\end{td})$.
\end{proof}
We recall that the Gelfand order $\A$ and its hereditary envelope $\B$ 
have a common ideal $I = \A e_\star \A = \B e_\star \B$.

\begin{prp}\label{prp:homology}
	In the notations above, 
	we consider the following conditions.
	\begin{enumerate}[label=$\mathsf{(\alph*)}$]	
		\item
		\label{H0a1} The complex $\CY$ is a projective resolution.
		\item
		\label{H0b} The complex $\CX$ is a projective resolution.
			\item
			\label{H0c}
			The complex $\CY$ satisfies $I H_n(\CY) = 0$ and $H_{-n}(\CY) = 0$ for all integers $n > 0$.
		\end{enumerate}
		Then the implications
		\ref{H0a1} $\Rightarrow$ \ref{H0b} $\Rightarrow$ \ref{H0c} hold true.
	\end{prp}
	\begin{proof}
		By the right diagram in \eqref{eq:pb2} there is a short exact sequence of complexes of $\A$-modules
		$$
		\begin{cd} 0 \ar{r} \& \CX \ar{r}{
				\begin{psmallmatrix} \alpha \\ \beta \end{psmallmatrix}	
			} \& 
			\CV \oplus \CY \ar{r}{\pi} \& \CW \ar{r}
			\& 0 
		\end{cd}
		$$
		$\pi =
		\begin{pmatrix} \vartheta^* & \eta_{\CY} \end{pmatrix}	$,
		and $\beta$ and $\vartheta^*$ are monomorphisms of complexes.
		The last exact sequence gives rise to a long exact homology sequence
		\begin{align}
			\label{eq:homology}
			\begin{cd}
				\ldots 
				\ar{r} \&
				\W_{i+1} \ar{r} \&
				H_i(\CX)
				\ar{r} \& 
				\V_i \oplus H_i(\CY) \ar{r}{H_i(\pi)}
				\&
				\W_i \ar{r} \& \ldots
			\end{cd}
		\end{align}
		with $H_i(\pi) = 	
		\begin{pmatrix}	
			\vartheta^*_i & H_i(\eta_{\CY})
		\end{pmatrix}$
		where we have used that each of the complexes $\CV$ and $\CW$ has zero differentials according to Lemma~\ref{lem:minimal}.
		\begin{itemize}	
			\item If \ref{H0a1} is satisfied, $\CY = (
			\begin{td} \Y_1 \ar[hookrightarrow]{r} \& \Y_0 \end{td})$ and because of the monomorphism $\beta$ we obtain that $\CX = (\begin{td}\X_1 \ar[hookrightarrow]{r}\& \X_0\end{td})$.  
			\item
			Assume 
			\ref{H0b}.
			Let $n > 0$. Since $H_n(\CX) = 0$ the  map $H_n(\eta_{\CY})$ yields a $\B$-linear embedding of $H_n(\CY)$ into the $\BI$-module $\W_n$, thus $I H_n(\CY) \subseteq I \W_n = 0$.
			This shows the last implication. \qedhere
		\end{itemize}
	\end{proof}
	
	We note the following characterization of gluing triples corresponding to projective resolution which follows from results in \cite[Section 3.3.1]{Gnedin}.
	\begin{rmk} \label{rmk:pr-la}
		The complex
		$\CX$
		is a projective resolution if and only if 
		condition \ref{H0c} holds, $2\, \dim \V_j = \dim \W_j$ for all $j \geq 2$ and
		$\im \vartheta^*_n \cap \im H_n(\eta_{\CY}) = 0$ for all $n \geq 1$. 
	\end{rmk}
	However, checking the last condition of the characterization in practice involves linear algebra computations or a truncation argument.
	Instead of using this characterization, we will
	show that any indecomposable complex   satisfying \ref{H0c} of Proposition~\ref{prp:homology}, which is a `candidate' for a projective resolution, does qualify as such by a reduction to a complex  
	satisfying \ref{H0a1}, and thus  \ref{H0b}.

	\subsection{Bands and strings of projective resolutions}

In this section, we leave the categorical language of gluing triples 
and make use of the explicit combinatorics of band and string complexes provided by Subsection~\ref{sec:gluing}.

The following description of the normalization of a band or string complex was already given implicitly.
\begin{lem}\label{lem:norm-complex}
Let $\omega$ be a band or string of $\chfXA$, and $\CX(\omega)$ the corresponding band respectively string complex. Then the complex $\B \oA \CX(\omega)$ 
is isomorphic to the direct sum of two-term complexes obtained from the polarized gluing diagram of $\omega$ 
from Subsection~\ref{subsec:polar}
by deleting all gluing edges, replacing $P_{\pm}$ by $\BP_\diamond$, and $P_\star$ by $\BP_\star$.
\end{lem}
\begin{proof}
As $\F(\CX(\omega)) \cong \F \Gl\Ind({\mu(\omega)}) \cong \Ind({\mu(\omega)})$ by Proposition~\ref{prp:gluing} and Theorem~\ref{thm:detect}, it follows that $\B \oA \CX(\omega)$ is isomorphic to the complex $\CY$ described in 
Remark~\ref{rmk:ind-triple}. The latter coincides with the description in the claim.
\end{proof}

\begin{thm}\label{thm:proj-res}
Let $\omega$ be a
band or string of $\chfXA$.
Let $\CX = \CX(\omega)$ denote the corresponding band or string complex and set $\CY = \B \oA \CX$.
\begin{enumerate}
	\item \label{res1}The complex $\CX$ is a projective resolution of a finite-dimensional $\A$-module.
	\item \label{res2}
	The complex $\CY$ has the form $(\begin{td} \Y_2 \ar{r} \& \Y_1 \ar{r} \& \Y_0 \end{td} )$ with $IH_1(\CY) = 0$ and $\Y_0 \neq 0$.
	\item \label{res3}
	The complex $\CY$ satisfies $\Y_0 \neq 0$ and each indecomposable summand of $\CY$ is isomorphic to the shift of the simple $\B$-module $\wt{S}_{\diamond}[1]$
	or to an indecomposable finite-dimensional $\B$-module, that is, to one of the minimal projective resolutions
	\begin{align*}
		(\begin{gd} \underset{2}{\BP_\star} \ar{r}{\cdot t}\& \underset{1}{\BP_{\diamond}}\end{gd})
		\qquad
		\text{or}
		\qquad
		(\begin{gd}
			\underset{1}{\BP_{u}} \ar{r}{\cdot t^{n}} \& \underset{0}{\BP_v} \end{gd})
	\end{align*}
	where $u,v \in \{\star,\diamond\}$ and $n \in \N_0$ with $n \neq 0$ if $(u,v) \neq (\diamond,\star)$.
	\item \label{res4} The word $w$ underlying $\omega$ satisfies the following conditions:
	\begin{itemize}
		\item The support of $\omega$ is given by $\{0,1\}$ or $\{0,1,2\}$.
		\item If some successive oriented numbers $(x_i,x_{i+1})$ in $w$ are right oriented, then
		$\omega$ is a string, $d=2$, 			$i=1$ and $w$ begins with
		$(\alpha,x_1) = (\star,\ora{1})$.
		\item If some successive oriented numbers $(x_{j-1},x_{j})$ in $w$ are left oriented, then
		$\omega$ is a string, $d^{\op}  =2$,
		$j=\ell$ and $w$ ends with $(x_\ell,\beta) = (\ola{1},\star)$.
	\end{itemize}
	\item \label{res5} 
	$\omega$ is is equivalent to a band or string from Table~\ref{tab:res-str-bands}. 
\end{enumerate}
\end{thm}

\begin{table}[h]
\centering
\caption{Strings and bands of projective resolutions}
\label{tab:res-str-bands}
$
\begin{array}{|l|l|}
	\hline
			\text{Usual strings }(d,w) 
		&
		\text{Special strings }(d,w,\varepsilon) \\
		\hline
		\quad (0,(\star,
		{\color{blue} \ola{q}_1,\ora{p}_1, \ldots 
			\ola{q}_{k},\ora{p}_{k},}
		\star)) 
		&
		\quad (0,(\star,
		{\color{blue} \ola{q}_1,\ora{p}_1, \ldots \ola{q}_{k},\ora{p}_{k},}
		\diamond),\varepsilon)
		\\
		\quad (1,
		(\star,
		{\color{blue} \ora{p}_1,\ola{q}_1, \ldots \ora{p}_{k},\ola{q}_{k},}
		\star)) 
		&
		\quad
		(1,(\star,
		{\color{blue} \ora{p}_1,\ola{q}_1, \ldots \ora{p}_k,\ola{q}_k,}
		\diamond),\varepsilon)
		\\
		\quad (1,(\star,
		{\color{blue} \ora{p}_1,\ola{q}_1, \ldots \ora{p}_{k},}
		\star)) 
		&
		\quad (1,(\star,
		{\color{blue} \ora{p}_1,\ola{q}_1, \ldots \ora{p}_k,}
		\diamond),\varepsilon) 
		\\
		(2,(\star,\ora{1},
		{\color{blue} \ora{p}_1,\ola{q}_1, \ldots \ora{p}_k,}
		\star)) 
		&
		\quad (0,(\star,
		{\color{blue} \ola{q}_1,\ora{p}_1, \ldots \ola{q}_{k},}
		\diamond),\varepsilon)  
		\\
		(2,(\star,\ora{1},
		{\color{blue} \ora{p}_1,\ola{q}_1, \ldots \ora{p}_k,\ola{q}_k,}
		\star)) 
		&
		(2,(\star,\ora{1},
		{\color{blue} \ora{p}_1,\ola{q}_1, \ldots \ora{p}_k,}
		\diamond),\varepsilon) 
		\\
		(2,(\star,\ora{1},
		{\color{blue} \ora{p}_1,\ola{q}_1, \ldots \ora{p}_k,\ola{q}_k,}
		\ola{1},\star)) 
		&
		(2,(\star,\ora{1},
		{\color{blue} \ora{p}_1,\ola{q}_1, \ldots \ora{p}_k,\ola{q}_k,}
		\diamond),\varepsilon) 
		\\
		\hline
		\hline
			\text{Bispecial strings }(d,\varepsilon_1,w,\varepsilon_2)
		&\text{Bands }(d,w,m,\lambda)
		\\
		\hline
			\quad (0,\varepsilon_1,(\diamond,
			{\color{blue} \ola{q}_1,\ora{p}_1, \ldots \ola{q}_{k},\ora{p}_{k},}
			\diamond),\varepsilon_2)
		&
		((\ora{p}_i, \ola{q}_i)_{i \in \Z}, m,\lambda)
		\\
			\quad
			(1,\varepsilon_1,(\diamond,
			{\color{blue} \ora{p}_1,\ola{q}_1, \ldots \ora{p}_k,\ola{q}_k,}
			\diamond),\varepsilon_2)
		&
		\\
			\quad (1,\varepsilon_1,(\diamond,
			{\color{blue} \ora{p}_1,\ola{q}_1, \ldots \ora{p}_k,}
			\diamond),\varepsilon_2) 
		&
		\\ 
		\hline
	\end{array}$
\end{table}
\begin{proof}
	The implication \eqref{res1} $\Rightarrow$ \eqref{res2} holds by Proposition~\ref{prp:homology} \ref{H0b} $\Rightarrow$ \ref{H0c}.
	The implication
	\eqref{res2} $\Rightarrow$ \eqref{res3} 
	follows from the well-known classification of indecomposable objects in $\DbRep{\B}$ which was recalled in Lemma~\ref{lem:ind-env}.
	To show
	\eqref{res3} $\Rightarrow$ \eqref{res4},
	let $\CY$ have the form in \eqref{res3}. 
	\begin{itemize}
		\item The restriction on the support is immediate from $\supp(\omega) =
		\supp \X = 
		\supp \Y$
		and the assumptions in \eqref{res3}, where the last equality holds by  Remark~\ref{rmk:support}.
		\item Let ($x_i$,$x_{i+1}$)=($\ora{n}_i$,$\ora{n}_{i+1}$) be a subsequence of $w$. By Lemma~\ref{lem:norm-complex} 
		the number $x_i$ corresponds to a direct summand
		$
		(\begin{td} \BP_{u} \ar{r}{t^{n_i}} \& \BP_{v}  \end{td})^\mu[d_i]
		$
		of $\CY$
		with certain $u,v \in \{\star,\diamond\}$, $n \in \N_0$, $d_i \geq 1$ and $\mu = 1$ if $\omega$ is a string respectively $\mu = m$ if $\omega$ is a band $(d,w,m,\lambda)$.
		\eqref{res3} implies that $d_i = 1$, $(u,v) = (\star,\diamond)$ and $n_i = 1$.
		As $\star$ must be an end of the gluing diagram it follows that $\omega$ is a string with $\alpha = \star$ and $d=2$.
		\item The second condition follows from application of the first to the equivalent band respectively string $\omega^{\op}$.
	\end{itemize}
	This shows that \eqref{res4} follows.
	
	\eqref{res4} $\Rightarrow$ \eqref{res5}:
	Assume that \eqref{res4} holds.
	\begin{itemize}
		\item Assume that $w$ does not contain a pair of successive equioriented numbers.
		Then the orientations in the  sequence
		of oriented numbers
		$(x_1,x_2,\ldots, x_{\ell})$ of $w$ are alternating. The support condition forces
		$d = 0$ if $x_1$ is left oriented, and
		$d=1$ if $x_1$ is right oriented.
		By passing to an equivalent string or band we may assume that $d \geq d^{\op}$. If $\omega$ is a band, we may replace $\omega$ by an equivalent band in which the 
		first oriented number $x_1$ of the underlying word $w$ is right oriented.
		These restrictions yield all strings with $d \in \{0,1\}$ and the family of bands in Table~\ref{tab:res-str-bands}.
		\item If $\omega$ contains a pair of successive equioriented numbers we may assume that these numbers are right-oriented by passing to an equivalent string. Then $\omega$ is a string with $d=2$, and $w$ begins with $(\star,\ora{1})$ followed by an
		alternating sequence $(\ora{n}_1, \ola{n}_2, \ldots \olra{n}_{\ell})$ and 
		ends with $(\ola{1},\star)$, $(\star)$ or $(\diamond)$. This yields the strings in Table~\ref{tab:res-str-bands} with $d=2$.
	\end{itemize}
	
	\eqref{res5} $\Rightarrow$ \eqref{res1}.
	Assume that $\omega$ is a string or a band which is equivalent to one from Table~\ref{tab:res-str-bands}.
	\begin{itemize}
		\item If $\omega$ is band or a string with $d \neq 2$, then $\CY$ is given by a projective resolution $(\begin{td}\Y_1 \ar{r} \& \Y_0\end{td})$, and thus $\CX$ is a projective resolution by Proposition~\ref{prp:homology}~\ref{H0a1} $\Rightarrow$ \ref{H0b}.  
		\item Otherwise, $\omega$ is equivalent to a usual string $(d,w)$ from Table~\ref{tab:res-str-bands} with $d=2$,
		or a special string $(d,w,\varepsilon)$ from Table~\ref{tab:res-str-bands} with $d=2$. 
		\begin{itemize}
			\item Let $\omega$ be a usual string $(2,w)$
			with $w = (\star,\ora{1},\ora{p}_1, \ola{q_1},\ldots\ora{p}_k, \ola{q}_k,\ola{1},\star)$, in particular, $d^{\op}=2$.
			Lemma~\ref{lem:qis}
			yields the quasi-isomorphism $\pi$ of complexes in the diagram below.
			$$
			\qquad
			\begin{cd}
				\CX \ar[->>]{d}[swap]{\pi}  \&		P_\star \ar{d} \ar{r}{
					\begin{psmallmatrix}
						t \\ -t 
					\end{psmallmatrix}	
				} \& P_{+-}  \ar[->>]{d}[swap]{
					\begin{psmallmatrix} \iota_+ & \iota_-
					\end{psmallmatrix}	
				} \ar{r}{M_1 t^{p_1}} \& 
				P_{+-} \ar[equal]{d} \ar[densely dotted,-]{rr}\&[-1cm] \& P_{+-} \ar[<-]{r}{M_{2k} t^{q_k}} \ar[equal]{d} \& P_{+-}  \ar[->>]{d}{
					\begin{psmallmatrix} \iota_+ & \iota_-
					\end{psmallmatrix}	
				}
				\ar[<-]{r}{
					\begin{psmallmatrix}
						t \\ -t 
					\end{psmallmatrix}	
				}
				\& P_{\star}   \ar{d}\\
				\CL \ar[hookrightarrow,dashed]{d}[swap]{\iota}\&
				0 \ar{r} \ar{d} \& I_\star \ar[dashed,hookrightarrow]{d}[swap]{\iota_\star}
				\ar{r}{M_1 t^{p_1}} \&
				P_{+-} \ar[equal]{d}	\ar[densely dotted,-]{rr}  \& \& P_{+-} \ar[<-]{r}{M_{2k} t^{q_k}} \ar[equal]{d}\& I_\star \ar[hookrightarrow,dashed]{d}{\iota_\star} \ar[<-]{r}\& 0  \ar{d} \\
				\CX' \&
				0\ar{r} \& P_\star \ar{r}{M_1 t^{p_1}} \& P_{+-} \ar[densely dotted,-]{rr} \& \& P_{+-} \ar[<-]{r}{M_{2k} t^{q_k}} \& P_{\star} \ar[<-]{r} \& 0
			\end{cd}
			$$
			In this diagram, $P_{+-} = P_+ \oplus P_-$, and the maps $\iota_+, \iota_-$ and $\iota_\star$ are natural inclusions
			with simple cokernels.
			Folding the morphisms between the bottom two rows yields a monomorphism $\iota$ of complexes  of $\A$-lattices.
			As $\coker \iota \cong S_\star^2[1]$ and $\CX$ 
			has finite-dimensional homology, so does the complex $\CX'$.
			By Lemma~\ref{lem:tri-fd}
			the two-term complex $\CY' = \B \oA \CX'$ has finite-dimensional homology as well.
			Since $\B$ is hereditary, $H_i(\CY') =0 $ for any $i \neq 0$.
			Because there are monomorphisms $\begin{td} \CL \ar[hookrightarrow]{r} \&\CX' \ar[hookrightarrow]{r} \& \CY' \end{td}$ of two-term complexes of $\A$-modules it follows that  $H_i(\CX) \cong H_i(\CL) = 0$ for any $i > 0$.	
			\item In the case of any other string with $d=2$ from Table~\ref{tab:res-str-bands}, we may use a simplification of the argument  above, in which the right ends of the line diagrams of $\CX$, $\CL$ and $\CX'$ are identical and 
			$\coker \iota \cong S_\star[1]$.
		\end{itemize}
		In both cases, we conclude that $H_i(\CX) =0$ for any $i > 0$.
	\end{itemize}
	This shows that all five conditions are equivalent.
\end{proof}

\subsection{Jordan-Hölder multiplicities via pullback data}

The next observation will be useful to compute the Jordan--Hölder multiplicities, that is, the dimension vectors, of indecomposable objects in $\Rep{\A}$.
\begin{lem}\label{lem:dim-formula}
Assume that $\CX$ is a minimal projective resolution of an object $M$ in $\Rep{\A}$.
Set $\CV = \AI \oA \CX$ and $\CY = \B \oA \CX$.
Then the dimension vector of the finite-dimensional $\A$-module $M$ is given by the formula
\begin{align}
\label{eq:dim-formula}
\uldim M = \uldim H_0(\CY) 
- \uldim H_1(\CY) +
\uldim \inv(\V_1) - \uldim \inv(\V_0)
\end{align}
where $H_1(\CY)$ and $H_0(\CY)$ are considered as $\A$-modules.
\end{lem}
\begin{proof}
There is a minimal triple $(\CV,\CY,\vartheta)$ isomorphic to $\F(\CX)$ because of Lemma~\ref{lem:minimal}.
According to Theorem~\ref{thm:detect} there is an isomorphism $\Gl\F(\CX) \cong \CX$.
By \eqref{eq:pb-complex} there is a commutative diagram with exact rows of complexes of $\A$-modules:
\begin{align*}
\begin{cd}
0 \ar{r} \&
\CX \ar[twoheadrightarrow]{d}[swap]{\beta} \ar{r}{\alpha} 		\arrow[phantom]{dr}[very near start,description]{\pullback}
\& \CY \ar[twoheadrightarrow]{d}{\eta_{\CY}} \ar{r} \& \coker \alpha  
\arrow[d, "\rotatebox{90}{$\sim$}"]
\ar{r} \& 0 \\
0\ar{r} \&	\CV 
\ar{r}{\vartheta^*}
\&
\CW \ar{r} \& \coker \vartheta^* \ar{r} \& 0
\end{cd}  
\end{align*}
Since
$\mathstrut_{\A} (\BI \otimes_{\AI} S_{\pm}) \cong S_+ \oplus S_-$, it follows that there are isomorphisms of complexes of $\A$-modules
$\coker \alpha \cong \coker \vartheta^* \cong \coker (\begin{td}\eta_{\CV}\colon \CV \ar{r} \& \BI \oAI \V\end{td}) \cong \inv(\CV)$, where the latter is the complex obtained from $\CV$ by interchanging $S_+$ and $S_-$ at all degrees.
Using that $\CV$ has zero differentials, the first row in the last diagram yields the long exact homology sequence:
\begin{align*}
\begin{cd}
\ldots 
H_{i+1}(\CY) \ar{r} \& \inv(\V_{i+1}) \ar{r} \&
H_i(\CX) \ar{r} \&
H_i(\CY) \ar{r} \&
\inv(\V_i) \ar{r}  \&  
\ldots
\end{cd}
\end{align*}
As $\CX$ is an $H_0$-complex,
this sequence 
implies 
the formula for the dimension vector.
\end{proof}	
We observe the following consequence, which does not require the notion of gluing triples.
\begin{cor}\label{cor:dim}
Let $\CX$ and $\CX'$ be projective resolutions of objects $M$ respectively $M'$ in $\Rep{\A}$.
Assume that there are isomorphisms
$\B \oA \CX \cong \B \oA \CX'$ 
in $\DbRep{\B}$ and $\AI \oA \CX \cong \AI \oA \CX'$ in $\DbRep{\AI}$.
Then $\uldim M = \uldim M'$.
\end{cor}
\begin{proof}
By Lemma~\ref{lem:minimal},
$\F(\CX)$ and	 $\FF(\CX')$
are isomorphic to minimal triples
$(\CV,\CY,\vartheta)$ and
$	(\CV', \CY',\vartheta')$, respectively.
Since $\CY \cong \CY'$ and $\CV \cong \CV'$, 
it follows that $\uldim H_i(\CY) = \uldim H_i(\CY')$ and $S_i \cong S'_i$ for any $ i \in \Z$. 
Lemma~\ref{lem:dim-formula}
implies that $\uldim M= \uldim M'$.
\end{proof}

\section{Band and string resolutions of the Gelfand quiver}
\label{sec:resolutions}

In this section, we
introduce bands and strings of $\Acat$
as combinatorial data which encode minimal projective resolutions. 
The main result of this section (Theorem~\ref{thm:prbij})
yields a bijection between bands and strings of $\Acat$ and indecomposable projective resolutions in $\DbRep{\A}$.
The proof is based on Theorem~\ref{thm:proj-res}
together with the explicit description of band and string complexes of Subsection~\ref{sec:gluing}.
For a band or string $\omega$ the corresponding resolution $\CX(\omega)$ 
can be viewed as gluing of projective resolutions of cyclic finite-dimensional $\A$-modules.
These projective resolutions are precisely the direct summands of the \emph{cyclification of $\CX(\omega)$} introduced in Subsection~\ref{subsec:cyclif}, which will be useful in the next sections.

\subsection{Motivation: Projective resolutions of cyclic modules}
\label{subsec:motiv-cyc}

To motivate the notions below, we consider first 
the `gluing diagrams'
\begin{align}
	\label{eq:cyc-diag1}
	\begin{gd}
		\& P_+ \ar[densely dotted,-]{d} \& \smash[b]{\underset{0}{P_\beta}} \\
		\underset{2}{P_\star}	 \ar{ru}{t}		\& \underset{1}{P_-} \ar{ru}{t^p} \& 
	\end{gd}
	&&  \text{and}
	&&
	\begin{gd}
		\& \smash[b]{\underset{0}{P_\beta}}\\
		\underset{1}{P_\alpha} \ar{ru}{t^p} \& 
	\end{gd}
\end{align}
where $\alpha, \beta \in \{\star,+,-\}$ is a vertex of the Gelfand quiver and $p \in \N_0$.
We require also that $p \neq 0 $ if $\beta \neq\star$ in the left, respectively if $(\alpha,\beta) \neq (\pm,\star)$ in the right diagram
in order to ensure that each multiplication by $t^p$ is well-defined and non-invertible.
We set $\alpha = \hat{\star}$ for 
the left diagram to indicate the appearance of a projective 
$P_\star$ at degree $2$.

A gluing diagram of the form above or its parameter triple $(\alpha,p,\beta)$ gives rise to
a minimal projective resolution, which already appeared in precisely the same form in~\eqref{eq:cyc-res},
\begin{align}\label{eq:cyc-res2}
	\CX(\alpha,p,\beta) = 
	\left\{
	\begin{array}{rl}
		\begin{cd} P_\star \ar{r}{
				\left(\begin{smallmatrix}
					\phantom{-}t \\ -t
				\end{smallmatrix}
				\right)
			} \& P_+ \oplus P_- \ar{r}{
				\left(\begin{smallmatrix}
					{t^{p}}&
					{t^{p}}	
				\end{smallmatrix}\right)	} \& P_\beta
		\end{cd} & \text{if }\alpha =\hat{\star}, \\
		\begin{cd} \mathstrut \& P_\alpha \ar{r}{t^{p}} \& P_\beta
		\end{cd}
		& \text{if }\alpha \in Q_0,
	\end{array} 
	\right.
\end{align}
of a finite-dimensional cyclic module $M(\alpha,p,\beta) = H_0(\CX(\alpha,p,\beta))$.

In this section, we will `glue' a projective resolution of an arbitrary 
finite-dimensional indecomposable $\A$-module from 
the above projective resolutions of the form 
$\CX(\alpha,p,\beta)$.
The basic idea in this gluing process is that we may `concatenate' ends of certain gluing diagrams 
by adding a gluing edge between them.
In this gluing process, we have to allow 
\emph{half-turns} of the diagrams above:
\begin{align}
	\label{eq:cyc-diag2}
	\begin{gd}
		\& P_- \ar[densely dotted,-]{d} \& \smash[b]{\underset{2}{P_\star}} \\
		\underset{0}{P_\beta}	 \ar[<-]{ru}{t}		\& \underset{1}{P_+} \ar[<-]{ru}{t^p} \& 
	\end{gd}
	&&  \text{and}
	&&
	\begin{gd}
		\& \smash[b]{\underset{1}{P_\alpha}}\\
		\underset{0}{P_\beta} \ar[<-]{ru}{t^p} \& 
	\end{gd}
\end{align}
An example of gluing of diagrams is given  as follows.
\begin{align*}
	\begin{gd}
		\& P_+ \ar[densely dotted,-]{d} \& \smash[b]{\underset{0}{P_+}} \\
		\underset{2}{P_\star}	 \ar{ru}{t}		\& \underset{1}{P_-} \ar{ru}{t^p} \& 
	\end{gd}
	&&
	\begin{gd}
		\& \smash[b]{\underset{1}{P_\star}}\\
		\underset{0}{P_-} \ar[<-]{ru}{t^q} \& 
	\end{gd}
	&&
	\Rightarrow
	&&
	\begin{gd}
		\& P_+ \ar[densely dotted,-]{d} \& {P_+} \ar[densely dotted,-]{d} \& \underset{1}{P_\star} \\
		\underset{2}{P_\star}	 \ar{ru}{t}		\& \underset{1}{P_-} \ar{ru}{t^p} \& \underset{0}{P_-} \ar[<-]{ru}{t^q} \& 
	\end{gd}
\end{align*}
In a formal language this gluing will be viewed as composition
of sequences $(\hat{\star},\ora{p},+)$ with $(-,\ola{q},\star)$
which results in the expression
$(\hat{\star},\ora{p},\ola{q},\star)$.
In a later step, the last diagram will obtain new arrows, and give rise to an indecomposable projective resolution.

			\subsection{Band resolutions}\label{sec:proj-res1}
			Next, we will modify the definitions
			of periodic 
			words
			and bands 
			of $\chfXA$ 
			in order to define counterparts for $\Acat$. For readers familiar with Subsection~\ref{subsec:strands}, we have  highlighted the modifications in the new definitions.
			Bands of $\Acat$ will be translated in certain diagrams.
			
			The goal of this is to define a projective resolution $\CX(\omega)$
			for each band $\omega$ of $\Acat$.
			The main technical step is to transform the underlying word $w$ into a datum $w^{\updownarrow}$ containing  a vertical arrow $\uparrow$ or $\downarrow$ between any two successive oriented numbers.

			\subsubsection{Periodic words and bands}
			
			As before, we will use the sets of formals symbols
			\begin{align*}
				\kX =  \ora{\N} \cup \ola{\N}
				\quad	\text{ with } \quad \ora{\N} = \{ \ora{n} \mid n \in \N \}
				\quad	\text{ and } \quad
				\ola{\N} = \{ \ola{n} \mid n \in \N \},
			\end{align*}
			For any oriented number $x$, flipping its arrow defines $x^t$.
			\begin{dfn}
				Let $w = (x_i)_{i \in \Z}$ be a sequence of positive oriented numbers
				{\color{blue} with alternating orientations}.
				\begin{enumerate}
					\item For any integer $j \in \Z$  the \emph{$j$-th shift of $w$} is given by 
					$w_j \colonequals (x_{i+j})_{i \in \Z}$.
					\item The sequence $w$ has period $\ell$ if $w=w_{\ell}$ and $\ell$ is the minimal number from $2\N$ 
					with this property.
					In this case, the sequence $w$ is called a \emph{periodic word} and we denote
					$\dot{w} = (x_1,x_2,\ldots x_{\ell-1}, x_{\ell})$.
					\item The \emph{opposite word of $w$} is defined by 
					$
					w^{\op} \colonequals (x^{t}_{1-i})_{i \in \Z}$.
					\item The word $w$ is \emph{symmetric} if there is an integer $j \in \Z$ such that $w_j = (w_j)^{\op}$, that is, $x_{j-p} = x^{t}_{j+1+p}$ for any $p \geq 0$.
				\end{enumerate}
			\end{dfn}
			We give a diagrammatic interpretation of these notions.
			For an explicit notation below, we assume that the periodic part $\dot{w}$ of $w$ begins with a right-oriented number, that is,
			$\dot{w} =(\ora{p}_1,\ola{q}_1,\ldots  \ora{p}_k,\ola{q}_k)$. To the sequence $\dot{w}$ we associate the  gluing diagram
			\begin{align} \label{eq:per-ab}
				\begin{gd}
					{\color{blue} \bt} \ar[color=blue,densely dotted,-]{d}\& \bt \ar[densely dotted,-]{d} \&\bt  \ar[densely dotted,-]{d} \&  \cdots   \&
					\bt \ar[densely dotted,-]{d}  \& \bt \ar[densely dotted,-]{d}  \&  {\color{blue} \bt} \ar[densely dotted,-,color=blue]{d}   \\
					{\color{blue} \bt}  \ar[->]{ru}[description]{p_1} \& \bt  \ar[<-]{ru}[description]{q_1} \& \bt \&  \cdots  \& \bt \ar[->]{ru}[description]{p_k} \& \bt \ar[<-]{ru}[description]{q_k} 	\& {\color{blue} \bt}
				\end{gd}
			\end{align}
			where the first and the last vertical edge are identified.
			The diagram of the first shift $w_1$ has periodic part $(\ola{q}_1, \ldots  \ora{p}_{k}, \ola{q}_{k}, \ora{p}_1)$ and its diagram can be obtained 
			by a cyclic permutation of the diagram above.
			The diagram of $w^{\op}$ has periodic part $(\ola{q}_{k}, \ora{p}_{k}, \ldots  \ola{q}_1, \ora{p}_1)$ and
			can be obtained by a half-turn of the diagram of $w$.
			The word $w$ is symmetric if and only if its diagram admits a midpoint of symmetry.

			\begin{dfn}
				\label{dfn:bd-sh}
				A \emph{band of $\Acat$} is a triple $(w,m,\lambda)$ given by a periodic word $w$ of $\Acat$, a number $m\in \N$ and a  parameter $\lambda \in \kk^*$ such that $\lambda \neq -1$ if $w$ is symmetric. 	
			\end{dfn}
			
			\begin{rmk}\label{rmk:trans1}
				Any band $\omega=(w,m,\lambda)$ of $\Acat$ gives rise to a band $\breve{\omega}=(d,w,m,\lambda)$ of $\chfXA$ with $d=1$ if the first letter of $w$ is right-oriented, and $d=0$ otherwise.
			\end{rmk}

			\subsubsection{From band diagrams to band resolutions}
			Throughout this subsection, let $\omega$ be a band $(w,m,\lambda)$ of $\Acat$
			such that 
			$\dot{w} =(\ora{p}_1,\ola{q}_1,\ldots  \ora{p}_k,\ola{q}_k)$.

			\paragraph{Polarized band diagram}
			Next, the diagram~\eqref{eq:per-ab} is transformed into the diagram of projective $\A$-modules
			\begin{align}\label{eq:per-ab2}
				\begin{gd}
					{{\color{blue} P_{+}^m}} \ar[color=blue,densely dotted,-]{d}\& P_{+}^m \ar[densely dotted,-]{d} \& P_{+}^m  \ar[densely dotted,-]{d} \&  \cdots   \&
					P_{+}^m \ar[densely dotted,-]{d}  \& P_{+}^m \ar[densely dotted,-]{d}  \&  {\color{blue} P_{+}^m} \ar[densely dotted,-,color=blue]{d}   \\
					\underset{1,\lambda}{{\color{blue} P_{-}^m}}  \ar[->]{ru}[description]{t^{p_1}} \& \underset{0}{P_{-}^m}  \ar[<-]{ru}[description]{t^{q_1}} \& \underset{1}{P_{-}^m} \&  \cdots  \& \underset{1}{P_{-}^m} \ar[->]{ru}[description]{t^{p_k}} \& \underset{0}{P_{-}^m} \ar[<-]{ru}[description]{t^{q_k}} 	\& \underset{1,\lambda}{{\color{blue} P_{-}^m}}
				\end{gd}
			\end{align}
			where the first and the last column have to be identified.

			\paragraph{Orienting the vertical edges}
			The periodic word $w = (x_i)_{i\in \Z}$
			gives rise to another periodic combinatorial sequence
			$w^{\updownarrow} = 
			(x_i \kappa_i)_{i \in \Z}$
			with intermediate arrows $\kappa_i \in \{\uparrow, \downarrow\}$ for each index $i \in \Z$,
			which are determined as follows.
			
			\begin{rmk}\label{rmk:kappa1}
				For any index  $i \in \Z$, 
				the arrow $\kappa_i$ is  $\uparrow$ if and only if $w$ does not have a midpoint at symmetry between $x_i$ and $x_{i+1}$,
				that is,  there exists
				a minimal number $j \in \N_0$ with $x_{i-j} \neq x^t_{i+1+j}$,
				and one of the following conditions is satisfied.
				\begin{align}\label{eq:or-ru} 
					(x_{i-j},x_{i+1+j})= 
					(\ora{p},\ola{q})\text{ with }p < q,\text{ or }
					(\ola{q},\ora{p})\text{ with }q > p.
				\end{align}
		\end{rmk}

		In what follows, we need to distinguish between the case of a band with a two-periodic underlying word and the remaining bands.
		
		\paragraph{Band words of period two}
		Let $\omega$ be a band $(w,m,\lambda)$ with $\dot{w} = (\ora{p},\ola{q})$, $m\in \N$ and $0 \neq \lambda \in \kk$ such that $\lambda \neq -1$ if $p = q$.
		We assume also that $p \geq q$.
		It then follows that $w^{\updownarrow} = (\ora{p} \downarrow \ola{q} \downarrow)$.
		Then the final gluing diagram 
		and the projective resolution $\CX(\omega)$ are given as follows.
		\begin{align}\label{eq:sh-bd-res}
			\begin{gd}
				{{\color{blue} P_{+}^m}} \ar[color=blue,densely dotted]{d}\& P_{+}^m \ar[densely dotted,->]{d} \ar[<-,color=red,dashed]{rd} \ar[<-,color=red]{r}{t^{q}}  \&  
				{\color{blue} P_{+}^m} \ar[densely dotted,->,color=blue]{d}   \\
				\underset{1,\lambda}{{\color{blue} P_{-}^m}}  \ar[->]{ru}[description]{t^{p}} \& \underset{0}{P_{-}^m}  \ar[<-]{ru}[description]{t^{q}} 	\& \ar[color=green]{l}{\JB_{\frac{1}{\lambda}} t^{q}}  \underset{1,\lambda}{{\color{blue} P_{-}^m}}
			\end{gd}
			\quad
			\begin{array}{ll}
				\Rightarrow &
				(\begin{cd} X_1 \ar{r}{\partial_1} \& X_0 \end{cd})
				\\
				&
				\text{with }
				\partial_1 = 
				\begin{pNiceArray}{cc}[first-row, last-col, columns-width = auto]
					\CodeBefore
					\cellcolor{red!25}{1-1}
					\cellcolor{blue!25}{2-1}
					\cellcolor{blue!25}{1-2}
					\cellcolor{red!25}{2-2}
					\Body
					P_-^m & P_+^m
					\\
					t^{p} + \JB_{\frac{1}{\lambda}} t^q &  t^{q} & P_+^m \\
					\JB_{\frac{1}{\lambda}} t^{q} & t^{q} 
					&P_-^m  \end{pNiceArray}
		\end{array}
	\end{align}
	As before, the first vertical edge in the gluing diagram has to be identified with the last vertical edge,
	and the dashed diagonal arrow in the second square is labeled by $\JB_{\frac{1}{\lambda}} t^q$.
	Such a band is peculiar, because
	the gluing diagram has an induced arrow with the same source and target as one of the initial arrows.
	The parallel arrows are combined to the sum of entries in the top left block matrix of the differential $\partial_1$.
	\paragraph{Band words of higher period}
	Next, we consider the family of bands $(w,m,\lambda)$
	such that the period $2k$ of the word $w$ is at least four.
	As before, we 
	assume that 
	$\dot{w}$ starts with a right-oriented arrow, that is,
	$\dot{w} =(\ora{p}_1,\ola{q}_1,\ldots  \ora{p}_k,\ola{q}_k)$.

	Remark~\ref{rmk:kappa1} yields a periodic sequence $w^{\updownarrow} = (\ora{p}_i \kappa_{2i-1} \ola{q}_i \kappa_{2i})_{i \in \Z}$
	with vertical arrows.
	These vertical arrows are used to  define binary coefficients for each integer $i \in \Z$:
	\begin{align}
		\label{eq:coeff}	
		\begin{array}{ll}
			(\alpha_i, \ol{\alpha}_i)
			= \begin{cases}
				(1,0) & \text{if $\kappa_{2i-1}$ is  $\downarrow$},\\
				(0,1) & \text{otherwise},
			\end{cases}
			&
			\gamma_i = \alpha_i \beta_i = 
			\begin{cases}
				1 & \text{if } (\kappa_{2i-1},\kappa_{2i}) = (\downarrow,\downarrow), \\
				0 & \text{otherwise},
			\end{cases}
			\\
			(\beta_i,\ol{\beta}_i) = 
			\begin{cases}
				(1,0) & \text{if $\kappa_{2i}$ is $\downarrow$}, \\
				(0,1) & \text{otherwise},
			\end{cases} &
			\ol{\gamma}_i = \ol{\beta}_i \ol{\alpha}_{i+1} = 
			\begin{cases}
				1 & \text{if } (\kappa_{2i},\kappa_{2i+1}) = (\uparrow,\uparrow), \\
				0 & \text{otherwise}.
			\end{cases}   
		\end{array}
	\end{align}
	\begin{rmk}\label{rmk:kap-imp}
		The following implications are true for any integer $i \in \Z$:
		\begin{align*}
			\begin{array}{cclclcl}
				p_i > q_i &\Rightarrow &
				(\alpha_i,\ol{\alpha}_i) = (1,0)
				\\
				p_i < q_i 
				&\Rightarrow&
				(\alpha_i,\ol{\alpha}_i) =(0,1) 
			\end{array}&&
			\begin{array}{cclclcl}
				q_i < p_{i+1} &\Rightarrow &
				(\beta_i,\ol{\beta}_i) =(1,0) 
				\\
				q_i > p_{i+1} &\Rightarrow&
				(\beta_i,\ol{\beta}_i) =(0,1) 
			\end{array}
		\end{align*}
	\end{rmk}
	Then $\CX(\omega) $ is given by a projective resolution of the form 
	\begin{align}\label{eq:l-bd-diff}
		(\begin{cd} \X_1 \ar{r}{\partial_1} \& \X_0 \end{cd})\text{ with 
			the 
			differential $\partial_1$ given by the block matrix}
	\end{align}
	\begin{align*}
		\scalebox{1.0}{$
			\begin{pNiceArray}{cccccccccc}[first-row, last-col
				, columns-width = auto
				]
				\CodeBefore
				\chessboardcolors{red!10}{blue!10}
				\foreach \i in {1,...,10}{
					\pgfmathtruncatemacro{\jmin}{max(1,\i-2)}
					\pgfmathtruncatemacro{\jmax}{min(10,\i+2)}
					\foreach \j in {\jmin,...,\jmax}{
						\pgfmathtruncatemacro{\sum}{mod(\i-\j,2)}
						\ifnum\sum=0
						\cellcolor{red!25}{\i-\j}
						\else
						\cellcolor{blue!25}{\i-\j}
						\fi
					}
				}
				\cellcolor{blue!25}{1-10}
				\cellcolor{red!25}{2-10}
				\cellcolor{red!25}{9-1}
				\cellcolor{blue!25}{10-1}
				\Body
				P_-^m & P_+^m & P_-^m & P_+^m &P_-^m &P_+^m &  &  & P_-^m & P_+^m 
				\\
				t^{p_1} & \alpha_1 t^{q_1} & \gamma_{1} t^{q_1} & 0 &\Cdots &\Cdots&\Cdots&\Cdots&0& \scalebox{0.8}{$\ol{\beta}_k \JB_{\lambda} t^{p_1}$} &P_+^m \\
				\ol{\alpha}_1 t^{p_1} & t^{q_1} & \beta_1 t^{q_1} & 0 &\Ddots&&& & 0 & 
				\scalebox{0.8}{$\ol{\gamma}_k  \JB_{\lambda} t^{p_1}$}
				&P_-^m  \\
				0 & \ol{\beta}_1 t^{p_2} & t^{p_2 }& \alpha_2 t^{q_2} & \gamma_2 t^{q_2} &\Ddots&&&0 & 0 &P_+^m  \\
				0 & \ol{\gamma}_1 t^{p_2} & \ol{\alpha}_2 t^{p_2} & t^{q_2} & \beta_2 t^{q_2} & 0 &\Ddots&& & \Vdots & P_-^m \\
				\Vdots		&\Ddots & 0 & \ol{\beta}_2 t^{p_3} & t ^{p_3}& \alpha_3 t^{q_3} & &  \Ddots &&& P_+^m\\
				\Vdots	&&\Ddots & \ol{\gamma}_2 t^{p_3} & \ol{\alpha}_3 t^{p_3} & t^{q_3} &  & \Ddots &  \Ddots &  \Vdots & P_-^m \\
				\Vdots &&&&		\Ddots & \Ddots &		
				\Ddots
				& 
				\Ddots
				&
				\Ddots
				&  0 &  \\
				0 & 0 &	&&& \Ddots & 	
				\Ddots
				& 
				\Ddots
				& 
				\Ddots
				&  0 &   \\
				\scalebox{0.8}{$\gamma_k   \JB_{\frac{1}{\lambda}} t^{q_k} $}
				&	0  	&&&& 0&		\Ddots & 
				\Ddots
				& t^{p_{k}}& \alpha_k t^{q_k} & P_+^m \\
				\scalebox{0.8}{$\beta_k  \JB_{\frac{1}{\lambda}} t^{q_k} $} & 0 & \Cdots && &&0 & 
				*
				& \ol{\alpha}_{k} t^{p_k} & t^{q_k} &	P_-^m \end{pNiceArray}$}
	\end{align*}
	In particular, the block matrix is `pentadiagonal' when its horizontal and vertical block indices are considered modulo the period length $2k$.
	Remark~\ref{rmk:kap-imp} yields an immediate description of the differential $\partial_1$ in case that any two successive numbers in the periodic word $w$ are distinct.

	\subsection{String resolutions}\label{sec:proj-res2}
	
	Next, we define 
	finite words 
	and strings 
	for $\Acat$. These have counterparts for $\DbRep{\A}$.
	Then we define a projective resolution $\CX(\omega)$ for each string $\omega$ of $\Acat$ by similar steps as in the previous subsection.

\subsubsection{Finite words and strings}

\begin{dfn}
A \emph{finite word} of $\Acat$
is given by  a sequence
$(\alpha, x_1, x_2, \ldots x_{\ell}, \beta)$ 
with ends $\alpha, \beta \in \{ {\color{blue} \wh{\star}}, \star,\diamond\}$, oriented numbers $x_2,\ldots x_{\ell-1} \in \kX$
and $x_1, x_{\ell} \in \kX_0$  such that 
the following conditions hold.
\begin{itemize}
\item {\color{blue} The sequence of orientations in $(x_1,x_2,\ldots x_{\ell})$ is alternating.}
\item {\color{blue} If $\alpha = \hat{\star}$, then $x_1$ is right-oriented. Similarly, if $\beta = \hat{\star}$, then $x_{\ell}$ is left-oriented.}
\item $x_i = \ola{0}$ only if $i=1$, $\alpha = \star$ and $w \neq (\star,x_1,\star)$.
\item $x_i = \ora{0}$ only if $i = \ell$, $\beta = \star$ and $w \neq (\star,x_1,\star)$.
\end{itemize}  
Any word $w$ of the form above gives rise to the following notions.
\begin{enumerate}
\item The \emph{opposite word} of the word $w$ is defined as 
$
w^{\op} \colonequals (\beta,  x^t_{\ell},  \ldots x_2^t, x^t_1, \alpha)$.
\item The word $w$ is \emph{asymmetric} if $w \neq w^{\op}$.
\item  An end $\alpha$ or $\beta$ of the word $w$ is called \emph{special}
if it is given by $\diamond$.
\item The word $w$ is \emph{usual}, \emph{special} respectively \emph{bispecial} if it has zero, one respectively two special ends.	
\end{enumerate}
\end{dfn}

A finite word  $(\alpha,\ora{p}_1,\ola{q}_1,\ora{p}_2,\ola{q}_2, \ldots  \ora{p}_{k}, \beta)$ will be depicted by a \emph{gluing diagram}
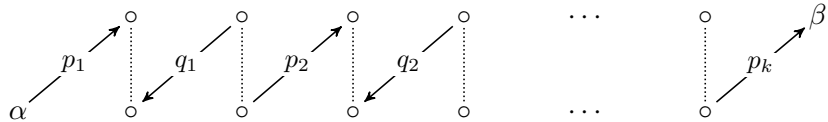
\begin{align}
\label{eq:fin-ab}
\begin{gd}
\& \bt \ar[densely dotted,-]{d} \&\bt  \ar[densely dotted,-]{d} \&\bt  \ar[densely dotted,-]{d} \&\bt  \ar[densely dotted,-]{d} \&  \cdots   \&
\bt \ar[densely dotted,-]{d}  \&   \beta \\
\alpha  \ar[->]{ru}[description]{p_1} \& \bt  \ar[<-]{ru}[description]{q_1} \& \bt \ar[->]{ru}[description]{p_2} \& \bt  \ar[<-]{ru}[description]{q_2} \& \bt \&  \cdots  \& \bt \ar[->]{ru}[description]{p_k} 	\& 
\end{gd}
\end{align}

\begin{dfn}
Let $w = (\diamond,x_1,x_2, \ldots x_{\ell},\diamond)$ be a bispecial word of $\Acat$.
The \emph{primitive root} $v$ of $w$
is the shortest possible subsequence
$v = (x_1,\ldots x_i)$ such that 
$v= (\diamond, v,v^{\op},  \ldots v, \diamond)$
or
$v = (\diamond, v,v^{\op},  \ldots v, v^{\op},\diamond)$,
where $v^{\op} = (x^t_i, \ldots  x^t_2, x^t_1)$.
In this case, the positive integer $m = \frac{\ell}{i}$ is called the \emph{multiplicity} of the primitive root in $w$. 
\end{dfn}

\begin{dfn}\label{dfn:str-sh}
There are three classes of
\emph{strings of $\Acat$ }
defined as follows.
\begin{enumerate}
\item A \emph{usual string} $w$ is given by a usual asymmetric word $w$.
\item  A \emph{special string} $(w, \varepsilon)$ is given by any special word $w$, and a sign $\varepsilon\in \{+,-\}$.
\item A \emph{bispecial string} $(\varepsilon_1,w, \varepsilon_2)$
is given by	any bispecial word $w$, and signs $\varepsilon_1, \varepsilon_2 \in \{+,-\}$.
\end{enumerate}
In each case, $w$ denotes  a finite word of $\Acat$.
\end{dfn}

\begin{rmk}\label{rmk:trans2}
Strings in the sense of Definition~\ref{dfn:str-sh}
translate into  
strings in the sense of Definition~\ref{dfn:strings}
as follows.	
\begin{itemize}
\item A finite word $w = (\alpha,x_1,x_2,\ldots x_{\ell},\beta)$ of $\Acat$ gives rise to a 
word $\breve{w}$ of $\chfXA$, which we define by three parts as follows.
\begin{enumerate}
	\item $\breve{w}$ begins $(\star,\ora{1})$ if $\alpha = \hat{\star}$, and $\breve{w}$ begins with $\alpha$ otherwise.
	\item The intermediate part of $\breve{w}$ is the interior part $(x_1,x_2,\ldots x_{\ell})$ of $w$.
	\item $\breve{w}$ ends with $(\ola{1},\star)$
	if $\beta = \hat{\star}$, and $\breve{w}$ ends with $\beta$ otherwise.
\end{enumerate}
\item Any string $\omega =w$, $(w,\varepsilon)$ or $(\varepsilon_1, w, \varepsilon_2)$ of $\Acat$ 
yields a string 
$\breve{\omega} = (d,\breve{w})$, $(d,\breve{w},\varepsilon)$ respectively $(d,\varepsilon_1, \breve{w},\varepsilon_2)$
of $\chfXA$,
where $d =2$ if the first two numbers of $\breve{w}$ are right-oriented, $d=1$ if $\breve{w}$ only the first number is right-oriented,
and $d=0$ otherwise.
\end{itemize}
\end{rmk}
Translating the strings in Table~\ref{tab:rep-str-bands}
via Remark~\ref{rmk:trans2}
yields the strings in Table~\ref{tab:res-str-bands}.
\begin{table}
\centering
\caption{Strings of $\Acat$}
\label{tab:rep-str-bands}

$
\setlength{\arraycolsep}{3pt}
\begin{array}{|l||l||l|}
\hline
\text{Usual strings } w & \text{Special strings } (w,\varepsilon) & \text{Bispecial strings } (\varepsilon_1,w,\varepsilon_2) \\
\hline
(\star, {\color{blue} \ora{p}_1,\ola{q}_1, \ldots \ora{p}_{k},} \star) 
& ((\star, {\color{blue} \ola{q}_1,\ora{p}_1 \ldots \ola{q}_{k}},\diamond),\varepsilon) 
& (\varepsilon_1,({\color{blue} \diamond, \ora{p}_1, \ola{q}_1, \ldots \ora{p}_k},\diamond),\varepsilon_2) \\
(\wh{\star}, {\color{blue} \ora{p}_1,\ola{q}_1, \ldots \ora{p}_{k}}, \star)
& ((\star, {\color{blue} \ora{p}_1,\ola{q}_1 \ldots \ora{p}_k},\diamond),\varepsilon) 
& (\varepsilon_1,({\color{blue} \diamond, \ola{q}_1,\ora{p}_1, \ldots \ola{q}_k, \ora{p}_{k}},\diamond),\varepsilon_2) \\
(\wh{\star}, {\color{blue} \ora{p}_1,\ola{q}_1, \ldots \ora{p}_{k}, \ola{q}_{k},} \star)
& ((\wh{\star}, {\color{blue} \ora{p}_1,\ola{q}_1 \ldots \ora{p}_k},\diamond),\varepsilon) 
& (\varepsilon_1,({\color{blue}\diamond, \ora{p}_1,\ola{q}_1, \ldots \ora{p}_k, \ola{q}_k},\diamond),\varepsilon_2) \\
(\star, {\color{blue} \ola{q}_1,\ora{p}_1, \ldots \ola{q}_{k}, \ora{p}_{k},} \star)
& ((\star, {\color{blue} \ola{q}_1,\ora{p}_1 \ldots \ola{q}_k, \ora{p}_k},\diamond),\varepsilon) & 
\\
\cline{3-3}
(\star, {\color{blue} \ora{p}_1, \ola{q}_1, \ldots \ora{p}_k, \ola{q}_{k},} \star)
& ((\star, {\color{blue} \ora{p}_1,\ola{q}_1 \ldots \ora{p}_k, \ola{q}_k},\diamond),\varepsilon) & 
\text{Bands }(w,m,\lambda)		
\\
\cline{3-3}
(\wh{\star}, {\color{blue} \ora{p}_1,\ola{q}_1, \ldots \ora{p}_k, \ola{q}_{k},} \wh{\star})
& ((\wh{\star}, {\color{blue} \ora{p}_1,\ola{q}_1 \ldots \ora{p}_k, \ola{q}_k},\diamond),\varepsilon) & 
((\ora{p}_i, \ola{q}_i)_{i \in \Z}, m, \lambda)
\\
\hline
\end{array}$
\end{table}

\subsubsection{From string diagrams to string resolutions}
Let $\omega$ be a string $w$, $(w,\varepsilon)$ or $(\varepsilon_1, w, \varepsilon_2)$ of $\Acat$.

\paragraph{Polarized diagrams for strings}

Let $\breve{\omega}$ denote the string defined by the tranlsation in Remark~\ref{rmk:trans2}. First, we consider the possible outcomes for the string $\breve{\omega}$.
We recall that  $\breve{\omega}$ contains a degree $d \in \{0,1,2\}$.
The string $\breve{\omega}$ gives rise to another degree $d^{\op} \in \{0,1,2\}$
given by
$d^{\op}=2$ if $\breve{w}$ ends with two left-oriented numbers,
$d^{\op}=1$ if $\breve{w}$ ends with one left-oriented number, and $d^{\op}=0$ if $\breve{w}$ does not end with a left-oriented number.
To obtain a uniform notation for all possibilities of $\breve{\omega}$,
we view $\breve{w}$ as being composed from three parts of the following table.
\begin{align}\label{eq:str-pos}
\begin{array}{rrcll}
d & \text{beginning}& \text{middle part} & \text{ending} & d^{\op} \\
2 & (\star,\ora{1}, \ora{p}_1, \, \tikzmark{a1}&  & \tikzmark{d1}\, \ola{q}_k,\ola{1},\star) & 2 \\
1 & (\alpha,\ora{p}_1, \, \tikzmark{a2}& \tikzmark{b} \,  \ola{q}_1, \ora{p}_2,\ldots  \ola{q}_{k-1},\ora{p}_k, \, \tikzmark{c} & \tikzmark{d2} \, \ola{q}_k,\beta) & 1\\
0 & (\alpha,\, \tikzmark{a3} & & \tikzmark{d3}\, \beta) & 0
\end{array}
\begin{tikzpicture}[remember picture, overlay]
\draw[-] ([yshift=0.5ex]pic cs:a1) -- ([yshift=0.5ex]pic cs:b);
\draw[-] ([yshift=0.5ex]pic cs:a2) -- ([yshift=0.5ex]pic cs:b);
\draw[-] ([yshift=0.5ex]pic cs:a3) -- ([yshift=0.5ex]pic cs:b);
\draw[-] ([yshift=0.5ex]pic cs:c) -- ([yshift=0.5ex]pic cs:d1);
\draw[-] ([yshift=0.5ex]pic cs:c) -- ([yshift=0.5ex]pic cs:d2);
\draw[-] ([yshift=0.5ex]pic cs:c) -- ([yshift=0.5ex]pic cs:d3);
\end{tikzpicture}
\end{align}
In particular, the word $\breve{w}$ may have any of the three beginnings, followed by a (possibly empty) middle part, and then any of the three endings.

In the next step, we introduce 
a pair of vertices  $(a_i,b_i)$ for each oriented number $\ora{p}_i$ in $\breve{w}$, and a pair of vertices $(v_i,u_i)$ for each oriented number $\ola{q}_i$ in $\breve{w}$
as well as certain pairs 
in cases $d=2$ and 
$d^{\op}=2$, as the table below indicates.
\begin{align}
\label{eq:str-pos2}
\begin{array}{rcl}
(\star,u_0 ,a_1,b_1, \, \tikzmark{a1'}&  & \tikzmark{d1'}\,  v_k, u_k, a_{k+1},\star)  \\
(a_1,b_1, \, \tikzmark{a2'}& \tikzmark{b'} \, v_1, u_1,a_2, b_2,
\ldots v_{k-1} u_{k-1}, a_k, b_k,  \, \tikzmark{c'} & \tikzmark{d2'} \, v_k,u_k) \\
(\emptyset\, \tikzmark{a3'} & & \tikzmark{d3'}\, 
\emptyset) 
\end{array}
\begin{tikzpicture}[remember picture, overlay]
\draw[-] ([yshift=0.5ex]pic cs:a1') -- ([yshift=0.5ex]pic cs:b');
\draw[-] ([yshift=0.5ex]pic cs:a2') -- ([yshift=0.5ex]pic cs:b');
\draw[-] ([yshift=0.5ex]pic cs:a3') -- ([yshift=0.5ex]pic cs:b');
\draw[-] ([yshift=0.5ex]pic cs:c') -- ([yshift=0.5ex]pic cs:d1');
\draw[-] ([yshift=0.5ex]pic cs:c') -- ([yshift=0.5ex]pic cs:d2');
\draw[-] ([yshift=0.5ex]pic cs:c') -- ([yshift=0.5ex]pic cs:d3');
\end{tikzpicture}
\end{align}
The presentation of the word $\breve{w}$ in \eqref{eq:str-pos} and its sequence of vertices  can be visualized by a diagram divided into three parts as shown below.
\begin{align}
\label{eq:str-ab2}
\scalebox{0.85}{$
\begin{array}{rcl}
	\text{beginning} & \text{middle part} & \text{ending} \\
	\begin{gd}
		\& P_{u_0} \ar[densely dotted,-]{d} \& P_{b_1}  \ar[densely dotted,-]{d} 
		\\
		\underset{2}{P_{\star}}  \ar[->]{ru}[description]{1} \& \underset{1}{P_{a_1}}  \ar[->]{ru}[description]{p_1} \& \mathstrut
	\end{gd}
	&
	\multirow{3}{*}{$
		\begin{gd}
			\\
			\& P_{u_1} \ar[densely dotted,-]{d}  
			\ar[phantom]{rd}[description]{\cdots}
			\ar[densely dotted,-]{d} 
			\&[1cm] 
			\mathstrut 
			\ar[densely dotted,-]{d} 
			\& \smash[b]{\underset{0}{P_{b_k}}}   \\
			\underset{0}{P_{v_1}} \ar[<-]{ru}[description]{q_1} 
			\&	
			\mathstrut\& 
			\underset{1}{P_{a_k}} \ar{ru}[description]{p_k} \&  
		\end{gd}$
	}
	&
	\begin{gd}
		\mathstrut \ar[densely dotted,-]{d}
		\& P_{u_k} \ar[densely dotted,-]{d}  \&  \smash[b]{\underset{2}{P_{\star}}} 
		\\
		\underset{0}{P_{v_k}} \ar[<-]{ru}[description]{q_k} \& \underset{1}{P_{a_{k+1}}} \ar[<-]{ru}[description]{1} 	\& 
	\end{gd}
	\\	
	\\
	\begin{gd}
		\&  \& P_{b_1}  \ar[densely dotted,-]{d} 
		\\
		\&	\underset{1}{P_{a_1}}  \ar[->]{ru}[description]{p_1} \& \mathstrut
	\end{gd}
	&
	&
	\begin{gd}
		\mathstrut \ar[densely dotted,-]{d}  \& \smash[b]{\underset{1}{P_{u_k}}}  \&   \\
		\underset{0}{P_{v_k}} \ar[<-]{ru}[description]{q_k} \& 	\& 
	\end{gd}
	\\
	\\
	\emptyset
	&&
	\emptyset
\end{array}$}
\end{align}
The vertices above are defined as follows.
\begin{itemize}[leftmargin=15pt]
\item Assume that $\omega$ is usual or special.
Any pair of intermediate vertices 
$(u_i,a_{i+1})$ 
or  $(b_i, v_i)$ 
is set to $(+,-)$.
\begin{itemize}
\item If $\omega$ is usual,
the remaining outer vertices are set to $\star$.
\item If $\omega$ is special $(w,\varepsilon)$ with a special left end, the first vertex $a_1$ respectively $v_1$ is set to $\varepsilon$, and the last vertex to $\star$.
\item  If $\omega$ is special with special right end, the first one is $\star$ and the last one $\varepsilon$.
\end{itemize}
\item 
Assume that $\omega$ is bispecial.
Let $r$ denote the primitive root of $w$ and $m\in \N$ its multiplicity, that is,
$w = (r,r^{\op},  \ldots r)$ if $m$ is odd, respectively
$w = (r,r^{\op},  \ldots r, r^{\op})$ if $m$ is even,
where $r^{\op} = (n_i, \ldots  n_1)$.
Any subsequence of vertices corresponding to a copy of the root $r$ is given by
$(\varepsilon_1, +,-, \ldots  +,-,\varepsilon_2)$, while  any subsequence of vertices corresponding to a copy of $r^{\op}$
is set to $(\ol{\varepsilon}_2,-,+,\ldots  -,+,\ol{\varepsilon}_1)$.
\end{itemize}
Similar to the case of bands, gluing edges connect  vertices of type $+$ and type $-$, but pairs $(-,+)$ when read from top to bottom may appear for bispecial strings.

\begin{ex}
Assume that $\omega$ is the bispecial string $(\varepsilon_1,w,\varepsilon_2)$ with
$w = (\ora{1},\ola{2},\ora{1},\allowbreak \ola{1},\ora{2},\ola{1},\allowbreak \ora{1},\ola{2},\ora{1},\ola{1},\ora{2},\ola{1})$.
Then $w$ has primitive root  $(\ora{1},\ola{2},\ora{1})$ with multiplicity four.
The sequence of vertices for the word $w$ 
is given by two copies of
	$(
	{ \varepsilon_1}, +,-,+,-,{\varepsilon_2}
	,
	\allowbreak
	{ \ol{\varepsilon}_2}, -,+,-,+,{ \ol{\varepsilon}_1})$.
\end{ex}

\paragraph{Orienting vertical edges}
The word $w$ 
yields a sequence of the form
$w^{\updownarrow} = (
x_1 \kappa_1 x_2 \kappa_2 \ldots  x_{\ell})$
with $\kappa_i \in \{\uparrow, \downarrow\}$ for each index $1 \leq i \leq \ell$
as follows.
\begin{dfn}
Let $w = (\alpha,x_1,x_2,\ldots x_\ell,\beta)$ be a finite word of $\Acat$. 
\begin{itemize}
\item Set $v = (y_1, y_2, \ldots y_{\ell})$  
where $y_1$ is obtained from $x_1$ by changing $n_1$ to $n_1 + \frac{1}{2}$ if $(\alpha,x_1) = (\star,\ola{n}_1)$
respectively to $n_1 - \frac{1}{2}$ if $(\alpha,x_1) = (\star,\ora{n}_1)$,
$y_i = x_i$ for $1 < i <\ell$,
and
$y_{\ell} = n_{\ell} - \frac{1}{2}$ if $(x_{\ell},\beta) = (\ola{n}_\ell,\star)$, respectively $n_{\ell} + \frac{1}{2}$
if $(x_{\ell},\beta) = (\ora{n}_\ell,\star)$.
\item Set $\overline{w} = (\star,v,\star)$ if $w$ is usual, 
$\overline{w} = (\star,v, v^{\op},\star)$
if $w$ is special and $\beta=\diamond$,
$\overline{w} = (\star,v^{\op},v,\star)$
if $w$ is special $\alpha=\diamond$,
and $\overline{w} = (y_i)_{i \in \Z}$ with $y_{\ell+1+a} = y_{\ell-a}^t$
for any $0 \leq a < \ell$ and $y_{2\ell+i} = y_i$ for any integer $i \in \Z$ 
if $w$ is bispecial.
\end{itemize}	
\end{dfn}
\begin{rmk}\label{rmk:kappa2}
The vertical arrows in $w^{\updownarrow}$ are determined as follows.
For each index $1 \leq i < \ell$
set $\kappa_{i}=\uparrow$ if
and only if one of the following occurs.
\begin{itemize}
\item 
$\overline{w} = (\hat{\star}, v, v^{\op} \ldots)$ where $v$ ends with $y_i$.
\item $\overline{w} =  
(\ldots y_{i-j}, v, v^{\op}, y_{i+1+j} \ldots)$
where $v$ ends with $y_i$,
and 
\eqref{eq:or-ru} holds
with the pair
$({x}_{i-j},{x}_{i+1+j})$			replaced by $(y_{i-j},y_{i+1+j})$.	
\end{itemize}
\end{rmk}
The oriented word $w^{\updownarrow}$
determines 
the first differential of the string complex $\CX(\omega)$ described below.

\paragraph{The first differential}
Similar to the gluing diagram,
the first differential $\partial_1$ of the string complex $\CX(\omega)$
can be divided into several parts.
More precisely, it is given by
\begin{align}
\label{eq:olap}
\partial_1 = 
\scalebox{0.5}{$
\begin{bNiceArray}[columns-width = auto, cell-space-limits=5pt, margin]{cccccccccccccccc}
\Block[draw=blue]{4-4}{\scalebox{2.0}{$\partial'_1$}} &&&& \Block[draw=blue]{2-9}<\large>{0}  &  &&& &&&& & \Block[draw=blue]{2-3}<\large>{0}&  &  \\
&&&& & &&& &&&& &&  &  \\
&&&
* 
\Block[draw=blue]{10-10}{\scalebox{2.0}{$\partial''_1$}} & * &* & 0 & \Cdots & &&&& 0 &\Block[draw=blue]{8-3}<\large>{0} && \\
&&&*&   &&\Ddots& \Ddots & &&&& \Vdots &&& \\
\Block[draw=blue]{8-3}<\large>{0}&&&*& &&&& &&&& &&& \\
&&& 0 &\Ddots &&&& &&&& &&& \\
&&& \Vdots &\Ddots &&&& &&&& &&& \\
&&&& &&&& &&&& &&& \\
&&&& &&&& &&&& 0 &&& \\
&&&& &&&& &&&& *&&& \\
&&&& &&&& &&&& * \Block[draw=blue]{4-4}{\scalebox{2.0}{$\partial'''_1$}} &&& \\
&&& 0 & \Cdots  &&&& &0&*&*& 
* 
&& &  \\
\Block[draw=blue]{2-3}<\large>{0}&& &\Block[draw=blue]{2-9}<\large>{0}& &&&& &&&& & && \\
&&&& &&&& &&&& &&& 
\CodeAfter
\tikz \draw [blue,very thick] (1-|1) rectangle (5-|5) ;
\tikz \draw [blue,very thick] (3-|4) rectangle (13-|14) ;
\tikz \draw [blue,very thick] (11-|13) rectangle (16-|17) ;
\end{bNiceArray}$}
\end{align}
where $\partial'_1$ 
and $\partial''_1$, as well as $\partial''_1$ and $\partial'''_1$ overlap in a $2$-by-$1$-minor.
Depending on the degree $d$ of $\omega$, 
the top left block $\partial'_1$ is given by one of the matrices:
\begin{align*}
\begin{array}{c}
d=2 \\
\scalebox{0.75}{$
\begin{bNiceArray}{cccc}[first-row, last-col
, columns-width = auto
]
\CodeBefore
\chessboardcolors{blue!25}{red!25}
\Body
P_{u_0} & P_{a_1} & P_{u_1} & P_{a_2}  & 
\\
t^{p_1} & t^{p_1} & \alpha_1 t^{q_1} & \gamma_{1} t^{q_1}  & P_{b_1} \\
\ol{\alpha}_1 t^{p_1} & \ol{\alpha}_1 t^{p_1} & t^{q_1} & \beta_1 t^{q_1} 
&P_{v_1} \\
0 & 0 & \ol{\beta}_1 t^{p_2} & { t^{p_2 } }& P_{b_2}  \\ 
0 & 0 & \ol{\gamma}_1 t^{p_2} & {\ol{\alpha}_2 t^{p_2 } }& P_{v_2}
\CodeAfter
\tikz \draw [blue] (3-|4) rectangle (5-|5) ; 
\end{bNiceArray}$}
\end{array}
&&
\begin{array}{c}
d=1\\
\scalebox{0.75}{$
\begin{bNiceArray}{ccc}[first-row, last-col
, columns-width = auto
]
\CodeBefore
\chessboardcolors{red!25}{blue!25}
\Body
P_{a_1} & P_{u_1} & P_{a_2} & 
\\
t^{p_1} & \alpha_1 t^{q_1} & \gamma_{1} t^{q_1}  & P_{b_1} \\
\ol{\alpha}_1 t^{p_1} & t^{q_1} & \beta_1 t^{q_1}  
&P_{v_1} \\
0 & \ol{\beta}_1 t^{p_2} &{ t^{p_2 } }& P_{b_2}  
\\ 
0 & \ol{\gamma}_1 t^{p_2} & {\ol{\alpha}_2 t^{p_2 } }& P_{v_2} 
\CodeAfter
\tikz \draw [blue] (3-|3) rectangle (5-|4) ;
\end{bNiceArray}$}
\end{array}
&&
\begin{array}{c}
d=0 \\
\\
\scalebox{0.75}{$
\begin{bNiceArray}{cc}[first-row, last-col
, columns-width = auto
]
\CodeBefore
\chessboardcolors{red!25}{blue!25}
\Body
P_{u_1} & P_{a_2}  & 
\\
t^{q_1} & \beta_1 t^{q_1} & P_{v_1} \\
\ol{\beta}_1 t^{p_2} & { t^{p_2 } }&  P_{b_2} \\
\ol{\gamma}_1 t^{p_2} & {\ol{\alpha}_2 t^{p_2 } }& P_{v_2} 
\CodeAfter
\tikz \draw [blue] (2-|2) rectangle (4-|3) ;
\end{bNiceArray}$}
\end{array}
\end{align*}
The central block $\partial''_1$ is given by the following pentadiagonal square matrix:
\begin{align}\label{eq:str-m2}
\scalebox{0.75}{$
\begin{bNiceArray}{cccccccccc}[first-row, last-col
, columns-width = auto
]
\CodeBefore
\chessboardcolors{red!10}{blue!10}
\foreach \i in {1,...,10}{
\pgfmathtruncatemacro{\jmin}{max(1,\i-2)}
\pgfmathtruncatemacro{\jmax}{min(10,\i+2)}
\foreach \j in {\jmin,...,\jmax}{
\pgfmathtruncatemacro{\sum}{mod(\i-\j,2)}
\ifnum\sum=0
\cellcolor{red!25}{\i-\j}
\else
\cellcolor{blue!25}{\i-\j}
\fi
}
}
\Body
P_{a_2} & P_{u_2} & P_{a_3} & P_{u_3} &&& P_{a_{k-2}} & P_{u_{k-2}} & P_{a_{k-1}} & P_{u_{k-1}} 
\\
{ t^{p_2 } } & \alpha_2 t^{q_2} & \gamma_{2} t^{q_2} & 0 &\Cdots &\Cdots&\Cdots&\Cdots&\Cdots& 0  &P_{b_2} \\
{\ol{\alpha}_2 t^{p_2}} & t^{q_2} & \beta_2 t^{q_2} & 0 &\Ddots&&& &  & 
&P_{v_2}  \\
0 & \ol{\beta}_2 t^{p_3} & t^{p_3}& \alpha_3 t^{q_3} & \Ddots &\Ddots&&& &  &P_{b_3}  \\
0 & \ol{\gamma}_2 t^{p_3} & \ol{\alpha}_3 t^{p_3} & t^{q_3} &\Ddots &\Ddots&\Ddots&& & \Vdots & P_{v_3} \\
\Vdots		&\Ddots &\Ddots &\Ddots & & \Ddots & &  \Ddots & \\
\Vdots	&&\Ddots &\Ddots &\Ddots &\Ddots &  & \Ddots &  \Ddots &  \Vdots & \\
\Vdots &&&&		\Ddots & \Ddots &		
t^{p_{k-2}}
& 
*
&
*
&  0 & P_{b_{k-2}}  \\
&  &	&&& \Ddots & 	
*
& 
t^{q_{k-2}} 
& 
*
&  0 & P_{v_{k-2}}  \\
&	  	&&&& &		0 & 
*
& t^{p_{k-1}}& 
\scalebox{0.7}{$
{ \alpha_{k-1} t^{q_{k-1}}}$}
& P_{b_{k-1}} \\
0 & \Cdots &  && &&0 & 
* 
& 
*
& { t^{q_{k-1}}} &	P_{v_{k-1}} 
\CodeAfter
\tikz \draw [blue] (1-|1) rectangle (3-|2) ;
\tikz \draw [blue] (9-|10) rectangle (11-|11) ;	 	 
\end{bNiceArray}$}	\end{align}
Similar to the first block,
the bottom right block $\partial'''_1$ has one of three forms, which is determined by $d^{\op}$ associated to $\omega$:
\begin{align*}
\begin{array}{c}
d^{\op}=2\\ 
\scalebox{0.7}{$
\begin{bNiceArray}{cccc}[first-row, last-col
, columns-width = auto
]
\CodeBefore
\chessboardcolors{red!25}{blue!25}
\Body
P_{u_{k-1}} & P_{a_k} & P_{u_{k}} & P_{a_{k+1}}  & 
\\
\scalebox{0.9}{${ \alpha_{k-1} t^{q_{k-1}}}$} & \scalebox{0.9}{$\gamma_{k-1} t^{q_{k-1}}$} & 0 & 0  
&P_{b_{k-1}} \\
{ t^{q_{k-1}}} & \scalebox{0.9}{$\beta_{k-1} t^{q_{k-1}}$} & 0 & 0  
&P_{v_{k-1}} \\
\scalebox{0.9}{$\ol{\beta}_{k-1} t^{p_k}$} &  t^{p_k} & \alpha_k t^{q_k} & {\alpha}_k t^{q_k} &  P_{b_k} \\
\scalebox{0.9}{$\ol{\gamma}_{k-1} t^{p_k}$} & \ol{\alpha}_k t^{p_k} & t^{q_k} & t^{q_k}& P_{v_k}
\CodeAfter
\tikz \draw [blue] (1-|1) rectangle (3-|2) ;	 
\end{bNiceArray}$}
\end{array}
&&
\begin{array}{c}
d^{\op}=1\\
\scalebox{0.7}{$
\begin{bNiceArray}{ccc}[first-row, last-col
, columns-width = auto
]
\CodeBefore
\chessboardcolors{red!25}{blue!25}
\Body
P_{u_{k-1}} & P_{a_k} & P_{u_k} &  
\\
\scalebox{0.9}{${ \alpha_{k-1} t^{q_{k-1}}}$} & \scalebox{0.9}{$\gamma_{k-1} t^{q_{k-1}}$} & 0   
&P_{b_{k-1}} \\
{ t^{q_{k-1}}} & \scalebox{0.9}{$\beta_{k-1} t^{q_{k-1}}$} & 0   
&P_{v_{k-1}} \\
\scalebox{0.9}{$\ol{\beta}_{k-1} t^{p_k}$} &  t^{p_k} & \alpha_k t^{q_k}  &  P_{b_k} \\
\scalebox{0.9}{$\ol{\gamma}_{k-1} t^{p_k}$} & \ol{\alpha}_k t^{p_k} & t^{q_k}&  P_{v_k}
\CodeAfter 
\tikz \draw [blue] (1-|1) rectangle (3-|2) ;
\end{bNiceArray}$}
\end{array}
&&
\setlength{\arraycolsep}{0pt}
\begin{array}{c}
d^{\op}=0 \\
\scalebox{0.7}{$
\begin{bNiceArray}{cc}[first-row, last-col
, columns-width = auto
]
\CodeBefore
\chessboardcolors{red!25}{blue!25}
\Body
P_{u_{k-1}} & P_{a_k} & 
\\
\scalebox{0.9}{${ \alpha_{k-1} t^{q_{k-1}}}$} & \scalebox{0.9}{$\gamma_{k-1} t^{q_{k-1}}$}  
&P_{b_{k-1}} \\
{ t^{q_{k-1}}} & \scalebox{0.9}{$\beta_{k-1} t^{q_{k-1}}$}   
&P_{v_{k-1}} \\
\scalebox{0.9}{$	\ol{\beta}_{k-1} t^{p_k } $}&  t^{p_k } & P_{b_k} \\
\CodeAfter
\tikz \draw [blue] (1-|1) rectangle (3-|2) ;
\end{bNiceArray}$}
\\
\mathstrut
\end{array}
\end{align*}
In case $d=d^{\op}=2$,
there are pairs of binary coefficients $(\alpha_i,\ol{\alpha}_i)$, $(\beta_i,\ol{\beta}_i) \in \{(1,0),(0,1)\}$ for each index $1 \leq i \leq k$
and coefficients $\gamma_{i}, \ol{\gamma}_i \in \{0,1\}$ with $1 \leq i < k$, which are defined by the prescriptions in \eqref{eq:coeff}.
In the general case, these coefficients are defined by the same prescriptions, for certain smaller index sets.

\paragraph{The second differential}
In case $d = 2$ or $d^{\op}=2$ the complex $\CX(\omega)$ has the form $(\begin{td} \X_2 \ar{r}{\partial_2} \& \X_1 \ar{r}{\partial_1} \& \X_0 \end{td})$
with the following three possibilities for $\X_2$ and $\partial_2$.
$$
\begin{array}{ccccc}
d = 2 = d^{\op} && d=2 \neq d^{\op} && d \neq 2 = d^{op} \\
\scalebox{0.75}{$	\begin{bNiceArray}{cc}[first-row, last-col, columns-width = auto
]
\CodeBefore
\chessboardcolors{blue!25}{red!25}
\Body
P_{\star} & P_{\star} & 
\\
\phantom{-}t & 0 & P_{u_0} \\
-t & 0&  P_{a_1} \\
0 & 0 &  P_{u_1} \\
\Vdots & \Vdots & \Vdots \\
& & \\
0 & 0 & P_{a_k}  \\
0 & -t & P_{u_k } \\
0 & \phantom{-}t &  P_{a_{k+1}}\\
\end{bNiceArray} $}
& &
\scalebox{0.75}{$
\begin{bNiceArray}{c}[first-row,last-col]
\CodeBefore
\chessboardcolors{blue!25}{red!25}
\Body
P_{\star} &  \\
\phantom{-}t & P_{u_0} \\
-t & P_{a_1}\\
0 & P_{u_1} \\
\Vdots & \mathstrut \\
\Vdots & \mathstrut \\
0 & \mathstrut
\end{bNiceArray}$}
&&
\scalebox{0.75}{$
\begin{bNiceArray}{c}[first-row,last-col]
\CodeBefore
\chessboardcolors{blue!25}{red!25}
\Body
P_{\star} & \mathstrut  \\
0 & \mathstrut \\
\Vdots & \mathstrut \\
\Vdots & \mathstrut \\
0 & P_{a_k} \\
-t & P_{u_k} \\
\phantom{-}t & P_{a_{k+1}}
\end{bNiceArray}$}
\end{array}
$$

\subsection{Bijection theorem}
Isomorphisms of certain band complexes or string complexes motivate the next notion.
\begin{dfn}
We consider the smallest equivalence relation on the set of bands and strings of $\Acat$ such that the following holds.
\begin{itemize}
\item Any usual string $w$ is equivalent to $w^{\op}$.
\item Any special  string $(w,\varepsilon)$ is equivalent to $(w^{\op},\varepsilon)$.
\item Any bispecial string $(\varepsilon_1, w,\varepsilon_2)$ is equivalent to $(\varepsilon_2, (v^{\op})^m,\varepsilon_1)$, 
where $v$ denotes the primitive root of $w$ and $m$ its multiplicity.
\item Any band $(w,m,\lambda)$ is equivalent to
$(w_{1}, m, \lambda)$ as well as  $(w^{\op},m,\lambda^{\mp})$ 
with $\lambda^\mp = \frac{1}{\lambda}$ if $w$ is asymmetric and $\lambda^{\mp} = \lambda$ otherwise.
\end{itemize}
\end{dfn}
Any band of $\Acat$ is equivalent to a band 
$(w,m,\lambda)$ with $\dot{w} =(\ora{p}_1,\ola{q}_1,\ldots  \ora{p}_k,\ola{q}_k)$.
Any string of $\Acat$ is equivalent to a string from
Table~\ref{tab:rep-str-bands}.

\begin{dfn}
For any band or string $\omega$ of $\Acat$ the \emph{band representation} respectively \emph{string representation} $M(\omega)$ is defined as the quiver representation of the $\A$-module $H_0(\CX(\omega))$, with the complex $\CX(\omega)$ defined by Subsection~\ref{sec:proj-res1}
respectively~\ref{sec:proj-res2}.
\end{dfn}

\begin{thm} \label{thm:prbij} The following statements hold.
\begin{enumerate}
\item \label{thm:prbijA} For any band or string $\omega$ of $\Acat$ the complex $\CX(\omega)$ is a minimal projective resolution of an indecomposable finite-dimensional $\A$-module.
\item \label{thm:prbijB}	The assignment $\begin{td} \omega \ar[mapsto]{r} \& M(\omega) \end{td}$ yields a bijection 
$$
\begin{cd}
\{\text{bands and strings of }\Acat  \}\!\big/_{\textstyle \approx}
\ar{r}{\sim}
\&
\ind \Acat
\end{cd}
$$
between the set of equivalence classes of bands and strings of $\Acat$ and the set of isomorphism classes of indecomposable objects in $\Acat$.
\end{enumerate}
\end{thm}
\begin{proof}
For any band or string $\omega$ of $\Acat$,
the complex $\CX(\omega)$ introduced in the previous subsection is precisely 
the complex $\CX(\breve{\omega})$ defined in Section~\ref{sec:gluing}
where $\breve{\omega}$ is the band or string translation of $\omega$. With previous results, there is a commutative diagram
\begin{align}\label{eq:prbij}
\begin{cd}
\{\text{bands and strings of }\Acat  \}\!\big/_{\textstyle \approx}
\ar[yshift=0pt]{r} \ar{d}{\text{Remarks }\ref{rmk:trans1},\, \ref{rmk:trans2}}[swap]{\breve{(-)}}
\&
\text{[}\Db{\Acat}\text{]}
\& \omega \ar[mapsto]{r} \ar[mapsto]{d}\& \CX(\omega)   \\
\{\text{bands and strings of }\chfXA  \}\!\big/_{\textstyle \approx} 
\ar[yshift=0pt]{r}{\sim}[swap]{\eqref{eq:bij-perf-fd2}}\&
\ind \Db{\Acat}  \ar[hookrightarrow]{u}\& \breve{\omega} \ar[mapsto]{r} \& \CX(\breve{\omega}) 
\ar[mapsto,dashed]{u}
\end{cd}
\end{align}	
where $[\DbRep{\A}]$ denotes the set of isomorphism classes of objects in $\DbRep{\A}$. 
Since
the bottom horizontal map is a bijection,
it follows that the complex $\CX(\omega)$ is indecomposable for any band or string $\omega$ of $\Acat$. 

It is straightforward to check that the 
left vertical map $\breve{(-)}$
is injective. Its image is given by the set of equivalence classes of bands and strings of $\chfXA$
whose complexes are projective resolutions by Theorem~\ref{thm:proj-res}.
Therefore, 
the codomain $\Db{\Acat}$ in diagram~\eqref{eq:prbij} may be replaced by the set of isomorphism classes of minimal projective resolutions of indecomposable objects of $\Acat$, and the remaining claims follow.
\end{proof}
The last result allows to recover the classification of objects in $\Rep{\A}$
which can be generated by only one element. 
\begin{cor}\label{cor:cyclic}
Table~\ref{tab:cyc-base} lists
all cyclic representations in $\Rep{\A}$
up to isomorphism.
\end{cor}
\begin{proof}
By Theorem~\ref{thm:prbij} any gluing diagram with only one indecomposable projective at degree zero is given by a diagram of one of the forms in~\eqref{eq:cyc-diag1}.
These gluing diagrams lead to the minimal projective resolutions in~\eqref{eq:cyc-res2}, 
which are quasi-isomorphic to the syzygy resolutions
of the form $(\begin{td} \syz(M) \ar[hookrightarrow]{r} \& P_0 \end{td})$ described in Table~\ref{tab:cyc-base}. 
\end{proof}

\subsection{Cyclification of a band or string resolution}
\label{subsec:cyclif}
For any band or string, we define a complex with the same Betti numbers 
as the associated band or string complex.

\begin{dfn}\label{dfn:cyc}
For any band or string $\omega$ of $\Acat$, 
the cyclification $\CX^{cyc}(\omega)$ is given by the 
the complex 
obtained from 
$\CX(\omega)$
by setting all orientation-dependent coefficients of its differential $\partial_1$ in \eqref{eq:l-bd-diff}
respectively~\eqref{eq:olap} to zero.
The module
$M^{cyc}(\omega) = H_0(\CX^{cyc}(\omega))$ is called the \emph{cyclification} of the representation $M(\omega)$.
\end{dfn}
We describe the cyclification resolution for all bands and strings up to equivalence.
\begin{enumerate}
\item Let $\omega$ be a band $(w,m,\lambda)$
such that $\dot{w} = (\ora{p}_1,\ola{q}_1,\ldots \ora{p}_{k},\ola{q}_{k})$.
Then $\CX^{cyc}(\omega) = \bigl(\bigoplus_{i=1}^k
\CX(-,p_i,+) \oplus \CX(+,q_i,-)\bigr)^m$.
\item Let $\omega$ be a string.
We recall that its underlying word translates into a finite word $\breve{w}$ as in \eqref{eq:str-pos} and gives rise to a sequence of vertices~\eqref{eq:str-pos2}.
Subwords of $\breve{w}$ induce projective resolutions of the form ~\eqref{eq:cyc-res2}
of finite-dimensional cyclic $\A$-modules as follows.
\begin{itemize}
\item A beginning $(\ora{1},\ora{p}_1)$ induces the resolution
$\CX(\hat{\star},p_1,b_1)$.
\item Similarly, an
end $(\ola{q}_k,\ola{1})$ yields 
$\CX(\hat{\star},q_k,a_{k+1})$.
\item 
Any of the remaining right-oriented numbers $\ora{p}_i$ gives rise to a summand $\CX(a_i,p_i,b_i)$, and any other left-oriented number $\ola{q}_i$ corresponds to a summand $\CX(u_i, q_i, v_i)$.
\end{itemize} 
The complex $\CX^{cyc}(\omega)$ is given by the direct sum of these  projective resolutions.
\end{enumerate}
In different terms, the cyclification $\CX^{cyc}(\omega)$
is the direct sum of projective resolutions
obtained from the gluing diagram \eqref{eq:per-ab2} respectively \eqref{eq:str-ab2} of $\omega$ 
by deleting each vertical edge between any two consecutive arrows if none of the arrows starts a projective degree two
and replacing each connected gluing diagram of one of the forms in \eqref{eq:cyc-diag1}
or \eqref{eq:cyc-diag2}
by the corresponding projective resolution in \eqref{eq:cyc-res2}. Therefore, $M^{cyc}(\omega)$
is a direct sum of cyclic finite-dimensional $\A$-modules.
This allows to view the complex $\CX(\omega)$ informally as `gluing' of the direct summands of $\CX^{cyc}(\omega)$ as  demonstrated in Subsection~\ref{subsec:motiv-cyc}. 

For an object $M$ of $\Acat$ 
we denote 
its \emph{dimension vector}
by 
$\uldim M = (\dim M_v)_{v \in \{-,\star,+\}}$.
\begin{rmk}\label{rmk:JH-mult}
As $\Acat$ is a finite length category, 
for each vertex $v \in Q_0$
the positive integer $\dim M_v$ is equal to  and can be interpreted as  
the multiplicity of the simple module $S_v$
in any Jordan--Hölder filtration of the $\A$-module $M$.
Therefore, the dimension vector 
records the \emph{Jordan--Hölder multiplicities} 
of $M$.
\end{rmk}
\begin{prp} \label{prp:cyc}
For any band or string $\omega$ of $\Acat$
the module $M^{cyc}(\omega)$ is a direct sum of cyclic $\A$-modules 
such that
$\uldim M(\omega) = \uldim M^{cyc}(\omega)$. 
\end{prp}
\begin{proof}
By definition, the complex $\CX^{cyc}(\omega)$ 
is a direct sum of complexes 
of the form \eqref{eq:cyc-res2}, which are minimal projective resolutions of cyclic modules.
As the complexes $\CX(\omega)$ and $\CX^{cyc}(\omega)$ have identical projective at each degree, it holds  that
$\AI \oA \CX(\omega) \cong \AI \oA \CX^{cyc}(\omega)$.
The complex  $\B \oB \CX(\omega)$ is isomorphic to $ \B \oB \CX^{cyc}(\omega)$ by Lemma~\ref{lem:norm-complex}.
Therefore,	Corollary~\ref{cor:dim} can be applied to deduce that $\uldim M(\omega) = \uldim M^{cyc}(\omega)$.
\end{proof}
The last result yields a useful tool to determine the dimension vector of the homology of any band or string resolution, which does not require a full description of the differentials.
\begin{ex}\label{ex:dim}
Let   
$\omega$
be the usual string $w=(\hat{\star},1,1,\star)$.
Then the projective resolution $\CX^{cyc}(\omega)$
is given by 
$(\begin{td} P_\star \ar{r}
\& P_+ \oplus P_- \ar{r}{
\left(\begin{smallmatrix}
	{t}&
	{t}	
\end{smallmatrix}\right)	} \& P_+
\end{td})$ 
$\oplus$
$(\begin{td} P_\star \ar{r}{t} \& P_- \end{td})$,
and its homology
$M^{cyc}(\omega)$ by $(P_+/\rad^2 P_+) \oplus S_-$.
Thus $\uldim M(\omega) = \uldim M^{cyc}(\omega) = 
\begin{pmatrix}
1, 1, 1 
\end{pmatrix}.$
\end{ex}

\section{Invariants and functorial properties of band and string representations}
\label{sec:invariants}

With the description of band and string resolutions in hand, 
we describe their behaviour under the involution $\inv$ and the contragredient duality $\cdual{(-)}$ in Theorem~\ref{thm:dualities}
and identify the self-dual resolutions (Corollary~\ref{cor:self-dual}). Moreover, we 
compute
projective and injective dimension, top, socle, dimension vector, and Euler characteristic of the indecomposable band or string representation $M(\omega)$ in $\Acat$
directly from the numerical data of $\omega$ (Theorem~\ref{thm:inv}, Table~\ref{tab:inv-str} and Corollary~\ref{cor:chi}).
In different terms, these invariants depend only on the much coarser cyclification $M^{cyc}(\omega)$, that is, 
on the `cyclic ingredients' in the gluing of $M(\omega)$.

\subsection{Functorial properties of indecomposable representations}
Previous considerations in \eqref{E:reflection},\eqref{E:dual} and \eqref{E:LieDuality} provide 
self-inverse equivalences
$$\inv,\dual{(-)}, \cdual{(-)}\colon \begin{cd} \Acat \ar{r}{\sim} \& \Acat \end{cd}$$
such that $\inv$ is covariant, whereas $\dual{(-)}$ and $\cdual{(-)} = \dual{(-)} \circ \inv = \inv \circ \dual{(-)}$ are contravariant.
The action of these functors on a representation $M$ corresponds 
to a horizontal flip, taking the $\kk$-linear dual $(-)^* = \Hom_{\kk}(-,\kk)$, respectively both operations
as the following diagram shows.
The goal of this subsection is to describe the operation of these three functors on bands and strings of $\Acat$.
\begin{align*}
\begin{array}{rcccrcc}
	M &=& \left[\begin{gqrep}{M_+}{M_\star}{M_-}{B_+}{B_-}{A_+}{A_-}\end{gqrep}\right] 
	&\begin{td}\mathstrut \ar[mapsto]{r}{\dual{(-)}} \& \mathstrut\end{td} &
	\dual{M} &=& 
	\left[
	\begin{gqrep}{M^*_+}{M^*_\star}{M^*_-}{A^*_+}{A_-^*}{B_+^*}{B^*_-}\end{gqrep}
	\right]
	\\
	&&
	\begin{td}\phantom{\cdot} \ar[mapsto]{d}{\inv} \\
		\phantom{\cdot} 
	\end{td}
	& & & & \begin{td}\phantom{\cdot} \ar[mapsto]{d}{\inv} \\
		\phantom{\cdot} 
	\end{td} \\
	\inv(M) &=& 
	\left[\begin{gqrep}{M_-}{M_\star}{M_+}{B_-}{B_+}{A_-}{A_+}\end{gqrep}\right] 
	&\begin{td}\mathstrut \ar[mapsto]{r}{\dual{(-)}} \& \mathstrut\end{td}  &
	\cdual{M} &=& 
	\left[\begin{gqrep}{M^*_-}{M^*_\star}{M^*_+}{A^*_-}{A_+^*}{B_-^*}{B^*_+}\end{gqrep}\right]
\end{array}
\end{align*}

For any $n \in \N_0$ we denote $n^{\pm} = n \pm 1$ below.
\begin{dfn}
Any word $w$ of $\Acat$
gives rise a word $\dual{w}$ of $\Acat$ 
as follows.
\begin{itemize}
	\item We consider certain exceptional words first.	
	\begin{itemize}
		\item If $w$ or $w^{\op} = (\star,\ora{1},\diamond)$, then $\dual{w} = w$. 
		\item If $w= (\star,\ora{p}_1,\star)$, then $\dual{w} = (\star,\ora{1},\ola{p}_1,\star)$, and vice versa.
		\item If $w=(\star,\ola{q}_1,\star)$, then $\dual{w} = 
		(\star,\ola{1}, \ora{q}_1,\star)$, and vice versa.
	\end{itemize}
	\item  Assume that $w$ is none of the previous cases.  
	\begin{itemize}
		\item If $w$ begins with $\star$ or $\wh{\star}$, we change the beginning of $w$ as follows.
		\begin{align*}
			\begin{array}{|r|rl|r|r|r|}
				\hline
				w & (\star,\ora{p}_1,\ldots)& p_1 > 1 
				& (\star,\ora{1},\ola{q}_1,\ldots) 
				& (\wh{\star},\ora{p}_1,\ldots) & (\star, \ola{q}_1, \ldots) 
				\\
				\hline
				\dual{w}&
				(\star,\ora{1},{\ola{p}}_1^-\ldots) & & (\star, {\ora{q}}_1^+, \ldots) 
				& (\star,\ola{p}_1^-,\ldots)&
				(\wh{\star},\ora{q}^+_1,\ldots)
				\\
				\hline
			\end{array}
		\end{align*}
		\item If $w$ ends with $\star$ or $\wh{\star}$, the end is changed similarly as follows.
		\begin{align*}
			\begin{array}{|r|ll|l|l|l|l|}
				\hline
				w &
				(\ldots \ola{q}_{k},\star)& q_k > 1 
				& (\ldots \ora{p}_{k},\ola{1},\star)& (\ldots  \ola{q}_k,\wh{\star})
				& (\ldots  \ora{p}_k,\star)
				\\
				\hline
				\dual{w} 
				& (\ldots \ora{q}_{k}^-,\ola{1},\star)& & (\ldots \ola{p}^+_{k},\star) 
				&
				(\ldots  \ora{q}^-_k,\star)
				&
				(\ldots  \ola{p}^+_k,\wh{\star})
				\\ \hline
			\end{array}
		\end{align*}
	\end{itemize}
	At last, any oriented number $x_i$ of the word $w$ which is not affected by the changes above, has to be replaced by $x_i^t$, that is, we flip the orientation of its arrow.
\end{itemize}
\end{dfn}
In particular, $\begin{td} w \ar[mapsto]{r} \& \dual{w} \end{td}$ defines an involution on the set of finite and periodic words of $\Acat$.
We recall that for a sign $\varepsilon = \pm$, we denote its opposite sign by $\ol{\varepsilon} = \mp$.
\begin{thm}\label{thm:dualities}
Let $\omega$ be a band or string of $\Acat$, and $M(\omega)$ the corresponding indecomposable object.
Then there are isomorphisms $\inv(M(\omega)) \cong M(\inv(\omega))$, $\dual{(M(\omega))} \cong M(\dual{\omega})$ and $\cdual{(M(\omega))} \cong M(\cdual{\omega})$ 
with $\inv(\omega)$, $\dual{\omega}$ and $\cdual{\omega}$ defined as follows.
\begin{itemize}
	\item If $\omega$ is a simple special string, that is, if $M(\omega) \cong S_{\varepsilon}$, 
	or if $\omega$ is equivalent to $((\star,\ora{1},\diamond),\varepsilon)$,
	then $
	\inv(\omega) = (w,\ol{\varepsilon})$,
	$\dual{\omega} = (\dual{w}, \varepsilon)$,
	and $\cdual{\omega} = (\dual{w}, \ol{\varepsilon})$.
	\item
	Otherwise 
	$\inv(\omega)$, $\dual{\omega}$ and $\cdual{\omega}$
	are given by the following table.
	\begin{align*}
		\begin{array}{|c||c|c|c|c|}
			\hline
			\omega & w & (w,\varepsilon) &  (\varepsilon_1, w,\varepsilon_2) & (w,m,\lambda) \\
			\hline \hline
			\inv(\omega) &
			w &  (w,\ol{\varepsilon}) &
			(\ol{\varepsilon}_1, w, \ol{\varepsilon}_2)
			&
			(w,m,\lambda) \\ \hline
			\dual{\omega} &
			\dual{w}
			&
			(\dual{w},\ol{\varepsilon})
			&
			(\ol{\varepsilon}_2,\dual{w},\ol{\varepsilon}_1)
			&
			(\dual{w},m,\lambda^{\mp}) \\ \hline
			\cdual{\omega} & \dual{w}
			&
			(\dual{w},{\varepsilon})
			&
			({\varepsilon}_2,\dual{w},{\varepsilon}_1)
			&
			(\dual{w},m,\lambda^{\mp}) \\ \hline
		\end{array}
	\end{align*}
\end{itemize}	
\end{thm}
\begin{proof}
A band or string $\omega$ of $\Acat$ has one of the forms from Table \ref{tab:rep-str-bands}. We denote by $\breve{\omega}$ its expanded notation in the sense of Remark~\ref{rmk:trans1} respectively~\ref{rmk:trans2}.
Proposition~\ref{prp:inv} 
implies that $\inv(\CPr(\breve{\omega})) \cong \CPr(\inv(\breve{\omega}))$ and 
$\inv(\breve{\omega}) = (\inv(\omega))'$.
Thus, $\inv(M(\omega)) \cong M(\inv(\omega))$.
Similarly, it can be verified that 
$\dual{(\breve{\omega})}$ is equivalent to the translation of $\dual{\omega}$ and the isomorphism 
$\dual{(M(\omega))} \cong M(\dual{\omega})$
follows from Theorem~\ref{thm:duals}.
The third isomorphism follows from the previous two using
that $\cdual{(-)} = \inv \circ \dual{(-)}$.
\end{proof}

\begin{cor}\label{cor:self-dual}
Let $M$ be an indecomposable object in $\Acat$
and $\omega$ the band or string of $\Acat$ such that $M \cong M(\omega)$.
Then the following statements hold.
\begin{enumerate}
	\item 
	It holds that	$\cdual{M} \cong M$ if and only if 
	$\omega$ is equivalent to one of the following:
	\begin{enumerate}
		\item a usual string $w = (\wh{\star},\ora{p}_1,v,\ora{p}_1,\star)$
		or $w = (\star,\ora{p}_1^+,
		v, \ora{p}_1, \ola{1}, \star)$
		with $p_1^+ =p_1 + 1$,
		\item a bispecial string $(\varepsilon_1, w,\varepsilon_1)$ with $w = (\ora{p}_1, v, \ora{p}_1)$,
		\item a band $(w,m,\lambda)$ 
		such that 
			$\dot{w} = (\ora{p}_{1}, v)$
			or
			$\dot{w}=
			(x_1,x_2,\ldots  x_k, x_1^t, x_2^t, \ldots x_k^t)$ with an odd number $k \geq 3$ 
			and $\lambda  \in \{ -1,1\}$ if $w$ is asymmetric,
	\end{enumerate}
	where
	$v= (x_2, \ldots x_{k}, x_{k+1}, x_{k} \ldots x_2)$ with a left-oriented number $x_2$ and
	$k \geq 1$.
	\item It holds that $\dual{M} \cong M $
	if and only if $\omega$ is equivalent to one of the following:
	\begin{enumerate}
		\item 
		a usual string or a band from the previous list, 
		\item
		the special string $(w, \varepsilon)$ with $w = (\star,\ora{1},\diamond)$,
		or, equivalently, $M \cong S_{+}$ or $S_-$,
		\item  a bispecial string $(\varepsilon_1,w, \ol{\varepsilon}_1)$ with  $w = (\ora{p}_1, v, \ora{p}_1)$ and $v$ as before.
	\end{enumerate}
\end{enumerate}
\end{cor}
\begin{proof}
The `if'-implication in both claims follow directly from Theorem~\ref{thm:dualities} and the equivalence conditions.
Vice versa, assume that $\cdual{M} \cong M$ or $\dual{M} \cong M$.
\begin{itemize}
	\item 
	If $\omega$ is a string with underlying word $w$, then $\dual{w}$ must have the same length as $w$. Together with Theorem~\ref{thm:dualities} and the equivalence conditions, this observation leads to the classes of strings in question.
	\item If $\omega$ is a band $(w,m,\lambda)$, then it is equivalent to $(\dual{w},m,\lambda^{\mp})$ by Theorem~\ref{thm:dualities}.
	In particular, there exists $j \in \Z$ such that
	$(w_{2j+1},m,\lambda^{\mp}) =
	(w^{\op},m,\lambda^{\mp})$ or $(w,m,\lambda)$. 
	\begin{itemize}
		\item
		In the first case, it follows that $x_{j-i} = x_{j+i}$ for any $i \in \Z$, which leads to the bands of the first type.
		\item In the second case, since $\dot{w}$ has period $2k$ for certain $k \in \N$, it follows that $w_{2j+1}   = w_k$ and $\lambda^{\mp} = \lambda$. Thus, $\omega$ is a band of the second type of the claim.
	\end{itemize}
\end{itemize}
This concludes the proof of both `only if'-implications.
\end{proof}
\begin{ex}\label{ex:dual-bd}
Any band $\omega = (w,m,\lambda)$ with an underlying word of period two yields a representation
$M = M(\omega)$ satisfying $\cdual{M} \cong \dual{M} \cong M$.
In case of the smallest possible band 
$\omega=(w,1,\lambda)$ 
with $\dot{w} = (\ora{1},\ola{1})$ and $\lambda \in \kk \backslash\{-1,0\}$
the band representation $M(\omega)$ is isomorphic to the representation $M_{\tilde{\lambda}}$ in \eqref{eq:sh-band}
with $\tilde{\lambda} = \frac{\lambda}{\lambda+1} \notin \{0,1\}$.
For this representation, an explicit isomorphism $\begin{td}M_{\tilde{\lambda}} \ar{r}{\sim} \& \dual{(M_{\tilde{\lambda}})} \end{td}$
is given by the matrices
$(\phi_-,\phi_\star,\phi_+)=
\left(\left(	\begin{smallmatrix}
	0 & 1 \vphantom{\tilde{\lambda}} \\
	1 & 0 \vphantom{\tilde{\lambda}}
\end{smallmatrix}\right),
\left(
\begin{smallmatrix}
	\tilde{\lambda} & \tilde{\lambda} \\
	\tilde{\lambda} & 1 
\end{smallmatrix}
\right),
\left(\begin{smallmatrix}
	0 & \tilde{\lambda} \\
	\tilde{\lambda} & 0 
\end{smallmatrix}\right) \right)$. 
\end{ex}

\subsection{Homological invariants of indecomposable representations}
\label{subsec:inv}
For any object $M$ of $\Acat$,
we denote by 
$\head M$ its top,
and by $\soc M$  its socle 
as in Subsection~\ref{subsec:cyc-mod-res}.
We recall that $\uldim M= (\dim M_v)_{v \in \{\star,+,-\}}$ 
records the Jordan--Hölder multiplicities of $M$ (Remark~\ref{rmk:JH-mult}).

\begin{thm}\label{thm:inv}
Let $M$ be an indecomposable object in $\Acat$.
Let  $\omega$ denote the band respectively the string of $\Acat$ such that $M \cong M(\omega)$. 
We use the following notation, suppressing the dependence from $\omega$.
\begin{enumerate}
	\item
	For any integer $i \in \Z$
	and any vertex $v \in Q_0$
	we denote by $\ct{i}{v} \in \N_0$ the multiplicity of the indecomposable projective $P_v$ in $\X_i(\omega)$.
	\item We denote by $\ctpr{1}{\star} \in \{0,1,2\}$ the sum of the number of beginnings of the form $(\star,\ora{p}_1)$ with $p_1 > 1$ and the number of endings $(\ola{q}_{\ell},\star)$ with $q_\ell > 1$ in the word $w$ underlying $\omega$.
	\item Denoting $\dot{w} = (\olra{n}_1, \olra{n}_2, \ldots  \olra{n}_{\ell})$, we set	$n = 
	\sum_{i=1}^{\ell} \mu n_i$ 
	with
	$\mu = m$
	if $\omega$ is a band $(w,m,\lambda)$ respectively
	$\mu = 1$ if $\omega$ is a string with underlying word $w$.  
\end{enumerate}
In this notation, the following equalities hold.
	\begin{align}
		\label{eq:hom-dim}
		\prdim M = 
		\begin{cases} 2 & \text{if }\ct{2}{\star} \neq 0, \\
			1 & \text{otherwise},
		\end{cases}
		&&
		\injdim M = 
		\begin{cases} 2 & \text{if }\ct{0}{\star} \neq 0, \\
			1 & \text{otherwise},
		\end{cases}
	\end{align}
	\begin{align*}
		\setlength{\arraycolsep}{2pt}	
		\begin{array}{rclccrclccrcl}
			\head M &\cong& S_-^{t_-} \oplus S_\star^{t_\star} \oplus S_+^{t_+}
			&&& \text{ with } 
			t_\pm &=& \ct{0}{\pm}, &&&     
			t_\star &=& \ct{0}{\star},
			\\
			\soc M &\cong& S_-^{u_-} \oplus S_\star^{u_\star} \oplus S_+^{u_+}
			&&& \text{ 
				with } u_\pm &=&   \ct{1}{\mp} + \ctpr{1}{\star} -\ct{2}{\star},  &&& 
			u_{\star} &=& \ct{2}{\star},
			\\ 
			\uldim M &=&  \begin{pmatrix} 
				n_-, n_\star, n_+
			\end{pmatrix} &&&
			\text{with }
			n_\pm &=& n 
			+  \ct{1}{\mp} - \ct{0}{\mp} -\ct{2}{\star}, &&&
			n_\star &=& n - \ct{1}{\star} +  \ct{0}{\star}.
		\end{array}
	\end{align*}
\end{thm}
\begin{proof}
	We recall that for $M$ as above, there  exists a band or string $\omega$ such that $\CP(\omega)$ is a minimal projective resolution of $M$ by Theorem~\ref{thm:prbij}.
	\begin{enumerate}
		\item The non-zero terms of the complex $\CP(\omega)$ 
		are already determined by its gluing diagram in 
		\eqref{eq:per-ab2} respectively \eqref{eq:str-ab2}.
		Since any gluing diagram of a projective resolution of length two may have only vertices of type $\star$ at degree $2$, the statement on $\prdim M$ follows.
		\item Since $\head M \cong \head \X_0(\omega)$ the formula for the top follows by counting indecomposable projective at degree zero in the gluing diagram 
		\eqref{eq:per-ab2} respectively \eqref{eq:str-ab2}.
	\end{enumerate}
	Note that the passage from the gluing diagram of $\omega$ to that of $\dual{\omega}$ admits the following interpretation in terms of gluing edges at degree one.
	\begin{itemize}
		\item Any gluing edge at degree one between arrows with support $\{0,1\}$ 
		gives rise to a gluing edge at degree zero for $\dual{\omega}$.
		\item Any gluing edge between an arrow with support $\{0,1\}$ and an arrow with support $\{1,2\}$ in the gluing diagram $\omega$
		has no corresponding gluing edge at degree zero in the gluing diagram of $\dual{\omega}$.
		\item Any proper subdiagram	of the form $\begin{gd} {P_\star} \ar{r}{t^n} \& P_{\varepsilon} \end{gd}$ with $n >1$
	in the gluing diagram of $w$ gives rise to a gluing edge at degree zero in $\dual{\omega}$.
\end{itemize}
Therefore, 
it follows that
	$\ct{0}{\star}(\dual{\omega}) = \ct{2}{\star}(\omega)$ and 
	$\ct{0}{\pm}(\dual{\omega}) = \ct{1}{\mp}(\omega)
	-\ct{2}{\star}(\omega)+\ctpr{1}{\star}(\omega)$.
	\begin{enumerate}
		\setcounter{enumi}{2}
		\item The claim on the injective dimension of $M$ follows 
		from
		$\injdim M = \prdim \dual{M} $, Theorem~\ref{thm:dualities}
		and the fact that
		$\ct{2}{\star}(\dual{\omega}) \neq 0$ if and only if $ \ct{0}{\star}(\omega) \neq 0$.
		\item 
	Using the previous equalities, the description of $\dual{(-)}$ in Theorem~\ref{thm:dualities}
	and the isomorphisms 
	$\soc M(\omega) \cong \head \dual{(M(\omega))} \cong \head M(\dual{\omega})$
	allow to derive
	the claims for the socle-multiplicities from the formulas for the top-multiplicities.
	\item We give two proofs to establish the claim on $\uldim M(\omega)$.
	\begin{itemize}
		\item 
		Let
		$\omega_1$, $\omega_2,$ $\ldots$ $\omega_{\ell}$ denote the strings with underlying words of length one such that  $\CPr^{cyc}(\omega) = \bigoplus_{j=1}^{\ell} \CPr(\omega_j)$.
		By Proposition~\ref{prp:cyc}
		it holds that $\uldim M(\omega) = \uldim M^{cyc}(\omega) = \sum_{j=1}^{\ell} \uldim M(\omega_j)$.
		Moreover, $\sum_{j=1}^{\ell} \ct{i}{v}(\omega_j)
		= \ct{i}{v}(\omega)$ for any $i \in \{0,1,2\}$ and any $v \in Q_0$. 
		Therefore, it is sufficient to verify the 
		formula for $\uldim M(\omega)$ assuming that $M(\omega)$ is cyclic.
		This is straightforward using the list of cyclic representations from Table~\ref{tab:cyc-base},
		which is complete by Proposition~\ref{prp:fd-cyc} or Corollary~\ref{cor:cyclic}.
		\item Alternatively, there is the following, conceptual argument, which shows how the formula for $\uldim M(\omega)$ was obtained in the first place.
		Set $\CY = \B \otimes_{\A} \CPr$ and $\CV = \AI \otimes_{\A} \CPr$.
		According to~\eqref{eq:dim-formula} it holds that
		$\uldim M(\omega) =  d(\CV) + d(\CY)$
		with 
		$$d(\CV) \colonequals \uldim \inv(\V_1) 
		- \uldim \inv(\V_0) =  \begin{pmatrix} \ct{1}{-} - \ct{0}{-}, 0, \ct{1}{+} - \ct{0}{+}  \end{pmatrix}
		$$
		and
		$d(\CY) \colonequals
		\uldim \mathstrut_{\A} H_0(\CY) 
		- \uldim \mathstrut_{\A} H_{1}(\CY)$, which can be determined as follows.
		The complex $\CY$ can be determined using the gluing diagram of $\omega$ together with Lemma~\ref{lem:norm-complex}.
		Since $\CPr(\omega)$ is a projective resolution,
		Theorem~\ref{thm:proj-res}
		restricts the possibilities
		for each indecomposable summand $\CY'$ of $\CY$
		to the following cases.
		\begin{align*}
					\begin{array}{|c||c|c|c|c|}
						\hline
						\text{subword}
						& \wh{\star} 
						& (\star,\ora{n}_i) \text{ or } (\ola{n}_i,\star) 
						& (\ora{n}_i,\star) \text{ or } (\star,\ola{n}_i) 
						& \text{other }(\olra{n}_i)
						\\ \hline
						\CY'
						& \begin{gd} \underset{2}{\BP_{\star}} \ar{r}{\cdot t} \& \underset{1}{\BP_{\diamond}} \end{gd}
						& \begin{gd} \underset{1}{\BP_{\star}} \ar{r}{\cdot t^{n_i}} \& \underset{0}{\BP_{\diamond}} \end{gd} 
						& \begin{gd} \underset{1}{\BP_{\diamond}} \ar{r}{\cdot t^{n_i}} \& \underset{0}{\BP_{\star}} \end{gd} 
						& \begin{gd} \underset{1}{\BP_{v}} \ar{r}{\cdot t^{n_i}} \& \underset{0}{\BP_{v}} \end{gd} 
						\\ \hline
							d(\CY')
						& \begin{pmatrix} -1 , 0 , -1 \end{pmatrix} 
						& \begin{pmatrix} n_i, n_i-1, n_i \end{pmatrix} 
						& \begin{pmatrix} n_i, n_i+1, n_i \end{pmatrix} 
						& \begin{pmatrix} n_i, n_i, n_i \end{pmatrix} 
						\\
						\hline
					\end{array}
				\end{align*} 
				The restriction from $\B$ to $\A$ preserves the dimension vector of any finite-dimensional $\B$-module.
				The multiplicities of the first three summands
				in the direct sum decomposition of $\CY$ are determined by $\ct{2}{\star}$, $\ct{1}{\star}$ and $\ct{0}{\star}$, respectively.
				It follows that
				\begin{align*}
					d(\CY) 
					&=  
					\begin{pmatrix} n - \ct{2}{\star},  n - \ct{1}{\star}+
						\ct{0}{\star}, n -\ct{2}{\star}  \end{pmatrix}
				\end{align*}
				which has to be checked separately
				for the case that $\CY$ is given by the last type of summand with $u = v =\star$.
				Using this expression for
				the summand $d(\CY)$ of $\uldim M(\omega)$
				yields the claimed formula.
				\qedhere
			\end{itemize}
		\end{enumerate}
	\end{proof}
	We note that \eqref{eq:hom-dim} follows also from Proposition~\ref{P:BoundsDimensions}.
	The invariants of the last theorem admit more explicit descriptions.
	\begin{cor} \label{cor:inv} The following statements hold.
		\begin{enumerate}
			\item \label{inv-band}
			Let $\omega$ be a band $(w,m,\lambda)$ of $\Acat$.
			Denoting $\ell$ and $n$ as in  Theorem~\ref{thm:inv},
			the band representation $M = M(\omega)$
			has the  invariants
			\begin{align*}
				\prdim M = \injdim M = 1,&&
				\head M \cong \soc M \cong (S_- \oplus S_+)^{\frac{m\ell}{2}}, 
				&&
				\uldim M =  \begin{pmatrix}
					m\,n, m\,n, m\,n
				\end{pmatrix}.
			\end{align*} 
			\item \label{inv-str} For any string $\omega$ of $\Acat$ up to involution, the top-, the socle- and the Jordan-Hölder-multiplicities of the string representation $M(\omega)$ are described in 
			Table~\ref{tab:inv-str},
			where $n$ denotes the sum of all numbers in the word underlying $\omega$.
		\end{enumerate}
	\end{cor}
	\begin{proof}
		The first claim follows directly from Theorem~\ref{thm:inv}.
		To show the second claim, let $\omega$ be a string.
		For any integer $i \in \{0,1,2\}$, 
		we denote by
		\begin{itemize}
			\item  
			$\alpha_{iv}$
			the number of indecomposable projectives $P_v$ in $\X_i(\omega)$
			which correspond to \emph{ends} of the gluing diagram of $\omega$ for a fixed vertex $v \in Q_0$,
			\item $\ctdeg{i} = \sum_{v \in Q_0} \alpha_{iv}$ the total number of ends at degree $i$,
			\item $\gamma_i$ 
			the number of gluing edges at degree $i$ in the gluing diagram of $\omega$.
		\end{itemize}
		We note that $\beta_{i \star} = \alpha_{i\star}$ 
		and
		$\beta_{i\pm} = \alpha_{i\pm} + \gamma_i$ 
		for any $i \in \{0,1,2\}$.
		A case-by-case analysis of the possible gluing diagrams allows to verify that the numbers of gluing edges are given by the formulas
		\begin{align} \label{eq:edge-count}
			\gamma_1 = \frac{1}{2} \left( \ell + \ctdeg{2} -
			\ctdeg{1} \right) && 
			\gamma_0 = \frac{1}{2} \left(\ell - \ctdeg{0} \right) 
		\end{align}
		where $\ell$ denotes the length of the word $\dot{w}$
		in case $\omega$ is a band respectively 
		the length of $w$ 
		if $\omega$ is a string.
		Combining 
		Theorem~\ref{thm:inv} with these observations
		facilitates the computations of the multiplicities in Table~\ref{tab:inv-str} for the string representation $M(\omega)$.
	\end{proof}

			\begin{longtable}{|l||ccc|ccc|ccc|}
				\caption{Top, socle and dimension vector of string representations}
				\label{tab:inv-str} \\
				\hline			\multicolumn{1}{|l||}{\text{Usual strings }$w$ } & \multicolumn{3}{c}{$\uldim\head M$} & \multicolumn{3}{|c|}{$\uldim\soc M$} & \multicolumn{3}{c|}{$\uldim M$} 
				\\ \hline
				$(\star,
				{\color{blue} \ora{p}_1,\ola{q}_1, \ldots \ora{p}_{k},}
				\star) 
				$&$
				k-1 $&$ 1 $&$ k-1
				$&$
				k' $&$ 0 $&$ k'
				$&$ n $&$ n $&$ n$
				\\ \hline
				$(\wh{\star},
				{\color{blue} \ora{p}_1,\ola{q}_1, \ldots \ora{p}_{k}},
				\star)
				$&$
				k-1 $&$ 1 $&$ k-1
				$&$
				k-1 $&$ 1 $&$ k-1
				$&$
				n $&$ n+1 $&$ n$
				\\ \hline
				$			(\wh{\star},
				{\color{blue} \ora{p}_1,\ola{q}_1, \ldots  \ora{p}_{k}, \ola{q}_{k},}
				\star)
				$&$
				k $&$ 0 $&$ k $&$ 	
				k'' $&$ 1 $&$ k''
				$&$ n-1 $&$ n-1 $&$ n-1$
				\\ \hline
				$					(\star,
				{\color{blue} \ola{q}_1,\ora{p}_1, \ldots \ola{q}_{k}, \ora{p}_{k},}
				\star)
				$&$
				k-1 $&$ 2 $&$ k-1
				$&$
				k $&$ 0 $&$ k
				$&$
				n+1 $&$ n+2 $&$ n+1$
				\\ \hline
				$					(\star,
				{\color{blue} \ora{p}_1, \ola{q}_1, \ldots \ora{p}_k, \ola{q}_{k},}
				\star)
				$&$ 
				k $&$ 0 $&$ k
				$&$
				k''' $&$ 0 $&$ k'''
				$&$
				n-1 $&$ n-2 $&$ n-1$
				\\ \hline
				$					(\wh{\star},
				{\color{blue} \ora{p}_1,\ola{q}_1, \ldots \ora{p}_k, \ola{q}_{k},}
				\wh{\star}) 
				$&$
				k $&$ 0 $&$ k
				$&$ 
				k-1 $&$ 2 $&$ k-1
				$&$
				n-1 $&$ n $&$ n-1$
				\\ \hline	\multicolumn{10}{|l|}{{
						\text{where $k' = k - \delta_{1,p_1}$, \quad $k'' = k - \delta_{1,q_{k}}$ \quad and \quad $k''' = k + 1 - \delta_{1,p_1}- \delta_{1,q_{k}}$. 
				}}}\\
				\hline
			\end{longtable}
			\begin{longtable}{|l||ccc|ccc|ccc|} 
				\hline			
				\multicolumn{1}{|l||}{\text{Special strings } $(w,\varepsilon)$ }
				& \multicolumn{3}{c}{$\uldim\head M$} & \multicolumn{3}{|c|}{$\uldim\soc M$} & \multicolumn{3}{c|}{$\uldim M$} 
				\\ \hline
				$					((\star,
				{\color{blue} \ola{q}_1,\ora{p}_1 \ldots \ola{q}_{k}}
				,\diamond),-) 
				$&$ k-1 $&$ 1 $&$ k-1 $&$ k-1 $&$ 0 $&$ k $&$ n $&$ n+1 $&$ n+1 $
				\\ \hline 
				$					((\star,
				{\color{blue} \ora{p}_1,\ola{q}_1 \ldots \ora{p}_k}
				,\diamond),-) 
				$ & $ k $&$ 0 $&$ k-1 $&$ k' $&$ 0 $&$ k' $&$ n $&$ n-1 $&$ n-1 $\\ \hline
				$					((\wh{\star},
				{\color{blue} \ora{p}_1,\ola{q}_1 \ldots \ora{p}_k},
				\diamond),-) 
				$&$ k $&$ 0 $&$ k-1 $&$ k-1 $&$ 1 $&$ k-1 $&$ n $&$ n $&$ n-1$ \\ \hline
				$					((\star,
				{\color{blue} \ola{q}_1,\ora{p}_1 \ldots \ola{q}_k, \ora{p}_k}
				,\diamond),-) 
				$&$ k $&$ 1 $&$ k-1 $&$ k $&$ 0 $&$ k $&$ n+1 $&$ n+1 $&$ n $\\ \hline
				$					((\star,
				{\color{blue} \ora{p}_1,\ola{q}_1 \ldots \ora{p}_k, \ola{q}_k}
				,\diamond),-)
				$&$ k $&$ 0 $&$ k $&$ k' $&$ 0 $&$ k'+1 $&$ n-1 $&$ n-1 $&$ n$ \\ \hline
				$					((\wh{\star},
				{\color{blue} \ora{p}_1,\ola{q}_1 \ldots \ora{p}_k, \ola{q}_k}
				,\diamond),-)
				$&$ k $&$ 0 $&$ k $&$ k-1 $&$ 1 $&$ k $&$ n-1 $&$ n $&$ n $\\
				\hline 
				\multicolumn{10}{|l|}{{
						\text{where $k' = k - \delta_{1,p_1}$.
				}}}\\
				\hline
			\end{longtable}
			\begin{longtable}{|l||ccc|ccc|ccc|}
				\hline
				\multicolumn{1}{|l||}{\text{Bispecial strings } $(\varepsilon_1,w,\varepsilon_2)$ }
				& \multicolumn{3}{c}{$\uldim\head M$} & \multicolumn{3}{|c|}{$\uldim\soc M$} & \multicolumn{3}{c|}{$\uldim M$} 
				\\ \hline
				$					(-,(\diamond,
				{\color{blue} \ora{p}_1, \ola{q}_1, \ldots  \ora{p}_k}
				,\diamond),-) $&$ 
				k $&$ 0 $&$ k-1 $&$ k-1 $&$ 0 $&$ k $&$ n $&$ n $&$ n $\\ \hline
				$					(-,(\diamond,
				{\color{blue} \ora{p}_1, \ola{q}_1, \ldots  \ora{p}_k}
				,\diamond),+)		$&$	 
				k-1 $&$ 0 $&$ k $&$ k-1 $&$ 0 $&$ k $&$ n+1 $&$ n $&$ n-1 $\\ \hline
				$					(-,(\diamond,
				{\color{blue} \ola{q}_1,\ora{p}_1, \ldots \ola{q}_k, \ora{p}_{k}}
				,\diamond),-)^1 $&$ 
				k+1 $&$ 0 $&$ k-1 $&$ k $&$ 0 $&$ k $&$ n-1 $&$ n $&$ n+1$ \\ \hline
				$					(-,(\diamond,
				{\color{blue} \ora{p}_1,\ola{q}_1, \ldots \ora{p}_k, \ola{q}_k}
				,\diamond),-)^1 
				$&$
				k $&$ 0 $&$ k $&$ k-1 $&$ 0 $&$ k+1 $&$ n+1 $&$ n $&$ n-1$ \\ \hline
				$			
				(\varepsilon_1,(\diamond,
				{\color{blue} x_1,y_1, \ldots x_k, y_k}
				,\diamond),\varepsilon_2)^2	 $&$ k $&$ 0 $&$ k $&$ k $&$ 0 $&$ k $&$ n $&$ n $&$ n $\\ \hline 
				\multicolumn{10}{|l|}{{
						\text{
							$\mathstrut^1$ assuming that
							the primitive root of $w$ has odd multiplicity.
				}}}
				\\
				\multicolumn{10}{|l|}{{
						\text{
							$\mathstrut^2$ assuming that $\varepsilon_1 \neq \varepsilon_2$ if  the primitive root of $w$ has odd multiplicity.
					}}
				}\\
				\hline
			\end{longtable}

							\addtocounter{equation}{-2}
							\begin{rmk} 
								We refer to the invariants in Theorem~\ref{thm:inv} as `numerical invariants' below and 
								note the following.
								\begin{enumerate}
									\item The statements on the projective and injective dimension recover Proposition~\ref{P:BoundsDimensions}.
									\item The number of generators of $M(\omega)$ as $\A$-module is given by $\ell \mu$.
									\item The modules $M(\omega)$ and $M^{cyc}(\omega)$ share the same numerical invariants.
									\item 
									Reordering the numbers without change of orientations in the word underlying a band or string $\omega$ does not change the numerical invariants of $M(\omega)$.
									\item  To any band $\omega = (w,m,\lambda)$ with $\dot{w} = (x_i)_{i=1}^{\ell}$ we may associate a bispecial string $\omega' = (-,w',+)$ setting $w' = (x_i)_{i=1}^{m\ell}$.
									Then the band representation $M(\omega)$ and the bispecial string representation $M(\omega')$
									have the same numerical invariants.
									It appears that $M(\omega)$ and $M(\omega')$ 
									can be distinguished only via the involution $\inv$; the band representation $M(\omega)$ is $\inv$-invariant, while  the bispecial string representation $M(\omega')$ is not.	
								\end{enumerate}
							\end{rmk}
							\begin{ex}
								In case of the band $\omega$ from Example~\ref{ex:dual-bd}, 
								the band representation $M(\omega)$ was given by $M_{\tilde{\lambda}}$ in \eqref{eq:sh-band}
								with $\tilde{\lambda} \neq 0,1$ while the associated bispecial string representation $M(\omega')$ is isomorphic to the degeneration $M_0$.
							\end{ex}

							\begin{lem}\label{lem:chi}
								The Euler characteristic of a pair $(M,N)$ of objects from $\Acat$
								is given by
								\begin{align*}
									\chi(M,N) =  
									\sum_{i=0}^2 (-1)^i \, \dim_{\kk} \Ext^i_{\A}(M,N)
									= x_+ y_+ + x_- y_- 
								\end{align*}
								where $x_{\pm} = \dim M_\pm - \dim M_{\star}$
								and $y_{\pm} = \dim N_\pm - \dim N_{\star}$. 
							\end{lem}
							\begin{proof}
								We recall that $\chi(M,N)$ 
								is initially defined by $\sum_{i \in \Z}
								(-1)^i \, \dim_{\kk} \Ext^i_{\A}(M,N)$.
								As $\A$ has global dimension two, the first equality follows.
								Since the Euler form is additive on short exact sequences, it holds also that 
								\begin{align*}
									\chi(M,N) = 
									\uldim M \
									\Phi_{\A} \ (\uldim N)^T
									\quad \text{ with }
									\Phi_{\A} = (\chi(S_i,S_j))_{i,j \in Q_0} = 
									\begin{pNiceArray}{ccc}[first-row, last-col]
										- & \star & + \\
										\phantom{-}1 & -1 & \phantom{-}0 & -\\
										-1 & \phantom{-}2 & -1 & \star \\
										\phantom{-}0 & -1 & \phantom{-}1 & +
									\end{pNiceArray}.
								\end{align*}
								This yields
								the remaining equality.
							\end{proof}

							\begin{cor}\label{cor:chi}
								For any indecomposable objects $M$,$N$ in $\Acat$
								the following statements hold.
								\begin{enumerate}
									\item
									It holds that
									$
									\lvert \chi(M,N)\rvert \leq 2$.
									\item Let $\omega$ be the band or string of $\Acat$ such that $M \cong M(\omega)$.
									Then								the following statements 		on 
									the Euler characteristic $\chi(M) \colonequals \chi(M,M)$		
									are true.									\begin{itemize}
										\item $\chi(M) =2$ if and only if 
										$\omega$ is one of the following;
										\begin{itemize}
											\item $\omega$ is a usual string $w$ with $(\ct{2}{\star},\ct{1}{\star},\ct{0}{\star})\in\{(1,1,0),(0,1,1)\}$,
											\item $\omega$ is a bispecial string $(\varepsilon_1, w,\varepsilon_2)$ 
											such that the length of $w$  is odd and $\varepsilon_1 \neq \varepsilon_2$,
											\item $\omega$ is a bispecial string $(\varepsilon_1, w,\varepsilon_2)$
											such that the length of $w$ is even, $\varepsilon_1 = \varepsilon_2$ and the primitive root of $w$ has odd multiplicity.
										\end{itemize} 
										\item
										$\chi(M) = 1$ if and only if $\omega$ is a special string.
										\item $\chi(M)=0$ if and only if $\omega$ is a band or another string than the ones above.
									\end{itemize}
								\end{enumerate}
							\end{cor}
							\begin{proof}
								Set $\ctalt{v} = \sum_{i=0}^2 (-1)^i \alpha_{iv}$
								for each $v \in Q_0$.
								In the notations of Lemma~\ref{lem:chi}
								and Theorem~\ref{thm:inv}, using \eqref{eq:edge-count},
								it follows that 
								$$x_{\pm} =  n_{\pm} - n_{\star} = \ct{1}{\mp} - \ct{0}{\mp} 
								- \ctalt{\star}
								= 
								\frac{1}{2}( -\ctalt{\star} \pm (\ctalt{+} - \ctalt{-}))   \in \{-1,0,1\}.
								$$ 
								Similarly, $y_{\pm} \in \{-1,0,1\}$. 
								The first claim follows from Lemma~\ref{lem:chi}. 
								The second claim can be verified using that $\chi(M) = x_+^2 + x_-^2$.
							\end{proof}

\section{Explicit band and string representations}
\label{sec:bases}
We recall that a nilpotent representation of the Gelfand quiver 
is given by a diagram of finite-dimensional $\kk$-linear spaces and $\kk$-linear maps
\begin{align}\label{eq:gqrep}
	\begin{gqrep}{M_+}{M_\star}{M_-}{B_+}{B_-}{A_+}{A_-}\end{gqrep}
	\quad
	\text{such that }B_- A_- = B_+ A_+ \text{ and }B_- A_- \text{ is nilpotent.}
\end{align}
In this section,
we give an explicit description of the indecomposable quiver representation
$M(\omega)$ associated to each band or string $\omega$ of $\Acat$.
This section requires notions from the previous one.
Except for an argument providing the dimension of $M(\omega)$, the proofs of this section are self-contained. 
Some of the examples include gluing diagrams with induced arrows. 
The precise rules for induction of arrows have been described in Subsection~\ref{subsec:ind-arr}, though reviewing these technicalities is not necessary when taking the band and string resolutions from Subsections~\ref{sec:proj-res1} and~\ref{sec:proj-res2} as given.
The explicit representation matrices of bands are described in Propositions~\ref{prp:sh-band}, \ref{prp:bands}, those for strings in Proposition~\ref{prp:str-base}.

\subsection{Bases of cyclic modules}
Let $\beta \in Q_0$. Let $g$ denote the idempotent generator of the indecomposable projective $\A$-module $P_{\beta}$.
As observed in Subsection~\ref{subsec:ind-proj},
a topological basis of the infinite-dimensional $\kk$-linear space $P_{\beta}$ 
is given by the monomials $t^n_v g$ with $(v,n)$
running through the index set
\begin{align*}
	I_{\beta} \colonequals
	(Q_0 \times \N) \cup 
	\begin{cases}
		\{(v,0) \mid v \in Q_0\} & \text{if } \beta = \star,\\
		\{(\beta,0)\}  & \text{otherwise}.
	\end{cases}
\end{align*}
The action of the four matrices on the basis vectors of the indecomposable projective $\A$-module $P_{\beta}$ is completely determined by the relations
\begin{align}
	\label{eq:std-rel}
	B_- t_-^p g &= t_\star^{p+1} g 
	&&
	(p \geq \delta_{\beta-} )
	&&&
	B_+ t_+^p g &= t_\star^{p+1} g&&
	(p \geq \delta_{\beta+})
	\\
	A_- t_\star^p g &= t_-^{p} g && (p\geq \delta_{\beta\star})
	&&&
	A_+ t_\star^p g &= t_+^p g
	&&
	(p \geq \delta_{\beta\star})
\end{align}
where $p \in \N$ and the  relations on the right are well-defined in any projective $P_b$.
We consider the following partial order relation on elements of the above chosen basis  of each indecomposable projective. 
\begin{dfn}\label{dfn:poset}
	For each vertex $u \in Q_0$ and  any two well-defined vectors $t_v^m$ and $t_w^n$ in $P_u$  we 
	set
	$t_v^m \leq t_w^n$ if any of the following equivalent conditions are satisfied.
	\begin{enumerate}
		\item $\A t_v^m \ni t_w^n$.
		\item $\A t_m^v \supseteq \A t_n^w$.
		\item $m < n$ or $m = n$ and $(v,w) = (\star,\pm)$ or $m=n$ and $v=w$.
		\item $t_v^m = t_w^n$ or there is a path in the diagram of $P_u$ in \eqref{eq:proj-diag} from $t_v^m$ to $t_w^n$. 
	\end{enumerate}
\end{dfn}

\begin{rmk}
	The proof of Proposition~\ref{prp:fd-cyc} implies that
	for each $u \in Q_0$
	there is an isomorphism of posets
	$$
	(\{t^n_v \mid (v,n) \in I_u\},\leq)
	\cong
	(\{L \subseteq P_u \mid L \text{ a cyclic submodule}\},\supseteq).
	$$
	This isomorphism  provides a natural interpretation of the  partial order introduced above.
\end{rmk}

Let $(\alpha,p,\beta)$ be a triple of a finite-dimensional cyclic $\A$-module, that is,
in addition to $\beta \in Q_0$ we fix
$\alpha \in Q_0 \cup \{\hat{\star}\}$, and $p\in \N_0$
with $p \neq 0$ if  
$(\alpha,\beta) \notin \{ (\hat{\star}, \star), (+,\star), (-,\star)\}$.

The monomial basis of $P_{\beta}$ described above 
gives rise to a basis of the finite-dimensional cyclic 
$\A$-module $M(\alpha,p,\beta) = P_{\beta}/ I(\alpha) e_{\beta}$, where $I(\alpha) = \A t^p_{\alpha}$ if $\alpha \in Q_0$ respectively $I(\hat{\star}) =  \A(t_+^p,t_-^p )$ otherwise.
Let $g \in M(\alpha,p,\beta)$ denote the residue class $x + \im \partial_1$.
\begin{dfn}
	For a triple  $(\alpha,p,\beta)$ as above,
	we 
	define a basis
	$\mathfrak{B}(t_{\alpha}^{p} g,\beta)$ of the cyclic $\A$-module $M(\alpha,p,\beta)$
	by 	
	\begin{align*}\mathfrak{B}(t_\alpha^{p} g,\beta) &=
		\begin{cases}
			\{ \ol{t^n_v g} \mid (v,n) \in I_{\beta}: t^{p}_+ \not\leq t^n_v \text{ and } t^{p}_- \not\leq t^n_v 
			\} & \text{if }a = \hat{\star},
			\\
			\{ \ol{t^n_v g} \mid (v,n) \in I_b: t^{p}_{\alpha} \not\leq t^n_v 
			\} & \text{otherwise},
		\end{cases}
	\end{align*}
	where $\ol{t^n_v g}$ denotes the residue class $t^n_v g + \im \partial_1$.
\end{dfn}

\begin{dfn}\label{dfn:std-rel}
	Given a triple $(\alpha,p,\beta)$ as above,
	let $(r_i)_{i \in Q_0}$ and $(s_i)_{i \in Q_0}$ denote the non-negative integers such
	that
	$$\mathfrak{B}(t_a^{p }g,b) = \{ (t^n_+ g)_{n=r_+}^{s_+}, (t^n_\star g)_{n=r_\star}^{s_\star}, (t^n_- g)_{n=r_-}^{s_-}\}.$$
	We will say that the \emph{standard relations} of the basis vectors in 
	$\mathfrak{B}(t_{\alpha}^{p}g,\beta)$
	are given by the equations:
	\begin{align*}
			B_- t_{-}^n g&= 
			t_{\star}^{n+1} g&
			&(r_- \leq n < s_\star)
			&&&
			B_+ t_{+}^n g&=
			t_{\star}^{n+1} g&
			&
			(r_+ \leq n < s_\star)
			\\
			A_- t_{\star}^n g &= 
			t_-^{n} g&
			&
			(r_\star \leq n \leq s_-)
			&&&
			A_+ t_{\star}^n g &= 
			t_{+}^{n} g&
			&
			(r_\star \leq n \leq s_+)&
		\end{align*}
	\end{dfn}

	\begin{dfn}\label{dfn:glu-vec}
		For a triple $(\alpha,p,\beta)$ as above, an element $z \in  \mathfrak{B}(t_{\alpha}^{p } g,\beta)$
		will be called a \emph{boundary vector}
		if
			$B_{\pm} z =0$ when
			$z \in e_{\pm} M$, respectively  
			$A_+ z = 0$ or $A_- z = 0$ when $z \in e_\star M$,
		where $B_+, B_-, A_+$ and $A_-$ denote the matrices in the quiver representation of $M = M(\alpha,p,\beta)$.
	\end{dfn}

	More explicitly, 
	for all possible triples $(\alpha,p,\beta)$,
	the associated projective presentation, quiver representation $M(\alpha,p,\beta)$, and its basis $\mathfrak{B}(t_{\alpha}^{p}g,\beta)$ together with the boundary vectors
	are described in Table~\ref{tab:cyc-base} on page~\pageref{tab:cyc-base}.
	In the next remark, we give some motivation for the subsequent sections.
	\begin{rmk}
		\label{rmk:base-motiv}
		Let $\omega$ be a band or string of $\Acat$.
		We recall that the corresponding band respectively string resolution admits the so-called cyclification $\CX^{cyc}(\omega)$ introduced in Definition~\ref{dfn:cyc}.
		In particular, there are canonical projections
		$$
		\begin{cd}
			M^{cyc}(\omega)\&  \ar[->>]{l}[swap]{\pi^{cyc}} X_0^{cyc}(\omega)  = X_0(\omega) \ar[->>]{r}{\pi} \& M(\omega)
		\end{cd}
		$$
		Let $\tilde{\mathfrak{B}}$ denote the set of monomial basis vectors $t^n_v g$ of the projective module $\X_0(\omega)$ such that $\pi^{cyc}(t^n_v g) \neq 0$.
		Then $\pi^{cyc}(\tilde{\mathfrak{B}})$ is a basis of 
		the direct sum of cyclic modules $M^{cyc}(\omega)$.
		We will show ultimately that $\pi(\tilde{\mathfrak{B}})$ is also
		a basis, for the $\kk$-linear space $M(\omega)$.
	\end{rmk}

\subsection{Explicit band and string representations}

To deduce that a certain
candidate set $\mathfrak{B}$ is a basis of $M(\omega)$ we will use the following simple observation.
\begin{lem}\label{lem:gen2}
Let $\omega$ be a band or string of $\Acat$ and $M = M(\omega)$ its band or string representation. A
subset	 $\mathfrak{B}$ of $M$ is a $\kk$-linear basis of $M$ if it satisfies the following conditions.
\begin{enumerate}
\item The set $\mathfrak{B}$ contains a set of $\A$-linear generators of $M$.
\item 		The subspace $U \colonequals \langle \mathfrak{B} \rangle_{\kk}$ is a $\A$-submodule of $M$,
or, equivalently, 
\begin{itemize}
\item 
any element $x \in \mathfrak{B}$ is contained in $e_{v} M$ for some vertex $v \in Q_0$.
\item 
any element $x_\pm \in e_{\pm} \mathfrak{B}$
satisfies
$B_{\pm} x_\pm \in  U$, and 
\item any element $x_\star \in e_\star \mathfrak{B}$ satisfies
$A_+ x_\star, A_- x_\star \in U$.
\end{itemize}
\item It holds that $\lvert \mathfrak{B} \rvert = \dim M^{cyc}(\omega)$.
\end{enumerate}
\end{lem}
\begin{proof}
As $M$ is finite-dimensional, so is $U$ and there is a number $n \in \N$ 
with $\rad^n U = 0$.
The assumptions on $\mathfrak{B}$ imply that  
$M = \A U = (\A/\rad^n \A) U = U/\rad^n U= U$.
The third  property and 
Proposition~\ref{prp:cyc} imply that $\lvert \mathfrak{B} \rvert= \dim M$, and thus $\mathfrak{B}$ is a basis.
\end{proof}
For any numbers $m, n \in \N_0$ we set $\delta^{m < n} = 1$ if $m < n$ and $\delta^{m < n} = 0$ otherwise, 
and define $\delta^{m > n}$ in a similar way. This notation will be used frequently below.
\subsection{Band representations}

In this subsection, we give explicit descriptions of quiver representation $M(\omega)$ for any band $\omega$ of $\Acat$. We will need to distinguish between two classes of bands.

\subsubsection{Bands with two-periodic words}
Let $\omega$ be a band $(w,m,\lambda)$ of $\Acat$
such that $w$ is a $2$-periodic word.
Using equivalence of bands, we may assume that
$\dot{w} = (\ora{p},\ola{q})$ with $p \geq q$.
The remaining parameters are given 
by any $m \in \N$ and $\lambda \in \kk^*$ with $\lambda \neq -1$ if $p=q$.
The band representation $M(\omega)$ 
admits the following description.
\begin{prp}\label{prp:sh-band}
For the band $\omega$ as above, 
set
$\mathfrak{B} = \bigcup_{j=1}^m \mathfrak{B}_j$ with 
\begin{align*}
\mathfrak{B}_{j} &= 
\{ 
(\mathbf{u}^n_j)_{n=0}^{p-1},
(\mathbf{v}^n_j)_{n=1}^{p},
(\mathbf{w}^n_j)_{n=1}^{p},
(\mathbf{x}^n_j)_{n=1}^{q},
(\mathbf{y}^n_j)_{n=1}^{q},
(\mathbf{z}^n_j)_{n=0}^{q-1}
\}.
\end{align*}
Then $\mathfrak{B}$ can be identified with a basis of the quiver representation $M(\omega)$ and
the matrices $B_-$, $B_+$, $A_-$ and $A_+$ 
of $M(\omega)$
with respect to this basis
are given by the following series of relations, valid for any index $1 \leq j \leq m$.
\begin{enumerate}
\item \label{shb-std} 
Standard relations:
\begin{align*}
\begin{array}{llcllcll}
e_- \mathbf{u}_j^n = \mathbf{u}_j^n
&
{\textstyle 
(0\leq n < p)}
&&
e_\star \mathbf{v}_j^n = \mathbf{v}_j^n
&
{\textstyle 
(1\leq n \leq p)}
&&
e_+ \mathbf{w}_j^n = \mathbf{w}_j^n
&
{\textstyle 
(1\leq n \leq p)}
\\
e_- \mathbf{x}_j^n = \mathbf{x}_j^n
&
{\textstyle 
(1\leq n \leq q)}
&&
e_\star \mathbf{y}_j^n = \mathbf{y}_j^n
&
{\textstyle 
(1\leq n \leq q)}
&&
e_+ \mathbf{z}_j^n = \mathbf{z}_j^n			
&
{\textstyle 
(0\leq n < q)}
\end{array}
\end{align*}
\begin{align*}
\begin{array}{llcllcllcll}
B_- \mathbf{u}^n_{j }= 
\mathbf{v}^{n+1}_{j}
&(0 \leq n < p)
&\qquad&
B_+ \mathbf{w}_j^n =
\mathbf{v}^{n+1}_{j}
&
(1 \leq n < p)
\\
A_- \mathbf{v}^n_{j} = 
\mathbf{u}^{n}_{j}
&
(1 \leq n < p)
&&
A_+ \mathbf{v}^n_{j} = 
\mathbf{w}^{n}_{j}
&
(0 \leq n \leq p)
\\
B_- \mathbf{x}^n_{j}=
\mathbf{y}^{n+1}_{j}
&
(0 \leq n < q)
&&
B_+ \mathbf{z}^n_{j}=
\mathbf{y}^{n+1}_{j}
&
(1 \leq n < q)
\\
A_- \mathbf{y}^n_{j} = 
\mathbf{x}^{n}_{j}
&
(1 \leq n \leq q)
&&
A_+ \mathbf{y}^n_{j}= 
\mathbf{z}^{n}_{j}
&
(1 \leq n < q)
\end{array}
\end{align*}
\item  \label{shb-glu} Gluing relations:
\begin{align*}
\begin{array}{lcl}
B_- \mathbf{x}^q_j = 	
- \delta^{p > q}
\mathbf{v}^{q+1}_j
&&
B_+ \mathbf{w}^p_j  = 0  
\\
A_- 
\mathbf{v}^p_j = 
\begin{cases}
-\mathbf{\hat{u}}^q_j
-\mathbf{\check{x}}^q_j& \text{if } p > q \\
-\frac{\lambda}{\lambda+1} \mathbf{x}^p_{j} 
-\sum_{i=1}^{j-1} (\frac{-1}{\lambda+1})^{j-i} \mathbf{x}^p_{i} 
& \text{if }p=q
\end{cases}
&&	
A_+ 
\mathbf{y}^q_j = -
\mathbf{w}^q_j
\end{array}
\end{align*}
where
$\mathbf{\hat{u}}^q_j \colonequals 
\lambda \mathbf{u}^q_j + \mathbf{u}_{j-1}^q$,
with $\mathbf{u}^q_0 \colonequals 0$, 
and $\mathbf{\check{x}}^q_j =
\frac{1}{\lambda} \mathbf{x}^q_j + \mathbf{x}_{j-1}^q
$ with $\mathbf{x}^q_0 \colonequals 0$.

\end{enumerate}
\end{prp}
The proof of the statement yields a template for the case of other band and arbitrary string representations. 
\begin{proof}
The band resolution $\CX(\omega)$ 
has the form $(\begin{td} \partial_1\colon \X_1(\omega) \ar{r} \& \X_0(\omega)\end{td})$ and 
was described in~\eqref{eq:sh-bd-res}.
For any index $1 \leq j \leq m$ let ${g}_j$ denote the idempotent generator 
of the $j$-th copy of $P_+$, 
and ${h}_j$ the idempotent generator of the $j$-th copy of $P_-$ in $\X_0(\omega) = P_+^m \oplus P_-^m$.
By definition, $M(\omega) = \coker \partial_1$.
For each index $1 \leq j \leq m$ 
and any number $0\leq n_1 < p$
respectively $1 \leq n_2 \leq p$
we identify each of the elements $
\mathbf{u}_j^{n_1}, \mathbf{v}_j^{n_2},	\mathbf{w}_j^{n_2}$
with the residue class $t^{n_1}_- {g}_j + \im \partial_1$, 
$t^{n_2}_\star {g}_j + \im \partial_1$, respectively
$t^{n_2}_+ {g}_j + \im \partial_1$, 
in $M(\omega)$.
Similarly,
$		\mathbf{x}_j^{n_1}, \mathbf{y}_j^{n_1},	\mathbf{z}_j^{n_2}$
are identified with residue classes of $t_-^{n_1} h_j$, $t_\star^{n_1} h_j$ respectively $t_+^{n_2} h_j$.
Next, we establish the relations in \eqref{shb-std} and \eqref{shb-glu}. 
\begin{enumerate}
\item 
For any $0 \leq n <p$ it holds that $B_+ t_+^n {g}_j = t_\star^{n_1} {g}_j$ in $\X_0(\omega)$ by definitions. 
Passing to the quotient $M(\omega)$ we obtain the first relation in \eqref{shb-std}, under the identifications above.
The remaining relations
in \eqref{shb-std} follow similarly.
\item 
Taking the differential in \eqref{eq:sh-bd-res} into account,
the $\A$-module $M(\omega)$ can be written as the quotient of $\bigoplus_{j=1}^{m} \A {g}_j $
modulo its $\A$-submodule generated by the set $\bigcup_{j=1}^{m} \{t_-^p {g}_{j} - R_{j}, 
t_+^q {h}_{j} + t_+^q {g}_j\}
$
with 
$R_{j} = 
- t_-^{q}  \hat{g}_j - t_-^{q} \check{h}_{j}  $ if  $p > q$ \text{respectively}
$R_j = 
-\frac{1}{\lambda+1} (t_-^p {g}_{j-1} + t_-^p \check{h}_{j})$ if $p=q$.
For any index $1 \leq j \leq m$
we denote by $e_j$ and $f_j$ the idempotent generators of the $j$-th copy of $P_-$ respectively $P_+$ in $\X_1(\omega)$, 
and define $\hat{f}_j$ similar to ${\hat{x}}_j$.
In the notation below, we consider all congruences modulo $\im \partial_1$.
We make two observations.
\begin{enumerate}[ref=(\alph*)]
\item \label{sh-bd-a}
It holds that $\partial_1 (t_\star^1 f_j) = t_{\star}^{q+1} ({g}_j + {h}_j)$, 
and thus $t_\star^{q+1} {h}_j \equiv - t_\star^{q+1} {g}_j$.
\item \label{sh-bd-b} It holds that
$\partial_1(t^1_\star e_j - t_\star^1 \hat{f}_j) = t_{\star}^{p+1} {g}_{j}$,
and thus
$t_{\star}^{p+1} {g}_{j} \equiv 0 $.
\end{enumerate}
For any index $ 1 \leq j \leq m$ it holds that
$A_- t_\star^p {g}_j = t_-^p {g}_j$, $B_+ t_+^p {g}_j = B_- t_-^p {g}_j$,
$A_+ t_\star^q {h}_j = t_+^q {h}_j$, $B_- t_-^q {h}_j = B_+ t_+^q {h}_j$ in $\X_0(\omega)$. 
Using $t_-^p {g}_j \equiv R_j$ and $t_+^q {h}_j \equiv - t_+^q {g}_j$, we obtain the following congruences in the quotient $M(\omega)$.
\begin{align*}
A_- t_\star^p {g}_j \equiv 
R_{j} &&
B_+ t_+^p {g}_j \equiv 
B_- R_{j} &&
A_+ t_\star^q {h}_j \equiv 
-t_+^q {g}_j &&
B_- t_-^q {h}_j 
\equiv -t_\star^{q+1} {g}_j
\end{align*}
Using the  observations~\ref{sh-bd-a} and~\ref{sh-bd-b}, and replacing 
${g}_j$ by ${x}_j$ and ${h}_j$ by ${y}_j$,
we obtain all relations in \eqref{shb-glu}
except the first relation in case $p =q$.

To show the remaining relation, let $p=q$. Set $\mu = \frac{-1}{\lambda+1}$, $a_j = t_-^p {g}_j$ and $b_j = t_-^p {h}_j$.
Since $a_{j-i} \equiv R_{j-i} \equiv -\frac{1}{\lambda+1} (t_-^p g_{j-i-1} + t_-^p \check{h}_{j-i})$ for any index $1 \leq i < j$, 
it follows that
\begin{align*}
\frac{1}{\mu} a_j \equiv 
a_{j-1} + \hat{b}_j
\equiv (\mu (a_{j-2}  +  \hat{b}_{j-1}) +  \hat{b}_j)
\equiv \ldots
\equiv  \sum_{i=1}^j \mu^{j-i} \hat{b}_i
\equiv \sum_{i=1}^{j-1} \mu^{j-i} b_i + \lambda b_j 
\end{align*}
Multiplying this congruence with $\mu$ 
induces the first relation of \eqref{shb-glu} in case $p = q$. 
\end{enumerate}
The right hand side of each relation in \eqref{shb-std} and \eqref{shb-glu} 
is contained in the span $\langle \mathfrak{B} \rangle_\kk$. Moreover, these relations describe completely the action of the operators $a_-, a_+, b_-$ and $b_+$ from $\A$ on the vectors from $\mathfrak{B}$.
Therefore, the assumptions of
Lemma~\ref{lem:gen2} are satisfied, and its application yields 
that $\mathfrak{B}$ is a basis of $M(\omega)$.
\end{proof}

\begin{ex}
For a band $\omega = (w,m,\lambda)$ of period two as above, 
set $n = p + q$.
\begin{itemize}
\item Assume that $m=1$.
Then $\uldim M(\omega) =(n,n,n)$
and the matrices of 
the quiver representation of $M(\omega)$ 
are given by 
\begin{align*}
\begin{array}{ll} 
B_- = J_0 - \delta^{p > q} E_{q+1,n},
&
{ B_+ = {\displaystyle \smash{\sum_{i=1, i \neq p+1}^n E_{i,i}}}},
\\
A_- = {\displaystyle \smash{\sum_{i=1, i \neq p}^n E_{i,i}}-}
\begin{cases}
\lambda (E_{q,n} + E_{n,p})  & \text{if }p>q,\\
\frac{\lambda}{\lambda+1} E_{n,p}& \text{if }p=q,
\end{cases}
&
A_+ = J_0 - E_{q+1,n} ,
\end{array}
\end{align*}
where $J_0$ denote the lower triangular Jordan block of size $n$ and eigenvalue $0$.
\item 
If $m > 1$, then each $\kk$-vector space has dimension $n = m(p+q)$ and each non-zero entry
of any of the four matrices has to be replaced by an identity matrix of size $m$, $\lambda$ by $J_{\lambda,m}$ and $-\frac{\lambda}{\lambda+1} E_{p,p}$
by the matrix $\sum_{j=1}^m \sum_{i=1}^{j-1} (\frac{-1}{\lambda+1})^{j+1-i} E_{ij} - \frac{\lambda}{\lambda+1} {\color{black} E_{n,j}}$ in case $p =q$.
\end{itemize}
\end{ex}
\subsubsection{Band words of higher period}
Next, we consider a band $\omega =(w,m,\lambda)$ of $\Acat$
with a $2k$-periodic word $w$, where $k > 1$.
In particular, $m \in \N$, $\lambda \in \kk^*$ with $\lambda \neq -1$. 
As in the previous section, we may assume that $\dot{w}$ is given by a sequence of numbers  $(\ora{p}_1,\ola{q}_1, \ldots, \ora{p}_k, \ola{q}_k)$ all of which are positive integers.
The band resolution $\CX(\omega)$ and its  differential $\partial_1$ were described in 
\eqref{eq:l-bd-diff}.
For the next statement, we denote
\begin{align}
\label{eq:nota-gen}
\X_0 = \bigoplus_{i=1}^{k} \X_{0,i} \quad\text{with } X_{0,i} = P_+^m \oplus P_-^m = \bigoplus_{j=1}^m \A g_{i,j} \oplus \bigoplus_{j=1}^m  \A h_{i,j},
\end{align}
that is, $g_{i,j}$ denotes the idempotent generator of the $j$-th copy of the summand $P_{+}$ in $\X_{0,i}$, and
$h_{i,j}$ denotes the idempotent generator of the $j$-th copy of $P_-$ in $\X_{0,i}$.

\begin{rmk}	
For the band $\omega$ as above, the cyclification $\CX^{cyc}(\omega)$ 
is given by 
$$
(
\begin{cd}
\X_1 \ar{r}{\partial_1^{cyc}} \& \X_0
\end{cd})
= 
\bigoplus_{i=1}^k
\big((\begin{cd} P_-^m \ar{r}{t^{p_i}} \& P_+^m \end{cd}) \oplus(\begin{cd} P_+^m \ar{r}{t^{q_i}} \& P_-^m \end{cd})\big)
$$
Remark~\ref{rmk:base-motiv}
put forward a candidate for a $\kk$-linear basis of $M(\omega)$ given by 
the image of the set 	
$		\tilde{\mathfrak{B}} = 	\bigcup_{i=1}^{k} \bigcup_{j=1}^m
\{
t_-^{p_i} g_{i,j}, t_+^{q_i} h_{i,j}
\}$ under the canonical projection
$\begin{td}	\X_0 \ar[twoheadrightarrow]{r} \& M(\omega).
\end{td}$
This motivates to consider the diagonal entries of $\partial_1$ in more detail.
\end{rmk}

To proceed, we need first an approximate description of the image of $\partial_1$.
\begin{lem}\label{lem:im-del}
For the band $\omega$ as above
and any indices $1 \leq i \leq k$ and $1\leq j \leq m$,
the image of the differential $\partial_1$ of the band resolution $\CX(\omega)$ contains the elements
\begin{itemize} 
\item $t^n_v g_{i,j} + t^n_v h_{i,j}$ for any $(v,n) \in Q_0 \times \N$
such that 
$t^n_v > t_-^{p_i}$ or $t^n_v > t_+^{q_i}$, as well as
\item $t^n_v g_{i,j} $ and $t^n_v h_{i,j}$ 
for any  $(v,n) \in Q_0 \times \N$ such that
$t^n_v > t_-^{p_i}$ and $t^n_v > t_+^{q_i}$.
\end{itemize}
\end{lem}
\begin{proof}
Similarly to \eqref{eq:nota-gen}, we may denote
$X_1(\omega) = \bigoplus_{i=1}^{k} X_{1,i}$ with $X_{1,i} = P_-^m \oplus P_+^m = \bigoplus_{j=1}^m \A e_{i,j} \oplus \bigoplus_{j=1}^m  \A f_{i,j}$.
Let $1 \leq i \leq k$ and $1 \leq j \leq m$.
It can be checked that
\begin{align*}
\partial_1(-\beta_i t_\star^1 f_{i,j} + t_\star^1 e_{i+1,j}) &= t_{\star}^{p_{i+1}+1} g_{i+1,j} + \ol{\alpha}_{i+1} t_{\star}^{p_{i+1}+1} h_{i+1,j} &&\text{if }i \neq k, \\
\partial_1(t_{\star}^1 f_{i,j} - \ol{\beta}_i t_{\star}^1 e_{i+1,j}) &= \alpha_i t_{\star}^{q_i + 1} g_{i,j} + t_{\star}^{q_i + 1} h_{i,j}&&\text{if }i \neq k,\\
\partial_1(-\beta_k t_\star^1 (\lambda f_{k,j} + f_{k,j-1}) + t_\star^1 e_{1,j}) &= t_{\star}^{p_{1}+1} g_{1,j} + \ol{\alpha}_{1} t_{\star}^{p_{1}+1} h_{1,j},&&  \text{where } f_{k,0} \colonequals 0,\\
\partial_1(t_{\star}^1 f_{k,j} - \ol{\beta}_k t_{\star}^1 (\tfrac{1}{\lambda} e_{1,j} + e_{1,j-1}) &= \alpha_k t_{\star}^{q_k + 1} g_{k,j} + t_{\star}^{q_k + 1} h_{k,j},&&
\text{where } e_{k,0} \colonequals 0.
\end{align*}
We set $a_i = \min(p_i,q_i)+1 $ and $b_i = \max(p_i,q_i)+1 $, and observe the following.
\begin{itemize}
\item If $(\alpha_i,\ol{\alpha}_i) = (1,0)$, then $p_i \geq q_i$ by Remark~\ref{rmk:kap-imp},
and 
$t_{\star}^{b_i} g_{i,j} $,
$t_{\star}^{a_i } g_{i,j} + t_{\star}^{a_i} h_{i,j} \in \im \partial_1$. 
\item Otherwise $(\alpha_i,\ol{\alpha}_i) =(0,1)$. Then $p_i \leq q_i$, and thus
$t_{\star}^{a_i} g_{i,j} + t_{\star}^{a_i} h_{i,j}$,
$t_{\star}^{b_i} h_{i,j} \in \im \partial_1$.
\end{itemize}
Since the image of $\partial_1$ is a $\A$-submodule, each of the two claims follows.
\end{proof} 

\begin{lem}\label{lem:long-bd}
For the band $\omega$ as above,
the following statements hold.
\begin{enumerate}
\item \label{eq:long-bd1}
There is an integer $i \in \Z$ such that $p_i \neq q_i$ or $\gamma_i = 0$.
\item 
Similarly, there is an integer $j \in \Z$ such that $p_j \neq q_{j}$ or $\ol{\gamma}_j = 0$.
\end{enumerate}
\end{lem}
\begin{proof}
We recall that the sequences $(p_i)_{i \in \Z}$, $(q_i)_{i \in \Z}$, $(\gamma_i)_{i \in \Z}$
$(\ol{\gamma}_i)_{i \in \Z}$ are $k$-periodic.
To show \eqref{eq:long-bd1}, assume that $p_i = q_i$ and $\gamma_i \neq 0$ for any $i \in \Z$. 
Then $\beta_i \neq 0$, which implies that $q_i \leq p_{i+1}$ by Remark~\ref{rmk:kap-imp}, for any $i \in \Z$.
As $w$ is $2k$-periodic with $k \geq 2$ by assumption, it follows that all numbers in $p_1 = q_1 \leq p_{2} = q_{2}\leq \ldots p_{k} = q_k \leq p_1$ are equal. Therefore, $w$ has period two, a contradiction, which shows the first claim.

The second claim is proven similarly, using that $\ol{\gamma}_j \neq 0$ only if $p_{j+1} \geq q_{j+1}$.
\end{proof}
From now, we omit the standard relations involving idempotents, assuming that any vectors with the letter $\mathbf{u}$ or $\mathbf{x}$ belong to $M_-(\omega)$,
vectors with letters $\mathbf{v}$ or $\mathbf{y}$ to $M_\star(\omega)$, and vectors with letters $\mathbf{w}$ or $\mathbf{z}$ to $M_+(\omega)$.  
\begin{prp}\label{prp:bands}
For a band $\omega$ as above, 
we set
$\mathfrak{B} \colonequals \bigcup_{i=1}^k \bigcup_{j=1}^m 
\mathfrak{B}_{i,j}$
with
\begin{align*}
\mathfrak{B}_{i,j} 
&= 
\{ 
(\mathbf{u}^n_{i,j})_{n=0}^{p_i-1},
(\mathbf{v}^n_{i,j})_{n=1}^{p_i},
(\mathbf{w}^n_{i,j})_{n=1}^{p_i},
(\mathbf{x}^n_{i,j})_{n=1}^{q_i},
(\mathbf{y}^n_{i,j})_{n=1}^{q_i},
(\mathbf{z}^n_{i,j})_{n=0}^{q_i-1}
\}.
\end{align*}
For formal reasons, given any integer $a\in \Z$, index $1\leq j\leq m$ 
and $0 \leq n < p_i$ respectively $1 \leq n \leq q_i$,
we set
\begin{itemize}
\item $\mathbf{\hat{u}}^n_{a,j} 
= \mathbf{\hat{u}}^n_{i,j}$ and $\mathbf{\check{x}}^n_{a,j}= \mathbf{\check{x}}^n_{i,j}$, where
$1 \leq i \leq k$ is the unique integer such that $a \equiv_k i$,
\item $\mathbf{\hat{u}}^n_{a,0} = \mathbf{\check{x}}^n_{a,0} = 0$,  
\item 
$\mathbf{\hat{x}}^n_{a,j} \colonequals 
\begin{cases}
\lambda {\mathbf{x}^n_{a,j}}+\mathbf{x}^n_{a,j-1} & \text{if }a \equiv_k 1, \\
\ \mathbf{x}^n_{a,j} & \text{otherwise},
\end{cases}$
\quad and \quad
$	\mathbf{\check{u}}^n_{a,j} \colonequals
\begin{cases}
\frac{1}{\lambda} \mathbf{u}^n_{a,j} + \mathbf{u}^n_{a,j-1} & \text{if }a \equiv_k 0,  \\
\ \mathbf{u}^n_{a,j}& \text{otherwise}.
\end{cases}$
\end{itemize}
Similarly we define $\mathbf{\hat{w}}_{a,j}^n$ and $\mathbf{\check{z}}_{a,j}^n$ for any $a \in \Z$, $1 \leq j \leq m$ and appropriate $n$.

Then $\mathfrak{B}$ can be identified with a $\kk$-linear basis of  the band representation $M(\omega)$  and the matrices $B_-$, $B_+$, $A_-$ and $A_+$ 
of $M(\omega)$
with respect to this basis are given by the following series of relations, valid for any indices $1 \leq i \leq k$ and $1 \leq j \leq m$.
\begin{enumerate}
\item  \label{bd-rel1}
Standard relations:
\begin{align*}
\begin{array}{llcllcllcll}
B_- \mathbf{u}^n_{i,j }= 
\mathbf{v}^{n+1}_{i,j}
&(1 \leq n < p_i)
&\qquad&
B_+ \mathbf{w}_{i,j}^n =
\mathbf{v}^{n+1}_{i,j}
&
(0 \leq n < p_i)
\\
A_- \mathbf{v}^n_{i,j} = 
\mathbf{u}^{n}_{i,j}&
(1 \leq n < p_i)&&
A_+ \mathbf{v}^n_{i,j} = 
\mathbf{w}^{n}_{i,j}
&
(1 \leq n \leq p_i)
\\
B_- \mathbf{x}^n_{i,j }=
\mathbf{y}^{n+1}_{i,j}
&
(0 \leq n < q_i)
&&
B_+ \mathbf{z}^n_{i,j }=
\mathbf{y}^{n+1}_{i,j}
&
(1 \leq n < q_i)
\\
A_- \mathbf{y}^n_{i,j} = 
\mathbf{x}^{n}_{i,j}
&
(1 \leq n \leq q_i)
&&
A_+ \mathbf{y}^n_{i,j} = 
\mathbf{z}^{n}_{i,j}
&
(1 \leq n < q_i)
\end{array}
\end{align*}
\item \label{bd-rel2} Gluing relations:
\begin{align*}
B_- \mathbf{x}^{q_i}_{i,j} &= 
-\delta^{p_i > q_i} \mathbf{v}^{q_i+1}_{i,j} \qquad \qquad
B_+ \mathbf{w}^{p_i}_{i,j} = 
- \delta^{p_i < q_i } \mathbf{y}^{p_i+1}_{i,j}
\\
A_- \mathbf{v}^{p_i}_{i,j} &=
-\ol{\alpha}_i { \mathbf{x}^{p_i}_{i,j}}
+
{\displaystyle
\sum_{a=1}^{b}} 
(-1)^a \beta_{i-a} \mathbf{\check{x}}^{q_{i-a}}_{i-a,j} 
+ (-1)^{b}  \gamma_{i-b} \mathbf{\hat{u}}^{q_{i-b}}_{i-b,j} 
\intertext{where 
$b \colonequals b(i) = \min \{ a \geq 1 \mid
\gamma_{i-a} = 0
\text{ or } 				p_{i-a} \neq q_{i-a}\}$, and}
A_+ \mathbf{y}^{q_i}_{i,j} &= 
-\alpha_i \mathbf{w}^{q_i}_{i,j}
- {\displaystyle \sum_{c=i}^{d}}  (-1)^{c-i} 	 \ol{\beta}_{c} \mathbf{\hat{w}}^{p_{c+1}}_{c+1,j}
- (-1)^{d-i}  \ol{\gamma}_{d} \mathbf{\check{z}}^{p_{d+1}}_{d+1,j}
\end{align*}
where 
$d \colonequals d(i)= \min \{  c \geq i \mid
\ol{\gamma}_{c} = 0
\text{ or } 				p_{c} \neq q_{c}\}$.
\end{enumerate}
\end{prp}
In the statement above, the parameters $b=b(i)$ and $d=d(i)$ are well-defined for any index $1 \leq i \leq k$ because of Lemma~\ref{lem:long-bd}.
\begin{proof}
Using the notation from \eqref{eq:nota-gen}, 
we define the elements $g_{a,j}$, $h_{a,j}$,
for any integer $a \in \Z$ and index $1 \leq j \leq m$.
Moreover,  we define $\hat{g}_{a,j}$ and $\check{h}_{a,j}$ 
by the same procedure as $\mathbf{\hat{u}}^n_{a,j}$
respectively $\mathbf{\check{x}}^n_{a,j}$ ignoring the power $n$.
We identify $\mathbf{u}^n_{a,j}$ with the residue class of $t^n_- g_{a,j}$ modulo the image of $\partial_1$, and, similarly,
$\mathbf{x}_{a,j}^n$ with $t^n_- h_{a,j} + \im \partial_1$
for any indices $1 \leq i\leq k$ and $1 \leq j\leq m$,
whenever the values of $(v,n)$
imply that
$\mathbf{u}_{a,j}^n$ respectively $\mathbf{x}_{a,j}^n$ is a well-defined element of $\mathfrak{B}$. Similar identifications are made with $\mathbf{u}$ replaced by $\mathbf{v}$ or $\mathbf{w}$,
or $\mathbf{x}$ replaced by $\mathbf{y}$ or $\mathbf{z}$.

Applying $\partial_1$ to the idempotent generators of $\X_1$ yield the following congruences in $M(\omega)$
\begin{align}
\label{eq:bd-gen1}
t^{p_i}_- g_{i,j} &\equiv 
-\ol{\alpha}_i t^{p_i}_- h_{i,j}	- \beta_{i-1} t^{q_{i-1}}_- \check{h}_{i-1,j} 
- \gamma_{i-1} t^{q_{i-1}}_- \hat{g}_{i-1,j} && 
\\
\label{eq:bd-gen2}
t^{q_i}_+ h_{i,j} &\equiv  -\alpha_i t^{q_i}_+ g_{i,j}
- \ol{\beta}_{i} t^{p_{i+1}}_+ \hat{g}_{i+1,j} 
- \ol{\gamma}_{i} t^{p_{i+1}}_+ \check{h}_{i+1,j}
\end{align}
where $1 \leq i \leq k$ and $1 \leq j \leq m$.
Next, we verify the relations. 	
\begin{enumerate}
\item
Under the identifications above, 
the standard relations in \eqref{bd-rel1} are induced from their counterparts in $X_0(\omega)$
with $g_{i,j}$ and $h_{i,j}$ instead of $x_{i,j}$ and $y_{i,j}$.
\item 
Equation \eqref{eq:bd-gen1} yields that
\begin{align}\label{eq:im-rel}
A_- t_{\star}^{p_i}g_{i,j}\equiv		{\color{black}t_-^{p_i} g_{i,j}} 
\equiv-\ol{\alpha}_i {\color{black}t^{p_i}_- h_{i,j}}
- \beta_{i-1} {\color{black} t^{q_{i-1}}_- \check{h}_{i-1,j}} 
- \gamma_{i-1} {\color{black} t^{q_{i-1}}_- \hat{g}_{i-1,j}} 
\end{align}
If $p_{i-1} = q_{i-1}$ and $\gamma_{i-1} \neq 0$, the last summand may be replaced further using relation \eqref{eq:bd-gen1}
with index $i-1$ instead of $i$. Iterating this procedure yields the first relation in \eqref{bd-rel2}.

To deduce the second relation in \eqref{bd-rel2},
we claim that
$v_{i-1} = 
\beta_{i-1} {\color{black} t^{q_{i-1}+1}_\star \check{h}_{i-1,j}}  +
\gamma_{i-1} {\color{black} t^{q_{i-1}+1}_{\star} \hat{g}_{i-1,j}} \in \im \partial_1$.
There are the following cases.	
\begin{itemize}
\item
If $\gamma_{i-1} = 1$, then $\beta_{i-1} = 1$, $p_{i-1} \geq q_{i-1}$ and
the first statement of Lemma~\ref{lem:im-del} implies that $v_{i-1} \in \im \partial$. 
\item If $\gamma_{i-1} = 0 \neq \beta_{i-1}$, then $\alpha_{i-1} = 0$, 
and thus $p_{i-1} \leq q_{i-1}$ and $v_{i-1} \in \im \partial$ by the second statement of Lemma~\ref{lem:im-del}.
\item Otherwise $\gamma_{i-1} = \beta_{i-1}=0$ and thus $v_{i-1} = 0$.
\end{itemize}
Thus, $v_{i-1} \equiv 0$ in $M(\omega)$.
Multiplying \eqref{eq:im-rel} with $B_-$  from the left yields
the second congruence in
$
B_+ t_+^{p_i}g_{i,j} \equiv
B_- {\color{black}t_-^{p_i} g_{i,j}}
\equiv
-\ol{\alpha}_i {\color{black}t^{p_i+1}_\star h_{i,j}}
- v_{i-1}
\equiv
-\ol{\alpha}_i {\color{black}t^{p_i+1}_\star h_{i,j}}.
$
Comparing $p_i$ to $q_i$, there are three possibilities.
\begin{itemize}
\item If $p_i < q_i$, Remark~\ref{rmk:kap-imp} implies that $\ol{\alpha}_i = 1$.
\item If $p_i = q_i$, then ${\color{black}t^{p_i+1}_\star h_{i,j}} \in \im \partial$ by the second claim in Lemma~\ref{lem:im-del}.
\item If $p_i > q_i$, then $\ol{\alpha}_i = 0$.
\end{itemize}
The last chain of congruences and these cases yield the second relation in \eqref{bd-rel2}.

The remaining two relations in \eqref{bd-rel2} can be derived by similar arguments using \eqref{eq:bd-gen2}.
\end{enumerate}
Let $\tilde{\mathfrak{B}}$ denote all monomial vectors in $\X_0(\omega)$ of the form $t^n_v g_{i,j}$ or $t^n_{v} h_{i,j}$ which map to basis vectors in $\mathfrak{B}$ 
under the canonical projection $\begin{td} \X_0(\omega)\ar[twoheadrightarrow]{r} \& M(\omega) \end{td}$.
Next, we make the following observations in which $1\leq i \leq k$ and $1 \leq j \leq m$ are arbitrary.
\begin{itemize}
\item If $\ol{\alpha}_i \neq 0$, 
Remark~\ref{rmk:kap-imp}
implies that  $p_i \leq q_i$, and thus $t_-^{p_i} h_{i,j} + \im \partial_1 \in \tilde{\mathfrak{B}}$. If also $p_i < q_i$, then $t_{\star}^{p_i + 1} h_{i,j} \in \tilde{\mathfrak{B}}$.
\item Similarly, $\alpha_i \neq 0$ implies that
$p_i \geq q_i$, and thus $t_+^{q_i} g_{i,j} \in \tilde{\mathfrak{B}}$. If, moreover, $p_i > q_i$, it follows that $t_{\star}^{q_i+1} g_{i,j} \in \tilde{\mathfrak{B}}$.
\item By definition of $\tilde{\mathfrak{B}}$ it holds that $t_+^{p_i} \hat{g}_{i,j}, t_-^{q_i} \check{h}_{i,j} \in \tilde{\mathfrak{B}}$.
\item If $\gamma_{i-b}\neq 0$ then $p_{i -b} > q_{i-b}$ by definition of 
$b = b(i)$ and Remark~\ref{rmk:kap-imp}, and thus $t_-^{q_{i-b}} g_{i-b,j} \in \tilde{\mathfrak{B}}$.
Similarly, $\ol{\gamma}_{i+d} \neq 0$ implies that
$t_+^{p_{i+d}} h_{i+d,j} \in \tilde{\mathfrak{B}}$. 
\end{itemize}
These observations imply that
the right hand side of each relation in \eqref{bd-rel1} and \eqref{bd-rel2} is contained in the span $\langle \mathfrak{B} \rangle_{\kk}$.
Because of this, Lemma~\ref{lem:gen2} yields
that 
$\mathfrak{B}$ can be identified with  a basis of $M(\omega)$, which completes the proof. 
\end{proof}

The last proposition can be used to 
describe the quiver representation of a band representation $M(\omega)$ explicitly, which is demonstrated next.

\begin{ex}
Let $\omega$ be the band representation $(w,m,\lambda)$ 
with
$\dot{w} = (\ora{2},\ola{1},\ora{2},\ola{3})$, $m=1$ and $\lambda \in \kk^*$. 
Then $w^{\updownarrow} = (2\downarrow 1 \downarrow 2 \uparrow 3 \uparrow)$.
The oriented gluing diagram of $\omega$ has the following form. We refer to Subsection~\ref{subsec:ind-arr} for details.
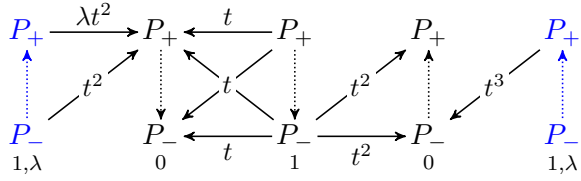
\begin{align*}
\begin{gd}
{{\color{blue} P_{+}}} \ar[color=blue,densely dotted,<-]{d} \ar{r}{\lambda t^2} \& P_{+} \ar[densely dotted,->]{d} \ar[<-]{r}{t} \ar[<-]{rd} \& P_{+}  \ar[densely dotted,->]{d} \&  
P_{+} \ar[densely dotted,<-]{d}  \&   {\color{blue} P_{+}} \ar[densely dotted,<-,color=blue]{d}   \\
\underset{1,\lambda}{{\color{blue} P_{-}}}  \ar[->]{ru}[description]{t^{2}} \& \underset{0}{P_{-}}  \ar[<-]{ru}[description]{t} 
\ar[<-]{r}[swap]{t}
\& \underset{1}{P_{-}} \ar[->]{ru}[description]{t^{2}} \ar{r}[swap]{t^2} \& \underset{0}{P_{-}} \ar[<-]{ru}[description]{t^{3}} 	\& \underset{1,\lambda}{{\color{blue} P_{-}}}
\end{gd}
\end{align*}
It follows that
the band resolution $\CX(\omega)$ 
and the basis $\mathfrak{B}$
provided by
Proposition~\ref{prp:bands}
are given as follows.		
\begin{align*}
\begin{array}{ccc}
\begin{array}{lcl}
\CX(\omega)
=(\begin{cd}\X_1 \ar{r}{\partial_1} \& \X_0 \end{cd}) 
&&
\mathfrak{B} = \bigcup_{i=1}^4 \mathfrak{B}_i\\
\text{where }&&
\text{where }
\\
\partial_1=
\begin{pNiceArray}{cccc}[ first-row, last-col, columns-width = auto]
	P_- & P_+ & P_- & P_+ & \\
	t^2 & t & t & \lambda t^2 & P_+\\
	0 & t & t & 0 & P_- \\
	0 & 0 & t^2 & 0 & P_+ \\
	0 & 0 & t^2 & t^3 & P_-
\end{pNiceArray}
&&
\setlength{\arraycolsep}{1pt} 
\begin{array}{rccr}
\mathfrak{B}_1 = \{& 
\mathbf{u}^1_{1},&
(\mathbf{v}^n_1)_{n=1}^2, &
(\mathbf{w}^n_1)_{n=0}^2 \},
\\
\mathfrak{B}_2 =\{&
(\mathbf{x}^n_{1})_{n=0}^1,&
\mathbf{y}^n_1\}, &
\\
\mathfrak{B}_3 = \{&
\mathbf{u}^1_2, &
(\mathbf{v}^n_2)_{n=1}^2,&
(\mathbf{w}^n_2)_{n=0}^2 \},
\\
\mathfrak{B}_4 = \{ &
(\mathbf{x}^n_{2})_{n=0}^3,&
(\mathbf{y}^n_2)_{n=1}^3,&
(\mathbf{z}^n_2)_{n=1}^2
\}.
\end{array}
\end{array}
\end{array}
\end{align*}
The module $M^{cyc}(\omega)$ is given by the direct sum of the following cyclic representations.
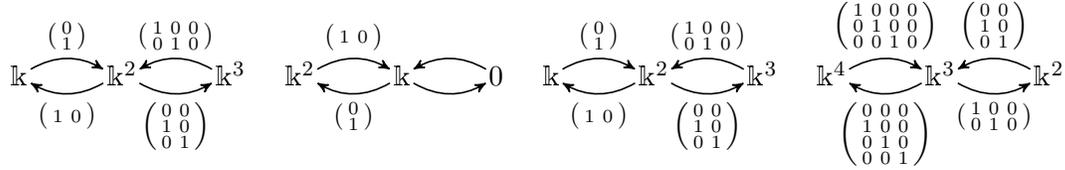
\begin{align*}
\begin{gqrep}{\kk^3}{\kk^2}{\kk}{
\begin{psmallmatrix} 
1 & 0 & 0 \\
0 & 1 & 0 \end{psmallmatrix}}{
\begin{psmallmatrix} 0\\ 1  \end{psmallmatrix}	
}{
\begin{psmallmatrix}0 & 0 \\ 1 & 0 \\ 0 & 1 \end{psmallmatrix}
}{
\begin{psmallmatrix} 1 & 0  \end{psmallmatrix}
}\end{gqrep}
&&
\begin{gqrep}{0}{\kk}{\kk^2}{}{
\begin{psmallmatrix} 1 & 0  \end{psmallmatrix}	
}{}{
\begin{psmallmatrix} 0 \\ 1  \end{psmallmatrix}
}\end{gqrep}&&
\begin{gqrep}{\kk^3}{\kk^2}{\kk}{
\begin{psmallmatrix} 
1 & 0 & 0 \\
0 & 1 & 0 \end{psmallmatrix}
}{
\begin{psmallmatrix} 0\\ 1  \end{psmallmatrix}
}{
\begin{psmallmatrix}0 & 0 \\ 1 & 0 \\ 0 & 1 \end{psmallmatrix}
}{
\begin{psmallmatrix} 1 & 0  \end{psmallmatrix}
}\end{gqrep}&&
\begin{gqrep}{\kk^2}{\kk^3}{\kk^4}{
\begin{psmallmatrix} 0 & 0  \\ 1 & 0 \\0 & 1   \end{psmallmatrix}
}{
\begin{psmallmatrix} 1 & 0 & 0 & 0 \\ 0 & 1 & 0 & 0 \\ 0 & 0 & 1 & 0  \end{psmallmatrix}
}{
\begin{psmallmatrix} 1 & 0 & 0 \\ 0 & 1 & 0  \end{psmallmatrix}
}{
\begin{psmallmatrix} 0 & 0 & 0 \\ 1 & 0 & 0 \\ 0 & 1 & 0 \\ 0 & 0 & 1  \end{psmallmatrix}
}\end{gqrep}		
\end{align*}
According to Proposition~\ref{prp:bands}, $M(\omega)$
admits the same basis as $M^{cyc}(\omega)$ with the standard relations corresponding to the matrix prescriptions above, together with 
the following `gluing relations'.
\begin{align*}
\begin{array}{llll}
B_- \mathbf{x}_1^1 = - \mathbf{v}^2_1	
&
B_+ \mathbf{w}^2_1 =0
&
A_- \mathbf{v}^2_1 = 0 
&
A_+ \mathbf{y}^1_1 = - \mathbf{w}^1_1 
\\
B_- \mathbf{x}^3_2 =0
&
B_+ \mathbf{w}^2_2 = -\mathbf{y}^3_2
&
A_- \mathbf{v}^2_2 = 
-\mathbf{x}^2_2 - \mathbf{x}^1_1 - \mathbf{u}^{1}_1 
&
A_+ \mathbf{y}^3_2 = -\lambda \mathbf{w}^{2}_1
\end{array}
\end{align*}
It follows that $M(\omega)$
is the representation with $\uldim M(\omega) = (8,8,8)$
and the matrices 
\begin{align*}
\scalebox{1.0}{$
\begin{array}{cccc}
B_- & B_+ & A_- & A_+ \\
\begin{pNiceArray}[small]{cc>{\color{green}}:c:ccccc}
0 &  & 0 &  &  & & &  \\
1 &  & -1 & &&&&\\ 
& 1 & 0 && &&& \\ 
&  &0& 0 &  &  \\
&&0& 1 &  &  \\ 
&&0&& 1 & 0 & 0 & 0 \\
&&0& & 0 & 1 & 0 & 0 \\
&&0& & 0 & 0 & 1 & 0 
\CodeAfter 
\tikz \draw [blue] (1-|1) rectangle (3-|2) rectangle (4-|4) rectangle (6-|5) rectangle (9-|9) ;
\end{pNiceArray},&
\begin{pNiceArray}[small]{cc>{\color{green}}:c:cc>{\color{green}}:ccc}
1 & 0 & 0 &   &    &0  &  &  \\
0 & 1 & 0 &  &&0 &&  \\
&  & 0 &  &  & 0&  &  \\
&&0& 1 & 0 & 0 &  & \\
& &0& 0 & 1 & 0 &  & \\
&&0&  &&0& 0 & 0 \\
&&0& &&0& 1 & 0 \\
&&0& &&-1& 0 & 1
\CodeAfter 
\tikz \draw [blue] (1-|1) rectangle (3-|4) rectangle (4-|4) rectangle (6-|7) rectangle (9-|9) ;
\end{pNiceArray},
&
\begin{pNiceArray}[small]{cccc>{\color{green}}:c:ccc}
1 & 0 &  &  & -1 &  &  &  \\
&  & 0 &  & 0 &  \\
&  & 1 & & -1  \\
&&& 1 & 0 &  \\
&& &  &0& 0 & 0 & 0 \\
&&& &0& 1 & 0 & 0 \\
&&& &-1& 0 & 1 & 0 \\
&&& &0& 0 & 0 & 1
\CodeAfter 
\tikz \draw [blue] (1-|1) rectangle (2-|3) rectangle (4-|4) rectangle (5-|6) rectangle (9-|9);
\end{pNiceArray},
&
\begin{pNiceArray}[small]{cc>{\color{green}}:c:cccc>{\color{green}}:c} 
0 & 0 & 0 & \Block{3-2}{} & & \Block{3-2}{} && 0\\
1 &  0 & -1 & &&&& 0 \\
0 & 1 & 0 & &&&& -\lambda \\
\Block{3-2}{}
&&0& 0 & 0 & \Block{3-2}{} && 0  \\
&&0& 1 & 0 & && 0  \\
&&0& 0 & 1 & && 0 \\
\Block{2-2}{} &&0& \Block{2-2}{}& & 1 & 0 & 0 \\
&&0& & & 0 & 1 & 0 
\CodeAfter 
\tikz \draw [blue] (1-|1) rectangle (4-|3) rectangle (4-|4) rectangle (7-|6) rectangle (9-|9);
\end{pNiceArray}.
\end{array}$}
\end{align*}
Now assume that $m \in \N$ is arbitrary.
In this case, 
$\uldim M(\omega) = (8m,8m,8m)$,  
and in each of the four matrices above
each zero entry has to be replaced by a square zero matrix of size $m$,
each
entry $\pm 1$ by a matrix $\pm \Id$ of size $m$,
and
$-\lambda$  by  $-\JB_{\lambda}$ where $\JB_{\lambda}$ denotes a Jordan block with eigenvalue $\lambda$ and size $m$.
\end{ex}

\subsection{String representations}
Let $\omega$ be a string of $\Acat$. We may assume that it is given by one of the strings in  Table~\ref{tab:rep-str-bands}.
Let $I'_w$ denote the set of all indices $i$ such that $\ora{p}_i$ is defined, and $I''_w$ 
the set of such indices with respect to the powers $\ola{q}_i$, that is, 
\begin{align*}
I'_w = \begin{cases}
\begin{array}{rl}
\{1,2,\ldots, k\} & \text{if }d > 0,\\
\phantom{1, }\{2,\ldots, k\} & \text{if }d =0, 
\end{array}
\end{cases}
&&
I''_w = 
\begin{cases}
\{1,\ldots,k-1, k\} & \text{if }d^{\op} > 0,\\
\{1,\ldots, k-1\} & \text{if }d^{\op} =0. 
\end{cases}
\end{align*} 
Moreover, since the bases of cyclic representations are described in Table~\ref{tab:cyc-base}, 
we assume that $M(\omega)$ is not cyclic, that is,
$I'_w \cap I''_w \neq \emptyset$.

Let $\CX(\omega)$ denote the associated string resolution.
We may write $\X_0(\omega) = \bigoplus_{i \in I'_w} P_{b_i} \oplus \bigoplus_{i \in I''_w} P_{v_i}$
with the vertices $b_i$ and $v_i$ defined in the comments following \eqref{eq:str-pos}.
For any $i \in I'_w$ 
let $g_i$ denote the idempotent generator of the projective summand $P_{b_i}$ in $\X_0(\omega)$.
Similarly, for $i \in I''_w$ let $h_i$ denote the idempotent generator of the projective $P_{v_i}$ in $\X_0(\omega)$.

The first differential $\partial_1$ of the string resolution~$\CX(\omega)$
was described by
\eqref{eq:olap} and its subsequent matrices.
In particular,
$\partial_1$ has a diagonal which induces 
monomial basis vectors $(t^{p_i}_{a_i} g_i)_{i \in I'_w}$ and $(t^{q_i}_{u_i} h_i)_{i\in I''_w}$ in $\X_0(\omega)$.
\begin{rmk}\label{rmk:kap-imp2}
For any integer $i \in I'_w \cap I''_w$ 
the following implications hold.
\begin{align*}
p_i > q_i \Rightarrow 
t^{p_i}_{a_i} > t^{q_i}_{u_i} \Rightarrow (\alpha_i,\ol{\alpha}_i) = (1,0) 
&&
q_i > p_i \Rightarrow
t^{q_i}_{u_i} > t^{p_i}_{a_i} \Rightarrow (\alpha_i,\ol{\alpha}_i) = (0,1)  
\end{align*}
These implications give a refinement of Remark~\ref{rmk:kap-imp}.
\end{rmk}
The vectors
$(t^{p_i}_{a_i} g_i)_{i \in I'_w}$ and $(t^{q_i}_{u_i} h_i)_{i\in I''_w}$
generate the image of $\partial_1^{cyc}$.
With
respect to the image of $\partial_1$, we have the following 
analogue of Lemma~\ref{lem:im-del} for the string $\omega$.

\begin{lem}\label{lem:im-del2}
In the setup above, the image of the differential $\partial_1$ contains the vectors 
\begin{enumerate}
\item $t^n_v h_1$  for any $(v,n) \in I_{v_1}$ with $t^n_v > t^{q_1}_{u_1} $ if $d=0$,
\item $t^n_v g_k$  for any $(v,n) \in I_{b_k}$ with $t^n_v > t^{p_k}_{a_k}$ if $d^{\op} = 0$,
\item and, for any index $i \in I'_w \cap I''_w$,
\begin{enumerate}
\item \label{im-del2a} $t^n_v g_{i} + t^n_v h_{i}$ for any $(v,n) \in Q_0 \times \N$  with $t^n_v > t^{p_i}_{a_i}$ or $t^n_v > t^{q_i}_{u_i}$, 
\item \label{im-del2b} $t^n_v g_{i} $ and $t^n_v h_{i}$ for any $(v,n) \in Q_0 \times \N$ with 
$t^n_v > t^{p_i}_{a_i}$ and $t^n_v > t^{q_i}_{u_i}$.
\end{enumerate} 
\end{enumerate}
\end{lem}
\begin{proof}
For any index $1 < i < k$
set $X{(i)} = P_{a_i} \oplus P_{u_i}$ 
and denote by $e_i$ and $f_i$ the corresponding idempotent generators.
Let  $X{(1)} = P_{u_0} \oplus P_{a_1} \oplus P_{u_1}$ if $d=2$,
$P_{a_1} \oplus P_{u_1}$ if $d=1$ and $P_{u_1}$ if $d=0$. 
We denote the generators of $X(1)$ by $\{f_0, e_1, f_1\}$ in the first, by $\{e_1, f_1\}$ in the second and by $\{f_1\}$ in the third case.
Similarly, $X{(k)}$ is set to $P_{a_k} \oplus P_{u_k} \oplus P_{a_{k+1}} = \langle e_k, f_k, e_{k+1} \rangle_{\A}$
if $d^{\op} =2$, to
$P_{a_k} \oplus P_{u_k} = \langle e_k, f_k \rangle_{\A}$ if $d^{\op}=1$
and to $P_{a_k} = \langle e_k \rangle_{\A}$ if $d^{\op} =0$.
With this notation, 
$X_1(\omega) = \bigoplus_{i=1}^{k} X{(i)}$.
We show the first two claims.
\begin{enumerate}
\item 
Assume that  $d=0$. Then $\partial_1( t_\star^1 f_1 - \ol{\beta}_1 t_\star^1 e_2) = t_\star^{{q_1 } + 1} h_1$. 
Note that $t_\star^{{q_1 } + 1}$ is the direct successor of $t^{q_1}_{u_1}$, that is, 
$t^n_v > t^{q_1}_{u_1}$
if and only if $t^n_v \geq t_\star^{{q_1 } + 1}$.
This shows the first claim.
\item 
If $d^{\op} = 0$, the second claim follows because $\partial_1(t_\star^1 e_k -\beta_{k-1} t_\star^1 f_{k-1}) = t_\star^{{p_k }+1} g_{k}$.
\end{enumerate}
Let $i \in I'_w \cap I''_w$ and $(v,n) \in Q_0 \times \N$. Next, we claim that $t^{n}_v g_i + \ol{\alpha}_i t^n_v h_i \in \im \partial_1$ 
if $t^n_v > t^{p_i}_{a_i}$. 
\begin{itemize}
\item
If $i =1$, then $d \in \{1,2\}$.
\begin{itemize}
\item If $d=2$ or $d=1$ and $a_1 = \pm$, then $\partial_1(t_{\star}^1 e_1) = t_\star^{p_1+1}g_1 + \ol{\alpha}_1 t_\star^{p_1+1} h_1$.
\item If $d=1$ and $a_1=\star$,
then 
$\partial_1(t_{\pm}^0 e_1) = t_\pm^{p_1} g_1 + \ol{\alpha}_1 t_\pm^{p_1} h_1$.
\end{itemize}	
\item
If $i \neq 1$, it holds that
$\partial_1(t_\star^1 e_{i}-\beta_{i-1} t_\star^1 f_{i-1}) =
t_{\star}^{p_{i}+1} g_{i} + \ol{\alpha}_{i} t_{\star}^{p_{i}+1} h_{i}$.
\end{itemize}
In each case above, the resulting linear combinations 
$t^{n}_v g_i + \ol{\alpha}_i t^n_v h_i$ contain all direct successors $t^n_v$ of $t^{p_i}_{a_i}$ in the partial order, which shows the claim.

Similarly, we claim that $\alpha_i t^n_v g_i + t^n_v h_i \in \im \partial_1$ if $t^n_v > t^{q_i}_{u_i}$. 
This follows from the following observations.
\begin{itemize}
\item If $i \neq k$, then 
$\partial_1(t_{\star}^1 f_{i} - \ol{\beta}_i t_{\star}^1 e_{i+1}) = \alpha_i t_{\star}^{q_i + 1} g_{i} + t_{\star}^{q_i + 1} h_{i}$.
\item If $i = k$, then $d^{\op} \in \{1,2\}$.
\begin{itemize}
\item If $d^{\op}=2$ or $d^{\op}=1$ and $a_k = \pm$, then $\partial_1(t_{\star}^1 f_k) = \alpha_k t_\star^{q_k+1} g_k +  t_\star^{q_k+1} h_k$.
\item If $d^{\op}=1$ and $a_k=\star$,
then 
$\partial_1(t_{\pm}^0 f_k) = \alpha_k t_\pm^{q_k}g_k +  t_\pm^{q_k} h_k$.
\end{itemize}	
\end{itemize}
Now we consider the parameter $\alpha_i$.
\begin{itemize}
\item Let $\alpha_i = 1$.   
Then $t^{p_i}_{a'_i} \not< t^{q_i}_{u_i}$ by Remark~\ref{rmk:kap-imp2}. Thus, $t^{p_i}_{a_i} \geq t^{q_i}_{u_i}$ or
$(t^{p_i}_{a_i},t^{q_i}_{u_i}) = (t_{\pm}^{p_i},t_{\mp}^{p_i})$.

\begin{itemize}
\item If $t^n_v > t^{p_i}_{a_i}$, then $t^n_v > t^{q_i}_{u_i}$ by the previous observation, and $t_{n}^{v} g_{i} \in \im \partial_1$ as $\ol{\alpha}_i=0$.
\item If $t^n_v > t^{q_i}_{u_i}$, then 
$t_{n}^{v} g_{i} + t_{n}^{v} h_{i} \in \im \partial_1$ by the second claim.
\end{itemize}
\item Otherwise $\alpha_i=0$. Then $\ol{\alpha}_i = 1$, which implies that $t^{p_i}_{a_i} \leq t^{q_i}_{u_i}$ or $(t^{p_i}_{a_i},t^{q_i}_{u_i}) = (t_{\pm}^{p_i}, t_{\mp}^{p_i})$.
\begin{itemize}
\item If $t^n_v > t^{q_i}_{u_i}$, then $t^n_v > t^{p_i}_{a_i}$ and $t^n_v h_i \in \im \partial_1$.
\item If $t^n_v > t^{p_i}_{a_i}$, then $t_{n}^{v} g_{i} + t_{n}^{v} h_{i} \in \im \partial_1$. 
\end{itemize}
\end{itemize}
As $\im \partial_1$ is a $\kk$-linear subspace, the remaining two claims follow in both cases.
\end{proof}

\begin{prp}\label{prp:str-base}
In the setup of the string $\omega$ above, 
set $\mathfrak{B} = \bigcup_{i \in I'_w} \mathfrak{B}'_i
\cup \bigcup_{i\in I''_w} \mathfrak{B}''_i$
with
\begin{align*}
\mathfrak{B}'_i = 
\{
(\mathbf{u}_{i}^{n})_{n=b_-}^{p_i-a_-},
(\mathbf{v}_{i}^{n})_{n=b_\star}^{p_i-a_\star},
(\mathbf{w}_{i}^{n})_{n=b_+}^{p_i-a_+}
\}, 
&&
\mathfrak{B}''_i = 
\{
(\mathbf{x}_{i}^{n})_{n=v_-}^{q_i-u_-},
(\mathbf{y}_{i}^{n})_{n=v_\star}^{q_i-u_\star},
(\mathbf{z}_{i}^{n})_{n=v_+}^{q_i-u_+}
\},
\end{align*}
where the limits of the sequences in $\mathfrak{B}'_i$ are determined by
\begin{align}
\label{eq:delta-lim}
(a_-,a_\star,a_+) = 
\begin{cases}
(0,0,1) & \text{if }a'_i = -,\\
(1,0,0) & \text{if }a'_i = +,\\
(1,1,1) & \text{if }a'_i = \star,\\
(1,0,1) & \text{if }a'_i = \hat{\star},
\end{cases}
&&
(b_-,b_\star,b_+) = 
\begin{cases}
(0,1,1) & \text{if }b_i=-, \\
(1,1,0) & \text{if }b_i = +, \\
(0,0,0) & \text{if }b_i=\star,
\end{cases}
\end{align}
and the limits of the sequences in $\mathfrak{B}''_i$ are defined analogously by the vertices $u_i$ and $v_i$.

Then $\mathfrak{B}$ can be identified with a basis of the string representation $M(\omega)$
and the matrices $B_-$, $B_+$, $A_-$ and $A_+$ of $M(\omega)$ with respect to the basis $\mathfrak{B}$ are given by the following 
prescriptions.
\begin{enumerate}
\item \label{str-std-rel} 
For each index $i \in I'_w$ the vectors 
in $\mathfrak{B}'_i$ satisfy the 
standard relations described in Definition~\ref{dfn:std-rel}.
The same is true for vectors from 
$\mathfrak{B}''_i$ for each
index $i \in I''_w$.
\item \label{eq:beg-rel}
The gluing vectors in $\mathfrak{B}'_1$ satisfy the following relations.
\begin{itemize}
\item If $d=2$, then
$A_- \mathbf{v}_{1}^{p_1} = -\ol{\alpha}_1 {\mathbf{x}_{1}^{p_1}}$
and
$A_+ \mathbf{v}_{1}^{p_1} = -\ol{\alpha}_1 {\mathbf{z}^{p_1}_{1}}$.
\item If $d=1$ and $a_1= \star$, then 
$B_+ \mathbf{w}_{1}^{p_1-1}=  B_- \mathbf{u}_{1}^{p_1-1} =
-\ol{\alpha}_1 \mathbf{y}_{1}^{p_1}$.
\item If $d=1$ and $a_1= -$, then
$A_{-} \mathbf{v}_{1}^{p_1} = -\ol{\alpha}_1 \mathbf{x}_{1}^{p_1}$ 
and 
$	B_{+} \mathbf{w}_{1}^{p_1} =
- 
\delta^{p_1 < q_1}
\mathbf{y}_{1}^{p_1+1}$.
\item If $d=1$ and $a_1= +$, then
$A_{+} \mathbf{v}_{1}^{p_1} =-\ol{\alpha}_1 \mathbf{z}_{1}^{p_1}$ 
and 
$	B_{-} \mathbf{u}_{1}^{p_1} =
- 
\delta^{p_1 < q_1}
\mathbf{y}_{1}^{p_1+1}
$.
\end{itemize}			
\item \label{str-iter-rel}
For any index $1 < i \leq k$,
the gluing vectors in $\mathfrak{B}'_i$
satisfy
\begin{align*}
\begin{cases}
B_{+} \mathbf{w}^{p_i}_{i} = 	-\delta^{p_i < q_i} \mathbf{y}^{p_i+1}_{i} 	\quad \text{and}		& \\
A_- \mathbf{v}^{p_i}_{i} =
-\ol{\alpha}_i \mathbf{x}^{p_i}_{i}
+
{\displaystyle
\sum_{a=1}^{b}} 
(-1)^{a} \beta_{i-a} \mathbf{x}^{q_{i-a}}_{i-a}
+ (-1)^{b}  \gamma_{i-b} \mathbf{u}^{q_{i-b}}_{i-b}
& \text{if }a_i = -,\\
B_{-} \mathbf{u}^{p_i}_{i} = 	-\delta^{p_i < q_i} \mathbf{y}^{p_i+1}_{i} \quad \text{and}	 & \\
A_+ \mathbf{v}^{p_i}_{i} =
-\ol{\alpha}_i \mathbf{z}^{p_i}_{i}
+
{\displaystyle
\sum_{a=1}^{b}} 
(-1)^{a} \beta_{i-a} \mathbf{z}^{q_{i-a}}_{i-a}
+ (-1)^{b}  \gamma_{i-b} \mathbf{w}^{q_{i-b}}_{i-b}	
& \text{if }a_i = +,
\end{cases}
\intertext{where 
$b \colonequals  b(i) =
\min \{ 1 \leq a \leq i \mid
\gamma_{i-a} = 0 \text{ or }
p_{i-a} \neq q_{i-a} \}$. \newline Similarly, for any index $1 \leq i < k$ 
the gluing vectors in $\mathfrak{B}''_i$ satisfy}
\begin{cases}
B_{-} \mathbf{x}^{q_i}_{i} =
-
\delta^{p_i > q_i}
\mathbf{v}^{q_i+1}_{i} \quad \text{and}		&\\
A_{+} \mathbf{y}^{q_i}_{i} =
-\alpha_i {\mathbf{w}^{q_i}_{i}}
- {\displaystyle \sum_{c=i}^{d}}  (-1)^{c-i} 	{\ol{\beta}_{c} \mathbf{w}^{p_{c+1}}_{c+1}} 
- (-1)^{d-i} {\ol{\gamma}_{d} \mathbf{z}^{p_{d+1}}_{d+1}}
& 
\text{if }u_i = +,\\
B_{+} \mathbf{z}^{q_i}_{i} =
-
\delta^{p_i > q_i}
\mathbf{v}^{q_i+1}_{i} \quad \text{and}		& \\
A_{-} \mathbf{y}^{q_i}_{i} =
-\alpha_i {\mathbf{u}^{q_i}_{i}}
- {\displaystyle \sum_{c=i}^{d}}  (-1)^{c-i} 	{\ol{\beta}_{c} \mathbf{u}^{p_{c+1}}_{c+1}} 
- (-1)^{d-i} {\ol{\gamma}_{d} \mathbf{x}^{p_{d+1}}_{d+1}}
&  \text{if }u_i = -,
\end{cases}
\end{align*}
where $
d \colonequals d(i)= \min \{  i \leq c \leq k \mid
\ol{\gamma}_{c} = 0 \text{ or }p_{c+1} \neq q_{c+1}\}$.
\item \label{str-end-rel} The gluing vectors in $\mathfrak{B}''_k$ satisfy the following relations.
\begin{itemize}
\item If $d^{\op}=2$, then
$A_- \mathbf{y}^{q_k}_{k} = -{\alpha}_k \mathbf{u}^{q_k}_k$
and
$A_+ \mathbf{y}^{q_k}_{k} = -{\alpha}_k \mathbf{w}^{q_k}_k$.
\item If $d^{\op}=1$ and $u_k = \star$,
then
$B_{-} \mathbf{x}^{q_k-1}_{k} = B_{+} \mathbf{z}^{q_k-1}_{k} =-\alpha_k \mathbf{v}^{q_k}_k$.
\item If $d^{\op}=1$ and $u_k = +$,
then $A_{+} \mathbf{y}^{q_k}_{k} = -{\alpha}_k \mathbf{w}^{q_k}_{k}$ 
and $B_- \mathbf{x}^{q_k}_{k} = - 
\delta^{p_k > q_k} 		 
\mathbf{v}^{q_k+1}_{\star,k}$.
\item If $d^{\op}=1$ and $u_k = -$,
then $A_{-} \mathbf{y}^{q_k}_{k} = -{\alpha}_k \mathbf{u}^{q_k}_{k}$ 
and $B_+ \mathbf{z}^{q_k}_{k} = - 
\delta^{p_k > q_k} 		 
\mathbf{v}^{q_k+1}_{\star,k}$.
\end{itemize}
\end{enumerate}	
\end{prp}
The underlying assumption that $M(\omega)$ is not cyclic merely 
ensures that the coefficients $\delta^{p_1 < q_1}$
and $\delta^{p_k > q_k}$ in the formulations above are well-defined.
As the proof of Proposition~\ref{prp:bands}
is quite similar to the next proof, we focus on the essential steps leaving out rather formal details.
\begin{proof}
For any $i \in I'_w$
and $b_- \leq n \leq p_i - a_-$
the vector $\mathbf{u}^n_i$ 
is identified with the residue class $\ol{t^n_- e_{b_i}} \in 
M(\omega) = \X_0(\omega)/\im \partial_1$.
Similarly, $\mathbf{v}_i^n = \ol{t^n_\star e_{b_i}}$ and $\mathbf{w}_i^n =
\ol{t^n_\star e_{b_i}}$ whenever the former vectors are defined.
In the same spirit, vectors of $\mathfrak{B}''_i$
are given by residue classes  $\ol{t^n_v h_i}$ with appropriate vertex $v$.
For any $i \in I'_w$.
we denote by $e_{a_i}$ the idempotent generator of the summand $P_{a_i}$ of $\X_1(\omega)$.
Next, we verify the relations.
\begin{enumerate}
\item 
With the identifications above, the standard relations are induced by their analogues in $X_0(\omega)$.
\item 	 Assume that $i = 1$. Then $d > 0$ and $t^{p_1}_{a_1} g_1 \equiv \partial_1(e_{a_1}) -\ol{\alpha}_1 t_{a_1}^{p_1}$ above.
\begin{itemize}
\item If $d=2$, then $A_+
t_\star^{p_1} g_1 = t_+^{p_1} g_1 \equiv 
\partial_1(e_{u_0}) - \ol{\alpha}_1 t_+^{p_1} h_1$ and $A_- t_\star^{p_1} g_1 = t^{p_1}_{a_1} g_1$.
\item If $d=1$ and $a_1 = \star$, then $B_+ t_+^{p_1} g_1 = B_- t_-^{p_1} g_1 = t^{p_1}_{a_1} g_1$. 
\item If $d=1$ and $a_1 = \pm$, then $A_{\mp} t_\star^{p_1} g_1 = t^{p_1}_{a_1} g_1$
and $B_{\mp} t^{p_1}_{a_1} g_1 = B_{\pm} t_{\pm}^{p_1} g_1 = - \ol{\alpha}_1 t_\star^{p_1+1} h_1$ which is  contained in $\im \partial_1$ by Lemma~\ref{lem:im-del2}~\eqref{im-del2b}
if $p_1 \leq q_1$. 
\end{itemize}
These considerations yield the relations in \eqref{eq:beg-rel}.
\item Assume that $i \neq 1$. In this case, $a_i \neq \star$, and thus we may denote $a_i = \mp$.  
Resolving $\partial_1( e_{a_i})$ for $t^{p_i}_{a_i} e_{a_i}$  yields
\begin{align}\label{eq:crit-rel}
{\color{black}t_{a_i}^{p_i} g_{i}} 
\equiv -\ol{\alpha}_i {\color{black}t^{p_i}_{a_i} h_{i}}
- \beta_{i-1} {\color{black} t^{q_{i-1}}_{a_i} h_{i-1}} 
- \gamma_{i-1} {\color{black} t^{q_{i-1}}_{a_i} {g}_{i-1}}.
\end{align}
The first and third relation in \eqref{str-iter-rel}
follows  
by resolving $\sum_{a=0}^{b-1} (-1)^{a} \partial_1( e_{i-a})$ for $t^{p_i}_{a_i} g_i$
and writing $t^{p_i}_{a_i} g_i$ as $A_{\mp} t_\star^{p_i} g_i$. 
The second and fourth relation is obtained by application of $B_{\mp}$
to \eqref{eq:crit-rel} and discarding summands by use of Lemma~\ref{lem:im-del2}~\eqref{im-del2a}.
\item The remaining relations in \eqref{str-iter-rel} and \eqref{str-end-rel}
follow by straightforward adaptations of the previous arguments.
\end{enumerate}	
We observe the validity of the following order relations.
\begin{enumerate}
\item 
For any $i \in I'_w $, the residue classes of $\alpha_i t^{q_i}_{u_i}, 
\ol{\beta}_{i-1} t^{p_i}_{u_{i-1}}, \delta^{p_i>q_i}t_\star^{q_i+1}$
are contained in $\mathfrak{B}'_i$ by Remark~\ref{rmk:kap-imp2}, whenever the expressions are well-defined and non-zero, as well as
$t^{p_i}_{a_i} \not> \gamma_{i-b} t_{a_i}^{q_{i-b}} g_{i-b}$
because $p_{i-a} < q_{i-a}$ if $\gamma_{i-b} = 1$ by definition of $b$.
\item 
Similarly, for any $i \in I''_w$
it holds that $t^{q_i}_{u_i} \not> \delta^{p_i < q_i} t_\star^{p_i+1}$, $t^{q_i}_{u_i} \not> \ol{\alpha}_i t^{p_i}_{a_i}$ by Remark~\ref{rmk:kap-imp2} and $t^{q_i}_{u_i} \not> \beta_i t^{q_i}_{a_{i+1}}$, and $t^{q_i}_{u_i} \not> \ol{\gamma}_{i+d} t_{u_i}^{p_{i+1+d}}$ whenever the right hand side is defined and non-zero.
\end{enumerate}
These observations imply that the right hand side of each relation
is contained in $\langle \mathfrak{B} \rangle_{\kk}$.
In particular, $\langle \mathfrak{B} \rangle_{\kk}$ 
is closed under the arrow operators of $\A$ and 
Lemma~\ref{lem:gen2} can be applied to deduce
that $\mathfrak{B}$ is a basis of $M(\omega)$, which  completes the proof.
\end{proof}

\begin{ex}
Consider the special string
$(w,-)$ with 
$w = (\hat{\star},\ora{p}_1,\ola{q}_1,\ora{p}_2,\ola{q}_2)$ where $p_1 \leq q_1 \geq p_2 \leq q_2$, which is a special case of the example in Subsection~\ref{subsec:sp-str}. 
The string resolution $\CX(\omega)$ 
and its first differential are given by 
\eqref{eq:sp-str-d} with
with $q_0=1$ and $\varepsilon = -$.
The representation $M^{cyc}(\omega)$
is the direct sum of the cyclic representations with their monomial bases depicted below.
\begin{align*}
\begin{array}{cc}
\begin{gqrep}{\kk^{p_1}}{\kk^{p_1}}{\kk^{p_1-1}}{\Id}{\begin{psmallmatrix} 0 \\ \Id  \end{psmallmatrix}}{J_0}{\begin{psmallmatrix} \Id & 0  \end{psmallmatrix}}
\end{gqrep} 	&
\begin{gqrep}{\kk^{q_1-1}}{\kk^{q_1}}{\kk^{q_1+1}}{\begin{psmallmatrix} 0 \\ \Id  \end{psmallmatrix}}{\begin{psmallmatrix} \Id & 0  \end{psmallmatrix}}{
	\begin{psmallmatrix}	\Id & 0 \end{psmallmatrix}
}{\begin{psmallmatrix} 0 \\ \Id  \end{psmallmatrix}}\end{gqrep}\\
\{(\mathbf{u}^n_{1})_{n=1}^{p_1-1}, (\mathbf{v}^n_{1})_{n=1}^{p_1}, (\mathbf{w}^n_{1})_{n=0}^{p_1-1}\}
&
\{(\mathbf{x}^n_{1})_{n=0}^{q_1}, (\mathbf{y}^n_{1})_{n=1}^{q_1}, (\mathbf{z}^n_{1})_{n=1}^{q_1-1}\}
\\
\\
\begin{gqrep}{\kk^{p_2+1}}{\kk^{p_2}}{\kk^{p_2-1}}{\begin{psmallmatrix} \Id & 0  \end{psmallmatrix}}{\begin{psmallmatrix} 0 \\ \Id  \end{psmallmatrix}}{\begin{psmallmatrix} 0 \\ \Id  \end{psmallmatrix}}{
	\begin{psmallmatrix}	\Id & 0 \end{psmallmatrix}}\end{gqrep}	
&
\begin{gqrep}{\kk^{q_2}}{\kk^{q_2}}{\kk^{q_2}}{J_0}{\Id}{\Id}{J_0}\end{gqrep}\\
\{(\mathbf{u}^n_{2})_{n=1}^{p_2-1}, (\mathbf{v}^n_{2})_{n=1}^{p_2}, (\mathbf{w}^n_{2})_{n=0}^{p_2}\}
&
\{(\mathbf{x}^n_{2})_{n=0}^{q_2-1}, (\mathbf{y}^n_{2})_{n=1}^{q_2}, (\mathbf{z}^n_{2})_{n=1}^{q_2}\}
\end{array}
\end{align*}
By Proposition~\ref{prp:str-base}, a basis of the string representation $M(\omega)$ is given by 
the union of the four bases above with the basis vectors satisfying the following gluing relations.
\begin{align*}
\begin{array}{llll}
A_+ \mathbf{v}_1^{p_1} = - \mathbf{z}_1^{p_1}  	&
A_+ \mathbf{y}_1^{q_1}  = - \mathbf{w}_2^{p_2}  - \mathbf{z}_2^{p_2} &
A_- \mathbf{v}_2^{p_2}  = - \mathbf{z}_2^{p_2}  & 
A_- \mathbf{y}_2^{q_2} = 0
\\
A_- \mathbf{v}_1^{p_1} = - \mathbf{x}_1^{p_1}  &
B_- \mathbf{x}_1^{q_1} = 0 &
B_+ \mathbf{w}_2^{p_2} = - \mathbf{y}_2^{p_2+1} &
B_+ \mathbf{z}_2^{q_2} = 0.
\end{array}
\end{align*}	
Setting $n = \sum_{i=1}^2 p_i + q_i$, the representation $M(\omega)$ is the quiver representation with
$\uldim M(\omega) = (n-1,n,n)$ 
and the matrices
\begin{align*}
	\begin{array}{lcl}
		B_- = 
		\begin{pNiceArray}[small,last-col,first-row,right-margin, xdots={horizontal-labels}]{ccc|ccc>{\color{green}}|[color=blue]c|ccc|ccc}
			\Hbrace{3}{p_1-1} & \Hbrace{4}{q_1+1} &  \Hbrace{3}{p_2-1} &  \Hbrace{3}{q_2} \\
			0 & \Cdots & 0 & & & & & & & &&&& \Vbrace{4}{p_1} \\
			1 & & & & & & & & & & &&& \\
			& \Ddots & & & & & & & & & &&& \\
			& & 1 & & & & & & & & && \\ \hline
			&&& 1 & && 0 & & & & & & & \Vbrace{3}{q_1}\\
			&&& & \Ddots & & \Vdots & &  & && \\
			&&& &&1 &0 & & & & & \\
			\hline
			&&&& &&& 0 & \Cdots & 0 & && & \Vbrace{4}{p_2} \\
			&&&& &&& 1 & && \\
			&&&& &&& & \Ddots &&  \\
			&&&& &&& & & 1 &\\ \hline
			&&&& &&& &&& 1 & & & \Vbrace{3}{q_2} \\
			&&&& &&& &&&  & \Ddots & & \\
			&&&& &&& &&&  & & 1 & 
		\end{pNiceArray} 					
		&&
		B_+ = 
		\begin{pNiceArray}[small,first-row,last-col,right-margin, xdots={horizontal-labels} ]{ccc|ccc|ccc>{\color{green}}|[color=blue]c|ccc>{\color{green}}|[color=blue]c}
			\Hbrace{3}{p_1} & \Hbrace{3}{q_1-1} &  \Hbrace{4}{p_2+1} &  \Hbrace{4}{q_2} \\
			1 & & & & & & & & & & & & &  & \Vbrace{3}{p_1} \\
			& \Ddots & & & & & & & & & & & &  &\\
			&&1 & & & & & & & & & & &    & \\
			\hline
			& & & 0 & \Cdots & 0 & & & & & & & &  & \Vbrace{4}{q_1} \\
			&&& 1 & & & & & & & & & &  & \\
			&&& & \Ddots & & & & & & & & & & \\
			&&& & & 1 & & & & & & & &  & \\ \hline
			&&&&&& 1 & & & 0 & & & & & \Vbrace{3}{p_2} \\
			&&&&&& & \Ddots && \Vdots & & & & &\\
			&&&&&& & & 1 & 0 & & & & & \\ \hline
			&&&&&& & &  & \Block{4-1}{\scriptstyle y_1}& 0 &\Cdots& 0 & 0 & \Vbrace{4}{q_2}\\[-10pt]
			&&&&&& & &  &  & 1 & \Ddots & \Vdots & \Vdots  \\
			&&&&&& & &  & &  & \Ddots &0 & 0 \\
			&&&&&& & &  & &  & &1 & 0 \\
		\end{pNiceArray} 
	\end{array}			
\end{align*}
\begin{align*}
	A_- = 
	\begin{pNiceArray}[small,last-col,right-margin,first-row, xdots={horizontal-labels}]{ccc>{\color{green}}|[color=blue]c|ccc|ccc>{\color{green}}|[color=blue]c|cccc}
		\Hbrace{4}{p_1} & \Hbrace{3}{q_1} &  \Hbrace{4}{p_2} &  \Hbrace{4}{q_2} \\
		1 & && 0 & & & & & & & & & & && \Vbrace{3}{p_1-1} \\
		& \Ddots & & \Vdots & & & & \\
		&&1 &0  \\
		\hline
		&&&& 0 & \Cdots & 0 &&&&&&&&& \Vbrace{4}{q_1+1} \\
		&&&{\scriptstyle y_5}& 1 & & \\
		&&&&  & \Ddots &  \\
		&&&& & & 1 \\ \hline
		&&&& &&& 1 & && 0 & & & & & \Vbrace{3}{p_2-1} \\
		&&&& &&& & \Ddots & & \Vdots & & & & \\
		&&&& &&& &&1 &0 &   \\
		\hline
		&&&& &&& 	&&&& 0 &\Cdots& \Cdots & 0 & \Vbrace{4}{q_2}\\[-10pt]
		&&&& &&& 	&&&& 1 &\Ddots&& \Vdots \\
		&&&& &&& 	&&&&  &\Ddots &&\Vdots  \\
		&&&& &&& 	&&&&  & &1 & 0 \\
	\end{pNiceArray} 
	&&
	A_+ = 
	\begin{pNiceArray}[small,last-col,right-margin,first-row, xdots={horizontal-labels}]{ccc>{\color{green}}|[color=blue]c|ccc>{\color{green}}|[color=blue]c|ccc|ccc}
		\Hbrace{4}{p_1} & \Hbrace{4}{q_1} &  \Hbrace{3}{p_2} &  \Hbrace{3}{q_2} \\
		0 &\Cdots& 0 & 0  & & & & & & &&& && \Vbrace{4}{p_1} \\[-10pt]
		1 &\Ddots&& \Vdots & & & & & && && &&   \\
		&\Ddots &0 & 0     & & & & && & && & & \\
		& &1 & 0           & & & & && & & &&  &\\ \hline
		&&& & 1 & & & 0 & & & & & & &\Vbrace{3}{q_1-1} \\
		&&& {\scriptstyle y_2}&   & \Ddots && \Vdots & & & && &&   \\
		&&&&   & & 1 & 0 & & & & & &  &\\ \hline
		&&&&& &&   & 0 & \Cdots & 0 & && & \Vbrace{4}{p_2+1}  \\
		&&&&& &&  & 1 & & & & & & \\
		&&&&& && {\scriptstyle y_3} &   & \Ddots & & && &  \\
		&&&&& &&&&& 1 & & & & \\ \hline
		&&&&& &&&&& & 1 & & & \Vbrace{3}{q_2} \\
		&&&&& &&{\scriptstyle y_4}&&& &&\Ddots &&  \\
		&&&&& &&&&& &&&1   \\
	\end{pNiceArray} 					
\end{align*}
with the canonical basis vectors
$y_1 = -e_{p_2+1}$, $y_2 = - e_{p_1}$,
$y_3 = - e_{p_2+1}$,
$y_4 = - e_{p_2}$ and $y_5 = - e_{p_1+1}$.
\end{ex}

\begin{ex}[Primitive root of length one]
Let $\omega$ be the bispecial string $(\varepsilon_1, w,\varepsilon_2)$
such that $w = (\ora{p}, \ola{p}, \ldots )$ has length $\ell$ for some $\ell \in \N$.
Equivalently, $w$ has root $(\ora{p})$ with multiplicity $\ell$.
In this case,
the projective resolution $\CX(\omega)$
is given by
$\begin{td} P(\varepsilon_1) \ar{r}{\partial_1} \& P(\varepsilon_2) \end{td}$ 
where $\partial_1 = \left[C t^{p}\right]$
and
$C = \Id_m + \sum_{j=2}^{m} E_{j-1,j} + \sum_{3 \leq j \leq m \text{ odd}} E_{j-2,j}$.

Set $\mathfrak{B} = \bigcup_{i \in I} \mathfrak{B}'_i \cup \bigcup_{j \in J} 
\mathfrak{B}''_j$
with
\begin{align*}
\mathfrak{B}'_i &\colonequals
\{(\mathbf{u}^{n}_i)_{n=\delta_{2,+}}^{p-\delta_{1,-}},
(\mathbf{v}^{n}_{i})_{n=1}^{p},
(\mathbf{w}^{n}_{i})_{n=\delta_{2,-}}^{p-\delta_{1,+}}\},
&& I = \left\{i \in \N \mid 1 \leq i \leq \lceil \frac{m}{2} \rceil \right\},
\\
\mathfrak{B}''_j &\colonequals
\{(\mathbf{x}^{n}_j)_{n=\delta_{2,-}}^{p-\delta_{1,+}},
(\mathbf{y}^{n}_{j})_{n=1}^{p},
(\mathbf{z}^{n}_{j})_{n=\delta_{2,+}}^{p-\delta_{1,-}}
\},
&&
J = \left\{j \in \N \mid 
1 \leq j \leq \lfloor \frac{m}{2} \rfloor\right\},
\end{align*}
where $\delta_{1,\pm} = \delta_{\varepsilon_1,\pm}$ and
$\delta_{2,\pm} = \delta_{\varepsilon_2, \pm}$.
According to Proposition~\ref{prp:str-base}
a basis of the string representation $M(\omega)$
is given by $\mathfrak{B}$. 
The matrices of $M(\omega)$ are given by 
the standard relations and the following gluing relations
for any indices $i \in I$ and $j \in J$.
\begin{align*}
\begin{cases}
\setlength{\arraycolsep}{5pt}
\begin{array}{lllll}
B_{-} \mathbf{x}^{p}_j = 0 
&
B_{+} \mathbf{w}^p_{i} = 0 
&
A_{-} \mathbf{v}^{p}_i = {\displaystyle \sum_{a=1}^{i-1}} (-1)^a \mathbf{x}^{p}_{i-a} 
&
A_{+} \mathbf{y}^{p}_j = - \mathbf{w}^{p}_j 
&
\text{if }\varepsilon_1= -,\\
B_{-} \mathbf{u}^p_{i} = 0 
&
B_{+} \mathbf{z}^{p}_j = 0 
&
A_{-} \mathbf{y}^{p}_j = - \mathbf{u}^{p}_j 
&
A_{+} \mathbf{v}^{p}_i = {\displaystyle \sum_{a=1}^{i-1}} (-1)^a \mathbf{z}^{p_1}_{i-a} 
&
\text{if }\varepsilon_1= +.
\end{array}
\end{cases}
\end{align*}
In particular, it holds that
$A_{\varepsilon_1} \mathbf{v}^{p}_1 = 0$ for each sign $\varepsilon_1 \in \{+,-\}$.
\end{ex}

At last, we note
that the basis $\mathfrak{B}$ in
each of the 
Propositions~\ref{prp:sh-band}, \ref{prp:bands}, and \ref{prp:str-base}
is precisely the basis $\pi(\tilde{\mathfrak{B}})$ proposed by Remark~\ref{rmk:base-motiv}.
In this sense, the basis $\mathfrak{B}$ of $M(\omega)$ can be identified with the monomial basis of $M^{cyc}(\omega)$
for any band or string $\omega$ of $\Acat$.

\newpage
\clearpage
\renewcommand{\thesection}{A}
\renewcommand{\theequation}{\arabic{equation}}
\setcounter{equation}{0}
\setlength{\arraycolsep}{2pt} 
\begin{table}
	\centering
	\caption{Cyclic representations and their standard bases}
	\label{tab:cyc-base}
	$
	\begin{array}{|c|c|c|c|}
		\hline
		(\alpha,p,\beta)
		&
		\begin{array}{c}
			\text{syzygy}\\
			\text{resolution}
		\end{array}
		&
		\text{quiver representation}
		&
		\begin{array}{c}
			\text{basis }\mathfrak{B}(t_\alpha^p x, \beta) \text{ and }
			\\
			\text{boundary vectors}
		\end{array}
		\\
		\hline
		\hline
		\begin{array}{c}
			(-,p,+)
			\\
			p \in \N
		\end{array}
		& \begin{cd}
			P_- \ar{r}{t^p} \& P_+
		\end{cd}
		&
		\begin{gqrep}{\kk^{p+1}}{\kk^{p}}{\kk^{p-1}}{H}{G}{G}{H}
		\end{gqrep}
		&
		\begin{array}{c}
			\{({\mathbf{u}}^n)_{n=1}^{p-1}, ({\mathbf{v}}^n)_{n=1}^{p}, ({\mathbf{w}}^n)_{n=0}^{p}\}\\
			{\mathbf{v}}^p, 
			{\mathbf{w}}^p
		\end{array}
		\\
		\hline
		\begin{array}{c}
			(+,p,-) 
			\\
			p \in \N
		\end{array}				
		& 
		\begin{cd}
			P_+ \ar{r}{t^p} \& P_-
		\end{cd}
		&
		\begin{gqrep}{\kk^{p-1}}{\kk^{p}}{\kk^{p+1}}{G}{H}{H}{G}
		\end{gqrep}
		&
		\begin{array}{c}
			\{({\mathbf{u}}^n)_{n=0}^{p}, ({\mathbf{v}}^n)_{n=1}^{p}, ({\mathbf{w}}^n)_{n=1}^{p-1}\}
			\\
			{\mathbf{u}}^p, {\mathbf{v}}^p
		\end{array}
		\\	
		\hline
		\begin{array}{c}
			(+,p,+) \\
			p \in \N
		\end{array}
		&
		\begin{cd}
			P_+ \ar{r}{t^p} \& P_+
		\end{cd}
		&
		\begin{gqrep}{\kk^{p}}{\kk^{p}}{\kk^{p}}{\Id}{J}{J}{\Id}
		\end{gqrep} 
		&
		\begin{array}{c}
			\{({\mathbf{u}}^n)_{n=1}^{p}, ({\mathbf{v}}^n)_{n=1}^{p}, ({\mathbf{w}}^n)_{n=0}^{p-1}\} 
			\\
			{\mathbf{u}}^p,
			{\mathbf{w}}^p
		\end{array}
		\\
		\hline
		\begin{array}{c}
			(-,p,-) 
			\\
			p \in \N
		\end{array}				
		&
		\begin{cd}
			P_- \ar{r}{t^p} \& P_-
		\end{cd}
		&
		\begin{gqrep}{\kk^{p}}{\kk^{p}}{\kk^{p}}{J}{\Id}{\Id}{J}
		\end{gqrep}
		&
		\begin{array}{c}
			\{({\mathbf{u}}^n)_{n=0}^{p-1}, ({\mathbf{v}}^n)_{n=1}^{p}, ({\mathbf{w}}^n)_{n=1}^{p}\}
			\\
			{\mathbf{v}}^p,
			{\mathbf{w}}^p
		\end{array}
		\\
		\hline
		\begin{array}{c}
			(\star,p,+)
			\\
			p \in \N
		\end{array}
		&
		\begin{cd}
			P_\star \ar{r}{t^p} \& P_{+}
		\end{cd}
		&
		\begin{gqrep}{\kk^{p}}{\kk^{p-1}}{\kk^{p-1}}{H}{J}{G}{\Id}
		\end{gqrep}
		&
		\begin{array}{c}
			\{({\mathbf{u}}^n)_{n=1}^{p-1}, ({\mathbf{v}}^n)_{n=1}^{p-1}, ({\mathbf{w}}^n)_{n=0}^{p-1}\}
			\\
			{\mathbf{u}}^{p-1},
			{\mathbf{w}}^{p-1}
		\end{array}
		\\
		\hline
		\begin{array}{c}
			(\star,p,-) 
			\\
			p \in \N
		\end{array}
		&
		\begin{cd}
			P_\star \ar{r}{t^p} \& P_{-}
		\end{cd}&
		\begin{gqrep}{\kk^{p-1}}{\kk^{p-1}}{\kk^{p}}{J}{H}{\Id}{G}
		\end{gqrep}
		&
		\begin{array}{c}
			\{({\mathbf{u}}^n)_{n=0}^{p-1}, ({\mathbf{v}}^n)_{n=1}^{p-1}, ({\mathbf{w}}^n)_{n=1}^{p-1}\}
			\\
			{\mathbf{u}}^{p-1},
			{\mathbf{w}}^{p-1}
		\end{array}
		\\
		\hline
		\begin{array}{c}				
			(+,p,\star) 
			\\
			p \in \N_0
		\end{array}
		&
		\begin{cd}
			P_+ \ar{r}{t^{p}} \& P_{\star}
		\end{cd}
		&
		\begin{gqrep}{\kk^{p}}{\kk^{p+1}}{\kk^{p+1}}{G}{J}{H}{\Id}
		\end{gqrep}
		&
		\begin{array}{c}
			\{({\mathbf{u}}^n)_{n=0}^{p}, ({\mathbf{v}}^n)_{n=0}^{p}, ({\mathbf{w}}^n)_{n=0}^{p-1}\}
			\\
			{\mathbf{u}}^p,
			{\mathbf{w}}^p
		\end{array}
		\\
		\hline
		\begin{array}{c}
			(-,p,\star) 
			\\
			p \in \N_0
		\end{array}
		&
		\begin{cd}
			P_- \ar{r}{t^{p}} \& P_{\star}
		\end{cd}
		&
		\begin{gqrep}{\kk^{p+1}}{\kk^{p+1}}{\kk^{p}}{J}{G}{\Id}{H}
		\end{gqrep}
		&
		\begin{array}{c}
			\{({\mathbf{u}}^n)_{n=0}^{p-1}, ({\mathbf{v}}^n)_{n=0}^{p}, ({\mathbf{w}}^n)_{n=0}^{p}\}
			\\
			{\mathbf{v}}^p, {\mathbf{w}}^p
		\end{array}
		\\
		\hline
		\begin{array}{c}
			(\star,p,\star)
			\\
			p \in \N
		\end{array}
		&
		\begin{cd}
			P_\star \ar{r}{t^p} \& P_{\star}
		\end{cd}
		& 
		\begin{gqrep}{\kk^p}{\kk^p}{\kk^p}{J}{J}{\Id}{\Id}
		\end{gqrep}
		&
		\begin{array}{c}
			\{({\mathbf{u}}^n)_{n=0}^{p-1}, ({\mathbf{v}}^n)_{n=0}^{p-1}, ({\mathbf{w}}^n)_{n=0}^{p-1}\}
			\\
			{\mathbf{u}}^{p-1}, {\mathbf{w}}^{p-1}
		\end{array}
		\\
		\hline
		\begin{array}{c}
			(\hat{\star},p,+) 
			\\
			p \in \N 
		\end{array}
		&
		\begin{cd}
			I_\star \ar{r}{t^{p}}
			\& P_+
		\end{cd}
		&
		\begin{gqrep}{\kk^{p}}{\kk^{p}}{\kk^{p-1}}{\Id}{G}{J}{H}
		\end{gqrep} 
		&
		\begin{array}{c}
			\{({\mathbf{u}}^n)_{n=1}^{p-1}, ({\mathbf{v}}^n)_{n=1}^{p}, ({\mathbf{w}}^n)_{n=0}^{p-1}\} 
			\\
			{\mathbf{v}}^p
		\end{array}
		\\
		\hline
		\begin{array}{c}
			(\hat{\star},p,-)
			\\
			p \in \N
		\end{array}
		&
		\begin{cd}
			I_\star \ar{r}{t^{p}}
			\& P_{-}
		\end{cd}
		&
		\begin{gqrep}{\kk^{p-1}}{\kk^{p}}{\kk^{p}}{G}{\Id}{H}{J}
		\end{gqrep} 
		&
		\begin{array}{c}
			\{({\mathbf{u}}^n)_{n=0}^{p-1}, ({\mathbf{v}}^n)_{n=1}^{p}, ({\mathbf{w}}^n)_{n=1}^{p-1}\}
			\\
			{\mathbf{v}}^p
		\end{array}
		\\
		\hline
		\begin{array}{c}
			(\hat{\star},p,\star)\\
			p \in \N_0
		\end{array}
		&
		\begin{cd}
			I_\star \ar{r}{t^{p}}
			\& P_{\star}
		\end{cd}
		& 
		\begin{gqrep}{\kk^{p}}{\kk^{p+1}}{\kk^{p}}{G}{G}{H}{H}
		\end{gqrep} 
		&
		\begin{array}{c}
			\{({\mathbf{u}}^n)_{n=0}^{p-1}, ({\mathbf{v}}^n)_{n=0}^{p}, ({\mathbf{w}}^n)_{n=0}^{p-1}\}
			\\
			{\mathbf{v}}^p
		\end{array}
		\\
		\hline
		\multicolumn{4}{|l|}{}
		\\[-10pt]
		\multicolumn{4}{|l|}{
			\text{with matrices
				$G = \begin{pNiceArray}[small]{cccc}
					0 & 0 & \Cdots & 0 \\
					\hline
					1&&& \\
					&1&&\\
					&&\Ddots&\\
					&&&1\\
				\end{pNiceArray}$,
				$H = \begin{pNiceArray}[small]{c|cccc}
					0 & 1&&& \\
					0& &1 &&\\
					\Vdots&&&\Ddots&\\
					0 &&&&1\\
				\end{pNiceArray}$  and
				$J =
				{
					\begin{pNiceArray}{ccccc}[small]
						0 & & &  &\\
						1 & 0 & & & \\
						&  & \Ddots & & \\
						& & \Ddots & \Ddots & \\
						&	& & 1 & 0 \\
				\end{pNiceArray}}$.
		}}
		\\
		\hline
	\end{array}
	$
\end{table}

\end{document}